%% file: main.tex
\documentclass[11pt,reqno]{amsart}

\usepackage{
  amsmath,
  amsfonts,
  amssymb,
  amsthm,
  amscd,
  graphicx,
  mathtools,
  etoolbox,
  enumitem,
  mathdots,
  fancybox,
  bbm,
  stackrel,
  mathrsfs
}
\usepackage[all]{xy}
\usepackage[dvipsnames]{xcolor}
\usepackage[colorlinks=true, linkcolor=black, citecolor=black, urlcolor=blue, breaklinks=true]{hyperref}

\def\BE#1{\textcolor[rgb]{.75,0.00,0.00}{[BE: #1]}}
\def\PT#1{\textcolor[rgb]{0, .5, .25}{[PT: #1]}}
\def\PTv2#1{\textcolor[rgb]{0, 0.5, 0.25}{}}

\makeatletter
\@namedef{subjclassname@2020}{%
  \textup{2020} Mathematics Subject Classification}
\makeatother

\usepackage{tikz}
\usetikzlibrary{decorations.markings,decorations.pathmorphing,decorations.pathreplacing}
\usetikzlibrary{calc}
\usetikzlibrary{cd}
\usetikzlibrary{patterns}
\usetikzlibrary{shapes}

\usepackage[dvips]{epsfig}
\usepackage{pinlabel}
\newcommand{\ig}[2]{\vcenter{\xy (0,0)*{\includegraphics[scale=#1]{fig/#2}} \endxy}}

\usepackage{zref-clever}
\zcRefTypeSetup{definition}{
  Name-sg = Definition,
  name-sg = Definition,
  Name-pl = Definitions,
  name-pl = Definitions
}

\zcRefTypeSetup{property}{
  Name-sg = Property,
  name-sg = Property,
  Name-pl = Properties,
  name-pl = Properties
}

\zcRefTypeSetup{example}{
  Name-sg = Example,
  name-sg = Example,
  Name-pl = Examples,
  name-pl = Examples
}

\zcRefTypeSetup{lemma}{
  Name-sg = Lemma,
  name-sg = Lemma,
  Name-pl = Lemmas,
  name-pl = Lemmas
}

\zcRefTypeSetup{conjecture}{
  Name-sg = Conjecture,
  name-sg = Conjecture,
  Name-pl = Conjectures,
  name-pl = Conjectures
}

\zcRefTypeSetup{corollary}{
  Name-sg = Corollary,
  name-sg = Corollary,
  Name-pl = Corollaries,
  name-pl = Corollaries
}

\zcRefTypeSetup{theorem}{
  Name-sg = Theorem,
  name-sg = Theorem,
  Name-pl = Theorems,
  name-pl = Theorems
}

\zcRefTypeSetup{remark}{
  Name-sg = Remark,
  name-sg = Remark,
  Name-pl = Remarks,
  name-pl = Remarks
}

\zcRefTypeSetup{claim}{
  Name-sg = Claim,
  name-sg = Claim,
  Name-pl = Claims,
  name-pl = Claims
}

\zcRefTypeSetup{prop}{
  Name-sg = Proposition,
  name-sg = Proposition,
  Name-pl = Propositions,
  name-pl = Propositions
}

\zcRefTypeSetup{defn}{
  Name-sg = Definition,
  name-sg = Definition,
  Name-pl = Definitions,
  name-pl = Definitions
}

\zcRefTypeSetup{equation}{
  Name-sg = ,
  name-sg = ,
  Name-pl = ,
  name-pl =
}

\zcRefTypeSetup{enumi}{
  Name-sg = ,
  name-sg = ,
  Name-pl = ,
  name-pl =
}

\newcommand{\cref}[1]{\zcref{#1}}
\newcommand{\Cref}[1]{\zcref[S]{#1}}

\newcommand\Z{\mathbb{Z}}
\newcommand\Q{\mathbb{Q}}
\newcommand\N{\mathbb{N}}
\DeclareMathOperator{\Aut}{Aut}
\DeclareMathOperator{\End}{End}

\newtheorem{theorem}{Theorem}[section]
\newtheorem{lemma}[theorem]{Lemma}
\newtheorem*{lemma*}{Lemma}
\newtheorem*{prop*}{Proposition} 
\newtheorem{corollary}[theorem]{Corollary}
\newtheorem{conjecture}[theorem]{Conjecture}
\theoremstyle{definition}
\newtheorem{definition}[theorem]{Definition}
\newtheorem{property}[theorem]{Property}
\newtheorem{remark}[theorem]{Remark}

\newtheorem{example}[theorem]{Example}

\newtheorem{defn}[theorem]{Definition} 
\numberwithin{equation}{section}
\allowdisplaybreaks

\AddToHook{env/lemma/begin}{%
  \zcsetup{countertype={theorem=lemma}}}

\AddToHook{env/corollary/begin}{%
  \zcsetup{countertype={theorem=corollary}}}

\AddToHook{env/conjecture/begin}{%
  \zcsetup{countertype={theorem=conjecture}}}

\AddToHook{env/definition/begin}{%
  \zcsetup{countertype={theorem=definition}}}

\AddToHook{env/property/begin}{%
  \zcsetup{countertype={theorem=property}}}

\AddToHook{env/remark/begin}{%
  \zcsetup{countertype={theorem=remark}}}

\AddToHook{env/claim/begin}{%
  \zcsetup{countertype={theorem=claim}}}

\AddToHook{env/example/begin}{%
  \zcsetup{countertype={theorem=example}}}

\AddToHook{env/prop/begin}{%
  \zcsetup{countertype={theorem=prop}}}

\AddToHook{env/defn/begin}{%
  \zcsetup{countertype={theorem=defn}}}

\usetikzlibrary{decorations.markings, arrows.meta, fit} 
\usepackage{svg}
\usepackage{scalerel,stackengine} 
\stackMath
\newcommand\reallywidehat[1]{%
\savestack{\tmpbox}{\stretchto{%
  \scaleto{%
    \scalerel*[\widthof{\ensuremath{#1}}]{\kern.1pt\mathchar"0362\kern.1pt}%
    {\rule{0ex}{\textheight}}
  }{\textheight}%
}{2.4ex}}%
\stackon[-6.9pt]{#1}{\tmpbox}%
}

\DeclareMathAlphabet{\mathcal}{OMS}{cmsy}{m}{n}
\DeclareMathOperator{\Ima}{Im}

\DeclareMathOperator{\rk}{rk}
\DeclareMathOperator{\Hom}{Hom}
\DeclareMathOperator{\grk}{grk}
\DeclareMathOperator{\Sym}{Sym}
\DeclareMathOperator{\rt}{root}
\DeclareMathOperator{\wt}{weight}
\DeclareMathOperator{\codim}{codim}

\DeclareMathOperator{\Id}{Id}
\DeclareMathOperator{\Ind}{Ind}
\DeclareMathOperator{\Res}{Res}
\DeclareMathOperator{\std}{std}
\DeclareMathOperator{\perm}{perm}
\DeclareMathOperator{\fin}{fin}
\DeclareMathOperator{\KM}{KM}
\DeclareMathOperator{\Rep}{Rep}

\DeclareMathOperator{\Webs}{Webs}
\DeclareMathOperator{\Semis}{Semis}

\DeclareMathOperator{\ch}{ch}

\DeclareMathOperator{\ext}{ext}
\DeclareMathOperator{\sph}{sph}
\DeclareMathOperator{\Abe}{Abe}
\DeclareMathOperator{\Span}{span}

\DeclareMathOperator{\rex}{rex}
\DeclareMathOperator{\GL}{GL}
\DeclareMathOperator{\PGL}{PGL}
\DeclareMathOperator{\SL}{SL}
\DeclareMathOperator{\op}{op}
\DeclareMathOperator{\pre}{pre}

\DeclareMathOperator{\rotation}{rot}
\DeclareMathOperator{\reversal}{rev}

\DeclareMathOperator{\start}{start}
\DeclareMathOperator{\finish}{end}
\DeclareMathOperator{\ELL}{ELL}
\DeclareMathOperator{\id}{id}
\newcommand{\expr}{\leftrightharpoons}
\newcommand{\ot}{\otimes}
\DeclareMathOperator{\rot}{sh}

\newcommand{\ula}{\underline{\lambda}}
\newcommand{\comm}[1]{} 
\newcommand{\kk}{\Bbbk}
\newcommand{\al}{\alpha}
\newcommand{\dd}{\partial}
\newcommand{\LL}{\mathcal{L}}
\newcommand{\SBim}{\mathbb{S}\text{Bim}}
\newcommand{\qbinom}{\genfrac{[}{]}{0pt}{}}

\newcommand{\hh}[1]{\widehat{#1}}

\newcommand{\SSBim}{\mathcal{S}\SBim}
\newcommand{\BSBim}{\mathbb{BS}\text{Bim}}
\newcommand{\SBSBim}{\mathcal{S}\mathbb{BS}\text{Bim}}
\newcommand{\DiagSBS}{\mathcal{D}}
\DeclareMathOperator{\BS}{BS}
\DeclareMathOperator{\bsh}{bs}
\newcommand{\mSSBim}{m\SSBim}
\newcommand{\mSBSBim}{m\SBSBim}
\newcommand{\dSBSBim}{\mc{D}}

\newcommand{\pdmSBSBim}{\tilde{\mc{D}}}
\newcommand{\dmSBSBim}{\mc{D}}
\newcommand{\C}{\mathbb{C}}
\newcommand{\HH}{\mathcal{H}}
\newcommand{\KL}{\underline{H}}
\newcommand{\triv}{\text{triv}}

\newcommand{\SBimod}[2]{{}_{#1}\mathcal{S}_{#2}}
\newcommand{\Hecke}[2]{{}^{#1}\HH^{#2}}
\newcommand{\HAB}[3]{{}^{#1}H^{#2}_{#3}}
\newcommand{\KLB}[3]{{}^{#1}\underline{H}^{#2}_{#3}}
\newcommand{\Bimod}[3]{{}_{#1}{#3}_{#2}}
\newcommand{\bModgr}[2]{\text{${#1}$-Mod$^{\text{gr}}$-${#2}$}}

\newcommand{\cwebs}{\Webs^{\Om}}
\newcommand{\mc}[1]{\mathcal{#1}}
\newcommand{\un}[1]{\underline{#1}}
\newcommand{\set}[1]{\{#1\}}

\graphicspath{{arxiv-figures/}}

\newcommand{\gr}[1]{\scalebox{0.64}{\color{gray}$#1$ \color{black}}} 
\newcommand{\dgr}[1]{\scalebox{0.6}{\color{gray}$\reallywidehat{#1}$\color{black}}} 
\newcommand{\blk}[1]{ \scalebox{0.75}{$#1$} } 
\newcommand{\blkcoeff}[1]{ \scalebox{0.85}{$#1$} } 
\newcommand{\mapstosize}[1]{\scalebox{#1}{$\mapsto$}}
\newcommand{\ddm}[2]{\dd^{\hh{#1}}_{\hh{#2}}}
\newcommand{\Del}[3]{\Delta^{\hh{#1}}_{\hh{#2}\;(#3)}}
\newcommand{\qRep}{\Rep_q^{\Om}}
\newcommand{\qFund}{\Fund_q^{\Om}}
\newcommand{\QGS}{\mathbb{GS}_\zeta}
\newcommand{\Rev}{\mathrm{R}}
\newcommand{\Dual}{\mathrm{D}}
\DeclareMathOperator{\scl}{scl}
\newcommand{\scalewebs}{\scl} 

\makeatletter
\newcommand{\superimpose}[2]{{%
  \ooalign{%
    \hfil$\m@th#1\@firstoftwo#2$\hfil\cr
    \hfil$\m@th#1\@secondoftwo#2$\hfil\cr
  }%
}}
\makeatother

\newcommand{\ddq}{\mathbin{\mathpalette\superimpose{{\dd}{\cdot}}}}
\newcommand{\etadot}{\dot{\eta}}
\newcommand{\procopdtq}{\mathbin{\mathpalette\superimpose{{\mu}{\cdot}}}}

\newcommand{\qdiag}{\mc{D}_q}
\newcommand{\zdiag}{\mc{D}_\zeta}
\newcommand{\HHH}{\mathbb{H}} 
\newcommand{\HHHH}{\widetilde{\mathbb{H}}} 
\newcommand{\HAbd}{\mathscr{H}}
\newcommand{\KRepom}{\mc{K}_{\Om}} 

\DeclareRobustCommand{\rchi}{{\mathpalette\irchi\relax}}
\newcommand{\irchi}[2]{\raisebox{\depth}{$#1\chi$}}
\newcommand{\inv}{\iota}
\def\Om{\Omega}
\newcommand{\AC}{\mathcal{A}}
\newcommand{\KC}{\mathcal{K}}
\newcommand{\OC}{\mathcal{O}}
\newcommand{\PC}{\mathcal{P}} 
\DeclareMathOperator{\aff}{aff}
\DeclareMathOperator{\GS}{GS}
\DeclareMathOperator{\Perv}{Perv}
\def\ZZ{\mathbb{Z}}
\def\CC{\mathbb{C}}

\newcommand{\bMod}[2]{\text{${#1}$-Mod-${#2}$}}

\DeclareMathOperator{\Fund}{Fund}
\newcommand{\gf}{\mathfrak{g}}
\newcommand{\slf}{\mathfrak{sl}}
\DeclareMathOperator{\Cat}{Cat}
\DeclareMathOperator{\gbimod}{gbimod}

\usepackage{array}
\newcolumntype{O}{>{\centering\arraybackslash\(}p{1in}<{\)}}
\newcolumntype{T}{>{\centering\arraybackslash\(}p{2in}<{\)}}
\newcolumntype{F}{>{\centering\arraybackslash\(}p{4in}<{\)}}
\tikzset {
  ->-/.style = {
      decoration={markings, mark=at position 0.5 with {\arrow{>}}},
      postaction={decorate}
    },
  ->--/.style = {
      decoration={markings, mark=at position 0.2 with {\arrow{>}}},
      postaction={decorate}
    },
  -<-/.style = {
      decoration={markings, mark=at position 0.5 with {\arrow{<}}},
      postaction={decorate}
    },
}

\title{Quantum Satake in type $A$: the general and generic case}

\author[]{Ben Elias}
\address{University of Oregon.}
\email{belias@uoregon.edu}

\author[]{Koppara Philip Thomas}
\address{University of Oregon.}
\email{pthomas@uoregon.edu}
\date{}

\begin{document}
\begin{abstract}
    In work of the first author, the geometric Satake equivalence was (non-rigorously) reinterpreted as an equivalence between two algebraically-defined $2$-categories, one built from representations of a Lie algebra, and one built using singular Soergel bimodules. It was then explained how to $q$-deform this equivalence in type $A$, replacing the special linear Lie algebra with its quantum group, and using singular Soergel bimodules for a deformed reflection representation. Both the algebraic reformulation and its $q$-deformation were only proven in types $A_1$ and $A_2$. In this paper, we prove the result in type $A_{n-1}$ for $n \ge 4$, while working generically: more precisely, we work over a field of characteristic zero, where $q$ is not a root of unity, and having adjoined an $n$-th root of $q$. Along the way we generalize certain results in the literature (e.g. the Soergel-Williamson categorification theorem, the Soergel conjecture for spherical elements, various symmetries) to the deformed reflection representation.
\end{abstract}
\maketitle
\input{Introduction}

\input{Combo}

\section{SSBim from Affine Type A Realization}
\input{SSBim_from_Affine_Type_A_Realization}

\section{Frobenius Data and Calculations with Demazure Operators}
\input{Frobenius_Data_and_Calculations_with_Demazure_Operators}
\section{The Diagrammatic Quantum Satake Functor}
\input{The_Diagrammatic_Quantum_Satake_Functor}

\section{Checking well-defined-ness of the functor}\label{sec well-defined-ness of the functor}

\input{Well-Definedness_of_the_Functor}
\section{Supplementary}\label{sec supplementary}

\input{Supplementary_Material}

\bibliographystyle{amsalpha} 
\bibliography{references}

\end{document}

%% file: Introduction.tex
\section{Introduction}
\subsection{Geometric Satake and numerical coincidences}

 In \cite{EQuantumI}, the geometric Satake equivalence \cite{LusztigGS, Ginz95, MirkVil} was reinterpreted algebraically. This reinterpretation makes sense in any type, but was only made explicit in type $A$, where we now focus our attention.

On one side of the equivalence is a 2-category built from the $\Q$-linear representation theory of $\mathfrak{sl}_n$, which we call $\Rep^{\Omega} \slf_n$. 
The objects of this $2$-category form a group $\Omega \cong \Z/n\Z$, which is isomorphic to the quotient of the weight lattice by the root lattice; the $1$-morphisms are representations whose weights live in the appropriate coset within $\Omega$, and $1$-morphism composition is tensor product of representations. Loosely speaking, $\Rep^{\Omega} \slf_n$ is just the monoidal category $\Rep \slf_n$ with additional bookkeeping.

On the other side of the equivalence is a 2-category of singular Soergel bimodules \cite{WillSingular} in type $\tilde{A}_{n-1}$, which can be thought of as an algebraic replacement for equivariant perverse sheaves on the affine Grassmannian. 
Crucial to note is that the vertices of the affine Dynkin diagram $\tilde{A}_{n-1}$ are in bijection with $\Omega$, and these vertices also parametrize copies of the finite Weyl group $S_n$ inside the affine Weyl group $W_{\aff}$ as (standard) parabolic subgroups. 
These finite parabolic subgroups are objects in the 2-category of singular Soergel bimodules $\SSBim$. The $1$-morphisms are particular bimodules over certain polynomial rings.

The geometric Satake equivalence becomes a (conjectural) $\Q$-linear 2-functor $\GS$ from $\Rep^{\Omega} \slf_n$ to $\SSBim$. 
Note that singular Soergel bimodules are a graded category whereas representations are not, so one cannot expect $\GS$ to be an equivalence. Instead it is a \emph{degree zero equivalence}: it is fully faithful onto the space of degree zero 2-morphisms, and there are no 2-morphisms of negative degree between the 1-morphisms in its image.
We call this the \emph{Soergelified (geometric) Satake equivalence}. We refer to it as conjectural for two reasons. First, the reformulation from traditional geometric Satake to Soergelified Satake was sketched but not made rigorous in \cite{EQuantumI}. Second, even were it made rigorous, geometric Satake would imply the existence of a $2$-functor $\GS$, but it would not be obvious how to make this $2$-functor explicit on the level of $2$-morphisms.

In \cite{EQuantumI} an explicit (potential) construction of the $2$-functor $\GS$ is given for all $n \ge 2$, and it is proven correct for $n=2, 3$. The construction uses webs, which are a generators-and-relations description of the monoidal category $\Rep \slf_n$. For $n=2$ webs form the well-known Temperley-Lieb category; for $n=3$ they were constructed by Kuperberg \cite{Kupe}; for $n>3$ they were constructed by Cautis-Kamnitzer-Morrison \cite{CKM} based on earlier work of Morrison. Concretely, \cite{EQuantumI} explains where the generating webs should go under $\GS$, but only checks the relations for $n=2,3$. In this paper we check the relations for $n > 3$.

Amongst the relations defining the web category, various signed binomial coefficients appear. For example, the value of a $1$-labeled circle in the web category is $(-1)^{n-1} \binom{n}{1} = (-1)^{n-1} n$, which is the categorical dimension of the standard representation of $\slf_n$. Meanwhile, the construction of singular Soergel bimodules takes as input the reflection representation $V$ of $W_{\aff}$, equipped with a choice of simple roots and coroots. 
The technical heart of Soergelified Satake is the realization that these binomial coefficients are also coded within $V$. The simplest example of this phenomenon is the fact that $-2$, the value of the Temperley-Lieb circle, agrees with the off-diagonal entry in the Cartan matrix for $\tilde{A}_1$.

The representation $V$ is spanned by the simple roots $\{\alpha_i\}_{i \in \Omega}$, and the action of $W_{\aff}$ is governed by the affine Cartan matrix. Let $R = R_{\aff}$ be the polynomial ring whose linear terms are $V$. For each simple reflection $s_i$, $i \in \Omega$, one defines Demazure operators $\partial_i \colon R \to R$,
\begin{equation} \partial_i(f) = \frac{f - s_i(f)}{\alpha_i}. \end{equation}
For example, the entries of the Cartan matrix are $\partial_i(\alpha_j)$; the number $-2$ from the previous paragraph is $\partial_1(\alpha_0)$ when $n=2$. More generally, one can realize binomial coefficients by applying certain sequences of Demazure operators to certain products of positive roots. For example, when $n=4$, one has
\begin{equation} \label{4intro} \partial_1 \partial_2 \partial_3 \left(\alpha_0 \cdot s_1(\alpha_0) \cdot s_2 s_1(\alpha_0) \right) = -4 = (-1)^{4-1} \binom{4}{1}. \end{equation}
The reader can guess the generalization of \eqref{4intro} to any $n \ge 2$. Similarly, 
\begin{equation} \label{42intro} \partial_2 \partial_3 \partial_1 \partial_2(\alpha_0 \cdot s_1(\alpha_0) \cdot s_3(\alpha_0) \cdot s_1 s_3(\alpha_0)) = \binom{4}{2}.
\end{equation}
To our knowledge, these observations and related numerical coincidences have not yet appeared in the literature.

The computations \eqref{4intro} and \eqref{42intro} arise when checking that the circle relations of the web category continue to hold in $\SSBim$ after applying $\GS$. The product of roots appears via a coproduct for a Frobenius extension between subrings of $R$, and the sequence of Demazure operators is the trace map for a Frobenius extension between different subrings of $R$. To check the remaining web relations requires more computations of a similar flavor: applying Demazure operators to tensors (of polynomials) built from Frobenius coproducts. These computations can be quite thorny.

Once one has constructed the functor $\GS$, one can prove the rest of the Soergelified Satake equivalence from scratch (i.e. without going through geometric Satake). It remains to show that $\GS$ is a degree zero equivalence, and here is a brief sketch. Williamson's Hom formula (\cref{thm-hom formula}) shows that the graded dimension of $2$-morphism spaces between singular Soergel bimodules is determined by a certain pairing on its Grothendieck group. The seminal work of Lusztig \cite{LusztigGS} determines the images of (characters of) irreducible representations under the Satake isomorphism (i.e. the decategorification of the geometric Satake equivalence).  These results combine to show that the dimensions of the $2$-morphism spaces are appropriate for $\GS$ to be a degree zero equivalence. In particular, if $\GS$ is faithful then it is also full (to degree zero) by a dimension count.

If we define our categories to be $\Q$-linear, then we can continue as follows. The semisimple 2-category $\Rep^{\Omega} \slf_n$ has no nontrivial monoidal ideals (see \cref{lemma monoidal ideals in qRep} for a proof, but this is an adaptation of a well-known argument), so the kernel of $\GS$ is either zero or everything. The latter is absurd, since identity maps go to identity maps. Hence $\GS$ is faithful. This proof will not work integrally or over fields of finite characteristic, where nontrivial monoidal ideals exist; we will have more to say about this problem later.


\subsection{The quantum deformation}

One side of the Satake equivalence has a well-known $q$-deformation: the representation theory of the quantum group $U_q(\slf_n)$. These can also be described with webs, and quantum binomial coefficients replace ordinary binomial coefficients. Interestingly, \cite{EQuantumI} introduced a $q$-deformed Cartan matrix of type $\tilde{A}_{n-1}$, leading to a deformations $V_q$ and $R_q$, and to a deformation of singular Soergel bimodules which currently has no analogue in geometry. The \emph{quantum (Soergelified) Satake equivalence} is the corresponding (conjectural) $\C(q)$-linear degree zero equivalence $\GS_q$ from $\Rep^{\Omega} U_q(\slf_n)$ to $q$-deformed singular Soergel bimodules. Here is the deformed Cartan matrix in type $\tilde{A}_3$.
\begin{equation} \label{eq:deformedq} \begin{bmatrix}
    2 &-1 &0 & -q^{-1} \\
    -1 &2 &-1 &0\\
    0 &-1 &2 &-q \\
    -q & 0&-q^{-1} &2
\end{bmatrix} \end{equation}
Repeating the computation of \eqref{4intro} for this Cartan matrix, one obtains
\begin{equation} \label{q4intro} \partial_1 \partial_2 \partial_3 \left(\alpha_0 \cdot s_1(\alpha_0) \cdot s_2 s_1(\alpha_0) \right) = -[4]_q. \end{equation}

Once again, \cite{EQuantumI} explicitly constructed $\GS_q$ for $n=2,3$ and proved it was well-defined. The construction of $\GS_q$ for $n \ge 4$ as a deformation of $\GS$ entails a choice of invertible scalar for each generating web, such that the web relations hold after applying $\GS_q$. In this paper we provide those scalars and check the relations, thus constructing $\GS_q$ for $n \ge 4$. We also prove that this $2$-functor is a degree zero equivalence over $\Q(q)$, finishing the work begun in \cite{EQuantumI}. The details can be found in \cref{thm main}. More accurately, we work over an extension of $\Q(q)$, see \S\ref{subsec:covintro}.

The technical heart of our proof is a number of thorny computations with Demazure operators acting on $R_q$. This requires even more care than one might expect, because the $q$-deformed Demazure operators no longer satisfy the braid relations, so one must keep careful track of Frobenius structures and transformations between reduced expressions in addition to the other algebraic work.

Our proof that $\GS_q$ is a degree zero equivalence over $\Q(q)$ follows the same outline as for $\GS$ in the previous section. However, Williamson did not prove his Hom formula in this level of generality. The representation $V_q$ is not \emph{balanced} (see \S\ref{subsec-balanced}), a technical condition without which one cannot pick out positive versus negative roots consistently. The Hom formula is one of several basic facts about Soergel bimodules which had not been proven for the $q$-deformation, a situation we rectify in this paper. This also allows us to use certain shortcuts: knowing the dimensions of morphism spaces allows us to bypass certain computations. 

Since our goal is to establish the quantum geometric Satake equivalence rigorously, and since the literature has a tendency to wave its hands at some technicalities, we will be concrete and fill in many details. So that later work can more easily establish the integral or finite characteristic versions of (quantum) Soergel Satake without redoing fundamental work, we prove our results over general base rings when possible. For example, we introduce reflection faithfulness (another technical condition) over general base rings, and show that the $q$-deformed reflection representation is reflection faithful. Along similar lines, the literature on geometric Satake tends to use the extended affine Weyl group, whereas we prefer using different parabolic subgroups in the affine Weyl group. There are some frustrations in matching the choices of convention (coming from right actions vs. left actions\footnote{The commutativity of the spherical Hecke algebra has obscured the effects of some conventional issues.}, positive vs. negative Iwahori subgroups, etcetera), so we redevelop some of the combinatorics from scratch. We provide a description of the combinatorial correspondence between dominant weights and double cosets in a fully explicit fashion, focusing on particular representatives of these cosets which are best adapted to the concept of a ``singular reduced expression'' as introduced in \cite{WillThesis, EKo}.

Another goal of this paper was originally to rigorously establish the relationship between Soergel Satake and geometric Satake (in the non-quantum setting). This is work in progress and should hopefully appear in an appendix to a future version of this paper. The relationship is sketched in \cite[\S 6]{EQuantumI} but no attempt is made at a rigorous proof there.  The reader who is curious about the relationship is encouraged to read the introduction of \cite{EQuantumI}, which instructs the geometrically-minded reader on how to skip most of the paper and efficiently reach \S 6.

\subsection{Remarks on integral forms}\label{subsec:youthfuloptimism}

In \cite{EQuantumI} it was promised that the general case of quantum Soergel Satake would be proven in a sequel. The twelve year time lapse between \cite{EQuantumI} and its sequel (the present paper) deserves some explanation. In fact, what \cite{EQuantumI} proposed is a stronger result about the integral forms of these categories, which we discuss now. We do not prove this stronger result.

\begin{remark} The first author apologizes for their youthful optimism; the intervening years have cured both these flaws, youthfulness and optimism. At the time of \cite{EQuantumI} the basic computations like \eqref{q4intro} had been done (and \eqref
{42intro} for $n=4$), but not the nastier ones. \end{remark}

By presenting a category by generators and relations, one makes it straightforward to define functors therefrom. Another aspect of a presentation is that it defines an \emph{integral form} of the category. The category one is presenting may originally be $\Q$-linear (or $\Q(q)$-linear), but if all the coefficients in all the relations live in a subring $A$ (like $\Z$ or $\Z[q,q^{-1}]$), then one can define an $A$-linear category by generators and relations instead. We think of an integral form as being ``correct'' if morphism spaces are free over $A$, and specialize in an interesting way (e.g. to finite characteristic, or setting $q$ to be a root of unity). It was proven in \cite{ELLCC} that Cautis-Kamnitzer-Morrison's webs are the correct integral form for $\Rep^{\Omega}\slf_n$, in that they specialize to describe the category of tilting modules over any field.

Meanwhile, it was already known by the authors of \cite{CKM} that many of their relations were redundant when working non-integrally, over $\Q(q)$. For example, a streamlined presentation can be found in work of Bigelow \cite{Bigelow} (which also fixes a sign error in \cite{CKM}), and we reiterate that Bigelow's presentation is knowingly the ``incorrect'' integral form.

Our construction of $\GS$ and $\GS_q$ factors through a diagrammatic description of $\SSBim$ which we denote $\DiagSBS$. While an integral form for ordinary Soergel bimodules was found in the early 2010s by Elias, Khovanov, and Williamson \cite{EKho, ECathedral, EWGr4sb}, constructing an integral form for singular Soergel bimodules remains a long-standing programme (one which the youthful optimist expected to have been concluded by now). 
A diagrammatic $2$-category was constructed in \cite{EWSFrob}, called $\DiagSBS(\Gamma)$ in the body of the paper, with a 2-functor to $\SSBim$ which is full (proven\footnote{in the case when the realization is balanced.}  in \cite{EKLP}) but known \emph{not} to be faithful; there are not enough relations in \cite{EWSFrob}. Let $\DiagSBS$ be the quotient of $\DiagSBS(\Gamma)$ by the kernel, so that $\DiagSBS$ maps fully faithfully to $\SSBim$. Some of the relations in this kernel are known from the literature (e.g. from \cite{ECathedral, ELosev}), but it is expected that more relations are required, especially to get the ``correct'' integral form. Thus $\DiagSBS$ is not very explicit in its current state, but is nonetheless a useful language to encode morphisms in $\SSBim$.

One of the ``shortcomings'' of this paper is that we do not prove our 2-functor $\GS_q$ is well-defined using only known relations in $\DiagSBS$, but resort to calculations in $\SSBim$. We cannot be certain, but feel there are not enough known relations to do the job, so we await an improvement in the available diagrammatic technology. Since the state of the art constrains us to work with $\SSBim$ and not to work with an integral form thereof, we avail ourselves of certain tricks. We use the Soergel-Williamson Hom formula (which specifies the dimension of certain morphism spaces) to simplify some computations. We avail ourselves of Bigelow's presentation.

We have been referring to the (not-explicitly-defined) $2$-category $\DiagSBS$ as an integral form for $\SSBim$, but this is abuse of terminology. It should be an integral form for the \emph{Hecke $2$-category}, which is supposed to categorify the Hecke algebroid in such a way that the Soergel-Williamson Hom formula holds. It is known that singular Soergel bimodules are \emph{not} well-behaved in some settings (e.g. finite characteristic in affine type) so they are only a good model for the Hecke $2$-category generically, not for its integral form. Abe \cite{Abe} has provided an algebraic model for the Hecke $2$-category, though not a presentation.

When the programme to present the Hecke $2$-category is complete, there should be a version of $\DiagSBS$ defined over $\Z[q,q^{-1}]$, with a presentation, such that morphism spaces are \textbf{free} over $\Z[q,q^{-1}]$ of the expected size. Showing freeness given a presentation is very difficult work! What it would take to upgrade our $2$-functor $\GS_q$ to a $2$-functor from Cautis-Kamnitzer-Morrison's presentation to the integral form of $\DiagSBS$? What remains is less difficult, though there are some subtleties.

We have defined the $2$-functor on objects and morphisms already, and need only check the relations. Most of the relations from \cite{CKM} have the form: take a diagram $X$, and rewrite it as a linear combination of elements in a basis $\mathbb{B}$ for this morphism space. Consider the collection of diagrams $\GS_q(\mathbb{B})$ in $\DiagSBS$. This collection is linearly independent, which can be proven after base change to $\Q(q)$ where it follows from our results here. If $\GS(X)$ could be rewritten as a linear combination of $\GS_q(\mathbb{B})$, then the coefficients must be as desired, again provable after base change to $\Q(q)$. However, it is not obvious that $\GS_q(\mathbb{B})$ is a basis, even though it has the desired size; it might only span a $\Z[q,q^{-1}]$-sublattice. What remains is not a computational problem, but an abstract question: is $\GS_q$ full?

Major advances in the aforementioned programme were recently made by the first author together with Ko, Libedinsky, and Patimo \cite{EKo, KELPDemazure, KELPBruhat, KELPLeibniz, EKLP}. In particular, in \cite{EKLP} a basis for 2-morphisms in $\SSBim$ is constructed combinatorially as the image of certain diagrams in $\DiagSBS(\Gamma)$, called the \emph{double leaves basis}\footnote{It is proven in \cite{EKLP} that the double leaves maps become a basis after applying the 2-functor to $\SSBim$, but this is done under the balanced assumption. More work would be needed to adapt their arguments to the $q$-deformation.}. This is expected to be a basis for the correct integral form of $\DiagSBS$, once that integral form is established. Meanwhile, a combinatorial basis for webs was constructed in \cite{ELLCC}, the \emph{double ladders basis}. One way to prove the fullness of $\GS$ is to prove the following.

\begin{conjecture} \label{conj:double} Under $\GS$ (or $\GS_q$), double ladders are sent to double leaves (up to invertible scalar in $\Z[q,q^{-1}]$). \end{conjecture}

Double ladders are built from elementary light ladders, which live in one-dimensional morphism spaces modulo lower terms. That the image of elementary light ladders agrees with a (non-elementary) light leaf modulo lower terms is thus not surprising, but in the examples we have computed they agree on the nose.

At first glance this conjecture seems straightforward, but it is surprisingly complex. We illustrate with an example, and though we have not yet introduced the diagrammatic categories in question, we expect the reader will get the point. Complete details on this example are found in \cref{subsec:LLappendix}, so pardon the incompleteness here.

A relatively simple form of elementary light ladder is exemplified by the following web. We also show where it goes under $\GS$.
\begin{equation} \label{nastyELL}\vcenter{\xy (0,0)*{\def\svgscale{0.15}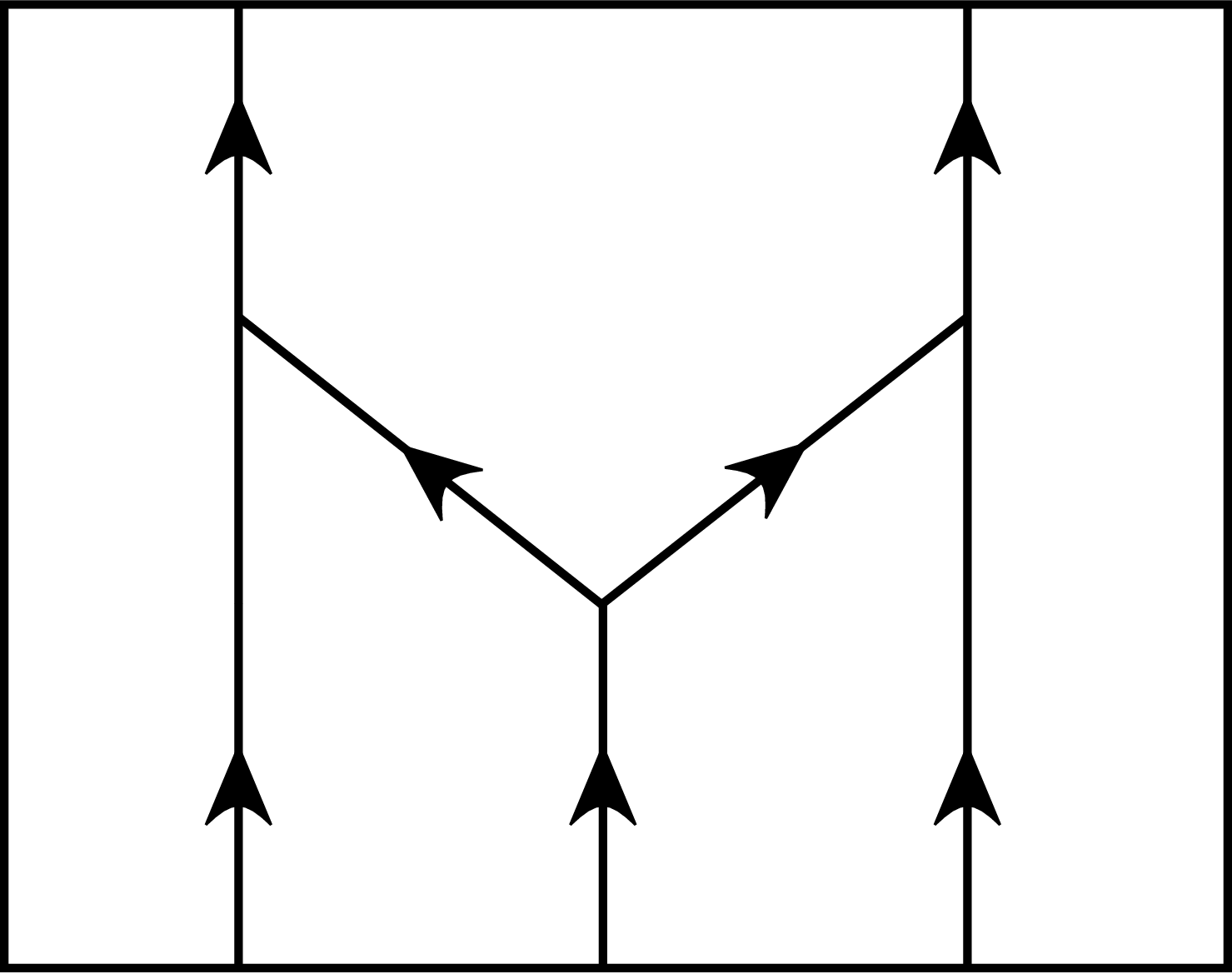} \endxy} \quad\mapstosize{1.5} \quad {
\labellist
\small\hair 2pt
 \pinlabel {$\blk{5}$} [ ] at 8 -5
 \pinlabel {$\blk{6}$} [ ] at 24 -5
 \pinlabel {$\blk{2}$} [ ] at 40 -5
 \pinlabel {$\blk{5}$} [ ] at 56 -5
 \pinlabel {$\blk{0}$} [ ] at 72 -5
 \pinlabel {$\blk{2}$} [ ] at 88 -5
 \pinlabel {$\blk{4}$} [ ] at 8 80
 \pinlabel {$\blk{6}$} [ ] at 24 80
 \pinlabel {$\blk{0}$} [ ] at 72 80
 \pinlabel {$\blk{4}$} [ ] at 88 80
 \pinlabel {$\dgr{0}$} [ ] at 90 45
 \pinlabel {$\dgr{6}$} [ ] at 2 45
\endlabellist
\centering
\ig{1}{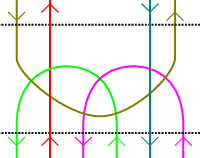}
} \end{equation}
\vspace{.2cm}
Meanwhile, the corresponding light leaf in $\DiagSBS$ is the following diagram.
\begin{equation} \label{nastyLL} \ig{1}{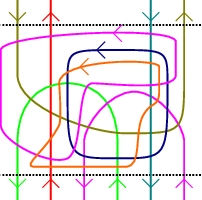} \end{equation}
It turns out that \eqref{nastyELL} and \eqref{nastyLL} are actually equal, though this requires a number of non-obvious diagrammatic transformations (e.g. ``switchback relations'').

\subsection{A change of variables} \label{subsec:covintro}

We wish to note a change of variables that is used throughout the paper. The $q$-deformed Cartan matrix for $\tilde{A}_n$ from \eqref{eq:deformedq} has the nice feature that it contains the usual Cartan matrix for $A_n$ inside, but it admits very little symmetry! Each parabolic subgroup of $\tilde{A}_n$ behaves a little bit differently. In the body of this paper, we do most computations for a different deformation of the usual $\tilde{A}_n$ Cartan matrix over a different formal variable $\zeta$, exemplified here for $n=3$.
\begin{equation} \begin{bmatrix}
    2 &-\zeta &0 &-\zeta^{-1} \\
    -\zeta^{-1} &2 &-\zeta &0 \\
    0 &-\zeta^{-1} &2 &-\zeta \\
    -\zeta & 0 &-\zeta^{-1} &2 
\end{bmatrix}\end{equation}
This Cartan matrix was introduced in \cite{EJY1}. It leads to deformations $V_{\zeta}$ and $R_{\zeta}$. We will also define a $2$-functor $\GS_{\zeta}$ from $\Rep^{\Omega} \slf_n$ to $\SSBim$ as defined over $R_{\zeta}$. 

When $\zeta^n=q^{-2}$, there is an isomorphism $V_q \cong V_{\zeta}$ which intertwines the root data up to scalar. This change of variables requires an extension of $\Q(q)$. The extension $\Q(q^{\frac{1}{n}})$ will suffice, and we use this for the introduction, see \cref{subsubsec quantum numbers} for more details.

The advantage of the $\zeta$-deformation is its obvious rotational symmetry (i.e. rotation of the Dynkin diagram, the action of $\Omega$; not rotation of diagrams), a symmetry we denote $\sigma$ in this paper, and a symmetry shared by colored webs. Diagrammatically, $\sigma$ will just change various labels (depicted as colors) in each diagram, without changing the underlying diagrams or rescaling them. The 2-functor $\GS_{\zeta}$ will preserve $\sigma$, making it easier to check relations.

\begin{remark} One could translate the symmetry $\sigma$ through the isomorphism with $V_q$ to obtain a symmetry on $V_q$ and on categories built from $V_q$. However, this symmetry would not just relabel things but would also multiply them by annoying scalars. \end{remark}

The disadvantage of the $\zeta$ deformation is that it obscures the quantum parameter $q$ which is most common in the literature. Knowing the isomorphism $V_q \cong V_{\zeta}$ is not quite enough to easily convert from one parameter to another. The most natural choice of Frobenius extension structures between $R_q$ and its subrings (all of which can be defined over $\Q(q)$ without needing an extension) differs from the most natural choice for $R_\zeta$. The very meaning of the diagrams we use (in $\DiagSBS$) depends on this choice of Frobenius extension structure. We distinguish the two diagrammatic categories as $\DiagSBS_\zeta$ and $\DiagSBS_q$. In \cref{appendix:comparisonqz} we explicitly describe the way in which one converts between these two conventions.

One might also hope to define $\GS_q$ over $\Q(q)$ rather than its extension $\Q(q^{\frac{1}{n}})$, and this is possible, but we do not attempt to make it explicit, because of the following obnoxious issue.

As noted earlier, diagrams in the web category are sent to scalar multiples of diagrams in $\DiagSBS$. We offer a particular preferred choice of these scalars in \cref{goodchoiceoflambdanu}. What precisely these scalars are is not intrinsic to geometric Satake! One could obtain a new $2$-functor $\GS'$ with different scalars by precomposing with an object-fixing autoequivalence of webs, one which rescales the generating morphisms. However, certain relationships between these scalars are required for the $2$-functor to be well-defined, see \eqref{scalarconditions}. When we define $\GS_{\zeta}$ in \cref{defn dQGS functor}, we assume implicitly that these scalars will be preserved by $\sigma$. Abandoning this assumption will lead to more bookkeeping.

If we naively try to define $\GS_q$ by using our preferred choice of scalars and pass through the equivalence between $\DiagSBS_{\zeta}$ and $\DiagSBS_q$, the resulting scalars will live in $\Q(q^{\frac{1}{n}})$ and not the subring $\Q(q)$. One must choose scalars which are not preserved by $\sigma$ in order to obtain a 2-functor defined over $\Q(q)$ (one must treat each parabolic subgroup a little differently). This is what was done for $n=3$ in \cite{EQuantumI}. Currently, the extra bookkeeping does not seem to be worth the effort, though we may spell this out in a future version of this paper.

\subsection{More symmetries}

Web categories have three ``natural'' symmetries we focus on: rotation of diagrams by 180 degrees, the horizontal flip, and the vertical flip. The $2$-category $\DiagSBS(\Gamma)$ from \cite{EWSFrob} also has these symmetries\footnote{In $\DiagSBS(\Gamma)$, each strand in a diagram has an orientation, but this orientation is cosmetic, a convention to help one keep track of a labeling on regions. After flipping a diagram, convention dictates that one will also need to invert orientations.}. Note that the composition of a horizontal and a vertical flip is the rotation by 180 degrees. The 2-functor $\GS_q$ appears to preserve these symmetries; a flipped web is sent to a flipped diagram (up to sign). However, as noted in \S\ref{subsec:youthfuloptimism}, we are only able to construct a $2$-functor to singular Soergel bimodules $\SSBim$; not to $\DiagSBS(\Gamma)$, nor to an explicit quotient thereof which a priori possesses the same symmetries.

The bimodule $2$-category $\SBSBim$ also has three symmetries in the literature: adjunction, swapping the left and right actions, and (Poincar\'{e}) duality (in the same order as the corresponding diagrammatic symmetries above). Duality is the most complex of these, an exposition for which can be found in \cite[Section 3.4]{BBE}. For an $(R^I, R^J)$-bimodule $B$, its right dual (ignoring grading shifts) would be $\Hom_{(-,R^I)}(B,R^I)$, where this represents the space of morphisms as a right $R^I$-module. Any singular Bott-Samelson bimodule is isomorphic to its right (or left) dual, and that isomorphism is defined by fixing a particular invariant bilinear form on the underlying bimodule. Crucially, these invariant bilinear forms can be defined using same Frobenius extension structures which are used to define adjunction. It is for this reason that ``rotation'' agrees with the composition of the ``horizontal'' and the ``vertical'' flip in the category of bimodules. This connection between bilinear forms and Frobenius structures is sometimes explicit and sometimes implicit in the literature, but what is seemingly absent from the literature is the statement that the 2-functor $\Phi$ from diagrams $\DiagSBS(\Gamma)$ to bimodules intertwines these three symmetries, especially in the unbalanced case. We clarify the situation in \cref{subsec reversal duality and phi}, and state the compatibility between $\GS_q$ and these symmetries in \cref{lemma symmetry of GS-zeta}. 



\subsection{Acknowledgments}

Both authors were supported by the first author's NSF grant DMS-2201387. We appreciate the support given to our research group by NSF RTG grant DMS-2039316. This work forms a substantial portion of the second author's doctoral thesis. The second author is extremely grateful to the first author, who as a doctoral advisor was a source of invaluable guidance, patience, and support. The second author also thanks Victor Ostrik for priceless mathematical discussions and encouragement throughout his doctoral studies.

%% file: arxiv-figures/ELLforBen2_svg-tex.eps_tex
\begingroup%
  \makeatletter%
  \providecommand\color[2][]{%
    \errmessage{(Inkscape) Color is used for the text in Inkscape, but the package 'color.sty' is not loaded}%
    \renewcommand\color[2][]{}%
  }%
  \providecommand\transparent[1]{%
    \errmessage{(Inkscape) Transparency is used (non-zero) for the text in Inkscape, but the package 'transparent.sty' is not loaded}%
    \renewcommand\transparent[1]{}%
  }%
  \providecommand\rotatebox[2]{#2}%
  \newcommand*\fsize{\dimexpr\f@size pt\relax}%
  \newcommand*\lineheight[1]{\fontsize{\fsize}{#1\fsize}\selectfont}%
  \ifx\svgwidth\undefined%
    \setlength{\unitlength}{709.94558601bp}%
    \ifx\svgscale\undefined%
      \relax%
    \else%
      \setlength{\unitlength}{\unitlength * \real{\svgscale}}%
    \fi%
  \else%
    \setlength{\unitlength}{\svgwidth}%
  \fi%
  \global\let\svgwidth\undefined%
  \global\let\svgscale\undefined%
  \makeatother%
  \begin{picture}(1,0.78882)%
    \lineheight{1}%
    \setlength\tabcolsep{0pt}%
    \put(0,0){\includegraphics[width=\unitlength]{ELLforBen2_svg-tex.eps}}%
    \put(0.85338918,0.37430824){\color[rgb]{0,0,0}\makebox(0,0)[lt]{\lineheight{1.25}\smash{\begin{tabular}[t]{l}\gr{0}\end{tabular}}}}%
    \put(0.79852686,0.04517085){\color[rgb]{0,0,0}\makebox(0,0)[lt]{\lineheight{1.25}\smash{\begin{tabular}[t]{l}\blk{2}\end{tabular}}}}%
    \put(0.07077373,0.37492017){\color[rgb]{0,0,0}\makebox(0,0)[lt]{\lineheight{1.25}\smash{\begin{tabular}[t]{l}\gr{6}\end{tabular}}}}%
    \put(0.50167657,0.04468765){\color[rgb]{0,0,0}\makebox(0,0)[lt]{\lineheight{1.25}\smash{\begin{tabular}[t]{l}\blk{3}\end{tabular}}}}%
    \put(0.20378734,0.04270029){\color[rgb]{0,0,0}\makebox(0,0)[lt]{\lineheight{1.25}\smash{\begin{tabular}[t]{l}\blk{1}\end{tabular}}}}%
    \put(0.20858517,0.71733325){\color[rgb]{0,0,0}\makebox(0,0)[lt]{\lineheight{1.25}\smash{\begin{tabular}[t]{l}\blk{2}\end{tabular}}}}%
    \put(0.79978727,0.71773939){\color[rgb]{0,0,0}\makebox(0,0)[lt]{\lineheight{1.25}\smash{\begin{tabular}[t]{l}\blk{4}\end{tabular}}}}%
  \end{picture}%
\endgroup%

%% file: Combo.tex
\section{Combinatorics of affine double cosets}
\label{sec:combo}

\begin{defn}
    Let $(W,S)$ be a Coxeter system. A subset $I \subset S$ is said to be \textit{finitary} if the parabolic subgroup $W_I \subset W$ is finite. In this case we use $w_I$ to denote the longest element of $W_I$, and set $\ell(I):= \ell(w_I).$
\end{defn}

\begin{defn} \label{defn:maximalcoset}
    A subset $I \subset S$ is \emph{maximal finitary} if it is finitary, and maximal in the poset of finitary subsets of $S$. If $I$ and $J$ are maximal, then a double coset $p \in W_I \backslash W / W_J$ is called a \emph{maximal double coset}\footnote{There is a Bruhat (partial) order on the set of double cosets, giving a different meaning to the adjective ``maximal.'' However, in an affine Weyl group, there are no maximal elements under the Bruhat order on double cosets Thus we hope this slightly ambiguous terminology causes no confusion.}.
\end{defn}

\begin{remark} When $W$ is finite, there is a unique maximal finitary subset of $S$, namely $S$ itself, and a unique maximal double coset in $W \backslash W / W$. It is for infinite Coxeter groups where maximal double cosets become interesting.
\end{remark}

The following construction is well-known, see \cite[Chapter 8.3]{BjornerBrenti}.

\begin{theorem} Let $S_{\Z}$ denote the permutations of $\Z$. Fix $n \ge 2$. Consider the subset of $S_{\Z}$ consisting of all permutations $w$ such that
\begin{equation} \label{periodic} w(k+n) = w(k)+n \quad \text{for all } k \in Z, \end{equation}
\begin{equation}\label{netzero} \sum_{i=1}^n w(i) = \sum_{i=1}^n i. \end{equation}
Then this subset is a subgroup isomorphic to the affine Weyl group $(W,S)$ for $S_n$. The isomorphism sends $s_i$ to the permutation which swaps $k$ and $k+1$ whenever $k \equiv i$ modulo $n$. As usual, we identify $S$ with  $\Omega =\Z/n\Z$ and consider indices (for simple reflections) modulo $n$.

Viewing an element $w \in W$ as living in $S_{\Z}$, the usual Coxeter length function can be described by
\begin{equation} \ell(w) = \# \{(i,j) \in \Z \times \Z \mid i < j, w(i) > w(j) \} / \sim, \end{equation} 
where $(i,j) \sim (i+n, j+n)$. Equivalence classes of pairs $(i,j)$ in this set are called \emph{inversions} of $w$, though abusing notation, we sometimes call the pair $(i,j)$ itself an inversion. \end{theorem}

The goal of this section is to understand the bijection(s) between maximal double cosets in $W$ and dominant weights which underlies geometric Satake, and to prove some statements about reduced expressions for these cosets, while making things as transparent as possible in the language of $S_{\Z}$.

Throughout this chapter we fix $n \ge 2$, and let $(W,S)$ be the affine Weyl group as above. For $a,b \in \Omega$ we use notation like $\hh{ab}$ to denote the finitary subset $S \setminus \{s_a, s_b\}$.

\begin{remark} Note that most approaches to geometric Satake in the literature use double cosets for the finite Weyl group within the extended affine Weyl group, rather than the collection of all maximal double cosets in the affine Weyl group, and the bijection between dominant weights and cosets is usually phrased in that language. A dictionary between these approaches is provided in \cref{subsec extended affine to affine}. \end{remark}

\subsection{Distinguished representatives of cosets}

We recall some standard facts about minimal and maximal representatives in cosets.

\begin{lemma} An element $w \in W$ is minimal (in the Bruhat order) within its coset $w W_{\hh{0}}$ if and only if
\begin{equation} w(1) < w(2) < \ldots < w(n).\end{equation}
Equivalently, $w$ has no inversions with a representative $(i,j)$ satisfying $1 \le i < j \le n$.

Similarly, $w$ is maximal in $w W_{\hh{0}}$ if $w(1) > w(2) > \ldots > w(n)$, and $(i,j)$ is an inversion for all $1 \le i < j \le n$.

Similarly, for all $k \in \Z$, $w$ is minimal in $w W_{\hh{k}}$ if
\[ w(k+1) < w(k+2) < \ldots < w(k+n), \]
and maximal if the order is reversed.

An element $w$ is minimal (resp. maximal) in $W_{\hh{0}} w$ if and only if $w^{-1}$ is minimal (resp. maximal) in $w^{-1} W_{\hh{0}}$. \end{lemma}

\begin{proof} This is straightforward, and left to the reader. \end{proof}

\begin{lemma} \label{lem:maxcrit} Each double coset $p \in W_{\hh{k}} \backslash W / W_{\hh{k'}}$ has a unique minimal element $\underline{p}$ (resp. maximal element $\overline{p}$) in the Bruhat order. If  $w \in p$ is minimal in both $w W_{\hh{k'}}$ and $W_{\hh{k}} w$, then $w = \underline{p}$, and similarly if $w$ is maximal in both right and left cosets, then $w = \overline{p}$. \end{lemma}

\begin{proof} This is well-known, see e.g. \cite[Lemma 2.12]{EKo}. \end{proof}

\begin{example} \label{ex:runningw0} Our favorite representative of a double coset $p \in W_{\hh{k}} \backslash W / W_{\hh{0}}$ will be $\overline{p} w_{\hh{0}}$. Note that $\overline{p} w_{\hh{0}}$ will be minimal in its right $W_{\hh{0}}$ coset. Here is a prototypical example of such an element.
\begin{equation} \label{runningw0} \overline{p} w_{\hh{0}} = {
\labellist
\tiny\hair 2pt
 \pinlabel {$1$} [ ] at 120 0
 \pinlabel {$4$} [ ] at 144 45
\endlabellist
\centering
\ig{1.5}{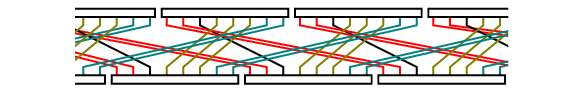}
}\end{equation}
In this example, $n=8$ and $k=3$. Sanity check: This permutation sends $3 \mapsto -2$ and $0 \mapsto 11$.

Here's a picture of the minimal element instead.
\begin{equation} \label{runningmin} \underline{p} = {
\labellist
\tiny\hair 2pt
 \pinlabel {$1$} [ ] at 120 0
 \pinlabel {$4$} [ ] at 144 45
\endlabellist
\centering
\ig{1.5}{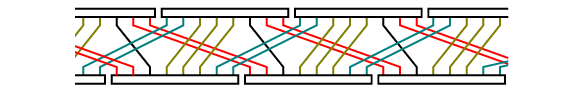}
}\end{equation}

The numbers $c_i$ from the following lemma are
\begin{equation} \vec{c} = (-2,-2,-1,0,0,0,1,1), \end{equation}
and the numbers $a_m$ from \eqref{eq:adef} are
\begin{equation} \vec{a} = (\ldots, 0, 0, 2, 1, \underline{3}, 2, 0, 0, \ldots), \end{equation}
with $a_0 = 3$. Both bits of numerology encode: two red strands, one black strand, three olive strands, and two teal strands. \end{example}

\begin{lemma}\label{lemma dist coset reps} Let $k \in \Z$. Any double coset $p \in W_{\hh{k}} \backslash W / W_{\hh{0}}$ contains a unique $w \in W$ such that
\begin{equation} \label{magicw} w(1) < w(2) < \ldots < w(n), \qquad w(i) \equiv k + i \text{ modulo } n.\end{equation} 
In fact, $w = \overline{p} w_{\hh{0}}$. For $1\leq i \leq n$ we let $c_i \in \Z$ be such that $w(i) = k + i + c_i n$.
Then
\begin{equation} \label{cprops} c_1 \le c_2 \le \ldots \le c_n, \qquad \sum c_i = -k. \end{equation}
\end{lemma}

For the purposes of this chapter we call $w$ the \emph{distinguished representative} of $p$.

\begin{proof} 
Let $y \in p$ be arbitrary. Define $c^y_i \in \Z$ to be the unique integer such that 
\[ k + c^y_i n < y(i) \le k + c^y_i n + n.\] Permutations $x \in W_{\hh{k}}$ will permute each block $\{k + c n + 1, k + cn + 2, \ldots, k + cn + n\}$ for $c \in \Z$, and thus $c^{xy}_i = c^y_i$ for each $i$. Precomposition with $z \in W_{\hh{0}}$ will preserve the image of the set $\{1, \ldots, n\}$, and thus the ordered set $(c^{yz}_i)$ is a permutation of $(c^y_i)$. Thus the multiset $\{c^y_i\}$ is independent of the choice of $y$ within the double coset $p$.

As for any element of $W$, the set $\{y(1), \ldots, y(n)\}$ descends bijectively to $\Z/n\Z$. Thus the set $\{y(i) - (k+c^y_i n)\}$ is equal to the set $\{1, \ldots, n\}$, and
\begin{equation} \sum_{i=1}^n y(i) - k - c^y_i n = \sum_{i=1}^n i.\end{equation}
Subtracting \eqref{netzero} we deduce that
\begin{equation} \sum_{i=1}^n c^y_i n = -kn, \qquad \text{whence } \sum_{i=1}^n c^y_i = -k.\end{equation}

For the rest of this proof we assume $y$ is a minimal length representative of $y W_{\hh{0}}$, so that $y(1) < y(2) < \ldots$. Obviously $c^y_1 \le c^y_2 \le \ldots$, and \eqref{cprops} will hold. Let $c_i = c^y_i$ for such $y$.

We have $y(1) < y(j)$ for all $j \in \Z_{> 1}$. If $y(1) \ne k + c_1 n + 1$, then for any $a$ with $k+c_1 n < a < y(1)$ we have $y^{-1}(a) \in \Z_{\le 0}$. Let $x = s_{k+1} s_{k+2} \cdots s_{y(1) - 1} \in W_{\hh{k}}$, or in cycle notation
\[ x = (k + c_1 n + 1, k + c_1 n + 2, \ldots, y(1)) \]
extended periodically. We have $xy(1) = k + c_1 n + 1$. We claim that $xy$ is also minimal in its right $W_{\hh{0}}$ coset. The inversion sets of $xy$ and of $y$ differ in at most $\ell(x)$ inversions, and for each $a$ as above the pair $(y^{-1}(a),1)$ is an inversion of $xy$ but not of $y$. Thus no other inversions are added (or subtracted), and $xy$ has no inversions with representatives of the form $(i,j)$ with $1 \le i < j \le n$, proving the claim. 

By replacing $y$ with $xy$ if necessary, we assume $y(1) = k + c_1 n + 1$ below.

We have $y(2) < y(j)$ for all $j \in \Z_{> 2}$ with $j \not\equiv 1$ modulo $n$. If $y(2) \ne k + c_2 n + 2$, then for any $a$ with $k+c_2 n + 1 < a < y(2)$ we have $y^{-1}(a) \in \Z_{\le 0}$. Let $x = s_{k+2} s_{k+3} \cdots s_{y(2) - 1} \in W_{\hh{k}}$. Note that $xy(1) = k + c_1 n + 1$ and $xy(2) = k + c_2 n + 2$. As before, one can argue that $xy$ is minimal in $xy W_{\hh{0}}$. By replacing $y$ with $xy$ if necessary, we can assume $y(2) = k + c_2 n + 2$.

Repeating this argument, one eventually finds an element $w \in p$ (equal to the final element $y$) which satisfies \eqref{magicw} and for which \eqref{cprops} holds. We now argue that, for any $w$ satisfying \eqref{magicw}, the element $\gamma = w w_{\hh{0}}$ will be maximal in its double coset. Since double cosets have a unique maximal element, this also proves the unicity of $w$.

Clearly $\gamma(1) > \gamma(2) > \ldots > \gamma(n)$, so $\gamma$ is maximal in its right $W_{\hh{0}}$ coset. Meanwhile, fix any $i, j$ with $1 \le i < j \le n$. Then $w^{-1}(k+i) = i - c_i n$ and $w^{-1}(k+j) = j - c_j n$, with $c_i \le c_j$. Conseqently, 
\[ \gamma^{-1}(k+i) = (n+1-i) - c_i n > (n+1-j) - c_j n = \gamma^{-1}(k+j). \]
Thus each $(k+i,k+j)$ is an inversion of $\gamma^{-1}$, proving $\gamma^{-1}$ is maximal in its right $W_{\hh{k}}$ coset, and $\gamma$ is maximal in $W_{\hh{k}} \gamma$. We conclude with \cref{lem:maxcrit}.
\end{proof}

\begin{corollary} \label{cor:bijection1} Fix $k \in \Z$. Double cosets $p \in W_{\hh{k}} \backslash W / W_{\hh{0}}$ are in bijection with sequences $(c_1, \ldots, c_n)$ of integers satisfying \eqref{cprops}. \end{corollary}

\begin{proof} By the previous lemma, a double coset $p$ has a unique representative $w$ giving rise to the integers $c_i$. Conversely, given integers $c_i$ satisfying \eqref{cprops}, let $w(i) = k+i + c_i n$, and extend this periodically to obtain a permutation of $\Z$. Clearly \eqref{netzero} is satisfied so $w \in W$, and \eqref{magicw} is satisfied so $w$ is the unique representative in its double coset. \end{proof}

This corollary and lemma did not assume that $0 \le k < n$, but applies for all $k \in \Z$. Since $W_{\hh{k}} = W_{\hh{k-n}}$, this gives multiple different ways of encoding a double coset $p$ in a sequence $(c_1, \ldots, c_n)$. In particular, shifting $k \mapsto k-n$ will shift $c_i \mapsto c_i + 1$, thus appropriately shifting $\sum c_i \mapsto \sum c_i + n$ so that \eqref{cprops} still holds. By enforcing that $0 \le k < n$, one obtains a unique sequence $(c_i)$ associated to each maximal double coset. We summarize as follows.

\begin{defn} Let 
\[ C = \{(c_1 \le c_2 \le \ldots \le c_n) \in \Z^n\}\] be the set of increasing sequences, and let $\sim$ be the equivalence relation where $(c_1, \ldots, c_n) \sim (c_1+1, \ldots, c_n+1)$. \end{defn}

\begin{corollary} \label{cor:twobijections} \cref{cor:bijection1} gives bijections between
\begin{equation} \coprod_{k=0}^{n-1} W_{\hh{k}} \backslash W / W_{\hh{0}} \leftrightarrow \{w \in W \mid \; \text{ \eqref{magicw} holds for some }k \} \leftrightarrow C/\sim. \end{equation}
For the last bijection, if $\sum_{i=1}^n c_i = -k$, then $w$ is the permutation sending $i \mapsto k + i + c_i n$, a formula which is constant on equivalence classes for $\sim$.
\end{corollary}

We wish to point out one other useful way to record the same information. Fix $\vec{c}$, and $k \in \Z$ so that \cref{cprops} hold, and let $p$ be the corresponding maximal coset, and $y \in p$ an arbitrary element. We define a sequence $\vec{a} = (a_m)_{m \in \Z}$ with 
\begin{equation}  \label{eq:adef} a_m := \# \{1 \le i \le n \mid c_i = m \} = \# \{0 \le i \le n \mid k-mn < y(i) \le k-mn+n\}.
\end{equation} That is, $a_m$ is the number of elements in the $0$-th ``block'' for $W_{\hh{0}}$ which are sent to the $m$-th ``block'' for $W_{\hh{k}}$, see the picture from \eqref{runningw0}.

\begin{defn} Let 
\[ A = \{\vec{a} = (\ldots, a_{-1}, a_0, a_1, \ldots) \in \Z_{\ge 0}^{\Z} \mid \sum a_m = n\},\] equipped with an equivalence relation $\sim$ which shifts all indices simultaneously. \end{defn}

\begin{lemma} \label{lem:CvsA} There is a bijection $C \to A$ which intertwines the operation $(c_i) \mapsto (c_i + 1)$ with the operation $(a_m)_{m \in \Z} \mapsto (a_{m-1})_{m \in \Z}$, thus intertwining the equivalence relations $\sim$. Note that $\sum_{i=1}^n c_i = \sum_{m \in \Z} m a_m$, a finite sum since all but finitely many values of $a_m$ are zero. \end{lemma}

\begin{proof} The inverse map is easy to describe: one sets $\vec{c}$ to be the concatenation of $a_m$ copies of $m$ (in increasing order). \end{proof}

\begin{example} \label{ex:runningcoset} Let $n=6$ and $\vec{c} = (-4,-3,-3,1,3,4) \in C$. Then 
\[ \vec{a} = (\ldots, 0,0,0,1,2,0,0,\underline{0},1,0,1,1,0,0,\ldots),\]
with the first nonzero entry being $a_{-4}$. If one adds $1$ to each $c_i$ to obtain $(-3, -2, -2, 2, 4, 5)$, the corresponding sequence $\vec{a}$ would look the same, except with first nonzero entry $a_{-3}$.
\end{example}

To summarize, we have a series of natural bijections
\begin{equation} \label{bijections} \coprod_{k=0}^{n-1} W_{\hh{k}} \backslash W / W_{\hh{0}} \leftrightarrow \{w \in W \mid \text{ \eqref{magicw} holds for some } k \} \leftrightarrow (C/\sim) \leftrightarrow (A/\sim) \leftrightarrow \Lambda^+_{\wt}, \end{equation}
with the final bijection appearing in the next section.
We will tacitly use these bijections below, i.e. if we fix $w \in W$ satisfying \eqref{magicw} then we will use the notation $c_i$ or $a_m$ freely, with the understanding that these numbers are well-defined only up to $\sim$, i.e. after fixing $k$.

\begin{lemma} Let $w$ satisfy \eqref{magicw} for some $k$. Then 
\begin{equation} \label{lengthformula} \ell(w) = \sum_{m < m' \in \Z} (m'-m) a_m a_{m'} = \sum_{1 \le i < j \le n} c_j - c_i. \end{equation}
\end{lemma}

\begin{proof} Each inversion has a unique representative $(h,j)$ with $1 \le j \le n$, and each inversion of $w$ must also have $h \le 0$. Let us count the inversions $(h,j)$ such that $1 \le j \le n$ and $c_j = m$ and $h = i - cn \le 0$ and $c_i = m'$.
We claim that there are zero such inversions if $m' \le m$, and $(m'-m) a_m a_m'$ such inversions if $m' > m$, corresponding to the $a_{m'}$ choices of $i$, the $a_m$ choices of $j$, and the choice of $c$ with $0 < c \le m'-m$. This is a relatively straightforward exercise, clarified by the following example. \end{proof}

\begin{example} Continuing \cref{ex:runningw0}, each cabled crossing of the red and olive strands contributes $2 \cdot 3 = a_{-2} \cdot a_0$ crossings. The number of cabled crossings of red and olive strands (up to periodicity) is $2 = 0-(-2)$. Similarly, each red-teal cabled crossing contributes $4 = a_{-2} \cdot a_1$ crossings, and there are $3 = 1 - (-2)$ red-teal cabled crossings (up to periodicity).
\end{example}

\subsection{Maximal cosets and dominant weights}

Let us set notation for weights.

\begin{defn} Let $\Lambda_{\wt}$ denote the weight lattice of $\slf_n$, which we view as the free $\Z$-module spanned by \emph{fundamental weights} $\varpi_i$ for $1 \le i \le n-1$.  The $\N$-linear span of $\{\varpi_i\}$ is the set of \emph{dominant weights}, and is denoted $\Lambda^+_{\wt}$.

We equip $\Lambda_{\wt}$ with a homomorphism $\rot \colon \Lambda_{\wt} \to \Omega$, sending $\varpi_i \mapsto i$.
\end{defn}

\begin{remark} The homomorphism $\rot$ realizes $\Omega$ as the quotient $\Lambda_{\wt}/\Lambda_{\rt}$. \end{remark}

\begin{lemma} \label{lem:cadditive} There is an additive map $\Lambda^+_{\wt} \to C$ sending
\begin{equation} \lambda = \sum_{i=1}^{n-1} g_i \varpi_i \mapsto \vec{c}(\lambda) := (0, g_1, g_1+g_2, \ldots), \qquad c(\lambda)_i := \sum_{j<i} g_j. \end{equation}
There is a map $C \to \Lambda^+_{\wt}$ sending
\begin{equation} \vec{c} = (c_1, \ldots, c_n) \mapsto \sum_{i=1}^{n-1} (c_{i+1} - c_i) \varpi_i. \end{equation}
These maps induce inverse bijections
\[ C/\sim \; \leftrightarrow \; \Lambda^+_{\wt}. \] \end{lemma}

\begin{proof} It is easy to verify that $\vec{c}(\lambda) \mapsto \lambda$, and that $\vec{c}(\lambda + \mu) = \vec{c}(\lambda) + \vec{c}(\mu)$. The image of $\Lambda^+_{\wt}$ consists of elements of $C$ with $c_1 = 0$, of which there is one in each equivalence class. \end{proof}

\begin{example} Continuing \cref{ex:runningcoset}, the coresponding weight is $\lambda = \varpi_1 + 4 \varpi_3 + 2 \varpi_4 + \varpi_5$. Note that $\vec{c}(\lambda) = (0,1,1,5,7,8)$ which is equivalent to $(-4,-3,-3,1,3,4)$. 
\end{example}

The additivity of this function $\Lambda^+_{\wt} \to C$ makes it our preferred representative. This gives rise for each $\lambda$ to a preferred choice of $k = -\sum c_i$, which will be equivalent modulo $n$ to $\rot(\lambda)$. We now state the bijection in these terms.

\begin{theorem} \label{thm:explicitpermutation} Let $\psi_0 \colon \Lambda^+_{\wt} \to \coprod_{k=0}^{n-1} W_{\hh{k}} \backslash W / W_{\hh{0}}$ be the bijection from \eqref{bijections}. Then a dominant weight $\lambda = \sum g_i \varpi_i$ is sent to the coset $\psi_0(\lambda) := p \in W_{\hh{\rot(\lambda)}} \backslash W / W_{\hh{0}}$ for which $w = \overline{p} w_{\hh{0}}$ satisfies
\begin{equation} w(i) = i + \sum_{1 \le j<i} g_j n + \sum_{1 \le j < n} g_j(j-n). \end{equation}
\end{theorem}

\begin{proof} We will define $w$ as in Corollary \ref{cor:bijection1}. For $\vec{c}(\lambda)$ we have $c_i = \sum_{1 \le j < i} g_j$ and
\begin{equation} k = - \sum_{i=1}^{n-1} c_i = - \sum_{i=1}^{n-1} \sum_{j<i} g_j = - \sum_{j=1}^{n-1} g_j (n-j) = \sum_{j=1}^{n-1} g_j(j-n).  \end{equation}
Clearly $k$ agrees with $\rot(\lambda) = \sum g_j j$ modulo $n$.
\end{proof}

We pause to discuss how the bijection behaves on fundamental weights.

\begin{defn}\label{defn cosets corresponding to fundamental weights}
For $0 < k < n$ let $p_{\varpi_k}$ be the $W_{\hh{k}} \backslash W / W_{\hh{0}}$ coset containing the identity. Let $w_{\varpi_k} = w_{\hh{k}} w_{\hh{0,k}}$ be its distinguished coset representative satisfying \eqref{magicw}. Pictured below is the case $n=8$ and $k=3$. 
\begin{equation} \underline{p}_{\varpi_k} = {
\labellist
\tiny\hair 2pt
 \pinlabel {$1$} [ ] at 120 0
 \pinlabel {$4$} [ ] at 144 45
\endlabellist
\centering
\ig{1.0}{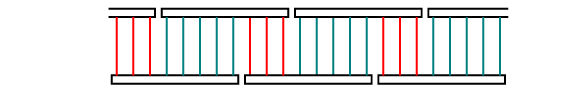}
}\end{equation}
\begin{equation} w_{\varpi_k} = {
\labellist
\tiny\hair 2pt
 \pinlabel {$1$} [ ] at 120 0
 \pinlabel {$4$} [ ] at 144 45
\endlabellist
\centering
\ig{1.0}{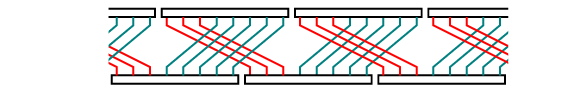}
}\end{equation}
\end{defn}

\begin{lemma}\label{lemma distinguished reps for fundamental weights} Under the bijections of \eqref{bijections}, $\psi_0(\varpi_k) = p_{\varpi_k}$. We have $w_{\varpi_k} = w_{\hh{k}} w_{\hh{0,k}}^{-1}$ and  $\ell(w_{\varpi_k}) = \ell(\hh{k}) - \ell(\hh{0,k}) = k (n-k)$. \end{lemma}

\begin{proof} This is straightforward. We note that $a_{-1} = k$ and $a_0 = n-k$, and $a_i = 0$ otherwise, so the length formula follows from \eqref{lengthformula}. \end{proof}

\begin{remark} Even though $w_I$ is an involution for all $I \subset S$ finitary, we sometimes write $w_I^{-1}$ rather than $w_I$ in compositions, to indicate how multiplication by $w_I$ affects the length of the expression. We could have written $\overline{p} w_{\hat{0}}^{-1}$ when discussing distinguished representatives above. \end{remark}

The following lemma states that addition of dominant weights also adds the lengths of their distinguished coset representatives. We will lift this to a statement about multiplying representatives in \S\ref{subsec:cosetcomp}.

\begin{lemma} \label{lem:lengthadds} Let $\lambda = \sum_{i=1}^{n-1} g_i \varpi_i \in \Lambda^+_{\wt}$, and let $w_{\lambda}$ be the corresponding element of $W$ satisfying \eqref{magicw}. Then
\begin{equation} \label{lengthadds} \ell(w_{\lambda}) = \sum_{i=1}^{n-1} g_i \ell(w_{\varpi_i}) = \sum_{i=1}^{n-1} g_i i(n-i). \end{equation} 
\end{lemma}

\begin{proof} When $\lambda = 0$, $\psi_0(0) \in W_{\hh{0}} \backslash W / W_{\hh{0}}$ is the coset containing the identity, and $w_{\lambda} = 1$ has length zero. This confirms the base case. If the result holds for $\lambda$, we prove it holds for $\lambda + \varpi_k$. Adding $\varpi_k$ will increase $c_i$ by one for all $i > k$. Therefore it will increase the length $\sum_{1 \le i < j \le n} c_j - c_i$ (see \eqref{lengthformula}) by one for each pair $(i,j)$ with $i \le k$ and $j > k$. There are $k(n-k)$ such pairs. \end{proof}

\subsection{Rotation}\label{subsec rotation}

Let $\omega \in S_{\Z} \setminus W$ be defined by $\omega(i) = i+1$. Conjugation by $\omega$ increases all indices in a permutation by one (in one-line notation, in cycle notation, etcetera). We have
\begin{equation} \omega s_i \omega^{-1} = s_{i+1} =: \sigma(s_i). \end{equation}
Conjugation by $\omega$ preserves $W$ and acts as a Dynkin diagram automorphism. We let $\sigma$ denote this automorphism of $W$ and of $\Omega$. Even though $\omega$ has infinite order, $\sigma^n = \id$ since $\omega^n$ centralizes $W$.

For $a \in \Z$, applying $\sigma^a$ to a double coset in $W_{\hh{k}} \backslash W / W_{\hh{0}}$ yields a double coset $p \in W_{\hh{k+a}} \backslash W / W_{\hh{a}}$. For such a double coset $p$, its distinguished representative $\overline{p} w_{\hh{a}}$ satisfies the rotation of \eqref{magicw}, namely
\begin{equation} \label{magicwa} w(a+1) < w(a+2) < \ldots < w(a+n), \qquad w(a+i) \equiv a+ k + i \text{ modulo } n.\end{equation} We let $c_i \in \Z$ be such that $w(a+ i) = a + k + i + c_i n$. This yields new bijections
\begin{equation} \label{bijectionsa} \coprod_{k=0}^{n-1} W_{\hh{k+a}} \backslash W / W_{\hh{a}} \leftrightarrow \{w \in W \mid \; \text{ \eqref{magicwa} holds for some } k \} \leftrightarrow (C/\sim) \leftrightarrow (A/\sim) \leftrightarrow \Lambda^+_{\wt}, \end{equation}
where the last two bijections are the same as before.

\begin{defn}\label{def:psia} For a weight $\lambda \in \Lambda^+_{\wt}$, let $\psi_a(\lambda)$ be the corresponding maximal double coset in $\coprod_{k=0}^{n-1} W_{\hh{k+a}} \backslash W / W_{\hh{a}}$. Equivalently, $\psi_a(\lambda) = \sigma^a(\psi_0(\lambda))$, where $\psi_0$ is the bijection from \cref{thm:explicitpermutation}. \end{defn}

\subsection{Composition} \label{subsec:cosetcomp}

\cref{lem:lengthadds} argues that the length of distinguished coset representatives is additive. In fact, this is because of a much stronger statement, that distinguished coset representatives compose as dominant weights add, and that this composition is reduced (the lengths add). However, when composing coset representatives, as for composing morphisms in a category, we should be careful that the sources and targets match appropriately. Given a coset $p \in W_I \backslash W / W_J$ and a coset $q \in W_J \backslash W / W_K$, we can multiply their distinguished representatives and hope to get a distinguished representative for a coset $p \star q$ in $W_I \backslash W / W_K$. We should not try to compose cosets in $W_{\hh{k}} \backslash W / W_{\hh{0}}$.

\begin{remark} In \cite{EKo} a category is introduced where the objects are finitary parabolic subsets and the morphisms are double cosets. Multiplication of distinguished representatives agrees with composition of cosets in this category only when the composition is reduced. To avoid introducing additional technical overhead from \cite{EKo}, we focus on reduced expressions in the next section. \end{remark}

\begin{theorem} \label{thm:weightsadd} Fix $a \in \Omega$, and $\lambda, \mu \in \Lambda^+_{\wt}$. Let $w_{\mu}$ be the distinguished representative of $\psi_a(\mu)$ (satisfying \eqref{magicwa}), let $w_{\lambda}$ be the distinguished representative of $\psi_{a+\rot(\mu)}(\lambda)$, and let $w_{\lambda + \mu}$ be the distinguished representative of $\psi_a(\lambda + \mu)$. Then
\begin{equation} w_{\lambda + \mu} = w_{\lambda} \cdot w_{\mu}, \qquad \ell(w_{\lambda + \mu}) = \ell(w_{\lambda}) + \ell(w_{\mu}). \end{equation}
\end{theorem}

\begin{proof} The statement about lengths is an immediate consequence of \cref{lem:lengthadds}, so we need only verify the equality of elements in $W$. Let $\lambda = \sum g_i \varpi_i$ and $\mu = \sum h_i \varpi_i$. Then by \cref{thm:explicitpermutation} we have
\begin{equation} w_{\mu}(a+i) = a+i + \sum_{1 \le j<i} h_j n + \sum_{1 \le j < n} h_j(j-n). \end{equation}
For brevity write $k = \sum_{1 \le j < n} h_j(j-n)$, noting that $k \equiv \rot(\mu)$. Thus
\begin{equation} w_{\lambda}((a + k) + i) = (a + k) + i + \sum_{1 \le j<i} g_j n + \sum_{1 \le j < n} g_j(j-n). \end{equation}
By periodicity, we deduce that
\begin{eqnarray} w_{\lambda}(a + k +i &+& \sum_{1 \le j<i} h_j n) \\ \nonumber &=& a + k + i + \sum_{1 \le j<i} g_j n + \sum_{1 \le j < n} g_j(j-n) + \sum_{1 \le j<i} h_j n \\ \nonumber &=&
a + i + \sum_{1 \le j<i} (g_j + h_j) n + \sum_{1 \le j < n} (g_j + h_j)(j-n). \end{eqnarray}
Hence $w_{\lambda}(w_{\mu}(a+i))$ agrees with $w_{\lambda + \mu}(a + i)$.
\end{proof}

In some sense, this proof boils down to the additivity of $\lambda \mapsto \vec{c}_{\lambda}$, see \cref{lem:cadditive}. 

\subsection{Reduced expressions}

Let us briefly recall the concept of reduced expressions for double cosets. We refer the reader to \cite{EKo} for a comprehensive exploration of this idea.

\begin{defn}
Fix a Coxeter system $(W,S)$. 
A  \emph{(singular) multistep expression} is the data of a sequence  of finitary subsets of $S$ of the form $I_0 \subset K_1 \supset I_1 \subset \ldots \subset K_{d} \supset I_d.$ We shall use the notation 
\begin{equation}\label{multistepprototype} IK_{\bullet}=[[I_0 \subset K_1 \supset I_1 \subset \ldots \subset K_{d} \supset I_d]]\end{equation}
to refer to singular multistep expressions.

A \emph{(singular) singlestep expression} is the data of a sequence of finitary subsets $I_0,I_1, \ldots, I_d$, where for each index $0 \le i < d$, the subsets $I_i$ and $I_{i+1}$ differ by the addition or subtraction of a single element. We use the notation 
\begin{equation} \label{singlestepprototype} I_{\bullet} = [I_0,I_1, \ldots, I_d].\end{equation} We also use notation like $[I_0 + s - t + u \ldots]$ to keep track of the simple reflection added or subtracted at each step. \end{defn}

\begin{example} In the affine Weyl group it is more common to keep track of which simple reflections are absent, rather than which are present. There is a singlestep expression
\begin{equation} \label{singlestepex} [\hh{2}, \hh{02}, \hh{024}, \hh{04}, \hh{034}, \hh{34}] = [\hh{2} - 0 - 4 + 2 - 3 + 0]. \end{equation}
Forgetting the order in which $0$ and $4$ are subtracted, we obtain the multistep expression
\begin{equation} \label{multistepex} [[\hh{2} \supset \hh{024} \subset \hh{04} \supset \hh{034} \subset \hh{34}]]. \end{equation}
This may not appear exactly like the form of \eqref{multistepprototype}, e.g. starting with $\supset$ rather than $\subset$, but keep in mind that equalities are permitted in \eqref{multistepprototype}. A more pedantic version of this multistep expression would be
\[ [[\hh{2} \subset \hh{2} \supset \hh{024} \subset \hh{04} \supset \hh{034} \subset \hh{34} \supset \hh{34}]]. \]
Meanwhile, the singlestep expression $[\hh{2} - 4 - 0 + 2 - 3 + 0]$ is distinct from the one in \eqref{singlestepex}, though it is adapted to the same multistep expression \eqref{multistepex}.
\end{example}

\begin{remark} \label{rmk:adapted} There are many multistep expressions adapted to a single step expression (allowing equalities willy nilly), and a finite number of single-step expressions adapted to each multistep expression (choosing the order in which simple reflections are added or subtracted). Ultimately, in the language of \cite{EKo}, all of these will be expressions for the same double coset, and if one of them is a reduced expression then they all are. We write $I_{\bullet} \expr IK_{\bullet}$ when a singlestep expression is adapted to a multistep expression. Later in this paper, singlestep expressions will take on the leading role as they parametrize objects in the diagrammatic 2-category $\DiagSBS$, but it is often more efficient to use multistep expressions for combinatorial definitions. \end{remark}

\begin{defn} \label{def:reducedexpression} Associated to a multistep expression $IK_{\bullet}$ as in \eqref{multistepprototype} we have the element
\begin{equation} \gamma = w_{K_1} \cdot (w_{I_1}^{-1}) \cdot w_{K_2} \cdot (w_{I_2}^{-1}) \cdots w_{K_{d-1}} \cdot (w_{I_{d-1}}^{-1}) \cdot w_{K_d}. \end{equation}
We say the expression is \emph{reduced} if
\begin{equation} \label{lengthviagamma} \ell(\gamma) = \ell(K_1) - \ell(I_1) + \ell(K_2) - \ell(I_2) + \cdots + \ell(K_d), \end{equation}
in which case $\gamma$ is the maximal element of the double coset $p = W_{I_0} \gamma W_{I_d}$, and we say $IK_{\bullet}$ is a reduced expression for $p$ and write $IK_{\bullet} \expr p$.
\end{defn}

\begin{lemma} \label{lem:otherreduced} Given a multistep expression $IK_{\bullet}$, let
\begin{equation} w = \gamma w_{I_d} = w_{K_1} \cdot (w_{I_1}^{-1}) \cdot w_{K_2} \cdot (w_{I_2}^{-1}) \cdots w_{K_d} \cdot (w_{I_d}^{-1}) = \prod_{i=1}^d w_{K_i} w_{I_i}^{-1}. \end{equation}
Then $IK_{\bullet}$ is a reduced expression if and only if 
\begin{equation}\label{lengthviaw} \ell(w) = \sum \ell(w_{K_i} w_{I_i}^{-1}), \qquad \textrm{ and $w$ is minimal in } w W_{I_d}. \end{equation} \end{lemma}

\begin{proof} For any element $w$, $w$ is minimal in the coset $w W_{I_d}$ if and only if $\gamma = w w_{I_d}$ is maximal in this coset, if and only if $\ell(w) + \ell(w_{I_d}) = \ell(\gamma)$. With this observation, it is straightforward to finish the proof (e.g. to go between \eqref{lengthviagamma} and \eqref{lengthviaw}). \end{proof}

To summarize, one should think that to each piece $[[K \supset I]]$ of a multistep expression one has the element $w_{K} w_{I}$, which is the ``distinguished representative'' of the double coset in $W_{K} \backslash W / W_{I}$ containing the identity. When the expression is reduced, then taking the product of these distinguished representatives over a multistep expression one obtains the distinguished representative $\overline{p} w_{I_d}$ of the coset $p$ being expressed, and the lengths add in this product.

\begin{example}\label{example rex for fund weights} The coset $p_{\varpi_k} = \psi_0(\varpi_k) \in W_{\hh{k}} \backslash W / W_{\hh{0}}$ has reduced expression $[[\hh{k} \supset \hh{0,k} \subset \hh{0}]] \expr [\hh{k} - 0 + k]$. Note that $\gamma=w_{\hh{k}}w_{\hh{0,k}}w_{\hh{0}}$ and  $w=w_{\hh{k}}w_{\hh{0,1}}$ in this case.
Conversely, for any $a \ne b \in \Omega$, $[[\hh{a} \supset \hh{a,b} \subset \hh{b}]]$ is a reduced expression for $\psi_b(\varpi_{a-b})$. \end{example}

\begin{example} For any $a \ne b \in \Omega$, $[[\hh{a} \supset \hh{a,b} \subset \hh{a}]]$ is not a reduced expression, as $w = w_{\hh{a}} w_{\hh{a,b}}$ is not minimal in $w W_{\hh{a}}$. \end{example}

\subsection{Reduced expressions for dominant weights}

\begin{defn} \label{defn:fundword} A \emph{fundword} is a finite sequence $(\varpi_{j_1}, \ldots, \varpi_{j_d})$ of fundamental weights. If $\sum_{r=1}^d \varpi_{j_r} = \lambda$ we call the fundword a \emph{reduced expression} or \emph{rex} for $\lambda$, and we often use the notation $\underline{\lambda} = (\varpi_{j_1}, \ldots, \varpi_{j_d})$. \end{defn}

\begin{defn} \label{def:GSfundword} Fix a fundword $\ula = (\varpi_{j_1}, \ldots, \varpi_{j_d})$
for $\lambda \in \Lambda^+_{\wt}$. A sequence of elements $a_0, a_1, \ldots, a_d \in \Omega$ is \emph{adapted} to $\ula$ if $a_{r-1} = a_r + j_r$ for all $1 \le r \le d$, in which we call this sequence a \emph{coloring}. Since the choice of a single $a_i$ determines the rest, and the first and last entries $a_0$ and $a_d$ have special importance, we denote the choice of $\ula$ together with a coloring as
\[ (a_0, \ula) \text{ or } (a_0,\ula, a_d) \text{ or } (\ula, a_d). \]
Such a choice is called a \emph{colored fundword}.
To a colored fundword, we define the corresponding singlestep expression $\GS(\underline{\lambda},a_d)$ as the concatenation of the expressions $[\hh{a_{r-1}} \supset \hh{a_{r-1}, a_r} \subset \hh{a_r}]$ for $1 \le r \le d$. \end{defn}

That is, we concatenate the reduced expressions for $p_{\varpi_{j_r}}$, appropriately rotated, to get the singular expression associated with $\underline{\lambda}$.

\begin{example} Let $\underline{\lambda} = (\varpi_2, \varpi_1, \varpi_3, \varpi_2)$ and $a_0 = 11$. Then
\begin{equation} \GS(11,\ula) = \GS(\ula,3) = [\hh{11} \supset \hh{9,11} \subset \hh{9} \supset \hh{9,8} \subset \hh{8} \supset \hh{5,8} \subset \hh{5} \supset \hh{3,5} \subset \hh{3}]. \end{equation}
For this example the value of $n$ is irrelevant so long as $n \ge 4$ (in order that $\varpi_3$ makes sense).
\end{example}

\begin{remark} Colored fundwords $(\ula, a_d)$ will parametrize the 1-morphisms of $\cwebs$, while singlestep expressions will parametrize 1-morphisms of $\DiagSBS$, and the 2-functor $\GS$ will send $(\ula, a_d) \mapsto \GS(\ula, a_d)$. This justifies the use of $\GS$ in our notation. \end{remark}

The following theorem is effectively just \cref{thm:weightsadd} with different terminology.

\begin{theorem} \label{thm:GSrex} For any rex $\ula$ for $\lambda \in \Lambda^+_{\wt}$, and any $a \in \Omega$, the expression $\GS(\ula,a)$ is a reduced expression for $\psi_a(\lambda)$. \end{theorem}

\begin{proof} Using \cref{lem:otherreduced} we are left with showing that the distinguished representative $w$ of $\psi_a(\lambda)$ (which is minimal in $w W_{\hh{a}}$) satisfies 
\begin{equation} \label{foobarb} w = \prod_{r=1}^{d} w_{\hh{a_{r-1}}} w_{\hh{a_{r-1} a_{r}}}, \end{equation}
and that the lengths add in this expression.
This is precisely what was proven in \cref{thm:weightsadd}. \end{proof}

\subsection{Extended affine Weyl groups}\label{subsec extended affine Weyl groups} We aim to match the combinatorics above with some standard conventions.

For a complex reductive algebraic group $G$ with a pinning $T\subset B \subset G$, we have the extended affine Weyl group defined as follows. 
Let $\rchi=\rchi_*(T)$ denote the cocharacter lattice\footnote{Because of various technicalities involving lattices and Langlands dual groups, we prefer to keep this notation distinct from e.g. $\Lambda_{\wt}$.} of $T$, $\rchi_{\rt} \subset \rchi$ the coroot lattice.
Let $W_f$ denote the (finite) Weyl group, thought of as a subgroup of $\Aut(\rchi)$.
Then $\rchi$ acts on $\rchi$ by translation, and we may define $W_{\ext,G}:=W_f\ltimes \rchi$, which naturally is a subgroup of $\Aut(\rchi)$. 
The affine Weyl group $W \cong W_f \ltimes \rchi_{\rt}$ is isomorphic to a subgroup of $W_{\ext, G}$.
For $\lambda \in \rchi$, we denote by $t_\lambda$ the corresponding translation element such that $t_\lambda(\mu)=\lambda + \mu$.

We shall only be interested in extended affine Weyl groups for two groups, namely, $\GL_n$ and $\PGL_n$, and we denote the corresponding extended Weyl groups as $W_{\ext}'$ and $W_{\ext}$ respectively. Note that the extended Weyl group of $\SL_n$ coincides with the affine Weyl group. 
For $\GL_n$, we explicitly choose $T$ to be the subgroup of diagonal matrices, and identify $\rchi=\rchi_*(T)$ with $\Z^n$ by choosing the basis $\set{\epsilon_i}_{1 \leq i \leq n}$ such that $\epsilon(t)$ is a diagonal matrix with $t$ in the $(i,i)$-th entry and $1$ in the other diagonal entries.
Let $\varpi_k:=(1,\ldots,1,0,\ldots,0) \in \Z^n$ for $1\leq k \leq n$, noting that $\varpi_n=(1,\ldots, 1)$. The cocharacter lattice for $\PGL_n$ is $\rchi/\langle \varpi_n \rangle$, and $W_{\ext} = W'_{\ext}/\langle t_{\varpi_n} \rangle$.

Let $W'$ denote the subgroup of $S_\Z$ generated by $W$ and $\omega$. The following result is well known and is left as an exercise to the reader.

\begin{lemma}\label{lemma extended affine of GLn}
There is an isomorphism of groups $W'_{\ext} \xrightarrow{\sim} W'$
which extends the inclusion $W \to S_{\Z}$. The inverse map sends $\omega \mapsto t_{\varpi_1}w_{\hh{0,1}}w_{\hh{0}}$, which explicitly acts on $\rchi$ as:
$$(\lambda_1, \ldots,\lambda_n)\mapsto (\lambda_n +1, \lambda_1, \ldots, \lambda_{n-1}).$$

Let $P' \subset S_{\Z}$ denote the subgroup of permutations $w$ for which $w(i) \equiv i$ modulo $n$ for all $i$, and let $P = P' \cap W$. Under the above isomorphism, the elements $t_\lambda \in W_{\ext}'$ correspond to the elements in $S_{\Z}$ which send $i \mapsto i + n\lambda_i$. This gives an isomorphism between $\rchi \subset W_{\ext}'$ and $P' \subset W'$, and an isomorphism between $\rchi_{\rt}$ and $P$.
\end{lemma}

\begin{remark} More generally, if we define $\lambda_{kn+i}:=\lambda_i-k$ for any $k \in \Z$ and $1\leq i \leq n$, we can define the inverse isomorphism to be 
$$w \mapsto \left((\lambda_1,\ldots,\lambda_n) \mapsto(\lambda_{w^{-1}(1)}, \ldots, \lambda_{w^{-1}(n)})\right)$$
for $w \in W'$. Note that $s_0$ acts explicitly on $\chi=\Z^n$ via
$$s_0(\lambda_1, \ldots, \lambda_n)=(\lambda_n +1, \lambda_2, \ldots, \lambda_{n-1}, \lambda_{1}-1).$$ \end{remark}

\begin{corollary}\label{formula for omega-k}
 Through the isomorphism in \cref{lemma extended affine of GLn}, we have that
$\omega^k$ is identified with $t_{\varpi_k}w_{\hh{0,k}}w_{\hh{0}}$ for $1\leq k \leq n-1$. Meanwhile, $\omega^n \mapsto t_{\varpi_n}$.
 \end{corollary}
\begin{proof} (Sketch) The easiest way to proceed is to show that both elements have the same action on $\rchi$ using the explicit description in \cref{lemma extended affine of GLn}.
\end{proof}

 
\begin{example}
When $n=8$, the image $t_{\varpi_3}$ in $S_{\Z}$ is explicitly given by 
\[ {
\labellist
\tiny\hair 2pt
 \pinlabel {$1$} [ ] at 120 0
\endlabellist
\centering
\ig{1.0}{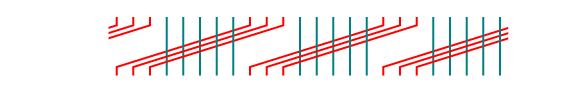}
}.\]
 \end{example}


For the rest of this section, we restrict ourselves to $W_{\ext}=W'/\langle\omega^n \rangle$.
Corresponding to each double coset $p \in W_{\hh{k}}\backslash W/W_{\hh{0}}$, we may associate to it the double coset $\omega^{-k}p \in W_{\hh{0}}\backslash W_{\ext}/W_{\hh{0}}.$

\begin{lemma}\label{lemma coset correspondence for extended affine spherical cosets}
    Let $\lambda \in \Lambda^+_{\wt}$ with $\rot(\lambda)=k$. 
    \begin{enumerate}
    \item $\omega^{-k}\psi_0(\lambda)$ coincides with the double coset $W_{\hh{0}}t_{-\lambda}W_{\hh{0}} \in W_{\hh{0}}\backslash W_{\ext}/W_{\hh{0}}$.
    \item The distinguished representative $w_{\lambda}$ of $\psi_0(\lambda)$ is given by $\omega^{k}t_{-\lambda}$.
    \end{enumerate}
\end{lemma}
\begin{remark}
The negative sign in $t_{-\lambda}$ in \cref{lemma coset correspondence for extended affine spherical cosets}(1) is perhaps unexpected. See \cref{subsec extended affine to affine} for relationships to other conventions in the literature.
\end{remark} 

\begin{proof}
    First, observe that (1) follows from (2): indeed, if the distinguished coset representative $w_\lambda$ equals $\omega^{k}t_{-\lambda}$, then 
    \begin{equation}\psi_0(\lambda)=W_{\hh{k}}\omega^{k}t_{-\lambda}W_{\hh{0}}=\omega^{k}W_{\hh{0}}t_{-\lambda}W_{\hh{0}}.\end{equation}
    
    Hence it suffices to show (2). Let $\lambda = \sum_{1 \leq i \leq n-1}g_i\varpi_i$.
    We proceed by induction on $\sum_{i}g_i$. If this sum is $0$, then we have that $\lambda=0$, $\vec{c}(\lambda)=(0,\ldots, 0)$ and the distinguished element is $t_0=1$.
    
    If $\lambda= \varpi_k$, for some $1\leq k  \leq n-1$, then we know $w_{\lambda}=w_{\hh{k}}w_{\hh{0,k}}$ (see \cref{defn cosets corresponding to fundamental weights}).
    Since \begin{equation} \omega^{-k}=\omega^{n-k}=t_{\varpi_{n-k}}w_{\hh{0,n-k}}w_{\hh{0}}\end{equation} in $W_{\ext}$, we have
    \begin{align}
    \omega^{-k}w_{\varpi_k}=\omega^{-k}w_{\hh{k}}w_{\hh{0,k}}&=w_{\hh{0}}w_{\hh{0,n-k}}\omega^{-k}\\ \nonumber
    &=w_{\hh{0}}w_{\hh{0,n-k}}t_{\varpi_{n-k}}w_{\hh{0,n-k}}w_{\hh{0}}\\ \nonumber 
    &=t_{w_{\hh{0}}w_{\hh{0,n-k}}(\varpi_{n-k})}=t_{-\varpi_k}
    \end{align}
    and the statement follows.
    Now suppose we know that the statement holds for $\lambda$, and want to show that it holds for $\lambda+\varpi_l$.
    Applying \cref{thm:weightsadd}, we have that 
    $$w_{\varpi_l+\lambda}=\omega^{\rot(\lambda)}w_{\hh{l}}w_{\hh{0,l}}\omega^{-\rot(\lambda)}\omega^{\rot(\lambda)}t_{-\lambda}=\omega^{\rot(\lambda)+l}t_{-\varpi_l}t_{-\lambda}=\omega^{\rot(\lambda)+l}t_{-(\lambda+\varpi_l)}$$
    as desired.
\end{proof}

%% file: SSBim_from_Affine_Type_A_Realization.tex
We recall the background material for setting up the algebraic category $\SSBim$ of singular Soergel bimodules from a given realization. We then recall the categorification results and the Williamson-Soergel Hom formula, and then show that these results are valid for the deformed affine type A realization we are working with, under suitable assumptions on the ground field. We also briefly discuss a reformulation of the classical Satake isomorphism, which is categorified by the the functor $\GS_\zeta$ constructed in \cref{subsec the functor}.

\subsection{Deformed affine realizations}
\label{subsec- realizations}


\begin{defn} (See \cite[Definition 3.1]{EWGr4sb}) Let $(W,S)$ be a Coxeter system and $\Bbbk$ a commutative domain. A \emph{realization} of $W$ over $\Bbbk$ is the data of a free
$\Bbbk$ module $\Lambda$ with a $W$-action, together with a set of roots $\{\al_s\} \subset \Lambda$ and coroots $\{\al_s^\vee\} \subset \Lambda^*=\Hom_{\kk}(V,\kk)$ indexed by $s \in S$. We require that the action of
$s \in S$ on $\Lambda$ is given by the formula \[ s(v) = v - \langle \al_s^\vee, v \rangle \al_s. \] We also require $\langle \al_s^\vee, \al_s \rangle = 2$ for all $s$, and $\langle
\al_s^\vee, \al_t \rangle = 0$ whenever $m_{st} = 2$. The \emph{Cartan matrix} of the realization is the $S \times S$ matrix with entries $\langle \al_s^\vee, \al_t \rangle$.
\end{defn}

Given a Coxeter system $(W,S)$, each subset $I \subset S$ defines a parabolic subgroup $W_I$ of $W$, and $(W_I,I)$ itself is a Coxeter system. 
Given a realization $(\Lambda,\Lambda^*,\{\al_s\}_{s\in S},\{\al_s^\vee\}_{s \in S})$ of $(W,S)$, we can then obtain the realization $(\Lambda,\Lambda^*, \set{\al_s}_{s \in I}, \set{\al_s^\vee}_{s\in I})$ of $(W_I,I)$, remembering only the roots and coroots for simple reflections in $I$, and $W_I$ acts on $\Lambda$ by restriction.

For a more detailed exposition on realizations, see \cite[Section 5.7]{EMTW}.

\begin{example}
    For $(W,S)=(S_n, \set{s_1, s_2, \ldots, s_{n-1}})$, we have the permutation realization (over $\kk$) with $\Lambda=\bigoplus_{i=1}^n\kk y_i$. 
    The roots are $ \al_{s_i} := y_i-y_{i+1}$ and $\al_{s_i}^\vee:= y_i^*-y_{i+1}^*$, where $\{y_i^*\}$ is the dual basis to $\{y_i\}$. This induces the usual permutation action on the set $\{y_i\}$.
\end{example}

We shall primarily focus on the affine Weyl group $W_{\aff}$ in Type A, associated to the symmetric group $S_n$. Viewing $W_{\aff}$ as a Coxeter group, it is generated by a set of simple reflections $S = \{s_i\}$ indexed by $\Om = \ZZ/n\ZZ$. We have $m_{s_i s_j} = 3$ if $j = i \pm 1$, and $m_{s_i
s_j} = 2$ otherwise. In this paper all indices will be considered modulo $n$. We say that $\{s_1, \ldots, s_{n-1}\}$ generate the finite Weyl group $S_n$, while $s_0$ is the \emph{affine reflection}. Because it leaves out $s_0$, we call the parabolic subgroup $(S_n,\{s_1, \ldots, s_{n-1}\})$ by the name $W_{\hh{0}}$. 

\begin{example}\cite[Definition 2.8, with $z$ replaced by $\zeta$]{EJY1} \label{defn:Vz} Let $\zeta$ be a formal variable, let $\kk'$ be a commutative domain, and let $\kk = \AC_{\zeta} := \kk'[\zeta,\zeta^{-1}]$. Let $\Lambda$ be the free $\AC_\zeta$-module $\AC_\zeta^{\oplus n}$ with basis $\{x_i\}_{i \in \Om}$. For all $i \in \Om$ define simple roots and coroots by
\begin{equation} \al_i = x_i - \zeta x_{i+1}, \qquad \al_i^\vee = x_i^* - \zeta^{-1} x_{i+1}^*. \end{equation}
This induces the action
\begin{equation} \label{zetarefrep} s_i(x_i) = \zeta x_{i+1}, \quad s_i(x_{i+1}) = \zeta^{-1} x_i, \quad s_i(x_j) = x_j \text{ if } j \notin \{i,i+1\}. \end{equation}
The Cartan matrix is
\begin{equation} \label{slnzeCartan} \left( \begin{array}{cccccc}
2 & -\zeta & 0 & \cdots & 0 & -\zeta^{-1} \\
-\zeta^{-1} & 2 & -\zeta & \cdots & 0 & 0 \\
0 & -\zeta^{-1} & 2 & \cdots &  & 0 \\
\vdots & \vdots & \vdots & \ddots & -\zeta & 0 \\
0 & 0 &   & -\zeta^{-1} & 2 & -\zeta \\
-\zeta & 0 & 0 & \cdots & -\zeta^{-1} & 2
\end{array} \right), \end{equation}
though in the special case $n=2$ the Cartan matrix is
\begin{equation} \label{sl2zeCartan} \left( \begin{array}{cc}
2 & -(\zeta+\zeta^{-1}) \\
-(\zeta+\zeta^{-1}) & 2 \end{array} \right). \end{equation}
Its determinant is $2 - \zeta^n - \zeta^{-n}$.
\end{example} 
\begin{defn}\label{defn-deformed affine realization}
When $\kk'=\Z$, we shall denote the realization from Example \ref{defn:Vz} by $\Lambda_{\zeta}$ and call it \emph{the (integral) deformed affine realization} of type $A_{n-1}$. 
For any commutative $\Z[\zeta,\zeta^{-1}]$-algebra $A$, we shall refer to the same realization but with $\Bbbk = A$ as the \emph{deformed affine realization over $A$}, and denote it by $\Lambda_A$. Note that $\Lambda_A=\Lambda_{\zeta}\otimes_{\Z[\zeta,\zeta^{-1}]} A.$
\end{defn}

The following lemma is an easy consequence of the definitions. 
\begin{lemma}\label{lemma-deformed affine restricted to Sn}
    The restriction of the (integral) deformed affine realization to $W_{\hh{0}} \cong S_n$ is isomorphic to the permutation realization $\bigoplus_{k=1}^n\Z[\zeta,\zeta^{-1}]y_k$, by sending $y_k \mapsto \zeta^{k-1}x_k$. 
    
\end{lemma}

The Dynkin diagram of $W_{\aff}$ admits dihedral symmetry, and these symmetries extend to $\Lambda_{\zeta}$, see \cite[Definition 2.10]{EJY1}.

\begin{defn}We have the \emph{rotation operator} \label{defn-rotation operator}$\sigma$, which is the Dynkin diagram automorphism
$\sigma: \Om \rightarrow \Om: i \mapsto i+1$. 
It acts on $\Lambda_{\zeta}$ by $\sigma(x_i) = x_{\sigma(i)}$ and $\sigma(\zeta) = \zeta$.
We also have the \emph{reflection operator} \label{defn-reflection operator} $\tau: \Om \rightarrow \Om: i \mapsto -i$. It acts on $\Lambda_{\zeta}$ by $\tau(x_i) = x_{1-i}$ and $\tau(\zeta) = \zeta^{-1}$. 
Note that $\tau$ does not preserve the roots; it satisfies $\tau(\alpha_i) = -\zeta^{-1} \alpha_{\tau(i)}$. \end{defn}

Note that $\tau$ fixes the affine vertex, and acts on the finite Weyl group by its traditional Dynkin automorphism (conjugation by the longest element). However, the isomorphism of Lemma \ref{lemma-deformed affine restricted to Sn} does not quite intertwine $\tau$ with the action of the longest element, being off by a multiplicative factor of $\zeta^{n-1}$.

\subsection{Balancedness}
\label{subsec-balanced}
We recall the concept of balancedness.

\begin{defn} (\cite[Definition A.4]{ECathedral}) Let $(W,S)$ be a Coxeter system equipped with a realization, and let $s, t \in S$. It is proven in \cite[\S A.2]{ECathedral} that, whenever $m_{st} = 2k+1$ is odd, one has
\begin{equation} \label{eq:lambdadef} (st)^k(\alpha_s) = \lambda \alpha_t, \quad (ts)^k(\alpha_t) = \lambda^{-1} \alpha_s, \end{equation}
for some scalar $\lambda \in \Bbbk^\times$. The realization is \emph{balanced for} $\{s,t\}$ if $\lambda = 1$. There is a separate definition for $m_{st}$ even that we do not recall here, and it is vacuous when $m_{st} = 2$. We say the realization is \emph{balanced} if it is balanced for all pairs $\{s,t\}$ with $m_{st} < \infty$. \end{defn}

When $m_{st} = 3$, being balanced for $\{s,t\}$ is equivalent to the corresponding entries in the Cartan matrix being $-1$. The deformed affine realization is not balanced.

Under standard notions of what ``positive roots'' should mean, two colinear positive roots should be equal, and both $\alpha_t$ and $(st)^k(\alpha_s)$ should be positive roots. These two requirements are in conflict when the realization is not balanced.

\subsection{Singular Soergel Bimodules from Realizations}
\label{subsec- SSBim}
In \cite{Soer07}, Soergel constructed a category $\SBim$ of bimodules to categorify the Hecke algebra, starting from a nice representation of a Coxeter group $(W,S)$. Williamson, in \cite{WillSingular}, generalized this construction to obtain a 2-category\footnote{In this paper, a 2-category refers to a weak 2-category/bicategory and not necessarily a strict 2-category.} $\SSBim$ of singular Soergel bimodules to categorify the Hecke algebroid\footnote{The Hecke algebroid is sometimes referred to as the Schur algebroid in the literature.}.  We recall these constructions over the next few sections.


Given a realization $(\Lambda,\Lambda^*,\{\al_s\},\{\al_s^\vee\})$ over a base (integral domain) $\kk$, we can construct the graded\footnote{All gradings in this paper are over $\Z$.} $\kk$-algebra $R:= \Sym_\kk \Lambda$ (with $\Lambda$ in degree 2). 
The algebra $R$ inherits a $W$-action from $\Lambda$ which preserves the grading. 
For $I\subset S$ we set 
\[ R^I:=R^{W_I}=\{ f \in R: s(f)=f \text{\; for all }s \in I \}.\]

\begin{defn}\label{def sbsbim}
Let $\SBSBim$, the $2$-category of \emph{singular Bott-Samelson bimodules}, be defined as follows. The objects are finitary subsets $I \subset S$, identified with the rings $R^I$. The morphism category $\Hom(I,J)$ will be a full subcategory of graded $(R^J,R^I)$-bimodules, itself denoted $\bModgr{R^J}{R^I}$. The $1$-morphisms in $\SBSBim$ are generated monoidally by the induction bimodule
$\Bimod{R^I}{R^J}{R^I}$ and the shifted restriction bimodule $\Bimod{R^J}{R^I}{R^I} (\ell(J)-\ell(I))$ for the extensions $R^J \subset R^I$ whenever $I \subset J$. \end{defn}

Before continuing we state some grading conventions. For a graded bimodule $M=\bigoplus_{k \in \Z}M_k$ with $k$-th graded component denoted $M_k$, the shifted bimodule $M(j)$ is obtained by modifying the grading so that $M(j)_k=M_{k+j}$.

For $M,N$ being $(R^I, R^J)$-bimodules, $\Hom_{(R^I,R^J)}(M,N)$ consists of bimodule morphisms $M\rightarrow N$ which preserve the grading. 
We also have the notion of a \emph{graded Hom space}. 
\begin{defn}\label{defn- graded Hom}
   The \emph{graded Hom space} $\Hom^\bullet(M,N)$ for $M, N \in \bModgr{R^I}{R^J}$ is defined as 
   $$\Hom^\bullet(M,N):=\bigoplus_{k\in \Z} \Hom_{(R^I,R^J)}(M,N(k)).$$
\end{defn}
For a Laurent polynomial $p=\sum p_iv^i \in \Z_{\geq0}[v,v^{-1}]$ with positive integer coefficients and $M \in \bModgr{R^I}{R^J}$, we define $M^{\oplus p} \in \bModgr{R^I}{R^J}$ to be
$$M^{\oplus p}:=\bigoplus_{i\in \Z}M(-i)^{\oplus p_i}.$$

\begin{defn}
    Let $A$ be a ring and $M$ be a graded $A$-module. Then we say that $M$ is a \emph{graded free} $A$-module if $M\cong A^{\oplus p}$ for a Laurent polynomial $p \in \Z_{\geq 0}[v,v^{-1}]$ with positive integer coefficients, in which case $p$ is called the \emph{graded rank} of $M.$
\end{defn}

Returning to $\SBSBim$, we can describe the 1-morphisms in a more combinatorial manner.

\begin{defn}
Associated to a singular multistep expression as in \eqref{multistepprototype}, we associate the \emph{Bott-Samelson bimodule}
$$\BS(IK_{\bullet}):= R^{I_0}\otimes_{R^{K_1}} \otimes R^{I_1}\otimes \ldots \otimes_{R^{K_d}}R^{I_d}(\sum_{i=1}^d \ell(K_i)-\ell(I_i)) \in \Hom_{\SBSBim}(I_d,I_0).$$
\end{defn}

\begin{example}
    For $S_4$ with simple reflections $\set{s,t,u}$, we have a multistep expression
    $$[[t \subset stu \supset s\subset st]].$$
    The corresponding Bott-Samelson bimodule is the graded $(R^t, R^{st})$-bimodule
    $$R^t \otimes_{R^{stu}}R^s(5).$$
\end{example}

\begin{defn} The $2$-category $\SSBim$ of \emph{singular Soergel bimodules} is the Karoubi envelope (of the graded, additive completion) of $\SBSBim$. For finitary subsets $I, J \subset S,$ we use the notation $\SBimod{I}{J}$ to refer to the category $\Hom_{\SSBim}(J,I)$. 
\end{defn}
Given two graded rings $R_1$ and $R_2$, the category $\bModgr{R_1}{R_2}$ is Karoubi closed. 
Since $\Hom_{\SBSBim}(J,I)$ is a full subcategory of $\bModgr{R^I}{R^J}$, we have that $\SBimod{I}{J}=\Hom_{\SSBim}(J,I)$ is also naturally a full subcategory of $\bModgr{R^I}{R^J}$.

We will see in \cref{thm-categorification} that (under certain assumptions) the indecomposable bimodules in $\SBimod{I}{J}$ are in bijection (up to isomorphism and grading shift) with double cosets $p \in W_I \backslash W / W_J$. Let $p$ correspond to the indecomposable bimodule $B_p$. Whenever $IK_{\bullet}$ is a reduced expression for $p$ as in \cref{def:reducedexpression}, $B_p$ will be a direct summand with multiplicity one inside $\BS(IK_{\bullet})$.



\subsection{The Hecke Algebroid}
\label{subsec-Hecke}

The Hecke algebroid associated to a Coxeter system $(W,S)$ is a relative version of the Hecke algebra associated to all pairs of finitary subsets $I, J \subset S$. One may think of it as a $\Z[v,v^{-1}]$-linear category with objects  indexed by the finitary subsets of $S$. 
We now go into the details of the construction.

Let $\HH$ denote the Hecke algebra associated to $(W,S)$. This is an algebra over $\Z[v,v^{-1}]$. Let $H_w$ denote the standard basis element and $\KL_w$ the Kazhdan-Lusztig basis element corresponding to $w \in W$. The Kazhdan-Lusztig basis is typically complicated, but not when $w = w_J$ for $J$ finitary, where we have
\[ \KL_{w_J} = \sum_{w \in W_J} v^{\ell(w_J) - \ell(w)} H_w. \]
Note that $\KL_{w_J}\KL_{w_J}=\pi(J)\KL_{w_J}$ for all finitary $J \subset S$, where $\pi(J)=\displaystyle \sum_{w\in W_J}v^{\ell(w_J)-2\ell(w)}.$ Thus $\KL_{w_J}$ is a quasi-idempotent in $\HH$. One should think of the Hecke algebroid as something like a (partial) Karoubi envelope associated to these quasi-idempotents.

For all pairs of finitary subsets $I,J \subset S,$ define:
\begin{align*}
    {}^I\HH &= \KL_{w_I}\HH, \quad \text{a right ideal,}\\
    \HH^J &= \HH\KL_{w_J}, \quad \text{a left ideal,}\\
    {^I}\HH^{J} &= {}^I\HH\cap\HH^J.
\end{align*}
Given finitary subsets $I, J, K$, we define a multiplication:
\begin{align*}
    \ast_J :{}^I\HH^J\times {}^J\HH^K \longrightarrow {}^I\HH^K:
    (h_1,h_2) \longmapsto h_1 \ast_J h_2 := \frac{1}{\pi(J)}h_1h_2.
\end{align*}
The operation $\ast_J$ is well-defined as the numerator is divisible by the denominator.
Indeed, for $h_1 \in \Hecke{I}{J}$ and $h_2 \in \Hecke{J}{K}$, we may write $h_1=\KL_{w_I}y_1=x_1\KL_{w_J}$ and $h_2=\KL_{w_J} y_2=x_2\KL_{w_K}$ for some $x_i, y_i \in \HH$, so that 
$$h_1\ast_J h_2=\frac{1}{\pi(J)}h_1h_2=\frac{1}{\pi(J)} x_1\KL_{w_J}\KL_{w_J}y_2=x_1\KL_{w_J}y_2=\KL_{w_I}y_1y_2=x_1x_2\KL_{w_K} \in \Hecke{I}{K}.$$

\begin{defn}\cite[Definition 2.3.1]{WillSingular}
    The Hecke algebroid $\HAbd$ is the following $\Z[v,v^{-1}]$-linear category:
    \begin{itemize}
        \item Objects are given by finitary subsets $I \subset S$;
        \item For finitary subsets $I, J \subset S$, $\Hom(J,I):={}^I\HH^J$, and composition ${}^I\HH^J\times {}^J\HH^K \rightarrow {}^I\HH^K$ is given by $\ast_J$.
    \end{itemize}
\end{defn}
As explained in \cite[section 2.3]{WillSingular}, $\Hecke{I}{J}$ has a standard basis $\{\HAB{I}{J}{p}: p\in W_I\backslash W/W_J\}$ and a Kazhdan-Lusztig (KL) basis $\{\KLB{I}{J}{p}:p \in W_I\backslash W/W_J \}$ where basis elements are indexed by double cosets $p \in W_I\backslash W/ W_J$. As usual, the KL basis is uniquely determined by two properties: self-duality under the bar involution (inherited from the Hecke algebra), and the change of basis formula
\begin{equation} \KLB{I}{J}{p} = \HAB{I}{J}{p} + \sum_{q < p} \textrm{coeff} \cdot \HAB{I}{J}{q} \end{equation}
where each coefficient lives in $v \Z[v]$. Here $<$ indicates the Bruhat order on double cosets.

When the parabolic subgroups $I$ and $J$ are understood, we write $H_p$ instead of $\HAB{I}{J}{p}$, and $\KL_{p}$ instead of $\KLB{I}{J}{p}$.

\begin{remark}
    When $(W,S)$ is the affine Weyl group associated to a simply connected complex semisimple group $G$, and $I$ is the finitary subset such that $W_I$ is the Weyl group associated to $G$, then $\End(I)=({}^I\HH^{I},\ast_I)$ coincides with the spherical Hecke algebra of $G$ (cf. \cref{lemma affine spherical Hecke in Hecke algebroid})
\end{remark}

\begin{remark}
One may also provide the following alternate description for the Hecke algebroid.
For a finitary subset $I\subset S$, let $\HH_I$ denote the Hecke algebra of $(W_I,I)$ thought of as a subalgebra of $\HH$, and let $\triv_I: \HH_I \rightarrow \Z[v,v^{-1}]$ denote the trivial representation (which sends the standard basis elements $H_{w}$ to $v^{-\ell(w)}$, and hence sends $\KL_{w_I}$ to $\pi(I)$). 
Then for finitary $I,J \subset S$, ${}^I\HH^J$ can be identified with $\Hom_{\HH}(\triv_J\otimes_{\HH_J}\HH, \triv_I\otimes_{\HH_I}\HH)$ and $\ast_J:{}^I\HH^J\times {}^J\HH^K \rightarrow {}^I\HH^K$ coincides with composition of (right) $\HH$-module morphisms via this identification.
\end{remark}

When $I \subset J$ one defines ${}^I H^J \in {}^I \HH^J$ to be equal to $\HAB{I}{J}{p}$, associated to the double coset $p$ containing the identity element. As $p$ is minimal in the Bruhat order we have $\HAB{I}{J}{p} = \KLB{I}{J}{p}$, and the underlying element in the Hecke algebra is $\KL_{w_J}$.  Define ${}^J H^I$ similarly. These are the \emph{standard generators}, as Williamson calls them, and they generate the Hecke algebroid. In \cite[Proposition 2.3.3]{WillSingular} one can find a rule for computing the action of standard generators on the standard basis.

Let $IK_{\bullet}$ be a multistep expression as in \eqref{multistepprototype}. We define the element $\bsh(IK_{\bullet}) \in {}^{I_0} \HH^{I_d}$ as the corresponding product of standard generators
\begin{equation} \bsh(IK_{\bullet}) = {}^{I_0} H^{K_1} \star_{K_1} {}^{K_1} H^{I_1} \star_{I_1} {}^{I_1} H^{K_2} \star_{K_2} \cdots \star_{K_d} {}^{K_d} H^{I_d}. \end{equation}
One can use a similar formula to define $\bsh(I_{\bullet})$ for singlestep expressions as in \eqref{singlestepprototype}. Whenever $I_{\bullet}$ is adapted to $IK_{\bullet}$ (see \cref{rmk:adapted}) one has $\bsh(I_{\bullet}) = \bsh(IK_{\bullet})$.

By iterating the rule from \cite[Proposition 2.3.3]{WillSingular}, one obtains an expansion of $\bsh(I_{\bullet})$ in the standard basis, of the form
\begin{equation} \label{eq:expandbs} \bsh([I_0,I_1,\ldots,I_d]) = \sum_{[p_0,p_1,\ldots,p_d]} \textrm{coeff} \cdot H_{p_d}. \end{equation}
The sum ranges over sequences of double cosets $p_i \in W_{I_0} \backslash W / W_{I_i}$ such that either $p_i \subset p_{i+1}$ or $p_i \supset p_{i+1}$, and such that $p_0$ contains the identity element. The coefficient is a certain polynomial in $\N[v,v^{-1}]$. One feature of this expansion is that, when $I_{\bullet}$ is a reduced expression for $p$, then the coefficient of $H_p$ will be $1$. We will not need the details except in the proof of \cref{lem:dim2}.

\subsection{Categorification results and Hom formula}
\label{subsec-categorification}
In section \ref{subsec- SSBim}, we defined (following \cite{WillSingular}) a 2-category constructed from a given realization of a Coxeter group $(W,S)$. 
Under nice enough assumptions (as in \cite{WillSingular}), or under a suitable modification of $\SSBim$ as in \cite{Abe} (see \cref{rmk-Abe SSBim}), the so-called Soergel-Williamson categorification theorem and Hom formula both hold. These assumptions do not hold for the deformed affine realization, so we will need to extend the known results. Let us first recall the results in the literature, and discuss assumptions afterwards. Throughout this section, we shall assume that the base $\kk$ is an infinite field.

\begin{theorem}[Categorification Theorem] \cite[Theorem 1, Theorem 2]{WillSingular}, \cite[Theorem 1.1]{Abe} \label{thm-categorification}
    If \cref{willassume} \hyperref[assumption 1]{(1)}, \hyperref[assumption 2]{(2)} and \hyperref[assumption 3]{(3)} hold, we have the following.
\begin{enumerate}
\item There exists a $\Z[v,v^{-1}]$-linear isomorphism ${}_{I}\ch_{J}\colon [\SBimod{I}{J}]\xrightarrow{\sim}{}^I{\mathcal{H}}^J$ such that the diagram
\[
\begin{tikzcd}\relax
[\SBimod{I}{J}]\times [\SBimod{J}{K}]\arrow[r,"\otimes_{J}"]\arrow[d,"{}_{I}{\ch}_{J}\times {}_{J}{\ch}_{K}"] &\relax [\SBimod{I}{K}]\arrow[d,"{}_{I}{\ch}_{K}"]\\
{}^{I}{\mathcal{H}}^{J}\times{}^{J}{\mathcal{H}}^{K}\arrow[r,"*_{J}"] & {}^{I}{\mathcal{H}}^{K}
\end{tikzcd}
\]
is commutative.
\item For each $p\in W_{I}\backslash W/W_{J}$ there exists an indecomposable object
$B_p \in \SBimod{I}{J}$
such that
${}_{I}{\ch}_{J}(B_p)\in {}^{I}{H}^{J}_{x} + \sum_{y < x}\Z[v,v^{-1}]{}^{I}{H}^{J}_{y}$.
Moreover, it is unique up to isomorphism.
\item For any indecomposable object $B\in \SBimod{I}{J}$ there exists $p\in W_{I}\backslash W/W_{J}$ and $n\in \Z$ such that
$B \simeq B_p(n)$.
\item If a singular multistep expression $IK_{\bullet}$ is a reduced expression for the double coset $p \in W_{I}\backslash W/W_{J}$, then $B_p$ is a summand (with multiplicity 1) of the Bott-Samelson bimodule $\BS(IK_{\bullet})$. All other summands are isomorphic to $B_q(n)$ for some $n \in \Z$ and some $q \le p$ in the Bruhat order on double cosets (see \cite[Proposition 3.3]{KELPBruhat}).
\end{enumerate}
\end{theorem}

Whenever $I \subset J$, $[[I \subset J]]$ is a reduced expression for the double coset $p$ containing the identity, and the induction bimodule $\BS(I \subset J) := \BS([[I \subset J]]) = R^I$ is indecomposable, so equals $B_p$. We have ${}_{I}{\ch}_{J}(B_p) = {}^I\KL^J_p = \KL_{w_J}$. Similar statements can be made for $[[J \supset I]]$ and the shifted restriction bimodule, which has character $\KL_{w_J}$ as well. The categorification theorem sends the generators of $\SSBim$ to the standard generators of the Hecke algebroid. Property (1) above then implies that the character of $\BS(I_{\bullet})$ is $\bsh(I_{\bullet})$ for any singular expression $I_{\bullet}$.

\begin{remark} For certain nice realizations (typically in characteristic zero) the \emph{Soergel conjecture} holds, which states that ${}_{I}{\ch}_{J}(B_p) = {}^I\KL^J_p$ for all $p$. Whether this holds is not relevant in this paper.
\end{remark}

Before stating the Williamson-Soergel Hom formula, we recall the standard pairing on the Hecke algebra $\HH$.
Let $\bar{( \cdot)}: \ZZ[v,v^{-1}] \rightarrow \ZZ[v, v^{-1}]$ denote the ring automorphism $v \mapsto v^{-1}$. 
This extends to the \emph{bar involution} or the \emph{Kazhdan-Lusztig involution} 
$$\bar{(\cdot)}: \HH \rightarrow \HH: H_w \mapsto H_{w^{-1}}^{-1}.$$
We also have the \emph{Kazhdan-Lustig anti-involution} $\omega\colon \HH\to \HH$  given by 
$$\omega(\sum_{w\in W}a_{w}H_{w}) := \sum_{w\in W}\overline{a_{w}}H_{w}^{-1}.$$
This is an anti-involution and $\omega(\KL_{w_I}) = \KL_{w_{I}}$.
We also have the \emph{standard trace} $\epsilon: \HH \rightarrow \ZZ[v,v^{-1}]$ which extracts the coefficient of $H_1$ when writing an element of $\HH$ in the standard basis. 
Note that $\epsilon(ab)=\epsilon(ba)$ (see, for example, \cite[Corollary 3.17]{EMTW}).
\begin{defn}\label{defn-standard form}
    The \emph{standard form} $(-,-):\HH \times \HH \rightarrow \ZZ[v,v^{-1}]$ is the sesquilinear form 
    $$(a,b):=\epsilon(\omega(a)b) = \epsilon(b\omega(a)).$$
    This form can be extended to the Hecke algebroid as $(-,-):\Hecke{I}{J} \times \Hecke{I}{J} \rightarrow \ZZ[v,v^{-1}]$, \begin{equation} (a,b):= v^{-\ell(J)}\epsilon(b*_J\omega(a)).\end{equation} Note that $\omega(\Hecke{I}{J}) = \Hecke{J}{I}$, so that $b *_J \omega(a)$ makes sense.\footnote{The factor of $v^{-\ell(J)}$ is not always included in the literature.}
\end{defn}

Just as for the Hecke algebra, the standard form on the Hecke algebroid has the following property, the \emph{asymptotic orthonormality} of the KL basis.

\begin{lemma} \label{lem:ortho} Let $p, q \in W_I \backslash W/W_J$. Then
\begin{equation} (\KLB{I}{J}{p}, \KLB{I}{J}{q}) \in \delta_{pq} + v \Z[v]. \end{equation}
\end{lemma}

Now we state the Hom formula. We recommend the reader see \cite[Theorem 4.16]{EKLP} for a reformulation which is often more useful in practice, one which also makes \cref{lem:ortho} more transparent.

\begin{theorem}[Hom formula]\cite[Theorem 7.4.1]{WillSingular}, \cite[Theorem 2.54]{Abe} \label{thm-hom formula}
If \cref{willassume} \hyperref[assumption 1]{(1)}, \hyperref[assumption 2]{(2)} and \hyperref[assumption 3]{(3)} hold, we have the following. 
Let $M_{1},M_{2}\in \SBimod{I}{J}$.
Then $\Hom_{\SBimod{I}{J}}(M_{1},M_{2})$ is a graded free $R^{I}$-module and we have
\[
\grk_{R^{I}} \Hom^\bullet_{\SBimod{I}{J}}(M_{1},M_{2}) = 
({}_{I}{\ch}_{J}(M_{1}), {}_{I}{\ch}_{J}(M_{2}))
\]
\end{theorem}

\begin{remark} The morphism space is also free as a right $R^J$-module with a different graded rank. \end{remark}

We now give an extended example, as the proof of the following lemma. This lemma will play a role in our proof of the well-definedness of the 2-functor $\GS$.

\begin{lemma} \label{lem:dim2} Let $(W,S)$ be the affine Weyl group using conventions from \cref{sec:combo}. Consider the fundword $\ula = (\varpi_{k-1}, \varpi_1)$ for some $2 \le k \le n$, and let $I_{\bullet} = \GS(\ula,0)$, see \cref{def:GSfundword}. Explicitly, we have
\begin{equation} I_{\bullet} = [\hh{k} \supset \hh{k1} \subset \hh{1} \supset \hh{01} \subset \hh{0}]. \end{equation}
Then $\dim \End^0(\BS(I_{\bullet})) = 2$. \end{lemma}

\begin{proof}
By the Hom formula, we should compute the pairing $(\bsh(I_{\bullet}), \bsh(I_{\bullet}))$, which yields the graded rank of the endomorphism space over $R^{\hh{k}}$. This does not agree with the graded dimension of the endomorphism space, but they have the same leading term (the smallest power of $v$) since $R^{\hh{k}}$ is concentrated in non-negative degree with a one-dimensional degree zero space. Thus if the pairing lives in $\Z[v]$, then $\dim \End^0(\BS(I_{\bullet}))$ agrees with the coefficient of $v^0$.

Our goal is to prove that
\begin{equation} \label{dim2bsdecomp} \bsh(I_{\bullet}) = \KL_{\psi_0(\varpi_1 + \varpi_{k-1})} + \KL_{\psi_0(\varpi_k)}. \end{equation}
By the almost orthonormality of the KL basis \cref{lem:ortho}, we deduce that 
\[ (\bsh(I_{\bullet}), \bsh(I_{\bullet})) \in 2 + v \Z[v], \]
whence we deduce that $\dim \End^0(\BS(I_{\bullet})) = 2$.

First assume $k < n$. We expand $\bsh(I_{\bullet})$ in the standard basis. We consider all sequences $[p_0 \supset p_1 \subset p_2 \supset p_3 \subset p_4]$ as in \eqref{eq:expandbs}. There is a unique choice of $p_0$ and $p_1$ and $p_2$: the underlying set of both $p_0$ and $p_1$ is $W_{\hh{k}}$, and $p_2 = W_{\hh{k}} 1 W_{\hh{1}}$. In fact one has
\begin{equation} \bsh([\hh{k} \supset \hh{k1} \subset \hh{1}]) = H_{p_2} \end{equation}
since this is a reduced expression and $p_2$ is minimal in the Bruhat order.

There are two double cosets in $W_{\hh{k}}\backslash W /W_{\hh{01}}$ contained in $p_2$, denoted $p_3 < p'_3$. They satisfy
\begin{equation} \underline{p}_3 = \underline{p}_2 = 1, \qquad \overline{p}_3' = \overline{p}_2 = w_{\hh{k}} w_{\hh{k1}}^{-1} w_{\hh{1}}. \end{equation} 
Finally, $p_3$ and $p_3'$ are contained respectively in $p_4$ and $p_4'$, with
\begin{equation} p_4 = \psi_0(\varpi_{k}), \qquad p_4' = \psi_0(\varpi_1 + \varpi_{k-1}). \end{equation}

As proven in \cref{thm:GSrex}, $I_{\bullet}$ is a reduced expression for $\psi_0(\varpi_1 + \varpi_{k-1})$, so the coefficient of $H_{p_4'}$ in this expansion will be $1$. Meanwhile, according to \cite[Proposition 2.3.3(1)]{WillSingular}, the coefficient of $H_{p_4}$ will be
\begin{equation} v^{\ell(\overline{p_2}) - \ell(\overline{p_3})} v^{\ell(\underline{p_3}) - \ell(\underline{p_4})} \frac{\pi(\hh{0k})}{\pi(\hh{01k})}, \end{equation}
where $\pi$ is the balanced Poincare polynomial of the parabolic subgroup. This ratio of Poincare polynomials is that for $S_k$ over $S_{k-1}$, and is equal to $[k] = v^{k-1} + v^{k-3} + \ldots + v^{-(k-1)}$. The difference $\ell(\overline{p_2}) - \ell(\overline{p_3})$ is equal to $k-1$, whereas the relevant minimal elements are equal to the identity. Thus we have
\begin{equation} \bsh(I_{\bullet}) = H_{p_4'} + v^{k-1} [k] H_{p_4} = H_{\psi_0(\varpi_1 + \varpi_{k-1})} + (1 + v^2 + \ldots + v^{2(k-1)}) H_{\psi_0(\varpi_k)}. \end{equation}

Now we re-express this in the KL basis. Note that $\bsh(I_{\bullet})$ is self-dual (as is any product of the self-dual generators), and $H_{\psi_0(\varpi_k)} = \KL_{\psi_0(\varpi_k)}$ since the double coset is minimal in the Bruhat order, so it is also self-dual. Thus the difference is self-dual, and the coefficient of $H_{\psi_0(\varpi_k)}$ lives in $v \Z[v]$. By the defining property of the KL basis we must have
\begin{equation} \KL_{\psi_0(\varpi_1 + \varpi_{k-1})} = \bsh(I_{\bullet}) - H_{\psi_0(\varpi_k)} = H_{\psi_0(\varpi_1 + \varpi_{k-1})} + (v^2 + \ldots + v^{2(k-1)}) H_{\psi_0(\varpi_k)}. \end{equation}
In particular, \eqref{dim2bsdecomp} holds as desired.

When $k=n$, replace $\varpi_k$ with the zero weight, so that $p_4 = \psi_0(0)$ is the double coset containing the identity. One should also interpret $\hh{0k}$ as $\hh{0}$ in this special case. Otherwise, the proof above works directly.
\end{proof}

\begin{remark} \label{rmk:dim2vsrepthry} On the other side of the Soergel Satake equivalence, the corresponding statement is
\begin{equation} \label{dim2inreptheory} \dim \End(L_{\varpi_{k-1}} \ot L_{\varpi_1}) = 2. \end{equation}
One typically proves this using Schur's lemma and the decomposition
\begin{equation}\label{1kdecompreptheory} L_{\varpi_{k-1}} \ot L_{\varpi_1} \cong L_{\varpi_{1} + \varpi_{k-1}} \oplus L_{\varpi_k}. \end{equation}
(Replace $\varpi_k$ by the zero weight if $k=n$.) The decategorified version of this decomposition is \eqref{dim2bsdecomp}.
Note however that \eqref{dim2inreptheory} holds even in finite characteristic where \eqref{1kdecompreptheory} might fail! Similarly, the decomposition
\begin{equation} \BS(I_{\bullet}) \cong B_{\psi_0(\varpi_1 + \varpi_{k-1})} \oplus B_{\psi_0(\varpi_k)} \end{equation}
need not hold in finite characteristic. We will not rely on such decompositions in this paper, only on the decategorified statements and their implications for the size of morphism spaces.
\end{remark}

\subsection{The Satake isomorphism in terms of the Hecke algebroid}\label{subsec Satake Hecke algebroid}
In this section, we rephrase the Satake isomorphism (which was categorified into the geometric Satake equivalence) as an isomorphism between the Grothendieck group of $\Rep(SL_n)$ and a version of the spherical Hecke algebra for $W_{\ext}$ constructed within the Hecke algebroid.
We then define $\Rep^{\Om}(SL_n)$ (or rather, its quantum analogue), and reinterpret the the Satake isomorphism as a functor into the Hecke algebroid.


We retain our notation from \cref{sec:combo}, so that $\Lambda_{\wt}$ denotes the weight lattice of $\mathfrak{sl}_n$, $\Lambda_{\wt}^+$ denotes the subset of dominant weights, $\rot: \Lambda_{\wt} \rightarrow \Om$ denotes the quotient map with respect to the root lattice, and $\psi_a$ denotes the bijection in \cref{def:psia}.
Let $\Z[\Lambda_{\wt}]$ denote the the group algebra of $\Lambda_{\wt}$.
It is well known that $[\Rep(\SL_n)]=\Z[\Lambda_{\wt}]^W$, the subring of $W$-invariants.
One one side of the Satake isomorphism is $[\Rep(\SL_n)]\otimes_{\Z} \Z[v^{\pm 1}]=\Z[v^{\pm 1}][\Lambda_{\wt}]^W.$

There is a decomposition $$\Rep(SL_n)=\bigoplus_{k \in \Om}\Rep(SL_n)_{k}$$
where $\Rep(SL_n)_{k}$ is the full subcategory of representations with weights lying entirely in the root coset $k$. 
This decomposition naturally induces an $\Om$-grading
\begin{equation}\label{eq grading on representation ring of sl_n}
   \Z[\Lambda_{\wt}]^W= \bigoplus_{k \in \Om}\Z[\Lambda_{\wt}]^W_k,
\end{equation}
which further induces an $\Om$-grading on $\Z[v^{\pm 1}][\Lambda_{\wt}]^W$.

The other side of the Satake isomorphism would then be the spherical Hecke algebra associated to the extended affine Weyl group $W_{\ext}$ of $\PGL_n$.
We do not discuss the spherical Hecke algebra here, but instead work within the Hecke algebroid associated to $W$.

Recall that we had the rotation operator $\sigma$ which was an automorphism of the Dynkin diagram (see \cref{defn-rotation operator}). 
Then $\sigma$ naturally acts on the Hecke algebroid, sending $\HAB{I}{J}{p}$ to $\HAB{\sigma(I)}{\sigma(J)}{\sigma(p)}$.

\begin{defn}\label{defn the algebra HHH}
Let 
\begin{equation}\label{eq HHH graded algebra}
    \HHH=\bigoplus_{0\leq k \leq n-1}\Hecke{\hh{k}}{\hh{0}}.
\end{equation} 
$\HHH$ is naturally a right module over $\Hecke{\hh{0}}{\hh{0}}$, but we define a $\Z[v,v^{-1}]$-algebra structure on $\HHH$ as follows.
For $f \in \Hecke{\hh{k}}{\hh{0}}$, $g \in \Hecke{\hh{l}}{\hh{0}}$, define
\begin{equation} \label{eq:badmultyuck}
    f\bullet g:= \sigma^l(f)\ast_{\hh{l}}g \in \Hecke{\hh{k+l}}{\hh{0}}
\end{equation}
and extend it linearly to $\HHH$. 
Then $\HHH$ is an $\Om$-graded algebra (with respect to the decomposition in \cref{eq HHH graded algebra}) and has $\Hecke{0}{0}$ as a subalgebra.
\end{defn}

The algebra $\HHH$ is in fact isomorphic to the spherical Hecke algebra associated to $W_{\ext}$, and hence may be used on the other side of the Satake isomorphism. Then Lusztig's work in \cite{LusztigGS} may be reformulated as follows; see \cref{subsec extended affine to affine} for further details.

\begin{theorem}\label{thm Soergel Satake mini}
    There is an isomorphism of $\Z[v^{\pm 1}]$-algebras
    $$\Z[v^{\pm 1}][\Lambda_{\wt}]^W \xrightarrow{\sim} \HHH$$
    which sends $\ch(L_\lambda)$ to $\KLB{\hh{\rot(\lambda)}}{\hh{0}}{\psi_0(\lambda)}$ for $\lambda \in \Lambda_{\wt}^+$, where $\ch(L_\lambda)$ is the character of the irreducible representation of highest weight $\lambda$.
    This isomorphism respects the $\Om$-grading on these algebras given by \cref{eq grading on representation ring of sl_n}, \cref{eq HHH graded algebra}.
\end{theorem}

The algebra structure on $\HHH$ may appear slightly artificial, hence we prefer to work directly with the Hecke algebroid, and do not discriminate between maximal finitary subsets of $S$.

\begin{defn}\label{defn KRepom}
    The category $\KRepom$ is defined as follows.
    \begin{itemize}
        \item  The objects of $\KRepom$ are elements of $\Om=\Z/n\Z$.
        \item The morphism space $\Hom(l,k)$ is the $\Z[\Lambda_{\wt}]^W_{k-l}$, spanned by the characters of representations whose weights live in the root coset $k-l \in \Om$. Composition is given by multiplication within $\Z[\Lambda_{\wt}]^W.$
    \end{itemize}
    This is naturally a $\Z$-linear category. 
    We denote by $\KRepom \otimes_{\Z}\Z[v^{\pm 1}]$ the $\Z[v^{\pm 1}]$-linear category with same objects, but with morphism spaces given by $\Hom_{\KRepom}(l,k)\otimes_{\Z}\Z[v^{\pm 1}]=\Z[v^{\pm 1}][\Lambda_{\wt}]^W_k$.
\end{defn}

The following theorem is equivalent \cref{thm Soergel Satake mini}, after translating the unnatural multiplication of \eqref{eq:badmultyuck} into the more natural multiplication in the Hecke algebroid $\HAbd$.

\begin{theorem}\label{thm Soergel Satake}
    There is a functor $\KRepom \rightarrow \HAbd$ which
    sends the object $k \in \Om$ to the maximal finitary subset $\hh{k}$,
    and sends $\ch(L_{\lambda})\in \Hom(l,k)$ to $\KLB{\hh{k}}{\hh{l}}{\psi_l(\lambda)} \in \Hecke{\hh{k}}{\hh{l}}$.
    The induced functor 
    $$\KRepom \otimes_{\Z}\Z[v^{\pm 1}] \rightarrow \HAbd$$ is fully faithful, with essential image the set of maximal finitary subsets of $S$.
\end{theorem}

We now discuss the categorification of $\KRepom$, which appears on one side of the Soergelified Satake equivalence.

In \cite[Section 4.2]{EQuantumI}, it is explained how to construct a 2-category from an \emph{$H$-graded-monoidal category} (\cite[Definition 4.5]{EQuantumI}) for a finite abelian group $H$. 
As seen before, the category $\Rep(SL_n)$ is naturally an $\Om$-graded-monoidal category, to which applying the construction yields a 2-category $\Rep^{\Om}$.
Similarly, for generic $q$, $\Rep(U_q(\mathfrak{sl}_n))$ is also an $\Om$-graded-monoidal category, to which applying the construction yields a 2-category $\qRep$. 
We recall the precise definition of $\qRep$ below, the definition of $\Rep^{\Om}$ is similar.

\begin{definition}\label{defn- qRep Omega}
    The 2-category $\qRep$ is defined as follows.
    \begin{itemize}
        \item The set of objects in the 2-category is $\Om$. 
        \item The Hom-category $\Hom(l,k)$ is the full subcategory 
        of $\Rep(U_q(\mathfrak{sl}_n))$ (defined over $\Q(q)$) consisting of representations whose weights live in the root coset $k-l$.
        \item Composition of 1-morphisms is given by tensor product of representations.
    \end{itemize}
\end{definition}

\begin{remark}
    In \cite{EQuantumI}, the set of objects is taken to be the set of ``removable'' vertices of the Dynkin diagram, which is a torsor over $\Om$. This set is then shown to be in canonical bijection with $\Om$, with the affine vertex identifying with $0$ (see \cite[Claim 4.13]{EQuantumI}).
\end{remark}

Recall that at the level of Grothendieck rings, $[\Rep(U_q(\mathfrak{sl}_n))]=[\Rep(SL_n)]=\Z[\Lambda_{\wt}]^W$ as $\Om$-graded rings. The following lemma is then a consequence of \cref{defn KRepom}, \cref{defn- qRep Omega} and this observation.

\begin{lemma}\label{lemma decategorification of Rep-omega}
    There are canonical isomorphisms $[\Rep^{\Om} ]=[\qRep]=\KRepom$.
\end{lemma}

In view of \cref{lemma decategorification of Rep-omega}, \cref{thm Soergel Satake} is the decategorified version of the geometric Satake functor when interpreted in terms of singular Soergel bimodules, or its quantization $\GS_\zeta$ thereof, constructed in \cref{subsec the functor}.

\subsection{Assumptions on a realization}


Now we discuss the assumptions made by Williamson in his proofs of these results, and how to weaken those assumptions.
First we recall the notion of reflection faithfulness for a representation $V$ of a Coxeter group $(W,I)$ as introduced in \cite{Soer07} and later used in \cite{WillSingular}.

\begin{defn}\label{defn-reff over field}
    Let $(W,S)$ be a Coxeter system, and $\kk$ be a field. A \emph{reflection faithful} representation of $W$ is a finite
dimensional representation $V$ of $W$ over $\kk$ such that:
\begin{itemize}
\item The representation is faithful;
\item 
We have $\codim V^w = 1$ if and only if $w$ is a reflection.
\end{itemize}
\end{defn}

\begin{lemma}\label{lemma reff implies (3')}
   Let char $\kk \neq 2$ and let $V$ be a reflection faithful representation of $(W,S)$. Then for any reflection $t \in W$, there exists a pair $(h_t,v_t)\in V^* \times V$ such that 
   \[t(\lambda)=\lambda - 2h_t(\lambda)v_t\]
   for all $\lambda\in V$. The pair $(h_t, v_t)$ is unique up to the equivalence relation $(h_t, v_t) \sim (a h_t, a^{-1} v_t)$ for units $a \in \kk^{\times}$.
\end{lemma}

\begin{proof}
When char $\Bbbk$ is not $2$, any involution acts diagonalizably with eigenvalues $1$ and $-1$.
When $V$ is a reflection faithful representation of $W$ and $t$ is a reflection, the second condition in \cref{defn-reff over field} then implies that $V$ has a one dimensional $-1$-eigenspace for $t$. Choose any non-zero element $v_t$ in the $-1$-eigenspace, and let $h_t \in V^*$ be the unique element which vanishes on $V^t$ and takes the value $1$ on $v_t$. By verifying the statement for $V^t$ and for $v_t$, one concludes that $t(\lambda)=\lambda - 2h_t(\lambda)v_t$. The uniqueness statement is straightforward.
\end{proof}

Of course, one should think of $(h_t, v_t)$ as a root and (half a) coroot attached to $t$.


Let $T$ denote the set of reflections in $W$. For $t,t'\in T$, fix some $h_t, h_{t'}$ as in the above lemma. If $h_t, h_{t'}$ are not linearly independent, then $V^t=\ker(h_t)=\ker(h_{t'})=V^{t'}$, so that $\codim V^{tt'}\leq 1$. 
Since $tt'$ is not a reflection, reflection faithfulness of $V$ implies that $V^{tt'}=V$, and we have that $t=t'$ by faithfulness of the representation $V$.
Hence we have the following lemma (\cite[Lemma 4.1.4]{WillSingular}):
\begin{lemma}\label{lemma coroots are linearly independent}
     Let char $\kk \neq 2$ and $V$ be a reflection faithful representation of $(W,S)$. The elements of $\set{h_t: t \in T} \subset V^*$ as above are pairwise linearly independent.
\end{lemma}

Theorems \ref{thm-categorification} and \ref{thm-hom formula} as proven in \cite{WillSingular} assume that the following properties hold.

\begin{property} \label{willassume} The following are properties of a finite-dimensional realization $\Lambda$ of $(W,S)$ over a field $\kk$. As usual, we regard $\Lambda^*$ as a $W$-module via the contragredient action, and let $R = \Sym_{\Bbbk}(\Lambda)$.

\begin{enumerate}
    \item $\Lambda^*$ is a reflection faithful representation of $(W,S)$, and  $\kk$ is an infinite field of characteristic $\ne 2$. \label{assumption 1}
    \item For all finitary $I \subset S$, $R$ is graded free over $R^I$, and one has an isomorphism of graded $R^I$-modules: 
    $$R \cong \Tilde{\pi}(I)R^I,$$
    where $\tilde{\pi}(I):= \displaystyle \sum_{w \in W_I}v^{-2\ell(w)}.$\label{assumption 2}
    \item One may choose a family of pairs $\{(h_t, v_t) \colon t \in T\}$ as in \cref{lemma reff implies (3')}, such that
    $$xh_s = h_t \text{ if } xsx^{-1} = t.$$
    \label{assumption 3}
\end{enumerate}
\end{property}

\begin{remark} The statement of \cref{willassume} \hyperref[assumption 1]{(1)} was phrased in terms of $\Lambda^*$ rather than $\Lambda$ to be consistent with Williamson's conventions. Regardless, when comparing with Lie theory, the roots (not the coroots) should be elements of $R$. It is not hard to verify that \hyperref[assumption 1]{(1)} holds for $\Lambda^*$ if and only if it holds for $\Lambda$. The same can be said for \hyperref[assumption 3]{(3)}, but we are not sure whether the same can be said for \hyperref[assumption 2]{(2)}. \end{remark}

\begin{remark} One can think of \cref{willassume} \hyperref[assumption 3]{(3)} as the statement that the realization has something like a ``non-integral root system.'' In \cite[Section 4.1]{WillSingular}, both \cref{willassume} \hyperref[assumption 1]{(1)} and \hyperref[assumption 2]{(2)} are stated as explicit assumptions, but \hyperref[assumption 3]{(3)} is merely asserted as a property which holds. This assertion is incorrect, so we view this property as a third assumption made by Williamson. \end{remark}

\begin{lemma} Property \hyperref[assumption 3]{(3)} fails for the deformed affine realization. \end{lemma}

\begin{proof} Suppose that \hyperref[assumption 3]{(3)} holds, and choose a family of elements $(h_t, v_t)$ for each reflection. Note that $h_{s_1}$ must be in the span of $x_1 - \zeta x_2$; by rescaling the entire family we can assume $h_{s_1} = x_1 - \zeta x_2$. Conjugating $s_1$ by $s_1 s_2$ yields, $s_2$, so we apply $s_1 s_2$ to $h_{s_1}$ and deduce that $h_{s_2} = \zeta x_1 - \zeta^2 x_2$. Continuing to apply $s_i s_{i+1}$ to $h_{s_i}$, we deduce that $h_{s_{n+1}} = \zeta^{n} x_1 - \zeta^{n+1} x_2$. But $s_{n+1} = s_1$, a contradiction. \end{proof}

Let us state a weaker property which a realization might possess. 

\begin{property} We continue the setup of \cref{willassume}.
\begin{enumerate}
    \item[(3')]\label{assumption 3'} 
    One may choose a family of pairs $\{(h_t, v_t) \colon t \in T\}$ as in \cref{lemma reff implies (3')}, such that
    $$xh_s \in \kk^\times \cdot h_t \text{ if } xsx^{-1} = t$$
    where $\kk^\times$ denotes the group of units in $\kk$.  
\end{enumerate}
\end{property}

\begin{lemma}  Property \hyperref[assumption 1]{(1)} implies Property \hyperref[assumption 3']{(3')}. \end{lemma}

\begin{proof} In fact, one can choose any pairs $(h_t, v_t)$ for each reflection $t$, and Property \hyperref[assumption 3']{(3')} will hold. By \cref{lemma reff implies (3')}, one need only show that $xh_s$ is a $-1$-eigenvector for $t$, which is obvious. \end{proof}

We shall see in section \ref{subsec-deformed realization} that our deformed affine realization satisfies assumptions \hyperref[assumption 1]{(1)}, \hyperref[assumption 2]{(2)}, and thus also \hyperref[assumption 3']{(3')}, when working over $\kk=\Q(\zeta)$ for an indeterminate $\zeta.$ 

Now we explain why $\hyperref[assumption 3']{(3')}$ is sufficient for the arguments in \cite{WillSingular} to go through.  
Property $\hyperref[assumption 3]{(3)}$ appears as eq $(4.1.3)$ in \cite{WillSingular} and is used only in the proofs of \cite[Definition/Proposition 4.3.2, Definition/Proposition 5.0.1, and Lemma 7.3.4]{WillSingular}. We recall these briefly below and explain exactly why $\hyperref[assumption 3']{(3')}$ is sufficient. Note that Williamson fixes the family $\{(h_t, v_t)\}$ satisfying $\hyperref[assumption 3]{(3)}$ once and for all; we too fix a family satisfying $\hyperref[assumption 3']{(3')}$.

\begin{prop*}\cite[Definition/Proposition 4.3.2]{WillSingular} \label{deprop-Wil 4.3.2}
Let $X \subset W$ be a finite subset. Consider the subspace
\begin{equation*}
R(X) = \left \{ f = (f_x) \in \bigoplus_{x \in X} R \; \middle | \; 
\begin{array}{c} f_x - f_{tx} \in (h_t) \\ \text{ for all } t \in T
\text{ and } x, tx \in X\end{array}  \right \} \subset \bigoplus_{x \in X} R.
\end{equation*}
Then $R(X)$ is a graded $k$-algebra under componentwise multiplication. 
It becomes an object of $\bMod{R}{R}$ if we define left and right actions of $r \in R$ via
\begin{align*}
 (rf)_x & = r \cdot f_x \\
 (fr)_x & = f_x\cdot x(r)
\end{align*}
for $f = (f_x) \in R(X)$ (where $\cdot$ on the RHS above is multiplication in $R$).
If a pair of subgroups $W_1, W_2 \subset W$ satisfy $W_1X = X = XW_2$ then $R(X)$ carries commuting left $W_1$- and right $W_2$-actions if we define
\begin{align*}
  (uf)_x &= u(f_{u^{-1}x}) & \text{for $u \in W_1$},\\
  (fv)_x &= f_{xv^{-1}} & \text{for $v \in W_2$}.
\end{align*}
\end{prop*}
\begin{remark}
    \cite[Proposition 4.3.2]{WillSingular} as stated also includes that if $X = \{ x \}$ is a singleton then $R(X) \cong R_x$, where $R_x$ denotes the standard bimodule as in \cite[Definition 4.2.1]{WillSingular}, and that if $X = \{x, y \}$ consists of two elements we write $R_{x,y}$ instead of $R(X)$. The former statement is a straightforward consequence of the definitions, and the latter is just notation.
\end{remark}
The proof of this proposition in \cite{WillSingular} only uses $\hyperref[assumption 3]{(3)}$ near the end of the proof to show that the $R(X)$ is stable under the left $W_1$-action. More explicitly, it is shown that if $x, tx \in X$ one has,
\begin{equation*}
  (wf)_x - (wf)_{tx} = w(f_{w^{-1}x} - f_{w^{-1}tx}) \in (w(h_{w^{-1}tw})) = (h_t)
\end{equation*}
where the last equality follows from $\hyperref[assumption 3]{(3)}.$ But the equality of the ideals $(w(h_{w^{-1}tw}))$ and $(h_t)$ also hold with the weaker assumption $\hyperref[assumption 3']{(3')}$. Thus the proposition holds under the assumption $\hyperref[assumption 3']{(3')}$ instead of $\hyperref[assumption 3]{(3)}.$
\vspace{5pt}

Now we proceed to discuss why \cite[Definition/Proposition 5.0.1]{WillSingular} still holds with $\hyperref[assumption 3']{(3')}$. We will continue using the setup in \cite[Definition/Proposition 4.3.2]{WillSingular} discussed above.
\begin{prop*}\cite[Definition/Proposition 5.0.1]{WillSingular}
Let $X, W_1, W_2 \subset W$ be as above.

\begin{enumerate}
\item For all reflections $t \in W_1$ there exists an operator
$f \mapsto \partial_tf$ on $R(X)$, the \emph{left Demazure operator to $t$}, uniquely determined by
\begin{equation*}
  f - tf = 2h_t(\partial_tf) \quad \text{for all $f \in R(X)$.}
\end{equation*}
This is a morphism in $\bMod{R^t}{R}$.

\item For all reflections $t \in W_2$ there exists an operator $f \mapsto f\partial_t$ on $R(X)$, the \emph{right Demazure operator to $t$}, uniquely determined by
\begin{equation*}
  f - ft = (f \partial_t)2h_t \quad \text{for all $f \in R(X)$.}
\end{equation*}
This is a morphism in $\bMod{R}{R^t}$.
\end{enumerate}
\end{prop*}
\begin{proof}
In \cite{WillSingular}, the left Demazure operator $f \mapsto \dd_tf$ is defined by
\begin{equation*}
    (\partial_tf)_x := \frac{f_x - tf_{tx}}{2h_t}, \qquad x \in X
  \end{equation*} 
and the right Demazure operator $f\mapsto f\dd_t$ is defined by
\begin{equation*}
  (f\partial_t)_x := \frac{f_x - f_{xt}}{2x(h_t)},\qquad x\in X.
\end{equation*}
The proof of part (1) of the proposition in \cite{WillSingular} (namely, that the left Demazure operator on $R(X)$ is well-defined, unique, and is a morphism in $\bMod{R^t}{R}$) is completely independent of the assumption $\hyperref[assumption 3]{(3)}$, so part (1) holds.
The explicit check that the right Demazure operator is well-defined as an operator on $R(X)$ is left to the reader in \cite{WillSingular}, so we explain the details of this check and show that it suffices to assume $\hyperref[assumption 3']{(3')}$ instead of $\hyperref[assumption 3]{(3)}$.  Unicity follows from the argument in \cite{WillSingular}.

Let us fix $f \in R(X)$ and $t\in W_2$. To see that $(f\dd_t)_x \in R$ (and not just in $R[x(h_t)^{-1}]$) for $x\in X$, we use that $f_x-f_{xt}=f_x-f_{xtx^{-1}x} \in (h_{xtx^{-1}})=(x(h_t))$. Note that this argument still holds with assumption $\hyperref[assumption 3']{(3')}$ instead of $\hyperref[assumption 3]{(3)}$, since all we need is the equality of the ideals $(h_{xtx^{-1}})$ and $(x(h_t))$. Now we are left to show that $f\dd_t$ in fact lives in $R(X)$ and not just in $\bigoplus_{x\in X}R$. 
In other words, we need to show for any reflection $t' \in W$ and $x \in X$ such that $t'x \in X$,
\begin{equation*}
    (f\dd_t)_x - (f\dd_t)_{t'x} \in (h_{t'}).
\end{equation*}
Adding and subtracting $t'((f\dd_t)_{t'x})$, we see that we need to show
\begin{equation}\label{eq 1 prop williamson demazure operators}
    (f\dd_t)_x-t'((f\dd_t)_{t'x}) - ((f\dd_t)_{t'x}-t'((f\dd_t)_{t'x})) \in (h_{t'}).
\end{equation}
Note that for any $g \in R$ which is $t'$-anti-invariant, $g$ vanishes on $V^{t'}$, so that $g \in (h_{t'})$ (since $(h_{t'})$ is the ideal of functions vanishing on the hyperplane $V^{t'} \subset V$).
In particular, $(f\dd_t)_{t'x}-t'((f\dd_t)_{t'x}) \in (h_{t'})$, so showing \cref{eq 1 prop williamson demazure operators} is equivalent to showing
\begin{equation}\label{eq 2 prop williamson demazure operators}
    (f\dd_t)_x-t'((f\dd_t)_{t'x}) \in (h_{t'}).
\end{equation}
Expanding using definitions and the fact that $t'$ is a reflection, and in particular an involution, we have 
\begin{align}
    (f\dd_t)_x-t'((f\dd_t)_{t'x}) &= \frac{f_x - f_{xt}}{2x(h_t)}-t'\left(\frac{f_{t'x} - f_{t'xt}}{2t'x(h_t)} \right) \nonumber \\
    &=\frac{f_x-t'(f_{t'x}) - (f_{xt}-t'(f_{t'xt}))}{2x(h_t)}. \label{eq-right demazure}
\end{align}

Let $X'=\set{x,xt,t'x,t'xt}$. 
Then $\set{1,t'}X'=X'=X'\set{1,t}$, and we have $\tilde{f}\in R(X')$ whose components are given by those of $f$.
We may then write the numerator above using the left Demazure operators (from part (1) of the proposition, for $\tilde{f} \in R(X')$) as
$$f_x-t'(f_{t'x}) - (f_{xt}-t'(f_{t'xt}))=2h_{t'}\left((\dd_{t'}\tilde{f})_x-(\dd_{t'}\tilde{f})_{xt}\right),$$
so that the numerator is in $(h_{t'}).$
So if $x(h_t)$ is coprime to $h_{t'}$, we have that the RHS of \cref{eq-right demazure} is in $(h_{t'})$ as desired.\\
Now suppose $x(h_t)$ is not coprime to $h_{t'}$. 
From assumption $\hyperref[assumption 3']{(3')}$, we have that $x(h_t) \in \kk^\times h_{xtx^{-1}}$, and reflection faithfulness implies that $\kk h_{xtx^{-1}}=\kk h_{t'}$ if and only if $xtx^{-1}=t'$ (see \cref{lemma coroots are linearly independent}).
Hence we have that $xtx^{-1}=t'$ in this case,
and $x(h_t)=\lambda h_{t'}$ for some $\lambda \in \kk^\times$. 
Then the RHS of \eqref{eq-right demazure} simplifies to $\lambda^{-1}\left((\dd_{t'}\tilde{f})_x-(\dd_{t'}\tilde{f})_{xt=t'x} \right)$ which is in $(h_{t'})$ since we already know that $\dd_{t'}\tilde{f} \in R(X')$ by part (1) of the proposition. 

The statement that $f \mapsto f\dd_t$ is a morphism in $\bMod{R}{R^t}$ follows from the observation that for $a\in R$ and $b \in R^t$,
\begin{align*}
(a(f\dd_t)b)_x=a\left(\frac{f_x - f_{xt}}{2x(h_t)} \right)x(b)= \frac{af_xx(b) - af_{xt}xt(b)}{2x(h_t)}=\frac{(afb)_x - (afb)_{xt}}{2x(h_t)}=((afb)\dd_t)_x.&\qedhere
\end{align*}
\end{proof}

Finally, we recall \cite[Lemma 7.3.4]{WillSingular} and show that it holds under assumption $\hyperref[assumption 3']{(3')}$ instead of $\hyperref[assumption 3]{(3)}$.
In this Lemma, $T$ is the set of reflections in $W$ and $p_-$ is the minimal coset representative of a $(W_I,W_J)$-coset $p$.
\begin{lemma*}\cite[Lemma 7.3.4]{WillSingular}
    The element 
    $$m_p=\prod_{t\in T: tp_-<p_-}h_t$$
    lies in $R^{I\cap p_-Jp_-^{-1}}$.
\end{lemma*}
\begin{proof}
    In the proof of the Lemma in \cite{WillSingular}, it is shown that for $s \in I \cap p_-Jp_-^{-1}$ and $t \in T$ such that $tp_-<p_-$, it follows that $stsp_-<p_-$. Note that this argument is purely within the Coxeter group $W$, and is independent of the representation $V$. With assumption $\hyperref[assumption 3]{(3)}$, it follows that $s(m_p)=m_p$ for all $s\in I\cap p_-Jp_-^{-1}$ since we can think of $m_p$ as a product of the $s$-invariant elements $h_th_{sts}$. We claim that $h_th_{sts}$ is still $s$-invariant with assumption \hyperref[assumption 3']{(3')}. Indeed, by \hyperref[assumption 3']{(3')}, we have that $s(h_{sts})=\lambda h_{t}$ for some $\lambda \in \kk^\times$. Then $h_{sts}=\lambda s(h_t)$, so that $s(h_t)=\lambda^{-1}h_{sts}.$ 
    It follows that 
    \begin{align*}s(h_th_{sts})=s(h_t)s(h_{sts})=\lambda^{-1}h_{sts}\lambda h_{t}=h_th_{sts}. &\qedhere 
    \end{align*}
\end{proof}

Hence we conclude that Theorems \ref{thm-categorification} and \ref{thm-hom formula} hold without requiring \hyperref[assumption 3]{(3)}.

\begin{theorem} Theorems \ref{thm-categorification} and \ref{thm-hom formula} hold when one only requires Properties \hyperref[assumption 1]{(1)} and \hyperref[assumption 2]{(2)} but not \hyperref[assumption 3]{(3)}. \end{theorem}

%

\subsection{The deformed affine Type A realization}
\label{subsec-deformed realization}
In this section, we discuss the deformed affine Type A realization (\cref{defn-deformed affine realization}) in more detail, and prove that Properties $\hyperref[assumption 1]{(1)}$ and $\hyperref[assumption 2]{(2)}$ from section \ref{subsec-categorification} hold when the ground field is $\Q(\zeta)$, so that we can use the categorification theorem and Hom formula from section \ref{subsec-categorification}.

In preparation for future work, we aim to establish this result in reasonable generality. For this purpose it will be helpful to think of realizations over integral domains that are not necessarily fields. Since both \cite{WillSingular}, \cite{Soer07} define ``reflection faithfulness" only over fields (see Definition \ref{defn-reff over field}), we give a generalized definition. There are many generalizations which will specialize to the definition in \cite{WillSingular, Soer07}, but this one is sufficient for our purposes.
\begin{defn}\label{defn-reff}
    $V$ is a \emph{reflection faithful} representation of $(W,S)$ over an integral domain $A$ if:
    \begin{enumerate}
        \item [(a)] $V$ is a faithful representation,
        \item [(b)] $V^w$ is free of finite rank over $A$ for all $w\in W$,
        \item [(c)] $V/V^w$ is free of finite rank over $A$ for all $w \in W$ (i.e., $V^w$ is a direct summand of $V$),
        \item [(d)] rk$_A(V/V^w)=1$ if and only if $w$ is a reflection in $W$.
    \end{enumerate}
\end{defn}

\begin{lemma}\label{reff flat base change}
    Reflection faithfulness is stable under flat base change, i.e., if $A\rightarrow B$ is a flat morphism of commutative domains and $V$ a reflection faithful representation of $(W,S)$ over $A$, then $V_B :=V\otimes_A B$ is a reflection faithful representation over $B$.
\end{lemma}
\begin{proof}
    Observe that for $w\in W$,
    \begin{align*}
        V^w &= \ker(V \xrightarrow{w-1}V), \\
        V/V^w &\cong\Ima(V\xrightarrow{w-1}V). \end{align*}
    Since $A \rightarrow B$ is flat, we have that for any map $\phi$, there are canonical isomorphisms $\ker(\phi)\otimes_A B \cong ker(\phi \otimes_A 1_B)$, $\Ima(\phi)\otimes_A B \cong \Ima(\phi \otimes_A 1_B)$. 
    We then have that
    \begin{align*}
        (V_B)^w &= \ker(V\otimes_A B \xrightarrow{(w-1)\otimes 1_B} V\otimes_A B) = V^w \otimes_A B,\\
        V_B/(V_B)^w &= \Ima(V\otimes_A B \xrightarrow{(w-1)\otimes 1_B} V\otimes_A B)=(V/V^w)\otimes_A B
    \end{align*}
    and (a), (b), (c) and (d) in Definition \ref{defn-reff} follow.
\end{proof}

The following lemma helps us compute fixed point sets. The product of a diagonal matrix and a permutation matrix is called a \emph{rescaled permutation matrix}.

\begin{lemma}\label{lemma rescaled permutation matrix fixed points} Let $\{x_1, \ldots, x_n\}$ be a basis of a free module $V$ over a commutative domain $\Bbbk$, and let $\rho \in GL(V)$ be a rescaled permutation matrix. Then both  the  fixed point set $V^{\rho}$ and the quotient $V/V^{\rho}$ are free. The rank of $V^{\rho}$ is the number of monodromy-free cycles of $\rho$ (as defined in the proof). \end{lemma}

\begin{proof} Let $\rho = w t$ where $w$ is a permutation matrix and $t$ a diagonal matrix. Let $\{1, \ldots, n\} = \coprod \mathcal{O}_k$ be the decomposition into cycles for $w$. As a module over $\rho$, $V = \bigoplus \Span \{x_i\}_{i \in \mathcal{O}_k}$. Each of the statements in the lemma can be proven individually for each cycle, which together proves the statement for the direct sum.
So, by restriction to a given cycle and by renaming our basis, we can assume that $w = (12 \cdots d)$ is a standard $d$-cycle, and that $t$ is the diagonal matrix with entries $(t_1, t_2, \ldots, t_d)$. 
Note that $d=1$ is allowed. Let $c = \prod_{i=1}^d t_i$. We call the cycle \emph{monodromy-free} if $c=1$. It is easy to verify that $\rho^d = c \cdot I$, which has no nonzero fixed points unless $c = 1$. Thus the fixed point set of $\rho$ is zero unless $c=1$.

Assume $c=1$. Then $\rho(x_1) = t_1 x_2$ and $\rho(x_2) = t_1 t_2 x_3$, etcetera. 
We claim that the vector
\[ \sum_{i=1}^d x_i(\prod_{j<i} t_j)\]
is a basis for the (rank one) $\rho$-fixed subspace, and that $\{x_2, \ldots, x_d\}$ descends to a basis for the quotient. We leave the straightforward argument to the reader. 
\end{proof}

Now we consider the deformed affine realization. For simplicity we fix $n \geq 3.$
We briefly recall the basics. 
Let $\Om =\Z/n\Z$, and $W = \tilde{S}_n$, with $S = S_{\aff} = \{s_i\}_{i \in \Om}$. Then
$\Lambda_{\zeta}$ is the realization of $(W,S)$ over $\Z[\zeta, \zeta^{-1}]$ with basis $\{x_i\}_{i \in \Om}$ with $\alpha_i=x_i - \zeta x_{i+1}$ and $\alpha_i^\vee= x_i^*- \zeta^{-1}x_{i+1}^*$, where $\{x_{i}^*\}$ is dual to $\{ x_i\}$. 
The action of $W$ was given by
$$ s_i(x_i)=\zeta x_{i+1}, \hspace{15pt} s_i(x_{i+1})=\zeta^{-1}x_i, \hspace{15pt} s_i(x_j)=x_j \;\text{  otherwise.}$$


Let us recall some well-known facts as a way to set notation and context.

\begin{lemma}\label{lemma reflections}
        We have an isomorphism $\Tilde{S}_n \xrightarrow{\sim} S_n \ltimes \Lambda_{\rt}$, realizing the affine Weyl group as affine transformations on the root lattice (of $\mathfrak{sl}_n$), where 
        \begin{align}
        &s_i \mapsto s_i \in S_n \text{ for } i\not=0 \\
        &s_0 \mapsto s_{\alpha_{long}}t_{-\alpha_{long}} = t_{\alpha_{long}} s_{\alpha_{long}}.
        \end{align}
When $\alpha$ is a root of $\mathfrak{sl}_n$, we let $s_{\alpha}$ denote the corresponding reflection in $S_n$ and $t_{\alpha}$ the corresponding translation (an element of $\Lambda_{\rt}$), and we let $t_{k\alpha} = t_{\alpha}^k$. Above, $\alpha_{long}$ is the longest root. Via this isomorphism, the reflections in $W$ get identified with
$\{s_{\alpha}t_{k\alpha}\}$ ranging over $k \in \Z$ and over all positive roots of $\mathfrak{sl}_n$.
\end{lemma}
    
\begin{proof}[Proof] 
The identification of $\tilde{S}_n$ with $S_n\ltimes\Lambda_{\rt}$ is standard and well-known; for example see  \cite[Chapter 4]{HumphreysCoxeter}. The classification of reflections is also well-known but not found in \cite{HumphreysCoxeter}; we sketch the proof.

For $\lambda \in \Lambda_{\rt}$, we shall use the notation $t_{\lambda}$ exclusively  to denote the corresponding element in $W$. Note that $S_n$ normalizes the subgroup $\{t_{\lambda}: \lambda \in \Lambda_{\rt}\}$. 
Explicitly, we have that for $w\in S_n$, $\lambda \in \Lambda_{\rt}$, $wt_{\lambda}=t_{w(\lambda)}w$ where $w(\lambda)$ refers to the usual action of $S_n$ on the root lattice. 
This will be used multiple times without explicit mention for the rest of this section.
    
Writing $w=w't_{\lambda}$ for $w' \in S_n$ and $\lambda \in \Lambda_{\rt}$, we compute that
       \begin{equation}
            ws_{\alpha}t_{k\alpha}w^{-1} = s_{w'(\alpha)}t_{(k-\langle \alpha^\vee, \lambda \rangle)w'(\alpha)}.
        \end{equation}

 In particular, any reflection $ws_iw^{-1}$ (for $w\in W,\; i \in \ZZ/n\ZZ$) is of the form $s_{\alpha}t_{k\al}$ for some $k \in \ZZ$ and some root $\al$ of $S_n$.
 Hence all reflections in $W$ are indeed of the desired form.
 Note that all the pairings here are for the geometric realization of $S_n$ coming from the root lattice of $\mathfrak{sl}_n$.

Suppose $n \ge 3$. Conversely, given a root $\al$ of $S_n$, we can find $w' \in S_n$ such that $w'(\al) = \al_1$ whence $w' s_{\al} (w')^{-1} = s_1$. Now we use the fact that $n \ge 3$, so that $s_1(\al_2) = \al_1 + \al_2$, and compute (for any $k \in \Z$) that
\begin{equation}  (t_{-k \al_2} w')(s_{\al} t_{k \al}) (t_{-k \al_2} w')^{-1} = t_{-k\alpha_2}s_1t_{k\alpha_1}t_{k\alpha_2} = s_1. \end{equation}
Thus any element of the form $s_{\al} t_{k \al}$ is a reflection.

The proof is different when $n=2$ since the reflections are not all conjugate, and we leave it to the reader. \end{proof}


\begin{lemma}\label{lemma-action matrix}
    Identify $GL_{\Z[\zeta, \zeta^{-1}]}(\Lambda_{\zeta})$ with $GL(n,\Z[\zeta, \zeta^{-1}])$ by fixing the basis $\{x_i\}$, and let $\rho: W \rightarrow GL(n,\Z[\zeta, \zeta^{-1}])$ be the map obtained from the realization $\Lambda_{\zeta}$. Then,
    \begin{enumerate}
        \item For $w \in S_n \subset W$, $\rho(w)$ is a rescaled permutation matrix with non-zero entries of the form $\zeta^k$ for $k \in \Z.$ 
        For a reflection (i.e., transposition) $w= (i,j) \in S_n$ with $i < j$, $\rho(w)$ is explicitly given by 
        \[\rho(w)=
        \left({\begin{array}{ccccccccccc}
        \Id_{i-1} &    & & &  \\
           &0  &0  & \zeta^{i-j} &\\
           &0    &\Id_{j-i-1} &0   &\\
           &\zeta^{j-i}   &0   &0   & \\
           &    &   &   &\Id_{n-j}\\
           \end{array} }  \right)
        \]
         where $\Id_k$ denotes the identity matrix of size $k$ (and unspecified entries of the matrix are all 0).
        \item Conjugation with $\rho(w), \; w \in S_n$ is an honest permutation on diagonal matrices, which permutes the diagonal entries according to the permutation $w$.
        \item For a positive root $\alpha$ associated to a transposition $(i,j)$ in $S_n$ with $i<j$, let $t_{\alpha}$ be the corresponding translation. Then we have that \[ \rho(t_{\alpha})= 
          \left({\begin{array}{ccccccccccc}
        \Id_{i-1}   \\
           &\zeta^{n} \\
           &    &\Id_{j-i-1}\\
           &    &   &\zeta^{-n} & \\
           &    &   &   &\Id_{n-j}\\
           \end{array} }  \right).
        \]
        More generally, for $\lambda=\sum_{i=1}^{n-1}\lambda_i\al_{s_i} \in \Lambda_{\rt}$, we have that $\rho(t_{\lambda})$ is a diagonal matrix given by
        \[ \rho(t_{\lambda})= 
          \left({\begin{array}{ccccccccccc}
        \zeta^{n\lambda_1} \\
            &\zeta^{n(\lambda_2-\lambda_1)}\\
            &   &\zeta^{n(\lambda_3-\lambda_2)} \\
            &   &   &\ddots \\
            &   &   &   &\zeta^{n(\lambda_{n-1}-\lambda_{n-2})}\\
            &   &   &   &   &\zeta^{-n\lambda_{n-1}}
           \end{array} }  \right).
        \]
        \item Any element of $W$ acts by a rescaled permutation matrix.
    \end{enumerate}
\end{lemma}
\begin{proof}
    (1) The fact that $\rho(s_i)$ coincides with the given matrix is immediate from the definition of the $W$-action in the deformed affine realization. 
    By writing any $w\in W$ as a composition of the simple reflections $s_i$ for $1\leq i \leq n-1$, we see that $\rho(w)$ is a generalized permutation matrix as desired. 
    
    Now, observe that as an $S_n-$module (for the finite Weyl group $S_n \subset W$), we can identify $\Lambda_{\zeta}$ with the permutation representation $\Z[\zeta,\zeta^{-1}][y_1, \ldots ,y_n]$ of $S_n$ over the ring $\Z[\zeta,\zeta^{-1}]$, by setting $y_k=\zeta^{k-1}x_k$ for $1\leq k \leq n$. 
    Then we see that for a transposition $w=(i,j) \in S_n$ with $i<j,$ 
    $$\rho(w)(x_k)=\rho(w)(\zeta^{1-k}y_k)= 
    \begin{cases}
    \zeta^{1-i}y_j=\zeta^{j-i}x_j \;\text{  if $k=i$}\\
    \zeta^{1-j}y_i=\zeta^{i-j}x_i \;\text{  if $k=j$}\\
    \zeta^{1-k}y_k=x_k \hspace{18pt}\text{  otherwise}
    \end{cases},$$
   so that
    \[\rho(w)=
        \left({\begin{array}{ccccccccccc}
        \Id_{i-1} &    & & &  \\
           &0  &0  & \zeta^{i-j} &\\
           &0    &\Id_{j-i-1} &0   &\\
           &\zeta^{j-i}   &0   &0   & \\
           &    &   &   &\Id_{n-j}\\
           \end{array} }  \right)
        .\]
    
    (2) follows from (1), since we can then write $\rho(w)$ as a permutation matrix times a diagonal matrix for any $w\in S_n$.
    
    (3) We know that $s_0=t_{\al_{long}}s_{\al_{long}}$ acts via the matrix
    \[\left({\begin{array}{ccccccccccc}
             
         0  &0 &\zeta\\
         0  &\Id_{n-2} &0\\
         \zeta^{-1} &0 &0  \\
           \end{array} }  \right).
        \]
    But $s_{\al_{long}}$ acts via  the matrix 
    \[\left({\begin{array}{ccccccccccc}
             
         0  &0 &\zeta^{1-n}\\
         0  &\Id_{n-2} &0\\
         \zeta^{n-1} &0 &0  \\
           \end{array} }  \right),
        \]
    so that 
    \[\rho(t_{\al_{long}})=\rho(s_0)\rho(s_{\al_{long}})= 
    \left({\begin{array}{ccccccccccc}
             
         \zeta^{n} &0 &0\\
         0  &\Id_{n-2} &0\\
         0 &0 &\zeta^{-n}  \\
           \end{array} }  \right).
        \]
        Now we use that any positive root $\alpha$ is $S_n$-conjugate to $\alpha_{long}$ so that there is some $w \in S_n$ such that $w(\alpha_{long})=\alpha$. 
        Then $t_\alpha=t_{w(\alpha_{long})}=wt_{\alpha_{long}}w^{-1}$. 
        By part (2), we see that $\rho(t_{\al})$ is then a diagonal matrix whose diagonal entries are given by permuting the diagonal entries of $\rho(t_{\al_{long}}).$ 
        If $\alpha$ corresponds to a transposition $(i,j)$ with $i<j$, then since $\alpha_{long}$ corresponds to the transposition $(1,n)$, we may choose $w=(1,i)(n,j) \in S_n$ so that $ws_{\al_{long}}w^{-1}=(i,j)=s_{\al}$ to get the desired description of $\rho(t_{\alpha})$ using part (2) of the lemma. The general description of $\rho(t_{\lambda})$ for $\lambda \in \Lambda_{\rt}$ follows.
        
    (4) Any element of $W$ is a product of an element of $S_n$ and a translation, and a product of rescaled permutation matrices is a rescaled permutation matrix. 
\end{proof}

\begin{lemma}\label{lemma freeness zeta infinite order} Let $\Lambda_A$ be a specialization of $\Lambda_{\zeta}$ as in \cref{defn-deformed affine realization}. For any element $y \in W$, the fixed point set $\Lambda_A^y$ and the quotient $\Lambda_A/\Lambda_A^y$ are both free modules over $A$.  The rank of $\Lambda_A^y$ over $A$ agrees with the rank of $\Lambda_{\zeta}^y$ over $\Z[\zeta, \zeta^{-1}]$ so long as $\zeta \in A^{\times}$ has infinite order. \end{lemma}

\begin{proof}
By \cref{lemma-action matrix}, $y$ acts by a rescaled permutation matrix. 
From \cref{lemma rescaled permutation matrix fixed points}, we deduce the freeness of $V_A^y$ and $V_A/V_A^y$. The entries of $t$ are powers of $\zeta$, so the element $c$ for any cycle is a power of $\zeta$. By the infinite order assumption, $c = 1$ in $A$ if and only if $c=1$ in $\Z[\zeta, \zeta^{-1}]$. \end{proof}


\begin{theorem}\label{reff thm} 
    The realization $\Lambda_{\zeta}$ is a reflection-faithful representation of $W$. 
    Let $A$ be any $\Z[\zeta, \zeta^{-1}]$-algebra for which the unit $\zeta \in A$ has infinite order. Then $\Lambda_A$ is reflection faithful as well.
    
\end{theorem}

\begin{proof}  

    Faithfulness of $V_A$ is an immediate consequence of \cref{lemma-action matrix} and \cref{lemma rescaled permutation matrix fixed points} (for $w\in S_n, \lambda \in \Lambda_{\rt}$, the only way $wt_\lambda$ acts as the identity matrix is if $w=1$ and $\lambda=0$, by the explicit description of the corresponding matrices in \cref{lemma-action matrix} and the argument in \cref{lemma rescaled permutation matrix fixed points}).

    We are left to show that $V_A$ is reflection faithful, i.e., we have to show that (b),(c) and (d) in Definition \ref{defn-reff} hold.
    \cref{lemma freeness zeta infinite order} shows us that (b) and (c) hold, so now we only have to show that (d) holds. By \cref{lemma rescaled permutation matrix fixed points} the fixed point set of any simple reflection $s$ has codimension $1$. The same is true for any reflection, as the fixed point set of $w s w^{-1}$ is $w(\Lambda_A^{s})$. 
    
    
    
    
    On the other hand, suppose we have $w=w't_{\lambda} \in W$ with $w'\in S_n, \; \lambda \in \Lambda_{\rt}$ such that  $\rk_A(\Lambda_A^w)=n-1$.
    By \cref{lemma rescaled permutation matrix fixed points}, $w'$ needs to have at least $n-1$ cycles, so either $w'$ is a reflection in $S_n$, or $w'$ is the identity. If $w'$ is the identity, so that $w = t_{\lambda}$, then we need exactly $n-1$ diagonal entries in $\rho(t_{\lambda})$ (see \cref{lemma-action matrix}) to be $1$, so that all but one cycle is monodromy-free. But $\rho(t_{\lambda})$ has determinant $1$, so the last entry is also $1$, a contradiction.
    
     Thus $w'$ is a reflection. Suppose $n \ge 3$. By further conjugating $w$ (and changing the meaning of $\lambda$) we may assume that $w = s_1 t_{\lambda}$. By \cref{lemma rescaled permutation matrix fixed points}, all $n-1$ cycles of $s_1$ must be monodromy-free. In particular, all diagonal entries in $\rho(t_{\lambda})$ must be $1$ except the first two, from which one deduces that $\lambda$ is a multiple of $\al_1$.  Thus $w = s_1 t_{k \al_1}$ for some $k \in \Z$, and $w$ is a reflection. 

    We leave to the reader the analogous proof when $n=2$.
\end{proof}

\begin{lemma}\label{lemma reff for V*} Let $A$ be any $\Z[\zeta, \zeta^{-1}]$-algebra for which the unit $\zeta \in A$ has infinite order. Then $\Lambda_A^*$ is reflection faithful as well. \end{lemma}

\begin{proof} The $W$-action on $\Lambda^*$ is the same as the $W$-action on $\Lambda$ up to replacing $\zeta$ with $\zeta^{-1}$ and vice versa. 
This follows immediately by checking that if $\{x_i^*\} \subset \Lambda^*$ denotes the dual basis to $\{x_i\} \subset \Lambda$, then 
$$
s_i(x_k^*)=
\begin{cases}
 \zeta^{-1}x_{i+1}^* \text{ if k=i}\\
 \zeta x_i^* \hspace{19pt}\text{ if k= i+1} \\
 x_k^* \hspace{23 pt}\text{ otherwise}.
\end{cases}
$$
(See \cite[Proposition 2.15]{EJY1} for the formal notion of a $\zeta-$antilinear isomorphism.)\end{proof}

Having concluded the discussion of reflection faithfulness, we discuss the other critical assumption.

\begin{lemma} \label{lem:assume2} Assumption  \hyperref[assumption 2]{(2)} holds for $\Lambda_{\zeta}$, and for any specialization $\Lambda_A$. \end{lemma}

\begin{proof} Fix a finitary subset $I\subset S$, and a maximal finitary subset $S'$ containing $I$. 
Then $S'$ is of the form $S\setminus\{k\}=\hh{k}$ for some $0\le k< n.$ 
We can then identify $\Lambda_A$ with the permutation representation $A \cdot \{y_1, \ldots y_n\}$ of  $W_{S'} \cong S_n$ by choosing $y_1=x_{k+1}, y_2=\zeta^{-1} x_{k+2}, \ldots, y_n=\zeta^{-(n-1)}x_{k+n}.$
It is well-known that Assumption \hyperref[assumption 2]{(2)} holds for the permutation representation of $S_n$ even over $\Z$ (and hence after base change to $A$). This is effectively a result of Demazure \cite{Demazure} (for a modern exposition in English see \cite[Theorem 4.2]{EKLP}). For example, one can find a basis for $R$ over $R^{\hh{0}}$ by applying Demazure operators to the polynomial $y_1^{n-1} y_2^{n-2} \cdots y_{n-1}$. \end{proof}

\begin{theorem}
Assumptions \hyperref[assumption 1]{(1)} and \hyperref[assumption 2]{(2)}  from section \ref{subsec-categorification} hold for the realization $\Lambda_{A}$ so long as $A$ is an infinite field of characteristic $\ne 2$, and $\zeta$ has infinite order in $A^{\times}$. Hence Theorems \ref{thm-categorification} and \ref{thm-hom formula} hold for $\SSBim$ when defined using $\Lambda_{A}$.
        
\end{theorem}
\begin{proof}
Assumption \hyperref[assumption 1]{(1)} follows from \cref{lemma reff for V*}.
Assumption \hyperref[assumption 2]{(2)} follows from \cref{lem:assume2}.
\end{proof}

\begin{remark}\label{rmk-Abe SSBim}
 Noriyuke Abe in \cite{Abe} generalizes the construction of \cite{WillSingular} by starting with weaker assumptions, including cases where $V$ is not a faithful representation of $(W,S)$ but is a faithful representation for all finitary subgroups, and allowing $\kk$ to be a complete local Noetherian integral domain.
 Let $Q$ denote the fraction field of $R$, and $Q^I$ the $W_I$-invariant subring. 
 Now we define $\SSBim_{\Abe}$  as follows: instead of just working with $(R^I,R^J)$-bimodules $M$, one has the extra data of certain $(Q^I,Q^J)$-bimodules $(M_Q^x)_{x\in W_I\backslash W/W_J}$ and isomorphisms $Q^I\otimes_{R^I}M\otimes_{R^J}Q^J\simeq \bigoplus_{x\in W_I\backslash W/W_J}M_Q^x $ which satisfy certain conditions. See \cite[section 2]{Abe} for further details of this construction. 
 It is shown then in \cite{Abe} that Theorems \ref{thm-categorification} and Theorems \ref{thm-hom formula} hold under his weaker assumptions for $\SSBim_{\Abe}$. 
 Abe is careful not to rely on the balanced assumption. We expect Abe's construction to be useful when studying specializations of $\Lambda_{\zeta}$ where $\zeta$ is a root of unity, as such specializations are no longer faithful.
\end{remark}

\subsection{Reverse Bott-Samelsons}\label{subsec reverse bott-samelsons}

The goal of this section is to prove that the 2-functor $\GS$ is essentially surjective after taking the Karoubi envelope. A brief pedagogical interlude on ordinary Bott-Samelson bimodules will set the stage. We work with a realization satisfying Assumptions \hyperref[assumption 1]{(1)} and \hyperref[assumption 2]{(2)}.

The traditional category of Bott-Samelson bimodules $\BSBim$ (these are $(R,R$)-bimodules) is the full subcategory of $\SBSBim(\emptyset, \emptyset)$ only containing objects associated to singular expressions of the form 
\begin{equation} \label{ordinaryexp} [[\emptyset \subset \{s_1\} \supset \emptyset \subset \{s_2\} \supset \cdots \subset \{s_d\} \supset \emptyset]] .\end{equation}
Thus $\BSBim$ and $\SBSBim(\emptyset, \emptyset)$ are not the same. 
However, these two monoidal categories have the same Karoubi envelope, which is not obvious, but follows from the classification of indecomposables in each category. This is stated as \cite[Theorem 24.12]{EMTW} but the proof is not explained, so let us elaborate.

\begin{theorem} \label{thm:sameKar1} The Karoubi envelopes of (the graded, additive completions of) $\BSBim$ and $\SBSBim(\emptyset, \emptyset)$ are equivalent. \end{theorem} 

\begin{proof} Since $\BSBim \subset \SBSBim(\emptyset, \emptyset)$ we need only prove that every indecomposable summand of a singular Bott-Samelson bimodule in $\SBSBim(\emptyset, \emptyset)$ is isomorphic to a summand of an ordinary Bott-Samelson bimodule. By \cref{thm-categorification}, the isomorphism classes of indecomposable objects in $\SBimod{I}{J}$ are in bijection with double cosets $p \in W_I \backslash W / W_J$, and appears as a direct summand of $\BS(IK_{\bullet})$ whenever $IK_{\bullet}$ is a (singular) reduced expression for $p$.

When $I = J = \emptyset$, $p = \{w\}$ is a singleton, and any ordinary reduced expression for $w$ gives rise to singular reduced expression as in \eqref{ordinaryexp} (though not all singular reduced expressions are of this form). Consequently, the indecomposable bimodule $B_p$ is a summand of an ordinary Bott-Samelson bimodule (it would traditionally be called $B_w$).
\end{proof}

Ultimately, the basic combinatorial fact underlying this theorem is: every double coset $p \in W_{\emptyset} \backslash W / W_{\emptyset}$ has a singular reduced expression of the form \eqref{ordinaryexp}.

Ordinary Bott-Samelson bimodules bounce between minimally singular and subminimally singular.
In this paper, we mostly deal with the other extreme. Recall the definition of a maximal double coset from \cref{defn:maximalcoset}.

\begin{defn}
Let $(W,S)$ be an affine Weyl group equipped with a realization satisfying Assumptions \hyperref[assumption 1]{(1)} and \hyperref[assumption 2]{(2)}.
The $2$-category of \emph{maximally singular Soergel bimodules}, $\mSSBim$ denote the 1-full sub-$2$-category\footnote{So the only objects of $\mSSBim$ are $\hh{s}$ for $s \in S$, and for each $s, t \in S$ we have ${}_{\hh{s}} \mSSBim_{\hh{t}} = \SBimod{\hh{s}}{\hh{t}}$.} of $\SSBim$ whose objects are the maximal finitary subsets $\hh{s}$ of $S$. 
The $2$-category of \emph{reverse Bott-Samelson bimodules} $\mSBSBim$ is the 2-full sub-$2$-category\footnote{So the 1-morphisms are restricted, but the 2-morphisms are the same.} of $\mSSBim$ containing 1-morphisms of the form
\begin{equation} \label{reverseexp} [[\hh{a_1} \supset \hh{a_1 a_2} \subset \hh{a_2} \supset \hh{a_2 a_3} \subset \cdots \supset \hh{a_{d-1} a_d} \subset \hh{a_d}]].\end{equation}
To avoid redundancy we assume that $a_j \ne a_{j+1}$ for all $1 \le j < d$.
\end{defn}

\begin{theorem}\label{thm mSBSBim vs mSSBim in type A} In affine type $A$, the Karoubi envelope of (the graded, additive completion of) $\mSBSBim$ is equivalent to $\mSSBim$. Moreover, every expression of the form \eqref{reverseexp} is a reduced expression. \end{theorem}

\begin{proof} Since $\mSBSBim \subset \mSSBim$, we need to prove that every indecomposable singular Bott-Samelson bimodule in $\SBimod{\hh{s}}{\hh{t}}$ is a summand of a reverse Bott-Samelson bimodule (for all $s, t \in S$). Just as in the proof of Theorem \ref{thm:sameKar1}, Williamson's categorification theorem reduces this problem to proving that every maximal double coset has a singular reduced expression of the form \eqref{reverseexp}. Every expression \eqref{reverseexp} is $\GS(\ula,a_d)$ for some $\ula$, whence it is reduced by \cref{thm:GSrex}. Fixing the parabolic subgroup $\hh{a_d}$ on the right, the bijections \eqref{bijectionsa} show that any maximal double coset has the form $\psi_{a_d}(\lambda)$ for some $\lambda \in \Lambda^+_{\wt}$, and \cref{thm:GSrex} provides many reduced expressions of the desired form (associated to rexes $\ula$ for $\lambda$). 
\end{proof}

\begin{remark} The Bott-Samelson bimodules for the reduced expressions from \cref{thm:GSrex} are the images under geometric Satake of tensor products of fundamental representations. The analogous statement to \cref{thm mSBSBim vs mSSBim in type A} on the other side of the Soergelified Satake equivalence is that $\Rep(\slf_n)$ is equivalent to the Karoubi envelope of $\Fund(\slf_n)$, the full subcategory monoidally generated by fundamental representations. The corresponding statement under geometric Satake states that the corresponding reverse (singular) Bott-Samelson variety gives a resolution of singularities of a Schubert variety associated to the corresponding dominant weight within the affine Grassmannian. \end{remark}

%% file: Frobenius_Data_and_Calculations_with_Demazure_Operators.tex
We recall the notion of Frobenius extensions, Demazure operators, and Frobenius hypercubes. These will be used in the next chapter to construct a diagrammatic category as in \cite{EWSFrob}. 
This chapter ends with various calculations involving Demazure operators for the deformed affine realization. These calculations will be used in section \ref{sec well-defined-ness of the functor}.
\subsection{Cubes of Frobenius Extensions} \label{subsec- cubes of frob}


\begin{defn} A \emph{graded (commutative) Frobenius extension of degree d} is an extension of commutative graded rings $\iota :  A \hookrightarrow B$, where $B$ is a finitely generated free graded $A$-module, equipped with an
$A$-linear map $\dd = \dd_A^B : B \to A$ of degree $-2d$, called the \emph{(Frobenius) trace}. In other words, $\dd$ is a homomorphism of graded $A$-modules $\dd: B \rightarrow A(-2d)$. The trace is required to be \emph{non-degenerate}, meaning the following map is an isomorphism:
$$B \rightarrow \Hom_A(B,A): b \mapsto (b'\mapsto \dd(bb')).$$

\end{defn}
\begin{remark}\label{rmk-non-degeneracy and dual basis}
The non-degeneracy of $\dd$ is equivalent
to requiring that $B$ admits a \emph{pair of dual bases}. This is a pair of bases
$\{x_i\}$ and $\{y_i\}$ for $B$ as an $A$-module, such that $\dd(x_i y_j)=\delta_{i j}$. The equivalence between these two notions is familiar from linear algebra. 
Within $\Hom_A(B,A)$ there is a dual basis $\{x_i^*\}$ to the basis $\{x_i\}$, and one takes $y_{i}$ to be the preimage of $x_{i}^*$ under $B \to \Hom_A(B,A)$. 
If $\dd$ is non-degenerate then any basis $\{x_{i}\}$ admits a unique dual basis.

If we start with a homogeneous basis $\{x_i\}\subset B$, then the dual basis $\{y_i\}$ is also homogeneous. After all, if $\{y_{i}\}$ is the dual basis to $\{x_{i}\}$, then the appropriate homogeneous components of $\{y_{i}\}$ are also dual to $\{x_{i}\}$, but the dual basis is unique.

\end{remark}

The data of a Frobenius extension also equips one with a comultiplication map $\Delta^B_A : B \to B \otimes_A B$, which is a $B$-bimodule map.  It sends $1 \mapsto \sum_\alpha x_i \otimes
y_i$, an element which is independent of the choice of dual bases. If $A \hookrightarrow B$ is a graded Frobenius extension of degree $d$, then the comultiplication map $\Delta^B_A$ has degree $2d$.  

\begin{defn}\label{defn-product coproduct element}
    The \emph{product-coproduct} element $\mu^B_A$ of a Frobenius extension $A \hookrightarrow B$ is defined to be 
    $$\mu_A^B:=m^B_A\circ\Delta^B_A (1)$$
    where $m^B_A:B\otimes_AB \rightarrow B$ denotes the multiplication map of $B$.
\end{defn}

\begin{lemma}
The element $\partial^B_A(\mu^B_A) \in A$ is a scalar, equal to the rank (not the graded rank) of $B$ as an $A$-module. \end{lemma}

\begin{proof} This follows directly from the definition of dual bases. \end{proof}

Being a Frobenius extension is a structure (the choice of $\dd$), not just a property. However, the choice of $\dd$ is unique up to a choice of a unit in $B$. 

\begin{lemma} \label{lem Frobenius structures unique up to unit}
Suppose $\dd, \dd': B \rightarrow A$ are two Frobenius traces. Then there exists a unit $b_0 \in B$ such that $\dd'(b) = \dd(b_0 b)$ for all $b \in B$. If $\mu$ and $\mu'$ represent $\mu^B_A$ for these two choices of Frobenius structure, then $\mu = b_0 \mu'$.
    
Thus the product-coproduct element is independent of the choice of Frobenius trace up to unit, and if $\mu = \mu'$ then $\dd = \dd'$. \end{lemma} 

\begin{proof} Via the isomorphism 
    $$B \rightarrow \Hom_A(B,A): b \mapsto (b' \mapsto \dd(bb')),$$ we have that $\dd' \in \Hom_A(B,A)$ corresponds to some $b_0\in B$, so that $\dd'(b)=\dd(b_0b)$ for any $b \in B$. We leave the reader to verify that $b_0$ is a unit, and the remaining statements.
\end{proof}


Of note is the case when $A \hookrightarrow B$ induces an equality $A^\times = B^\times$, because then $\dd' = a_0 \dd$ for $a_0 \in A^\times$. We say $\dd'$ is a rescaling of $\dd$. For example, in this paper it happens frequently that $A$ and $B$ are $\kk$-algebras and $A^\times=B^\times=\kk^\times$.

Note that if $A \subset B$ and $B \subset C$ are Frobenius
extensions, then $A \subset C$ is a Frobenius extension as well, with trace $\dd^C_A = \dd^B_A \dd^C_B$. It is proven in \cite[(2.5)]{EWSFrob} that \begin{equation} \label{eq:mumult} \mu^C_A = \mu^C_B \mu^B_A. \end{equation}

\begin{defn} A \emph{hypercube of (graded) Frobenius extensions} or a \emph{Frobenius hypercube} will be the following datum. \begin{itemize} \item
A finite set $\Gamma$. We also use $\Gamma$ to designate the entire datum. We consider the hypercube with vertices labelled by subsets of
$\Gamma$. An edge in this hypercube corresponds to $I \setminus \gamma \subset I$ for some $\gamma \in I \subset \Gamma$, and parallel
edges correspond to the same $\gamma$. \item A (contravariant) assignment of graded rings $R^I$ to vertices in the cube, so that $I \subset J
\implies R^J \subset R^I$. \item For each edge, a trace map $\dd^{I \setminus \gamma}_I : R^{I \setminus \gamma} \to R^I$ making $R^I \hookrightarrow R^{I
\setminus \gamma}$ into a graded Frobenius extension. 

\end{itemize}
We call this hypercube \emph{compatible} if, for every square $I \subset
J,J' \subset K$ in the hypercube, we have $\dd^J_K \dd^I_J = \dd^{J'}_K \dd^I_{J'}$. In this case,
there is a well-defined map $\dd^I_J : R^I \to R^J$ for every $I \subset J$, which endows the extension $R^J \subset R^I$ with the
structure of a Frobenius extension. We assume henceforth that every hypercube of Frobenius extensions is compatible.

We usually refer to $R^{\emptyset}$ as just $R$. Similarly we write $\dd_I := \dd^\emptyset_I$.


A \emph{partial hypercube of Frobenius extensions} or a \emph{partial Frobenius hypercube} is a finite set $\Gamma$ and a non-empty collection $\mathcal{F}$ of subsets of $\Gamma$ (i.e. vertices in the hypercube), together with the data of graded rings and Frobenius extensions as above but only for vertices within the subset $\mathcal{F}$. We require that $J \in \mathcal{F}$ and $I \subset J$ implies $I \in \mathcal{F}$.  Note that $\emptyset$ is always in $\mc{F}$. For the purposes of this paper, we shall use the term Frobenius hypercube to also to refer to partial Frobenius hypercubes.
\end{defn}

For a Frobenius hypercube and $I \subset J \in \mathcal{F}$, let $\mu^I_J \in R^I$ denote the product-coproduct element for the Frobenius extension $R^J \subset R^I$, and $\mu_I = \mu^{\emptyset}_I$. By \eqref{eq:mumult} we have
\begin{equation} \label{eq:mumult2} \mu^I_J = \frac{\mu_J}{\mu_I}.\end{equation}
It is straightforward to see that \eqref{eq:mumult2} is equivalent to compatibility.

For sake of simplicity, assume that all rings in the hypercube are $\Bbbk$-algebras and that their unit groups agree with $\Bbbk^\times$.  By \cref{lem Frobenius structures unique up to unit}, each individual Frobenius extension $I \subset J$ in the cube is uniquely determined up to a scalar, and this scalar is visible in the choice of $\mu^I_J$. When defining a new Frobenius hypercube structure by rescaling the various elements $\mu^I_J$, constraints are placed on these scaling factors by \eqref{eq:mumult2}. Ultimately, one can freely rescale the elements $\mu_I \in R$ (resp. the trace maps $\partial_I$) for each $I \in \mathcal{F}$, and there will be a unique way to extend this rescaling to the entire hypercube so that compatibility will be maintained.

Fix a Coxeter system $(W,S)$ with a realization, and let $R$ be as in \cref{subsec- SSBim}. 
Let $\Gamma = S$, and $\mathcal{F}$ be the collection of finitary subsets of $S$. Under some assumptions on the realization, there will be a Frobenius hypercube with rings $R^I$ associated to each $I \in \mathcal{F}$. For example, when $W$ is a Weyl group with its geometric realization over $\mathbb{R}$, there is a standard Frobenius hypercube structure where $\mu_I$ is the product of the positive roots for the parabolic subgroup $W_I$. However, as we have noted in \S\ref{subsec-balanced}, some realizations do not have a well-behaved notion of positive roots. Consequently, the choice of product-coproduct elements $\mu_I$, or equivalently the choice of trace maps $\partial_I$, will be subtle. To construct the trace maps, we need to study Demazure operators.
 

\subsection{Demazure operators}
\label{subsec-demazures}

For the rest of this section we fix a Coxeter system $(W,S)$, with a realization $(\Lambda,  \Lambda^\vee,\{\alpha_s\},\{\alpha_s^\vee\})$ over an integral domain $\kk$. We have the graded $\kk$-algebra $R:= \Sym_{\kk} \Lambda$, where $\Lambda$ is in degree 2, with its natural $W$-action.

We define \emph{the Demazure operator} or \emph{the divided difference operator} $\dd_s$ associated to a simple reflection $s\in S$ to be a $\kk$-linear operator on $R$:
$$\dd_s(f)= \frac{f-s(f)}{\al_s} \text{\hspace{20 pt}for $f\in R$}.$$
Note that $\dd_s$ may also be defined as the unique $\kk-$linear operator on $R$ satisfying the following properties:
\begin{itemize}
    \item $\dd_s|_{\Lambda}$ coincides with the functional defined by $\al_s^\vee$, i.e., $\dd_s(\lambda)=\langle\al_s^\vee,\lambda\rangle$ for $\lambda \in \Lambda$.
    \item $\dd_s$ satisfies the twisted Leibniz rule:
    \begin{equation}\label{eq twisted leibniz}
        \dd_s(fg)=\dd_s(f)g+s(f)\dd_s(g).
    \end{equation}
\end{itemize}

\begin{lemma}\label{lemma demazure operators coxeter relations}
    (\cite[Claim A.7]{ECathedral}) The Demazure operators $\dd_s, \dd_t$ associated to simple reflections $s,t \in S$ satisfy the following relations.
    \begin{enumerate}
        \item $\dd_s\dd_s=0$.
        \item For $m_{st} < \infty$, 
    \begin{equation} \label{eq:rescaled braid} \lambda \dd_s \dd_t \dd_s \cdots = \dd_t \dd_s \dd_t \cdots \end{equation} (both expressions alternate and have length $m_{st}$) for some unit $\lambda \in \kk^\times$. Moreover, $\lambda=1$ when $m_{st}$ is even.
    \end{enumerate}
The last property states that the braid relations hold up to scalar. This unit $\lambda$ agrees with the unit $\lambda$ in \eqref{eq:lambdadef}, so that $\lambda = 1$ when the realization is balanced.
\end{lemma}

\begin{remark} Let us fill in more details of the sketched proof from \cite[Claim A.7]{ECathedral}. We can write each $\dd_s$ (and compositions thereof) as a $Q$-linear combination of elements of $W$, where $Q$ is the fraction field of $R$. For example, when $m_{st} = 3$, we have
\begin{equation} \dd_s \dd_t \dd_s(f) = \frac{(-1)^3 sts(f)}{\alpha_s s(\alpha_t) st(\alpha_s)} + \ldots, \quad \dd_t \dd_s \dd_t(f) =\frac{(-1)^3 sts(f)}{\alpha_t t(\alpha_s) ts(\alpha_t)} + \ldots.\end{equation}
More generally, the coefficient of $w_{s,t}$ in $\dd_s \dd_t \cdots$ (the leading term) will be $(-1)^{m_{st}}$ times the reciprocal of the product of the $s$-aligned choice of roots as in \cite[Definition A.6]{ECathedral}, whereas for $\dd_t \dd_s \cdots$ one instead uses the $t$-aligned choice of roots. One can compare these leading coefficients using \eqref{eq:lambdadef}, and verify that \eqref{eq:rescaled braid} holds for the leading terms. With more combinatorial work one could verify this for all terms. \end{remark}

\begin{remark}
When $W_{s,t}$ acts faithfully and the Chevalley-Shephard-Todd theorem holds, one can bypass this combinatorial work by arguing that any two $R^{s,t}$-linear operators of degree $-2m$ are colinear, and that elements of $W_{s,t}$ are linearly independent over $Q$.

Also, note that both sides of \eqref{eq:rescaled braid} are zero when the action of $W_{s,t}$ is not faithful. If $W_{s,t}$ acts trivially, all Demazure operators are zero. Otherwise the action will factor through a smaller dihedral group, where \eqref{eq:rescaled braid} holds by induction, and then $\dd_s \dd_s=0$ implies the statement.
\end{remark}



\begin{definition}\label{defn-demazure operator associated to an expression}
For an expression $\un{w}=(s_1,\ldots, s_d)$ of $w\in W$, we define the Demazure operator $\dd_{\un{w}}$ associated to it by setting
$$\dd_{\un{w}}:= \dd_{s_1}\ldots \dd_{s_d}.$$
\end{definition}
When the expression $\un{w}$ is not reduced, \cref{lemma demazure operators coxeter relations} implies that $\dd_{\un{w}}$ is 0. \cref{lemma demazure operators coxeter relations}(2),(3)  together with Matsumoto's theorem implies that  $\dd_{\un{w}}$ is independent of the particular choice of a reduced expression $\un{w}$ of $w$ up to a unit in $\kk.$ We pause to note a conjecture, which to our knowledge is not in the literature at this level of generality (when the braid relations hold only up to scalar).

\begin{conjecture} If $W$ acts faithfully and $\un{w}$ is a reduced expression, then $\dd_{\un{w}}$ is nonzero. \end{conjecture}

We now record how Demazure operators behave for the deformed affine realization. Until the end of this subsection, $(W,S)$ is the affine Weyl group, where $S$ is in bijection with $\Omega = \Z/n\Z$, and $R_{\zeta}$ denotes the graded ring $\Sym_{\Z[\zeta,\zeta^{-1}]}(\Lambda_{\zeta})$ with $\Lambda_{\zeta}$ in degree 2.



\begin{lemma}\label{lemma demazure operators computation deformed affine}
For $i \in \Om$, let $\dd_i$ denote the corresponding Demazure operator on $R_{\zeta}$.
\begin{enumerate}
    \item For $i \in \Om$, $\dd_i(x_j)=\begin{cases}
        1 &\text{ if $j=i$,}\\
        -\zeta^{-1} &\text{ if $j=i+1$,}\\
        0 &\text{ otherwise.}
    \end{cases}$
    \item For $i,j \in \Om$ such that $s_is_j=s_js_i$, $\dd_i \dd_j= \dd_j \dd_i$.
    \item For $i \in \Om$, $\zeta \dd_{i}\dd_{i+1}\dd_{i}=\dd_{i+1}\dd_i\dd_{i+1}.$
    \end{enumerate}
\end{lemma}
\begin{proof}
    This is straightforward from \cref{lemma demazure operators coxeter relations}. The most interesting part is (3), where the unit $\zeta = \lambda$ agrees with \eqref{eq:lambdadef} when $s = s_i$ and $t = s_{i+1}$, given that
    \begin{equation} \zeta s_i(\alpha_{i+1}) = s_{i+1}(\alpha_i). \end{equation}
    Alternatively, one can compute the unit directly by evaluating both sides on $x_i^2x_{i+1}$.

\end{proof}

Recall that we had defined the rotation operator $\sigma: \Omega \rightarrow \Omega: i \mapsto i+1$  and extended it to an operator on $\Lambda_{\zeta}$ such that $\sigma(\zeta)=\zeta$ and $\sigma(x_i)=x_{\sigma(i)}=x_{i+1}$ (see \cref{defn-rotation operator}). Hence $\sigma$ naturally extends to a ring automorphism $\sigma: R\xrightarrow{\sim} R$ which preserves the grading.

\begin{lemma} We have
\begin{equation} \label{eq rotational symmetry deformed affine} \sigma \circ \dd_i \circ \sigma^{-1} (f) = \dd_{\sigma(i)}(f). \end{equation}
\end{lemma}

\begin{proof} One has $\sigma(s_{k}(x_{j}))=s_{k+1}(x_{j+1})=s_{\sigma(k)}(\sigma(x_j))$ for any $k, j\in \Om$, so that $\sigma\circ s_k \circ \sigma^{-1}=s_{\sigma(k)}$ as operators on $R$.
It follows that
\[ \sigma\dd_i \sigma^{-1}(f)=\sigma \left(\frac{\sigma^{-1}(f)-s_i(\sigma^{-1}(f))}{x_i-\zeta x_{i+1}} \right)
=\frac{f-s_{i+1}(f)}{x_{i+1}-\zeta x_{i+2}}=\dd_{\sigma(i)} (f). \]
\end{proof}
We refer to this as \emph{rotational symmetry} or \emph{color symmetry} in the later sections.

The action of $\sigma$ on $S=\Om$ also extends naturally to an action on the set of words in $S$ and on $W$. 
Given an expression $\un{w}:=(s_{i_1}, \ldots s_{i_k})$ for $w=s_{i_1}\ldots ,s_{i_d}$, we have $\sigma(\un{w}):=(s_{\sigma(i_1)}, \ldots ,s_{\sigma(i_d)})$ as an expression for $\sigma(w):=s_{\sigma(i_1)}\ldots s_{\sigma(i_d)}$. 
Then \cref{eq rotational symmetry deformed affine} implies that 
\begin{equation}
    \sigma\dd_{\un{w}}=\dd_{\sigma(\un{w})}\sigma.
\end{equation}

\begin{remark}\label{rmk symmetry for rotation operator}
    In the special case when $f \in R$ is such that $\dd_{\un{w}}(f)$ is of degree 0, we have that $\sigma(\dd_{\un{w}}(f))=\dd_{\un{w}}(f)=\dd_{\sigma(\un{w})}(\sigma(f)).$
\end{remark}

Recall that we had defined the reflection operator $\tau: \Om \rightarrow \Om: k \mapsto -k$, and extended it to an operator on $\Lambda_{\zeta}$ 
such that $\tau(\zeta)=\zeta^{-1}$ and $\tau(x_i)=x_{n-i+1}$ (see \cref{defn-reflection operator}). This extends to a ring automorphism $\tau:R \xrightarrow{\sim}R$.

\begin{lemma} We have
\begin{equation} \label{eq:reflection symmetry for reflection operator} \tau \circ \dd_i \circ \tau^{-1}(f) = - \zeta \dd_{\tau(i)}(f). \end{equation}
\end{lemma}

\begin{proof}
It is straightforward to check that $\tau(s_k(x_i))=s_{n-k}(\tau(x_i))$ for any $i, k\in \Om$, i.e., $\tau\circ s_k\circ\tau^{-1}=s_{\tau(k)}$.
It follows then that 
\[\tau\dd_i \tau^{-1}(f)= \tau \left( \frac{\tau^{-1}(f)-s_{i}(\tau^{-1}(f))}{x_i-\zeta x_{i+1}}
\right) = \frac{f-s_{n-i}(f)}{x_{n-i+1}-\zeta^{-1}x_{n-i}}=-\zeta \dd_{n-i}(f). \] 
\end{proof}

The action of $\tau$ on $S=\Om$ extends naturally to an action on the set of words in $S$, and on $W$. 
Given an expression $\un{w}:=(s_{i_1}, \ldots s_{i_d})$ for $w=s_{i_1}\ldots s_{i_d}$, we have $\tau(\un{w}):=(s_{\tau(i_1)}, \ldots s_{\tau(i_d)})$ as an expression for $\tau(w):=s_{\tau(i_1)}\ldots s_{\tau(i_d)}$. 
Then we have that 
\begin{equation}
    \tau\dd_{\un{w}}=(-\zeta)^{d}\dd_{\tau(\un{w})}\tau.
\end{equation}

\subsection{Frobenius hypercube realizations and the deformed affine setting}
\label{subsec-frob realizations}
\label{subsec-frob realization deformed affine}

The representation of $W$ produces the ring $R$; the additional data provided by a realization (roots and coroots) pins down the Frobenius trace $R \to R^s$. Let us pin down all the Frobenius traces.

\begin{defn}\label{defn frobenius hypercube realization} A \emph{Frobenius (hypercube) realization} or \emph{Frealization} is a realization equipped with a choice of a compatible partial Frobenius hypercube, with $\Gamma = S$, $\mathcal{F}$ the collection of finitary subsets of $S$, and $R^I$ the subring of $W_I$-invariants in $R$ for $I\in \mc{F}$.
In other words, we have a
Frobenius structure $\dd^J_I$ on $R^I \subset R^J$ for each finitary pair $J \subset I \subset
S$ such that $\dd^J_I \dd^K_J = \dd^K_I$ whenever $K \subset J$.
(In particular, we assume that such a Frobenius structure exists.) 
We also require that 
\begin{itemize}
\item The Frobenius trace $R \to R^s$ must simply be $\dd_s$, as determined by the realization. 
\item The Frobenius trace $\dd_I : R \to R^I$ must agree with $\dd_{\un{w}_I}$, up to a unit in $\kk$, for some reduced expression $\un{w}_I$ of $w_I$.
\end{itemize}
\end{defn}

\begin{remark} The second condition, that $\dd_I$ and $\dd_{\un{w}_I}$ agree up to unit, is automatic by \cref{lem Frobenius structures unique up to unit} so long as $\dd_{\un{w}_I}$ is nondegenerate. A prerequisite for $\dd_{\un{w}_I}$ to be nondegenerate is for it to be surjective as a map $R \to R^I$, which is called \emph{generalized Demazure surjectivity} \cite[Definition 3.5]{EKLP}.\end{remark}

\begin{remark} It seems as though the choice of $\dd_s$ is fixed whereas the choice of $\dd_I$ is an invertible scalar by \cref{lem Frobenius structures unique up to unit}. In truth, $\dd_s$ can also be rescaled by rescaling the choice of roots and coroots (i.e. changing the underyling realization). An alternate approach would be to omit the realization from the data, and recover it from the Frobenius extension structure for $R^s \subset R$. \end{remark}

\begin{example}\label{example nice realizations are frobenius}
If the realization $\Lambda$ is balanced (see \S\ref{subsec-balanced}) then we unambiguously define $\dd_w$ to be $\dd_{\un{w}}$ for any reduced expression $\un{w}$ of $w$. 
In this case, we can then define $\dd^\emptyset_I:= \dd_{w_I}$ for a finitary $I\subset S$, and 
$$\dd^J_I:=\dd_{w_Iw_J^{-1}}$$
for $J \subset I \subset S$. 
Compatibility will hold. If the realization is also faithful and satisfies generalized Demazure surjectivity, then $\dd^J_I$ is nondegenerate, and this data yields a Frealization. This is proven in \cite[Theorem 4.3]{EKLP} under these exact assumptions on the realization (see \cite[Section 3.1]{EKLP}).
\end{example}




Our goal for the rest of this section is to define a Frobenius realization extending the deformed affine realization. Once again, $(W,S)$ is the type $A$ affine Weyl group with our standard conventions.

\begin{defn}
We say that a nonempty subset $I \subset S$ is \emph{connected} if $I=\set{i, i+1, \ldots i+k}$ for some integers $i$ and $k$ (considered modulo $n$). If $I\subset S$ is not connected, we say it is \emph{disconnected}. 
A \emph{connected component} of a disconnected subset is a maximal connected subset.
\end{defn}

\begin{defn}\label{defn deformed affine frobenius realization}
We make the following choices for the deformed affine realization $\Lambda_{\zeta}$ (see Definition \ref{defn-deformed affine realization}) defined over $\Z[\zeta,\zeta^{-1}]$, letting $R = R_{\zeta}$.
\begin{itemize}
    \item We fix a system of \textbf{positive roots} $\Phi_I^+$ for each subset $I\subset S$. 
    \\If $I=\{s_1, s_2, \ldots , s_{d-1}\}$, we choose our positive roots to be 
    $$x_i - \zeta^{j-i}x_j \hspace{20pt} \text{for } 1\leq i<j\leq d.$$
    If $I$ is connected and $I = \sigma^l(J)$ for some $l \in \Z$, then we set $\Phi^+_I = \sigma^l(\Phi^+_J)$.
    \\ For $I$ disconnected, we take a disjoint union of the positive roots for its connected components.
    \\ Note that for $I \subset J$ finitary, $\Phi^+_I \subset \Phi^+_J$.
    \item Now for finitary $I \subset S$, we define $\mu_I$ to be 
    $$\mu_I:= \prod_{\al \in \Phi_I^+} \al.$$
    We will use these $\mu_I$ when defining the Frobenius trace map $\dd_I$ so that the $\mu_I$ coincide with the product-coproduct element of the Frobenius extension $R^I \hookrightarrow R$ thus defined.
    We also define $$\mu_I^J:= \prod_{\al \in \Phi_I^+\setminus \Phi_J^+}\al\;\;,$$
     and this will coincide with the product-coproduct element for the Frobenius extension $R^I \hookrightarrow R^J$. 
     Note that for $K\subset  J \subset I$, we have that $\mu^K_I=\mu^J_I\mu^K_J.$
    \item Now we fix Frobenius trace maps as follows.
    \\For $I \subset S$ finitary and any reduced expression $\un{w}_I$ of the longest element $w_I$, we define $\dd_I$ to be the (unique) rescaling of $\dd_{\un{w}_I}$ so that $\dd_I(\mu_I)= |W_I|.$  Moreover, for $J \subset I \subset S$, we define $\dd_I^J$ to be the unique rescaling of $\dd_{\un{w_Iw_J^{-1}}}$ (for any reduced expression $\un{w_Iw_J^{-1}}$ of $w_Iw_J^{-1}$) such that $\dd^J_I(\mu_I^J)=|W_I|/|W_J|$.
    
\end{itemize}
\end{defn}

In fact, in \cref{subsec:stdrex} we will choose a particular reduced expression $\un{w}_I^{\std}$ such that $\dd_I = \dd_{\un{w}_I^{\std}}$ on the nose.

\begin{lemma}\label{lemma deformed affine frobenius realization is frobenius} The above definition yields a well-defined Frealization. \end{lemma}

\begin{proof} 
Fix $k \in \Omega$. Combining \cref{lemma-deformed affine restricted to Sn} with the symmetry $\sigma$, we can identify the restriction of the deformed affine realization to $W_{\hh{k}}$ with the permutation representation  $V_{\perm}=\bigoplus_{i=1}^n \Z[\zeta,\zeta^{-1}] y_i$ of $S_n$ via 
\begin{equation}\label{eq identifying deformed affine restricted to W-k-hat with permutation}
y_1=x_{k+1}, y_2=\zeta x_{k+2}, \ldots ,y_{n-1}=\zeta^{n-2} x_{k-1}, y_n=\zeta^{n-1} x_k.     
\end{equation}
Through this identification, the simple reflections $s_{k+i}$ (for $1 \leq i \leq n-1)$ act as simple transpositions $t_i$ swapping $y_i$ and $y_{i+1}$.
The permutation realization has simple roots $\alpha_{t_i} = y_i-y_{i+1}$. The identification does not preserve the roots and coroots, and we have
\begin{equation} \alpha_{k+i} \mapsto \zeta^{-(i-1)} \alpha_{t_i}, \quad \dd_{k+i} \mapsto \zeta^{i-1} \dd_{t_i}. \end{equation}
The other roots in $\Phi^+_{\hh{k}}$ are also sent to scalar multiples of positive roots in $\{y_i - y_j\}_{i > j}$. To disambiguate from $\dd_{k+i}$, we write $\tilde{\dd}_{k+i}$ for $\dd_{t_i}$.

Let $I \subset \hat{k}$ be any parabolic subgroup, and let $\un{w}_I$ be any reduced expression for $w_I$. The permutation realization (over $\Z$ or $\Z[\zeta, \zeta^{-1}]$) is a balanced realization and satisfies generalized Demazure surjectivity, so the discussion of \cref{example nice realizations are frobenius} applies. Thus the composition $\tilde{\dd}_{\un{w}_I}$ is a non-degenerate map $R \to R^I$, and it is well-known (see \cite{Demazure}) that the corresponding product-coproduct element is the product $\tilde{\mu}_I$ of the positive roots for $I$ (a subset of $\{y_i - y_j\}_{i < j}$). Since $\dd_{\un{w}_I}$ agrees with $\tilde{\dd}_{\un{w}_I}$ up to a power of $\zeta$, then $\dd_{\un{w}_I}$ is also nondegenerate. By \cref{lem Frobenius structures unique up to unit}, its product-coproduct element agrees up to a power of $\zeta$ with $\tilde{\mu}_I$, and thus also with $\mu_I$.

Consequently, $\dd_{\un{w}_I}(\mu_I) = \zeta^l |W|$ for some $l \in \Z$. We can define $\dd_I$ as above, and it will be non-degenerate. Its product-coproduct element will be a scalar multiple of $\mu_I$, and the normalization condition $\dd_I(\mu_I) = |W_I|$ implies that it is $\mu_I$ on the nose. Moreover, $\dd_I$ is independent of the choice of $\un{w}_I$, again by \cref{lem Frobenius structures unique up to unit}. Similar statements can be made about $J \subset I \subset \hat{k}$ and the relative longest element $w_Iw_J^{-1}$.

We check the other requirements of \cref{defn frobenius hypercube realization}. Since $\mu_s = \al_s$, we have  $\dd^{\emptyset}_s =\dd_s$ as required. The compatibility condition \eqref{eq:mumult2} follows from the definition. \end{proof}

\begin{remark}\label{rmk dual basis from permutation realization} 
The identification \cref{eq identifying deformed affine restricted to W-k-hat with permutation} may be used for constructing dual basis for the Frobenius extensions in the deformed affine realization, by rescaling from the permutation realization. In \cref{subsec-dual basis} we pick dual bases explicitly.
\end{remark}

{
}

Before continuing to explore this Frealization, we check some technical properties which were introduced in \cite{EWSFrob} to ensure proper behavior of its diagrammatic category.

\begin{defn}\cite[Definition 1.7]{EWSFrob}\label{defn-condition star} We say that a Frobenius hypercube satisfies condition $\star$ if, for every $I \ne J \in \mc{F}$ with $I \cup J \in \mc{F}$, setting $K = I \cap J$ and $L = I \cup J$, we may choose a basis $\{x_\alpha\}$ of $R^K$ over $R^I$ such that $x_\alpha \in R^J$. \end{defn}
Note that in \cite{EWSFrob}, this is stated a bit differently, where it is required to have a pair of dual bases for $R^K$ over $R^I$ where one of the bases lives in $R^J$. 
But if we have any basis of $R^K$ over $R^I$ already in $R^J$, then we can find its dual basis by Remark \ref{rmk-non-degeneracy and dual basis}. The point to emphasize is that typically the dual basis will not live within $R^J$.

\begin{defn}\cite[Definition 1.8]{EWSFrob} \label{defn- no mu zero deivisors} We say that the Frobenius hypercube $\Gamma$ has \emph{no $\mu$-zero
    divisors} if $\mu_I \in R$ is not a zero divisor for all $I$ in the hypercube $\Gamma$.
\end{defn}

\begin{defn}\cite[Definition 1.10]{EWSFrob}\label{defn-condition R3} We say that a Frobenius hypercube satisfies condition R3, if for every $I \in \mc{F}$ and distinct  $i, j, k \in \Gamma \setminus I$ with $I \sqcup\set{i,j,k} \in \mc{F}$, we have that ${\mu_{Iij} \mu_{Iik} 
\mu_{Ijk} \mu_I} \mid {\mu_{Iijk} \mu_{Ii} \mu_{Ij}\mu_{Ik}}$.
\end{defn}

\begin{remark}
Note that in \cite{EWSFrob}, the phrasing for condition R3 is slightly different. Translated to our notation, \cite[Definition 1.10]{EWSFrob} asks that
${\mu^I_{Iij} \mu^I_{Iik} 
\mu^I_{Ijk}} \mid {\mu^I_{Iijk} \mu^I_{Ii} \mu^I_{Ij}\mu^I_{Ik}}$,
which is equivalent to the condition mentioned above, when there are no $\mu$-zero divisors.

As mentioned in \cite[Remark 1.11]{EWSFrob}, when $R$ is a UFD and 
 $\mu^I_{Iij}, \mu^I_{Ijk}, \mu^I_{Iik}$
 are relatively prime, the condition R3 automatically holds. 
\end{remark}

\begin{lemma}\label{lemma deformed affine frobenius satisfies extra conditions}
   The deformed affine realization equipped with the Frobenius structure in \cref{defn deformed affine frobenius realization} satsfies the conditions in 
   \cref{defn-condition star}, \cref{defn- no mu zero deivisors}, \cref{defn-condition R3}.
\end{lemma}

\begin{proof} Clearly this Frealization has no $\mu$-zero divisors, since each ring $R^I$ is a domain and each $\mu_I$ is nonzero. Given that each $\mu_I$ is defined as a product of a certain set of roots, condition R3 follows immediately from the inclusion-exclusion principle.

For condition $\star$, we use similar arguments as in the proof of \cref{lemma deformed affine frobenius realization is frobenius}, passing to the permutation realization by restricting to a maximal finitary parabolic subset containing $I \cup J$. The property $\star$ for the permutation realization over any commutative domain follows from \cite[Theorem 4.15, Example 4.17]{EKLP}. The assumptions for \cite[Theorem 4.15]{EKLP} are that the realization is faithful, balanced, and satisfies generalized Demazure surjectivity (see \cite[Section 3.1]{EKLP}), all of which are known to hold for the permutation realization.

\end{proof}


\subsection{Choosing reduced expressions} \label{subsec:stdrex}

Now we find the promised reduced expression for $w_I$, for which the corresponding operator $\dd_{\underline{w}_I}$ agrees with $\dd_I$ on the nose. However, when $J \subset I$ it will typically be impossible to find a reduced expression for $w_I w_J^{-1}$ for which the corresponding operator agrees with $\dd^J_I$ on the nose, see \cref{rmk:wontworkforIJ}.



\begin{lemma}\label{lemma- demazure standard rex}
    We have $\dd_{\{1,2, \ldots ,k\}}=\dd_1\dd_2\ldots \dd_k \dd_{\{1,2,\ldots , k-1\}}$ for any $k<n.$ \\Hence for the reduced expression
    $$\un{w}^{\std}_{\{1,2, \ldots,k\}}:=(s_1, s_2, \ldots, s_k, s_1,s_2, \ldots , s_{k-1}, \ldots, s_1,s_2,s_1)$$
    of $w_{\{1,2, \ldots,k\}}$, we have that
    $\dd_{\{1,2, \ldots ,k\}}=\dd_{\un{w}^{\std}_{\{1,2, \ldots,k\}}}$.
\end{lemma}
\begin{proof}
    It suffices to show that $\dd_1\dd_2\ldots \dd_k \dd_{\{1,2,\ldots , k-1\}}$ satisfies the normalization condition, i.e., 
    $$\dd_1\dd_2\ldots \dd_k \dd_{\{1,2,\ldots , k-1\}}(\mu_{\{1,\ldots,k\}})= (k+1)!\;$$
    where $\mu_{\{1,2,\ldots , k\}}$ is as defined in \cref{defn deformed affine frobenius realization}.
    
    \noindent We proceed by induction on $k$. The statement is true for the case $k=1$, since $\dd_1(x_1-\zeta x_2)=2$. For the general case, observe that    
    $$\mu_{\{1,\ldots,k\}} = \mu_{\{1,\ldots,k-1\}}(x_1-\zeta^kx_{k+1})(x_2-\zeta^{k-1}x_{k+1})\ldots (x_k -\zeta x_{k+1}).$$
    Since $(x_1-\zeta^kx_{k+1})(x_2-\zeta^{k-1}x_{k+1})\ldots (x_k -\zeta x_{k+1})$ is invariant for $W_{\{1,\ldots, k-1\}}$, we have
     \begin{align*}
     \dd_1\dd_2\ldots &\dd_k \dd_{\{1,2,\ldots , k-1\}}(\mu_{\{1,\ldots,k\}}) \\
     &= \dd_1 \dd_2 \ldots \dd_k\left(
     (x_1-\zeta^kx_{k+1})(x_2-\zeta^{k-1}x_{k+1})\ldots (x_k -\zeta x_{k+1})
          \dd_{\{1,2,\ldots , k-1\}}(\mu_{\{1,2,\ldots , k-1\}})
     \right)\\
     &=\dd_1\dd_2\ldots \dd_k\left((x_1-\zeta^kx_{k+1})(x_2-\zeta^{k-1}x_{k+1})\ldots (x_k -\zeta x_{k+1})\right)k!,
     \end{align*}
     so it suffices to show that
     \begin{equation} 
     \dd_1\dd_2\ldots \dd_k\left((x_1-\zeta^kx_{k+1})(x_2-\zeta^{k-1}x_{k+1})\ldots (x_k -\zeta x_{k+1})\right)=k+1.\label{eq 1 lemma-demazure standard rex}
     \end{equation}
    Observe that for any $s \in S$ $$\dd_s(f\al_s)= \frac{f\al_s+s(f)\al_s}{\al_s}=f+s(f).$$
    Taking $s=s_k$, we may simplify the LHS of \cref{eq 1 lemma-demazure standard rex} to
    \begin{align}
     \dd_1\dd_2\ldots &\dd_{k-1}\left(
     \begin{array}{c}
     (x_1-\zeta^kx_{k+1})(x_2-\zeta^{k-1}x_{k+1})\ldots (x_{k-1} -\zeta^2 x_{k+1})\vspace{5 pt}\\
     +\quad s_k((x_1-\zeta^k x_{k+1})(x_2-\zeta^{k-1}x_{k+1})\ldots (x_{k-1} -\zeta^2 x_{k+1}))
     \end{array}
     \right)  
     \label{eq 2 lemma-demazure standard rex}
    \end{align}
    For the first term, using that $x_{k+1}$ is invariant under $s_i$ for $1 \le i \le k-1$, we get
    \begin{align*}\dd_1\dd_2\ldots \dd_{k-1}&((x_1-\zeta^kx_{k+1})(x_2-\zeta^{k-1}x_{k+1})\ldots (x_{k-1} -\zeta^2 x_{k+1}))\\
    &=\dd_1\dd_2\ldots \dd_{k-1}(x_1x_2\ldots x_{k-1}) + x_{k+1}(\text{something of negative degree})\\
    &=1+0=1.
    \end{align*}
    For the second term, by the inductive hypothesis,
    $$\dd_1\dd_2\ldots \dd_{k-1}\left((x_1-\zeta^{k-1} x_{k})(x_2-\zeta^{k-2}x_{k})\ldots (x_{k-1} -\zeta x_{k})\right)=k.$$
    So the RHS of \cref{eq 2 lemma-demazure standard rex} further simplifies to $1+k$ as desired.
\end{proof}


\begin{defn}\label{defn-standard rex}
We define a \emph{standard rex}  for an arbitrary $I \subset S$ as follows.

When $I = \{1, 2, \ldots, k\}$ we define $\un{w}_I^{\std}$ as in the previous lemma, i.e.
\begin{equation}  \un{w}_{\{1, 2, \ldots,k\}}^{\std}:=(s_1,s_2,\ldots ,s_k)\bullet \un{w}^{\std}_{\{1,2, \ldots,k-1\}} \end{equation}
    where $\bullet$ is used for concatenation. Obviously $\un{w}_{\{1\}}^{\std}:=(s_1)$, the unique reduced expression.

When $I = \sigma^{l}(\{1, \ldots, k\})$, let $\un{w}_I^{\std} := \sigma^l(\un{w}_{\{1, \ldots, k\}}^{\std})$.

For $I$ disconnected,  let $I_1, \ldots, I_k$ denote the connected components of $I$ (in any order). Then define
\begin{equation} \un{w}_I^{\std} := \un{w}_{I_1}^{\std}\bullet \un{w}_{I_2}^{\std} \bullet \ldots \bullet \un{w}_{I_k}^{\std}. \end{equation}

\end{defn}

\begin{remark}
When $I$ is disconnected, the chosen order on components does change the expression $\un{w}_I^{\std}$ but not the operator $\dd_{\un{w}_I^{\std}}$, since $\dd_i$ and $\dd_j$ commute whenever $i$ and $j$ belong to different components.

\end{remark}

    Now we explore other reduced expressions.
    
\begin{defn} \label{defn-rex weight} For a finitary $I\subset S$, we define a function  $\mc{L}_I$ called the \emph{rex weight} on the set of reduced expressions of $w_I$. 
    First suppose $I=\set{1,2, \ldots, k}$ for some $k<n$.
    Let $\mathcal{L}_I'$ be the function on the collection of reduced expressions for $w_I$ such that 
    $$\mathcal{L}_I'((s_{i_1},s_{i_2},\ldots, s_{i_m}))=\sum_{j=1}^{m}i_j,$$
     and let $\mathcal{L}_I$ be the function such that 
    $$\LL_I((s_{i_1},s_{i_2},\ldots, s_{i_m}))= \LL_I'((s_{i_1},s_{i_2},\ldots, s_{i_m}))-\LL_I'(\un{w}_I^{\std}).$$
    In other words, $\LL_I$ counts the number of times we need to apply the the braid move $(s_i,s_{i+1},s_i)=(s_{i+1},s_i,s_{i+1})$ (with signs, i.e., +1 when going from $(s_i,s_{i+1},s_i)$ to $(s_{i+1},s_i,s_{i+1})$, and -1 the other way around) to get from the standard rex to the given reduced expression.

    For an arbitrary finitary $I \subset S$, we then define $\LL_I$ as follows.
    If $I$ is such that $\sigma^l(I)=\set{1,2, \ldots, k}$, we set $$\LL_I(\un{w}_I):=\LL_{\sigma^l(I)}(\sigma^l(\un{w}_I)).$$
    If $I$ is disconnected, with connected components $I_1, \ldots ,I_m$, we set
    $$\LL_I(\un{w}_I):=\sum_{i=1}^m\LL_{I_i}(\un{w}_{I_i})$$
    where $\un{w}_{I_j}$ are chosen such that
    $\un{w}_{I_1}\bullet \ldots \bullet \un{w}_{I_m}$ is a reduced expression of $w_I$ obtained from $\un{w}_I$ by only using braid moves
    of the form $(s_i,s_j)\sim(s_j,s_i)$.
\end{defn}

\begin{defn}\label{defn-rex weight2}
    For a subset $J$ of $I$ and a reduced expression $\un{w}_I^J$ of $w_Iw_J^{-1}$, we define the \emph{rex weight} $\LL_I^J(\un{w}_I^J)$ of $\un{w}_I^J$ to be $\LL_I^J(\un{w}_I^J):=\LL_I(\un{w}_I^J \bullet \un{w}_J^{\std})$.
\end{defn}
    
    Now \cref{lemma demazure operators computation deformed affine} (3) and Lemma \ref{lemma- demazure standard rex} quickly imply the following.
    
\begin{corollary}\label{corollary-normalization}
    For any $I\subset S$ finitary, and any reduced expression $(s_{i_1},s_{i_2},\ldots, s_{i_m})$ for $w_I$, we have
    $$\partial_I=\zeta^{-\LL_I((s_{i_1},s_{i_2},\ldots, s_{i_m})}\partial_{s_{i_1}}\ldots \partial_{s_{i_m}}. $$
    In particular, we have $\dd_I=\dd_{\un{w}_I^{\std}}.$

    Hence for $J\subset I$, and a reduced expression $\un{w}_I^J$ of $w_Iw_J^{-1}$, we have that
    $$\dd_I^J=\zeta^{-\LL_I^J(\un{w}_I^J)}\dd_{\un{w}_I^J}$$
\end{corollary}
\begin{remark} \label{rmk:wontworkforIJ}
    In the above setup, if $J\neq \emptyset$, it is not always true that there is a nice reduced expression $\un{w}_I^J$ of $w_Iw_J^{-1}$ so that $\dd_I^J=\dd_{\un{w}_I^J}$.
    For example, consider $J={2} \subset I=\set{1,2}$ for any $n>2$ (so that these are finitary). 
    Then $\dd_J=\dd_{s_2}$ and $\dd_I=\dd_{(s_1,s_2,s_1)}=\dd_{s_1}\dd_{s_2}\dd_{s_1}$. 
    But  $w_Iw_J^{-1}=s_2s_1$ has a unique reduced expression $(s_2,s_1)$ in this case, and $\dd_I^J=\zeta^{-1}\dd_{s_2}\dd_{s_1}$.
\end{remark}

\subsection{Calculations needed for geometric Satake}
\label{subsec-further frob calculations}

This section consists of various technical identities involving the Frealization structure on the deformed affine realization. We do not attempt to motivate each computation here; they will each be used in the verification of the diagrammatic Quantum Geometric Satake functor in \cref{sec well-defined-ness of the functor}. After \cref{subsubsec quantum numbers} and \cref{subsubsec-notation in further calculations.} the reader is welcome to skip this section.

\subsubsection{Quantum Numbers}\label{subsubsec quantum numbers}
For $n \in \Z_{\geq 1}$, we shall use the notation $[n]_q$ to refer to the \emph{n-th quantum number} explicitly given by 
$$[n]_q=\frac{q^n-q^{-n}}{q-q^{-1}}=q^{n-1}+q^{n-3}+\ldots q^{3-n}+q^{1-n} \in \Z[q,q^{-1}].$$
We also set $[0]_q=0$.
Thus $[1]_q=1, [2]_q=q+q^{-1}$. 
We also have explicit formulas for product of quantum numbers, for example 
$$[2]_q[n]_q=[n+1]_q+[n-1]_q.$$
More generally, for $m \geq n$, we have 
$$[m]_q[n]_q=[n]_q[m]_q=[m+n-1]_q+[n+m-3]_q+\ldots +[m-n+1]_q.$$
We often drop the subscript $q$ for convenience.
For $k \in \Z_{\geq 0}$, we define
$[k]!:= [k][k-1]\ldots [1]$.
We then define the \emph{quantum binomial coefficients} $\qbinom{n}{k}$ to be
$$\qbinom{n}{k}:= \frac{[n]!}{[n-k]! [k]!}.$$


For the calculations below we introduce a formal variable $q$ satisfying 
\begin{equation} \label{zetaq} \zeta^n = q^{-2}. \end{equation} We introduce $q$ so that we can view certain scalars as quantum numbers, and thus connect our calculations to calculations arising for quantum groups. 
More precisely, we work over 
\begin{equation}\label{eq coefficient ring with q and zeta}
    A_{\Z}:= \begin{cases}
\Z[q^{\pm 1}, \zeta^{\pm 1}]/(q^{-2}-\zeta^n) &\text{if $n$ odd.}\\
\Z[q^{\pm 1}, \zeta^{\pm 1}]/(q^{-1}-\zeta^{n/2}) &\text{if $n$ even.}
    \end{cases}
\end{equation}
Every time we use the variable $q$ it will be through one of the following two equalities:
\begin{equation} \sum_{i=0}^{k} \zeta^{-in} = q^k [k+1], \end{equation}
\begin{equation} \label{thisoneisworse} \zeta^{-\frac{n(n-1)}{2}} = q^{n-1}. \end{equation}
There is no need to introduce $q$, as all scalars we use live in the ring $\Z[\zeta,\zeta^{-1}]$. Indeed, we only use the subring of $\Z[q,q^{-1}]$ generated by $q^{\pm 2}$ and $q^{n-1}$, all interpreted within $\Z[\zeta,\zeta^{-1}]$ via \eqref{zetaq} or \eqref{thisoneisworse}. Not all quantum numbers in $q$ live in this subring, but $[n]_q$ does.


\subsubsection{Some notation} \label{subsubsec-notation in further calculations.}
For $k,l,m \in S$, we continue using $\hh{k}$ to denote the (maximal finitary) subset $S \setminus \set{s} \subset S$, and $\hh{klm}$ to denote the subset $S \setminus \set{k,l,m}$, etc. 
We then have the corresponding Frobenius trace maps $\dd_{\hh{k}}$, $\dd_{\hh{klm}}$ and the product-coproduct elements $\mu_{\hh{k}}$, $\mu_{\hh{klm}}$.

For finitary $I \subset S$, $i,j \in I$ with $i\neq j$, we set
\begin{equation} \label{eq mu I I minus j I minus i}
  \mu_{I}^{I\setminus \set{i}, I\setminus \set{j}}:=  \frac{\mu_I\mu_{I\setminus \set{i,j}}}{\mu_{I\setminus \set{i}}\mu_{I\setminus \set{j}}}=\frac{\mu_I^{I\setminus \set{i}}}{\mu_{I\setminus \set{j}}^{I\setminus \set{i,j}}}.
\end{equation}
By the inclusion-exclusion principle, this is the product of all the positive roots in $\Phi^+_I$ which are not in $\Phi^+_{I \setminus \set{i}}$ or $\Phi^+_{I \setminus \set{j}}$.
For example, 
\begin{equation} \label{mu01k} \mu_{\hh{1}}^{\hh{01},\hh{1k}}= \frac{\mu_{\hh{1}}^{\hh{01}}}{\mu_{\hh{1k}}^{\hh{01k}}}= \frac{\prod_{2\leq i \leq n}(x_i - \zeta^{n-i}x_1)}{\prod_{k+1\leq i \leq n}(x_i - \zeta^{n-i}x_1)}=\prod_{2\leq i \leq k}(x_i - \zeta^{n+1-i}x_1).
\end{equation}

We also use Sweedler notation: having a $\Delta^J_{I}{}_{(1)}$ and $\Delta^J_{I}{}_{(2)}$ in the same expression indicates a sum over a dual basis $\set{e_i}, \set{f_i}$ for the extension $R^I \subset R^J$, i.e., we substitute $e_i$ for $\Delta^J_{I}{}_{(1)}$ and $f_i$ for $\Delta^J_{I}{}_{(2)}$ and take a sum over $i$.
For example, $\Delta^J_{I}{}_{(1)}\otimes \Delta^J_{I}{}_{(2)} = \sum_i e_i \ot f_i$ is by definition the image of $1$ under the coproduct map  $\Delta_I^J$ (see \cref{subsec- cubes of frob}).

\begin{remark}
    Using Sweedler notation, the expression in \cref{eq mu I I minus j I minus i} is also equal to 
    \begin{equation}  \mu_{I}^{I\setminus \set{i}, I\setminus \set{j}} =  \dd_{I\setminus \set{j}}^{I\setminus \set{i,j}}(\Delta^{I\setminus \set{i}}_{I,(1)})\Delta^{I\setminus \set{i}}_{I,(2)} = \dd_{I\setminus \set{i}}^{I\setminus \set{i,j}}(\Delta^{I\setminus \set{j}}_{I,(1)})\Delta^{I\setminus \set{j}}_{I,(2)}. \end{equation} This is a consequence of \cref{lemma deformed affine frobenius satisfies extra conditions}.
\end{remark}

\subsubsection{Dual bases}
\label{subsec-dual basis}

In this section, we construct explicit dual bases for certain extensions, which are necessary for further computations.

\begin{lemma} Given dual bases for $R^I \subset R^J$, we get dual bases for $R^{\sigma(I)} \subset R^{\sigma(J)}$ by applying $\sigma$ to the bases. \end{lemma}

\begin{proof} Immediate from \cref{rmk symmetry for rotation operator}. \end{proof}

\begin{lemma}\label{lemma dual basis left ext} The sets 
   $$ \{(-1)^{k-i}h_i(x_1)\}_{0\leq i \leq k},\;\{e_{k-i}(\zeta x_2, \zeta^2 x_3, \ldots, \zeta^{k}x_{k+1}) \}_{0\leq i \leq k}$$ form dual bases for the Frobenius extensions $R^{\{1,2, \ldots, k\}} \hookrightarrow R^{\{2, \ldots, k\}}$ for $k<n$. 
 
\end{lemma}
\begin{proof}
    We want to show:
    \begin{equation}\label{eq 1 dual basis left ext} 
        \dd_{\set{1,2, \ldots,k}}^{\set{2,3, \ldots, k}}\left((-1)^{k-i}h_i(x_1)e_{k-j}(\zeta x_2, \zeta^2 x_3, \ldots, \zeta^{k}x_{k+1})\right)=\delta_{i,j}.
    \end{equation}
    By \cref{corollary-normalization}, we have that 
    \begin{equation} \label{eq 2 dual basis left ext}
        \dd_{\set{1,2, \ldots,k}}^{\set{2,3, \ldots, k}}=\zeta^{-l}\dd_k \dd_{k-1} \ldots \dd_1
    \end{equation}
    where $l$ is the rex weight of $(s_k, s_{k-1}, \ldots s_1)$.
    Explicitly, from \cref{defn-rex weight}, 
    $$l=\LL_{\set{1,2,\ldots k}}((s_k,s_{k-1},\ldots, s_1)\bullet (s_2, s_3, \ldots, s_k, s_2, \ldots, s_{k-1}, \ldots, s_2, s_3, s_2) ).$$
    Comparing with the standard rex $\un{w}_{\set{1,2,\ldots,k}}=(s_1,s_2, \ldots s_k, s_1, \ldots ,s_{k-1}, \ldots, s_1, s_2, s_1)$, we see that
    $l=k(k-1)/2$.
    Hence \cref{eq 1 dual basis left ext} becomes
    \begin{equation}\label{eq 3 dual basis left ext}
        \zeta^{-k(k-1)/2}\dd_k \dd_{k-1} \ldots \dd_1\left((-1)^{k-i}h_i(x_1)e_{k-j}(\zeta x_2, \zeta^2 x_3, \ldots, \zeta^{k}x_{k+1})\right)=\delta_{i,j}.
    \end{equation}
    Now let's proceed to prove \cref{eq 3 dual basis left ext}. We use induction on $k$ for $1\leq k \leq n-1$. The base case is easy and left as an exercise to the reader. So now consider $k>1$.\\
    We first consider the case $i=0$. 
    In this case, if $j>0$, due to degree reasons, we get 
    $$\dd_k \dd_{k-1} \ldots \dd_1\left((-1)^{k}e_{k-j}(\zeta x_2, \zeta^2 x_3, \ldots, \zeta^{k}x_{k+1})\right)=0.$$
    When $j=0$, we have
    \begin{align}
        \nonumber&\dd_k \dd_{k-1} \ldots \dd_1\left((-1)^{k}e_{k}(\zeta x_2, \zeta^2 x_3, \ldots, \zeta^{k}x_{k+1})\right)\\
        &\nonumber= (-1)^k\dd_k \ldots \dd_1
        \left(
         \zeta x_2 e_{k-1}(\zeta^2x_3,\ldots, \zeta^kx_{k+1})
        +e_{k}(\zeta^2x_3, \ldots \zeta^k x_{k+1})        
        \right) \\
        &\nonumber=(-1)^k\dd_k\ldots\dd_2\left(\dd_1(\zeta x_2)\;e_{k-1}(\zeta^2x_3, \ldots, \zeta^kx_{k+1}) + 0\right)\\
        &=(-1)^{k-1}\dd_k\ldots \dd_2\left(\zeta^{k-1}e_{k-1}(\zeta x_3, \ldots, \zeta^{k-1} x_{k+1})\right)\label{eq 4 dual basis left ext}\\
        &=(-1)^{k-1}\zeta^{k-1}\dd_{k-1}\ldots \dd_1\left(e_{k-1}(\zeta x_2, \ldots, \zeta^{k-1} x_{k})\right) \label{eq 5 dual basis left ext} \\
        &=\zeta^{k-1} \zeta^{(k-1)(k-2)/2}=\zeta^{k(k-1)/2} \label{eq 6 dual basis left ext}
    \end{align}
    and we have that \cref{eq 3 dual basis left ext} holds. Note that in the above computation, we used symmetry with respect to the rotation operator $\sigma$, (see \cref{rmk symmetry for rotation operator}) to go from \eqref{eq 4 dual basis left ext} to \eqref{eq 5 dual basis left ext}, and induction hypothesis to go from \eqref{eq 5 dual basis left ext} to \eqref{eq 6 dual basis left ext}.

    Now suppose $i>0$. Then we have: 
    \begin{align}
        \nonumber(-1)^{k-i}&\dd_k\ldots \dd_2 \dd_1\left(x_1^ie_j(\zeta x_2, \ldots , \zeta^k x_{k+1}\right) \\
        \nonumber&= (-1)^{k-i}\dd_k \ldots \dd_2\dd_1 \left(x_1^i\zeta x_2 e_{j-1}(\zeta^2x_3, \ldots, \zeta^k x_{k+1})+x_1^ie_j(\zeta^2 x_3, \ldots , \zeta^k x_{k+1}) \right)\\
       \nonumber &= (-1)^{k-i}\dd_k \ldots \dd_2\left(x_1 \zeta x_2 \dd_1(x_1^{i-1})e_{j-1}(\zeta^2 x_3, \ldots , \zeta^k x_{k+1}) + \dd_1(x_1^i)e_j(\zeta^2 x_3, \ldots , \zeta^k x_{k+1})
        \right)\\
        \nonumber&= (-1)^{k-i} \dd_k \ldots \dd_2\left(
        \begin{array}{cc}
         &\displaystyle x_1 \zeta x_2 \sum_{l'=0}^{i-2} x_1^{i-2-l'}(\zeta x_2)^{l'}e_{j-1}(\zeta^2 x_3, \ldots \zeta^k x_{k+1})    \\
          + &\displaystyle \sum_{l=0}^{i-1}x_1^{i-1-l}(\zeta x_2)^le_j(\zeta^2x_3, \ldots, \zeta^k x_{k+1})   
        \end{array}\right)\\
        &= (-1)^{k-i}\left(
        \begin{array}{cc}
        &\displaystyle \sum_{l'=0}^{i-2} x_1^{i-1-l'}\dd_k \ldots \dd_2((\zeta x_2)^{l'+1}e_{j-1}(\zeta^2 x_3, \ldots \zeta^k x_{k+1}))    \\
          + &\displaystyle \sum_{l=0}^{i-1}x_1^{i-1-l}\dd_k \ldots \dd_2((\zeta x_2)^le_j(\zeta^2x_3, \ldots, \zeta^k x_{k+1}))
        \end{array}\right)
    \end{align}
    Since $i,j \leq k$, we have that $i-1, j-1 \leq k-1$. 
    An argument using rotational symmetry and induction hypothesis shows then that
    \begin{align*}\dd_k \ldots \dd_2((\zeta x_2)^{l'+1}e_{j-1}(\zeta^2 x_3, \ldots \zeta^k x_{k+1}))&=\dd_{k-1} \ldots \dd_1((\zeta x_1)^{l'+1}e_{j-1}(\zeta^2 x_2, \ldots, \zeta^k x_{k}))\\
    &=\zeta^{l'+j}\dd_{k-1} \ldots \dd_1((x_1)^{l'+1}e_{j-1}(\zeta x_2, \ldots ,\zeta^{k-1} x_{k}))\\
    &=\zeta^{k-1}(-1)^{k-l'-2}\zeta^{(k-1)(k-2)/2}\delta_{j-1,k-l'-2}\\
    &=\zeta^{k(k-1)/2}(-1)^{j-1}\delta_{j-1,k-l'-2},
    \end{align*}
    and
    \begin{align*}
     \dd_k \ldots \dd_2((\zeta x_2)^le_j(\zeta^2x_3, \ldots, \zeta^k x_{k+1}))&=\dd_{k-1} \ldots \dd_1((\zeta x_1)^le_j(\zeta^2x_2, \ldots, \zeta^k x_{k}))  \\
     &= \zeta^{j+l} \dd_{k-1} \ldots \dd_1(( x_1)^le_j(\zeta x_2, \ldots, \zeta^{k-1} x_{k}))\\
     &=\zeta^{j+l}(-1)^{k-1-l}\zeta^{(k-1)(k-2)/2}\delta_{j, k-1-l}\\
     &=\zeta^{k(k-1)/2}(-1)^{j}\delta_{j,k-1-l}.
    \end{align*}
Hence RHS of \eqref{eq 6 dual basis left ext} simplifies to
\begin{align*}
    &(-1)^{k-i}\left(
    \begin{array}{c}
    \displaystyle \sum_{l'=0}^{i-2} x_1^{i-1-l'}\zeta^{k(k-1)/2}(-1)^{j-1}\delta_{j-1,k-l'-2}    \\
    + \displaystyle \sum_{l=0}^{i-1}x_1^{i-1-l}\zeta^{k(k-1)/2}(-1)^{j}\delta_{j,k-1-l}
     \end{array}   \right)\\
    &=(-1)^{k-i}\zeta^{k(k-1)/2}(-1)^{j} \left(
    \displaystyle \sum_{l=0}^{i-1}x_1^{i-1-l}\delta_{j,k-1-l}
    -\displaystyle \sum_{l'=0}^{i-2} x_1^{i-1-l'}\delta_{j-1,k-l'-2}
    \right)\\
    &=(-1)^{k-i-j}\zeta^{k(k-1)/2}\delta_{j,k-i}=\zeta^{k(k-1)/2}.
\end{align*}
This implies \eqref{eq 3 dual basis left ext}.
\end{proof}

\begin{remark}\label{rmk explicit dual basis permutation realization}
    \cref{lemma dual basis left ext} with $\zeta$ specialized to $1$ provides explicit dual basis for the corresponding extension for the permutation realization of $S_n$ over $\Z$. 
\end{remark}


\begin{corollary}\label{corollary dual basis right ext} The sets
     $$ \{(-1)^{i}\zeta^{k(k+1)/2}h_i(x_{k+1})\}_{0\leq i \leq k},\;\{e_{k-i}(\zeta^{-1} x_k, \zeta^{-2} x_{k-1}, \ldots, \zeta^{-k}x_1) \}$$ form dual bases for the Frobenius extensions $R^{\{1,2, \ldots, k\}} \hookrightarrow R^{\{1,2, \ldots, k-1\}}$ for $k<n$.
\end{corollary}
\begin{proof}
    One can do a similar explicit proof as in the proof of \cref{lemma dual basis left ext} to show this. 
    Alternatively, one may use \cref{rmk dual basis from permutation realization} and \cref{rmk explicit dual basis permutation realization} to see this as a consequence of \cref{lemma dual basis left ext}. 
    One may also see this as a consequence of \cref{lemma dual basis left ext} using symmetry of our Frobenius realization under the rotation \eqref{eq rotational symmetry deformed affine} and reflection \eqref{eq:reflection symmetry for reflection operator}.
\end{proof}

\subsubsection{Calculations}\label{subsubsec calculations}

\begin{lemma}\label{lemma demazure adjacent} Let $0\leq k<n$. Then we have the following equalities (whenever all subsets appearing are finitary):  
    \begin{enumerate}
        \item[(a)]   \hfill$ \dd_{\{1,2, \ldots ,k\}}^{\{1,2, \ldots ,k-1\}}\left((x_1-\zeta^{-(n-k)}x_{k+1})\ldots (x_{k}-\zeta^{-(n-1)}x_{k+1})\right)=q^{k}[k+1]_q$ \hfill \null
        \item[(b)] \hfill $\dd^{\hh{0,-1}}_{\hh{0}}(\mu^{\hh{0,-1}}_{\hh{-1}}) =(-1)^{n-1}[n]_q$. \hfill \null
        \item[(c)]  \hfill $\dd_{\{1,2, \ldots ,k-1\}}^{\{2, \ldots ,k-1\}}\left((x_1-\zeta^{-(n-k+1)}x_k)\ldots (x_1-\zeta^{-(n-1)}x_2)\right)=q^{k-1}[k]_q.$ \hfill \null\label{lemma demazure for bigon}
        \item[(d)] \hfill $ \dd_{\hh{0}}^{\hh{0,1}}\left(\mu^{\hh{0,1}}_{\hh{1}} \right)= (-1)^{n-1}[n]_q$. \hfill \null
    \end{enumerate}
\end{lemma}
\begin{proof}   
  To prove part (a), observe that 
  \begin{align*}(x_1-\zeta^{-(n-k)}x_{k+1})\cdots &(x_{k}-\zeta^{-(n-1)}x_{k+1})\\
  &= (\zeta^{-(n-k)})\cdots (\zeta^{-(n-1)})(\zeta^{n-k}x_1-x_{k+1})\cdots (\zeta^{n-1}x_k-x_{k+1}) \\
  &=\zeta^{-kn+\frac{k(k+1)}{2}}\sum_{0\leq i \leq k}(-1)^ih_i(x_{k+1})e_{k-i}(\zeta^{n-1}x_k, \ldots, \zeta^{n-k}x_1)\\
  &=q^{2k}\sum_{0 \leq i \leq k}(-1)^i\zeta^{\frac{k(k+1)}{2}}h_{i}(x_{k+1})\zeta^{n(k-i)}e_{k-i}(\zeta^{-1}x_k, \ldots, \zeta^{-k}x_1).
  \end{align*}
  Applying $\dd_{\{1,2, \ldots ,k\}}^{\{1,2, \ldots ,k-1\}}$ to both sides, using \cref{corollary dual basis right ext}, we see that
  \begin{align*}
    \dd_{\{1,2, \ldots ,k\}}^{\{1,2, \ldots ,k-1\}}\left((x_1-\zeta^{-(n-k)}x_{k+1})\ldots (x_{k}-\zeta^{-(n-1)}x_{k+1})\right)  
    =q^{2k}\sum_{0 \leq i \leq k}q^{-2(k-i)}=q^{k}[k+1]_q.
  \end{align*}
  
  Now we show that (a) implies (b).
    Observe that from \cref{defn deformed affine frobenius realization} and \cref{lemma deformed affine frobenius realization is frobenius}, we have that
    $$\mu_{\{0,1, \ldots ,n-2\}}^{\{1,2, \ldots ,n-2\}}=(x_n-\zeta x_1)\ldots (x_n-\zeta^{n-1} x_{n-1}).$$
    So (c) reduces to showing
    $$\dd_{\{1,2, \ldots ,n-1\}}^{\{1,2, \ldots ,n-2\}}((x_n-\zeta x_1)\ldots (x_n-\zeta^{n-1} x_{n-1}))=(-1)^{n-1}[n]_q.$$
    We can rewrite the L.H.S above as:
    \begin{align}       
   LHS &= (-\zeta)(-\zeta^2)\ldots (-\zeta^{n-1})\dd_{\{1,2, \ldots ,n-1\}}^{\{1,2, \ldots ,n-2\}}((x_1-\zeta^{-1}x_n)\ldots (x_{n-1}-\zeta^{-n+1}x_n)) \nonumber\\
   &= (-1)^{n-1}\zeta^{\frac{n(n-1)}{2}}\dd_{\{1,2, \ldots ,n-1\}}^{\{1,2, \ldots ,n-2\}}((x_1-\zeta^{-1}x_n)\ldots (x_{n-1}-\zeta^{-n+1}x_n))\nonumber\\
   &= (-1)^{n-1}q^{-(n-1)}\dd_{\{1,2, \ldots ,n-1\}}^{\{1,2, \ldots ,n-2\}}((x_1-\zeta^{-1}x_n)\ldots (x_{n-1}-\zeta^{-n+1}x_n)).\nonumber
   \end{align}
   Hence (b) follows from (a), taking $k$ to be $n-1$.

   Now we prove (c).
   LHS may be written as
    \begin{align*}
    &\dd_{\{1,2, \ldots ,k-1\}}^{\{2, \ldots ,k-1\}}\left((x_1-\zeta^{-(n-k+1)}x_k)\ldots (x_1-\zeta^{-(n-1)}x_2)\right) \\
    &\quad= \dd_{\{1,2, \ldots ,k-1\}}^{\{2, \ldots ,k-1\}}\left(\sum_{0 \leq i \leq k-1}h_i(x_1)e_{k-1-i}(-\zeta^{-n+1}x_2, -\zeta^{-n+2} x_3, \ldots, -\zeta^{-n+(k-1)}x_k)\right)\\
    &\quad=\dd_{\{1,2, \ldots ,k-1\}}^{\{2, \ldots ,k-1\}}\left(\sum_{0 \leq i \leq k-1}(-1)^{k-1-i}\zeta^{-n(k-1-i)}h_i(x_1)e_{k-1-i}(\zeta x_2, \zeta^{2} x_3, \ldots, \zeta^{-(k-1)}x_k)\right)\\
    &\quad = \sum_{0 \leq i \leq k-1} \zeta^{-n(k-1-i)} = \sum_{0 \leq i \leq k-1}q^{2(k-1-i)}=q^{(k-1)}[k]_q.
    \end{align*}
    where we used \cref{lemma dual basis left ext} for the first equality in the last line above.
    For (d), we have from \cref{defn deformed affine frobenius realization} that
    $$\mu^{\hh{0,1}}_{\hh{1}}= (x_2-\zeta^{n-1}x_1)(x_3-\zeta^{n-2}x_1)\ldots (x_n-\zeta x_1)=(-1)^{n-1}\zeta^{\frac{n(n-1)}{2}}\cdot(x_1-\zeta^{-1}x_{n})\ldots (x_1-\zeta^{-(n-1)}x_2).$$
    Noting that $\zeta^{n(n-1)/2}= q^{-(n-1)},$ (d) then follows from (c) when $k=n.$  
\end{proof}


\begin{corollary}\label{corollary- bigon bursting}
For $k \not= 0,1$, we have:
    $$\dd_{\hh{0,k}}^{\hh{0,1,k}}\left(\mu_{\hh{1}}^{\hh{0,1},\hh{1,k}}\right)=(-1)^{k-1}\zeta^{-k(k-1)/2}q^{-(k-1)}[k]_q.$$
    \begin{proof}
        Rewriting the LHS, we get
        \begin{align}
        \dd_{\hh{0,k}}^{\hh{0,1,k}}\left(\mu_{\hh{1}}^{\hh{0,1},\hh{1,k}}\right)&=\dd_{\hh{0,k}}^{\hh{0,1,k}}\left( (x_2-\zeta^{n-1}x_1)(x_3-\zeta^{n-2}x_1)\ldots (x_k-\zeta^{n+1-k}x_1)\right). \label{eq 5.7.1}
        \end{align}
        We also have that 
        \[\dd_{\hh{0,k}}^{\hh{0,1,k}}=\dd_{\set{1,\ldots, k-1}\sqcup\set{k+1,\ldots, n-1}}^{\set{2,\ldots, k-1}\sqcup \set{k+1, \ldots, n-1}}=\dd_{\set{1,\ldots, k-1}}^{\set{2,\ldots, k-1}}.\]
       Hence we may simplify \eqref{eq 5.7.1} further:
       \begin{align}
        \dd_{\hh{0,k}}^{\hh{0,1,k}}\left(\mu_{\hh{1}}^{\hh{0,1},\hh{1,k}}\right)&=\dd_{\set{1,\ldots, k-1}}^{\set{2,\ldots, k-1}}\left( (x_2-\zeta^{n-1}x_1)(x_3-\zeta^{n-2}x_1)\ldots (x_k-\zeta^{n+1-k}x_1)\right) \nonumber\\
        &= (-1)^{k-1}\zeta^{(n-1)+(n-2)+\ldots + (n+1-k)}
        \\&\hspace{20pt}\dd_{\set{1,\ldots, k-1}}^{\set{2,\ldots, k-1}}\left((x_1-\zeta^{-(n-k+1)}x_k)\ldots (x_1-\zeta^{-(n-1)}x_2)\right)\nonumber\\
        &=(-1)^{k-1}\zeta^{(k-1)n-\frac{(k-1)k}{2}}q^{(k-1)}[k]_q \label{eq 5.7.2}\\
        &=(-1)^{k-1}(q^{-2})^{k-1}\zeta^{-\frac{(k-1)k}{2}}q^{(k-1)}[k]_q = (-1)^{k-1}\zeta^{-\frac{(k-1)k}{2}}q^{-(k-1)}[k]_q \label{eq 5.7.3}
        \end{align}
        where we used Lemma \ref{lemma demazure for bigon} for the equality in \eqref{eq 5.7.2} and the fact that $\zeta^n=q^{-2}$ for the first equality in \eqref{eq 5.7.3}. 
    \end{proof}
\end{corollary}

\begin{lemma}\label{lemma square flop 1 tensor 1}
    The following equality holds in $\Bimod{\hh{1,k+1}}{}{R^{\hh{1,k+1}}}\otimes_{\hh{1}}\Bimod{}{\hh{0,1}}{R^{\hh{0,1}}}$ for $2\leq k\leq n-1$.
    
    $$\ddm{1,2, k+1}{1,k+1}\left(\ddm{1,2,k+1}{2,k+1}(\Del{1,2}{2}{1})\Del{1,2}{1}{1}\right) \otimes \ddm{0,1,2}{0,1}\left(\ddm{0,1,2}{0,2}(\Del{1,2}{2}{2})\Del{1,2}{1}{2} \right) =[k-1]_q(-1)^k\zeta^{\frac{-(k-1)k}{2} -1}q^{-k}(1\otimes 1)$$
\end{lemma}
\begin{proof}
    Let us start by analyzing the LHS. 
    Observe that we need to make explicit choices for the dual basis for the extensions $R^{\hh{2}}\hookrightarrow R^{\hh{1,2}}$ and $R^{\hh{1}}\hookrightarrow R^{\hh{1,2}}$ to make the computations involving $\Del{1,2}{1}{.}$ and $\Del{1,2}{2}{.}$. 
    This is precisely the content of \cref{lemma dual basis left ext} and \cref{corollary dual basis right ext}; we choose to let $\Del{1,2}{2}{1}$ run over $\set{e_{n-1-i}(\zeta^{-1}x_1, \zeta^{-2}x_n, \ldots , \zeta^{-(n-1)}x_3)}$ and $\Del{1,2}{2}{2}$ run over $ \set{(-1)^i\zeta^{n(n-1)/2}h_i(x_2)}$.
    Similarly, we choose $\Del{1,2}{1}{1}$ run over $\set{(-1)^{n-1-j}h_{j}(x_2)}$ and $\Del{1,2}{1}{2}$ run over $\set{e_{n-1-j}(\zeta x_3, \zeta^2 x_4, \ldots \zeta^{n-1} x_1)}.$
    So we may write the LHS as:
    \begin{align}
    \sum_{0\leq i,j \leq n-1} \left( \begin{array}{c}
    \ddm{1,2, k+1}{1,k+1}\left(\ddm{1,2,k+1}{2,k+1}\left(e_{n-1-i}(\zeta^{-1}x_1, \zeta^{-2}x_n, \ldots , \zeta^{-(n-1)}x_3)\right)(-1)^{n-1-j}h_{j}(x_2)\right) \vspace{5 pt}\\ \otimes \ddm{0,1,2}{0,1}\left(\ddm{0,1,2}{0,2}((-1)^i\zeta^{n(n-1)/2}h_i(x_2))e_{n-1-j}(\zeta x_3, \zeta^2 x_4, \ldots \zeta^{n-1} x_1) \right) \label{eq 1 lemma square flop 1 tensor 1}
    \end{array}\right)
    \end{align}
    Note that the sum above should give a degree 0 term in $\Bimod{\hh{1,k+1}}{}{R^{\hh{1,k+1}}}\otimes_{\hh{1}}\Bimod{}{\hh{1,2}}{R^{\hh{1,2}}}$. This implies that we need only consider tuples $(i,j)$ where both the tensor factors are of degree $0$.
    This happens precisely when $i-j=0$, i.e., when $i=j$.
    So we now just analyze the terms in \eqref{eq 1 lemma square flop 1 tensor 1} for $i=j.$
    
    Now let us simplify the second tensor factor on the LHS. This factor is completely independent of $k$.
    We have that $\ddm{0,1,2}{0,2}=\dd^{\set{3, 4, \ldots , n-1}}_{\set{1}\sqcup \set{3,4, \ldots, n-1}}=\dd_1$, so we have that the $i=0$ term is 0. 
    Now we can simplify the second term in \eqref{eq 1 lemma square flop 1 tensor 1} to be
    \begin{align} \label{floobarl}
    &(-1)^i\zeta^{n(n-1)/2} \ddm{0,1,2}{0,1}\left( \dd_1(h_i(x_2))e_{n-1-i}(\zeta x_3, \zeta^2 x_4, \ldots , \zeta^{n-1}x_1) \right)\\ \nonumber
    &=(-1)^i\zeta^{n(n-1)/2} \ddm{0,1,2}{0,1}\left( (-\zeta^{-1}) h_{i-1}(x_2, \zeta^{-1}x_1)e_{n-1-i}(\zeta x_3, \zeta^2 x_4, \ldots , \zeta^{n-1}x_1)
    \right)\\ \nonumber
    &=(-1)^{i+1}\zeta^{\frac{n(n-1)}{2} -1} \ddm{0,1,2}{0,1}\left(
    \sum_{l=0}^{i-1}(\zeta^{-1}x_1)^{i-1-l}h_l(x_2)e_{n-1-i}(\zeta x_3, \zeta^2 x_4, \ldots , \zeta^{n-1}x_1)
    \right)\end{align}
    Now $x_1 \in W_{\hh{0,1}}$ so it can be pulled out of the operator $\dd^{\hh{0,1,2}}_{\hh{0,1}}$. If this happens then $\dd^{\hh{0,1,2}}_{\hh{0,1}}$ applied to the remainder would be zero, for degree reasons. Consequently we can ignore all terms above where $x_1$ appears as a factor, both those terms where $l \ne i-1$, and those terms within $e_{n-1-i}$ where $\zeta^{n-1} x_1$ is chosen as a factor. Consequently \eqref{floobarl} simplifies to
     \begin{align}
    &(-1)^{i+1}\zeta^{\frac{n(n-1)}{2} -1} \ddm{0,1,2}{0,1}\left(
    h_{i-1}(x_2)
    e_{n-1-i}(\zeta x_3, \zeta^2 x_4, \ldots , \zeta^{n-2}x_0)
    \right)\\ \nonumber
    &=(-1)^{i+1}\zeta^{\frac{n(n-1)}{2}-1} (-1)^{n-2-(i-1)}
    =(-1)^{n}\zeta^{\frac{n(n-1)}{2}-1}.
    \end{align}
    The first equality above used \cref{lemma dual basis left ext}.

    Hence \eqref{eq 1 lemma square flop 1 tensor 1} becomes
    \begin{align}
    \sum_{1\leq i \leq n-1} \left( \begin{array}{c}
    \ddm{1,2, k+1}{1,k+1}\left(\ddm{1,2,k+1}{2,k+1}\left(e_{n-1-i}(\zeta^{-1}x_1, \zeta^{-2}x_n, \ldots , \zeta^{-(n-1)}x_3)\right)(-1)^{n-1-i}h_{i}(x_2)\right) \vspace{5 pt}\\ \otimes (-1)^{n}\zeta^{\frac{n(n-1)}{2}-1} \label{eq 2 lemma square flop 1 tensor 1}
    \end{array}\right)
    \end{align}

    Now let us analyze the first tensor factor.
    We may write the Frobenius trace maps $\ddm{1,2, k+1}{1,k+1}$ and $\ddm{1,2,k+1}{2,k+1}$ explicitly using \cref{corollary-normalization} as
    \[\ddm{1,2,k+1}{2,k+1}=\dd^{\set{3, 4, \ldots , k} \sqcup \set{k+2, \ldots, n-1, n}}_{\set{3,4, \ldots k}\sqcup \set{k+2,\ldots,n,1}}=\dd_{k+2}\ldots \dd_n \dd_1, \]
    \[ \ddm{1,2,k+1}{1,k+1}=\dd^{\set{3, 4, \ldots , k} \sqcup \set{k+2, \ldots, n-1, n}}_{\set{2,3 \ldots k}\sqcup \set{k+2,\ldots,n-1, n}}=\zeta^{-\frac{(k-1)(k-2)}{2}}\dd_k \ldots \dd_3 \dd_2=\dd^{\set{3,4, \ldots, k}}_{\set{2, 3, \ldots, k}}.\]

    The first tensor factor in \eqref{eq 2 lemma square flop 1 tensor 1} (for a given $i\geq0$) then becomes:
    \begin{align}
    &\dd^{\set{3,4, \ldots, k}}_{\set{2, 3, \ldots, k}} \left(
    \dd_{k+2}\ldots \dd_n \dd_1 (e_{n-1-i}(\zeta^{-1}x_1, \zeta^{-2}x_n, \ldots \zeta^{-(n-1)}x_3))(-1)^{n-1-i}h_i(x_2)    
    \right). \label{eq 3 lemma square flop 1 tensor 1}
    \end{align}
    Every term in $e_{n-1-i}(\zeta^{-1}x_1, \zeta^{-2}x_n, \ldots \zeta^{-(n-1)}x_3))$ which does not include $\zeta^{-1} x_1$ is invariant under $s_1$ and thus killed by $\dd_1$. Thus
    \begin{align} \dd_1(e_{n-1-i}&(\zeta^{-1}x_1, \zeta^{-2}x_n, \ldots \zeta^{-(n-1)}x_3)) = e_{n-2-i}(\zeta^{-2}x_n, \ldots \zeta^{-(n-1)}x_3) \dd_1(\zeta^{-1} x_1) = \nonumber \\
    & \zeta^{-1} e_{n-2-i}(\zeta^{-2}x_n, \ldots \zeta^{-(n-1)}x_3). \end{align}
    Similarly, applying $\dd_n$ to $e_{n-2-i}(\zeta^{-2}x_n, \ldots \zeta^{-(n-1)}x_3)$ will kill any term not including $\zeta^{-2} x_n$, and so forth. We deduce that
       \begin{align}
       \nonumber &\dd_{k+2}\ldots \dd_n \dd_1 (e_{n-1-i}(\zeta^{-1}x_1, \zeta^{-2}x_n, \ldots \zeta^{-(n-1)}x_3))\\ \nonumber
       & = \zeta^{-1-2-\ldots -(n-k)} e_{k-1-i}(\zeta^{-(n-k+1)}x_{k+1}, \ldots, \zeta^{-(n-2)}x_4, \zeta^{-(n-1)}x_3) \\ 
        &= \zeta^{-\frac{(n-k)(n-k+1)}{2}}e_{k-1-i}(\zeta^{-(n-k+1)}x_{k+1}, \ldots, \zeta^{-(n-2)}x_4, \zeta^{-(n-1)}x_3). \label{eq left tensor factor square flop}
    \end{align}
    Thus \eqref{eq 3 lemma square flop 1 tensor 1} simplifies to
    \begin{align}
    \notag  &\dd^{\set{3,4, \ldots, k}}_{\set{2, 3, \ldots, k}} \left(
     \zeta^{-\frac{(n-k)(n-k+1)}{2}}e_{k-1-i}(\zeta^{-(n-k+1)}x_{k+1}, \ldots, \zeta^{-(n-2)}x_4, \zeta^{-(n-1)}x_3)
    (-1)^{n-1-i}h_i(x_2)  
    \right) \\
    \notag&=\zeta^{-\frac{(n-k)(n-k+1)}{2}} (-1)^{n-1-i} \dd^{\set{3,4, \ldots, k}}_{\set{2, 3, \ldots, k}} \left( \zeta^{-n(k-1-i)}e_{k-1-i}(\zeta x_3, \zeta^2 x_4, \ldots, \zeta^{k-1}x_{k+1})h_i(x_2)
    \right)\\
    &= \begin{cases}
    \begin{array}{cc}
    \zeta^{-\frac{(n-k)(n-k+1)}{2}}\zeta^{-n(k-1-i)} (-1)^{n-1-i} (-1)^{k-1-i}  &\text{if $0\leq i \leq k-1$} \\
    0  &\text{otherwise} 
    \end{array}
    \end{cases} \label{eq 4 lemma square flop 1 tensor 1}\\
    &=\begin{cases}
    \begin{array}{cc}
    \zeta^{-\frac{(n-k)(n-k+1)}{2}-n(k-1-i)} (-1)^{n-k}  &\text{if $0\leq i \leq k-1$}\\
    0  &\text{otherwise}
    \end{array}
    \end{cases}
    \label{eq 5 lemma square flop 1 tensor 1}  
    \end{align}
    where we used \cref{lemma dual basis left ext} for the equality in \eqref{eq 4 lemma square flop 1 tensor 1}.
    Using \eqref{eq 5 lemma square flop 1 tensor 1} to simplify \eqref{eq 2 lemma square flop 1 tensor 1}, we get that the expression in \eqref{eq 2 lemma square flop 1 tensor 1} becomes
    \begin{align}\notag\sum_{1\leq i \leq k-1}\zeta^{-\frac{(n-k)(n-k+1)}{2}-n(k-1-i)}& (-1)^{n-k} \otimes (-1)^{n}\zeta^{\frac{n(n-1)}{2}-1} \\
   \notag &= (-1)^k \zeta^{-\frac{k(k-1)}{2}-1}\sum_{i=1}^{k-1} \zeta^{ni} \otimes 1\\
    &=(-1)^k\zeta^{-\frac{k(k-1)}{2}-1}q^{-k}[k-1]_q 1\otimes 1 \label{eq 6 lemma square flop 1 tensor 1}
    \end{align}
    as desired.
\end{proof}

\begin{lemma}\label{lemma square flop 1 tensor x1k}
The following equality holds in $\Bimod{\hh{1,k+1}}{}{R^{\hh{1,k+1}}}\otimes_{\hh{1}}\Bimod{}{\hh{0,1}}{R^{\hh{0,1}}}$, for $2\leq k \leq n-2$.
    \begin{align*}
    &\ddm{1,2, k+1}{1,k+1}\left(\ddm{1,2,k+1}{2,k+1}(\Del{1,2}{2}{1})\Del{1,2}{1}{1}\right) \otimes \ddm{0,1,2}{0,1}\left(\ddm{0,1,2}{0,2}(\Del{1,2}{2}{2}x_1^k)\Del{1,2}{1}{2} \right) \\
    =&[k-1]_q(-1)^k\zeta^{\frac{-(k-1)k}{2} -1}q^{-k}(1\otimes x_1^k)
    \;+\;\zeta^{k-1}\ddm{0,1,k+1}{1,k+1}\left( \ddm{0,1,k+1}{0, k+1}(x_1^k)\Del{0,1}{1}{1}\right) \otimes \Del{0,1}{1}{2} 
    \end{align*}
\end{lemma}
\begin{proof}
    We choose $\Del{1,2}{1}{.}$ and $\Del{1,1}{2}{.}$  as in the proof of \cref{lemma square flop 1 tensor 1}, so that the LHS may be written explicitly as 
      \begin{align}
    \sum_{0\leq i,j \leq n-1} \left( \begin{array}{c}
    \ddm{1,2, k+1}{1,k+1}\left(\ddm{1,2,k+1}{2,k+1}\left(e_{n-1-i}(\zeta^{-1}x_1, \zeta^{-2}x_n, \ldots , \zeta^{-(n-1)}x_3)\right)(-1)^{n-1-j}h_{j}(x_2)\right) \vspace{5 pt}\\ \otimes \ddm{0,1,2}{0,1}\left(\ddm{0,1,2}{0,2}((-1)^i\zeta^{n(n-1)/2}h_i(x_2)x_1^k)e_{n-1-j}(\zeta x_3, \zeta^2 x_4, \ldots \zeta^{n-1} x_1) \right) \label{eq 1 lemma square flop 1 tensor x1k}
    \end{array}\right)
    \end{align}
    Let us now simplify the second tensor factor. As in the proof of \cref{lemma square flop 1 tensor 1}, for a given $i,j$, it equals
    \begin{align}
    &(-1)^i\zeta^{n(n-1)/2} \ddm{0,1,2}{0,1}\left( \dd_1(h_i(x_2)x_1^k)e_{n-1-j}(\zeta x_3, \zeta^2 x_4, \ldots , \zeta^{n-1}x_1) \right)\label{eq 2 lemma square flop 1 tensor x1k}
    \end{align}
    \begin{enumerate}
        \item[Case 1:]  Suppose $i=k$. Then $\dd_1(h_i(x_2)x_1^k)=0$, so that there is no contribution to \eqref{eq 1 lemma square flop 1 tensor x1k} from this term.
    \item[Case 2:] Suppose $i>k$. We will not end up needing this case (the other tensor factor vanishes) but we record the simplification without proof for future reference. All the terms with $j \ne i$ vanish, and the term with $j=i$ is equal to $(-1)^n \zeta^{\frac{n(n-1)}{2} - 1} x_1^k$.

    \item[Case 3:]  Suppose $0\leq i< k$. Then \eqref{eq 2 lemma square flop 1 tensor x1k} simplifies to 
    \begin{align*}
        (-&1)^i\zeta^{n(n-1)/2} \ddm{0,1,2}{0,1}\left( \dd_1(h_i(x_2)x_1^k)e_{n-1-j}(\zeta x_3, \zeta^2 x_4, \ldots , \zeta^{n-1}x_1) \right)\\
        &=(-1)^i\zeta^{n(n-1)/2} \ddm{0,1,2}{0,1}\left((x_1x_2)^i \dd_1(x_1^{k-i})e_{n-1-j}(\zeta x_3, \zeta^2 x_4, \ldots , \zeta^{n-1}x_1) \right)\\
        &= (-1)^i\zeta^{n(n-1)/2}\ddm{0,1,2}{0,1}\left((x_1x_2)^i h_{k-i-1}(x_1, \zeta x_2)e_{n-1-j}(\zeta x_3, \zeta^2 x_4, \ldots , \zeta^{n-1}x_1) \right)\\
        &=(-1)^i\zeta^{n(n-1)/2}\ddm{0,1,2}{0,1}\left(\sum_{l=0}^{k-1-i}x_1^{i+(k-1-i-l)} \zeta^l x_2^{i+l} e_{n-1-j}(\zeta x_3, \zeta^2 x_4, \ldots , \zeta^{n-1}x_1) \right)\\
        &=(-1)^i\zeta^{n(n-1)/2}\sum_{l=0}^{k-1-i}x_1^{i+(k-1-i-l)} \ddm{0,1,2}{0,1}\left(\zeta^l x_2^{i+l} e_{n-1-j}(\zeta x_3, \zeta^2 x_4, \ldots , \zeta^{n-1}x_1) \right)\\
        &=(-1)^i\zeta^{n(n-1)/2}\sum_{l=0}^{k-1-i}x_1^{k-1-l} \ddm{0,1,2}{0,1}\left(\zeta^l x_2^{i+l}
        \left(\begin{array}{c}
    e_{n-1-j}(\zeta x_3, \zeta^2 x_4, \ldots , \zeta^{n-2}x_0) \vspace{5pt}\\
    + \zeta^{n-1}x_1e_{n-2-j}(\zeta x_3, \zeta^2 x_4, \ldots , \zeta^{n-2}x_0)
    \end{array}\right) \right)
    \end{align*}
    But $i+l\leq i+(k-1-i) =k-1 \leq n-2$, so we can use \cref{lemma dual basis left ext} to simplify the above sum to
    \begin{align}
     \notag (-1)^i\zeta^{n(n-1)/2}\sum_{l=0}^{k-1-i}x_1^{k-1-l}\zeta^l (-1)^{n-2-(i+l)}\left(\delta_{i+l, n-2-(n-1-j)}+\zeta^{n-1}x_1 \delta_{i+l, n-2-(n-2-j)} \right)\\
     \notag =\begin{cases}
         \begin{array}{cl}
           0   & \text{if $j<i$, or $j>k$} \\
           (-1)^i\zeta^{n(n-1)/2}x_1^{i}\zeta^{k-1-i}(-1)^{n-k+1}    & \text{if $j=k$}\\
           (-1)^i\zeta^{n(n-1)/2}x_1^{k+i-j}\zeta^{j-1-i}(-1)^{n-j+1}(1-\zeta^{n})    &\text{if $i<j<k$}\\
           (-1)^i\zeta^{n(n-1)/2}x_1^{k}(-1)^{n-2-j}\zeta^{n-1}    &\text{if $i=j$}
         \end{array}.
     \end{cases}\\
     =\begin{cases}
         \begin{array}{cl}
           0   & \text{if $j<i$, or $j>k$} \\
           (-1)^{n-1+(k-i)}\zeta^{\frac{n(n-1)}{2}-1+(k-i)}x_1^{i}   & \text{if $j=k$}\\
           (-1)^{n-1+(j-i)}\zeta^{\frac{n(n-1)}{2}-1+(j-i)}x_1^{k-(j-i)}(1-q^{-2})    &\text{if $i<j<k$}\\
           (-1)^{n}\zeta^{\frac{n(n-1)}{2}-1}q^{-2}x_1^k    &\text{if $i=j$}
         \end{array}.\label{eq 3.5 lemma square flop 1 tensor x1k}
     \end{cases}
    \end{align}  
    \end{enumerate}

    Now we look at the first tensor factor (for a given $i,j$):
    \begin{align}
        \ddm{1,2, k+1}{1,k+1}\left(\ddm{1,2,k+1}{2,k+1}\left(e_{n-1-i}(\zeta^{-1}x_1, \zeta^{-2}x_n, \ldots , \zeta^{-(n-1)}x_3)\right)(-1)^{n-1-j}h_{j}(x_2)\right) \label{eq 4 lemma square flop 1 tensor x1k}
    \end{align}
     From \eqref{eq left tensor factor square flop} in the proof of \cref{lemma square flop 1 tensor 1}, we have that
     \begin{align*}
         &\dd_{k+2}\ldots \dd_n \dd_1 (e_{n-1-i}(\zeta^{-1}x_1, \zeta^{-2}x_n, \ldots \zeta^{-(n-1)}x_3))\\
         &= \zeta^{-\frac{(n-k)(n-k+1)}{2}}e_{k-1-i}(\zeta^{-(n-k+1)}x_{k+1}, \ldots, \zeta^{-(n-2)}x_4, \zeta^{-(n-1)}x_3),
     \end{align*}
     so that the expression in \eqref{eq 4 lemma square flop 1 tensor x1k} can be written as
     \begin{align}
       \ddm{1,2, k+1}{1,k+1}\left(\zeta^{-\frac{(n-k)(n-k+1)}{2}}e_{k-1-i}(\zeta^{-(n-k+1)}x_{k+1}, \ldots, \zeta^{-(n-2)}x_4, \zeta^{-(n-1)}x_3)(-1)^{n-1-j}h_{j}(x_2)\right) \label{eq 5 lemma square flop 1 tensor x1k}  .
     \end{align}
     Hence we need only consider the case when $0\leq i\leq k-1$,  for which the second tensor factor is explicitly given in Case 3 above (eq \eqref{eq 3.5 lemma square flop 1 tensor x1k}), so that we may further assume that $i\leq j \leq k$.\\
     Since $\ddm{1,2,k+1}{1,k+1}=\dd^{\set{3,4, \ldots, k}}_{\set{2, 3, \ldots, k}}$, by \cref{lemma dual basis left ext}, we have that for $j\leq k-1$, 
     \eqref{eq 5 lemma square flop 1 tensor x1k} is $0$ unless $j=i$.
     Hence the computation from \eqref{eq 5 lemma square flop 1 tensor 1} in the proof of \cref{lemma square flop 1 tensor 1} tells us that the first tensor factor is 
     \begin{align}
    \begin{array}{cc}
    \notag \zeta^{-\frac{(n-k)(n-k+1)}{2}-n(k-1-i)} (-1)^{n-k} \delta_{i,j}  &\text{if $0\leq i,j \leq k-1$,}
    \end{array}
     \end{align}
     and the corresponding term in \eqref{eq 1 lemma square flop 1 tensor x1k} is 
     \begin{align}
         \zeta^{-\frac{(n-k)(n-k+1)}{2}-n(k-1-i)} (-1)^{n-k}\delta_{i,j} \otimes (-1)^{n}\zeta^{\frac{n(n-1)}{2}-1}q^{-2}x_1^k \\
         =(-1)^{k}q^{-2} \zeta^{\frac{n(n-1)}{2}-\frac{(n-k)(n-k+1)}{2}-n(k-1-i)-1}\delta_{i,j} \;\;\;1\otimes x_1^k
         \label{eq 6 lemma square flop 1 tensor x1k} &&\text{if $0\leq i,j \leq k-1$.}
     \end{align}
     
    Now we look at the case when $j=k$, in which case \eqref{eq 5 lemma square flop 1 tensor x1k} can be simplified to
    \begin{align*}
       &\ddm{1,2, k+1}{1,k+1}\left(\zeta^{-\frac{(n-k)(n-k+1)}{2}}e_{k-1-i}(\zeta^{-(n-k+1)}x_{k+1}, \ldots, \zeta^{-(n-2)}x_4, \zeta^{-(n-1)}x_3)(-1)^{n-1-j}h_{k}(x_2)\right) \\
       &=(-1)^{n-1-k}\zeta^{-\frac{(n-k)(n-k+1)}{2}}\dd^{\set{3,4, \ldots, k}}_{\set{2, 3, \ldots, k}}\left(\zeta^{-n(k-1-i)}e_{k-1-i}(\zeta^{k-1}x_{k+1}, \ldots, \zeta^{2}x_4, \zeta x_3)x_2 h_{k-1}(x_2)\right)\\
       &=(-1)^{n-1-k}\zeta^{-\frac{(n-k)(n-k+1)}{2}-n(k-1-i)}\dd^{\set{3,4, \ldots, k}}_{\set{2, 3, \ldots, k}}\left(
       \left(\begin{array}{ll}
       e_{k-i}(\zeta^{k-1}x_{k+1}, \ldots, \zeta^{2}x_4, \zeta x_3,x_2) \\
       -e_{k-i}(\zeta^{k-1}x_{k+1}, \ldots, \zeta^{2}x_4, \zeta x_3)
       \end{array}\right)x_2^{k-1}
       \right).
    \end{align*}
     Since $\dd^{\set{3,4, \ldots, k}}_{\set{2, 3, \ldots, k}}(e_{k-i}(\zeta^{k-1}x_{k+1}, \ldots, \zeta^{2}x_4, \zeta x_3)x_2^{k-1})=0$ for $i\leq k-1$ (by \cref{lemma dual basis left ext}) and $e_{k-i}(\zeta^{k-1}x_{k+1}, \ldots, \zeta^{2}x_4, \zeta x_3,x_2) \in R^{\set{2,3, \ldots ,k}}$, we may further simplify the term above to
     \begin{align}
     \notag &(-1)^{n-1-k}\zeta^{-\frac{(n-k)(n-k+1)}{2}-n(k-1-i)}e_{k-i}(\zeta^{k-1}x_{k+1}, \ldots, \zeta^{2}x_4, \zeta x_3,x_2)\dd^{\set{3,4, \ldots, k}}_{\set{2, 3, \ldots, k}}(
       x_2^{k-1})\\
       &=(-1)^{n-1-k}\zeta^{-\frac{(n-k)(n-k+1)}{2}-n(k-1-i)}e_{k-i}(\zeta^{k-1}x_{k+1}, \ldots, \zeta^{2}x_4, \zeta x_3,x_2),  
     \end{align}
     Using the explicit description of the second tensor factor in \eqref{eq 3.5 lemma square flop 1 tensor x1k}, we have that the corresponding term in \eqref{eq 1 lemma square flop 1 tensor x1k} is:
     \begin{align}
      \notag  &(-1)^{n-1-k}\zeta^{-\frac{(n-k)(n-k+1)}{2}-n(k-1-i)}e_{k-i}(\zeta^{k-1}x_{k+1}, \ldots, \zeta x_3,x_2)\otimes (-1)^{n-1+(k-i)}\zeta^{\frac{n(n-1)}{2}-1+(k-i)}x_1^{i} \\
      &=(-1)^{i}\zeta^{\frac{n(n-1)}{2}-\frac{(n-k)(n-k+1)}{2}-n(k-1-i)-1+(k-i)}e_{k-i}(\zeta^{k-1}x_{k+1}, \ldots,  \zeta x_3,x_2)\otimes x_1^{i} \label{eq 7 lemma square flop 1 tensor x1k}.
     \end{align}
     
    In summary, the only non-zero terms satisfy $0\leq i \leq k-1$ and either $j=i$ or $j=k$. From \eqref{eq 6 lemma square flop 1 tensor x1k}, \eqref{eq 7 lemma square flop 1 tensor x1k}, the expression in \eqref{eq 1 lemma square flop 1 tensor x1k} simplifies to
    \begin{align*}
        &\sum_{\substack{0\leq i \leq k-1\\(j=i)}} (-1)^{k}q^{-2} \zeta^{\frac{n(n-1)}{2}-\frac{(n-k)(n-k+1)}{2}-n(k-1-i)-1} \;\;\;1\otimes x_1^k\hspace{50 pt} 
        \\
        +&\sum_{\substack{0\leq i \leq k-1\\(j=k)}} (-1)^{i}\zeta^{\frac{n(n-1)}{2}-\frac{(n-k)(n-k+1)}{2}-n(k-1-i)-1+(k-i)}e_{k-i}(x_2, \zeta x_3, \ldots, \zeta^{k-1}x_{k+1})\otimes x_1^{i}
    \end{align*}
    The expression above simplifies to
    \begin{align}
       \notag  &\;\;\;\sum_{i=0}^{k-1} (-1)^kq^{-2}\zeta^{-\frac{k(k-1)}{2}-1}\zeta^{ni} 1\otimes x_1^k\\
     \notag   &+\sum_{0\leq i \leq k-1} (-1)^{i}\zeta^{-\frac{k(k-1)}{2}-1+k-i}
        \zeta^{ni}e_{k-i}(x_2, \zeta x_3, \ldots, \zeta^{k-1}x_{k+1})\otimes x_1^{i}\\
    =&\;\;\;\notag (-1)^k\zeta^{-\frac{k(k-1)}{2}-1}q^{-k}[k-1]_q\; 1\otimes x_1^k\\
    &+\sum_{0\leq i \leq k} (-1)^{i}\zeta^{-\frac{k(k-1)}{2}-1+k-i}
        q^{-2i} e_{k-i}(x_2, \zeta x_3, \ldots, \zeta^{k-1}x_{k+1})\otimes x_1^{i}. \label{eq 8 lemma square flop 1 tensor x1k}   
    \end{align}
    Comparing the expression above (which is equal the one in \eqref{eq 1 lemma square flop 1 tensor x1k}, which was equal to the LHS of the equation in \cref{lemma square flop 1 tensor x1k}) with the RHS of the equation in \cref{lemma square flop 1 tensor x1k},
    we see that it suffices to show 
    \begin{align}
       \notag \zeta^{k-1}\ddm{0,1,k+1}{1,k+1}&\left( \ddm{0,1,k+1}{0, k+1}(x_1^k)\Del{0,1}{1}{1}\right) \otimes \Del{0,1}{1}{2}\\
        &=\sum_{0\leq i \leq k} (-1)^{i}\zeta^{-\frac{k(k-1)}{2}-1+k-i}
        q^{-2i} e_{k-i}(x_2, \zeta x_3, \ldots, \zeta^{k-1}x_{k+1})\otimes x_1^{i}.\label{eq 9 lemma square flop 1 tensor x1k}
    \end{align}
    Let us make the LHS of the equation above explicit by choosing $
    \Del{0,1}{1}{2}, \Del{0,1,}{1}{1}$ to run over $\set{(-1)^ih_i(x_1)\zeta^{\frac{n(n-1)}{2}}}$, $\set{e_{n-1-i}(\zeta^{-1}x_n, \ldots \zeta^{-(n-1)}x_2)}$ respectively for $0 \leq i \leq n-1$ (using \cref{corollary dual basis right ext} for explicit dual basis). 
    Moreover, 
    \[\ddm{0,1,k+1}{1,k+1}=\dd^{\set{2,3, \ldots,k}\sqcup \set{k+2, \ldots, n-1}}_{\set{2,3, \ldots, k}\sqcup\set{k+2, \ldots, n-1,n}}=\dd^{\set{k+2, \ldots,n-1}}_{\set{k+2,\ldots,n}} ,\]
    \[\ddm{0,1,k+1}{0, k+1}=\dd^{\set{2,3,\ldots,k}\sqcup\set{k+2,\ldots, n-1}}_{\set{1,2,\ldots, k}\sqcup \set{k+2,\ldots,n-1}}=\dd^{\set{2,3,\ldots,k}}_{\set{1,2,\ldots, k}}.\]
    Hence, we have (from \cref{lemma dual basis left ext}) that
    \begin{align*}
        \ddm{0,1,k+1}{0,k+1}(x_1^k)=\dd^{\set{2,3,\ldots,k}}_{\set{1,2,\ldots, k}}(x_1^k)=1.
    \end{align*}
    Hence  the LHS of \eqref{eq 9 lemma square flop 1 tensor x1k} simplifies to
    \begin{align}
        \notag &\zeta^{k-1}\ddm{0,1,k+1}{1,k+1}\left(\Del{0,1}{1}{1}\right) \otimes \Del{0,1}{1}{2}\\
        &=\sum_{0\leq i\leq n-1}\zeta^{k-1}\dd^{\set{k+2, \ldots,n-1}}_{\set{k+2,\ldots,n}}\left( e_{n-1-i}(\zeta^{-1}x_n, \ldots ,\zeta^{-(n-1)}x_2) \right) \otimes (-1)^ih_i(x_1)\zeta^{\frac{n(n-1)}{2}} \label{eq 10 lemma square flop 1 tensor x1k}
    \end{align}
    Since
    \begin{align*}
        &\;\;\;\;\;\;e_{n-1-i}(\zeta^{-1}x_n, \ldots ,\zeta^{-(n-1)}x_2)\\
        &=\sum_{l=0}^{n-1-i}e_l(\zeta^{-1}x_n, \zeta^{-2}x_{n-1}, \ldots ,\zeta^{-(n-k+1)}x_{k+2})e_{n-1-i-l}(\zeta^{-(n-k)}x_{k+1}, \ldots ,\zeta^{-(n-2)}x_3,\zeta^{-(n-1)}x_2),
    \end{align*}
    we have that
    \begin{align*}
      &\dd^{\set{k+2, \ldots,n-1}}_{\set{k+2,\ldots,n}}\left( e_{n-1-i}(\zeta^{-1}x_n, \ldots ,\zeta^{-(n-1)}x_2) \right)  \\
      &=\sum_{l=0}^{n-1-i}e_{n-1-i-l}(\zeta^{-(n-k)}x_{k+1}, \ldots ,\zeta^{-(n-1)}x_2)\dd^{\set{k+2, \ldots,n-1}}_{\set{k+2,\ldots,n}}\left( e_l(\zeta^{-1}x_n, \ldots ,\zeta^{-(n-k+1)}x_{k+2})\right)\\
      &=\sum_{l=0}^{n-1-i}e_{n-1-i-l}(\zeta^{-(n-k)}x_{k+1}, \ldots ,\zeta^{-(n-1)}x_2)\zeta^{-\frac{(n-k)(n-k-1)}{2}}\delta_{l, n-k-1} \hspace{20pt}\text{(Using \cref{corollary dual basis right ext})}\\
      &=\zeta^{-\frac{(n-k)(n-k-1)}{2}}e_{k-i}(\zeta^{-(n-k)}x_{k+1}, \ldots ,\zeta^{-(n-1)}x_2).
    \end{align*}
    Hence the RHS of \eqref{eq 10 lemma square flop 1 tensor x1k} can be simplified to
    \begin{align*}
      &\sum_{0\leq i\leq k}\zeta^{k-1}\zeta^{-\frac{(n-k)(n-k-1)}{2}}e_{k-i}(\zeta^{-(n-k)}x_{k+1}, \ldots ,\zeta^{-(n-1)}x_2) \otimes (-1)^ih_i(x_1)\zeta^{\frac{n(n-1)}{2}}  \\
      &=\sum_{0 \leq i \leq k}(-1)^i\zeta^{\frac{n(n-1)}{2}-\frac{(n-k)(n-k-1)}{2}+k-1}\zeta^{-(n-1)(k-i)}e_{k-i}(\zeta^{k-1} x_{k+1}, \ldots ,x_2) \otimes x_1^i\\
     &=\sum_{0\leq i \leq k} (-1)^{i}\zeta^{-\frac{k(k-1)}{2}-1+k-i}
        q^{-2i} e_{k-i}(x_2, \zeta x_3, \ldots, \zeta^{k-1}x_{k+1})\otimes x_1^{i}
    \end{align*}
  which is precisely the RHS of \eqref{eq 9 lemma square flop 1 tensor x1k}.
    \end{proof}
    
\begin{lemma}\label{lemma square flop 1 tensor x1k for k=n-1}
The following equality holds in $\Bimod{\hh{0,1}}{}{R^{\hh{0,1}}}\otimes_{\hh{1}}\Bimod{}{\hh{1,2}}{R^{\hh{1,2}}}$.
    \begin{align*}
    &\ddm{0,1,2}{0,1}\left(\ddm{0,1,2}{0,2}(\Del{1,2}{2}{1})\Del{1,2}{1}{1}\right) \otimes \ddm{0,1,2}{0,1}\left(\ddm{0,1,2}{0,2}(\Del{1,2}{2}{2}x_1^{n-1})\Del{1,2}{1}{2} \right) \\
    =& (-1)^{n-1}\zeta^{\frac{-(n-1)(n-2)}{2}-1}q^{-(n-1)}\left([k-1]_q(1\otimes x_1^{n-1})
    \;+\; (-1)^{n-1}\ddm{0,1}{0}(x_1^{n-1})\Del{0,1}{1}{1} \otimes \Del{0,1}{1}{2} \right)
    \end{align*}
\end{lemma}
\begin{proof}
    The LHS simplifies exactly as in \cref{lemma square flop 1 tensor x1k} for $k=n-1$ to the expression in \cref{eq 8 lemma square flop 1 tensor x1k}, so that
    \begin{align}
    \text{LHS} =&\;\;\;\notag (-1)^{n-1}\zeta^{-\frac{(n-1)(n-2)}{2}-1}q^{-(n-1)}[n-2]_q\; 1\otimes x_1^{n-1}\\
    &+\sum_{0\leq i \leq n-1} (-1)^{i}\zeta^{-\frac{(n-1)(n-2)}{2}+n-2-i}
        q^{-2i} e_{n-1-i}(x_2, \zeta x_3, \ldots, \zeta^{n-2}x_{n})\otimes x_1^{i}.
    \end{align}
    Comparing with the RHS, we see that it suffices to show
    \begin{align}
        \notag &\zeta^{\frac{-(n-1)(n-2)}{2}-1}q^{-(n-1)}\ddm{0,1}{0}(x_1^{n-1})\Del{0,1}{1}{1} \otimes \Del{0,1}{1}{2} \\
        &\quad = \sum_{0\leq i \leq n-1} \left(
        (-1)^{i}\zeta^{-\frac{(n-1)(n-2)}{2}+n-2-i}
        q^{-2i}\; e_{n-1-i}(x_2, \zeta x_3, \ldots, \zeta^{n-2}x_{n})\otimes x_1^{i}\right)
    \end{align}
    which is equivalent to showing
    \begin{align}\label{eq 1 lemma square flop 1 tensor x1k for k=n-1}
       q^{-(n-1)}\ddm{0,1}{0}(x_1^{n-1})\Del{0,1}{1}{1} \otimes \Del{0,1}{1}{2} &= \sum_{0\leq i \leq n-1} \left(
        \begin{array}{l}(-1)^{i}\zeta^{n-1-i}
        q^{-2i}\vspace{5 pt}\\
        \quad e_{n-1-i}(x_2, \zeta x_3, \ldots, \zeta^{n-2}x_{n})\otimes x_1^{i}\end{array}\right)
    \end{align}
    By \cref{lemma dual basis left ext}, we know that $\ddm{0,1}{0}(x_1^{n-1})=1$. 
    From \cref{corollary dual basis right ext}, we know that we may choose
    $\Del{0,1}{1}{1}$ to run over $\set{e_{n-1-i}(\zeta^{-1} x_{n}, \zeta^{-2}x_{n-1},\ldots, \zeta^{-(n-1)}x_2 )}_{0 \leq i \leq n-1}$ and $\Del{0,1}{1}{2}$ to run over $\set{(-1)^i\zeta^{\frac{n(n-1)}{2}}x_1^i}_{0\leq i \leq n-1}$.
    So we may expand the LHS of \cref{eq 1 lemma square flop 1 tensor x1k for k=n-1} as
    \begin{align}
        \nonumber q^{-(n-1)}\sum_{0\leq i \leq n-1}&e_{n-1-i}(\zeta^{-1} x_{n}, \zeta^{-2}x_{n-1},\ldots, \zeta^{-(n-1)}x_2 )\otimes (-1)^i\zeta^{\frac{n(n-1)}{2}}x_1^i\\
       \nonumber &=\sum_{0 \leq i\leq n-1}(-1)^i\zeta^{\frac{n(n-1)}{2}}q^{-(n-1)}e_{n-1-i}(\zeta^{-1} x_{n}, \zeta^{-2}x_{n-1},\ldots, \zeta^{-(n-1)}x_2 )\otimes x_1^i\\
       \notag &=\sum_{0\leq i \leq n-1}(-1)^i\zeta^{n(n-1)}\zeta^{-(n-1)(n-1-i)}e_{n-1-i}(\zeta^{n-2} x_{n}, ,\ldots, \zeta x_3, x_2 ) \otimes x_1^i\\
       \notag  &=\sum_{0 \leq i \leq n-1} (-1)^i\zeta^{(n-1)(1+i)}e_{n-1-i}(x_2, \zeta x_3, \ldots, \zeta^{n-2}x_{n} ) \otimes x_1^i\\
       \notag &=\sum_{0 \leq i \leq n-1} (-1)^i\zeta^{n-1-i}q^{-2i}e_{n-1-i}(x_2, \zeta x_3, \ldots, \zeta^{n-2}x_{n} ) \otimes x_1^i
    \end{align}
    which is precisely the RHS of \cref{eq 1 lemma square flop 1 tensor x1k for k=n-1}.
\end{proof}

%% file: The_Diagrammatic_Quantum_Satake_Functor.tex
In this section, we introduce the combinatorial and diagrammatic replacements for both sides of the Quantum Geometric Satake functor: the category $\cwebs$ of webs, describing $\Fund_q^{\Om} \subset \Rep^{\Om}(U_q(\mathfrak{sl}_n))$, and the category $\DiagSBS$ of singular Soergel diagrams, describing $\SBSBim \subset \SSBim$. We then  describe the Quantum Geometric Satake $2$-functor in \cref{subsec the functor} and address fullness and faithfulness in \cref{subsec fullness and faithfulness}, while also showing that it categorifies the reformulated Satake isomorphism. 
The proof that this $2$-functor is well-defined will be postponed to \cref{sec well-defined-ness of the functor}. See \cref{thm main} for a summary of the key results.

\subsection{The fundamental subcategory}\label{subsec the fundamental subcategory}

\begin{defn} Let $\Fund_q = \Fund_q(\mathfrak{sl}_n)$) denote the full monoidal subcategory of $\Rep_q = \Rep(U_q(\mathfrak{sl}_n))$ monoidally generated by the fundamental representations. It is a $\Q(q)$-linear category for a formal variable $q$ (though not an additive category, as we have not included direct sums). \end{defn}

The objects of $\Fund$ correspond to fundwords, in the sense of \cref{defn:fundword}. The Karoubi envelope of $\Fund_q$ is $\Rep_q$. Had we defined the analogous category over specializations of $\Z[q,q^{-1}]$ to other fields, the Karoubi envelope of $\Fund_q$ would agree with the category of tilting modules over that field.

As explained in \cref{subsec Satake Hecke algebroid}, $\Rep_q$ is an $\Om$-graded monoidal category, corresponding to which we constructed a 2-category $\qRep$ in \cref{defn- qRep Omega}. Being a full monoidal subcategory of $\Rep_q$, $\Fund_q$ is an $\Om$-graded-monoidal category as well. Now we define the corresponding sub-2-category of $\qRep$.


\begin{definition}\label{defn- qFund Omega}
 The $2$-category $\qFund$ is defined as follows:
  \begin{itemize}
      \item The set of objects is $\Omega$.
      \item The Hom-category $\Hom(l,k)$ is defined as $\Fund_q^{k-l} := \Rep(U_q(\mathfrak{sl}_n))^{k-l} \cap \Fund_q$.
\end{itemize}
      
\end{definition}

In other words, $\qFund$ is a 2-full sub-2-category of $\qRep$ where the 1-morphisms are generated monoidally by the fundamental representations. The $1$-morphisms of $\qFund$ correspond to colored fundwords in the sense of \cref{def:GSfundword}.

\subsection{Colored \texorpdfstring{$\mathfrak{sl_n}$}{sl\_n}-Webs}\label{subsec-colored sln webs}


A diagrammatic presentation for $\Fund_q(\mathfrak{sl}_n)$ over $\C(q)$ was first given by Cautis--Kamnitzer--Morrison in \cite{CKM}. In this paper, we use a different presentation from work of Bigelow \cite{Bigelow}, where the base field is also assumed to be $\C(q)$.  More accurately, both \cite{CKM} and \cite{Bigelow} describe a category which is larger than $\Fund_q$, having as generating objects both the fundamental representations and their duals. However, $\Fund_q$ is equivalent to this larger category (the inclusion functor is essentially surjective), as the duals of fundamental representations are isomorphic to other fundamental representations. Diagrammatically, these isomorphisms are drawn as tags.

In \cite[Proposition 1]{Bigelow}, Bigelow shows that the web relations in \cite{CKM} follow from this smaller collection of relations over $\C(q)$, and the proof only requires that the quantum integers are invertible. In particular, the proofs in \cite{CKM} and \cite{Bigelow} both work over $\Q(q)$.

\begin{figure}[h]
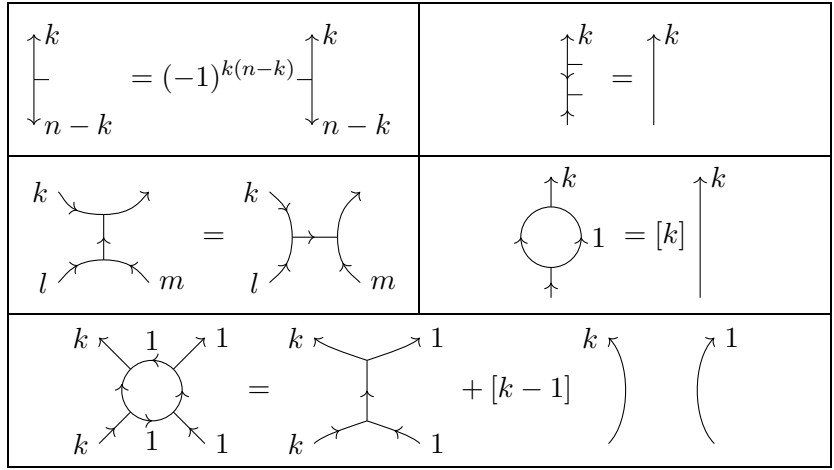

  \[
    \begin{array}{|T|T|}
      \hline
      {
        \phantom{ \tikz[x=1mm, y=1mm, baseline=-0.5ex]{\draw (0,-10) -- (0,10)}} 
        \tikz[x=1mm,
  y=1mm,
  baseline=-0.5ex] {
          \draw[->] (0,0) -- (0,6) node [right] {$k$};
          \draw[->] (0,0) -- (0,-6) node [right] {$n - k$};
          \draw (0,0) -- (2,0);
        }
        = (-1)^{k(n - k)} \tikz[x=1mm, y=1mm, baseline=-0.5ex]{
          \draw[->] (0,0) -- (0,6) node [right] {$k$};
          \draw[->] (0,0) -- (0,-6) node [right] {$n - k$};
          \draw (0,0) -- (-2,0);
        }
      } &
      {
        \tikz[x=1mm,
         y=1mm,
        baseline=-0.5ex] {
          \draw[->-] (0,-6) -- (0,-2);
          \draw (0,-2) -- (2,-2);
          \draw[->-] (0,2) -- (0,-2);
          \draw (0,2) -- (2,2);
          \draw[->] (0,2) -- (0,6) node [right] {$k$};
        }
        = \, \tikz[  x=1mm,
  y=1mm,
  baseline=-0.5ex] {
          \draw[->] (0,-6) -- (0,6) node [right] {$k$};
        }
      } \\ \hline
      {
        \phantom{ \tikz[x=1mm, y=1mm, baseline=-0.5ex] {\draw (0,-10) -- (0,10)}} 
        \tikz[x=1mm, y=1mm, baseline=-0.5ex] {
          \draw[->-] (-6,-6) node [left] {$l$}
            .. controls (-5,-5) and (-4,-3) .. (0,-3);
          \draw[->-] (6,-6) node [right] {$m$}
            .. controls (5,-5) and (4,-3) .. (0,-3);
          \draw[->-] (0,-3) -- (0,3);
          \draw[->-] (-6,6) node [left] {$k$}
            .. controls (-5,5) and (-5,3) .. (0,3);
          \draw[->] (0,3) .. controls (4,3) and (5,5) .. (6,6);
        }
        = \tikz[x=1mm, y=1mm, baseline=-0.5ex] {
          \draw[->-] (-6,-6) node [left] {$l$}
            .. controls (-5,-5) and (-3,-4) .. (-3,0);
          \draw[->-] (-6,6) node [left] {$k$}
            .. controls (-5,5) and (-3,4) .. (-3,0);
          \draw[->-] (-3,0) -- (3,0);
          \draw[->-] (6,-6) node [right] {$m$}
            .. controls (5,-5) and (3,-4) .. (3,0);
          \draw[->] (3,0) .. controls (3,4) and (5,5) .. (6,6);
        }
      } &
      {
        \tikz[x=1mm, y=1mm, baseline=-0.5ex] {
          \draw[->-] (0,-8) -- (0,-4);
          \draw[->-] (0,-4) arc (270:90:4);
          \draw[->-] (0,-4) arc (-90:90:4) node[pos=0.5,right]{$1$};
          \draw[->] (0,4) -- (0,8) node [right] {$k$};
        }
        = [k] \: \tikz[x=1mm, y=1mm, baseline=-0.5ex] {
          \draw[->] (0,-8) -- (0,8) node [right] {$k$};
        }
      } \\ \hline
    \multicolumn{2}{|F|}
    {
      \phantom{ \tikz[x=1mm, y=1mm, baseline=-0.5ex] {\draw (0,-10) -- (0,10)}} 
      \tikz[x=1mm, y=1mm, baseline=-0.5ex] {
        \draw[->-] (45:4) arc (45:135:4) node [pos=0.5, above]{$1$};
        \draw[-<-] (135:4) arc (135:225:4);
        \draw[->-] (225:4) arc (225:315:4) node [pos=0.5, below]{$1$};
        \draw[->-] (315:4) arc (315:405:4);
        \draw[->-] (7,-7) node [right] {$1$} -- (315:4);
        \draw[->] (45:4) -- (7,7) node [right] {$1$};
        \draw[->] (135:4) -- (-7,7) node [left] {$k$};
        \draw[->-] (-7,-7) node [left] {$k$} -- (225:4);
      }
      =
      \tikz[x=1mm, y=1mm, baseline=-0.5ex] {
        \draw[->-] (-7,-7) node [left] {$k$}
          .. controls (-6,-6) and (-3,-5) .. (0,-4);
        \draw[->-] (7,-7) node [right] {$1$}
          .. controls (6,-6) and (3,-5) .. (0,-4);
        \draw[->-] (0,-4) -- (0,4);
        \draw[->] (0,4) .. controls (-3,5) and (-6,6)
          .. (-7,7) node [left] {$k$};
        \draw[->] (0,4) .. controls (3,5) and (6,6)
          .. (7,7) node [right] {$1$};
      }
      + [k - 1] \tikz[x=1mm, y=1mm, baseline=-0.5ex] {
        \draw[->] (-7,-7) .. controls (-4,-4) and (-4,4)
          .. (-7,7) node [left] {$k$};
        \draw[->] (7,-7) .. controls (4,-4) and (4,4)
          .. (7,7) node [right] {$1$};
      }
    } \\ \hline
\end{array}
\]
\caption{Web relations. Some labels are omitted, but can be deduced.}
\label{fig:webrel}
\end{figure}

Bigelow's relations are those shown in \cref{fig:webrel}, up to taking mirror images and arrow reversals. While the generating morphisms in Bigelow's presentation look the same as in Cautis--Kamnitzer--Morrison's presentation, the meaning of some of the generating morphisms is different; they have different conventions for the ``tag.'' In particular, Bigelow's tags are defined in such a way that his webs admit both reflection symmetry and arrow reversal symmetry. The webs in \cite{CKM} are invariant under simultaneously reflecting and reversing arrows, but not under either individual operation.

\begin{remark} Bigelow ``fixes'' a subtle sign error in \cite{CKM}. The relations explicitly written in \cite{CKM} are correct, but the paper incorrectly asserts that the orientation-reversed relations also hold; this fails for equations (2.7) and (2.8) therein, being off by a sign, as Bigelow observed.  \end{remark}

\begin{definition}\label{defn-webs}
    The category of (Bigelow-style) $\mathfrak{sl}_n$-webs, denoted $\Webs$, is defined as follows.
    \begin{itemize}
        \item The objects are generated monoidally by oriented points labelled by $i\in\set{1,2, \ldots, n-1}$.
        \item The morphisms are generated ($\Q(q)$-linearly) by clockwise and counter-clockwise cups and caps, along with bivalent vertices (tags) and trivalent vertices as shown below:
        \[
            \begin{array}{OOOO}
        \phantom{ \tikz[x=1mm, y=1mm, baseline=-0.5ex]{\draw (0,-10) -- (0,10)}} 
         \tikz[x=1mm,
         y=1mm,
         baseline=-0.5ex] {
          \draw[->-] (0,-5) arc (0:180:5) node [pos=0.5, above] {$k$};
        } 
        &
        \tikz[x=1mm,
         y=1mm,
         baseline=-0.5ex] {
          \draw[-<-] (0,-5) arc (0:180:5) node [pos=0.5, above] {$k$};
        } 
        &\tikz[x=1mm,
         y=1mm,
         baseline=-0.5ex] {
          \draw[->-] (0,5) arc (180:360:5) node [pos=0.5, below] {$k$};
        } 
        &\tikz[x=1mm,
         y=1mm,
         baseline=-0.5ex] {
          \draw[-<-] (0,5) arc (180:360:5) node [pos=0.5, below] {$k$};
        } 
        \\
        \phantom{ \tikz[x=1mm, y=1mm, baseline=-0.5ex]{\draw (0,-10) -- (0,10)}} 
        \tikz[x=1mm,
         y=1mm,
         baseline=-0.5ex] {
          \draw[->-] (0,0) -- (0,6) node [above] {$k$};
          \draw[->-] (0,0) -- (0,-6) node [below] {$n - k$};
          \draw (0,0) -- (2,0);
        } 
        & 
        \tikz[x=1mm,
         y=1mm,
         baseline=-0.5ex] {
          \draw[-<-] (0,0) -- (0,6) node [above] {$k$};
          \draw[-<-] (0,0) -- (0,-6) node [below] {$n - k$};
          \draw (0,0) -- (2,0);
        }
        &
        \tikz[x=1mm, y=1mm, baseline=-0.5ex] {
        \draw[->-] (-6,-6) node [below] {$k$}
          .. controls (-5,-5) and (-3,-3) .. (0,0);
        \draw[->-] (6,-6) node [below] {$l$}
          .. controls (5,-5) and (3,-3) .. (0,0);
        \draw[->-] (0,0) -- (0,6) node[above]{$k+l$};
        }
        &
        \tikz[x=1mm, y=1mm, baseline=-0.5ex] {
        \draw[->-] (0,-6) node[below]{$k+l$} -- (0,0);
        \draw[-<-] (-6,6) node [above] {$k$} -- (0,0);
        \draw[-<-] (6,6) node [above] {$l$} -- (0,0);
      }
            \end{array}
            \]
        
        \item The cups and caps satisfy the isotopy relations, making the category pivotal (see \cite[Chapter 2]{Morr07} for further details.)
        \item  The morphisms satisfy the relations in \cref{fig:webrel}, and their mirror images and arrow reversals.
        As in \cite{Bigelow}, we shall refer to these relations in \cref{fig:webrel} as ``tag switching", ``tag cancellation”, ``I = H", ``bursting a digon" and ``bursting a square" respectively (from top to bottom, left to right).
        \item It is convenient to also allow strands labeled $0$ or $n$, with the following convention for deleting them. Any trivalent vertex that had a strand labeled $0$ simply ceases to be a vertex. Any trivalent vertex that had a strand labeled $n$ becomes a bivalent vertex, with a tag on the side of the deleted strand.
        
    \end{itemize}
\end{definition}


\begin{theorem}\label{thm webs present Fund}
    There is a $\Q(q)$-linear equivalence $\Fund_q \to \Webs$ which is not surjective on objects, but is essentially surjective. There is a fully faithful $\Q(q)$-linear functor $\Webs \rightarrow \Rep_q$. The composition of these functors is the inclusion functor $\Fund_q \subset \Rep_q$.
\end{theorem}

\begin{proof}  An equivalence between Cautis-Kamnitzer-Morrison's web category and the full subcategory of $\Rep_q$ generated by fundamental representations and their duals was proven in \cite[Theorem 3.3.1]{CKM}. We merely restrict the inverse equivalence to $\Fund_q$, and translate the result through Bigelow's equivalence. \end{proof}

\begin{remark} The analogous result for any $\Z[q,q^{-1}]$-algebra, when using Cautis-Kamnitzer-Morrison's webs, was proven in \cite{ELLCC}. \end{remark}


This diagrammatic description of $\Fund_q(\mathfrak{sl}_n)$ naturally provides a diagrammatic description of the associated 2-category $\qFund$.

\begin{definition}
The 2-category of \emph{colored $\mathfrak{sl_n}-webs$}, denoted $\cwebs$ is defined as follows:
\begin{itemize}
    \item The objects are elements of $\Om$. 
    These objects label the regions in diagrams. The region labels will be colored in gray.
    \item For each $a \in \Om$ and $i \in \set{1,2,\ldots, n-1}$, there are 1-morphisms $i \in \Hom(a, a+i)$ and $i^* \in \Hom(a, a-i)$ which generate the category monoidally.
    We draw the identity $2$-morphism of $i \in \Hom(a,a+i)$ as an upward arrow, and that of $i^* \in \Hom(a, a-i)$ as a downward arrow as shown below.
    \begin{equation*}
        \gr{a+i}\uparrow^{\blk{i}} \gr{a}, \qquad \qquad
        \gr{a-i}\downarrow^{\blk{i}} \gr{a}
    \end{equation*}
    \item The 2-morphisms are generated by clockwise and counter-clockwise cups and caps, along with bivalent vertices (tags) and trivalent vertices as shown below:
    \[\begin{array}{cccc}
      \vcenter{\xy (0,0)*{\def\svgscale{0.15}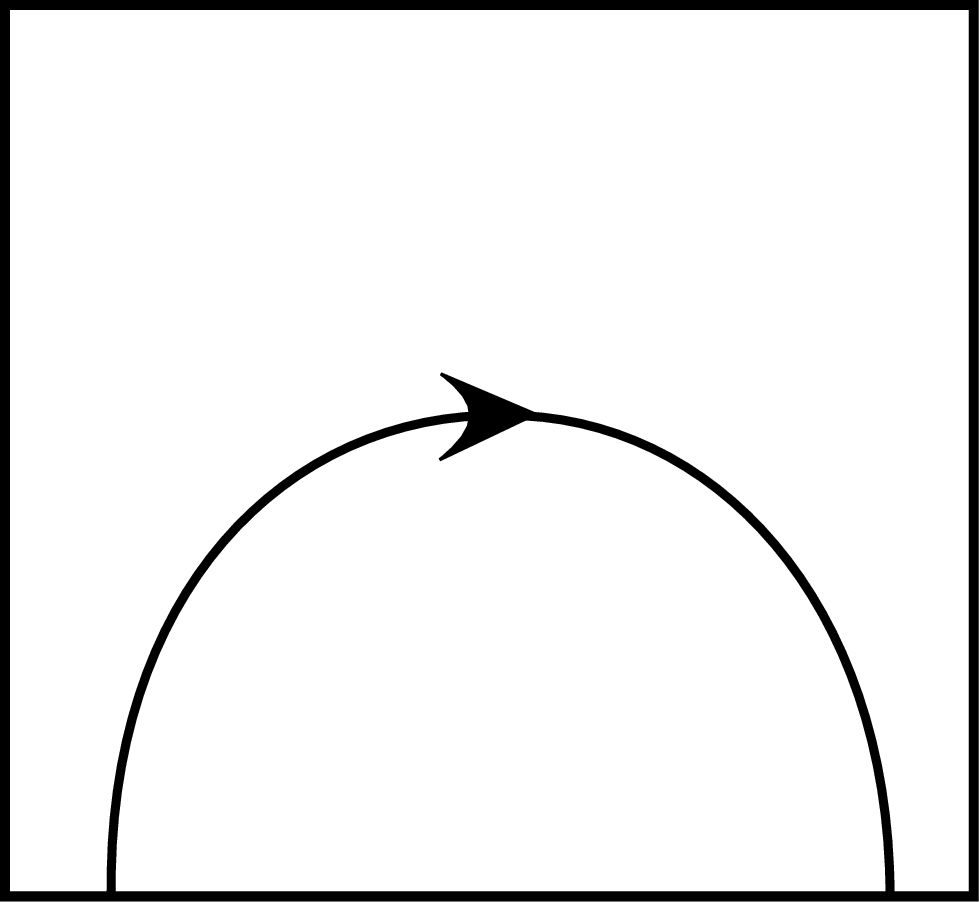} \endxy}   &\vcenter{\xy (0,0)*{\def\svgscale{0.15}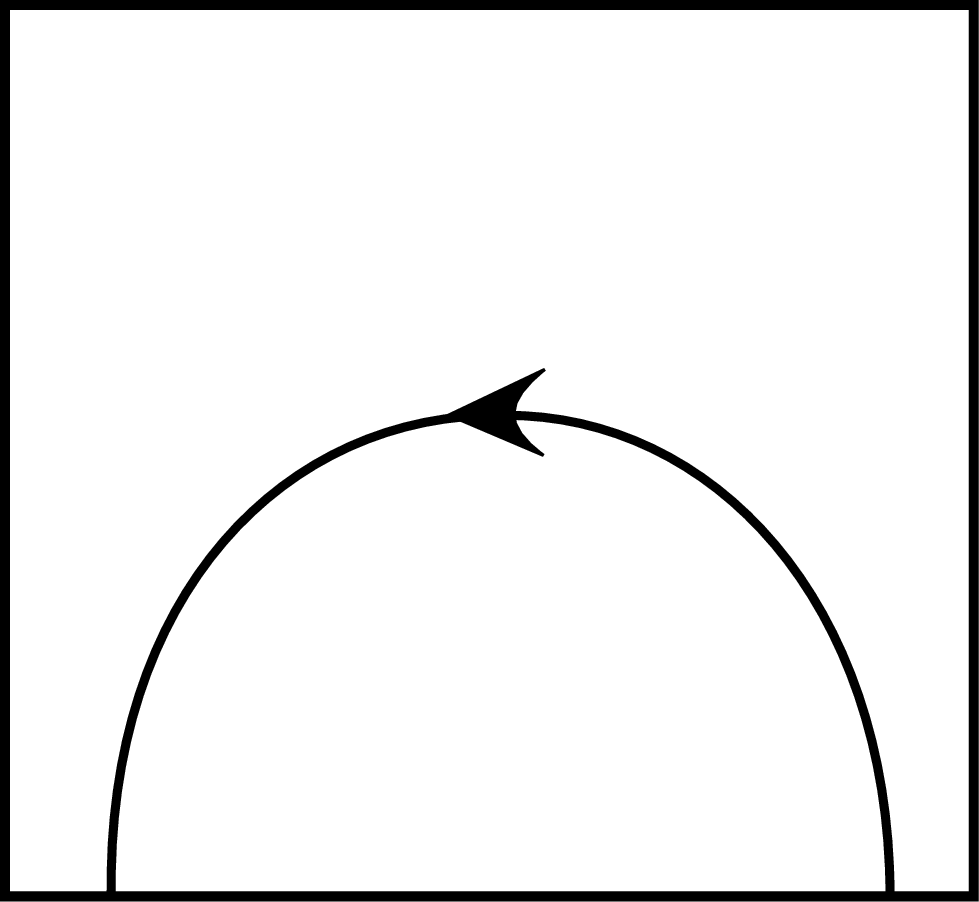} \endxy} &\vcenter{\xy (0,0)*{\def\svgscale{0.15}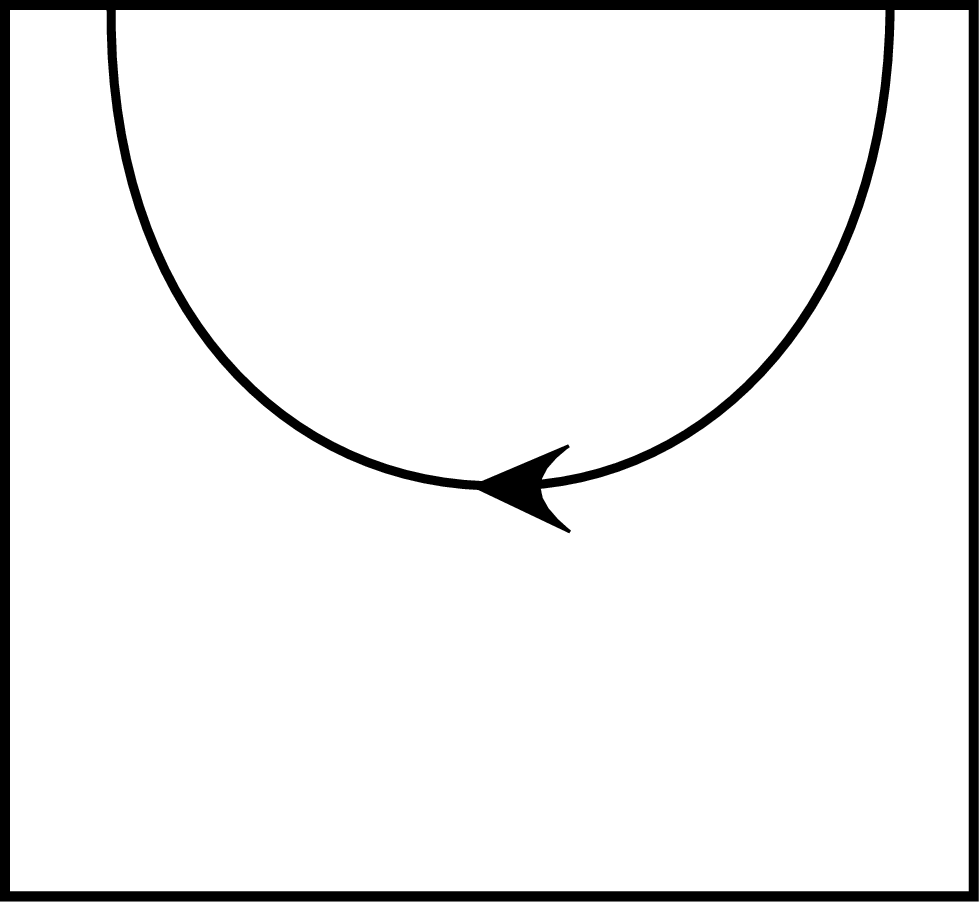} \endxy} &\vcenter{\xy (0,0)*{\def\svgscale{0.15}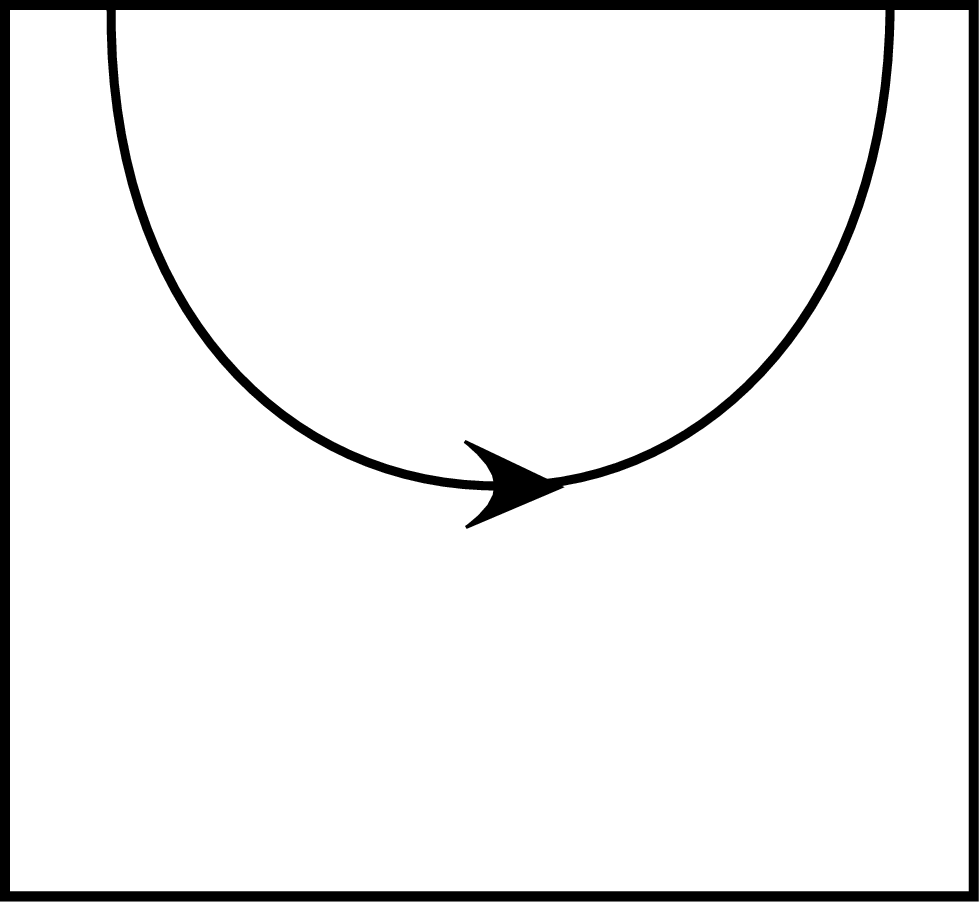} \endxy} 
    \end{array}
    \]
    \[\begin{array}{cccc}
        \vcenter{\xy (0,0)*{\def\svgscale{0.15}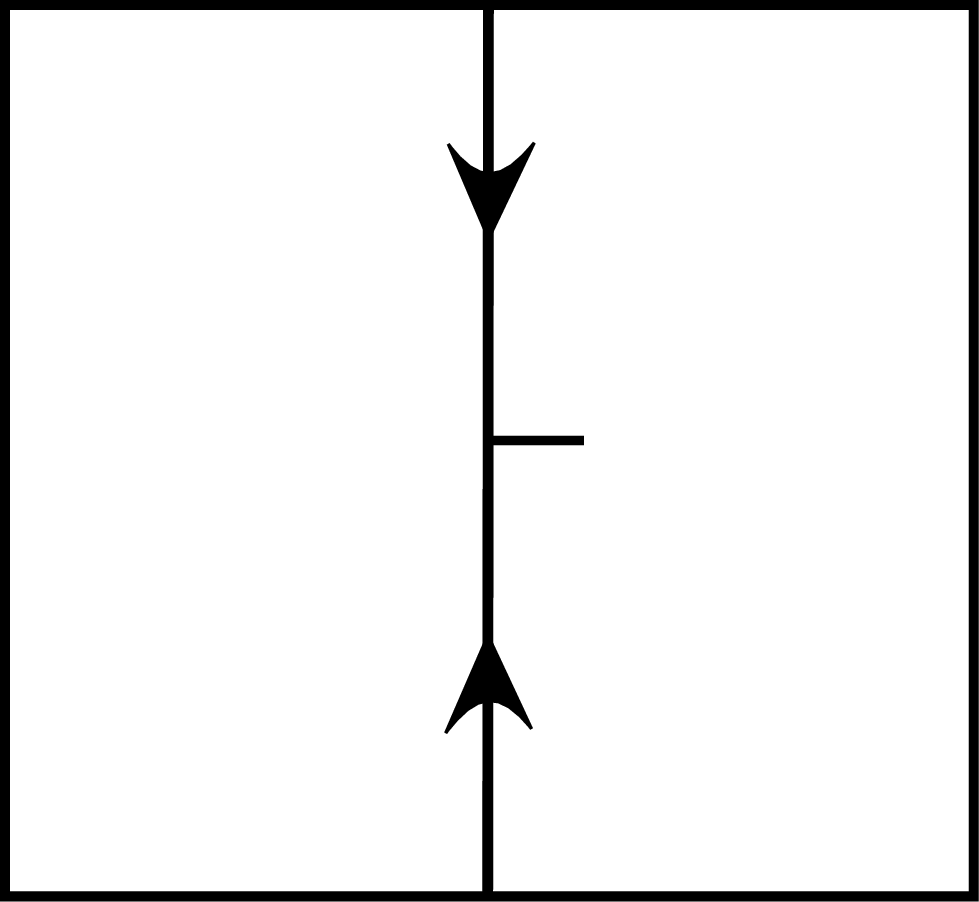} \endxy} &\vcenter{\xy (0,0)*{\def\svgscale{0.15}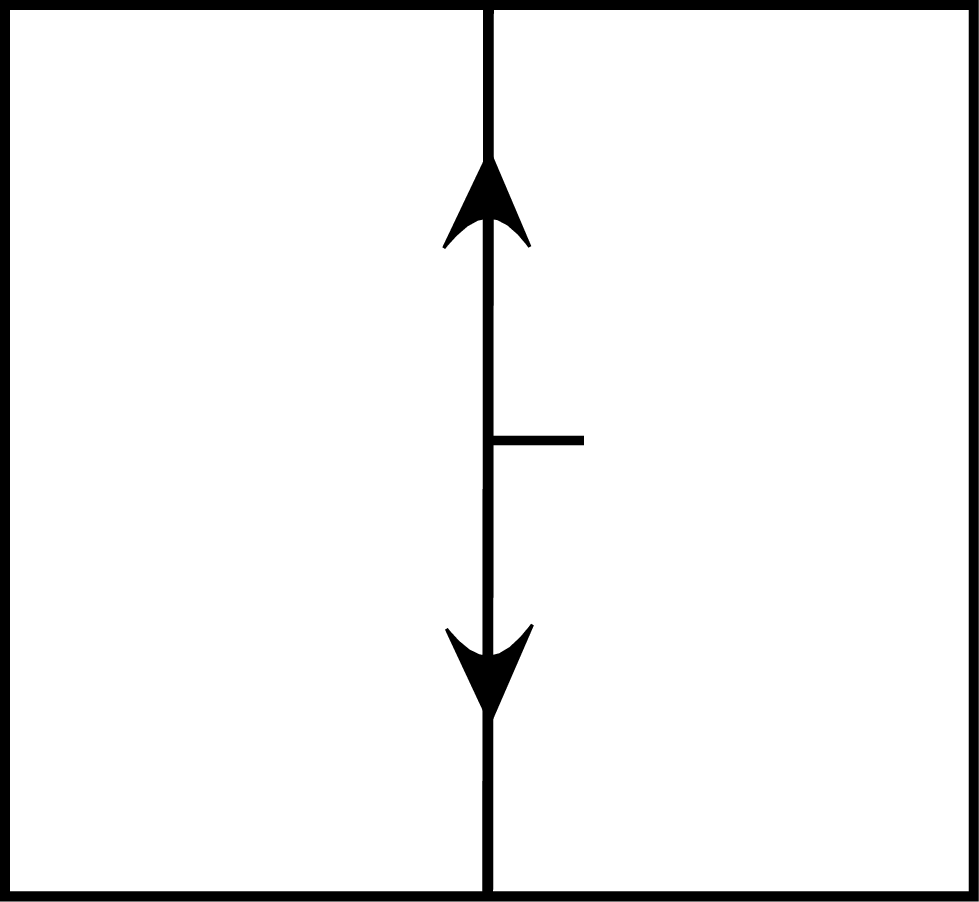} \endxy}  &\vcenter{\xy (0,0)*{\def\svgscale{0.15}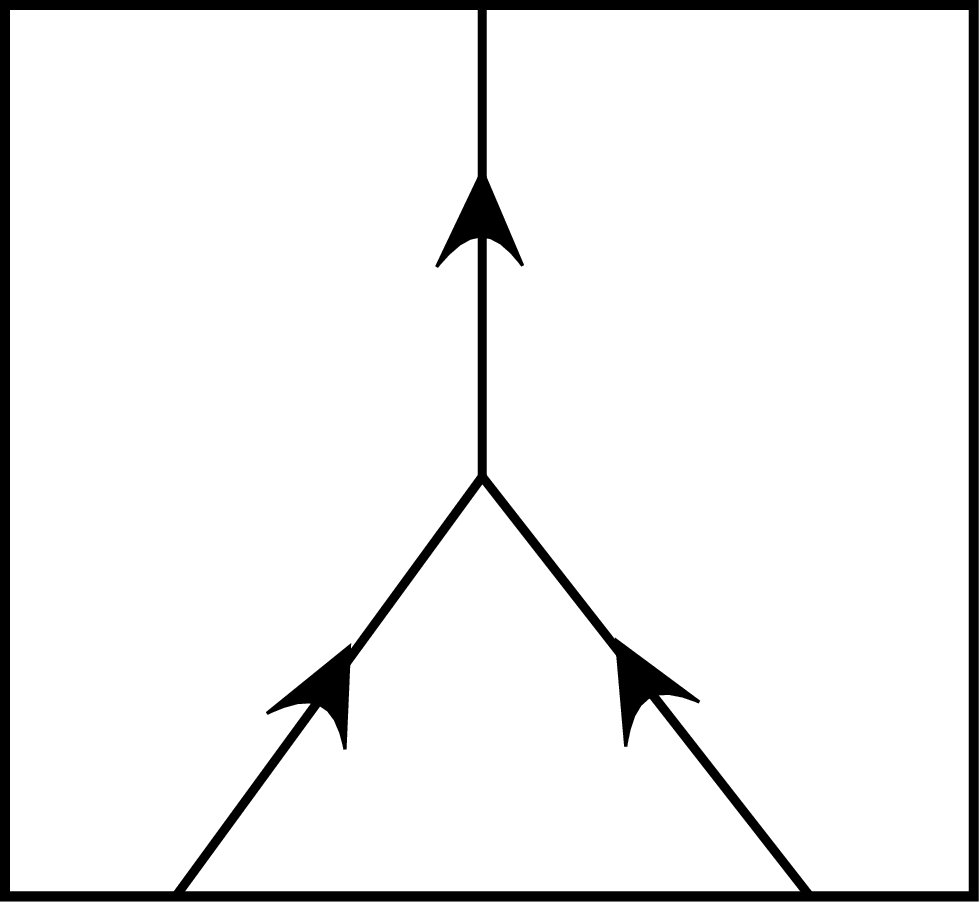} \endxy} &\vcenter{\xy (0,0)*{\def\svgscale{0.15}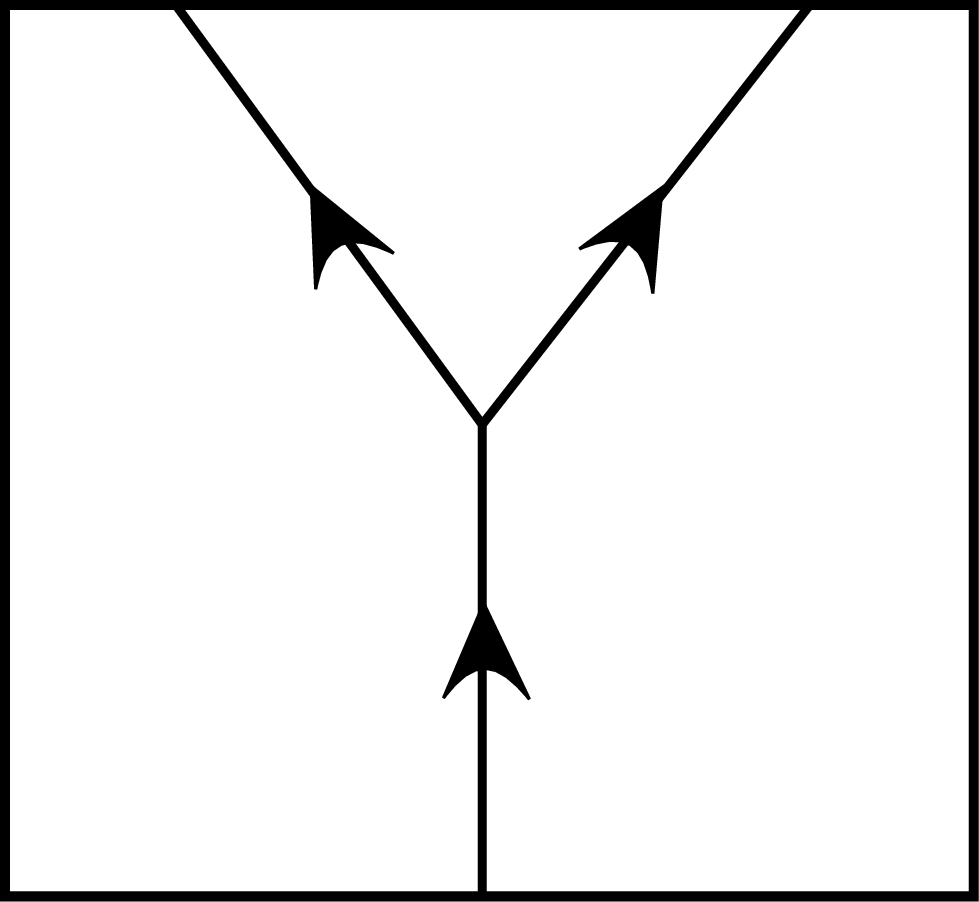} \endxy}
    \end{array}\]
    \item The 2-morphisms satisfy the relations in \cref{fig:webrel} and their mirror images and arrow reversals, for any allowed labeling of the regions (an upward arrow labelled $k$ should have a region $a+k$ to the left and $a$ to the right and so on). The cups and caps satisfy the isotopy relations as in \cref{isotopy}.
    \item Same conventions as in \cref{defn-webs} for deleting strands labeled $0$ or $n$.
\end{itemize}
\end{definition}

\begin{theorem}\label{theorem cwebs present qFund}
    There is a $\Q(q)$-linear $2$-functor $\Fund^\Omega_q \to \cwebs$ which is an equivalence; it is not surjective on $1$-morphisms, but is essentially surjective. There is a fully faithful $\Q(q)$-linear $2$-functor $\cwebs \rightarrow \Rep^\Omega_q$. The composition of these functors is the inclusion functor $\Fund^\Omega_q \subset \Rep^\Omega_q$.
\end{theorem}

\begin{defn} \label{defn:symmetries of webs} We define $2$-functors $\sigma, \rotation, \Upsilon_W, \reversal \colon \cwebs \to \cwebs$, the basic symmetries of this category.

On objects/region labels $\sigma$ sends $a \mapsto a+1$. Aside from this relabling, diagrams remain unchanged. This is a monoidal, covariant equivalence.

On objects $\rotation$ sends $a \mapsto a$, and it rotates all diagrams by 180 degrees. This is an antimonoidal, contravariant equivalence.

On objects $\Upsilon_W$ sends $a \mapsto -a$, and it flips each diagram vertically (i.e. reflecting across a horizontal axis). This is a monoidal, contravariant equivalence.

On objects $\reversal$ sends $a \mapsto -a$, and it reverses all the orientations in each diagram. This is a monodial, covariant equivalence.
\end{defn}

\subsection{Diagrammatics for Singular Bott-Samelson Bimodules}
\label{subsec-diagrammatics for frob hypercube}
\label{subsec-diagrammatic singular bott-samelson bimodules}

In \cite{EWSFrob}, it was explained how to construct a diagrammatic category $\DiagSBS(\Gamma)$, an algebraic category $\mc{C}(\Gamma)$ and a functor $\Phi: \DiagSBS(\Gamma) \rightarrow \mc{C}(\Gamma)$, starting from a Frobenius hypercube $\Gamma$. We describe these constructions for the (partial graded) Frobenius hypercube $\Gamma$ associated to the deformed affine Frealization in \cref{defn deformed affine frobenius realization}. The algebraic category obtained via this construction coincides with the category $\SBSBim$ described in \cref{subsec- SSBim}\footnote{To be precise, $C(\Gamma)$ is a sub-2-category of $\SBSBim$ which is only essentially full on 1-morphisms; see the discussion after \cite[Definition 1.3]{EWSFrob}.}, hence the rest of this section is primarily focused on describing the diagrammatic category $\DiagSBS(\Gamma)$ associated to the deformed affine Frealization and the functor $\Phi: \DiagSBS(\Gamma) \rightarrow \SBSBim$. 
There are certain relations in \cite{EWSFrob} that require the additional assumptions from \cref{defn-condition star}, \cref{defn- no mu zero deivisors}, and \cref{defn-condition R3}. By \cref{lemma deformed affine frobenius satisfies extra conditions}, these assumptions hold for the deformed affine Frealization, so these relations  are included when defining $\DiagSBS(\Gamma)$.

{


}

The construction is as follows.
Let $\mc{F}$ denote the set of finitary subsets of $S$.
The regions of $\mc{D}(\Gamma)$ are labeled by the finitary subsets $I \in \mc{F}$.
We typically label the regions in grey in small font, for example $\gr{I}$ for $I \in \mc{F}$ etc.
The functor $\Phi$ sends the region labeled by $I$ to the object $R^I$ in $\SBSBim$ (defined in \cref{subsec- SSBim}).

For $I \subset \Gamma$, $i,j\in \Gamma \setminus I$, we shall use the notation $Ii$ to refer to $I\sqcup \set{i}$, $Iij$ to refer to $I\sqcup \set{i,j}$ etc.
The generating 1-morphisms are $\Ind_{Ii}^{I}$ and $\Res_{Ii}^{I}$ for each $Ii \in \mc{F}$. 
We think of these generators as colored oriented points, where we fix a color to each $i \in S$.
The orientation will be upward for $\Ind_{Ii}^{I}$ with the region $Ii$ to the right (and $I$ to the left), while orientation is downward for $\Res_{Ii}^{I}$ with region $Ii$ to the left.
It follows that the space of 1-morphisms 
$\Hom(I,J)$ in $\DiagSBS(\Gamma)$ is exactly the set of singlestep expressions as in \cref{singlestepprototype}, with $I_0=J$ and $I_d=I$.

The functor $\Phi$ sends $\Ind_{Ii}^{I}$ to the induction bimodule $\Bimod{R^{I}}{R^{Ii}}{R^{I}}$ and $\Res_{Ii}^{I}$ to the 
shifted restriction bimodule $\Bimod{R^{Ii}}{R^{I}}{R^{I}} (\ell(Ii)-\ell(I))$.
More generally, $\Phi$ sends a singlestep expression $I_\bullet$ to the singular Bott-Samelson bimodule $\BS(I_\bullet)$.

The 2-morphisms are generated monoidally by oriented cups and caps (we have 4 of these for each $Ii$, $I$ in $\mc{F}$ ), polynomials $f \in R^I$ in a region labeled $I$, and crossings (upward and downward) for $Iij, Ij, Ii$ in $\mc{F}$ as shown below (where $i, j$ are colored blue and red respectively). The images of these generating 2-morphisms after applying $\Phi$ are given below the respective diagrams (suppressing grading shifts, and using standard canonical isomorphisms for convenience). 

\begin{equation}\label{diag generating cups and caps}
\hspace{30 pt}
{\begin{tikzpicture}[
    scale=1,
    line width=0.8pt,
    arrow_style/.style={
        decoration={
            markings,
            mark=at position 0.3 with {\arrow{Stealth[black]}}, 
            mark=at position 0.8 with {\arrow{Stealth[black]}}  
        },
        postaction={decorate}
    }
]
\draw[color={rgb,255:red,0; green,0; blue,192}, arrow_style] (-0.8,0) arc (180:0:0.8);
\node at (0,0.3) {$\gr{Ii}$};
\node at (0,1.5) {$\gr{I}$};
\draw (-1,0) rectangle (1,1.9);
\end{tikzpicture}
}
\hspace{20 pt}
{\begin{tikzpicture}[
    scale=1,
    line width=0.8pt,
    arrow_style/.style={
        decoration={
            markings,
            mark=at position 0.3 with {\arrow{Stealth[black]}}, 
            mark=at position 0.8 with {\arrow{Stealth[black]}}  
        },
        postaction={decorate}
    }
]
\draw[color={rgb,255:red,0; green,0; blue,192}, arrow_style] (0.8,0) arc (0:180:0.8);
\node at (0,1.5) {$\gr{Ii}$};
\node at (0,0.4) {$\gr{I}$};

\draw (-1,0) rectangle (1,1.9);
\end{tikzpicture}}
\hspace{20 pt}
{\begin{tikzpicture}[
    scale=1,
    line width=0.8pt,
    arrow_style/.style={
        decoration={
            markings,
            mark=at position 0.3 with {\arrow{Stealth[black]}}, 
            mark=at position 0.8 with {\arrow{Stealth[black]}}  
        },
        postaction={decorate}
    }
]
\draw[color={rgb,255:red,0; green,0; blue,192}, arrow_style] (-0.8,1.9) arc (180:360:0.8);
\node at (0,0.3) {$\gr{Ii}$};
\node at (0,1.5) {$\gr{I}$};

\draw (-1,0) rectangle (1,1.9);
\end{tikzpicture}}
\hspace{20 pt}
{\begin{tikzpicture}[
    scale=1,
    line width=0.8pt,
    arrow_style/.style={
        decoration={
            markings,
            mark=at position 0.3 with {\arrow{Stealth[black]}}, 
            mark=at position 0.8 with {\arrow{Stealth[black]}}  
        },
        postaction={decorate}
    }
]
\draw[color={rgb,255:red,0; green,0; blue,192}, arrow_style] (0.8,1.9) arc (360:180:0.8);
\node at (0,1.5) {$\gr{Ii}$};
\node at (0,0.4) {$\gr{I}$};
\draw (-1,0) rectangle (1,1.9);
\end{tikzpicture}}
\hspace{20 pt}
{\begin{tikzpicture}[
    scale=1,
    line width=0.8pt,
]
\draw (-1,0) rectangle (1,1.9);

\node[draw, rectangle] (f) at (0,1) {$f$};
\node at (0, 0.3) {$\gr{I}$};
\end{tikzpicture}}
\end{equation}

\hfill
\begin{tikzcd}
    R^{I} \\
    R^{I}  \otimes_{R^{Ii}} R^{I}  \arrow[u, "m"]
\end{tikzcd}
\hfill
\begin{tikzcd}
    R^{Ii} \\
    R^{I} \arrow[u, "\dd"]
\end{tikzcd}
\hfill 
\begin{tikzcd}
    R^{I} \\
    R^{Ii}  \arrow[u, "\iota"]
\end{tikzcd}
\hfill
\begin{tikzcd}
    R^{I}  \otimes_{R^{Ii}} R^{I}  \\
    R^{I} \arrow[u, "\Delta"]
\end{tikzcd}
\hfill
\begin{tikzcd}[column sep=0em]
     R^I  &f\\
   R^I \arrow[u] &1 \arrow[u, mapsto]
\end{tikzcd}.
\hfill
\vspace{15pt}

\begin{equation}\label{diag generating crossings}
\hfill
{
\begin{tikzpicture}[
    scale=1,
    line width=0.8pt,
    arrow_style/.style={
        decoration={
            markings,
            mark=at position 0.3 with {\arrow{Stealth[black]}}, 
            mark=at position 0.8 with {\arrow{Stealth[black]}}  
        },
        postaction={decorate}
    }
]

\draw[color={rgb,255:red,0; green,0; blue,192}, arrow_style] (1,-1) -- (-1,1);
\draw[red, arrow_style] (-1,-1) -- (1,1);

\node at (0, 0.7) {$\gr{Ii}$};
\node at (0.8, 0) {$\gr{Iij}$};
\node at (0, -0.7) {$\gr{Ij}$};
\node at (-1, 0) {$\gr{I}$};
\node at (1.7,0) {};
\node at (-1.7,0) {};
\node[draw, inner sep=0pt, fit=(current bounding box)] {};

\comm{
    \node (bbox) [fit=(current bounding box), inner sep=0pt] {};

    \draw[black] (bbox.north west) -- (bbox.north east); 
    \draw[black] (bbox.south west) -- (bbox.south east); 
}
\end{tikzpicture}
}
\hspace{100 pt}
{
\begin{tikzpicture}[
    scale=1,
    line width=0.8pt,
    arrow_style/.style={
        decoration={
            markings,
            mark=at position 0.3 with {\arrow{Stealth[black]}}, 
            mark=at position 0.8 with {\arrow{Stealth[black]}}  
        },
        postaction={decorate}
    }
]

\draw[color={rgb,255:red,0; green,0; blue,192}, arrow_style] (-1,1) -- (1,-1);
\draw[red, arrow_style] (1,1) -- (-1,-1);

\node at (0, -0.7) {$\gr{Ii}$};
\node at (-1, 0) {$\gr{Iij}$};
\node at (0, 0.7) {$\gr{Ij}$};
\node at (1, 0) {$\gr{I}$};
\node at (1.7,0) {};
\node at (-1.7,0) {};

\node[draw, inner sep=0pt, fit=(current bounding box)] {};

\comm{
    \node (bbox) [fit=(current bounding box), inner sep=0pt] {};

    \draw[black] (bbox.north west) -- (bbox.north east); 
    \draw[black] (bbox.south west) -- (bbox.south east); 
}
\end{tikzpicture}
}
\hfill \null
\end{equation}

\begin{tikzcd}[column sep=0em]
    R^{I}  \otimes_{R^{Ii}} R^{Ii} \otimes_{R^{Iij}} R^{Iij}  &f\otimes 1 \otimes 1\\
    R^{I}  \otimes_{R^{Ij}} R^{Ij} \otimes_{R^{Iij}} R^{Iij}  \arrow[u] &f \otimes 1 \otimes 1 \arrow[u, mapsto]
\end{tikzcd}
\hfill
\begin{tikzcd}[column sep=0em]
     R^{Iij} \otimes_{R^{Iij}} R^{Ij}  \otimes_{R^{Ij}}R^{I} &1\otimes 1 \otimes f\\
   R^{Iij} \otimes_{R^{Iij}} R^{Ii}  \otimes_{R^{Ii}}R^{I} \arrow[u] &1 \otimes 1 \otimes f \arrow[u, mapsto]
\end{tikzcd}. \vspace{10 pt}

Note that while the upward and downward crossings are degree zero morphisms, the cups, caps and polynomials need not be. The clockwise cup and cap in the diagrams above have positive degree, given by $\ell(Ii)-\ell(I)$, and the counterclockwise cup and cap have negative degree given by $\ell(I)-\ell(Ii)$. The degree of a polynomial $f$ is determined by its degree as an element of $R$.

The first set of relations in $\mathcal{D}(\Gamma)$ are the isotopy relations (for valid region labelings): 

 \begin{equation} \label{isotopy} \vcenter{\xy (0,0)*{\def\svgscale{1}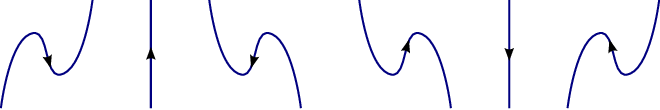} \endxy}, \end{equation} 
\begin{equation} \label{crosscyclic} \vcenter{\xy (0,0)*{\def\svgscale{1}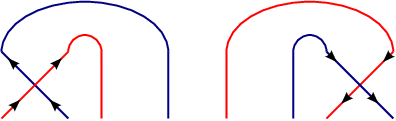} \endxy}, \end{equation}
\begin{equation}\label{crosscyclic2} \vcenter{\xy (0,0)*{\def\svgscale{1}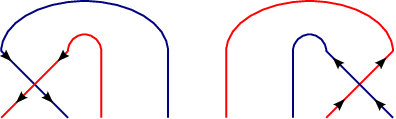} \endxy}. \end{equation}

We also want polynomials $f$ to behave as elements of $R$, so we have the following relation in $\DiagSBS(\Gamma)$ (for valid region labels).
\begin{equation} \label{polymult} \vcenter{\xy (0,0)*{\def\svgscale{1}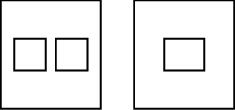} \endxy} \end{equation}

\begin{defn}\label{defn diagrammatic precategory}
    We define $\DiagSBS^{\pre}$ to be the diagrammatic category with the same generators as that of $\DiagSBS(\Gamma)$, with the isotopy relations in \cref{isotopy}, \cref{crosscyclic} and \cref{crosscyclic2} imposed, along with \cref{polymult}.
\end{defn}

The following are examples of 2-morphisms in $\DiagSBS^{\pre}$ (as well as in $\DiagSBS(\Gamma)$ and $\DiagSBS$, once defined). As in \cref{sec:combo}, we use $\hh{a}$ to denote the maximal finitary subset $S \setminus \set{a}$.
A label $k$ on a strand corresponds to the simple reflection $s_k$.

\[\begin{array}{cccc}
  \vcenter{\xy (0,0)*{\def\svgscale{0.15}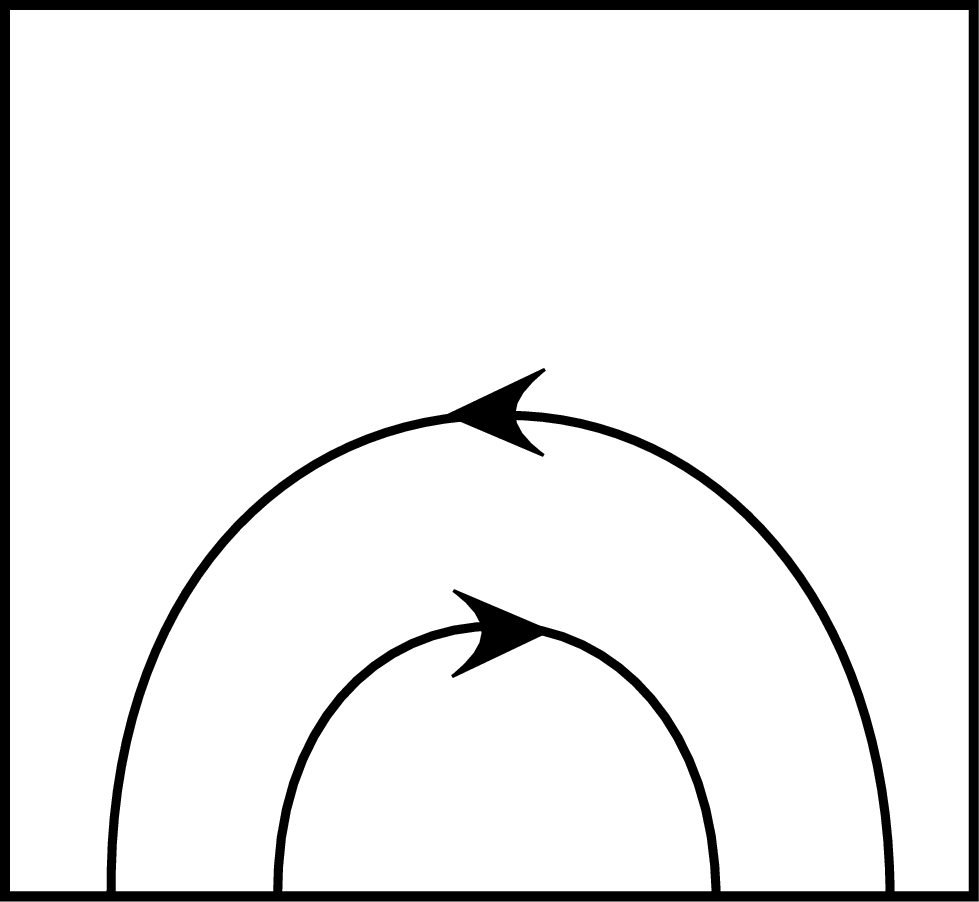} \endxy}    &\vcenter{\xy (0,0)*{\def\svgscale{0.15}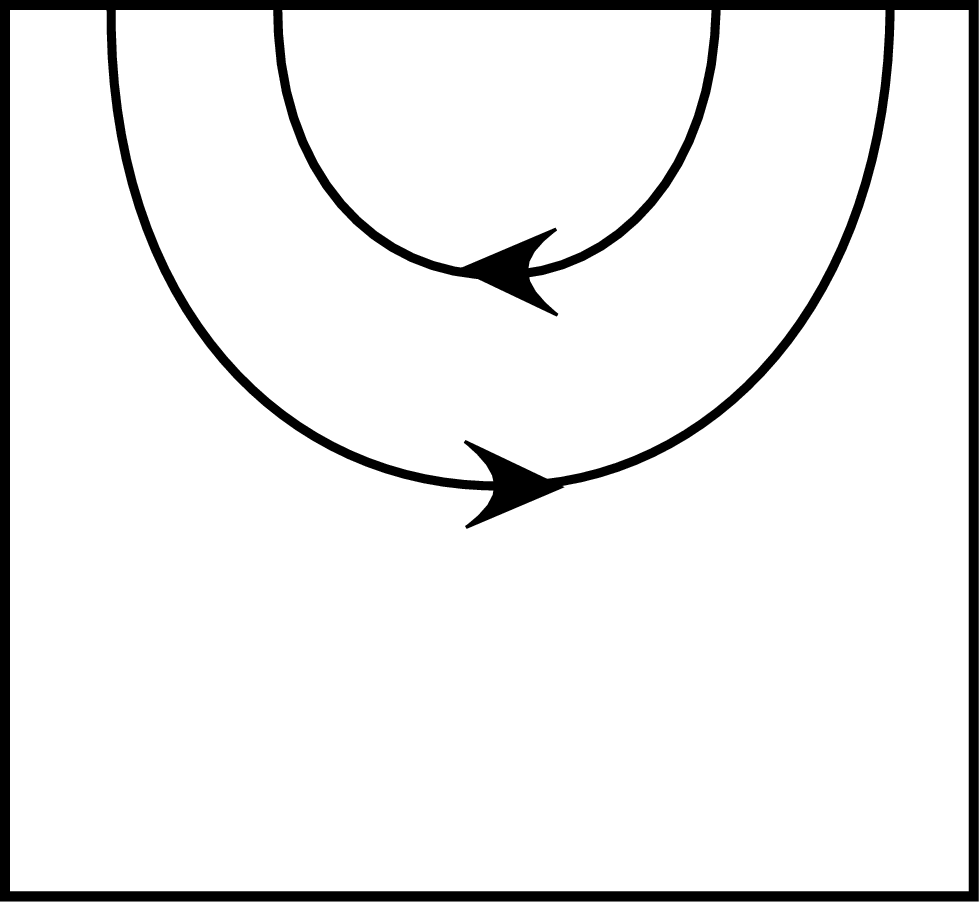} \endxy} &\vcenter{\xy (0,0)*{\def\svgscale{0.15}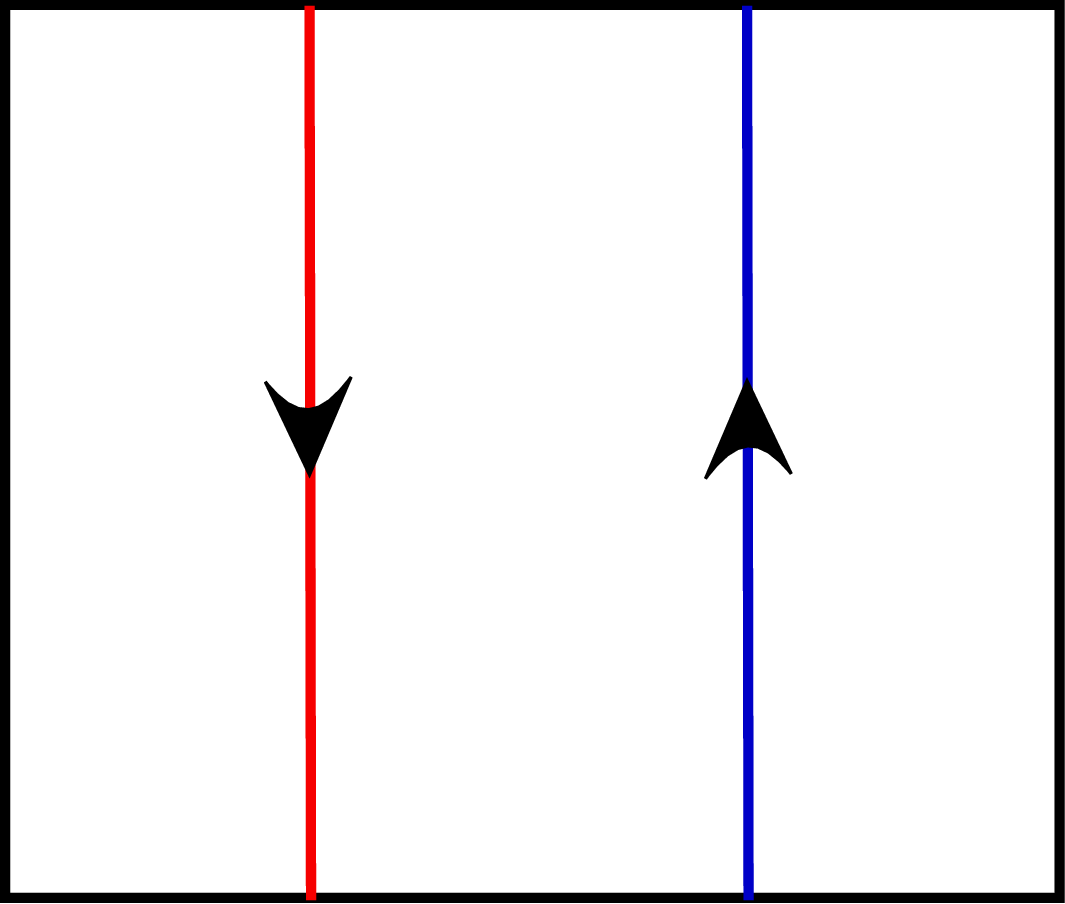} \endxy} 
  &\vcenter{\xy (0,0)*{\def\svgscale{0.15}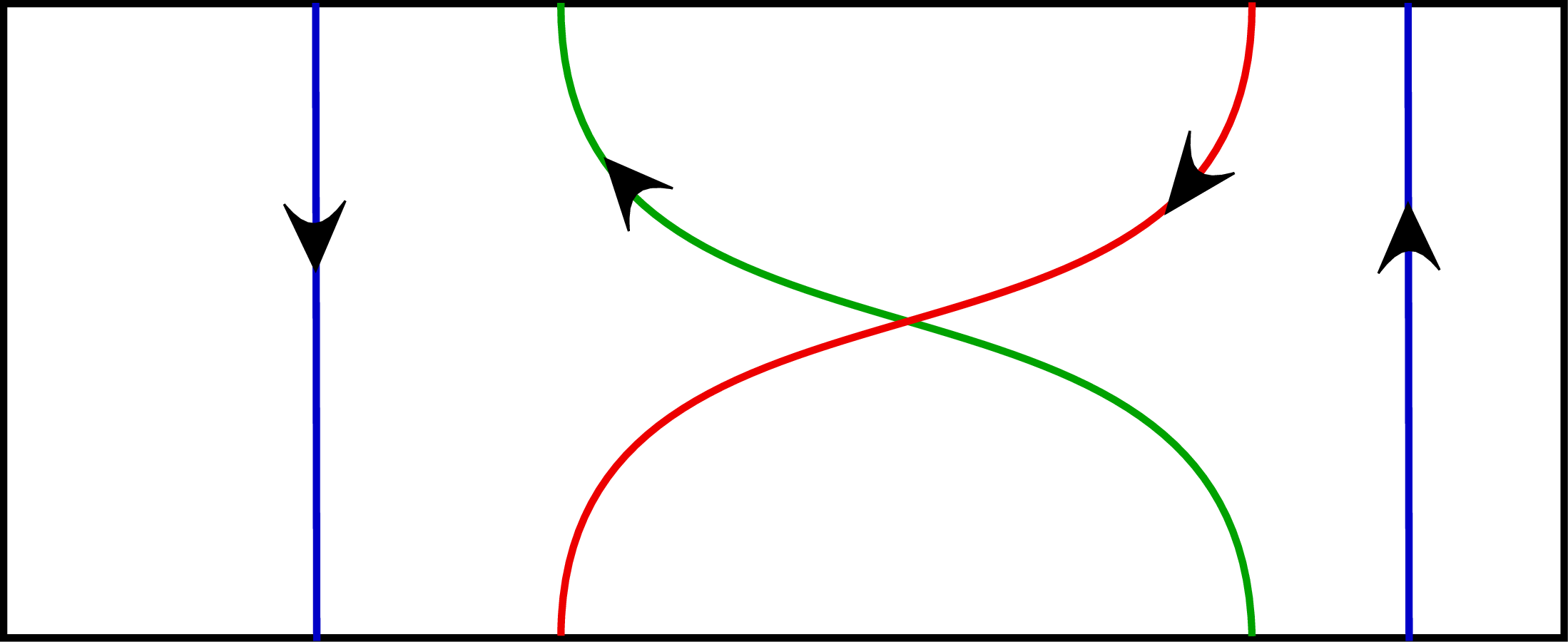} \endxy}
\end{array} 
\]
Note that we sometimes omit the strand and region labels, when these are already determined by the other labels. 
We will also simplify diagrams by drawing them up to isotopy. 
We have sideways crossings obtained by rotating the upward or downward crossings (which is well defined after imposing the isotopy relations) as shown below, for any $Iij \in \mc{F}$. 
\begin{equation}
    {
\begin{tikzpicture}[baseline=-0.5ex,
    scale=1,
    line width=0.8pt,
    arrow_style/.style={
        decoration={
            markings,
            mark=at position 0.3 with {\arrow{Stealth[black]}}, 
            mark=at position 0.8 with {\arrow{Stealth[black]}}  
        },
        postaction={decorate}
    }
]

\draw[color={rgb,255:red,0; green,0; blue,192}, arrow_style] (-1,1) -- (1,-1);
\draw[red, arrow_style] (-1,-1) -- (1,1);

\node at (0, 0.7) {$\gr{I}$};
\node at (0.8, 0) {$\gr{Ii}$};
\node at (0, -0.7) {$\gr{Iij}$};
\node at (-1, 0) {$\gr{Ij}$};
\node at (1.3,0) {};
\node at (-1.3,0) {};
\node[draw, inner sep=0pt, fit=(current bounding box)] {};

\comm{
    \node (bbox) [fit=(current bounding box), inner sep=0pt] {};

    \draw[black] (bbox.north west) -- (bbox.north east); 
    \draw[black] (bbox.south west) -- (bbox.south east); 
}
\end{tikzpicture}
} := 
    {\begin{tikzpicture}[baseline=-0.5ex,
    scale=1,
    line width=0.8pt,
    arrow_style/.style={
        decoration={
            markings,
            mark=at position 0.3 with {\arrow{Stealth[black]}}, 
            mark=at position 0.8 with {\arrow{Stealth[black]}}  
        },
        postaction={decorate}
    }
]
\draw[red, arrow_style] (.5,-1) to (-.5,1);
\draw[color={rgb,255:red,0; green,0; blue,192}, arrow_style] (-1.1,1) to[out=-90,in=-135,looseness=2] (0,0) to [out=45,in=90,looseness=2] (1.1,-1);
\node at (0.7,0.7) {$\gr{Ii}$};
\node at (-.7,-.7) {$\gr{Ij}$};
\node at (.5,-.2) {$\gr{Iij}$};
\node at (-.5,.2) {$\gr{I}$};
\node at (1.3,0) {};
\node at (-1.3,0) {};
\node[draw, inner sep=0pt, fit=(current bounding box)] {};
\end{tikzpicture}} 
    =
    {\begin{tikzpicture}[baseline=-0.5ex, 
    scale=1,
    line width=0.8pt,
    arrow_style/.style={
        decoration={
            markings,
            mark=at position 0.3 with {\arrow{Stealth[black]}}, 
            mark=at position 0.8 with {\arrow{Stealth[black]}}  
        },
        postaction={decorate}
    }
]
\draw[color={rgb,255:red,0; green,0; blue,192}, arrow_style] (.4,1) to (-.4,-1);
\draw[red, arrow_style] (-1.1,-1) to[out=90,in=135,looseness=2] (0,0) to [out=-45,in=-90,looseness=2] (1.1,1);
\node at (0.7,-.7) {$\gr{Ii}$};
\node at (-.7,.7) {$\gr{Ij}$};
\node at (.5,.3) {$\gr{I}$};
\node at (-.5,-.3) {$\gr{Iij}$};
\node at (1.3,0) {};
\node at (-1.3,0) {};
\node[draw, inner sep=0pt, fit=(current bounding box)] {};
\end{tikzpicture}}.
\end{equation}
\begin{equation}
    {
\begin{tikzpicture}[baseline=-0.5ex,
    scale=1,
    line width=0.8pt,
    arrow_style/.style={
        decoration={
            markings,
            mark=at position 0.3 with {\arrow{Stealth[black]}}, 
            mark=at position 0.8 with {\arrow{Stealth[black]}}  
        },
        postaction={decorate}
    }
]

\draw[color={rgb,255:red,0; green,0; blue,192}, arrow_style] (1,-1) -- (-1,1);
\draw[red, arrow_style] (1,1) -- (-1,-1);

\node at (0, 0.7) {$\gr{Iij}$};
\node at (0.8, 0) {$\gr{Ij}$};
\node at (0, -0.7) {$\gr{I}$};
\node at (-1, 0) {$\gr{Ii}$};
\node at (1.3,0) {};
\node at (-1.3,0) {};
\node[draw, inner sep=0pt, fit=(current bounding box)] {};

\comm{
    \node (bbox) [fit=(current bounding box), inner sep=0pt] {};

    \draw[black] (bbox.north west) -- (bbox.north east); 
    \draw[black] (bbox.south west) -- (bbox.south east); 
}
\end{tikzpicture}
} := 
    {\begin{tikzpicture}[baseline=-0.5ex,
    scale=1,
    line width=0.8pt,
    arrow_style/.style={
        decoration={
            markings,
            mark=at position 0.3 with {\arrow{Stealth[black]}}, 
            mark=at position 0.8 with {\arrow{Stealth[black]}}  
        },
        postaction={decorate}
    }
]
\draw[red, arrow_style] (-.5,1) to (.5,-1);
\draw[color={rgb,255:red,0; green,0; blue,192}, arrow_style] (1.1,-1) to[out=90,in=45,looseness=2] (0,0) to [out=-135,in=-90,looseness=2] (-1.1,1);
\node at (0.7,0.7) {$\gr{Ij}$};
\node at (-.7,-.7) {$\gr{Ii}$};
\node at (.5,-.2) {$\gr{I}$};
\node at (-.5,.2) {$\gr{Iij}$};
\node at (1.3,0) {};
\node at (-1.3,0) {};
\node[draw, inner sep=0pt, fit=(current bounding box)] {};
\end{tikzpicture}} 
    =
    {\begin{tikzpicture}[baseline=-0.5ex, 
    scale=1,
    line width=0.8pt,
    arrow_style/.style={
        decoration={
            markings,
            mark=at position 0.3 with {\arrow{Stealth[black]}}, 
            mark=at position 0.8 with {\arrow{Stealth[black]}}  
        },
        postaction={decorate}
    }
]
\draw[color={rgb,255:red,0; green,0; blue,192}, arrow_style] (-.5,-1) to (.5,1);
\draw[red, arrow_style] (1.1,1) to[out=-90,in=-45,looseness=2] (0,0) to [out=135,in=90,looseness=2] (-1.1,-1);
\node at (0.7,-.7) {$\gr{Ij}$};
\node at (-.7,.7) {$\gr{Ii}$};
\node at (.5,.3) {$\gr{Iij}$};
\node at (-.5,-.3) {$\gr{I}$};
\node at (1.3,0) {};
\node at (-1.3,0) {};
\node[draw, inner sep=0pt, fit=(current bounding box)] {};
\end{tikzpicture}}.
\end{equation}

Unlike the upward and downward crossings, these are not always of degree 0. The degree of the sidewards crossings given above are $\ell(Iij)+\ell(I)-\ell(Ii)-\ell(Ij)$. A mnemonic for the degree of the sideways crossings is ``big + small - middle - middle". We will henceforth use sidewards crossings in our diagrams and draw diagrams up to isotopy without explicit mention.


For example, the diagram in the RHS of the equation below will appear multiple times in later sections.
\begin{equation}\label{eq trivalent split}\vcenter{\xy (0,0)*{\def\svgscale{0.15}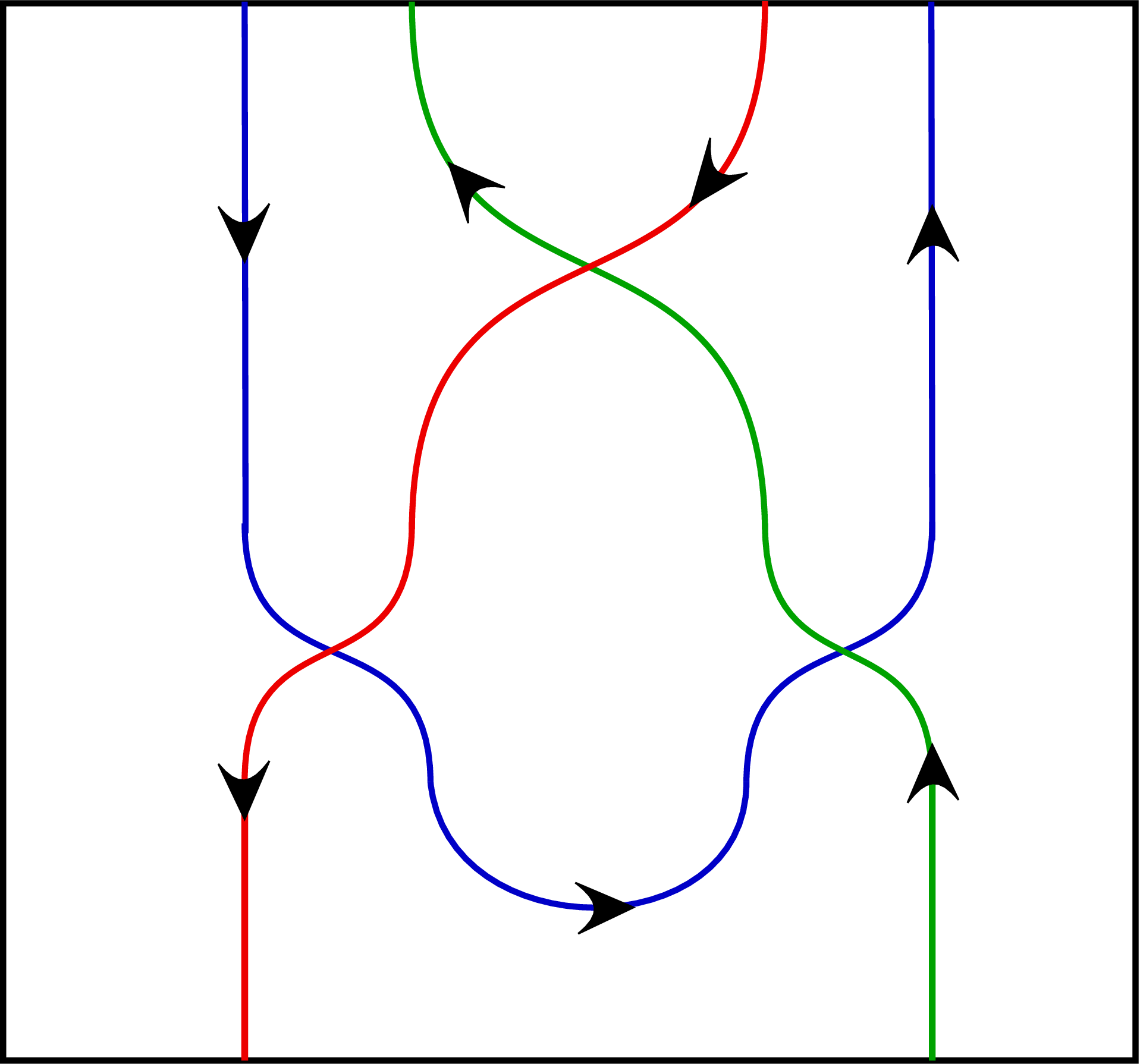} \endxy}= \vcenter{\xy (0,0)*{\def\svgscale{0.2}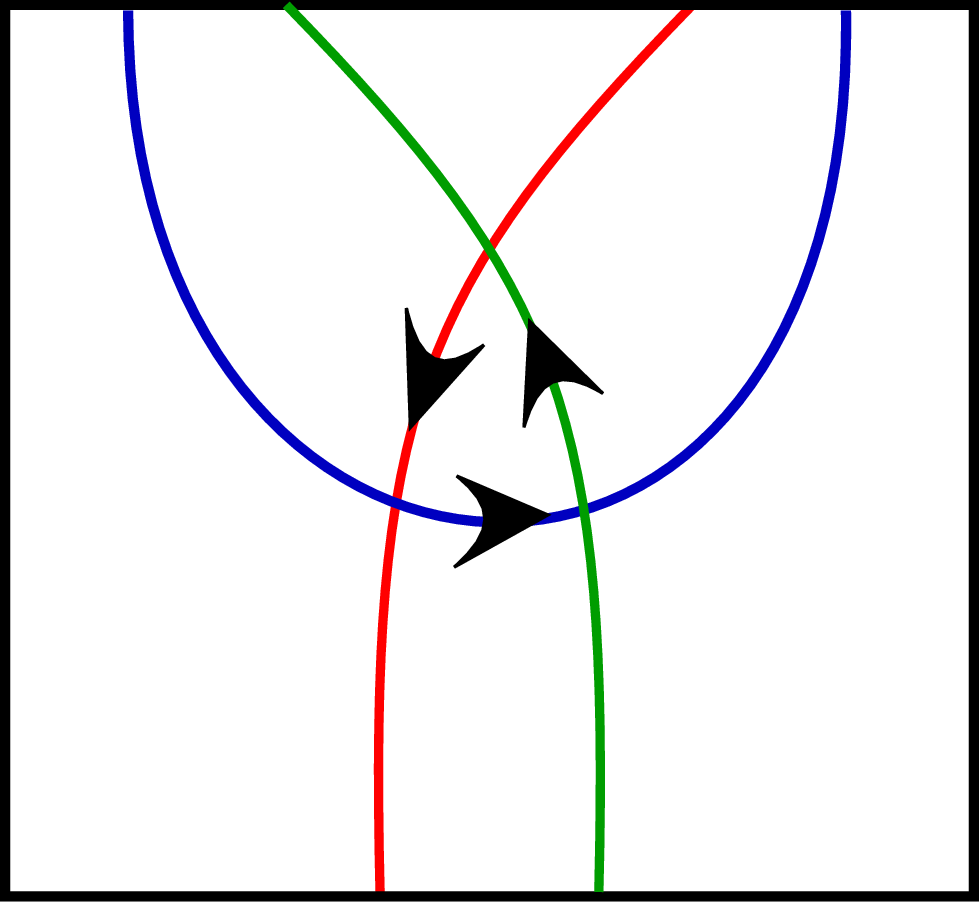} \endxy}\end{equation}
\begin{lemma}\label{lemma image of trivalent has degree 0}
    The 2-morphism in \cref{eq trivalent split} is of degree 0.
\end{lemma}
\begin{proof}
The rotation operator $\sigma$ (\cref{defn-rotation operator}) preserves lengths of parabolic subgroups, hence without loss of generality, we may assume $a=0$. 
The sidewards crossing in \cref{eq trivalent split} then has degree $\ell(\hh{l})+\ell(\hh{0,l,k+l})-\ell(\hh{0, l})-\ell(\hh{l, k+l})$ and the counter-clockwise cup has degree $\ell(\hh{0, l, k+l})-\ell(\hh{0, k+l})$. Note that the parabolic subgroups $W_{\hh{0,l,k+l}}, W_{\hh{0, l}}$ etc. are just a product of symmetric groups $S_l\times S_k \times S_{n-k-l}, S_l \times S_{n-l}$ etc., so that the net degree is 
\begin{align*}
    &\ell(\hh{l})+2\ell(\hh{0,l,k+l})-\ell(\hh{0, l})-\ell(\hh{l, k+l})-\ell(\hh{0, k+l})\\
    &= \binom{n}{2} + 2\left(\binom{l}{2}+\binom{k}{2}+\binom{n-k-l}{2}\right) - \left(\binom{l}{2}+\binom{n-l}{2}\right)\\
    &\qquad-\left(\binom{k}{2}+\binom{n-k}{2}\right)-\left(\binom{n-k-l}{2}+\binom{k+l}{2}\right)
    \\
    &=
    \binom{n}{2}+\binom{l}{2}+\binom{k}{2}+\binom{n-k-l}{2}-\binom{n-l}{2}-\binom{n-k}{2}-\binom{k+l}{2}=0. \qedhere
\end{align*}
\end{proof}

We now continue discussing the relations in $\DiagSBS(\Gamma)$.
All future relations hold when rotated or with the colors switched, but not with the orientations reversed. 

Consider $I, Ii \in \mc{F}$. 
We have the following relations associated to the Frobenius extension $(R^{Ii}\hookrightarrow R^I, \dd=\dd^I_{Ii})$.

\begin{equation} \label{polyslide} \vcenter{\xy (0,0)*{\def\svgscale{1}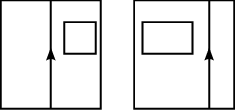} \endxy} \end{equation}
\begin{equation} \label{cccirc} \vcenter{\xy (0,0)*{\def\svgscale{1}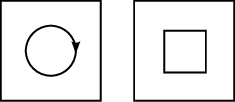} \endxy} \end{equation}
\begin{equation} \label{ccwcirc} \vcenter{\xy (0,0)*{\def\svgscale{1}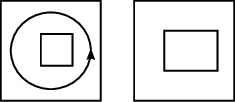} \endxy} \end{equation}
\begin{equation} \label{Bsplitting} \vcenter{\xy (0,0)*{\def\svgscale{1}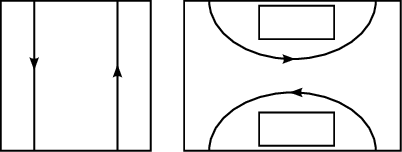} \endxy} \end{equation}
In \cref{Bsplitting}, we use the Sweedler notation from \cref{subsubsec-notation in further calculations.}

Now consider $I, Ii, Iij \in \mc{F}$. 
We then have the following \emph{Reidemeister II relations}.

\begin{equation} \label{R2oriented} \vcenter{\xy (0,0)*{\def\svgscale{0.8}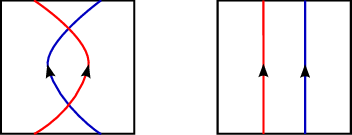} \endxy} \end{equation}
\begin{equation} \label{R2unoriented1} \vcenter{\xy (0,0)*{\def\svgscale{0.8}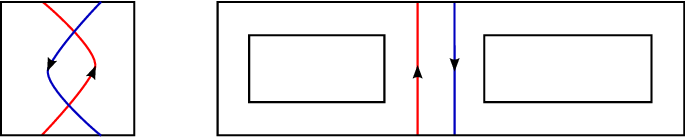} \endxy} \end{equation}
\begin{equation} \label{R2unoriented2} \vcenter{\xy (0,0)*{\def\svgscale{0.8}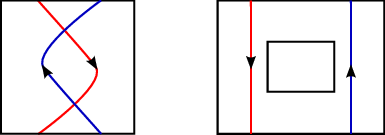} \endxy} .\end{equation}
The polynomial in the box in \cref{R2unoriented2} is  $\mu_{Iij}^{Ii, Ij}$ (see \cref{subsubsec-notation in further calculations.}).
The RHS of \cref{R2unoriented1} may be written in several equivalent ways using the following identity specialized at $f=1$:
\begin{align}\label{eq R2 unoriented equivalent formulas}
\notag \Delta^{Ii}_{Iij, (1)} \otimes  \dd^I_{Ij}(f\Delta^{Ii}_{Iij, (2)})= \dd^I_{Ii}(f\Delta^I_{Iij,(1)})\otimes \dd^I_{Ij}(\Delta^I_{Iij,(2)})&=\dd^I_{Ii}(\Delta^I_{Iij,(1)})\otimes \dd^I_{Ij}(f\Delta^I_{Iij,(2)})\\
= \dd^I_{Ii}(f\Delta^{Ij}_{Iij,(1)})\otimes \Delta^{Ij}_{Iij, (2)} & \qquad\text{for $f \in R^I$}.
\end{align}
The equation \eqref{eq R2 unoriented equivalent formulas} is a consequence of \cite[(2.7), Section 2.2]{EWSFrob}.

We mention a special instance where \cref{R2unoriented1}, \cref{R2unoriented2} become easier.
When $i$ and $j$ are in different connected components of $Iij \in \mc{F}$, we have that $\Phi_{Iij}^+ \setminus \Phi_{Ii}^+=\Phi_{Ij}^+ \setminus \Phi_I^+$, so that $\mu_{Iij}^{Ii, Ij}=\frac{\mu_{Iij}^{Ii}}{\mu_{Ij}^I}=1.$
We also have that $\dd^{I}_{I\sqcup \set{j}}=\dd^{I\sqcup \set{i}}_{I\sqcup \set{i,j}}$ in this case (for example, by using \cref{corollary-normalization}), so that all but one term on the RHS of \cref{R2unoriented1} vanish.
Hence both the LHS and RHS of the relations \cref{R2unoriented1}\;\;, \cref{R2unoriented2} are degree zero maps in this case, and we have the following \emph{distant Reidemeister II} relations (shortened as \emph{distant RII}).

\begin{align}\label{R2 distant}
 \vcenter{\xy (0,0)*{\def\svgscale{0.15}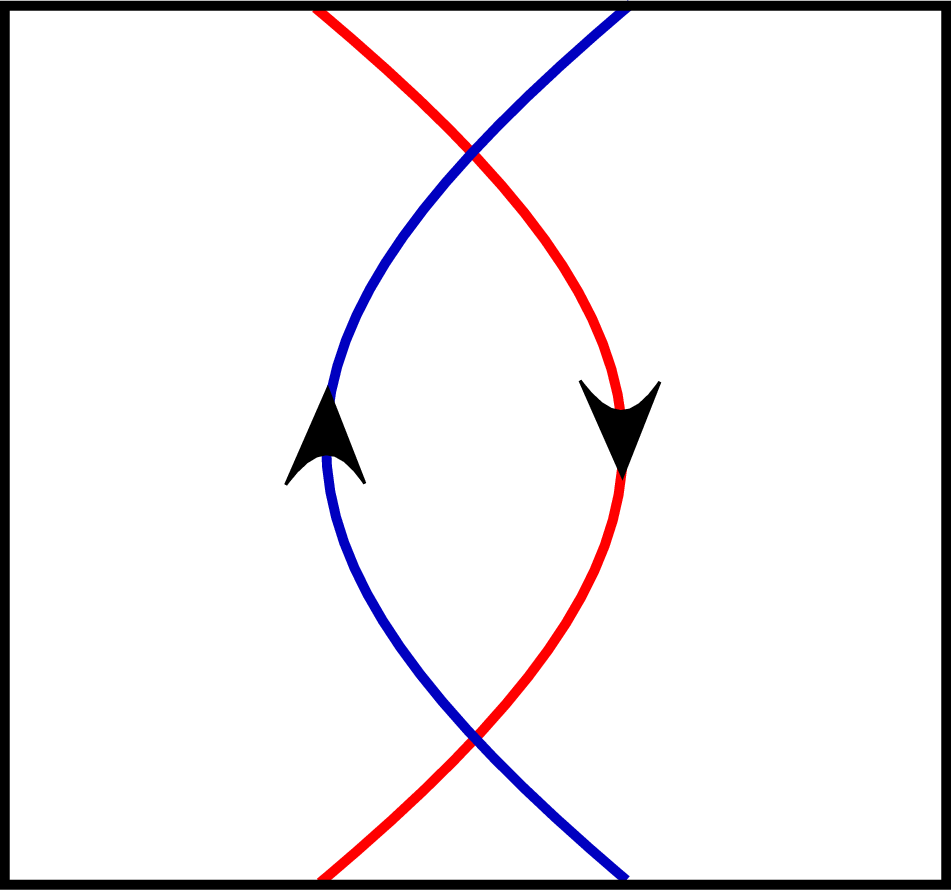} \endxy}= \vcenter{\xy (0,0)*{\def\svgscale{0.15}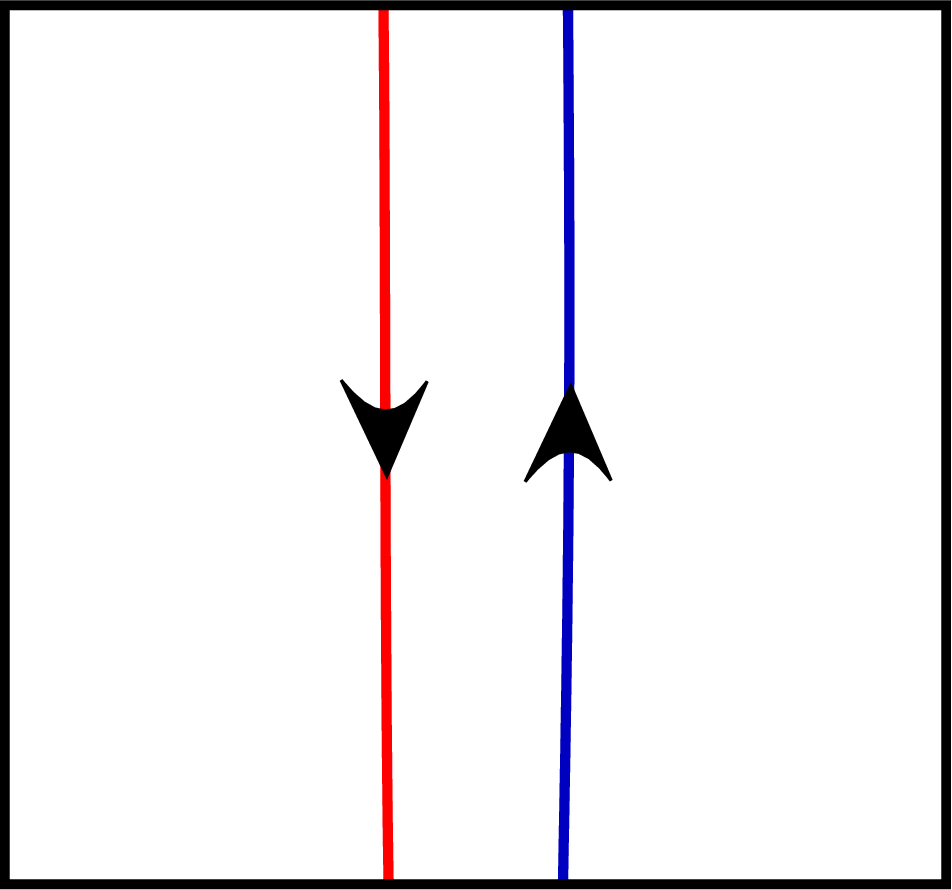} \endxy}, \qquad\vcenter{\xy (0,0)*{\def\svgscale{0.15}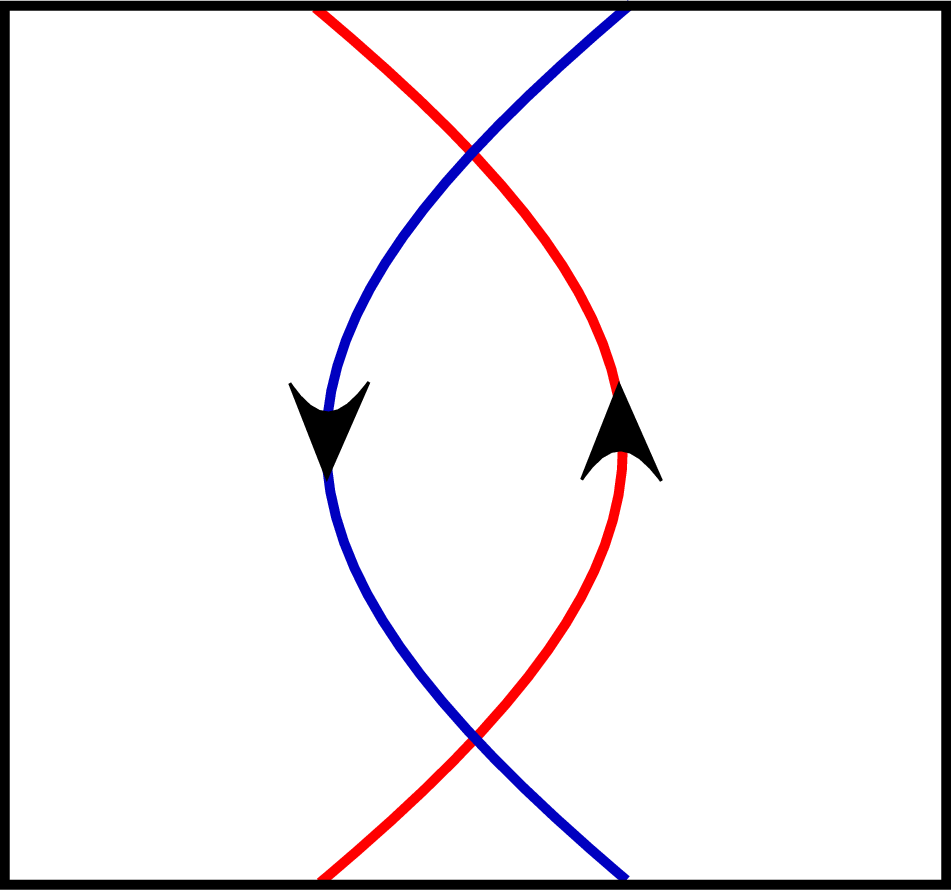} \endxy}= \vcenter{\xy (0,0)*{\def\svgscale{0.15}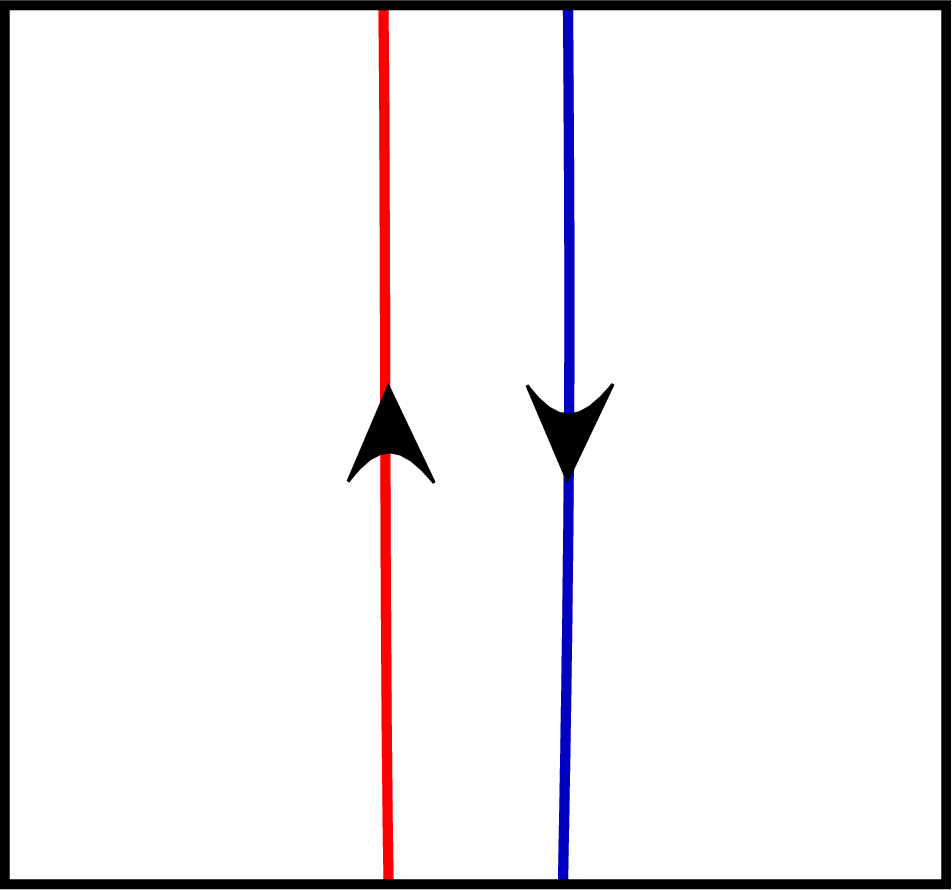} \endxy}.   
\end{align}

Now consider $I \in F$ and $i,j,k \in \Gamma$ such that $Iijk \in \mc{F}$. 
We have the following \emph{Reidemeister III relation}.

\begin{equation} \label{R3oriented} \vcenter{\xy (0,0)*{\def\svgscale{0.9}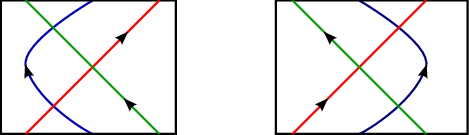} \endxy} .\end{equation}



We also have the \emph{hard Reidemeister III relation}:
\begin{equation} \label{R3unoriented}
  \vcenter{\xy (0,0)*{\def\svgscale{0.9}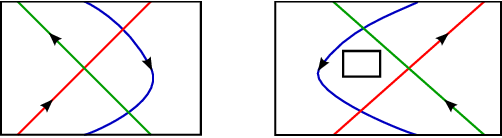} \endxy}. \end{equation} The polynomial in the
box is $\frac{\mu_{Iijk} \mu_{Ii} \mu_{Ij}\mu_{Ik}}{\mu_{Iij} \mu_{Iik} 
\mu_{Ijk} \mu_I}$, which is a genuine polynomial by assumption R3 (\cref{defn-condition R3}).

We mention a special case of \cref{R3unoriented}, which will be used in \cref{subsec- I=H}. 
When $i,j,k$ are not in the same component of $Iijk$, we have the following relation which we refer to as the \emph{distant Reidemeister III} relation, usually shortened to \emph{distant RIII}.
\begin{align} \label{distantR3}
\vcenter{\xy (0,0)*{\def\svgscale{0.15}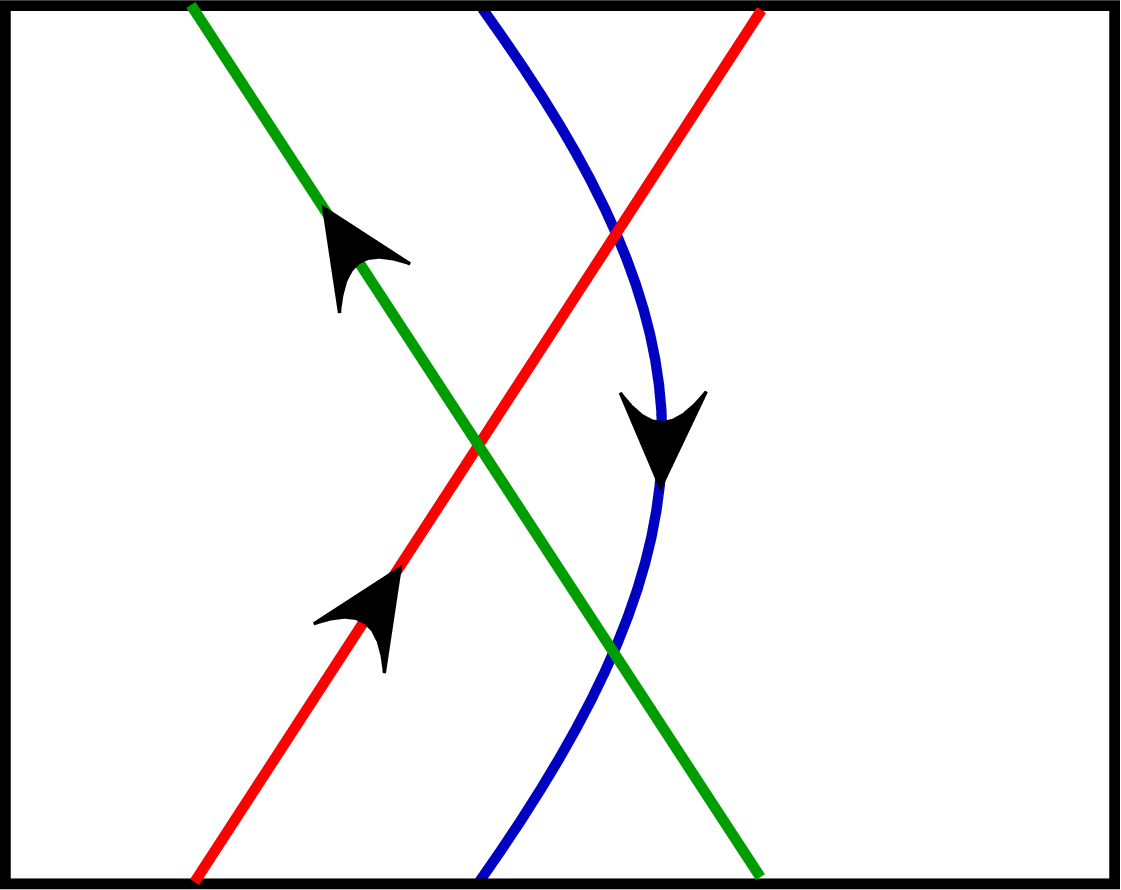} \endxy}=\vcenter{\xy (0,0)*{\def\svgscale{0.15}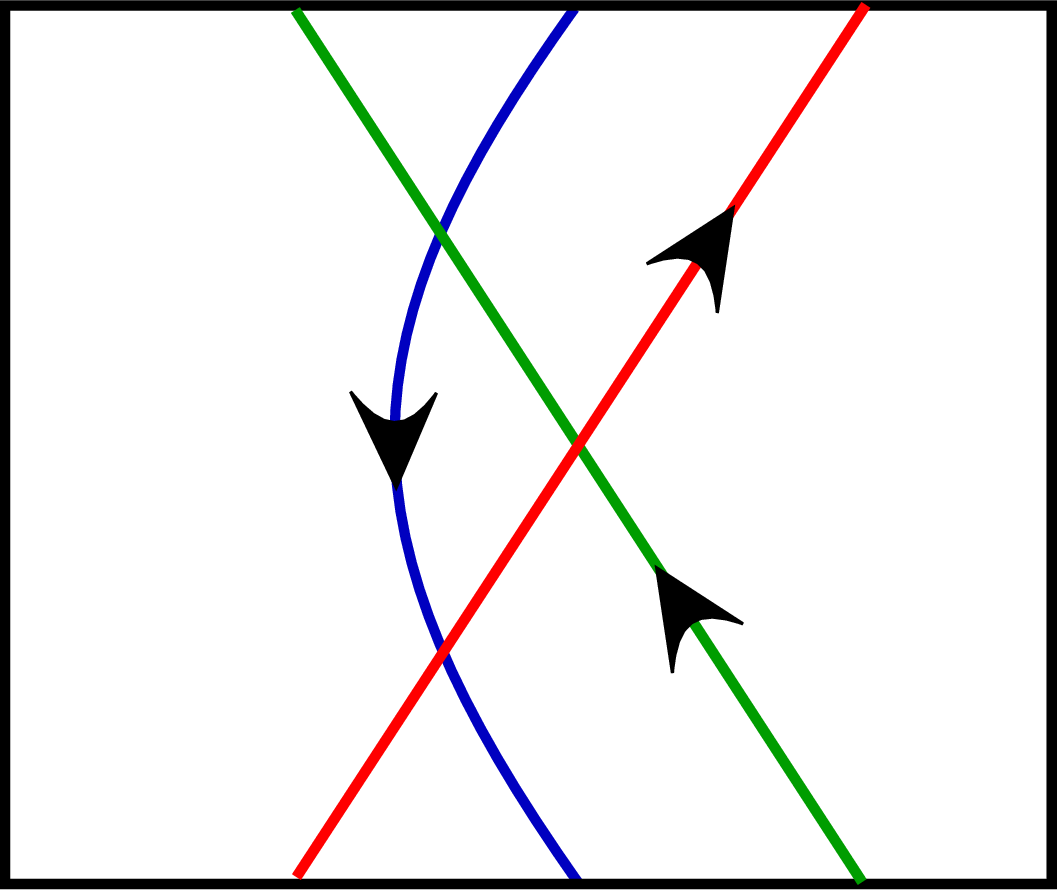} \endxy}.
\end{align}
We leave it as an exercise to the reader to check that the polynomial in the box in the RHS of \cref{R3unoriented} is 1 in this case, so that \cref{distantR3} is indeed just a special case of \cref{R3unoriented}.




The 2-functor $\Phi: \DiagSBS(\Gamma)$ has a kernel, and we define $\dSBSBim$ to be the 2-category obtained by killing the kernel, so that $\Phi$ factors through a faithful functor $\DiagSBS \rightarrow \SBSBim$, which we again denote by $\Phi$.

\PTv2{Still retaining this comment in case we want to talk a bit more about the kernel.}

Recall that $\sigma$ and $\tau$ denoted automorphisms of the affine Dynkin diagram, which extended to automorphisms of $\Lambda_\zeta$ and $R$. 
We use these automorphisms to define certain symmetries of $\DiagSBS^{\pre}$ and $\SBSBim$, which descend also to $\DiagSBS(\Gamma)$ and $\DiagSBS$.

\begin{defn} \label{defn:symmetries of diagsbs} 




We define 2-endofunctors $\sigma, \rotation, \Upsilon_S$ on $\DiagSBS^{\pre}$ and $\SBSBim$ as follows. 

On objects/region labels $\sigma$ sends $I \mapsto \sigma(I)$, e.g. it sends $\hh{a} \mapsto \hh{a+1}$. On $\DiagSBS$, this necessarily changes the labeling on edges by applying $\sigma$ to the set of simple reflections, and changes floating polynomials by applying $\sigma$ on polynomials. 
Aside from relabeling and changing polynomials, diagrams remain unchanged. 
On $\SBSBim$, the action of $\sigma$ on $\Hom(I,J)$ is obtained from the equivalence of graded $(R^I,R^J)$-bimodule categories induced by the isomorphism of graded algebras $R^I\otimes R^J\xrightarrow[\sigma \otimes \sigma]{\sim}R^{\sigma(I)}\otimes R^{\sigma(J)}$.
This is a monoidal, covariant equivalence.

On objects $\rotation$ sends $I \mapsto I$. On $\DiagSBS^{\pre}$, it rotates all diagrams by 180 degrees and leaves floating polynomials unchanged. Note that this 2-functor coincides with combining a horizontal flip and a vertical flip in either order.
On $\SBSBim$, $\rotation$ is obtained from adjunctions, but may also be thought of as a composition of natural symmetries corresponding to horizontal and vertical flips; see \cref{subsec reversal duality and phi} (in particular \cref{def rotation by 180 degrees}, \cref{corollary phi intertwines rotation} and \cref{rmk rotation and flips and adjunctions}) for the details. 
This is an antimonoidal, contravariant equivalence.

On objects $\Upsilon_S$ sends $I \mapsto \tau(I)$, e.g. it sends $\hh{a} \mapsto \hh{-a}$. This necessarily also acts on the edge labeling and on floating polynomials by $\tau$. It flips each diagram vertically (i.e. reflecting across a horizontal axis), reverses arrows to be consistent with region labels, and acts on any floating polynomials by $\tau$. 
On $\SBSBim$, $\Upsilon_S$ acts on 1-morphism categories via 
$$\Hom(I,J)^{op} \xrightarrow[D]{} \Hom(I,J) \xrightarrow[\tau_{\#}]{} Hom(\tau(I),\tau(J)),$$ where
$$D: \Hom(I,J)^{\op}\xrightarrow{\sim} \Hom(I,J): B \mapsto \Hom^\bullet_{R^J}(B, R^J)(2\ell(I)-2\ell(J)) $$ 
denotes the duality functor (see \cref{def duality symmetry D algebraic}),
and 
$$\tau_{\#}:\Hom(I,J) \xrightarrow{\sim} Hom(\tau(I),\tau(J))$$
is induced by the isomorphism of graded algebras $R^I\otimes R^J \xrightarrow[\tau \otimes \tau]{\sim}R^{\tau(I)}\otimes R^{\tau(J)}.$ 
This is a monoidal, contravariant equivalence. 

The 2-functor $\Phi: \DiagSBS^{\pre} \to \SBSBim$ intertwines the 2-equivalences $\sigma, \rotation, \Upsilon_S$, so that these 2-equivalences descend to $\DiagSBS(\Gamma)$ and $\DiagSBS$. 
Note here that ``intertwine" is not just a condition, but the data of certain natural isomorphisms (cf. \cref{thm phi commutes with dual and rev}, \cref{corollary phi intertwines rotation}).
We denote these 2-functors on $\mc{D}(\Gamma), \DiagSBS$ again by $\sigma, \rotation, \Upsilon_S$.
\end{defn}

\begin{remark} Unlike the category of webs, it would be absurd in $\DiagSBS$ to reverse the orientations of the strands without flipping, or vice versa. \end{remark}

\subsection{The functor \texorpdfstring{$\GS_{\zeta}$}{GS\_zeta}}\label{subsec the functor}

Define $\cwebs$ and $\dmSBSBim$ over the base ring $A=A_{\Z}\otimes_{\Z} \Q$,  
where $A_{\Z}$ is as in \cref{eq coefficient ring with q and zeta}, i.e.,
$A_{\Z}= \begin{cases}
\Z[q^{\pm 1}, \zeta^{\pm 1}]/(q^{-2}-\zeta^n) &\text{if $n$ odd;}\\
\Z[q^{\pm 1}, \zeta^{\pm 1}]/(q^{-1}-\zeta^{n/2}) &\text{if $n$ even.}
    \end{cases}$
    
\begin{defn}\label{defn dQGS functor}
Let $(\lambda, \nu)$ denote a collection of scalars $\lambda_k, \nu_k \in A^\times$ for various $1 \le k < n$, and scalars $\lambda_{k,l}, \nu_{k,l} \in A^{\times}$ for various $1 \le k,l < n$ with $k+l < n$. Attempt to define an $A$-linear $2$-functor $\GS_{\zeta} = \GS_{\zeta}(\lambda, \nu): \cwebs \rightarrow \dmSBSBim$ as follows. 

On objects, we have the bijection 
\begin{equation}\label{eq dQGS on objects}a \mapsto \hh{a}\end{equation}
for any 
$1\leq a \leq n$, where $\hh{a}$ is the maximal finitary subset $S\setminus \set{a} \subset S$.

On generating 1-morphisms, we have
\begin{align}\label{eq dQGS on 1-morphisms}
\notag \gr{a+k}\uparrow^{\blk{k}} \gr{a} \mapsto \dgr{a+k} \downarrow \dgr{a, a+k}\uparrow \dgr{a}    \\
\gr{a-k}\downarrow^{\blk{k}} \gr{a} \mapsto \dgr{a-k} \downarrow \dgr{a, a-k}\uparrow \dgr{a}    
\end{align}

On generating two morphisms, we have:
\begin{align}\label{diag-dQGS on 2-morphisms}
\notag&\vcenter{\xy (0,0)*{\def\svgscale{0.15}\input{arxiv-figures/cap_webs_c_svg-tex.eps_tex}} \endxy} \mapstosize{2} \vcenter{\xy (0,0)*{\def\svgscale{0.15}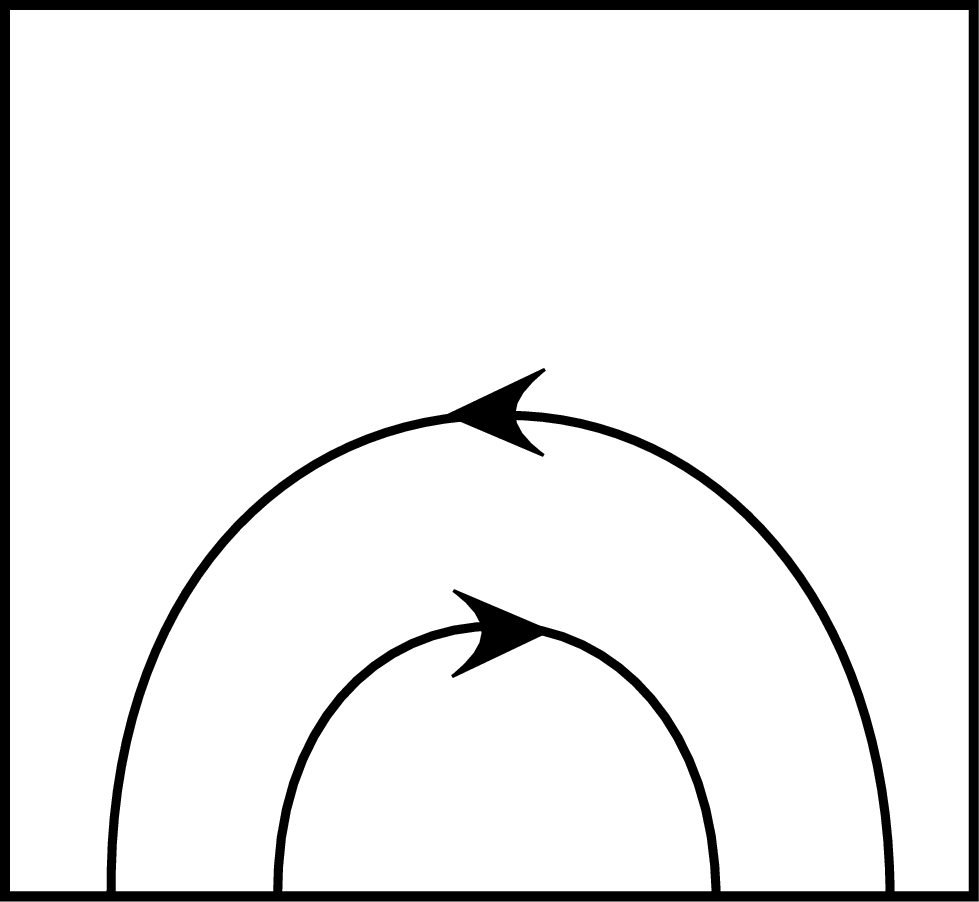} \endxy} (-1)^{k(n-k)}\hspace{10 pt}
&&\vcenter{\xy (0,0)*{\def\svgscale{0.15}\input{arxiv-figures/cap_webs_cc_svg-tex.eps_tex}} \endxy} \mapstosize{2} \vcenter{\xy (0,0)*{\def\svgscale{0.15}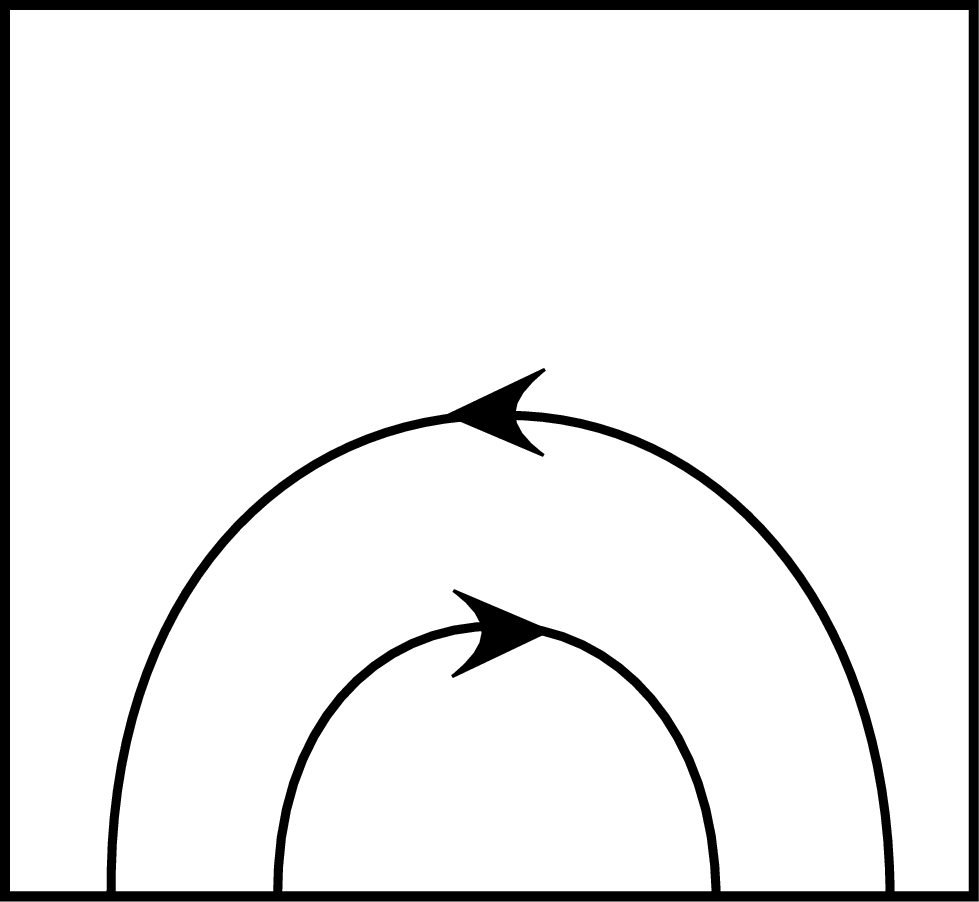} \endxy} \vspace{5 pt}\\
\notag&\vcenter{\xy (0,0)*{\def\svgscale{0.15}\input{arxiv-figures/cup_webs_c_svg-tex.eps_tex}} \endxy} \mapstosize{2}  \vcenter{\xy (0,0)*{\def\svgscale{0.15}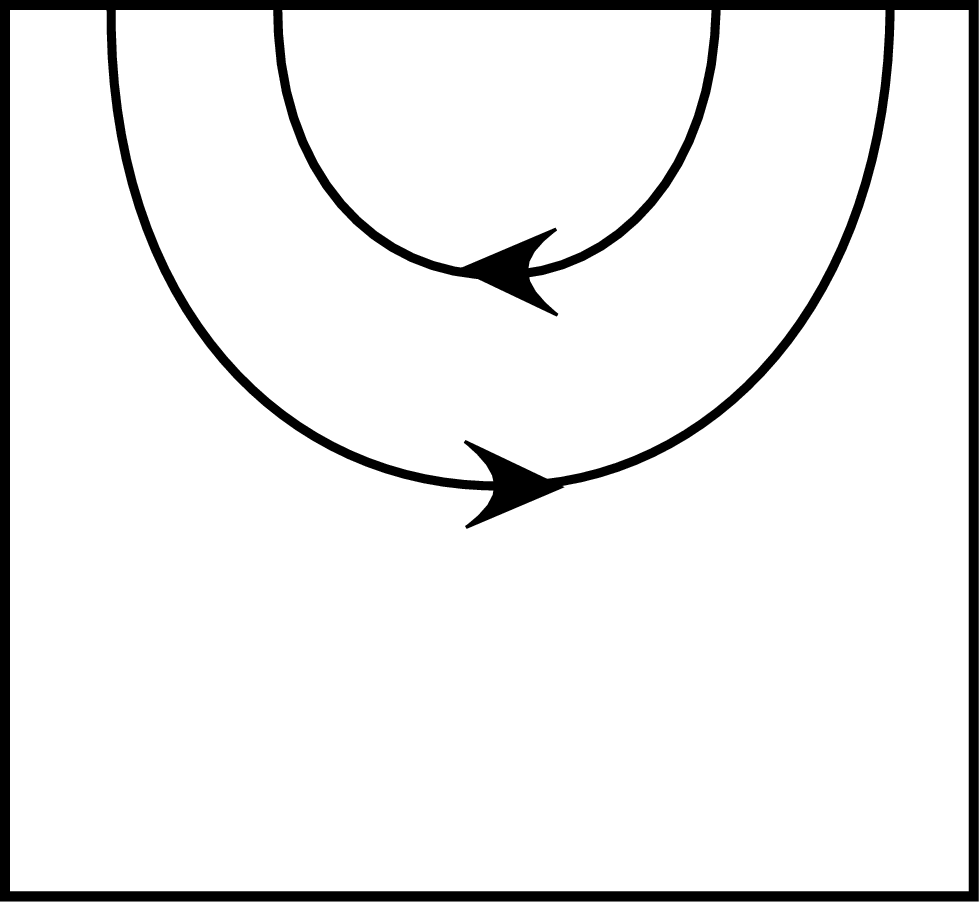} \endxy}
&&\vcenter{\xy (0,0)*{\def\svgscale{0.15}\input{arxiv-figures/cup_webs_cc_svg-tex.eps_tex}} \endxy} \mapstosize{2} \vcenter{\xy (0,0)*{\def\svgscale{0.15}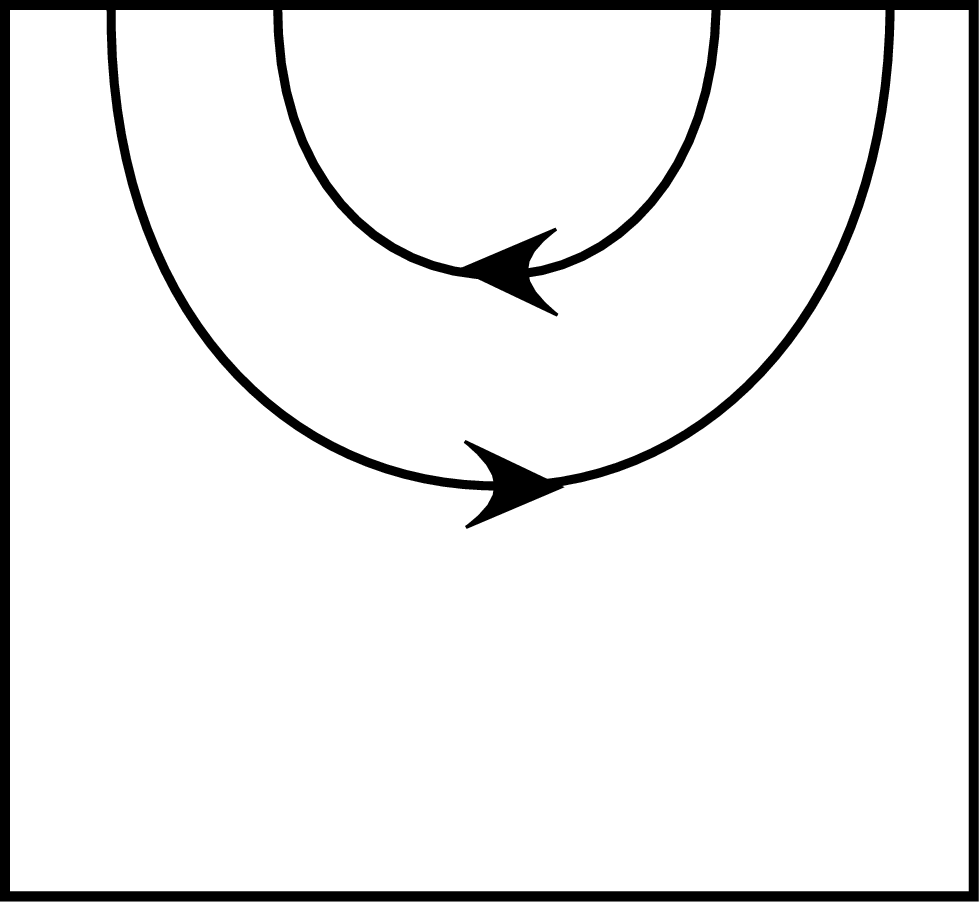} \endxy} (-1)^{k(n-k)}\vspace{5pt}\\
\notag&\vcenter{\xy (0,0)*{\def\svgscale{0.15}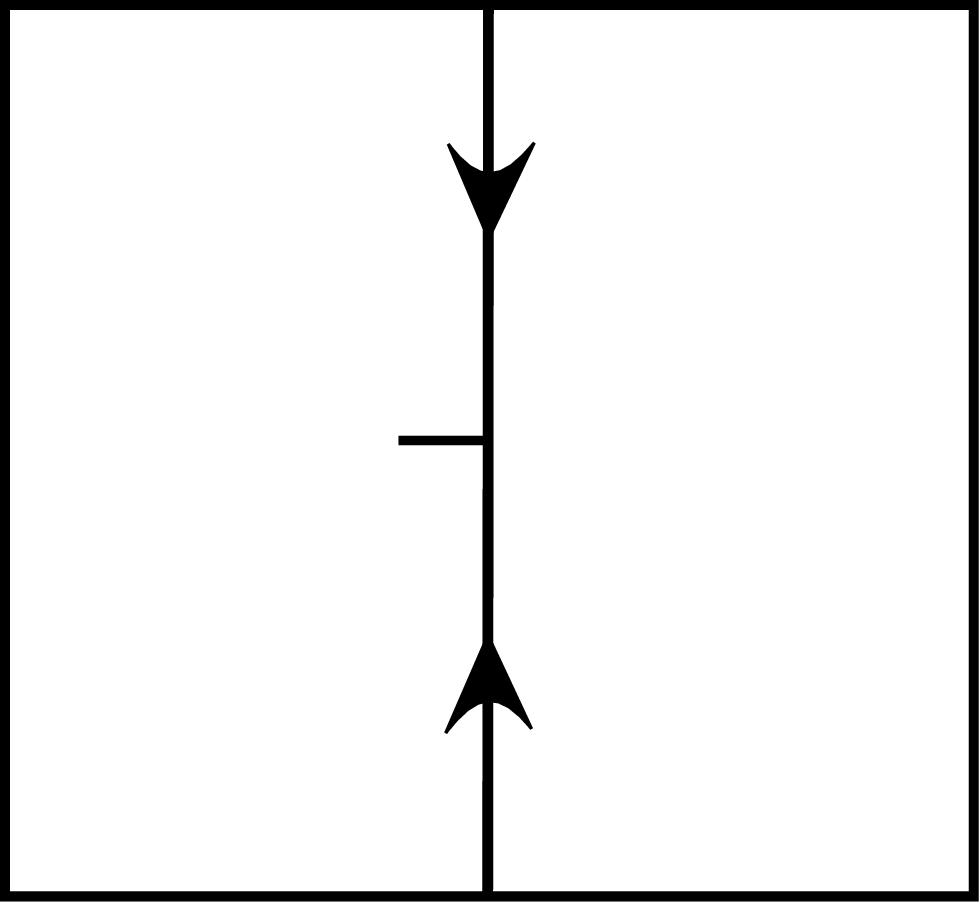} \endxy}  \mapstosize{2} \vcenter{\xy (0,0)*{\def\svgscale{0.15}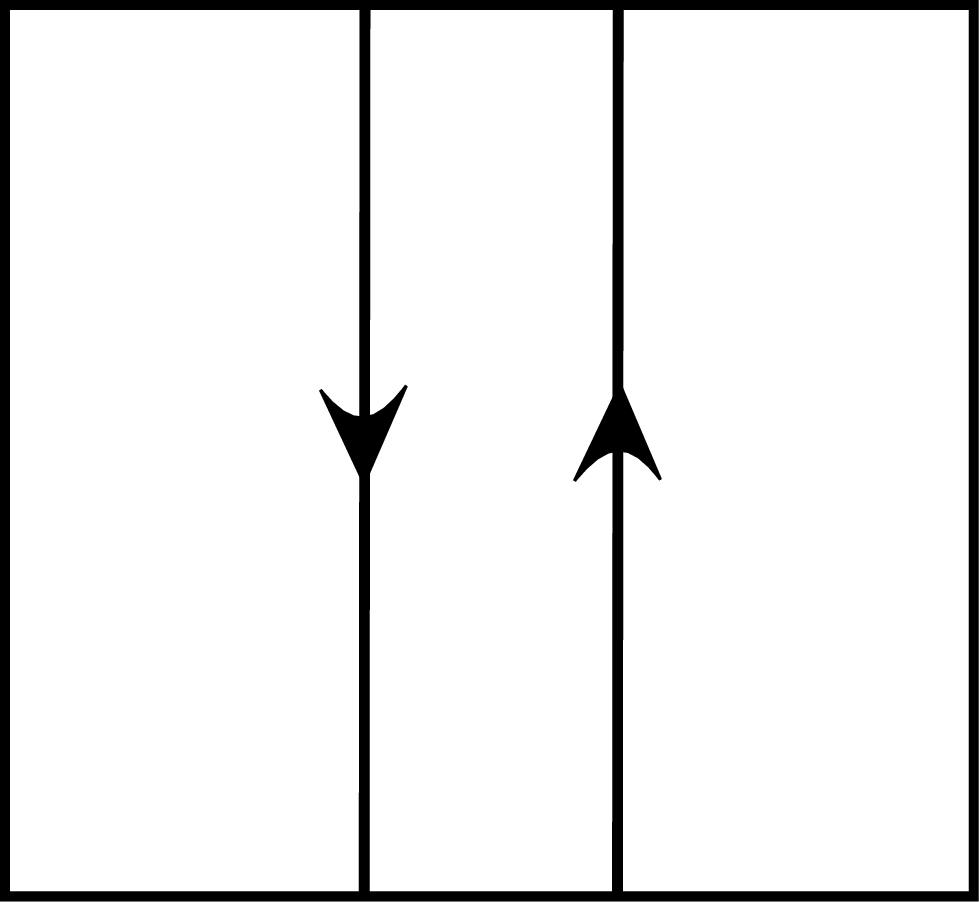} \endxy}\; \lambda_k 
&&\vcenter{\xy (0,0)*{\def\svgscale{0.15}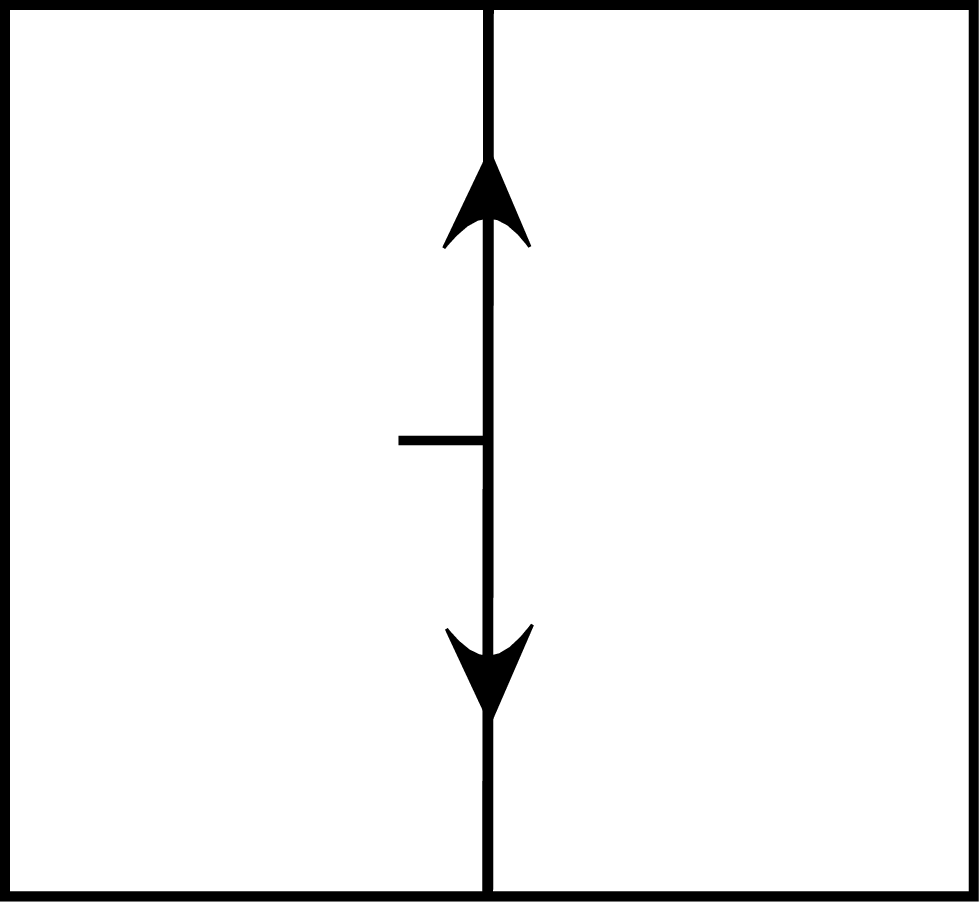} \endxy} \mapstosize{2} \vcenter{\xy (0,0)*{\def\svgscale{0.15}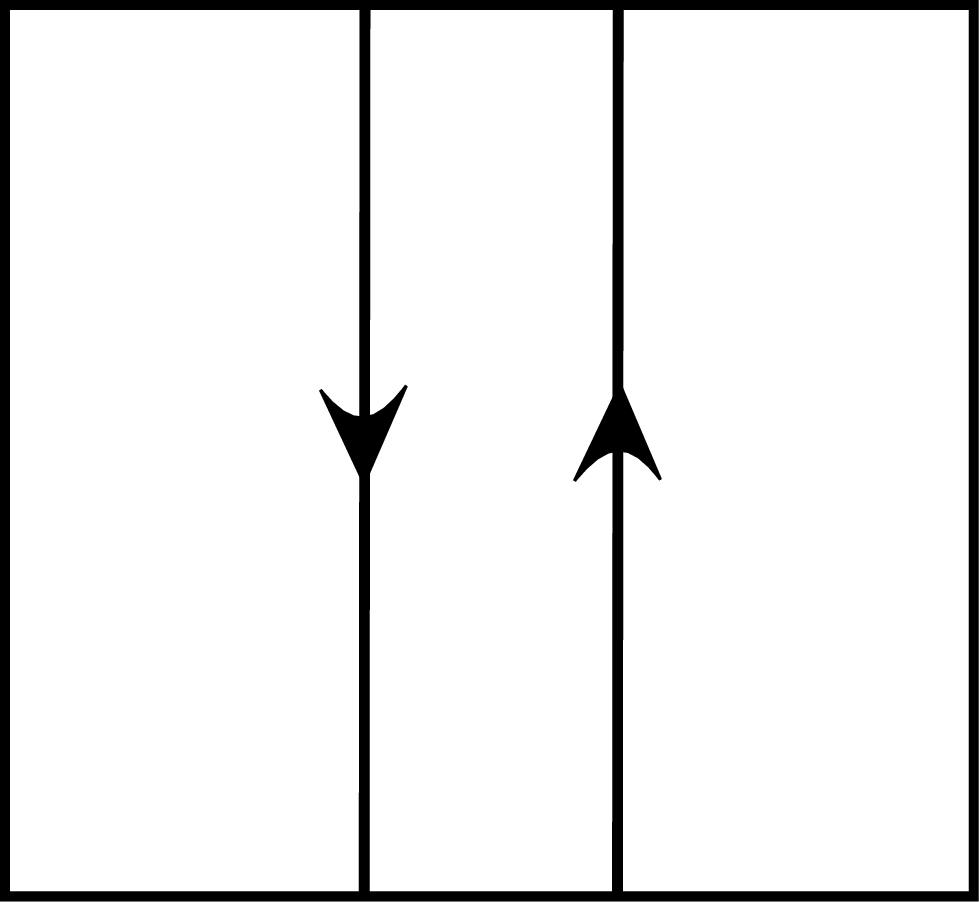} \endxy} \;\nu_k\vspace{5pt}\\
&\vcenter{\xy (0,0)*{\def\svgscale{0.15}\input{arxiv-figures/trivalent_vertex_up_svg-tex.eps_tex}} \endxy}    \mapstosize{2} \vcenter{\xy (0,0)*{\def\svgscale{0.15}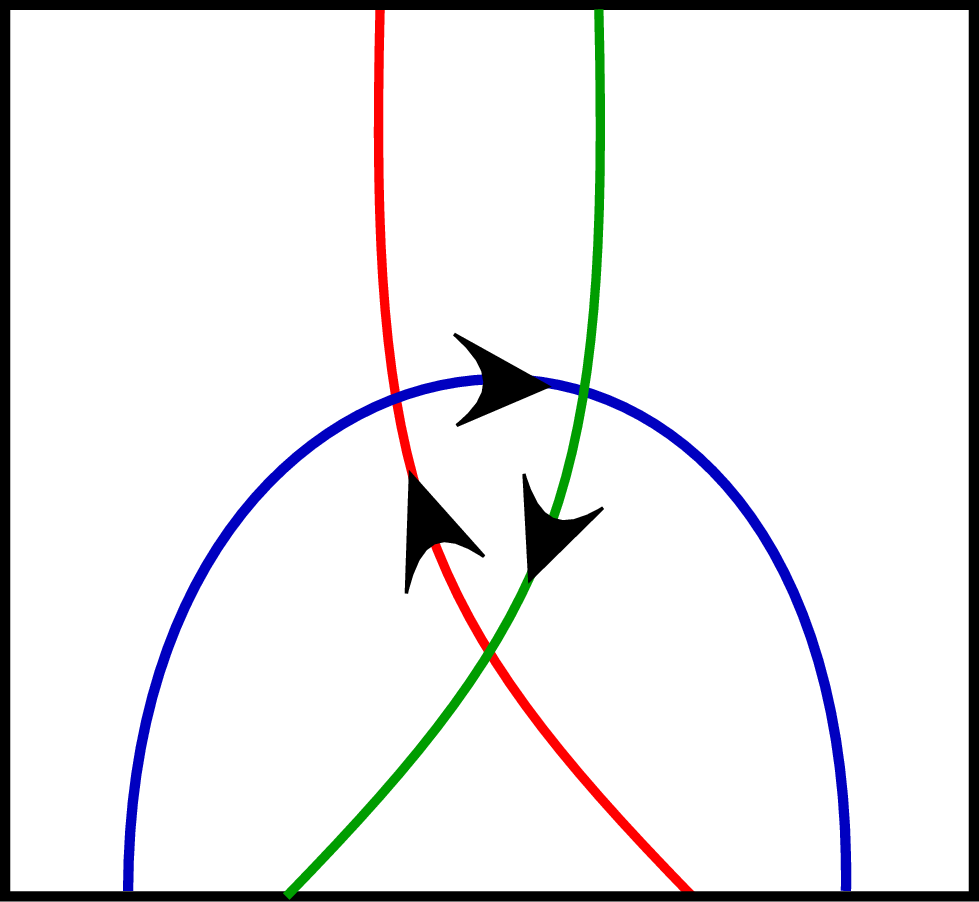} \endxy}\; \lambda_{k,l} 
&&\vcenter{\xy (0,0)*{\def\svgscale{0.15}\input{arxiv-figures/trivalent_vertex_down_svg-tex.eps_tex}} \endxy}  \mapstosize{2} \vcenter{\xy (0,0)*{\def\svgscale{0.15}\input{arxiv-figures/image_of_down_trivalent_svg-tex.eps_tex}} \endxy}\; \nu_{k,l}.
\end{align}
\end{defn}

\begin{defn}\label{defn convention for coefficients of trivalent with labels n} We introduce the following convention. Let $\lambda_{k,n-k}:=\lambda_k$ and $\nu_{k,n-k}:=(-1)^{k(n-k)}\nu_{n-k}$. \end{defn}

\begin{remark}\label{rmk compatibility of convention for lambdanu-k,(n-k) as lambdanu-k}
This convention is compatible with the convention that a tag can be thought of as a trivalent vertex with one strand (replacing the tag) having label $n$. 
\end{remark}

\begin{theorem} \label{thm dQGS functor}
The $2$-functor $\GS_{\zeta}(\lambda, \nu)$ is well-defined if and only if the scalars $(\lambda, \nu)$ satisfy the relations below for all $k, l, m$. By convention, if the indices are not within the range $\{1, \ldots, n-1\}$ then omit the relation. 
\begin{subequations} \label{scalarconditions}
\begin{gather}
   \lambda_k= \lambda_{n-k}, \quad \nu_k=\nu_{n-k} \label{cond-tag switching}\\ 
    \lambda_k= \nu_k^{-1}\label{cond-tag cancellation}\\ 
    \lambda_{k,l}\lambda_{k+l,m}= \lambda_{k,l+m}\lambda_{l,m}, \quad \nu_{k, l+m}\;\nu_{l,m}=\nu_{k+l,m}\nu_{k,l} \label{cond-I=H relation}\\  
    \lambda_{k-1,1}\nu_{k-1,1}= (-1)^{k-1}\zeta^{\frac{k(k-1)}{2}}q^{k-1}. \label{cond-bigon killing}
\end{gather}
\end{subequations}
\end{theorem}

\begin{proof}
In \cref{sec well-defined-ness of the functor}, we show that the conditions \cref{cond-tag switching}-\cref{cond-bigon killing} are the only constraints implied by the web relations in \cref{fig:webrel} (along with their mirror images and arrow reversals.)
\end{proof}

\begin{example} \label{goodchoiceoflambdanu}
One can check that 
\begin{equation} \lambda_{k,l}:=\zeta^{\frac{kl(k+l)}{2}}(-q)^{kl}, \quad \nu_{k,l}=1, \quad \lambda_{k}=\nu_k = (-1)^{k(n-k)} \end{equation} is a solution to the above equations.
Note that this example is compatible with the convention in \cref{defn convention for coefficients of trivalent with labels n}.
\end{example}

Having postponed this proof to the next chapter, we spend the remainder of this chapter discussing the scalars $(\lambda, \nu)$.


\begin{lemma}\label{lemma coefficient symmetry}
   Conditions \cref{cond-tag switching}-\cref{cond-bigon killing} imply that \begin{align*}
       \lambda_{k,l}=\lambda_{l,k}, \quad \quad \nu_{k,l}=\nu_{l,k}
   \end{align*}
   for any $1\leq k,l \leq n-1$ such that $k+l \leq n$.
\end{lemma}
\begin{proof}
    Condition \cref{cond-tag switching} implies that the lemma holds when $k+l=n$, so we  need only prove the lemma for $k+l \leq n-1$.

    First we show that $\lambda_{1,m}=\lambda_{m,1}$ for any $1 \leq m < n-1$ using induction on $m$. 
    The base case $m=1$ is trivial, so suppose now that $m>1$. 
    From condition \cref{cond-I=H relation}, we have that
    $$\lambda_{1,m-1}\lambda_{m,1}=\lambda_{1,m}\lambda_{m-1,1}.$$ 
    But $\lambda_{1,m-1}=\lambda_{m-1,1}$ by induction hypothesis, so $\lambda_{m,1}=\lambda_{1,m}$ as desired.
    
    Now we discuss the general case. 
    Condition \cref{cond-I=H relation} implies that
    for any $1\leq k \leq n-1$, $1<l \leq n-1$ such that $k+l \leq n-1$, we have $\lambda_{k,l-1}\lambda_{k+l-1,1}=\lambda_{k,l}\lambda_{l-1,1}$, so that $\dfrac{\lambda_{k,l+1}}{\lambda_{k,l}}=\dfrac{\lambda_{k+l,1}}{\lambda_{l,1}}$.
    We then iteratively have that
    \begin{align*}
     \frac{\lambda_{k,l}}{\lambda_{k,1}}= \frac{\lambda_{k+1,1}}{\lambda_{1,1}}\frac{\lambda_{k+2,1}}{\lambda_{2,1}}\ldots \frac{\lambda_{k+l-1,1}}{\lambda_{l-1,1}}
    \end{align*}
    so that 
    \begin{align}
    \lambda_{k,l}=\lambda_{k,1}\frac{\lambda_{k+1,1}}{\lambda_{1,1}}\frac{\lambda_{k+2,1}}{\lambda_{2,1}}\ldots \frac{\lambda_{k+l-1,1}}{\lambda_{l-1,1}}.\label{eq 1 lemma coefficient symmetry}
    \end{align}
    Similarly, using that
    $\lambda_{1,l-1}\lambda_{l,k}=\lambda_{1,k+l-1}\lambda_{l-1,k}$, we have that 
    \begin{align}
      \lambda_{l,k} =\lambda_{1,k}\frac{\lambda_{1,k+1}}{\lambda_{1,1}}\frac{\lambda_{1,k+2}}{\lambda_{1,2}}\ldots \frac{\lambda_{1,k+l-1}}{\lambda_{1,l-1}} \label{eq 2 lemma coefficient symmetry}.
    \end{align}
    The RHS of \cref{eq 1 lemma coefficient symmetry} and \cref{eq 2 lemma coefficient symmetry} are equal by the special case of the lemma when one of the indices is 1, so we have that $\lambda_{k,l}=\lambda_{l,k}$ as desired.
\end{proof}

\begin{lemma}\label{lemma symmetry of GS-zeta}
The 2-functor $\GS_\zeta$ intertwines the symmetries $\sigma$ and $\Upsilon_W$ of $\cwebs$ with the symmetries $\sigma$ and $\Upsilon_S$ of $\DiagSBS$ respectively. It intertwines $\rotation$ with $\rotation$, up to sign. More precisely, we have 
\begin{equation}\label{eq GS-zeta intertwines symmetries}
\GS_\zeta \circ \sigma=\sigma \circ \GS_\zeta, \quad \GS_\zeta \circ \Upsilon_W=\Upsilon_S \circ \GS_\zeta,  \quad \rotation \circ \GS_{\zeta} = \GS_{\zeta} \circ \rotation \circ \scalewebs,
\end{equation} 
where $\scalewebs$ is the object-fixing automorphism of $\cwebs$ that multiplies each web by $(-1)^{k(n-k)}$ for each upward-oriented $k$ in the source or target (downward-oriented strands have no effect)\footnote{For example, $\scalewebs$ multiplies all the generators in Definition 5.4 by $(-1)^{k(n-k)}$, except for the trivalent vertices which it multiplies by $(-1)^{k(n-k) + l(n-l) + (k+l)(n-k-l)} = 1$. The 2-functor $\scalewebs$ is isomorphic to the identity 2-functor; the natural isomorphism between them is multiplication by the appropriate sign on each object (a factor of $(-1)^{k(n-k)}$ for each upward-oriented $k$-labeled strand). Hence $\rotation\circ\GS_\zeta\cong\GS_\zeta \circ \rotation$, so that $\GS_\zeta$ also weakly intertwines with rotation.}.
\end{lemma}

\begin{proof}
    It suffices to check this on the generating 2-morphisms, under the assumptions in \cref{thm dQGS functor}. This check is an easy exercise and is left to the reader, where \cref{lemma coefficient symmetry} may be used to check the intertwining of $\Upsilon_W$ and $\Upsilon_S$.

\end{proof}

\begin{remark} The symmetries are intertwined not because of any application of relations, but already at the level of the diagrams themselves. 
To elaborate, one could define a ``free pivotal'' version $\Webs^{\Om, \pre}$ of webs which have isotopy relations but no other relations. 
One could define a version $\widetilde{\GS}_{\zeta}: \Webs^{\Om, \pre} \rightarrow \DiagSBS^{\pre}$ of $\GS_{\zeta}$; this functor would still involve scalars $\lambda$ and $\nu$ satisfying \eqref{scalarconditions}. 
Then \cref{lemma symmetry of GS-zeta} applies to $\widetilde{\GS}_{\zeta}$ as well, and descends to the corresponding result for $\GS_{\zeta}$. We mention this point because, as is common in the field, we use symmetries to check the well-definedness of $\GS_{\zeta}$ itself. We check certain relations, and then claim others hold by applying $\sigma$ or $\Upsilon$ or rotation. There is no fallacy in applying these symmetries before knowing of the existence of the $2$-functor, because we are really applying the symmetries to $\widetilde{\GS}_{\zeta}$. \end{remark}

\begin{remark} One might hope for a symmetry of $\DiagSBS$ corresponding to the arrow-reversing symmetry $\reversal$ of $\cwebs$. It would act on objects and $1$-morphisms by $\tau$. As evidenced by \eqref{diag-dQGS on 2-morphisms}, the inward-pointing tag (with label $k$ on bottom) is sent to an identity map times $\lambda_k$, and any autoequivalence of $\DiagSBS$ sends identity maps to identity maps. Applying $\reversal$ we obtain the outward-pointing tag, which goes to the corresponding identity map times $\nu_k$. So an intertwining symmetry of $\DiagSBS$ is only possible if $\lambda_k = \nu_k$ for all $k$.  Since $\lambda_k = \nu_k^{-1}$, we have $\lambda_k = \nu_k = \pm 1$. This constraint holds for \cref{goodchoiceoflambdanu}.

However, unlike all the other symmetries, this autoequivalence of $\DiagSBS$ would need to rescale diagrams (rather than merely manipulating them). As evidenced by \eqref{diag-dQGS on 2-morphisms}, a sign $(-1)^{k(n-k)}$ appears in the comparison of the different orientations of a $k$-labeled cap or cup. Worse still, $\lambda_{k,l}$ and $\nu_{k,l}$ are swapped in the comparison for trivalent vertices, and for \cref{goodchoiceoflambdanu} this would entail rescaling some diagrams with both signs and powers of $\zeta$. It is possible that such a symmetry exists.

\PTv2{I did try for more than an hour, but couldn't figure it out. Let's not add it in v1?}
\end{remark}

\subsection{Fullness and Faithfulness}\label{subsec fullness and faithfulness}
We now work over the fraction field of the base ring $A$ in \cref{subsec the functor}. The fraction field is a finite extension of $\Q(q)$ which we denote $\Q(q,\zeta)$.
Throughout this section, we use $\GS_\zeta$ to refer to the 2-functor $\GS_\zeta(\lambda, \nu)$ in \cref{subsec the functor} for any valid choice of coefficients (see \cref{thm dQGS functor}), and reuse the same notation for the functor (as well as the underlying categories) obtained after extending scalars along the localization map $A \hookrightarrow \Q(q,\zeta)$.

Recall from \cref{subsec-colored sln webs} that we have a factorization $\qFund \xrightarrow{\sim} \cwebs \rightarrow \qRep$ of the canonical inclusion of $\qFund$ into $\qRep$, and that the Karoubi envelope of $\qFund$ coincides with $\qRep$.
The 2-functor $\GS_\zeta: \cwebs \rightarrow \DiagSBS \hookrightarrow \SBSBim$ factors through $m\SBSBim$. We denote the induced map $\cwebs \rightarrow m\SBSBim$ again by $\GS_\zeta$.
Recall from \cref{subsec reverse bott-samelsons} that the Karoubi envelope of (the graded, additive completion of) $m\SBSBim$ coincides with the full sub-category $m\SSBim$ of $\SSBim$.
The 2-functor 
$$\qFund \xrightarrow{\sim} \cwebs \xrightarrow{\GS_\zeta} m\SBSBim$$
then canonically extends to the respective Karoubi envelopes, and we obtain an additive 2-functor $\QGS: \qRep \rightarrow m\SSBim$.

\begin{lemma}\label{lemma GS functor sends simples to indecomposables}
    The extended 2-functor $\QGS$ categorifies the reformulated Satake isomorphism in \cref{thm Soergel Satake}.
    For any $\lambda \in \Lambda_{\wt}^+,$and $a\in \Om$, $\QGS$ sends $L_\lambda \in \Hom_{\qRep}(a,a+\rot(\lambda))$ to the indecomposable singular Soergel bimodule $B_{\psi_a(\lambda)}$ in $\SBimod{a+\rot(\lambda)}{a}$. 
    In particular, 
\begin{equation} \label{thistoogoestoKL} {}_{\hh{a+\rot(\lambda)}}\ch_{\hh{a}}(B_{\psi_a(\lambda)})=\KLB{\hh{a+\rot(\lambda)}}{\hh{a}}{\psi_a(\lambda)}\end{equation}
    so that Soergel's conjecture holds for the maximally singular Soergel bimodules when constructed from the deformed affine realization over $\Q(q,\zeta)$.
\end{lemma}

\begin{proof}
By construction, $\GS_\zeta$ sends a colored fundword $(\ula, a_d)$ (which represents a 1-morphism in $\cwebs$) to the singlestep expression $\GS(\ula, a_d)$ (see \cref{def:GSfundword}), which corresponds to the bimodule $\BS(\GS(\ula, a_d))$ in $m\SBSBim$. In particular, $\GS_\zeta$ sends the generating object $((\varpi_k),a)$ to  
\[ \BS(\GS((\varpi_k),a)) = \BS([\hh{a+k}\supset \hh{a, a+k}\subset \hh{a}]) = R^{\hh{a,a+k}}(\ell(\hh{a,a+k})-\ell(\hh{a})) \]
where $R^{\hh{a,a+k}}$ is viewed as a $(R^{\hh{a+k}}, R^{\hh{a}})$-bimodule. Since $R^{\hh{a+k}}$ and $R^{\hh{a}}$ generate $R^{\hh{a,a+k}}$ as a ring (see \cref{defn-condition star}, \cref{lemma deformed affine frobenius satisfies extra conditions}), an endomorphism of this bimodule is determined by where it sends $1$, thus the endomorphism ring is $R^{\hh{a,a+k}}$ itself. Hence this bimodule is indecomposable. As $\GS((\varpi_k),a)$ is a reduced expression for $\psi_a(\varpi_k)$ (cf. \cref{example rex for fund weights}), we have
\[ \BS(\GS((\varpi_k),a)) \cong B_{\psi_a(\varpi_k)}.\]
It is also a straightforward computation in the Grothendieck group that its character is equal to both the standard and the Kazhdan-Lusztig basis element associated to the double coset $\psi_a(\varpi_k)$ (these are equal because $\psi_a(\varpi_k)$ is minimal in the Bruhat order on $W_{\hh{a+k}}\backslash W /W_{\hh{a}}$). That is,
\[ {}_{\hh{a+k}}\ch_{\hh{a}}(B_{\psi_a(\varpi_k)})=\KLB{\hh{a+k}}{\hh{a}}{\psi_a(\varpi_k)}=\HAB{\hh{a+k}}{\hh{a}}{\psi_a(\varpi_k)}.\]
Thus the decategorified functor $[\QGS]:[\qRep]=\KRepom \rightarrow \HAbd$ induced by $\QGS$ coincides with the reformulated Satake isomorphism in \cref{thm Soergel Satake}, since they agree on $\ch(L_{\varpi_k}) \in \KRepom(a,a+k)$ for any $1 \leq k\leq n-1$, $a \in \Om$, which generate $[\qRep]$ as an algebroid.

Since $\QGS$ categorifies the reformulated Satake isomorphism, it follows that 
\begin{equation} \label{itgoestoKL} {}_{\hh{a+\rot(\lambda)}}\ch_{\hh{a}}(\QGS(L_\lambda))=\KLB{\hh{a+\rot(\lambda)}}{\hh{a}}{\psi_a(\lambda)} \end{equation} for all $\lambda \in \Lambda^+_{\wt}$ such that $L_\lambda \in \Hom_{\qRep}(a, a+\rot(\lambda))$. The Hom-formula (\cref{thm-hom formula}) then implies that the degree zero endomorphism space $\End^0(\QGS(L_\lambda))$ is one dimensional, so it has no non-trivial idempotents, whence $\QGS(L_\lambda)$ is indecomposable.

Thus the Williamson-Soergel categorification theorem (\cref{thm-categorification}(2)) along with \cref{itgoestoKL} implies that $\QGS(L_\lambda) \cong B_{\psi_a(\lambda)}$. One then deduces \eqref{thistoogoestoKL}  from \eqref{itgoestoKL}, verifying the Soergel conjecture for this double coset. Since all double cosets in $W_{\hh{b}}\backslash W/W_{\hh{a}}$ are of the form $\psi_a(\lambda)$ for some $\lambda \in \Lambda^+_{\wt}$, we have the desired statement regarding Soergel's conjecture.
\end{proof}

We now proceed to address fullness and faithfulness of the 2-functor $\QGS$.

\begin{lemma}
   The images of the 2-morphisms in $\qRep$ via $\QGS$ have degree zero.
\end{lemma}
\begin{proof}
    It suffices to show that the images of the generators of $\qFund$ are of degree zero. For the images of cups (resp. caps), the degrees of the clockwise and counter-clockwise cups (resp. caps) cancel each other so that the images are of degree zero.
    The images of tags, being just identity maps up to a scalar, have degree zero.
    Finally, that the images of the trivalent vertices are of degree zero follows similarly, using \cref{lemma image of trivalent has degree 0} (and its analogue after applying $\Upsilon_S$).
\end{proof}

By a 2-ideal of $\qRep$, we mean a sub-2-category of $\qRep$ which is closed under both horizontal (left and right) and vertical (top and bottom) compositions by arbitrary (composable) morphisms in $\qRep$.

\begin{lemma}\label{lemma monoidal ideals in qRep}
    The semisimple 2-category $\qRep$ has no proper non-zero 2-ideals.
\end{lemma}

\begin{proof}
   This is an adaptation of a well-trodden argument, and relies on the fact that we work over a generic field. We use the notation $L_0^{(b)}$ to refer to the trivial representation in the 1-morphism category $\Hom_{\qRep}(b,b)$ for any $b \in \Om$.
   
   Let $\mc{I}$ be a non-zero 2-ideal.  
   Since $\qRep$ is a semisimple 2-category and $\mc{I}$ is a non-trivial 2-ideal, we have that $\mc{I}$ contains the identity 2-morphism $\id_L$ of some object $L$ in $\Hom_{\qRep}(a,a')$ for some $a,a' \in \Om$. 
   Then $L^*\otimes L$ in $\Hom_{\qRep}(a,a)$ has the trivial representation $L_0^{(a)}$ as a direct summand, so that $\id_{L_0^{(a)}}$ is in $\mc{I}$.
   For any $b \in \Om$, we then have $\id_{L_{\varpi_{b-a}}\otimes L_0^{(a)} \otimes L_{\varpi_{a-b}}}$ is a 2-morphism in $\mc{I}$. 
   Since $L_0^{(b)}$ is a summand of $L_{\varpi_{b-a}}\otimes L_0^{(a)} \otimes L_{\varpi_{a-b}}$, it follows that $\id_{L_0^{(b)}}$ is in $\mc{I}$ for any $b \in \Om$. For any object $M$ in $\Hom_{\qRep}(a,a')$ for any $a,a' \in \Om$, it follows that $\id_M = \id_M \ot \id_{L_0^{(a)}}$ is in $\mc{I}$. Containing the identity 2-morphisms of all 1-morphisms, one concludes that $\mc{I}$ contains all $2$-morphisms as well.
\end{proof}

\begin{lemma}
    The 2-functor $\QGS$ is faithful on 2-morphisms, i.e., the maps induced on the 2-morphism spaces are injective.
\end{lemma}

\begin{proof}
    Let $\ker (\QGS)$ denote the 2-ideal consisting of 2-morphisms $f$ such that $\QGS(f)=0$.
    Since $\QGS$ sends the identity 2-morphism of $L_{\varpi_k} \in \Hom_{\qRep}(0,k)$ to the identity 2-morphism of the corresponding indecomposable singular Soergel bimodule $B_{\psi_0(\varpi_k)}$, we have that these identity 2-morphisms are not in $\ker(\QGS)$, so that $\ker(\QGS)$ is a proper ideal. \cref{lemma monoidal ideals in qRep} then  implies that $\ker(\QGS)$ is trivial, so that $\QGS$ is faithful on 2-morphisms\footnote{An alternative proof which ignores the monoidal structure goes as follows. Using \cref{lemma GS functor sends simples to indecomposables} one deduces that the kernel does not contain the identity map of any simple $1$-morphism. From the semisimplicity of $\qRep$, every non-zero ideal contains the identity map of some simple $1$-morphism.}.
\end{proof}


\begin{lemma}\label{lemma GS  on 2-morphism spaces}
    The maps induced by $\QGS: \qRep \rightarrow m\SBSBim$ on the 2-morphism spaces are isomorphisms onto the space of degree zero 2-morphisms.
    Moreover, there are no negative degree 2-morphisms between 1-morphisms in the image of $\QGS$.
\end{lemma}
\begin{proof}
    We have already shown that the induced maps on 2-morphism spaces are injective. 
    Hence to show that the induced maps at the 2-morphism level are isomorphisms, it suffices to show that the 2-morphism spaces have the same (finite) dimensions.
    Since $\Rep_q$ is semisimple, we have that every 1-morphism in $\qRep$ is a direct sum of simple representations. 
    Given $N=\bigoplus_{\lambda}L_{\lambda}^{\oplus m_\lambda}$ and $M=\bigoplus_\lambda L_{\lambda}^{\oplus l_\lambda}$ in $\Hom_{\qRep}(a,b)$ (with only finitely many $m_\lambda, l_\lambda$ being non-zero, and sum running over dominant weights $\lambda$ such that $\rot(\lambda)=b-a$), additivity of $\QGS$ along with \cref{lemma GS functor sends simples to indecomposables} implies that $\QGS(N)=\bigoplus_\lambda B_{\psi_a(\lambda)}^{\oplus m_\lambda}$ and $\QGS(M)=\bigoplus_{\lambda} B_{\psi_a(\lambda)}^{\oplus l_\lambda}$.
    From \cref{lemma GS functor sends simples to indecomposables}, we have that $\ch(B_p)=\KLB{\hh{b}}{\hh{a}}{p}$ for any $p \in W_{\hh{b}}\backslash W/W_{\hh{a}}$.
    Using the Hom formula (\cref{thm-hom formula}), we then have that 
    $$\grk \Hom^\bullet_{\SBimod{\hh{b}}{\hh{a}}}(\QGS(N),\QGS(M))
    = \left(\sum_{\lambda} m_\lambda \KLB{\hh{b}}{\hh{a}}{\psi_a(\lambda)}, \sum_{\lambda'} l_{\lambda'} \KLB{\hh{b}}{\hh{a}}{\psi_a(\lambda')}\right)
    =\sum_{\lambda, \lambda'}m_\lambda l_{\lambda'}\delta_{\lambda, \lambda'} + vr$$
    for some $r \in \Z[v]$. 
    This implies that there are no negative degree morphisms from $\QGS(N)$ to $\QGS(M)$, and that $\dim \Hom^0(\QGS(N),\QGS(M))=\sum_{\lambda} m_\lambda l_\lambda$, 
    which (by Schur's Lemma) equals
    $\dim \Hom_{\qRep}(N,M)$.
\end{proof}

Recall the following notion of a \emph{degree zero equivalence}.
\begin{definition}\cite[Definition 2.14]{EQuantumI}
    A $2$-functor $\mc{G}$ from a $\kk$-linear $2$-category $\mc{A}$ to a $\kk$-linear graded $2$-category $\mc{B}$ is a \emph{degree-zero equivalence} if the following properties
hold. \begin{enumerate}[label=(\alph*)] 
\item $\mc{G}$ induces a bijection between the objects. 
\item $\mc{G}$ induces an isomorphism $\Hom_{\mc{A}}(M,N) \to \Hom_{\mc{B}}^0(\mc{G}(M),\mc{G}(N))$ to the space of
morphisms of degree $0$ (for 1-morphisms $M, N$ in $\mc{A}$ with same source and target). 
\item $\Hom_{\mc{B}}^k(\mc{G}(M),\mc{G}(N))=0$ for all $k<0$ (for 1-morphisms $M$ and $N$ in $\mc{A}$). 
\item Every $1$-morphism in $\mc{B}$ is isomorphic to a direct sum of grading shifts of 1-morphisms of the form $\mc{G}(M)$
for $1$-morphisms $M$ in the Karoubi envelope of $\mc{A}$; i.e. $\mc{G}$ is \emph{essentially surjective up to grading shift}. \end{enumerate}
\end{definition}

The following theorem is then an immediate consequence of \cref{lemma GS functor sends simples to indecomposables}, \cref{lemma GS  on 2-morphism spaces} and the Williamson-Soergel categorification theorem.

\begin{theorem}\label{thm dQGS functor is a degree zero equivalence}
    The 2-functors $\GS_\zeta: \cwebs \rightarrow m\SSBim$, $\QGS: \qRep \rightarrow m\SSBim$ are degree zero equivalences.
\end{theorem}

\begin{remark}
    The definition of a degree zero equivalence, especially conditions (b) and (c), may remind one of the inclusion of the abelian heart of a t-structure into a (monoidal) dg-category (provided the monoidal structure restricts nicely to the heart). Here is an explicit geometric example relevant to us.
    Let $\mc{O}, \mc{K}$ respectively denote the formal power series ring over $\CC$ and its fraction field. Fix sheaf coefficients to be a field $\kk$ of characteristic zero, and let $\Semis_{\SL_n({\mc{O}})\times \SL_n({\mc{O})}}(SL_n({\mc{K}}))$ denote the category of ${\SL_n({\mc{O}})\times \SL_n({\mc{O})}}$-equivariant semisimple complexes of perverse sheaves on the loop group $\SL_n(\mc{K})$, with $\Hom^{\bullet}(A,B)$ given by  $\bigoplus_i \Hom_{D^b}(A,B[i])$ (where $D^b$ is the ${\SL_n({\mc{O}})\times \SL_n({\mc{O})}}$-equivariant constructible derived category on the loop group).
    Then the inclusion of the perverse heart
    \begin{equation} \label{eq inclusion perverse heart} \Perv_{\SL_n({\mc{O}})\times \SL_n({\mc{O})}}(SL_n({\mc{K}})) \hookrightarrow \Semis_{\SL_n({\mc{O}})\times \SL_n({\mc{O})}}(SL_n({\mc{K}}))\end{equation}
    is a degree zero equivalence (considering both sides as monoidal categories with the convolution product).
    While there are no non-trivial extensions between perverse sheaves within the abelian heart, there may be non-vanishing positive degree morphisms between them in the equivariant derived category.
    From usual geometric Satake equivalence, the LHS of \cref{eq inclusion perverse heart} is equivalent to $\Rep_{\kk}(\PGL_n)$.
    The RHS of \cref{eq inclusion perverse heart} may moreover be identified with $ \SSBim(\hh{0},\hh{0})$ constructed from the Kac-Moody realization, so that the inclusion $\Rep(\PGL_n)\hookrightarrow \SSBim(\hh{0},\hh{0})$ so obtained is the restriction of the (non-quantized, Soergelified) Geometric Satake 2-functor $\mathbb{GS}: \Rep^{\Om} \rightarrow m\SSBim$ to the full subcategory with a single object $0 \in \Om$. 
\end{remark}

Summarizing the results from \cref{thm dQGS functor}, \cref{goodchoiceoflambdanu}, \cref{lemma symmetry of GS-zeta}, \cref{thm dQGS functor is a degree zero equivalence}, we have the following theorem.

\begin{theorem}\label{thm main}
Let $A$ denote the ring $\begin{cases}
\Q[q^{\pm 1}, \zeta^{\pm 1}]/(q^{-2}-\zeta^n) &\text{if $n$ odd;}\\
\Q[q^{\pm 1}, \zeta^{\pm 1}]/(q^{-1}-\zeta^{n/2}) &\text{if $n$ even.}
    \end{cases}$.
\begin{enumerate}
   \item Over the base ring $A$, there is a 2-functor 
   $$\GS_\zeta: \cwebs \rightarrow \DiagSBS$$ explicitly defined by \cref{eq dQGS on objects}, \cref{eq dQGS on 1-morphisms}, \cref{diag-dQGS on 2-morphisms} for any choice of scalars $\lambda, \nu$ satisfying conditions \cref{cond-tag switching}-\cref{cond-bigon killing},  with \cref{goodchoiceoflambdanu} offering a concrete choice.
    
  \item   The 2-functor $\GS_\zeta$ intertwines the symmetries $\sigma$, $\rotation$, and $\Upsilon_W$ of $\cwebs$ with the symmetries $\sigma$, $\rotation$, and $\Upsilon_S$ of $\DiagSBS$.
    
  \item  After passing to the fraction field $\Q(q,\zeta)$ of $A$, the 2-functor $\GS_\zeta$ naturally extends to a 2-functor $$\QGS: \qRep \rightarrow m\SSBim.$$
  \item Over the fraction field $\Q(q,\zeta)$ of $A$, the 2-functors 
    $\GS_\zeta: \cwebs \rightarrow m\SSBim$ and $\QGS: \qRep \rightarrow m\SSBim$ are degree zero equivalences, and $\QGS$ categorifies the reformulated Satake isomorphism in \cref{thm Soergel Satake}.
  \item Soergel's conjecture holds for the indecomposable maximally singular Soergel bimodules constructed from the deformed affine realization over the fraction field $\Q(q,\zeta)$ of $A$, i.e.,
  $${}_{\hh{b}}\ch_{\hh{a}}(B_{p})=\KLB{\hh{b}}{\hh{a}}{p}$$
  for any $p \in W_{\hh{b}} \backslash W/ W_{\hh{a}}.$
  
    \end{enumerate}
\end{theorem}

%% file: arxiv-figures/cap_webs_c_svg-tex.eps_tex
\begingroup%
  \makeatletter%
  \providecommand\color[2][]{%
    \errmessage{(Inkscape) Color is used for the text in Inkscape, but the package 'color.sty' is not loaded}%
    \renewcommand\color[2][]{}%
  }%
  \providecommand\transparent[1]{%
    \errmessage{(Inkscape) Transparency is used (non-zero) for the text in Inkscape, but the package 'transparent.sty' is not loaded}%
    \renewcommand\transparent[1]{}%
  }%
  \providecommand\rotatebox[2]{#2}%
  \newcommand*\fsize{\dimexpr\f@size pt\relax}%
  \newcommand*\lineheight[1]{\fontsize{\fsize}{#1\fsize}\selectfont}%
  \ifx\svgwidth\undefined%
    \setlength{\unitlength}{469.76640848bp}%
    \ifx\svgscale\undefined%
      \relax%
    \else%
      \setlength{\unitlength}{\unitlength * \real{\svgscale}}%
    \fi%
  \else%
    \setlength{\unitlength}{\svgwidth}%
  \fi%
  \global\let\svgwidth\undefined%
  \global\let\svgscale\undefined%
  \makeatother%
  \begin{picture}(1,0.92106259)%
    \lineheight{1}%
    \setlength\tabcolsep{0pt}%
    \put(0,0){\includegraphics[width=\unitlength]{cap_webs_c_svg-tex.eps}}%
    \put(0.38669729,0.26482665){\color[rgb]{0,0,0}\makebox(0,0)[lt]{\begin{minipage}{0.31158253\unitlength}\raggedright \gr{a-k}\end{minipage}}}%
    \put(0.46909927,0.82361498){\color[rgb]{0,0,0}\makebox(0,0)[lt]{\begin{minipage}{0.32445784\unitlength}\raggedright \gr{a}\end{minipage}}}%
    \put(0.74978098,0.49915726){\color[rgb]{0,0,0}\makebox(0,0)[lt]{\begin{minipage}{0.24463079\unitlength}\raggedright \blk{k}\end{minipage}}}%
  \end{picture}%
\endgroup%

%% file: arxiv-figures/cap_webs_cc_svg-tex.eps_tex
\begingroup%
  \makeatletter%
  \providecommand\color[2][]{%
    \errmessage{(Inkscape) Color is used for the text in Inkscape, but the package 'color.sty' is not loaded}%
    \renewcommand\color[2][]{}%
  }%
  \providecommand\transparent[1]{%
    \errmessage{(Inkscape) Transparency is used (non-zero) for the text in Inkscape, but the package 'transparent.sty' is not loaded}%
    \renewcommand\transparent[1]{}%
  }%
  \providecommand\rotatebox[2]{#2}%
  \newcommand*\fsize{\dimexpr\f@size pt\relax}%
  \newcommand*\lineheight[1]{\fontsize{\fsize}{#1\fsize}\selectfont}%
  \ifx\svgwidth\undefined%
    \setlength{\unitlength}{469.76640848bp}%
    \ifx\svgscale\undefined%
      \relax%
    \else%
      \setlength{\unitlength}{\unitlength * \real{\svgscale}}%
    \fi%
  \else%
    \setlength{\unitlength}{\svgwidth}%
  \fi%
  \global\let\svgwidth\undefined%
  \global\let\svgscale\undefined%
  \makeatother%
  \begin{picture}(1,0.92106259)%
    \lineheight{1}%
    \setlength\tabcolsep{0pt}%
    \put(0,0){\includegraphics[width=\unitlength]{cap_webs_cc_svg-tex.eps}}%
    \put(0.47167426,0.7618135){\color[rgb]{0,0,0}\makebox(0,0)[lt]{\lineheight{1.25}\smash{\begin{tabular}[t]{l}\gr{a}\end{tabular}}}}%
    \put(0.40214764,0.11032296){\color[rgb]{0,0,0}\makebox(0,0)[lt]{\lineheight{1.25}\smash{\begin{tabular}[t]{l}\gr{a+k}\end{tabular}}}}%
    \put(0.73948062,0.43993075){\color[rgb]{0,0,0}\makebox(0,0)[lt]{\lineheight{1.25}\smash{\begin{tabular}[t]{l}\blk{k}\end{tabular}}}}%
  \end{picture}%
\endgroup%

%% file: arxiv-figures/cup_webs_c_svg-tex.eps_tex
\begingroup%
  \makeatletter%
  \providecommand\color[2][]{%
    \errmessage{(Inkscape) Color is used for the text in Inkscape, but the package 'color.sty' is not loaded}%
    \renewcommand\color[2][]{}%
  }%
  \providecommand\transparent[1]{%
    \errmessage{(Inkscape) Transparency is used (non-zero) for the text in Inkscape, but the package 'transparent.sty' is not loaded}%
    \renewcommand\transparent[1]{}%
  }%
  \providecommand\rotatebox[2]{#2}%
  \newcommand*\fsize{\dimexpr\f@size pt\relax}%
  \newcommand*\lineheight[1]{\fontsize{\fsize}{#1\fsize}\selectfont}%
  \ifx\svgwidth\undefined%
    \setlength{\unitlength}{469.76640848bp}%
    \ifx\svgscale\undefined%
      \relax%
    \else%
      \setlength{\unitlength}{\unitlength * \real{\svgscale}}%
    \fi%
  \else%
    \setlength{\unitlength}{\svgwidth}%
  \fi%
  \global\let\svgwidth\undefined%
  \global\let\svgscale\undefined%
  \makeatother%
  \begin{picture}(1,0.92106259)%
    \lineheight{1}%
    \setlength\tabcolsep{0pt}%
    \put(0,0){\includegraphics[width=\unitlength]{cup_webs_c_svg-tex.eps}}%
    \put(0.48197448,0.22105048){\color[rgb]{0,0,0}\makebox(0,0)[lt]{\lineheight{1.25}\smash{\begin{tabular}[t]{l}\gr{a}\end{tabular}}}}%
    \put(0.40472267,0.741213){\color[rgb]{0,0,0}\makebox(0,0)[lt]{\lineheight{1.25}\smash{\begin{tabular}[t]{l}\gr{a-k}\end{tabular}}}}%
    \put(0.76523123,0.46310633){\color[rgb]{0,0,0}\makebox(0,0)[lt]{\lineheight{1.25}\smash{\begin{tabular}[t]{l}\blk{k}\end{tabular}}}}%
  \end{picture}%
\endgroup%

%% file: arxiv-figures/cup_webs_cc_svg-tex.eps_tex
\begingroup%
  \makeatletter%
  \providecommand\color[2][]{%
    \errmessage{(Inkscape) Color is used for the text in Inkscape, but the package 'color.sty' is not loaded}%
    \renewcommand\color[2][]{}%
  }%
  \providecommand\transparent[1]{%
    \errmessage{(Inkscape) Transparency is used (non-zero) for the text in Inkscape, but the package 'transparent.sty' is not loaded}%
    \renewcommand\transparent[1]{}%
  }%
  \providecommand\rotatebox[2]{#2}%
  \newcommand*\fsize{\dimexpr\f@size pt\relax}%
  \newcommand*\lineheight[1]{\fontsize{\fsize}{#1\fsize}\selectfont}%
  \ifx\svgwidth\undefined%
    \setlength{\unitlength}{469.76640848bp}%
    \ifx\svgscale\undefined%
      \relax%
    \else%
      \setlength{\unitlength}{\unitlength * \real{\svgscale}}%
    \fi%
  \else%
    \setlength{\unitlength}{\svgwidth}%
  \fi%
  \global\let\svgwidth\undefined%
  \global\let\svgscale\undefined%
  \makeatother%
  \begin{picture}(1,0.92106259)%
    \lineheight{1}%
    \setlength\tabcolsep{0pt}%
    \put(0,0){\includegraphics[width=\unitlength]{cup_webs_cc_svg-tex.eps}}%
    \put(0.48712469,0.21847554){\color[rgb]{0,0,0}\makebox(0,0)[lt]{\lineheight{1.25}\smash{\begin{tabular}[t]{l}\gr{a}\end{tabular}}}}%
    \put(0.3995726,0.78498907){\color[rgb]{0,0,0}\makebox(0,0)[lt]{\lineheight{1.25}\smash{\begin{tabular}[t]{l}\gr{a+k}\end{tabular}}}}%
    \put(0.79355696,0.47598168){\color[rgb]{0,0,0}\makebox(0,0)[lt]{\lineheight{1.25}\smash{\begin{tabular}[t]{l}\blk{k}\end{tabular}}}}%
  \end{picture}%
\endgroup%

%% file: arxiv-figures/tag_right_c_svg-tex.eps_tex
\begingroup%
  \makeatletter%
  \providecommand\color[2][]{%
    \errmessage{(Inkscape) Color is used for the text in Inkscape, but the package 'color.sty' is not loaded}%
    \renewcommand\color[2][]{}%
  }%
  \providecommand\transparent[1]{%
    \errmessage{(Inkscape) Transparency is used (non-zero) for the text in Inkscape, but the package 'transparent.sty' is not loaded}%
    \renewcommand\transparent[1]{}%
  }%
  \providecommand\rotatebox[2]{#2}%
  \newcommand*\fsize{\dimexpr\f@size pt\relax}%
  \newcommand*\lineheight[1]{\fontsize{\fsize}{#1\fsize}\selectfont}%
  \ifx\svgwidth\undefined%
    \setlength{\unitlength}{469.76640848bp}%
    \ifx\svgscale\undefined%
      \relax%
    \else%
      \setlength{\unitlength}{\unitlength * \real{\svgscale}}%
    \fi%
  \else%
    \setlength{\unitlength}{\svgwidth}%
  \fi%
  \global\let\svgwidth\undefined%
  \global\let\svgscale\undefined%
  \makeatother%
  \begin{picture}(1,0.92106259)%
    \lineheight{1}%
    \setlength\tabcolsep{0pt}%
    \put(0,0){\includegraphics[width=\unitlength]{tag_right_c_svg-tex.eps}}%
    \put(3.03676008,-1.61145364){\color[rgb]{0,0,0}\makebox(0,0)[lt]{\begin{minipage}{0.47689347\unitlength}\end{minipage}}}%
    \put(0.66480394,0.54550833){\color[rgb]{0,0,0}\makebox(0,0)[lt]{\begin{minipage}{0.16995404\unitlength}\end{minipage}}}%
    \put(2.92970235,-1.27081548){\color[rgb]{0,0,0}\makebox(0,0)[lt]{\begin{minipage}{0.70074114\unitlength}\end{minipage}}}%
    \put(0.70911357,0.49100215){\color[rgb]{0,0,0}\makebox(0,0)[lt]{\begin{minipage}{0.12875304\unitlength}\raggedright \gr{a}\end{minipage}}}%
    \put(0.1660838,0.47037437){\color[rgb]{0,0,0}\makebox(0,0)[lt]{\lineheight{1.25}\smash{\begin{tabular}[t]{l}\gr{a+k}\end{tabular}}}}%
    \put(0.53448945,0.3239251){\color[rgb]{0,0,0}\makebox(0,0)[lt]{\begin{minipage}{0.09010165\unitlength}\raggedright \blk{k}\end{minipage}}}%
    \put(0.51769813,0.6887245){\color[rgb]{0,0,0}\makebox(0,0)[lt]{\begin{minipage}{0.41508106\unitlength}\raggedright \blk{n-k}                           \end{minipage}}}%
    \put(0.51522089,0.32254877){\color[rgb]{0,0,0}\makebox(0,0)[lt]{\begin{minipage}{0.09010165\unitlength}\end{minipage}}}%
  \end{picture}%
\endgroup%

%% file: arxiv-figures/tag_right_d_svg-tex.eps_tex
\begingroup%
  \makeatletter%
  \providecommand\color[2][]{%
    \errmessage{(Inkscape) Color is used for the text in Inkscape, but the package 'color.sty' is not loaded}%
    \renewcommand\color[2][]{}%
  }%
  \providecommand\transparent[1]{%
    \errmessage{(Inkscape) Transparency is used (non-zero) for the text in Inkscape, but the package 'transparent.sty' is not loaded}%
    \renewcommand\transparent[1]{}%
  }%
  \providecommand\rotatebox[2]{#2}%
  \newcommand*\fsize{\dimexpr\f@size pt\relax}%
  \newcommand*\lineheight[1]{\fontsize{\fsize}{#1\fsize}\selectfont}%
  \ifx\svgwidth\undefined%
    \setlength{\unitlength}{469.76640848bp}%
    \ifx\svgscale\undefined%
      \relax%
    \else%
      \setlength{\unitlength}{\unitlength * \real{\svgscale}}%
    \fi%
  \else%
    \setlength{\unitlength}{\svgwidth}%
  \fi%
  \global\let\svgwidth\undefined%
  \global\let\svgscale\undefined%
  \makeatother%
  \begin{picture}(1,0.92106259)%
    \lineheight{1}%
    \setlength\tabcolsep{0pt}%
    \put(0,0){\includegraphics[width=\unitlength]{tag_right_d_svg-tex.eps}}%
    \put(0.7034299,0.45538123){\color[rgb]{0,0,0}\makebox(0,0)[lt]{\lineheight{11}\smash{\begin{tabular}[t]{l}\gr{a}\end{tabular}}}}%
    \put(0.16524195,0.45795626){\color[rgb]{0,0,0}\makebox(0,0)[lt]{\lineheight{11}\smash{\begin{tabular}[t]{l}\gr{a-k}\end{tabular}}}}%
    \put(0.54120103,0.19272484){\color[rgb]{0,0,0}\makebox(0,0)[lt]{\lineheight{11}\smash{\begin{tabular}[t]{l}\blk{k}\end{tabular}}}}%
    \put(0.54120099,0.70258714){\color[rgb]{0,0,0}\makebox(0,0)[lt]{\lineheight{11}\smash{\begin{tabular}[t]{l}\blk{n-k}\end{tabular}}}}%
  \end{picture}%
\endgroup%

%% file: arxiv-figures/trivalent_vertex_up_svg-tex.eps_tex
\begingroup%
  \makeatletter%
  \providecommand\color[2][]{%
    \errmessage{(Inkscape) Color is used for the text in Inkscape, but the package 'color.sty' is not loaded}%
    \renewcommand\color[2][]{}%
  }%
  \providecommand\transparent[1]{%
    \errmessage{(Inkscape) Transparency is used (non-zero) for the text in Inkscape, but the package 'transparent.sty' is not loaded}%
    \renewcommand\transparent[1]{}%
  }%
  \providecommand\rotatebox[2]{#2}%
  \newcommand*\fsize{\dimexpr\f@size pt\relax}%
  \newcommand*\lineheight[1]{\fontsize{\fsize}{#1\fsize}\selectfont}%
  \ifx\svgwidth\undefined%
    \setlength{\unitlength}{469.76640848bp}%
    \ifx\svgscale\undefined%
      \relax%
    \else%
      \setlength{\unitlength}{\unitlength * \real{\svgscale}}%
    \fi%
  \else%
    \setlength{\unitlength}{\svgwidth}%
  \fi%
  \global\let\svgwidth\undefined%
  \global\let\svgscale\undefined%
  \makeatother%
  \begin{picture}(1,0.92106259)%
    \lineheight{1}%
    \setlength\tabcolsep{0pt}%
    \put(0,0){\includegraphics[width=\unitlength]{trivalent_vertex_up_svg-tex.eps}}%
    \put(0.70803369,0.49197824){\color[rgb]{0,0,0}\makebox(0,0)[lt]{\lineheight{1.25}\smash{\begin{tabular}[t]{l}\gr{a}\end{tabular}}}}%
    \put(0.39840212,0.1363856){\color[rgb]{0,0,0}\makebox(0,0)[lt]{\lineheight{1.25}\smash{\begin{tabular}[t]{l}\gr{a+l}\end{tabular}}}}%
    \put(0.04522852,0.49197824){\color[rgb]{0,0,0}\makebox(0,0)[lt]{\lineheight{1.25}\smash{\begin{tabular}[t]{l}\gr{a+k+l}\end{tabular}}}}%
    \put(0.65391049,0.26421107){\color[rgb]{0,0,0}\makebox(0,0)[lt]{\lineheight{1.25}\smash{\begin{tabular}[t]{l}\blk{l}\end{tabular}}}}%
    \put(0.21455832,0.24282155){\color[rgb]{0,0,0}\makebox(0,0)[lt]{\lineheight{1.25}\smash{\begin{tabular}[t]{l}\blk{k}\end{tabular}}}}%
    \put(0.51451399,0.71694498){\color[rgb]{0,0,0}\makebox(0,0)[lt]{\lineheight{1.25}\smash{\begin{tabular}[t]{l}\blk{k+l}\end{tabular}}}}%
  \end{picture}%
\endgroup%

%% file: arxiv-figures/trivalent_vertex_down_svg-tex.eps_tex
\begingroup%
  \makeatletter%
  \providecommand\color[2][]{%
    \errmessage{(Inkscape) Color is used for the text in Inkscape, but the package 'color.sty' is not loaded}%
    \renewcommand\color[2][]{}%
  }%
  \providecommand\transparent[1]{%
    \errmessage{(Inkscape) Transparency is used (non-zero) for the text in Inkscape, but the package 'transparent.sty' is not loaded}%
    \renewcommand\transparent[1]{}%
  }%
  \providecommand\rotatebox[2]{#2}%
  \newcommand*\fsize{\dimexpr\f@size pt\relax}%
  \newcommand*\lineheight[1]{\fontsize{\fsize}{#1\fsize}\selectfont}%
  \ifx\svgwidth\undefined%
    \setlength{\unitlength}{469.76640848bp}%
    \ifx\svgscale\undefined%
      \relax%
    \else%
      \setlength{\unitlength}{\unitlength * \real{\svgscale}}%
    \fi%
  \else%
    \setlength{\unitlength}{\svgwidth}%
  \fi%
  \global\let\svgwidth\undefined%
  \global\let\svgscale\undefined%
  \makeatother%
  \begin{picture}(1,0.92106259)%
    \lineheight{1}%
    \setlength\tabcolsep{0pt}%
    \put(0,0){\includegraphics[width=\unitlength]{trivalent_vertex_down_svg-tex.eps}}%
    \put(0.74431863,0.39653035){\color[rgb]{0,0,0}\makebox(0,0)[lt]{\lineheight{1.25}\smash{\begin{tabular}[t]{l}\gr{a}\end{tabular}}}}%
    \put(0.0597425,0.39169242){\color[rgb]{0,0,0}\makebox(0,0)[lt]{\lineheight{1.25}\smash{\begin{tabular}[t]{l}\gr{a+k+l}\end{tabular}}}}%
    \put(0.40565913,0.78356984){\color[rgb]{0,0,0}\makebox(0,0)[lt]{\lineheight{1.25}\smash{\begin{tabular}[t]{l}\gr{a+l}\end{tabular}}}}%
    \put(0.70561478,0.65294408){\color[rgb]{0,0,0}\makebox(0,0)[lt]{\lineheight{1.25}\smash{\begin{tabular}[t]{l}\blk{l}\end{tabular}}}}%
    \put(0.16134038,0.67471497){\color[rgb]{0,0,0}\makebox(0,0)[lt]{\lineheight{1.25}\smash{\begin{tabular}[t]{l}\blk{k}\end{tabular}}}}%
    \put(0.53870391,0.23687657){\color[rgb]{0,0,0}\makebox(0,0)[lt]{\lineheight{1.25}\smash{\begin{tabular}[t]{l}\blk{k+l}\end{tabular}}}}%
  \end{picture}%
\endgroup%

%% file: arxiv-figures/isotopy_svg-tex.eps_tex
\begingroup%
  \makeatletter%
  \providecommand\color[2][]{%
    \errmessage{(Inkscape) Color is used for the text in Inkscape, but the package 'color.sty' is not loaded}%
    \renewcommand\color[2][]{}%
  }%
  \providecommand\transparent[1]{%
    \errmessage{(Inkscape) Transparency is used (non-zero) for the text in Inkscape, but the package 'transparent.sty' is not loaded}%
    \renewcommand\transparent[1]{}%
  }%
  \providecommand\rotatebox[2]{#2}%
  \newcommand*\fsize{\dimexpr\f@size pt\relax}%
  \newcommand*\lineheight[1]{\fontsize{\fsize}{#1\fsize}\selectfont}%
  \ifx\svgwidth\undefined%
    \setlength{\unitlength}{316.85077963bp}%
    \ifx\svgscale\undefined%
      \relax%
    \else%
      \setlength{\unitlength}{\unitlength * \real{\svgscale}}%
    \fi%
  \else%
    \setlength{\unitlength}{\svgwidth}%
  \fi%
  \global\let\svgwidth\undefined%
  \global\let\svgscale\undefined%
  \makeatother%
  \begin{picture}(1,0.1644351)%
    \lineheight{1}%
    \setlength\tabcolsep{0pt}%
    \put(0,0){\includegraphics[width=\unitlength]{isotopy_svg-tex.eps}}%
    \put(0.15876021,0.09226161){\color[rgb]{0,0,0}\makebox(0,0)[lt]{\begin{minipage}{0.02799467\unitlength}\raggedright $=$\end{minipage}}}%
    \put(0.2739154,0.09226161){\color[rgb]{0,0,0}\makebox(0,0)[lt]{\begin{minipage}{0.02722614\unitlength}\raggedright $=$\end{minipage}}}%
    \put(0.70313956,0.09072457){\color[rgb]{0,0,0}\makebox(0,0)[lt]{\begin{minipage}{0.03106876\unitlength}\raggedright $=$\end{minipage}}}%
    \put(0.81445229,0.09072456){\color[rgb]{0,0,0}\makebox(0,0)[lt]{\begin{minipage}{0.03030017\unitlength}\raggedright $=$\end{minipage}}}%
  \end{picture}%
\endgroup%

%% file: arxiv-figures/crosscyclic_svg-tex.eps_tex
\begingroup%
  \makeatletter%
  \providecommand\color[2][]{%
    \errmessage{(Inkscape) Color is used for the text in Inkscape, but the package 'color.sty' is not loaded}%
    \renewcommand\color[2][]{}%
  }%
  \providecommand\transparent[1]{%
    \errmessage{(Inkscape) Transparency is used (non-zero) for the text in Inkscape, but the package 'transparent.sty' is not loaded}%
    \renewcommand\transparent[1]{}%
  }%
  \providecommand\rotatebox[2]{#2}%
  \newcommand*\fsize{\dimexpr\f@size pt\relax}%
  \newcommand*\lineheight[1]{\fontsize{\fsize}{#1\fsize}\selectfont}%
  \ifx\svgwidth\undefined%
    \setlength{\unitlength}{201.49999237bp}%
    \ifx\svgscale\undefined%
      \relax%
    \else%
      \setlength{\unitlength}{\unitlength * \real{\svgscale}}%
    \fi%
  \else%
    \setlength{\unitlength}{\svgwidth}%
  \fi%
  \global\let\svgwidth\undefined%
  \global\let\svgscale\undefined%
  \makeatother%
  \begin{picture}(1,0.325062)%
    \lineheight{1}%
    \setlength\tabcolsep{0pt}%
    \put(0,0){\includegraphics[width=\unitlength]{crosscyclic_svg-tex.eps}}%
    \put(0.4801489,0.16253113){\color[rgb]{0,0,0}\makebox(0,0)[lt]{\lineheight{0}\smash{\begin{tabular}[t]{l}=\end{tabular}}}}%
  \end{picture}%
\endgroup%

%% file: arxiv-figures/crosscyclic2_svg-tex.eps_tex
\begingroup%
  \makeatletter%
  \providecommand\color[2][]{%
    \errmessage{(Inkscape) Color is used for the text in Inkscape, but the package 'color.sty' is not loaded}%
    \renewcommand\color[2][]{}%
  }%
  \providecommand\transparent[1]{%
    \errmessage{(Inkscape) Transparency is used (non-zero) for the text in Inkscape, but the package 'transparent.sty' is not loaded}%
    \renewcommand\transparent[1]{}%
  }%
  \providecommand\rotatebox[2]{#2}%
  \newcommand*\fsize{\dimexpr\f@size pt\relax}%
  \newcommand*\lineheight[1]{\fontsize{\fsize}{#1\fsize}\selectfont}%
  \ifx\svgwidth\undefined%
    \setlength{\unitlength}{189.0053227bp}%
    \ifx\svgscale\undefined%
      \relax%
    \else%
      \setlength{\unitlength}{\unitlength * \real{\svgscale}}%
    \fi%
  \else%
    \setlength{\unitlength}{\svgwidth}%
  \fi%
  \global\let\svgwidth\undefined%
  \global\let\svgscale\undefined%
  \makeatother%
  \begin{picture}(1,0.299902)%
    \lineheight{1}%
    \setlength\tabcolsep{0pt}%
    \put(0,0){\includegraphics[width=\unitlength]{crosscyclic2_svg-tex.eps}}%
    \put(0.48996091,0.17080448){\color[rgb]{0,0,0}\makebox(0,0)[lt]{\begin{minipage}{0.04232681\unitlength}\raggedright =\end{minipage}}}%
  \end{picture}%
\endgroup%

%% file: arxiv-figures/polymult_svg-tex.eps_tex
\begingroup%
  \makeatletter%
  \providecommand\color[2][]{%
    \errmessage{(Inkscape) Color is used for the text in Inkscape, but the package 'color.sty' is not loaded}%
    \renewcommand\color[2][]{}%
  }%
  \providecommand\transparent[1]{%
    \errmessage{(Inkscape) Transparency is used (non-zero) for the text in Inkscape, but the package 'transparent.sty' is not loaded}%
    \renewcommand\transparent[1]{}%
  }%
  \providecommand\rotatebox[2]{#2}%
  \newcommand*\fsize{\dimexpr\f@size pt\relax}%
  \newcommand*\lineheight[1]{\fontsize{\fsize}{#1\fsize}\selectfont}%
  \ifx\svgwidth\undefined%
    \setlength{\unitlength}{112.80001339bp}%
    \ifx\svgscale\undefined%
      \relax%
    \else%
      \setlength{\unitlength}{\unitlength * \real{\svgscale}}%
    \fi%
  \else%
    \setlength{\unitlength}{\svgwidth}%
  \fi%
  \global\let\svgwidth\undefined%
  \global\let\svgscale\undefined%
  \makeatother%
  \begin{picture}(1,0.46808507)%
    \lineheight{1}%
    \setlength\tabcolsep{0pt}%
    \put(0,0){\includegraphics[width=\unitlength]{polymult_svg-tex.eps}}%
    \put(0.53546108,0.21631178){\color[rgb]{0,0,0}\rotatebox{-180}{\makebox(0,0)[lt]{\lineheight{0}\smash{\begin{tabular}[t]{l}=\end{tabular}}}}}%
    \put(0.10992907,0.21631178){\color[rgb]{0,0,0}\makebox(0,0)[lt]{\lineheight{0}\smash{\begin{tabular}[t]{l}$f$\end{tabular}}}}%
    \put(0.2801419,0.21631178){\color[rgb]{0,0,0}\makebox(0,0)[lt]{\lineheight{0}\smash{\begin{tabular}[t]{l}$g$\end{tabular}}}}%
    \put(0.73404263,0.21631178){\color[rgb]{0,0,0}\makebox(0,0)[lt]{\lineheight{0}\smash{\begin{tabular}[t]{l}$fg$\end{tabular}}}}%
  \end{picture}%
\endgroup%

%% file: arxiv-figures/cap_ssbim_svg-tex.eps_tex
\begingroup%
  \makeatletter%
  \providecommand\color[2][]{%
    \errmessage{(Inkscape) Color is used for the text in Inkscape, but the package 'color.sty' is not loaded}%
    \renewcommand\color[2][]{}%
  }%
  \providecommand\transparent[1]{%
    \errmessage{(Inkscape) Transparency is used (non-zero) for the text in Inkscape, but the package 'transparent.sty' is not loaded}%
    \renewcommand\transparent[1]{}%
  }%
  \providecommand\rotatebox[2]{#2}%
  \newcommand*\fsize{\dimexpr\f@size pt\relax}%
  \newcommand*\lineheight[1]{\fontsize{\fsize}{#1\fsize}\selectfont}%
  \ifx\svgwidth\undefined%
    \setlength{\unitlength}{469.76640848bp}%
    \ifx\svgscale\undefined%
      \relax%
    \else%
      \setlength{\unitlength}{\unitlength * \real{\svgscale}}%
    \fi%
  \else%
    \setlength{\unitlength}{\svgwidth}%
  \fi%
  \global\let\svgwidth\undefined%
  \global\let\svgscale\undefined%
  \makeatother%
  \begin{picture}(1,0.92106259)%
    \lineheight{1}%
    \setlength\tabcolsep{0pt}%
    \put(0,0){\includegraphics[width=\unitlength]{cap_ssbim_svg-tex.eps}}%
    \put(0.47167429,0.85709077){\color[rgb]{0,0,0}\makebox(0,0)[lt]{\begin{minipage}{0.3502084\unitlength}\raggedright \dgr{a}\end{minipage}}}%
    \put(0.39442245,0.15667402){\color[rgb]{0,0,0}\makebox(0,0)[lt]{\begin{minipage}{0.51501237\unitlength}\raggedright \dgr{a-k}\end{minipage}}}%
    \put(0.65450371,0.22877576){\color[rgb]{0,0,0}\makebox(0,0)[lt]{\lineheight{1.25}\smash{\begin{tabular}[t]{l}\blk{a}\end{tabular}}}}%
    \put(0.73175552,0.45023107){\color[rgb]{0,0,0}\makebox(0,0)[lt]{\lineheight{1.25}\smash{\begin{tabular}[t]{l}\blk{a-k}\end{tabular}}}}%
    \put(0.33777112,0.43478071){\color[rgb]{0,0,0}\makebox(0,0)[lt]{\begin{minipage}{0.56136353\unitlength}\raggedright \dgr{a,a-k}\end{minipage}}}%
  \end{picture}%
\endgroup%

%% file: arxiv-figures/cup_ssbim_svg-tex.eps_tex
\begingroup%
  \makeatletter%
  \providecommand\color[2][]{%
    \errmessage{(Inkscape) Color is used for the text in Inkscape, but the package 'color.sty' is not loaded}%
    \renewcommand\color[2][]{}%
  }%
  \providecommand\transparent[1]{%
    \errmessage{(Inkscape) Transparency is used (non-zero) for the text in Inkscape, but the package 'transparent.sty' is not loaded}%
    \renewcommand\transparent[1]{}%
  }%
  \providecommand\rotatebox[2]{#2}%
  \newcommand*\fsize{\dimexpr\f@size pt\relax}%
  \newcommand*\lineheight[1]{\fontsize{\fsize}{#1\fsize}\selectfont}%
  \ifx\svgwidth\undefined%
    \setlength{\unitlength}{469.76640848bp}%
    \ifx\svgscale\undefined%
      \relax%
    \else%
      \setlength{\unitlength}{\unitlength * \real{\svgscale}}%
    \fi%
  \else%
    \setlength{\unitlength}{\svgwidth}%
  \fi%
  \global\let\svgwidth\undefined%
  \global\let\svgscale\undefined%
  \makeatother%
  \begin{picture}(1,0.92106259)%
    \lineheight{1}%
    \setlength\tabcolsep{0pt}%
    \put(0,0){\includegraphics[width=\unitlength]{cup_ssbim_svg-tex.eps}}%
    \put(0.40987283,0.26225155){\color[rgb]{0,0,0}\makebox(0,0)[lt]{\lineheight{1.25}\smash{\begin{tabular}[t]{l}\dgr{a+k}\end{tabular}}}}%
    \put(0.44077354,0.82876513){\color[rgb]{0,0,0}\makebox(0,0)[lt]{\lineheight{1.25}\smash{\begin{tabular}[t]{l}\dgr{a}\end{tabular}}}}%
    \put(0.40214764,0.54035823){\color[rgb]{0,0,0}\makebox(0,0)[lt]{\lineheight{1.25}\smash{\begin{tabular}[t]{l}\dgr{a,a+k}\end{tabular}}}}%
    \put(0.82703284,0.54550839){\color[rgb]{0,0,0}\makebox(0,0)[lt]{\lineheight{1.25}\smash{\begin{tabular}[t]{l}\blk{a}\end{tabular}}}}%
    \put(0.69312957,0.71803743){\color[rgb]{0,0,0}\makebox(0,0)[lt]{\lineheight{1.25}\smash{\begin{tabular}[t]{l}\blk{a+k}\end{tabular}}}}%
  \end{picture}%
\endgroup%

%% file: arxiv-figures/eg_mssbim_svg-tex.eps_tex
\begingroup%
  \makeatletter%
  \providecommand\color[2][]{%
    \errmessage{(Inkscape) Color is used for the text in Inkscape, but the package 'color.sty' is not loaded}%
    \renewcommand\color[2][]{}%
  }%
  \providecommand\transparent[1]{%
    \errmessage{(Inkscape) Transparency is used (non-zero) for the text in Inkscape, but the package 'transparent.sty' is not loaded}%
    \renewcommand\transparent[1]{}%
  }%
  \providecommand\rotatebox[2]{#2}%
  \newcommand*\fsize{\dimexpr\f@size pt\relax}%
  \newcommand*\lineheight[1]{\fontsize{\fsize}{#1\fsize}\selectfont}%
  \ifx\svgwidth\undefined%
    \setlength{\unitlength}{510.98442416bp}%
    \ifx\svgscale\undefined%
      \relax%
    \else%
      \setlength{\unitlength}{\unitlength * \real{\svgscale}}%
    \fi%
  \else%
    \setlength{\unitlength}{\svgwidth}%
  \fi%
  \global\let\svgwidth\undefined%
  \global\let\svgscale\undefined%
  \makeatother%
  \begin{picture}(1,0.84831313)%
    \lineheight{1}%
    \setlength\tabcolsep{0pt}%
    \put(0,0){\includegraphics[width=\unitlength]{eg_mssbim_svg-tex.eps}}%
    \put(0.82294175,0.57317524){\color[rgb]{0,0,0}\makebox(0,0)[lt]{\lineheight{1.25}\smash{\begin{tabular}[t]{l}\dgr{a}\end{tabular}}}}%
    \put(0.71925198,0.19866105){\color[rgb]{0,0,0}\makebox(0,0)[lt]{\lineheight{1.25}\smash{\begin{tabular}[t]{l}\blk{a+k}\end{tabular}}}}%
    \put(0.03272138,0.56796721){\color[rgb]{0,0,0}\makebox(0,0)[lt]{\lineheight{1.25}\smash{\begin{tabular}[t]{l}\dgr{a+k}\end{tabular}}}}%
    \put(0.20601122,0.19297938){\color[rgb]{0,0,0}\makebox(0,0)[lt]{\lineheight{1.25}\smash{\begin{tabular}[t]{l}\blk{a}\end{tabular}}}}%
    \put(0.35420712,0.33644077){\color[rgb]{0,0,0}\makebox(0,0)[lt]{\lineheight{1.25}\smash{\begin{tabular}[t]{l}\dgr{a,a+k}\end{tabular}}}}%
  \end{picture}%
\endgroup%

%% file: arxiv-figures/image_of_trivalent_top_section_svg-tex.eps_tex
\begingroup%
  \makeatletter%
  \providecommand\color[2][]{%
    \errmessage{(Inkscape) Color is used for the text in Inkscape, but the package 'color.sty' is not loaded}%
    \renewcommand\color[2][]{}%
  }%
  \providecommand\transparent[1]{%
    \errmessage{(Inkscape) Transparency is used (non-zero) for the text in Inkscape, but the package 'transparent.sty' is not loaded}%
    \renewcommand\transparent[1]{}%
  }%
  \providecommand\rotatebox[2]{#2}%
  \newcommand*\fsize{\dimexpr\f@size pt\relax}%
  \newcommand*\lineheight[1]{\fontsize{\fsize}{#1\fsize}\selectfont}%
  \ifx\svgwidth\undefined%
    \setlength{\unitlength}{1054.81886894bp}%
    \ifx\svgscale\undefined%
      \relax%
    \else%
      \setlength{\unitlength}{\unitlength * \real{\svgscale}}%
    \fi%
  \else%
    \setlength{\unitlength}{\svgwidth}%
  \fi%
  \global\let\svgwidth\undefined%
  \global\let\svgscale\undefined%
  \makeatother%
  \begin{picture}(1,0.40922131)%
    \lineheight{1}%
    \setlength\tabcolsep{0pt}%
    \put(0,0){\includegraphics[width=\unitlength]{image_of_trivalent_top_section_svg-tex.eps}}%
    \put(0.92546673,0.1208334){\color[rgb]{0,0,0}\makebox(0,0)[lt]{\lineheight{1.25}\smash{\begin{tabular}[t]{l}\dgr{a}\end{tabular}}}}%
    \put(0.52064265,0.34790192){\color[rgb]{0,0,0}\makebox(0,0)[lt]{\lineheight{1.25}\smash{\begin{tabular}[t]{l}\dgr{a+l}\end{tabular}}}}%
    \put(0.6903706,0.19766973){\color[rgb]{0,0,0}\makebox(0,0)[lt]{\lineheight{1.25}\smash{\begin{tabular}[t]{l}\dgr{a,a+l}\end{tabular}}}}%
    \put(0.21444421,0.1953761){\color[rgb]{0,0,0}\makebox(0,0)[lt]{\lineheight{1.25}\smash{\begin{tabular}[t]{l}\dgr{a+l,a+k+l}\end{tabular}}}}%
    \put(0.41903522,0.04514388){\color[rgb]{0,0,0}\makebox(0,0)[lt]{\lineheight{1.25}\smash{\begin{tabular}[t]{l}\dgr{a,a+l,a+k+l}\end{tabular}}}}%
    \put(0.01719279,0.12427392){\color[rgb]{0,0,0}\makebox(0,0)[lt]{\lineheight{1.25}\smash{\begin{tabular}[t]{l}\dgr{a+k+l}\end{tabular}}}}%
  \end{picture}%
\endgroup%

%% file: arxiv-figures/trivalent_split_svg-tex.eps_tex
\begingroup%
  \makeatletter%
  \providecommand\color[2][]{%
    \errmessage{(Inkscape) Color is used for the text in Inkscape, but the package 'color.sty' is not loaded}%
    \renewcommand\color[2][]{}%
  }%
  \providecommand\transparent[1]{%
    \errmessage{(Inkscape) Transparency is used (non-zero) for the text in Inkscape, but the package 'transparent.sty' is not loaded}%
    \renewcommand\transparent[1]{}%
  }%
  \providecommand\rotatebox[2]{#2}%
  \newcommand*\fsize{\dimexpr\f@size pt\relax}%
  \newcommand*\lineheight[1]{\fontsize{\fsize}{#1\fsize}\selectfont}%
  \ifx\svgwidth\undefined%
    \setlength{\unitlength}{919.06029191bp}%
    \ifx\svgscale\undefined%
      \relax%
    \else%
      \setlength{\unitlength}{\unitlength * \real{\svgscale}}%
    \fi%
  \else%
    \setlength{\unitlength}{\svgwidth}%
  \fi%
  \global\let\svgwidth\undefined%
  \global\let\svgscale\undefined%
  \makeatother%
  \begin{picture}(1,0.93429174)%
    \lineheight{1}%
    \setlength\tabcolsep{0pt}%
    \put(0,0){\includegraphics[width=\unitlength]{trivalent_split_svg-tex.eps}}%
    \put(0.85086438,0.35017766){\color[rgb]{0,0,0}\makebox(0,0)[lt]{\lineheight{1.25}\smash{\begin{tabular}[t]{l}\dgr{a}\end{tabular}}}}%
    \put(0.45073627,0.86481626){\color[rgb]{0,0,0}\makebox(0,0)[lt]{\lineheight{1.25}\smash{\begin{tabular}[t]{l}\dgr{a+l}\end{tabular}}}}%
    \put(0.02165141,0.35544259){\color[rgb]{0,0,0}\makebox(0,0)[lt]{\lineheight{1.25}\smash{\begin{tabular}[t]{l}\dgr{a+k+l}\end{tabular}}}}%
  \end{picture}%
\endgroup%

%% file: arxiv-figures/image_of_down_trivalent_svg-tex.eps_tex
\begingroup%
  \makeatletter%
  \providecommand\color[2][]{%
    \errmessage{(Inkscape) Color is used for the text in Inkscape, but the package 'color.sty' is not loaded}%
    \renewcommand\color[2][]{}%
  }%
  \providecommand\transparent[1]{%
    \errmessage{(Inkscape) Transparency is used (non-zero) for the text in Inkscape, but the package 'transparent.sty' is not loaded}%
    \renewcommand\transparent[1]{}%
  }%
  \providecommand\rotatebox[2]{#2}%
  \newcommand*\fsize{\dimexpr\f@size pt\relax}%
  \newcommand*\lineheight[1]{\fontsize{\fsize}{#1\fsize}\selectfont}%
  \ifx\svgwidth\undefined%
    \setlength{\unitlength}{469.83444011bp}%
    \ifx\svgscale\undefined%
      \relax%
    \else%
      \setlength{\unitlength}{\unitlength * \real{\svgscale}}%
    \fi%
  \else%
    \setlength{\unitlength}{\svgwidth}%
  \fi%
  \global\let\svgwidth\undefined%
  \global\let\svgscale\undefined%
  \makeatother%
  \begin{picture}(1,0.92107402)%
    \lineheight{1}%
    \setlength\tabcolsep{0pt}%
    \put(0,0){\includegraphics[width=\unitlength]{image_of_down_trivalent_svg-tex.eps}}%
    \put(0.78104139,0.33567976){\color[rgb]{0,0,0}\makebox(0,0)[lt]{\lineheight{1.25}\smash{\begin{tabular}[t]{l}\dgr{a}\end{tabular}}}}%
    \put(0.02236649,0.33310501){\color[rgb]{0,0,0}\makebox(0,0)[lt]{\lineheight{1.25}\smash{\begin{tabular}[t]{l}\dgr{a+k+l}\end{tabular}}}}%
    \put(0.39996336,0.79654896){\color[rgb]{0,0,0}\makebox(0,0)[lt]{\lineheight{1.25}\smash{\begin{tabular}[t]{l}\dgr{a+l}\end{tabular}}}}%
  \end{picture}%
\endgroup%

%% file: arxiv-figures/polyslide_svg-tex.eps_tex
\begingroup%
  \makeatletter%
  \providecommand\color[2][]{%
    \errmessage{(Inkscape) Color is used for the text in Inkscape, but the package 'color.sty' is not loaded}%
    \renewcommand\color[2][]{}%
  }%
  \providecommand\transparent[1]{%
    \errmessage{(Inkscape) Transparency is used (non-zero) for the text in Inkscape, but the package 'transparent.sty' is not loaded}%
    \renewcommand\transparent[1]{}%
  }%
  \providecommand\rotatebox[2]{#2}%
  \newcommand*\fsize{\dimexpr\f@size pt\relax}%
  \newcommand*\lineheight[1]{\fontsize{\fsize}{#1\fsize}\selectfont}%
  \ifx\svgwidth\undefined%
    \setlength{\unitlength}{112.80000195bp}%
    \ifx\svgscale\undefined%
      \relax%
    \else%
      \setlength{\unitlength}{\unitlength * \real{\svgscale}}%
    \fi%
  \else%
    \setlength{\unitlength}{\svgwidth}%
  \fi%
  \global\let\svgwidth\undefined%
  \global\let\svgscale\undefined%
  \makeatother%
  \begin{picture}(1,0.46808512)%
    \lineheight{1}%
    \setlength\tabcolsep{0pt}%
    \put(0,0){\includegraphics[width=\unitlength]{polyslide_svg-tex.eps}}%
    \put(0.53546103,0.2163118){\color[rgb]{0,0,0}\rotatebox{-180}{\makebox(0,0)[lt]{\lineheight{0}\smash{\begin{tabular}[t]{l}=\end{tabular}}}}}%
    \put(0.25929395,0.11409821){\color[rgb]{0,0,0}\makebox(0,0)[lt]{\lineheight{0}\smash{\begin{tabular}[t]{l}\gr{Ii}\end{tabular}}}}%
    \put(0.08444401,0.11409821){\color[rgb]{0,0,0}\makebox(0,0)[lt]{\lineheight{0}\smash{\begin{tabular}[t]{l}\gr{I}\end{tabular}}}}%
    \put(0.31560289,0.28014183){\color[rgb]{0,0,0}\makebox(0,0)[lt]{\lineheight{0}\smash{\begin{tabular}[t]{l}$f$\end{tabular}}}}%
    \put(0.92260938,0.10992882){\color[rgb]{0,0,0}\makebox(0,0)[lt]{\lineheight{0}\smash{\begin{tabular}[t]{l}\gr{Ii}\end{tabular}}}}%
    \put(0.6857806,0.10992882){\color[rgb]{0,0,0}\makebox(0,0)[lt]{\lineheight{0}\smash{\begin{tabular}[t]{l}\gr{I}\end{tabular}}}}%
    \put(0.64982271,0.28206906){\color[rgb]{0,0,0}\makebox(0,0)[lt]{\lineheight{0}\smash{\begin{tabular}[t]{l}$f$\end{tabular}}}}%
  \end{picture}%
\endgroup%

%% file: arxiv-figures/cccirc_svg-tex.eps_tex
\begingroup%
  \makeatletter%
  \providecommand\color[2][]{%
    \errmessage{(Inkscape) Color is used for the text in Inkscape, but the package 'color.sty' is not loaded}%
    \renewcommand\color[2][]{}%
  }%
  \providecommand\transparent[1]{%
    \errmessage{(Inkscape) Transparency is used (non-zero) for the text in Inkscape, but the package 'transparent.sty' is not loaded}%
    \renewcommand\transparent[1]{}%
  }%
  \providecommand\rotatebox[2]{#2}%
  \newcommand*\fsize{\dimexpr\f@size pt\relax}%
  \newcommand*\lineheight[1]{\fontsize{\fsize}{#1\fsize}\selectfont}%
  \ifx\svgwidth\undefined%
    \setlength{\unitlength}{113.20706089bp}%
    \ifx\svgscale\undefined%
      \relax%
    \else%
      \setlength{\unitlength}{\unitlength * \real{\svgscale}}%
    \fi%
  \else%
    \setlength{\unitlength}{\svgwidth}%
  \fi%
  \global\let\svgwidth\undefined%
  \global\let\svgscale\undefined%
  \makeatother%
  \begin{picture}(1,0.43106852)%
    \lineheight{1}%
    \setlength\tabcolsep{0pt}%
    \put(0,0){\includegraphics[width=\unitlength]{cccirc_svg-tex.eps}}%
    \put(0.47091386,0.20714137){\color[rgb]{0,0,0}\makebox(0,0)[lt]{\lineheight{0}\smash{\begin{tabular}[t]{l}=\end{tabular}}}}%
    \put(0.28763179,0.04559565){\color[rgb]{0,0,0}\makebox(0,0)[lt]{\lineheight{0}\smash{\begin{tabular}[t]{l}\gr{I}\end{tabular}}}}%
    \put(0.18020071,0.18834673){\color[rgb]{0,0,0}\makebox(0,0)[lt]{\lineheight{0}\smash{\begin{tabular}[t]{l}\gr{Ii}\end{tabular}}}}%
    \put(0.72433657,0.18890354){\color[rgb]{0,0,0}\makebox(0,0)[lt]{\lineheight{0}\smash{\begin{tabular}[t]{l}$\mu^I_{Ii}$\end{tabular}}}}%
    \put(0.89029628,0.0475301){\color[rgb]{0,0,0}\makebox(0,0)[lt]{\lineheight{0}\smash{\begin{tabular}[t]{l}\gr{I}\end{tabular}}}}%
  \end{picture}%
\endgroup%

%% file: arxiv-figures/ccwcirc_svg-tex.eps_tex
\begingroup%
  \makeatletter%
  \providecommand\color[2][]{%
    \errmessage{(Inkscape) Color is used for the text in Inkscape, but the package 'color.sty' is not loaded}%
    \renewcommand\color[2][]{}%
  }%
  \providecommand\transparent[1]{%
    \errmessage{(Inkscape) Transparency is used (non-zero) for the text in Inkscape, but the package 'transparent.sty' is not loaded}%
    \renewcommand\transparent[1]{}%
  }%
  \providecommand\rotatebox[2]{#2}%
  \newcommand*\fsize{\dimexpr\f@size pt\relax}%
  \newcommand*\lineheight[1]{\fontsize{\fsize}{#1\fsize}\selectfont}%
  \ifx\svgwidth\undefined%
    \setlength{\unitlength}{112.80000767bp}%
    \ifx\svgscale\undefined%
      \relax%
    \else%
      \setlength{\unitlength}{\unitlength * \real{\svgscale}}%
    \fi%
  \else%
    \setlength{\unitlength}{\svgwidth}%
  \fi%
  \global\let\svgwidth\undefined%
  \global\let\svgscale\undefined%
  \makeatother%
  \begin{picture}(1,0.43262408)%
    \lineheight{1}%
    \setlength\tabcolsep{0pt}%
    \put(0,0){\includegraphics[width=\unitlength]{ccwcirc_svg-tex.eps}}%
    \put(0.47261321,0.20788887){\color[rgb]{0,0,0}\makebox(0,0)[lt]{\lineheight{0}\smash{\begin{tabular}[t]{l}=\end{tabular}}}}%
    \put(0.01773048,0.03900695){\color[rgb]{0,0,0}\makebox(0,0)[lt]{\lineheight{0}\smash{\begin{tabular}[t]{l}\gr{Ii}\end{tabular}}}}%
    \put(0.08703257,0.13667821){\color[rgb]{0,0,0}\makebox(0,0)[lt]{\lineheight{0}\smash{\begin{tabular}[t]{l}\gr{I}\end{tabular}}}}%
    \put(0.61034082,0.05017712){\color[rgb]{0,0,0}\makebox(0,0)[lt]{\lineheight{0}\smash{\begin{tabular}[t]{l}\gr{Ii}\end{tabular}}}}%
    \put(0.71985793,0.19503537){\color[rgb]{0,0,0}\makebox(0,0)[lt]{\lineheight{0}\smash{\begin{tabular}[t]{l}$\partial(f)$\end{tabular}}}}%
    \put(0.21631184,0.19503537){\color[rgb]{0,0,0}\makebox(0,0)[lt]{\lineheight{0}\smash{\begin{tabular}[t]{l}$f$\end{tabular}}}}%
  \end{picture}%
\endgroup%

%% file: arxiv-figures/Bsplitting_svg-tex.eps_tex
\begingroup%
  \makeatletter%
  \providecommand\color[2][]{%
    \errmessage{(Inkscape) Color is used for the text in Inkscape, but the package 'color.sty' is not loaded}%
    \renewcommand\color[2][]{}%
  }%
  \providecommand\transparent[1]{%
    \errmessage{(Inkscape) Transparency is used (non-zero) for the text in Inkscape, but the package 'transparent.sty' is not loaded}%
    \renewcommand\transparent[1]{}%
  }%
  \providecommand\rotatebox[2]{#2}%
  \newcommand*\fsize{\dimexpr\f@size pt\relax}%
  \newcommand*\lineheight[1]{\fontsize{\fsize}{#1\fsize}\selectfont}%
  \ifx\svgwidth\undefined%
    \setlength{\unitlength}{192.80001149bp}%
    \ifx\svgscale\undefined%
      \relax%
    \else%
      \setlength{\unitlength}{\unitlength * \real{\svgscale}}%
    \fi%
  \else%
    \setlength{\unitlength}{\svgwidth}%
  \fi%
  \global\let\svgwidth\undefined%
  \global\let\svgscale\undefined%
  \makeatother%
  \begin{picture}(1,0.37759334)%
    \lineheight{1}%
    \setlength\tabcolsep{0pt}%
    \put(0,0){\includegraphics[width=\unitlength]{Bsplitting_svg-tex.eps}}%
    \put(0.39626558,0.20954371){\color[rgb]{0,0,0}\makebox(0,0)[lt]{\lineheight{0}\smash{\begin{tabular}[t]{l}=\end{tabular}}}}%
    \put(0.92323653,0.18464736){\color[rgb]{0,0,0}\makebox(0,0)[lt]{\lineheight{0}\smash{\begin{tabular}[t]{l}\gr{Ii}\end{tabular}}}}%
    \put(0.86514518,0.32572617){\color[rgb]{0,0,0}\makebox(0,0)[lt]{\lineheight{0}\smash{\begin{tabular}[t]{l}\gr{I}\end{tabular}}}}%
    \put(0.65723668,0.30716395){\color[rgb]{0,0,0}\makebox(0,0)[lt]{\lineheight{0}\smash{\begin{tabular}[t]{l}$\Delta^I_{Ii} {}_{(1)}$\end{tabular}}}}%
    \put(0.85630379,0.01736703){\color[rgb]{0,0,0}\makebox(0,0)[lt]{\lineheight{0}\smash{\begin{tabular}[t]{l}\gr{I}\end{tabular}}}}%
    \put(0.65821916,0.04040171){\color[rgb]{0,0,0}\makebox(0,0)[lt]{\lineheight{0}\smash{\begin{tabular}[t]{l}$\Delta^I_{Ii} {}_{(2)}$\end{tabular}}}}%
    \put(0.31327799,0.23029052){\color[rgb]{0,0,0}\makebox(0,0)[lt]{\lineheight{0}\smash{\begin{tabular}[t]{l}\gr{Ii}\end{tabular}}}}%
    \put(0.02282162,0.23029052){\color[rgb]{0,0,0}\makebox(0,0)[lt]{\lineheight{0}\smash{\begin{tabular}[t]{l}\gr{Ii}\end{tabular}}}}%
    \put(0.16804986,0.23029052){\color[rgb]{0,0,0}\makebox(0,0)[lt]{\lineheight{0}\smash{\begin{tabular}[t]{l}\gr{I}\end{tabular}}}}%
  \end{picture}%
\endgroup%

%% file: arxiv-figures/R2oriented_svg-tex.eps_tex
\begingroup%
  \makeatletter%
  \providecommand\color[2][]{%
    \errmessage{(Inkscape) Color is used for the text in Inkscape, but the package 'color.sty' is not loaded}%
    \renewcommand\color[2][]{}%
  }%
  \providecommand\transparent[1]{%
    \errmessage{(Inkscape) Transparency is used (non-zero) for the text in Inkscape, but the package 'transparent.sty' is not loaded}%
    \renewcommand\transparent[1]{}%
  }%
  \providecommand\rotatebox[2]{#2}%
  \newcommand*\fsize{\dimexpr\f@size pt\relax}%
  \newcommand*\lineheight[1]{\fontsize{\fsize}{#1\fsize}\selectfont}%
  \ifx\svgwidth\undefined%
    \setlength{\unitlength}{169.21245384bp}%
    \ifx\svgscale\undefined%
      \relax%
    \else%
      \setlength{\unitlength}{\unitlength * \real{\svgscale}}%
    \fi%
  \else%
    \setlength{\unitlength}{\svgwidth}%
  \fi%
  \global\let\svgwidth\undefined%
  \global\let\svgscale\undefined%
  \makeatother%
  \begin{picture}(1,0.38425099)%
    \lineheight{1}%
    \setlength\tabcolsep{0pt}%
    \put(0,0){\includegraphics[width=\unitlength]{R2oriented_svg-tex.eps}}%
    \put(0.01880664,0.18680802){\color[rgb]{0,0,0}\makebox(0,0)[lt]{\lineheight{1.25}\smash{\begin{tabular}[t]{l}\gr{I}\end{tabular}}}}%
    \put(0.15352452,0.188185){\color[rgb]{0,0,0}\makebox(0,0)[lt]{\lineheight{1.25}\smash{\begin{tabular}[t]{l}\gr{Ij}\end{tabular}}}}%
    \put(0.27883487,0.188185){\color[rgb]{0,0,0}\makebox(0,0)[lt]{\lineheight{1.25}\smash{\begin{tabular}[t]{l}\gr{Iij}\end{tabular}}}}%
    \put(0.16316323,0.02110436){\color[rgb]{0,0,0}\makebox(0,0)[lt]{\lineheight{1.25}\smash{\begin{tabular}[t]{l}\gr{Ii}\end{tabular}}}}%
    \put(0.16293193,0.3357572){\color[rgb]{0,0,0}\makebox(0,0)[lt]{\lineheight{1.25}\smash{\begin{tabular}[t]{l}\gr{Ii}\end{tabular}}}}%
    \put(0.7805287,0.21962337){\color[rgb]{0,0,0}\makebox(0,0)[lt]{\lineheight{1.25}\smash{\begin{tabular}[t]{l}\gr{Ii}\end{tabular}}}}%
    \put(0.6526957,0.21870845){\color[rgb]{0,0,0}\makebox(0,0)[lt]{\lineheight{1.25}\smash{\begin{tabular}[t]{l}\gr{I}\end{tabular}}}}%
    \put(0.90010011,0.22008543){\color[rgb]{0,0,0}\makebox(0,0)[lt]{\lineheight{1.25}\smash{\begin{tabular}[t]{l}\gr{Iij}\end{tabular}}}}%
    \put(0.47170448,0.18785983){\color[rgb]{0,0,0}\makebox(0,0)[lt]{\lineheight{1.25}\smash{\begin{tabular}[t]{l}$=$\end{tabular}}}}%
  \end{picture}%
\endgroup%

%% file: arxiv-figures/R2unoriented1_svg-tex.eps_tex
\begingroup%
  \makeatletter%
  \providecommand\color[2][]{%
    \errmessage{(Inkscape) Color is used for the text in Inkscape, but the package 'color.sty' is not loaded}%
    \renewcommand\color[2][]{}%
  }%
  \providecommand\transparent[1]{%
    \errmessage{(Inkscape) Transparency is used (non-zero) for the text in Inkscape, but the package 'transparent.sty' is not loaded}%
    \renewcommand\transparent[1]{}%
  }%
  \providecommand\rotatebox[2]{#2}%
  \newcommand*\fsize{\dimexpr\f@size pt\relax}%
  \newcommand*\lineheight[1]{\fontsize{\fsize}{#1\fsize}\selectfont}%
  \ifx\svgwidth\undefined%
    \setlength{\unitlength}{328.90000194bp}%
    \ifx\svgscale\undefined%
      \relax%
    \else%
      \setlength{\unitlength}{\unitlength * \real{\svgscale}}%
    \fi%
  \else%
    \setlength{\unitlength}{\svgwidth}%
  \fi%
  \global\let\svgwidth\undefined%
  \global\let\svgscale\undefined%
  \makeatother%
  \begin{picture}(1,0.1984182)%
    \lineheight{1}%
    \setlength\tabcolsep{0pt}%
    \put(0,0){\includegraphics[width=\unitlength]{R2unoriented1_svg-tex.eps}}%
    \put(0.38604046,0.09392284){\color[rgb]{0,0,0}\makebox(0,0)[lt]{\lineheight{0}\smash{\begin{tabular}[t]{l}$\Delta^{Ij}_{Iij, (1)}$\end{tabular}}}}%
    \put(0.73256283,0.09442347){\color[rgb]{0,0,0}\makebox(0,0)[lt]{\lineheight{0}\smash{\begin{tabular}[t]{l}$\dd_{Ij}^I\Delta^{Ii}_{Iij, (2)}$\end{tabular}}}}%
    \put(0.02507575,0.15623451){\color[rgb]{0,0,0}\makebox(0,0)[lt]{\lineheight{1.25}\smash{\begin{tabular}[t]{l}\gr{Ii}\end{tabular}}}}%
    \put(0.09379838,0.09455277){\color[rgb]{0,0,0}\makebox(0,0)[lt]{\lineheight{1.25}\smash{\begin{tabular}[t]{l}\gr{I}\end{tabular}}}}%
    \put(0.13185476,0.03578045){\color[rgb]{0,0,0}\makebox(0,0)[lt]{\lineheight{1.25}\smash{\begin{tabular}[t]{l}\gr{Ij}\end{tabular}}}}%
    \put(0.41611396,0.16577768){\color[rgb]{0,0,0}\makebox(0,0)[lt]{\lineheight{1.25}\smash{\begin{tabular}[t]{l}\gr{Ii}\end{tabular}}}}%
    \put(0.78969519,0.16589417){\color[rgb]{0,0,0}\makebox(0,0)[lt]{\lineheight{1.25}\smash{\begin{tabular}[t]{l}\gr{Ij}\end{tabular}}}}%
    \put(0.08550665,0.01004715){\color[rgb]{0,0,0}\makebox(0,0)[lt]{\lineheight{1.25}\smash{\begin{tabular}[t]{l}\gr{Iij}\end{tabular}}}}%
    \put(0.08551235,0.17414394){\color[rgb]{0,0,0}\makebox(0,0)[lt]{\lineheight{1.25}\smash{\begin{tabular}[t]{l}\gr{Iij}\end{tabular}}}}%
    \put(0.61533559,0.16454259){\color[rgb]{0,0,0}\makebox(0,0)[lt]{\lineheight{1.25}\smash{\begin{tabular}[t]{l}\gr{Iij}\end{tabular}}}}%
    \put(0.24472318,0.10706476){\color[rgb]{0,0,0}\makebox(0,0)[lt]{\begin{minipage}{0.03404139\unitlength}\raggedright $=$\end{minipage}}}%
  \end{picture}%
\endgroup%

%% file: arxiv-figures/R2unoriented2_svg-tex.eps_tex
\begingroup%
  \makeatletter%
  \providecommand\color[2][]{%
    \errmessage{(Inkscape) Color is used for the text in Inkscape, but the package 'color.sty' is not loaded}%
    \renewcommand\color[2][]{}%
  }%
  \providecommand\transparent[1]{%
    \errmessage{(Inkscape) Transparency is used (non-zero) for the text in Inkscape, but the package 'transparent.sty' is not loaded}%
    \renewcommand\transparent[1]{}%
  }%
  \providecommand\rotatebox[2]{#2}%
  \newcommand*\fsize{\dimexpr\f@size pt\relax}%
  \newcommand*\lineheight[1]{\fontsize{\fsize}{#1\fsize}\selectfont}%
  \ifx\svgwidth\undefined%
    \setlength{\unitlength}{187.34771442bp}%
    \ifx\svgscale\undefined%
      \relax%
    \else%
      \setlength{\unitlength}{\unitlength * \real{\svgscale}}%
    \fi%
  \else%
    \setlength{\unitlength}{\svgwidth}%
  \fi%
  \global\let\svgwidth\undefined%
  \global\let\svgscale\undefined%
  \makeatother%
  \begin{picture}(1,0.34656823)%
    \lineheight{1}%
    \setlength\tabcolsep{0pt}%
    \put(0,0){\includegraphics[width=\unitlength]{R2unoriented2_svg-tex.eps}}%
    \put(0.42674669,0.17331071){\color[rgb]{0,0,0}\makebox(0,0)[lt]{\lineheight{1.25}\smash{\begin{tabular}[t]{l}$=$\end{tabular}}}}%
    \put(0.03843118,0.26004791){\color[rgb]{0,0,0}\makebox(0,0)[lt]{\lineheight{1.25}\smash{\begin{tabular}[t]{l}\gr{Ii}\end{tabular}}}}%
    \put(0.16920414,0.01985282){\color[rgb]{0,0,0}\makebox(0,0)[lt]{\lineheight{1.25}\smash{\begin{tabular}[t]{l}\gr{I}\end{tabular}}}}%
    \put(0.15923645,0.30203904){\color[rgb]{0,0,0}\makebox(0,0)[lt]{\lineheight{1.25}\smash{\begin{tabular}[t]{l}\gr{I}\end{tabular}}}}%
    \put(0.2459734,0.07518805){\color[rgb]{0,0,0}\makebox(0,0)[lt]{\lineheight{1.25}\smash{\begin{tabular}[t]{l}\gr{Ij}\end{tabular}}}}%
    \put(0.14055453,0.1659285){\color[rgb]{0,0,0}\makebox(0,0)[lt]{\lineheight{1.25}\smash{\begin{tabular}[t]{l}\gr{Iij}\end{tabular}}}}%
    \put(0.73837333,0.2766852){\color[rgb]{0,0,0}\makebox(0,0)[lt]{\lineheight{1.25}\smash{\begin{tabular}[t]{l}\gr{I}\end{tabular}}}}%
    \put(0.57957785,0.27801953){\color[rgb]{0,0,0}\makebox(0,0)[lt]{\lineheight{1.25}\smash{\begin{tabular}[t]{l}\gr{Ii}\end{tabular}}}}%
    \put(0.9171853,0.27535062){\color[rgb]{0,0,0}\makebox(0,0)[lt]{\lineheight{1.25}\smash{\begin{tabular}[t]{l}\gr{Ij}\end{tabular}}}}%
  \end{picture}%
\endgroup%

%% file: arxiv-figures/R2_easy_distant_ssbim_LHS_svg-tex.eps_tex
\begingroup%
  \makeatletter%
  \providecommand\color[2][]{%
    \errmessage{(Inkscape) Color is used for the text in Inkscape, but the package 'color.sty' is not loaded}%
    \renewcommand\color[2][]{}%
  }%
  \providecommand\transparent[1]{%
    \errmessage{(Inkscape) Transparency is used (non-zero) for the text in Inkscape, but the package 'transparent.sty' is not loaded}%
    \renewcommand\transparent[1]{}%
  }%
  \providecommand\rotatebox[2]{#2}%
  \newcommand*\fsize{\dimexpr\f@size pt\relax}%
  \newcommand*\lineheight[1]{\fontsize{\fsize}{#1\fsize}\selectfont}%
  \ifx\svgwidth\undefined%
    \setlength{\unitlength}{456.45952739bp}%
    \ifx\svgscale\undefined%
      \relax%
    \else%
      \setlength{\unitlength}{\unitlength * \real{\svgscale}}%
    \fi%
  \else%
    \setlength{\unitlength}{\svgwidth}%
  \fi%
  \global\let\svgwidth\undefined%
  \global\let\svgscale\undefined%
  \makeatother%
  \begin{picture}(1,0.93509924)%
    \lineheight{1}%
    \setlength\tabcolsep{0pt}%
    \put(0,0){\includegraphics[width=\unitlength]{R2_easy_distant_ssbim_LHS_svg-tex.eps}}%
    \put(0.13035334,0.76012396){\color[rgb]{0,0,0}\makebox(0,0)[lt]{\lineheight{1.25}\smash{\begin{tabular}[t]{l}\gr{Ii}\end{tabular}}}}%
    \put(0.73171395,0.14112416){\color[rgb]{0,0,0}\makebox(0,0)[lt]{\lineheight{1.25}\smash{\begin{tabular}[t]{l}\gr{Ij}\end{tabular}}}}%
  \end{picture}%
\endgroup%

%% file: arxiv-figures/R2_easy_distant_ssbim_RHS_svg-tex.eps_tex
\begingroup%
  \makeatletter%
  \providecommand\color[2][]{%
    \errmessage{(Inkscape) Color is used for the text in Inkscape, but the package 'color.sty' is not loaded}%
    \renewcommand\color[2][]{}%
  }%
  \providecommand\transparent[1]{%
    \errmessage{(Inkscape) Transparency is used (non-zero) for the text in Inkscape, but the package 'transparent.sty' is not loaded}%
    \renewcommand\transparent[1]{}%
  }%
  \providecommand\rotatebox[2]{#2}%
  \newcommand*\fsize{\dimexpr\f@size pt\relax}%
  \newcommand*\lineheight[1]{\fontsize{\fsize}{#1\fsize}\selectfont}%
  \ifx\svgwidth\undefined%
    \setlength{\unitlength}{456.45952739bp}%
    \ifx\svgscale\undefined%
      \relax%
    \else%
      \setlength{\unitlength}{\unitlength * \real{\svgscale}}%
    \fi%
  \else%
    \setlength{\unitlength}{\svgwidth}%
  \fi%
  \global\let\svgwidth\undefined%
  \global\let\svgscale\undefined%
  \makeatother%
  \begin{picture}(1,0.93456902)%
    \lineheight{1}%
    \setlength\tabcolsep{0pt}%
    \put(0,0){\includegraphics[width=\unitlength]{R2_easy_distant_ssbim_RHS_svg-tex.eps}}%
    \put(0.12505309,0.77337466){\color[rgb]{0,0,0}\makebox(0,0)[lt]{\lineheight{1.25}\smash{\begin{tabular}[t]{l}\gr{Ii}\end{tabular}}}}%
    \put(0.73966426,0.14377427){\color[rgb]{0,0,0}\makebox(0,0)[lt]{\lineheight{1.25}\smash{\begin{tabular}[t]{l}\gr{Ij}\end{tabular}}}}%
  \end{picture}%
\endgroup%

%% file: arxiv-figures/R2_distant_ssbim_LHS_svg-tex.eps_tex
\begingroup%
  \makeatletter%
  \providecommand\color[2][]{%
    \errmessage{(Inkscape) Color is used for the text in Inkscape, but the package 'color.sty' is not loaded}%
    \renewcommand\color[2][]{}%
  }%
  \providecommand\transparent[1]{%
    \errmessage{(Inkscape) Transparency is used (non-zero) for the text in Inkscape, but the package 'transparent.sty' is not loaded}%
    \renewcommand\transparent[1]{}%
  }%
  \providecommand\rotatebox[2]{#2}%
  \newcommand*\fsize{\dimexpr\f@size pt\relax}%
  \newcommand*\lineheight[1]{\fontsize{\fsize}{#1\fsize}\selectfont}%
  \ifx\svgwidth\undefined%
    \setlength{\unitlength}{456.45952739bp}%
    \ifx\svgscale\undefined%
      \relax%
    \else%
      \setlength{\unitlength}{\unitlength * \real{\svgscale}}%
    \fi%
  \else%
    \setlength{\unitlength}{\svgwidth}%
  \fi%
  \global\let\svgwidth\undefined%
  \global\let\svgscale\undefined%
  \makeatother%
  \begin{picture}(1,0.93509924)%
    \lineheight{1}%
    \setlength\tabcolsep{0pt}%
    \put(0,0){\includegraphics[width=\unitlength]{R2_distant_ssbim_LHS_svg-tex.eps}}%
    \put(0.13035334,0.76012396){\color[rgb]{0,0,0}\makebox(0,0)[lt]{\lineheight{1.25}\smash{\begin{tabular}[t]{l}\gr{Ii}\end{tabular}}}}%
    \put(0.73171395,0.14112416){\color[rgb]{0,0,0}\makebox(0,0)[lt]{\lineheight{1.25}\smash{\begin{tabular}[t]{l}\gr{Ij}\end{tabular}}}}%
  \end{picture}%
\endgroup%

%% file: arxiv-figures/R2_distant_ssbim_RHS_svg-tex.eps_tex
\begingroup%
  \makeatletter%
  \providecommand\color[2][]{%
    \errmessage{(Inkscape) Color is used for the text in Inkscape, but the package 'color.sty' is not loaded}%
    \renewcommand\color[2][]{}%
  }%
  \providecommand\transparent[1]{%
    \errmessage{(Inkscape) Transparency is used (non-zero) for the text in Inkscape, but the package 'transparent.sty' is not loaded}%
    \renewcommand\transparent[1]{}%
  }%
  \providecommand\rotatebox[2]{#2}%
  \newcommand*\fsize{\dimexpr\f@size pt\relax}%
  \newcommand*\lineheight[1]{\fontsize{\fsize}{#1\fsize}\selectfont}%
  \ifx\svgwidth\undefined%
    \setlength{\unitlength}{456.45952739bp}%
    \ifx\svgscale\undefined%
      \relax%
    \else%
      \setlength{\unitlength}{\unitlength * \real{\svgscale}}%
    \fi%
  \else%
    \setlength{\unitlength}{\svgwidth}%
  \fi%
  \global\let\svgwidth\undefined%
  \global\let\svgscale\undefined%
  \makeatother%
  \begin{picture}(1,0.93456902)%
    \lineheight{1}%
    \setlength\tabcolsep{0pt}%
    \put(0,0){\includegraphics[width=\unitlength]{R2_distant_ssbim_RHS_svg-tex.eps}}%
    \put(0.12505309,0.77337466){\color[rgb]{0,0,0}\makebox(0,0)[lt]{\lineheight{1.25}\smash{\begin{tabular}[t]{l}\gr{Ii}\end{tabular}}}}%
    \put(0.73966426,0.14377427){\color[rgb]{0,0,0}\makebox(0,0)[lt]{\lineheight{1.25}\smash{\begin{tabular}[t]{l}\gr{Ij}\end{tabular}}}}%
  \end{picture}%
\endgroup%

%% file: arxiv-figures/R3oriented_svg-tex.eps_tex
\begingroup%
  \makeatletter%
  \providecommand\color[2][]{%
    \errmessage{(Inkscape) Color is used for the text in Inkscape, but the package 'color.sty' is not loaded}%
    \renewcommand\color[2][]{}%
  }%
  \providecommand\transparent[1]{%
    \errmessage{(Inkscape) Transparency is used (non-zero) for the text in Inkscape, but the package 'transparent.sty' is not loaded}%
    \renewcommand\transparent[1]{}%
  }%
  \providecommand\rotatebox[2]{#2}%
  \newcommand*\fsize{\dimexpr\f@size pt\relax}%
  \newcommand*\lineheight[1]{\fontsize{\fsize}{#1\fsize}\selectfont}%
  \ifx\svgwidth\undefined%
    \setlength{\unitlength}{228.24983835bp}%
    \ifx\svgscale\undefined%
      \relax%
    \else%
      \setlength{\unitlength}{\unitlength * \real{\svgscale}}%
    \fi%
  \else%
    \setlength{\unitlength}{\svgwidth}%
  \fi%
  \global\let\svgwidth\undefined%
  \global\let\svgscale\undefined%
  \makeatother%
  \begin{picture}(1,0.28455576)%
    \lineheight{1}%
    \setlength\tabcolsep{0pt}%
    \put(0,0){\includegraphics[width=\unitlength]{R3oriented_svg-tex.eps}}%
    \put(0.61533454,0.14232796){\color[rgb]{0,0,0}\makebox(0,0)[lt]{\lineheight{1.25}\smash{\begin{tabular}[t]{l}\gr{I}\end{tabular}}}}%
    \put(0.72048248,0.22995131){\color[rgb]{0,0,0}\makebox(0,0)[lt]{\lineheight{1.25}\smash{\begin{tabular}[t]{l}\gr{Ik}\end{tabular}}}}%
    \put(0.72048248,0.05470481){\color[rgb]{0,0,0}\makebox(0,0)[lt]{\lineheight{1.25}\smash{\begin{tabular}[t]{l}\gr{Ii}\end{tabular}}}}%
    \put(0.80608986,0.14294627){\color[rgb]{0,0,0}\makebox(0,0)[lt]{\lineheight{1.25}\smash{\begin{tabular}[t]{l}\gr{Iik}\end{tabular}}}}%
    \put(0.91325377,0.14232796){\color[rgb]{0,0,0}\makebox(0,0)[lt]{\lineheight{1.25}\smash{\begin{tabular}[t]{l}\gr{Iijk}\end{tabular}}}}%
    \put(0.01949619,0.14232796){\color[rgb]{0,0,0}\makebox(0,0)[lt]{\lineheight{1.25}\smash{\begin{tabular}[t]{l}\gr{I}\end{tabular}}}}%
    \put(0.24731676,0.14232796){\color[rgb]{0,0,0}\makebox(0,0)[lt]{\lineheight{1.25}\smash{\begin{tabular}[t]{l}\gr{Iijk}\end{tabular}}}}%
    \put(0.10711944,0.14232796){\color[rgb]{0,0,0}\makebox(0,0)[lt]{\lineheight{1.25}\smash{\begin{tabular}[t]{l}\gr{Ij}\end{tabular}}}}%
    \put(0.1772181,0.21242672){\color[rgb]{0,0,0}\makebox(0,0)[lt]{\lineheight{1.25}\smash{\begin{tabular}[t]{l}\gr{Ijk}\end{tabular}}}}%
    \put(0.1772181,0.03718002){\color[rgb]{0,0,0}\makebox(0,0)[lt]{\lineheight{1.25}\smash{\begin{tabular}[t]{l}\gr{Iij}\end{tabular}}}}%
    \put(0.08796866,0.01058398){\color[rgb]{0,0,0}\makebox(0,0)[lt]{\lineheight{1.25}\smash{\begin{tabular}[t]{l}\gr{Ii}\end{tabular}}}}%
    \put(0.09200039,0.25693722){\color[rgb]{0,0,0}\makebox(0,0)[lt]{\lineheight{1.25}\smash{\begin{tabular}[t]{l}\gr{Ik}\end{tabular}}}}%
    \put(0.81538974,0.25553956){\color[rgb]{0,0,0}\makebox(0,0)[lt]{\lineheight{1.25}\smash{\begin{tabular}[t]{l}\gr{Ijk}\end{tabular}}}}%
    \put(0.80794458,0.01360793){\color[rgb]{0,0,0}\makebox(0,0)[lt]{\lineheight{1.25}\smash{\begin{tabular}[t]{l}\gr{Iij}\end{tabular}}}}%
    \put(0.4501404,0.14173673){\color[rgb]{0,0,0}\makebox(0,0)[lt]{\lineheight{1.25}\smash{\begin{tabular}[t]{l}$=$\end{tabular}}}}%
  \end{picture}%
\endgroup%

%% file: arxiv-figures/R3unoriented_svg-tex.eps_tex
\begingroup%
  \makeatletter%
  \providecommand\color[2][]{%
    \errmessage{(Inkscape) Color is used for the text in Inkscape, but the package 'color.sty' is not loaded}%
    \renewcommand\color[2][]{}%
  }%
  \providecommand\transparent[1]{%
    \errmessage{(Inkscape) Transparency is used (non-zero) for the text in Inkscape, but the package 'transparent.sty' is not loaded}%
    \renewcommand\transparent[1]{}%
  }%
  \providecommand\rotatebox[2]{#2}%
  \newcommand*\fsize{\dimexpr\f@size pt\relax}%
  \newcommand*\lineheight[1]{\fontsize{\fsize}{#1\fsize}\selectfont}%
  \ifx\svgwidth\undefined%
    \setlength{\unitlength}{240.9041943bp}%
    \ifx\svgscale\undefined%
      \relax%
    \else%
      \setlength{\unitlength}{\unitlength * \real{\svgscale}}%
    \fi%
  \else%
    \setlength{\unitlength}{\svgwidth}%
  \fi%
  \global\let\svgwidth\undefined%
  \global\let\svgscale\undefined%
  \makeatother%
  \begin{picture}(1,0.27013876)%
    \lineheight{1}%
    \setlength\tabcolsep{0pt}%
    \put(0,0){\includegraphics[width=\unitlength]{R3unoriented_svg-tex.eps}}%
    \put(0.04566375,0.13552889){\color[rgb]{0,0,0}\makebox(0,0)[lt]{\lineheight{1.25}\smash{\begin{tabular}[t]{l}\gr{Ij}\end{tabular}}}}%
    \put(0.11807285,0.21854951){\color[rgb]{0,0,0}\makebox(0,0)[lt]{\lineheight{1.25}\smash{\begin{tabular}[t]{l}\gr{Ijk}\end{tabular}}}}%
    \put(0.11807285,0.05250846){\color[rgb]{0,0,0}\makebox(0,0)[lt]{\lineheight{1.25}\smash{\begin{tabular}[t]{l}\gr{Iij}\end{tabular}}}}%
    \put(0.20590703,0.13317076){\color[rgb]{0,0,0}\makebox(0,0)[lt]{\lineheight{1.25}\smash{\begin{tabular}[t]{l}\gr{Iijk}\end{tabular}}}}%
    \put(0.22604819,0.23996731){\color[rgb]{0,0,0}\makebox(0,0)[lt]{\lineheight{1.25}\smash{\begin{tabular}[t]{l}\gr{Ik}\end{tabular}}}}%
    \put(0.23194343,0.01104689){\color[rgb]{0,0,0}\makebox(0,0)[lt]{\lineheight{1.25}\smash{\begin{tabular}[t]{l}\gr{Ii}\end{tabular}}}}%
    \put(0.31537744,0.13303945){\color[rgb]{0,0,0}\makebox(0,0)[lt]{\lineheight{1.25}\smash{\begin{tabular}[t]{l}\gr{Iik}\end{tabular}}}}%
    \put(0.5777849,0.13434982){\color[rgb]{0,0,0}\makebox(0,0)[lt]{\lineheight{1.25}\smash{\begin{tabular}[t]{l}\gr{Ij}\end{tabular}}}}%
    \put(0.70513589,0.07919245){\color[rgb]{0,0,0}\makebox(0,0)[lt]{\lineheight{1.25}\smash{\begin{tabular}[t]{l}\gr{I}\end{tabular}}}}%
    \put(0.79886522,0.21841821){\color[rgb]{0,0,0}\makebox(0,0)[lt]{\lineheight{1.25}\smash{\begin{tabular}[t]{l}\gr{Ik}\end{tabular}}}}%
    \put(0.81546946,0.03577292){\color[rgb]{0,0,0}\makebox(0,0)[lt]{\lineheight{1.25}\smash{\begin{tabular}[t]{l}\gr{Ii}\end{tabular}}}}%
    \put(0.89259455,0.13294216){\color[rgb]{0,0,0}\makebox(0,0)[lt]{\lineheight{1.25}\smash{\begin{tabular}[t]{l}\gr{Iik}\end{tabular}}}}%
    \put(0.45017889,0.13539759){\color[rgb]{0,0,0}\makebox(0,0)[lt]{\lineheight{1.25}\smash{\begin{tabular}[t]{l}$=$\end{tabular}}}}%
    \put(0.70471633,0.24280182){\color[rgb]{0,0,0}\makebox(0,0)[lt]{\lineheight{1.25}\smash{\begin{tabular}[t]{l}\gr{Ijk}\end{tabular}}}}%
    \put(0.70568239,0.0117074){\color[rgb]{0,0,0}\makebox(0,0)[lt]{\lineheight{1.25}\smash{\begin{tabular}[t]{l}\gr{Iij}\end{tabular}}}}%
  \end{picture}%
\endgroup%

%% file: arxiv-figures/R3_distant_ssbim_LHS_svg-tex.eps_tex
\begingroup%
  \makeatletter%
  \providecommand\color[2][]{%
    \errmessage{(Inkscape) Color is used for the text in Inkscape, but the package 'color.sty' is not loaded}%
    \renewcommand\color[2][]{}%
  }%
  \providecommand\transparent[1]{%
    \errmessage{(Inkscape) Transparency is used (non-zero) for the text in Inkscape, but the package 'transparent.sty' is not loaded}%
    \renewcommand\transparent[1]{}%
  }%
  \providecommand\rotatebox[2]{#2}%
  \newcommand*\fsize{\dimexpr\f@size pt\relax}%
  \newcommand*\lineheight[1]{\fontsize{\fsize}{#1\fsize}\selectfont}%
  \ifx\svgwidth\undefined%
    \setlength{\unitlength}{537.50792083bp}%
    \ifx\svgscale\undefined%
      \relax%
    \else%
      \setlength{\unitlength}{\unitlength * \real{\svgscale}}%
    \fi%
  \else%
    \setlength{\unitlength}{\svgwidth}%
  \fi%
  \global\let\svgwidth\undefined%
  \global\let\svgscale\undefined%
  \makeatother%
  \begin{picture}(1,0.79364966)%
    \lineheight{1}%
    \setlength\tabcolsep{0pt}%
    \put(0,0){\includegraphics[width=\unitlength]{R3_distant_ssbim_LHS_svg-tex.eps}}%
    \put(0.72490643,0.37190408){\color[rgb]{0,0,0}\makebox(0,0)[lt]{\lineheight{1.25}\smash{\begin{tabular}[t]{l}\gr{Ijk}\end{tabular}}}}%
    \put(0.12195058,0.37319425){\color[rgb]{0,0,0}\makebox(0,0)[lt]{\lineheight{1.25}\smash{\begin{tabular}[t]{l}\gr{Ii}\end{tabular}}}}%
  \end{picture}%
\endgroup%

%% file: arxiv-figures/R3_distant_ssbim_RHS_svg-tex.eps_tex
\begingroup%
  \makeatletter%
  \providecommand\color[2][]{%
    \errmessage{(Inkscape) Color is used for the text in Inkscape, but the package 'color.sty' is not loaded}%
    \renewcommand\color[2][]{}%
  }%
  \providecommand\transparent[1]{%
    \errmessage{(Inkscape) Transparency is used (non-zero) for the text in Inkscape, but the package 'transparent.sty' is not loaded}%
    \renewcommand\transparent[1]{}%
  }%
  \providecommand\rotatebox[2]{#2}%
  \newcommand*\fsize{\dimexpr\f@size pt\relax}%
  \newcommand*\lineheight[1]{\fontsize{\fsize}{#1\fsize}\selectfont}%
  \ifx\svgwidth\undefined%
    \setlength{\unitlength}{507.26599205bp}%
    \ifx\svgscale\undefined%
      \relax%
    \else%
      \setlength{\unitlength}{\unitlength * \real{\svgscale}}%
    \fi%
  \else%
    \setlength{\unitlength}{\svgwidth}%
  \fi%
  \global\let\svgwidth\undefined%
  \global\let\svgscale\undefined%
  \makeatother%
  \begin{picture}(1,0.84096501)%
    \lineheight{1}%
    \setlength\tabcolsep{0pt}%
    \put(0,0){\includegraphics[width=\unitlength]{R3_distant_ssbim_RHS_svg-tex.eps}}%
    \put(0.71614996,0.41076893){\color[rgb]{0,0,0}\makebox(0,0)[lt]{\lineheight{1.25}\smash{\begin{tabular}[t]{l}\gr{Ijk}\end{tabular}}}}%
    \put(0.13707689,0.40975135){\color[rgb]{0,0,0}\makebox(0,0)[lt]{\lineheight{1.25}\smash{\begin{tabular}[t]{l}\gr{Ii}\end{tabular}}}}%
  \end{picture}%
\endgroup%

%% file: arxiv-figures/cap_ssbim_c_svg-tex.eps_tex
\begingroup%
  \makeatletter%
  \providecommand\color[2][]{%
    \errmessage{(Inkscape) Color is used for the text in Inkscape, but the package 'color.sty' is not loaded}%
    \renewcommand\color[2][]{}%
  }%
  \providecommand\transparent[1]{%
    \errmessage{(Inkscape) Transparency is used (non-zero) for the text in Inkscape, but the package 'transparent.sty' is not loaded}%
    \renewcommand\transparent[1]{}%
  }%
  \providecommand\rotatebox[2]{#2}%
  \newcommand*\fsize{\dimexpr\f@size pt\relax}%
  \newcommand*\lineheight[1]{\fontsize{\fsize}{#1\fsize}\selectfont}%
  \ifx\svgwidth\undefined%
    \setlength{\unitlength}{469.76640848bp}%
    \ifx\svgscale\undefined%
      \relax%
    \else%
      \setlength{\unitlength}{\unitlength * \real{\svgscale}}%
    \fi%
  \else%
    \setlength{\unitlength}{\svgwidth}%
  \fi%
  \global\let\svgwidth\undefined%
  \global\let\svgscale\undefined%
  \makeatother%
  \begin{picture}(1,0.92106259)%
    \lineheight{1}%
    \setlength\tabcolsep{0pt}%
    \put(0,0){\includegraphics[width=\unitlength]{cap_ssbim_c_svg-tex.eps}}%
    \put(0.47167429,0.85709077){\color[rgb]{0,0,0}\makebox(0,0)[lt]{\begin{minipage}{0.3502084\unitlength}\raggedright \dgr{a}\end{minipage}}}%
    \put(0.39442245,0.15667402){\color[rgb]{0,0,0}\makebox(0,0)[lt]{\begin{minipage}{0.51501237\unitlength}\raggedright \dgr{a-k}\end{minipage}}}%
    \put(0.33777112,0.43478071){\color[rgb]{0,0,0}\makebox(0,0)[lt]{\begin{minipage}{0.56136353\unitlength}\raggedright \dgr{a,a-k}\end{minipage}}}%
  \end{picture}%
\endgroup%

%% file: arxiv-figures/cap_ssbim_cc_svg-tex.eps_tex
\begingroup%
  \makeatletter%
  \providecommand\color[2][]{%
    \errmessage{(Inkscape) Color is used for the text in Inkscape, but the package 'color.sty' is not loaded}%
    \renewcommand\color[2][]{}%
  }%
  \providecommand\transparent[1]{%
    \errmessage{(Inkscape) Transparency is used (non-zero) for the text in Inkscape, but the package 'transparent.sty' is not loaded}%
    \renewcommand\transparent[1]{}%
  }%
  \providecommand\rotatebox[2]{#2}%
  \newcommand*\fsize{\dimexpr\f@size pt\relax}%
  \newcommand*\lineheight[1]{\fontsize{\fsize}{#1\fsize}\selectfont}%
  \ifx\svgwidth\undefined%
    \setlength{\unitlength}{469.76640848bp}%
    \ifx\svgscale\undefined%
      \relax%
    \else%
      \setlength{\unitlength}{\unitlength * \real{\svgscale}}%
    \fi%
  \else%
    \setlength{\unitlength}{\svgwidth}%
  \fi%
  \global\let\svgwidth\undefined%
  \global\let\svgscale\undefined%
  \makeatother%
  \begin{picture}(1,0.92106259)%
    \lineheight{1}%
    \setlength\tabcolsep{0pt}%
    \put(0,0){\includegraphics[width=\unitlength]{cap_ssbim_cc_svg-tex.eps}}%
    \put(0.47167429,0.85709077){\color[rgb]{0,0,0}\makebox(0,0)[lt]{\begin{minipage}{0.3502084\unitlength}\raggedright \dgr{a}\end{minipage}}}%
    \put(0.39442245,0.15667402){\color[rgb]{0,0,0}\makebox(0,0)[lt]{\begin{minipage}{0.51501237\unitlength}\raggedright \dgr{a+k}\end{minipage}}}%
    \put(0.33777112,0.43478071){\color[rgb]{0,0,0}\makebox(0,0)[lt]{\begin{minipage}{0.56136353\unitlength}\raggedright \dgr{a,a+k}\end{minipage}}}%
  \end{picture}%
\endgroup%

%% file: arxiv-figures/cup_ssbim_c_svg-tex.eps_tex
\begingroup%
  \makeatletter%
  \providecommand\color[2][]{%
    \errmessage{(Inkscape) Color is used for the text in Inkscape, but the package 'color.sty' is not loaded}%
    \renewcommand\color[2][]{}%
  }%
  \providecommand\transparent[1]{%
    \errmessage{(Inkscape) Transparency is used (non-zero) for the text in Inkscape, but the package 'transparent.sty' is not loaded}%
    \renewcommand\transparent[1]{}%
  }%
  \providecommand\rotatebox[2]{#2}%
  \newcommand*\fsize{\dimexpr\f@size pt\relax}%
  \newcommand*\lineheight[1]{\fontsize{\fsize}{#1\fsize}\selectfont}%
  \ifx\svgwidth\undefined%
    \setlength{\unitlength}{469.76640848bp}%
    \ifx\svgscale\undefined%
      \relax%
    \else%
      \setlength{\unitlength}{\unitlength * \real{\svgscale}}%
    \fi%
  \else%
    \setlength{\unitlength}{\svgwidth}%
  \fi%
  \global\let\svgwidth\undefined%
  \global\let\svgscale\undefined%
  \makeatother%
  \begin{picture}(1,0.92106259)%
    \lineheight{1}%
    \setlength\tabcolsep{0pt}%
    \put(0,0){\includegraphics[width=\unitlength]{cup_ssbim_c_svg-tex.eps}}%
    \put(0.47939945,0.20817522){\color[rgb]{0,0,0}\makebox(0,0)[lt]{\lineheight{1.25}\smash{\begin{tabular}[t]{l}\dgr{a}\end{tabular}}}}%
    \put(0.41759802,0.79271426){\color[rgb]{0,0,0}\makebox(0,0)[lt]{\lineheight{1.25}\smash{\begin{tabular}[t]{l}\dgr{a-k}\end{tabular}}}}%
    \put(0.36352177,0.49400718){\color[rgb]{0,0,0}\makebox(0,0)[lt]{\lineheight{1.25}\smash{\begin{tabular}[t]{l}\dgr{a,a-k}\end{tabular}}}}%
  \end{picture}%
\endgroup%

%% file: arxiv-figures/cup_ssbim_cc_svg-tex.eps_tex
\begingroup%
  \makeatletter%
  \providecommand\color[2][]{%
    \errmessage{(Inkscape) Color is used for the text in Inkscape, but the package 'color.sty' is not loaded}%
    \renewcommand\color[2][]{}%
  }%
  \providecommand\transparent[1]{%
    \errmessage{(Inkscape) Transparency is used (non-zero) for the text in Inkscape, but the package 'transparent.sty' is not loaded}%
    \renewcommand\transparent[1]{}%
  }%
  \providecommand\rotatebox[2]{#2}%
  \newcommand*\fsize{\dimexpr\f@size pt\relax}%
  \newcommand*\lineheight[1]{\fontsize{\fsize}{#1\fsize}\selectfont}%
  \ifx\svgwidth\undefined%
    \setlength{\unitlength}{469.76640848bp}%
    \ifx\svgscale\undefined%
      \relax%
    \else%
      \setlength{\unitlength}{\unitlength * \real{\svgscale}}%
    \fi%
  \else%
    \setlength{\unitlength}{\svgwidth}%
  \fi%
  \global\let\svgwidth\undefined%
  \global\let\svgscale\undefined%
  \makeatother%
  \begin{picture}(1,0.92106259)%
    \lineheight{1}%
    \setlength\tabcolsep{0pt}%
    \put(0,0){\includegraphics[width=\unitlength]{cup_ssbim_cc_svg-tex.eps}}%
    \put(0.47939945,0.21332538){\color[rgb]{0,0,0}\makebox(0,0)[lt]{\lineheight{1.25}\smash{\begin{tabular}[t]{l}\dgr{a}\end{tabular}}}}%
    \put(0.40472267,0.77468885){\color[rgb]{0,0,0}\makebox(0,0)[lt]{\lineheight{1.25}\smash{\begin{tabular}[t]{l}\dgr{a+k}\end{tabular}}}}%
    \put(0.37639703,0.50173237){\color[rgb]{0,0,0}\makebox(0,0)[lt]{\lineheight{1.25}\smash{\begin{tabular}[t]{l}\dgr{a,a+k}\end{tabular}}}}%
  \end{picture}%
\endgroup%

%% file: arxiv-figures/tag_left_c_svg-tex.eps_tex
\begingroup%
  \makeatletter%
  \providecommand\color[2][]{%
    \errmessage{(Inkscape) Color is used for the text in Inkscape, but the package 'color.sty' is not loaded}%
    \renewcommand\color[2][]{}%
  }%
  \providecommand\transparent[1]{%
    \errmessage{(Inkscape) Transparency is used (non-zero) for the text in Inkscape, but the package 'transparent.sty' is not loaded}%
    \renewcommand\transparent[1]{}%
  }%
  \providecommand\rotatebox[2]{#2}%
  \newcommand*\fsize{\dimexpr\f@size pt\relax}%
  \newcommand*\lineheight[1]{\fontsize{\fsize}{#1\fsize}\selectfont}%
  \ifx\svgwidth\undefined%
    \setlength{\unitlength}{469.76640848bp}%
    \ifx\svgscale\undefined%
      \relax%
    \else%
      \setlength{\unitlength}{\unitlength * \real{\svgscale}}%
    \fi%
  \else%
    \setlength{\unitlength}{\svgwidth}%
  \fi%
  \global\let\svgwidth\undefined%
  \global\let\svgscale\undefined%
  \makeatother%
  \begin{picture}(1,0.92106259)%
    \lineheight{1}%
    \setlength\tabcolsep{0pt}%
    \put(0,0){\includegraphics[width=\unitlength]{tag_left_c_svg-tex.eps}}%
    \put(0.70342989,0.44250588){\color[rgb]{0,0,0}\makebox(0,0)[lt]{\lineheight{1.25}\smash{\begin{tabular}[t]{l}\gr{a}\end{tabular}}}}%
    \put(0.10086545,0.44250588){\color[rgb]{0,0,0}\makebox(0,0)[lt]{\lineheight{1.25}\smash{\begin{tabular}[t]{l}\gr{a+k}\end{tabular}}}}%
    \put(0.53090066,0.26225155){\color[rgb]{0,0,0}\makebox(0,0)[lt]{\lineheight{1.25}\smash{\begin{tabular}[t]{l}\blk{k}\end{tabular}}}}%
    \put(0.52060039,0.63821059){\color[rgb]{0,0,0}\makebox(0,0)[lt]{\lineheight{1.25}\smash{\begin{tabular}[t]{l}\blk{n-k}\end{tabular}}}}%
  \end{picture}%
\endgroup%

%% file: arxiv-figures/image_of_tag_c_svg-tex.eps_tex
\begingroup%
  \makeatletter%
  \providecommand\color[2][]{%
    \errmessage{(Inkscape) Color is used for the text in Inkscape, but the package 'color.sty' is not loaded}%
    \renewcommand\color[2][]{}%
  }%
  \providecommand\transparent[1]{%
    \errmessage{(Inkscape) Transparency is used (non-zero) for the text in Inkscape, but the package 'transparent.sty' is not loaded}%
    \renewcommand\transparent[1]{}%
  }%
  \providecommand\rotatebox[2]{#2}%
  \newcommand*\fsize{\dimexpr\f@size pt\relax}%
  \newcommand*\lineheight[1]{\fontsize{\fsize}{#1\fsize}\selectfont}%
  \ifx\svgwidth\undefined%
    \setlength{\unitlength}{469.76640848bp}%
    \ifx\svgscale\undefined%
      \relax%
    \else%
      \setlength{\unitlength}{\unitlength * \real{\svgscale}}%
    \fi%
  \else%
    \setlength{\unitlength}{\svgwidth}%
  \fi%
  \global\let\svgwidth\undefined%
  \global\let\svgscale\undefined%
  \makeatother%
  \begin{picture}(1,0.92106259)%
    \lineheight{1}%
    \setlength\tabcolsep{0pt}%
    \put(0,0){\includegraphics[width=\unitlength]{image_of_tag_c_svg-tex.eps}}%
    \put(0.73948071,0.46310633){\color[rgb]{0,0,0}\makebox(0,0)[lt]{\lineheight{1.25}\smash{\begin{tabular}[t]{l}\dgr{a}\end{tabular}}}}%
    \put(0.05966447,0.45795617){\color[rgb]{0,0,0}\makebox(0,0)[lt]{\lineheight{1.25}\smash{\begin{tabular}[t]{l}\dgr{a+k}\end{tabular}}}}%
  \end{picture}%
\endgroup%

%% file: arxiv-figures/tag_left_d_svg-tex.eps_tex
\begingroup%
  \makeatletter%
  \providecommand\color[2][]{%
    \errmessage{(Inkscape) Color is used for the text in Inkscape, but the package 'color.sty' is not loaded}%
    \renewcommand\color[2][]{}%
  }%
  \providecommand\transparent[1]{%
    \errmessage{(Inkscape) Transparency is used (non-zero) for the text in Inkscape, but the package 'transparent.sty' is not loaded}%
    \renewcommand\transparent[1]{}%
  }%
  \providecommand\rotatebox[2]{#2}%
  \newcommand*\fsize{\dimexpr\f@size pt\relax}%
  \newcommand*\lineheight[1]{\fontsize{\fsize}{#1\fsize}\selectfont}%
  \ifx\svgwidth\undefined%
    \setlength{\unitlength}{469.76640848bp}%
    \ifx\svgscale\undefined%
      \relax%
    \else%
      \setlength{\unitlength}{\unitlength * \real{\svgscale}}%
    \fi%
  \else%
    \setlength{\unitlength}{\svgwidth}%
  \fi%
  \global\let\svgwidth\undefined%
  \global\let\svgscale\undefined%
  \makeatother%
  \begin{picture}(1,0.92106259)%
    \lineheight{1}%
    \setlength\tabcolsep{0pt}%
    \put(0,0){\includegraphics[width=\unitlength]{tag_left_d_svg-tex.eps}}%
    \put(0.65450371,0.445081){\color[rgb]{0,0,0}\makebox(0,0)[lt]{\lineheight{1.25}\smash{\begin{tabular}[t]{l}\gr{a}\end{tabular}}}}%
    \put(0.08799013,0.44250588){\color[rgb]{0,0,0}\makebox(0,0)[lt]{\lineheight{1.25}\smash{\begin{tabular}[t]{l}\gr{a-k}\end{tabular}}}}%
    \put(0.53605082,0.18242471){\color[rgb]{0,0,0}\makebox(0,0)[lt]{\lineheight{1.25}\smash{\begin{tabular}[t]{l}\blk{k}\end{tabular}}}}%
    \put(0.53347574,0.74378808){\color[rgb]{0,0,0}\makebox(0,0)[lt]{\lineheight{1.25}\smash{\begin{tabular}[t]{l}\blk{n-k}\end{tabular}}}}%
  \end{picture}%
\endgroup%

%% file: arxiv-figures/image_of_tag_d_svg-tex.eps_tex
\begingroup%
  \makeatletter%
  \providecommand\color[2][]{%
    \errmessage{(Inkscape) Color is used for the text in Inkscape, but the package 'color.sty' is not loaded}%
    \renewcommand\color[2][]{}%
  }%
  \providecommand\transparent[1]{%
    \errmessage{(Inkscape) Transparency is used (non-zero) for the text in Inkscape, but the package 'transparent.sty' is not loaded}%
    \renewcommand\transparent[1]{}%
  }%
  \providecommand\rotatebox[2]{#2}%
  \newcommand*\fsize{\dimexpr\f@size pt\relax}%
  \newcommand*\lineheight[1]{\fontsize{\fsize}{#1\fsize}\selectfont}%
  \ifx\svgwidth\undefined%
    \setlength{\unitlength}{469.76640848bp}%
    \ifx\svgscale\undefined%
      \relax%
    \else%
      \setlength{\unitlength}{\unitlength * \real{\svgscale}}%
    \fi%
  \else%
    \setlength{\unitlength}{\svgwidth}%
  \fi%
  \global\let\svgwidth\undefined%
  \global\let\svgscale\undefined%
  \makeatother%
  \begin{picture}(1,0.92106259)%
    \lineheight{1}%
    \setlength\tabcolsep{0pt}%
    \put(0,0){\includegraphics[width=\unitlength]{image_of_tag_d_svg-tex.eps}}%
    \put(0.73948071,0.46310633){\color[rgb]{0,0,0}\makebox(0,0)[lt]{\lineheight{1.25}\smash{\begin{tabular}[t]{l}\dgr{a}\end{tabular}}}}%
    \put(0.05966447,0.45795617){\color[rgb]{0,0,0}\makebox(0,0)[lt]{\lineheight{1.25}\smash{\begin{tabular}[t]{l}\dgr{a-k}\end{tabular}}}}%
  \end{picture}%
\endgroup%

%% file: arxiv-figures/image_of_up_trivalent_svg-tex.eps_tex
\begingroup%
  \makeatletter%
  \providecommand\color[2][]{%
    \errmessage{(Inkscape) Color is used for the text in Inkscape, but the package 'color.sty' is not loaded}%
    \renewcommand\color[2][]{}%
  }%
  \providecommand\transparent[1]{%
    \errmessage{(Inkscape) Transparency is used (non-zero) for the text in Inkscape, but the package 'transparent.sty' is not loaded}%
    \renewcommand\transparent[1]{}%
  }%
  \providecommand\rotatebox[2]{#2}%
  \newcommand*\fsize{\dimexpr\f@size pt\relax}%
  \newcommand*\lineheight[1]{\fontsize{\fsize}{#1\fsize}\selectfont}%
  \ifx\svgwidth\undefined%
    \setlength{\unitlength}{469.83444011bp}%
    \ifx\svgscale\undefined%
      \relax%
    \else%
      \setlength{\unitlength}{\unitlength * \real{\svgscale}}%
    \fi%
  \else%
    \setlength{\unitlength}{\svgwidth}%
  \fi%
  \global\let\svgwidth\undefined%
  \global\let\svgscale\undefined%
  \makeatother%
  \begin{picture}(1,0.92107402)%
    \lineheight{1}%
    \setlength\tabcolsep{0pt}%
    \put(0,0){\includegraphics[width=\unitlength]{image_of_up_trivalent_svg-tex.eps}}%
    \put(0.75271986,0.52579907){\color[rgb]{0,0,0}\makebox(0,0)[lt]{\lineheight{1.25}\smash{\begin{tabular}[t]{l}\dgr{a}\end{tabular}}}}%
    \put(0.02751587,0.52579898){\color[rgb]{0,0,0}\makebox(0,0)[lt]{\lineheight{1.25}\smash{\begin{tabular}[t]{l}\dgr{a+k+l}\end{tabular}}}}%
    \put(0.39223929,0.03918291){\color[rgb]{0,0,0}\makebox(0,0)[lt]{\lineheight{1.25}\smash{\begin{tabular}[t]{l}\dgr{a+l}\end{tabular}}}}%
  \end{picture}%
\endgroup%

%% file: Well-Definedness_of_the_Functor.tex
In this section, we show that the functor $\GS_{\zeta}$ defined in \cref{subsec the functor} is well-defined, i.e., the web relations in \cref{fig:webrel} hold in $\dmSBSBim$ after applying $\GS_{\zeta}$.

\subsection{Tag Switching Relations}\label{subsec- tag switching relation}
For any valid region labels, we have the tag switching relation in $\cwebs$ given below (up to mirror images and arrow reversals):
\begin{align}
    {
        \phantom{ \tikz[x=1mm, y=1mm, baseline=-0.5ex]{\draw (0,-10) -- (0,10)}} 
        \tikz[x=1mm,
  y=1mm,
  baseline=-0.5ex] {
          \draw[->] (0,0) -- (0,6) node [right] {$n-k$};
          \draw[->] (0,0) -- (0,-6) node [right] {$k$};
          \draw (0,0) -- (2,0);
        }
        = (-1)^{k(n - k)} \tikz[x=1mm, y=1mm, baseline=-0.5ex]{
          \draw[->] (0,0) -- (0,6) node [right] {$n-k$};
          \draw[->] (0,0) -- (0,-6) node [right] {$k$};
          \draw (0,0) -- (-2,0);
        }
      }. \label{diag-tag switching}\end{align}
We have that
\[ \vcenter{\xy (0,0)*{\def\svgscale{0.15}\input{arxiv-figures/tag_right_d_svg-tex.eps_tex}} \endxy} = \vcenter{\xy (0,0)*{\def\svgscale{0.15}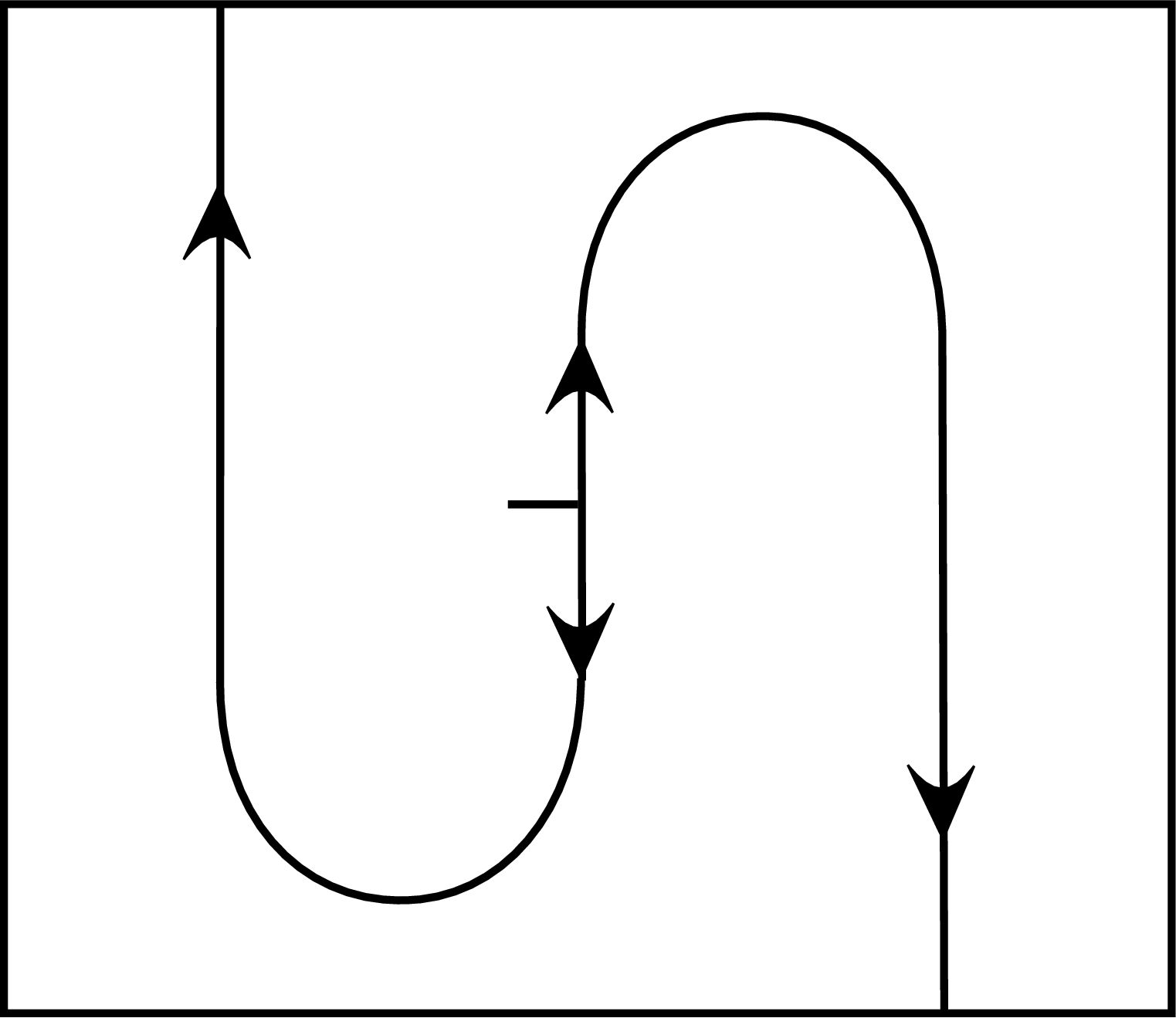} \endxy}   \]
The LHS of \cref{diag-tag switching} is sent via $\GS_{\zeta}$ to
$$(-1)^{k(n-k)} \nu_{n-k}\; \vcenter{\xy (0,0)*{\def\svgscale{0.15}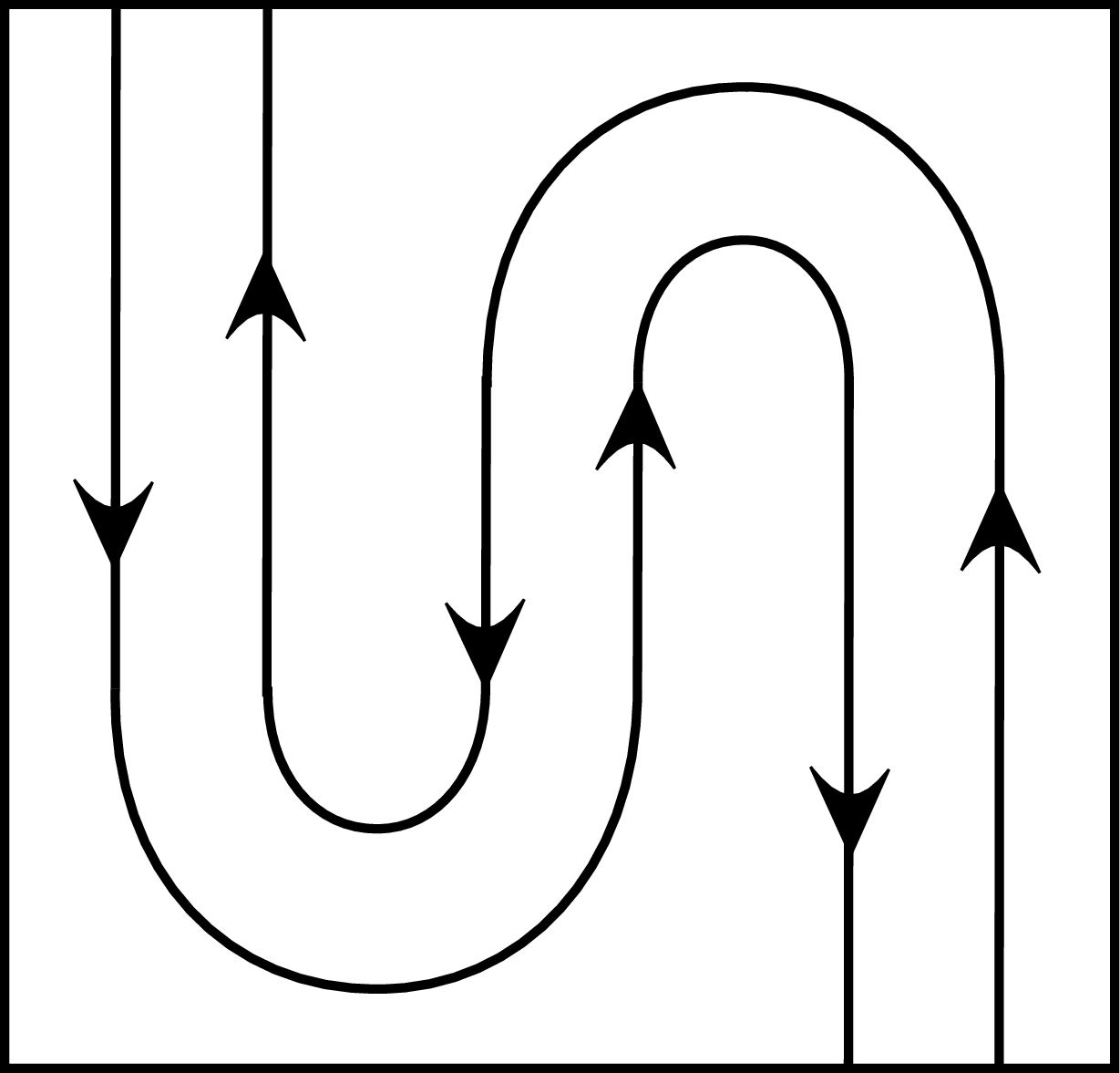} \endxy}= (-1)^{k(n-k)}\nu_{n-k}\; \vcenter{\xy (0,0)*{\def\svgscale{0.15}\input{arxiv-figures/image_of_tag_d_svg-tex.eps_tex}} \endxy}.$$
where the last equality is from the isotopy relations.
The RHS of \cref{diag-tag switching} is sent via $\GS_{\zeta}$ to 
$$(-1)^{k(n-k)} \nu_{k}\; \vcenter{\xy (0,0)*{\def\svgscale{0.15}\input{arxiv-figures/image_of_tag_d_svg-tex.eps_tex}} \endxy}.$$
Hence, we see that \cref{diag-tag switching} holds after applying $\GS_{\zeta}$ if and only if $\nu_k=\nu_{n-k}$.

Similarly, for the relation \cref{diag-tag switching} with reverse arrows, we have the condition that $\lambda_k=\lambda_{n-k}$.
These conditions are exactly those in \cref{cond-tag switching}.
There are no non-trivial mirror images to consider in this case.

\subsection{Tag Cancellation}\label{subsec-tag cancellation}
For any valid region labels, we have the tag cancellation relation in $\cwebs$ given below (up to mirror images and arrow reversals):
\begin{align}
    {
        \tikz[x=1mm,
         y=1mm,
        baseline=-0.5ex] {
          \draw[->-] (0,-6) -- (0,-2);
          \draw (0,-2) -- (2,-2);
          \draw[->-] (0,2) -- (0,-2);
          \draw (0,2) -- (2,2);
          \draw[->] (0,2) -- (0,6) node [right] {$k$};
        }
        = \, \tikz[  x=1mm,
  y=1mm,
  baseline=-0.5ex] {
          \draw[->] (0,-6) -- (0,6) node [right] {$k$};
        }
      }.\label{diag-tag cancellation}
\end{align}
$\GS_{\zeta}$ applied to the LHS is just $\lambda_k\;\nu_{n-k}$ times $\GS_{\zeta}$ applied to the LHS. So this relation imposes the constraints $\lambda_{k} \;\nu_{n-k}=1$, which along with \cref{cond-tag switching} gives us that $\lambda_{k}=\nu_{k}^{-1}$, which is precisely \cref{cond-tag cancellation}.
It is easy to check that the mirror images and arrow reversals also yield the same conditions \cref{cond-tag cancellation}.

\subsection{I=H relation}\label{subsec- I=H}
For valid region labels, we have the I=H relations in $\cwebs$ given below (up to mirror images and arrow reversals):
\begin{gather}{
        \phantom{ \tikz[x=1mm, y=1mm, baseline=-0.5ex] {\draw (0,-10) -- (0,10)}} 
        \tikz[x=1mm, y=1mm, baseline=-0.5ex] {
          \draw[->-] (-6,-6) node [left] {$l$}
            .. controls (-5,-5) and (-4,-3) .. (0,-3);
          \draw[->-] (6,-6) node [right] {$m$}
            .. controls (5,-5) and (4,-3) .. (0,-3);
          \draw[->-] (0,-3) -- (0,3);
          \draw[->-] (-6,6) node [left] {$k$}
            .. controls (-5,5) and (-5,3) .. (0,3);
          \draw[->] (0,3) .. controls (4,3) and (5,5) .. (6,6);
        }
        = \tikz[x=1mm, y=1mm, baseline=-0.5ex] {
          \draw[->-] (-6,-6) node [left] {$l$}
            .. controls (-5,-5) and (-3,-4) .. (-3,0);
          \draw[->-] (-6,6) node [left] {$k$}
            .. controls (-5,5) and (-3,4) .. (-3,0);
          \draw[->-] (-3,0) -- (3,0);
          \draw[->-] (6,-6) node [right] {$m$}
            .. controls (5,-5) and (3,-4) .. (3,0);
          \draw[->] (3,0) .. controls (3,4) and (5,5) .. (6,6);
        }
      }.
\end{gather}
Adding a clockwise cap on the top left, this relation is equivalent to
\begin{align}\label{diag-I=H relation}
    \vcenter{\xy (0,0)*{\def\svgscale{0.15}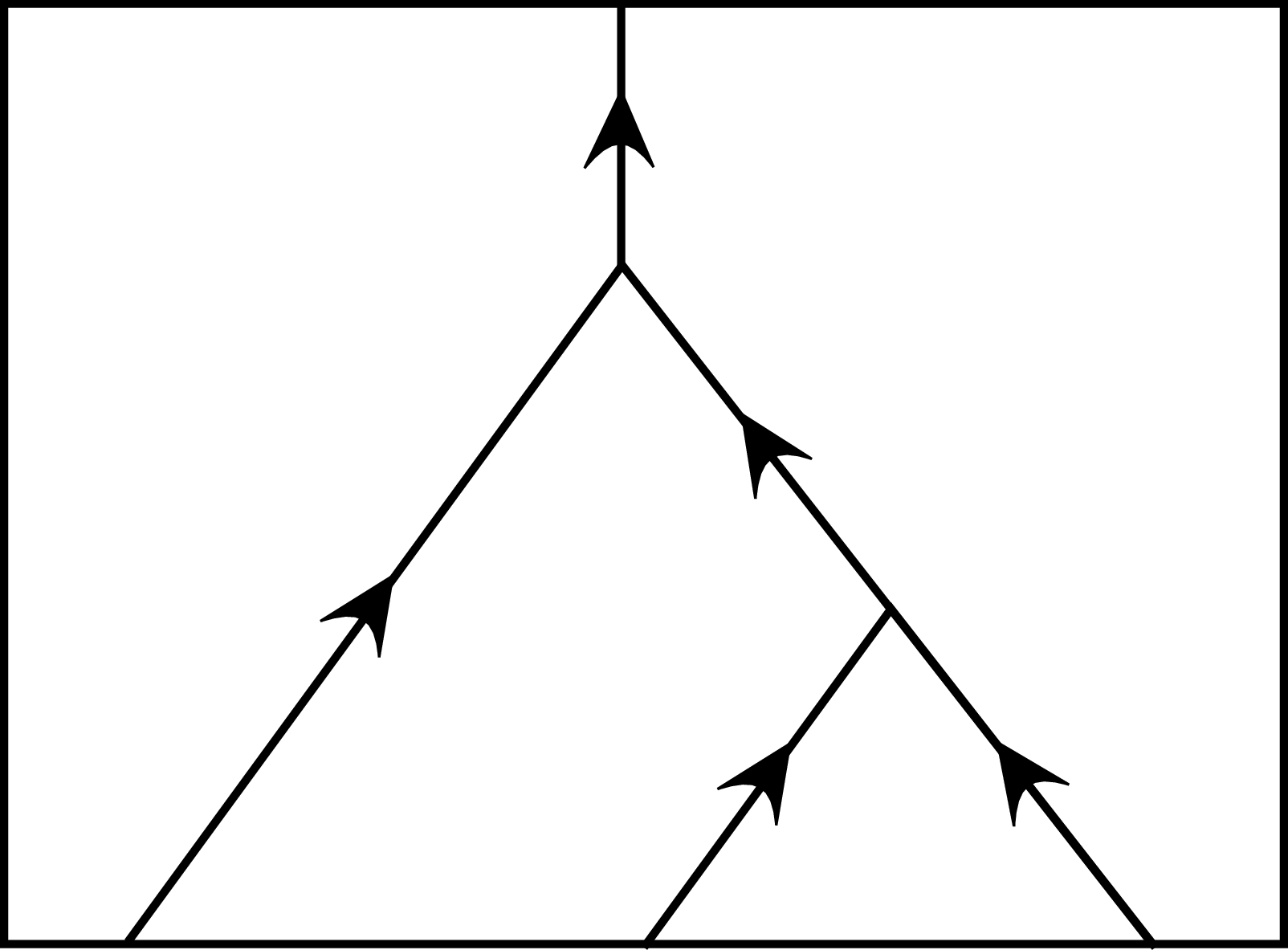} \endxy} =\vcenter{\xy (0,0)*{\def\svgscale{0.15}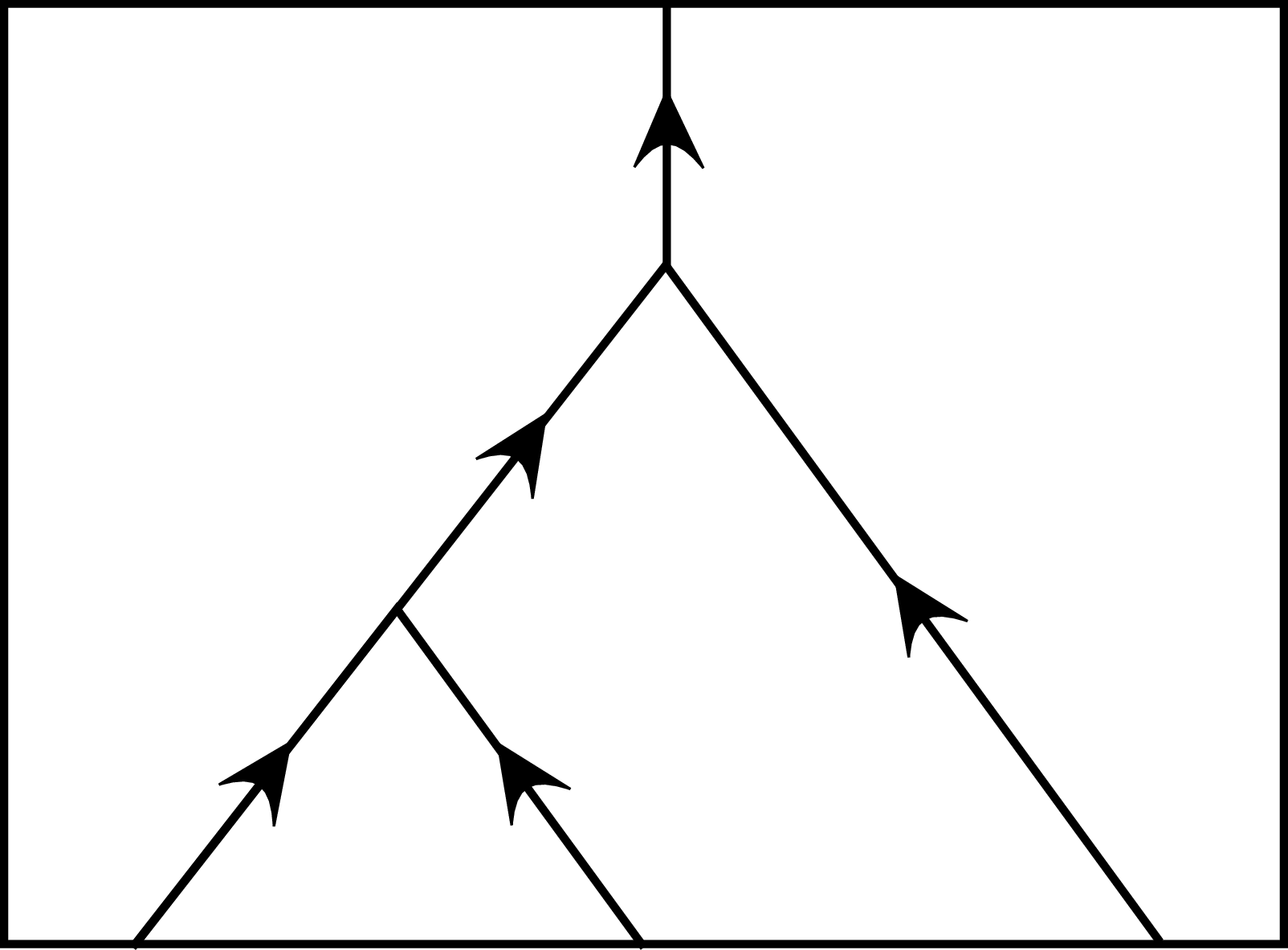} \endxy}.
\end{align}

So we desire the following equality in $\dmSBSBim$, for $k+l+m<n$:
\begin{align}\label{diag-image of I=H}
    \lambda_{k,l+m}\;\lambda_{l,m}\vcenter{\xy (0,0)*{\def\svgscale{0.16}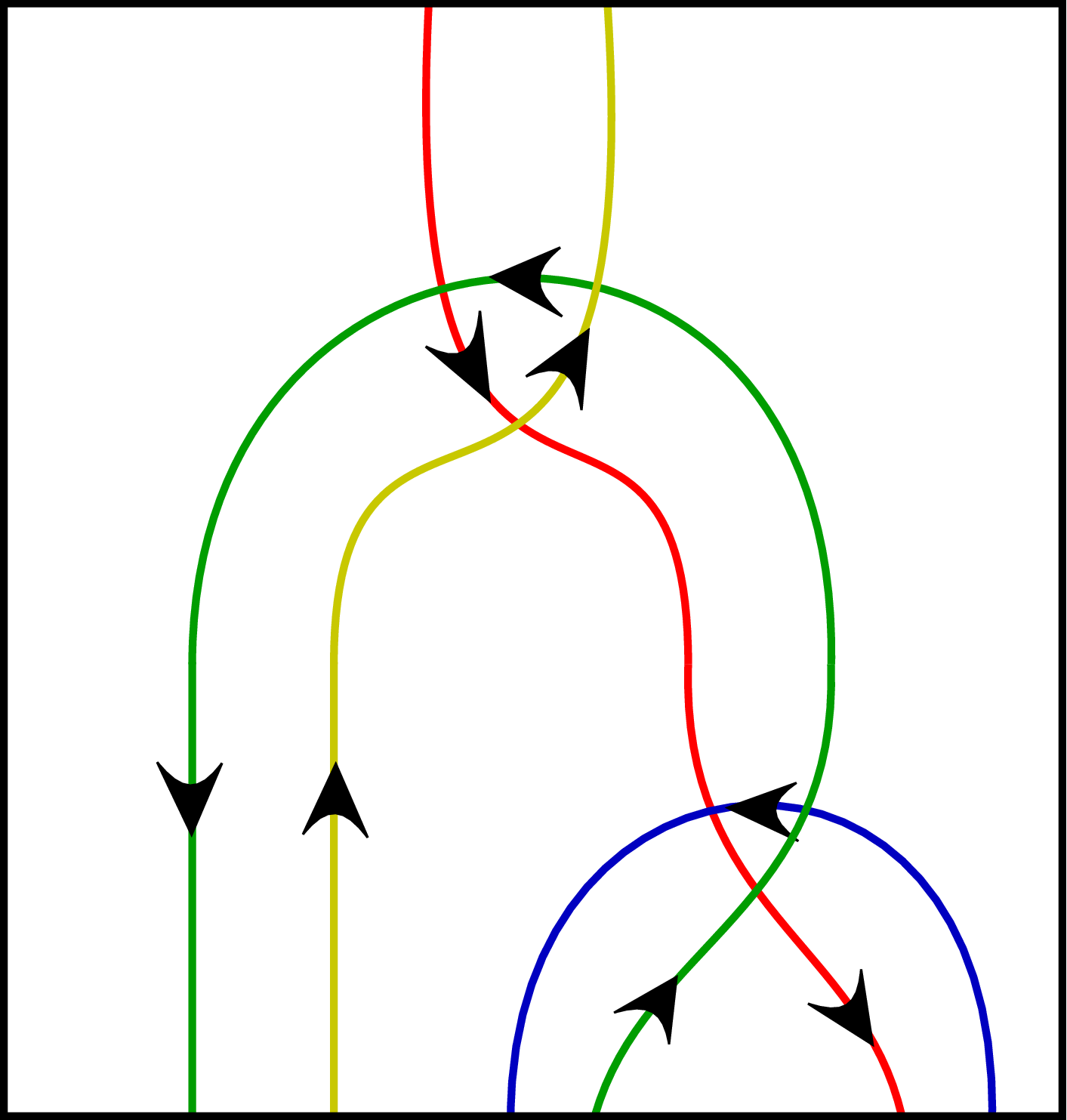} \endxy} =\lambda_{k,l}\;\lambda_{k+l,m}\vcenter{\xy (0,0)*{\def\svgscale{0.16}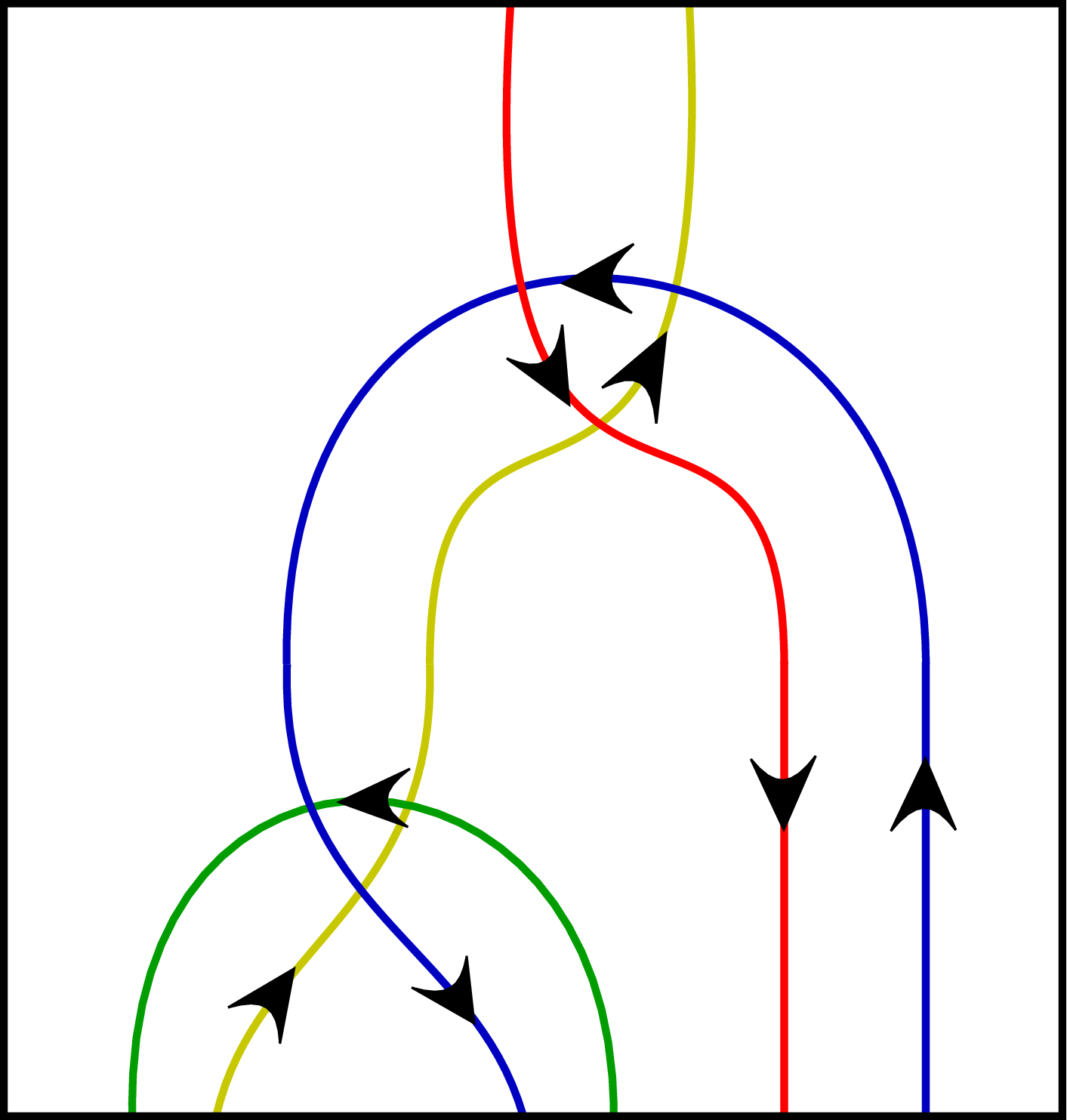} \endxy}
\end{align}\vspace{5pt}\\
and the following one for $k+l+m=n$:
\begin{align}\label{diag-image of I=H k+l+m=n}
   \lambda_k \lambda_{l,m} \vcenter{\xy (0,0)*{\def\svgscale{0.16}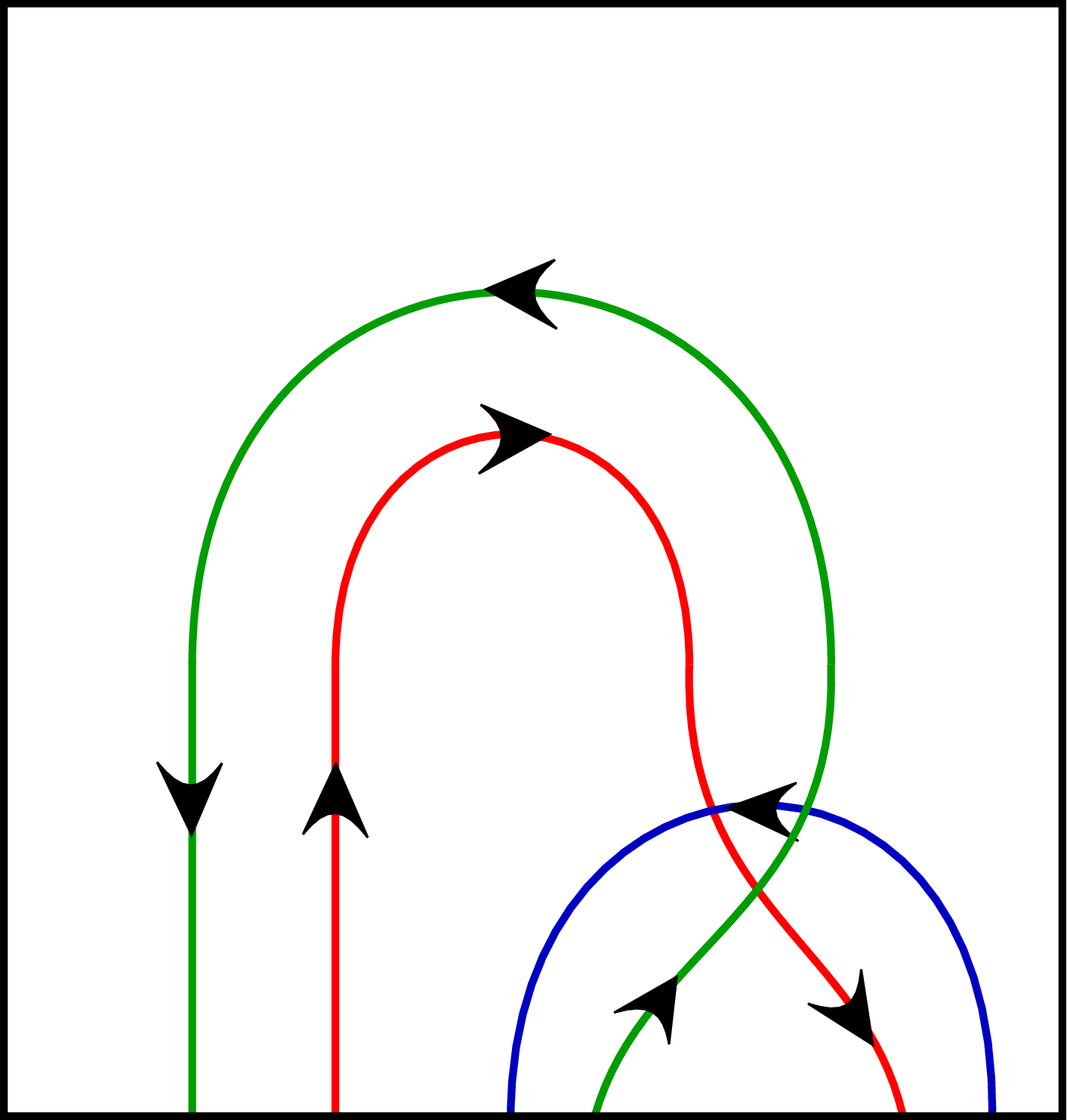} \endxy}= \lambda_{k,l} \lambda_{m} \vcenter{\xy (0,0)*{\def\svgscale{0.16}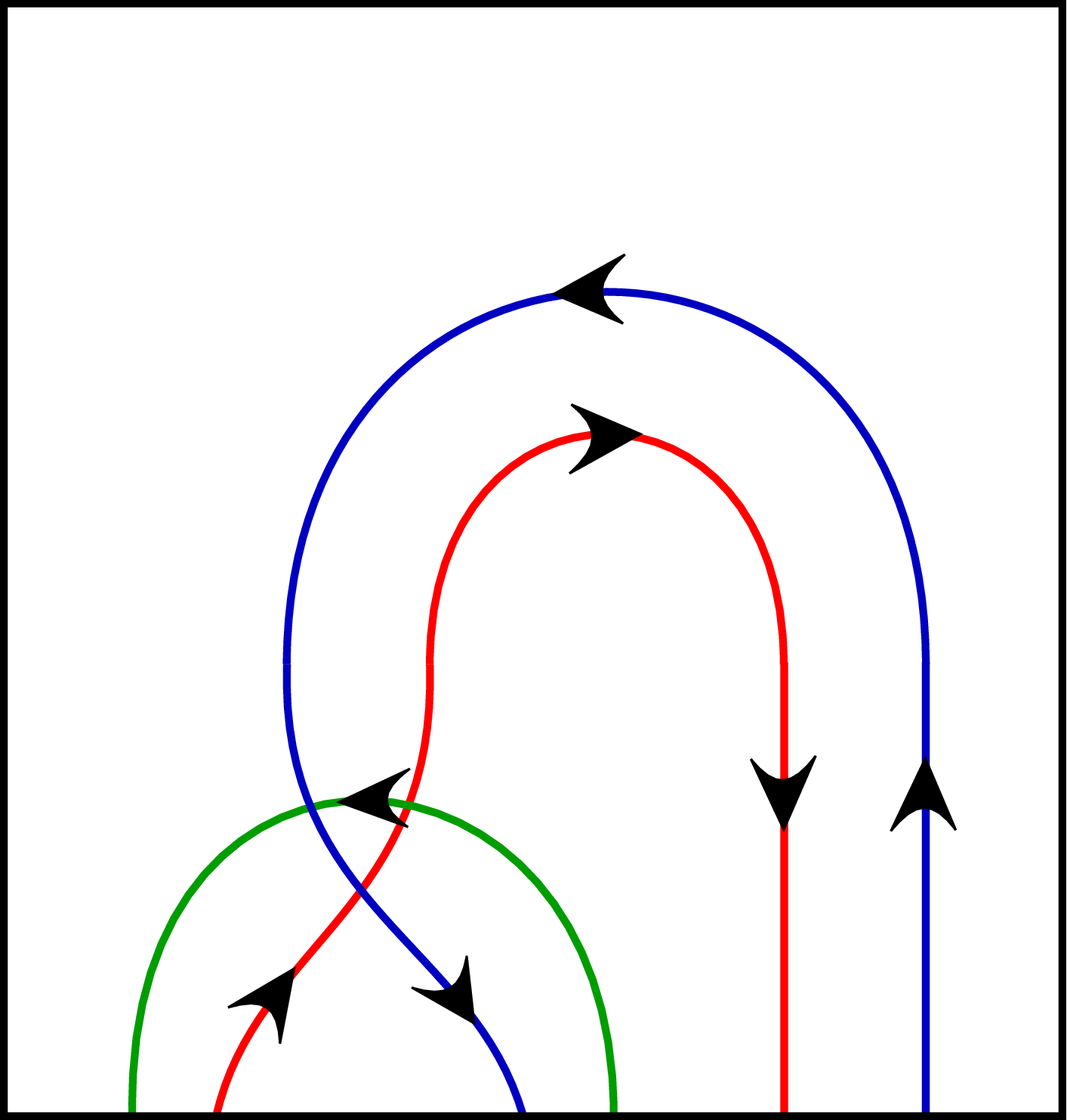} \endxy}.
\end{align}

If the diagrams (without coefficients) appearing on either side of $\cref{diag-image of I=H}, \cref{diag-image of I=H k+l+m=n}$ were equal, then the desired relation is equivalent to the equality of the coefficients, which is precisely \cref{cond-I=H relation} (with the convention in \cref{defn convention for coefficients of trivalent with labels n} for $k+l+m=n$ case).
For \cref{diag-image of I=H k+l+m=n}, the equality of the diagrams follows from isotopy. For \cref{diag-image of I=H}, we prove it in the following lemma.


\begin{lemma}\label{lemma associativity ssbim}
\[ \vcenter{\xy (0,0)*{\def\svgscale{0.16}\input{arxiv-figures/image_of_I=H_LHS_svg-tex.eps_tex}} \endxy}=\vcenter{\xy (0,0)*{\def\svgscale{0.16}\input{arxiv-figures/image_of_I=H_RHS_svg-tex.eps_tex}} \endxy}\] in $\pdmSBSBim$, hence in $\dmSBSBim$.    
\end{lemma}
\begin{proof}
By color symmetry, we may assume $a=0$. 
With this choice, red represents $0$, blue $m$, green $l+m$ and yellow $k+l+m$. Note that when red and green are absent, blue and yellow are in different connected components and vice versa.

Using the distant Reidemeister II relation, we may cross the blue and yellow strands within the red-green bigon, as shown below.
\begin{equation}
    \vcenter{\xy (0,0)*{\def\svgscale{0.15}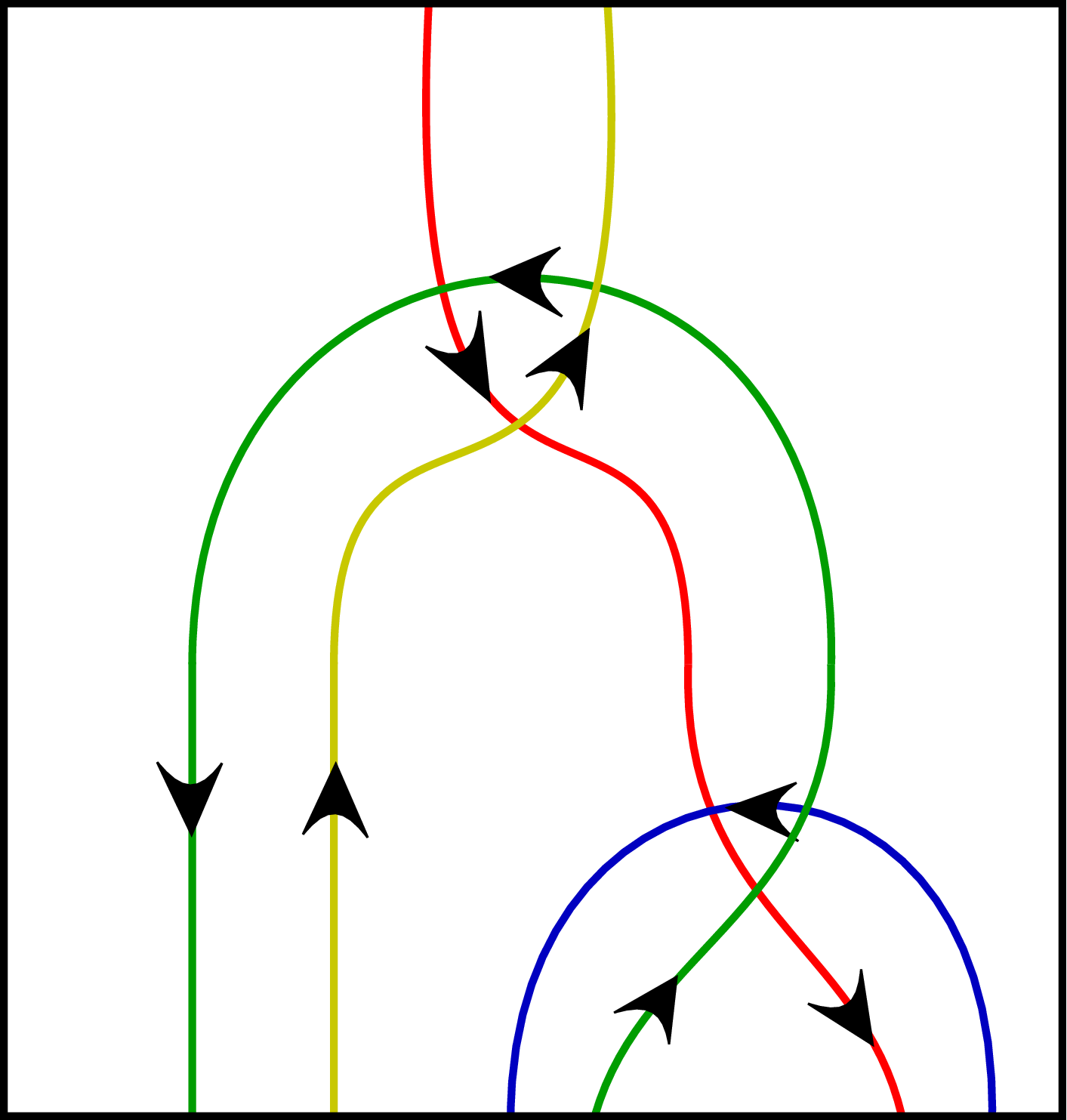} \endxy} = \vcenter{\xy (0,0)*{\def\svgscale{0.15}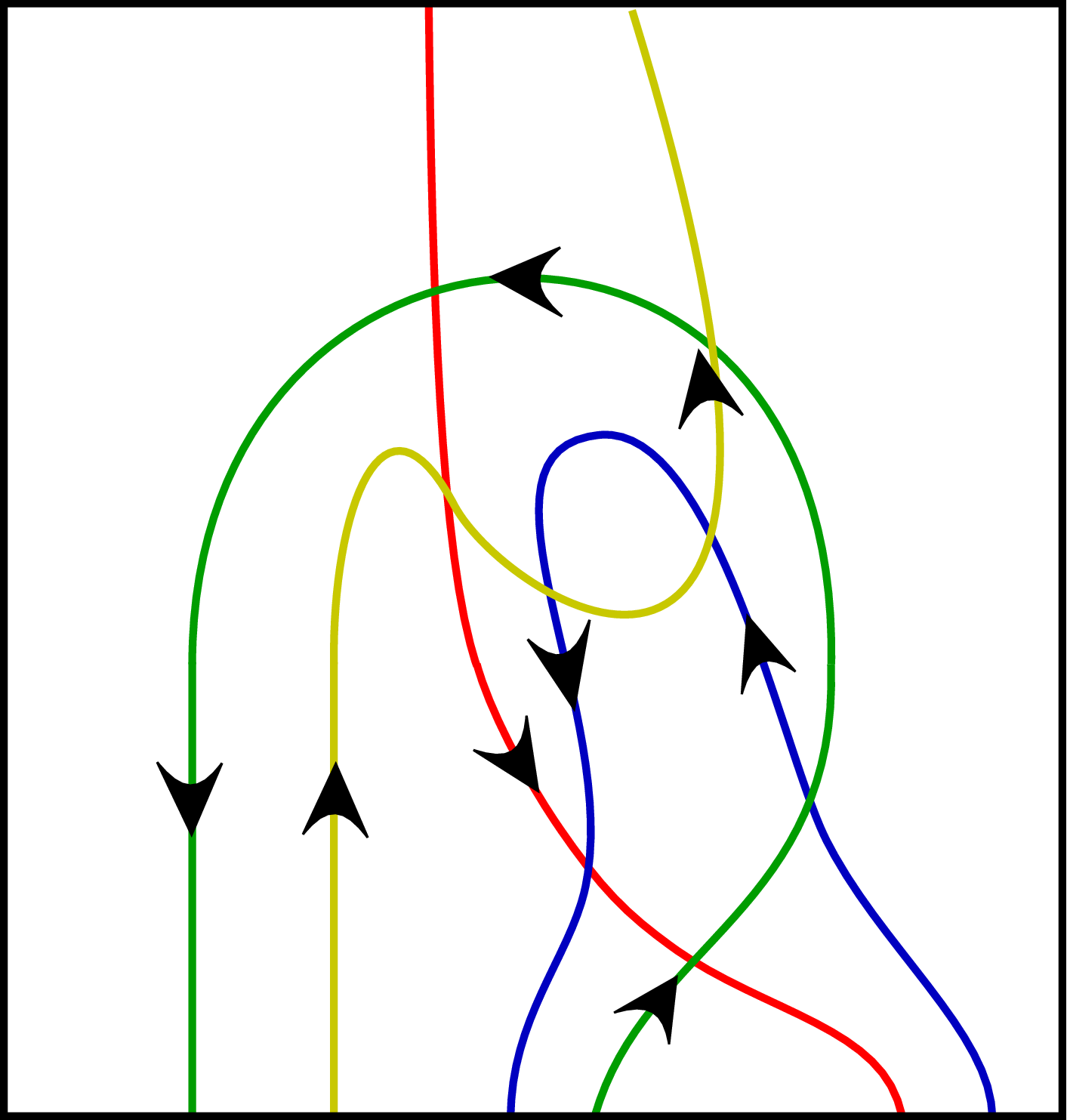} \endxy}
\end{equation}.

Then we do two oriented RIII relations to bring the yellow-blue crossings out of the red-green bigon.
\begin{equation}
 \vcenter{\xy (0,0)*{\def\svgscale{0.15}\input{arxiv-figures/associativity_ssbim_1_svg-tex.eps_tex}} \endxy} = \vcenter{\xy (0,0)*{\def\svgscale{0.15}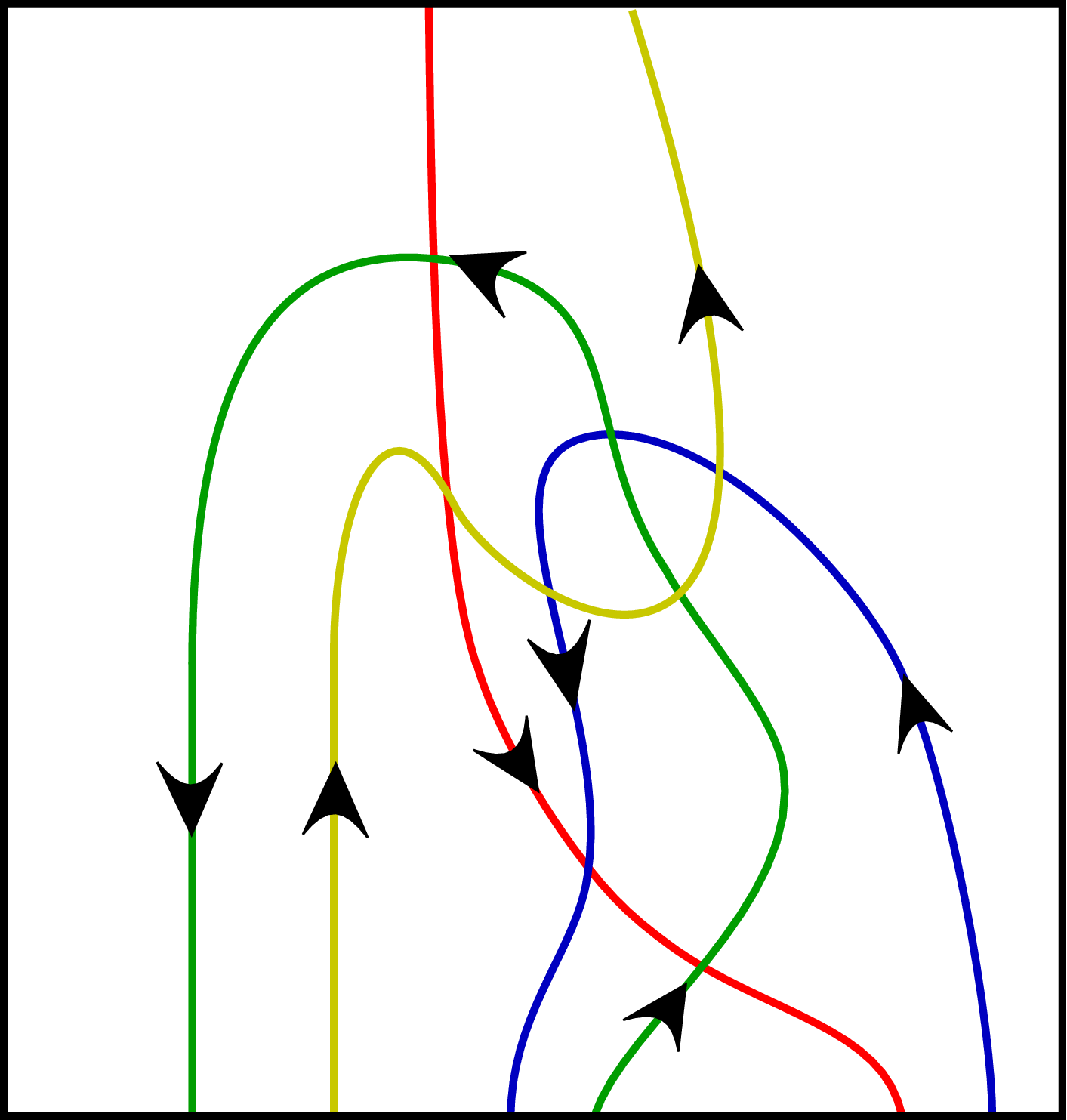} \endxy}=   \vcenter{\xy (0,0)*{\def\svgscale{0.15}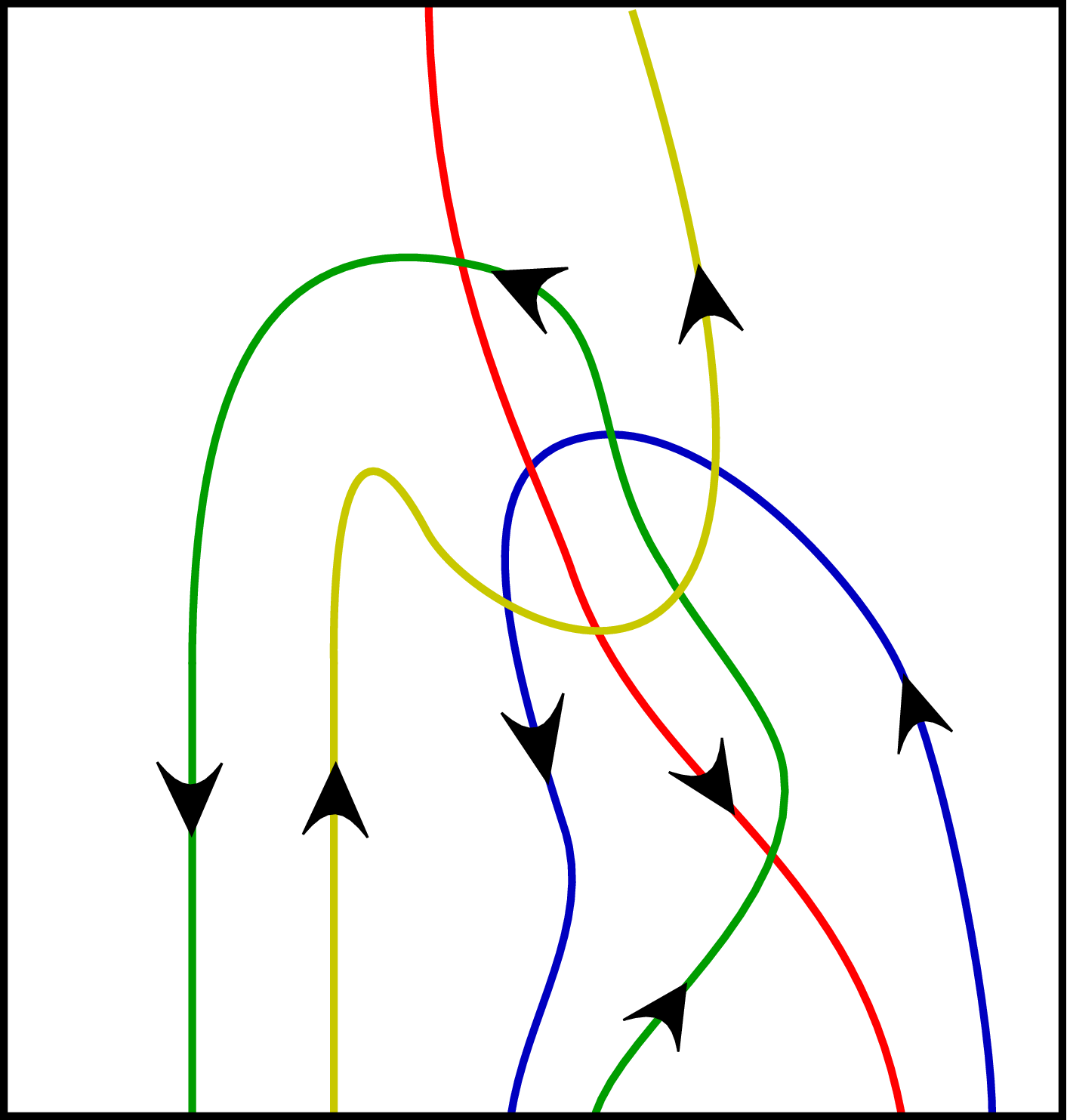} \endxy}
\end{equation}
Now two more oriented RIII will bring the green-red crossings into the blue-yellow bigon.
\begin{equation}
    \vcenter{\xy (0,0)*{\def\svgscale{0.15}\input{arxiv-figures/associativity_ssbim_3_svg-tex.eps_tex}} \endxy}=\vcenter{\xy (0,0)*{\def\svgscale{0.15}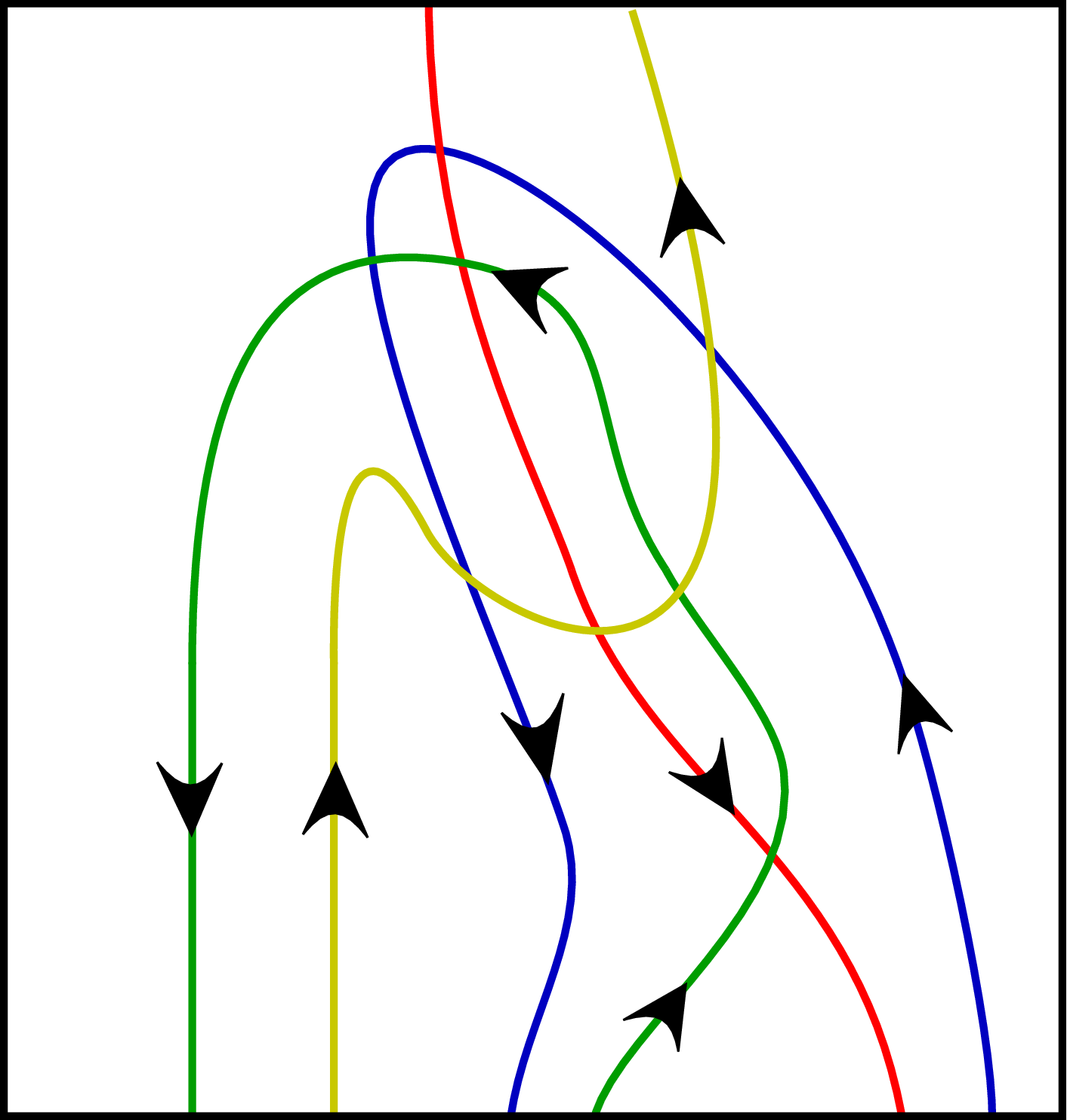} \endxy}=\vcenter{\xy (0,0)*{\def\svgscale{0.15}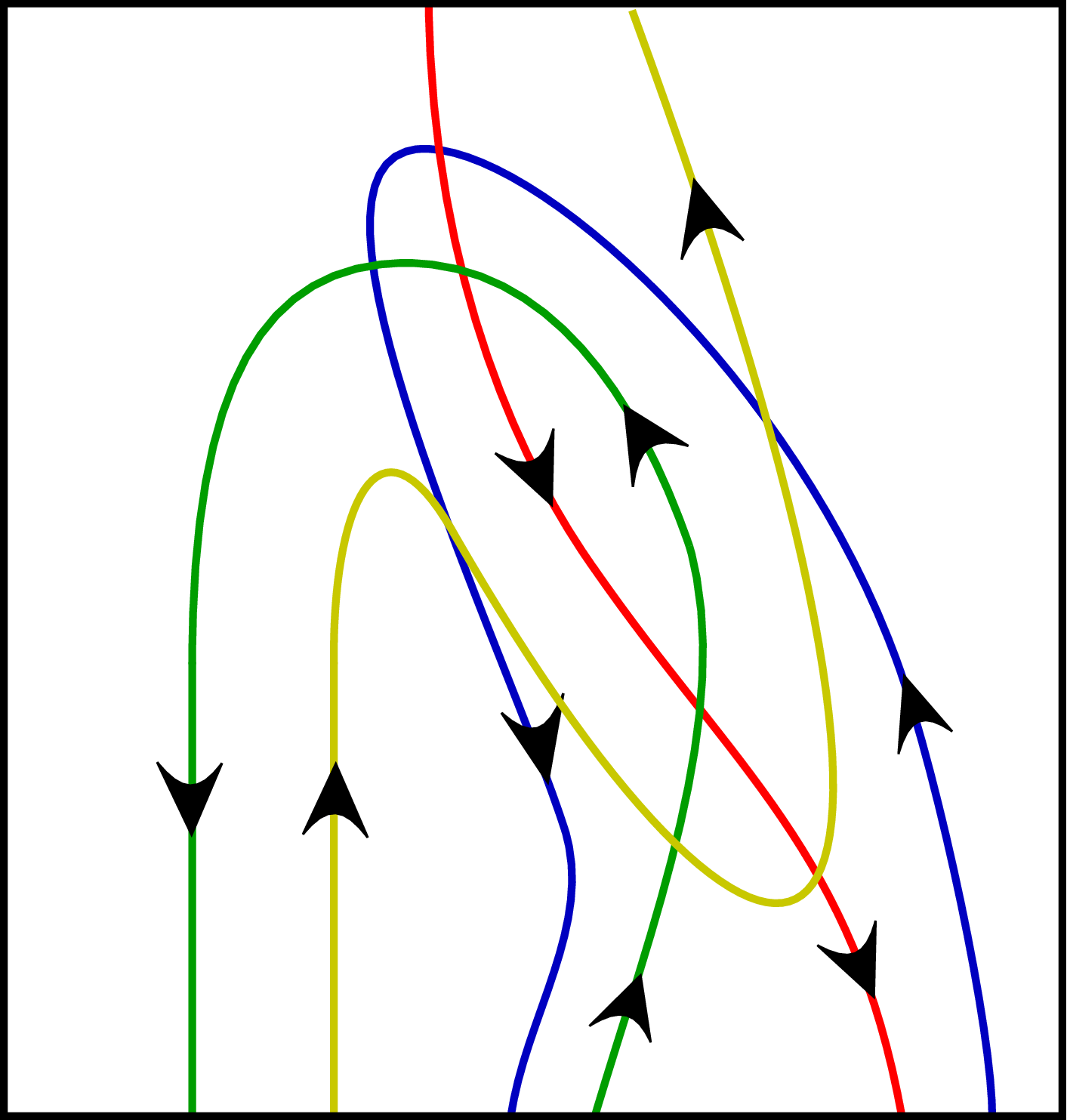} \endxy}
\end{equation}

Finally, we do a distant RII to uncross the green and red strands within the blue-yellow bigon.
\begin{equation}
  \vcenter{\xy (0,0)*{\def\svgscale{0.15}\input{arxiv-figures/associativity_ssbim_5_svg-tex.eps_tex}} \endxy}=\vcenter{\xy (0,0)*{\def\svgscale{0.15}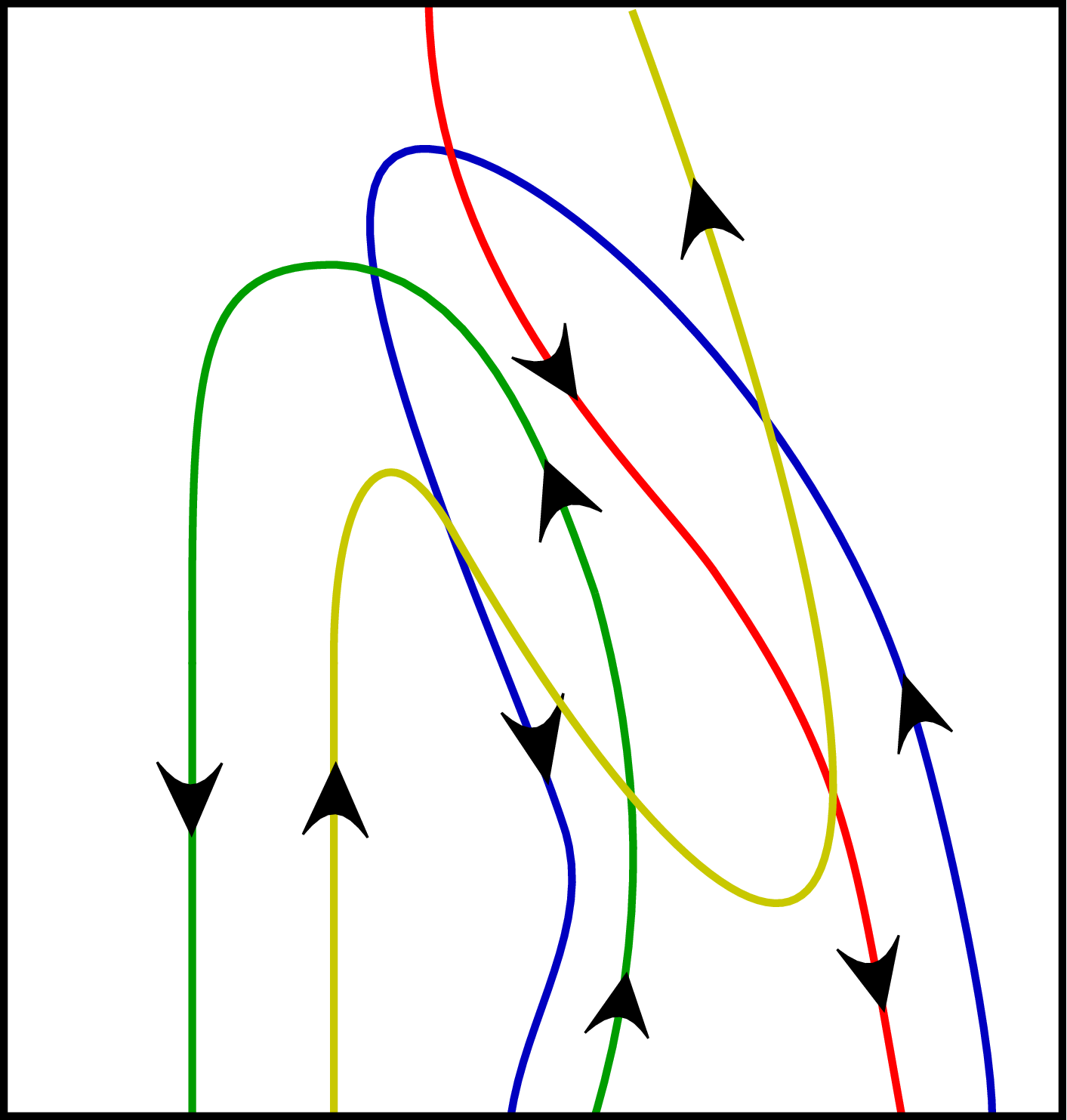} \endxy}.  
\end{equation}
The last term above is the RHS in \cref{lemma associativity ssbim} and we are done.
\end{proof}

Similarly, the relation obtained by taking the arrow reversals of \cref{diag-I=H relation} implies $\nu_{k, l+m}\;\nu_{l,m}=\nu_{k+l,m}\;\nu_{k,l}$.

\subsection{Relations and symmetries}

Above we checked specific versions of tag cancellation and $I=H$, and discussed some (but not all) of their mirror images and arrow reversals. This was enough to force the constraints on the scalars $\lambda$ and $\nu$ which were needed for \cref{lemma symmetry of GS-zeta}. Thanks to this lemma we do not need to check rotations, color changes via $\sigma$, or versions of relations obtained via $\Upsilon$.

For the remaining relations in $\cwebs$, namely \cref{diag-digon bursting} and \cref{diag-square flop} below, we only check one version of the relation. 

However, we must also verify the mirror images (resp. arrow reversals) of all the relations in $\cwebs$, and this does not correspond to a symmetry in $\DiagSBS$. Thankfully, the remaining mirror images and arrow reversals can be 
obtained alternatively from ones we have checked using only rotation and $\Upsilon$.

\subsection{Bursting a digon}\label{subsec-bigon bursting}
For valid region labels, we have the following ``bursting a digon" relation in $\cwebs$ (up to mirror images and arrow reversals).
\begin{align}\label{diag-digon bursting}{
        \tikz[x=1mm, y=1mm, baseline=-0.5ex] {
          \draw[->-] (0,-8) -- (0,-4);
          \draw[->-] (0,-4) arc (270:90:4);
          \draw[->-] (0,-4) arc (-90:90:4) node[pos=0.5,right]{$1$};
          \draw[->] (0,4) -- (0,8) node [right] {$k$};
        }
        = [k] \: \tikz[x=1mm, y=1mm, baseline=-0.5ex] {
          \draw[->] (0,-8) -- (0,8) node [right] {$k$};
        }
      }
      .\end{align}

For $k \neq 1$ or $n$, the LHS of the above diagram with a region label of $a$ on the right is sent to the following diagram via $\GS_{\zeta}$: 
\begin{align}\label{diag-image of bigon}
    \lambda_{k-1,1}\;\nu_{k-1,1}\vcenter{\xy (0,0)*{\def\svgscale{0.15}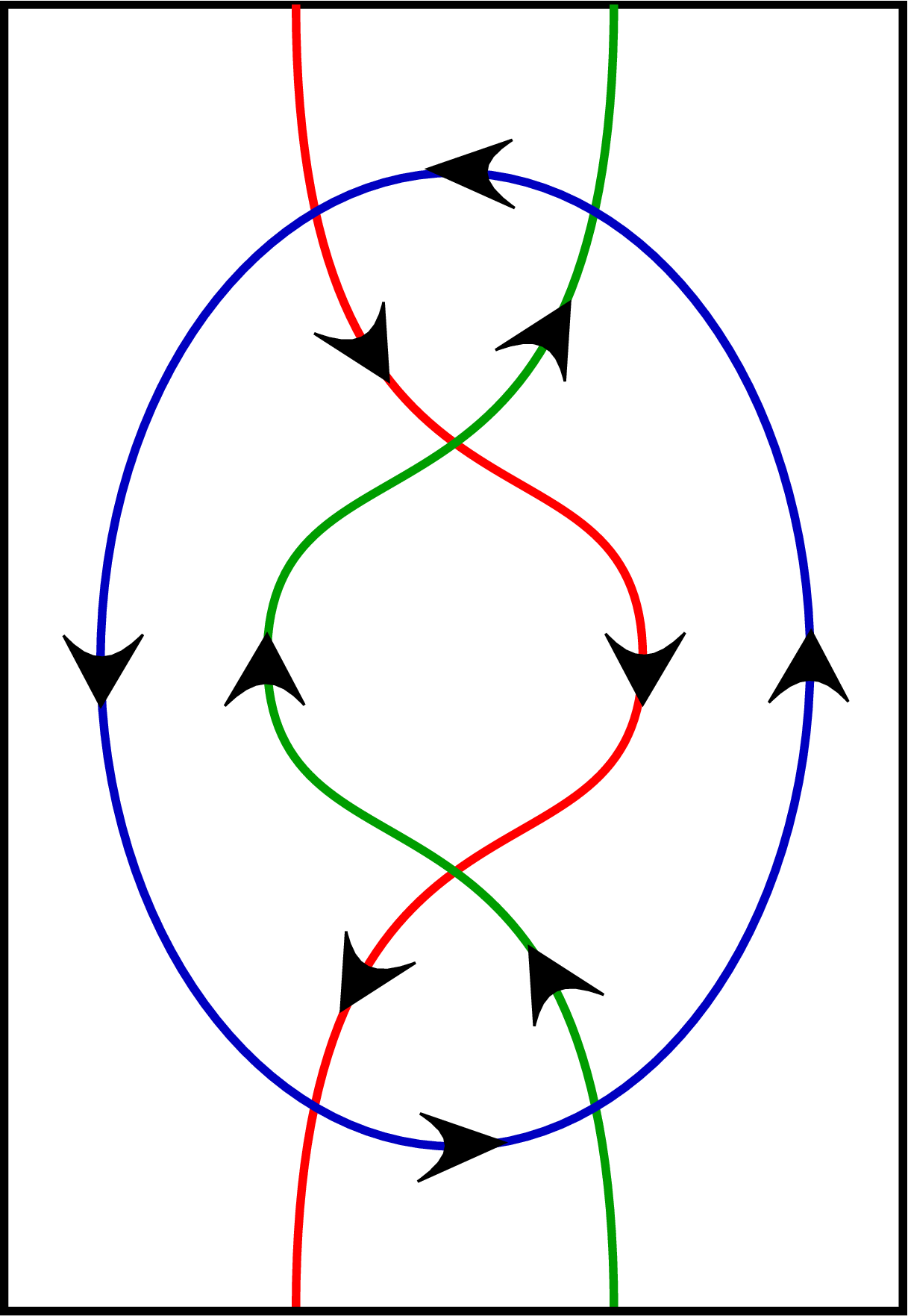} \endxy}\;.
\end{align}

By color symmetry, we may assume that $a=0$. 
Using \cref{R2unoriented2}, \cref{R2oriented} and \cref{ccwcirc}, we can simplify the diagram in \cref{diag-image of bigon} (without the coefficients) to 
\begin{align}\label{eq-bigon burst simplification}
    \vcenter{\xy (0,0)*{\def\svgscale{0.15}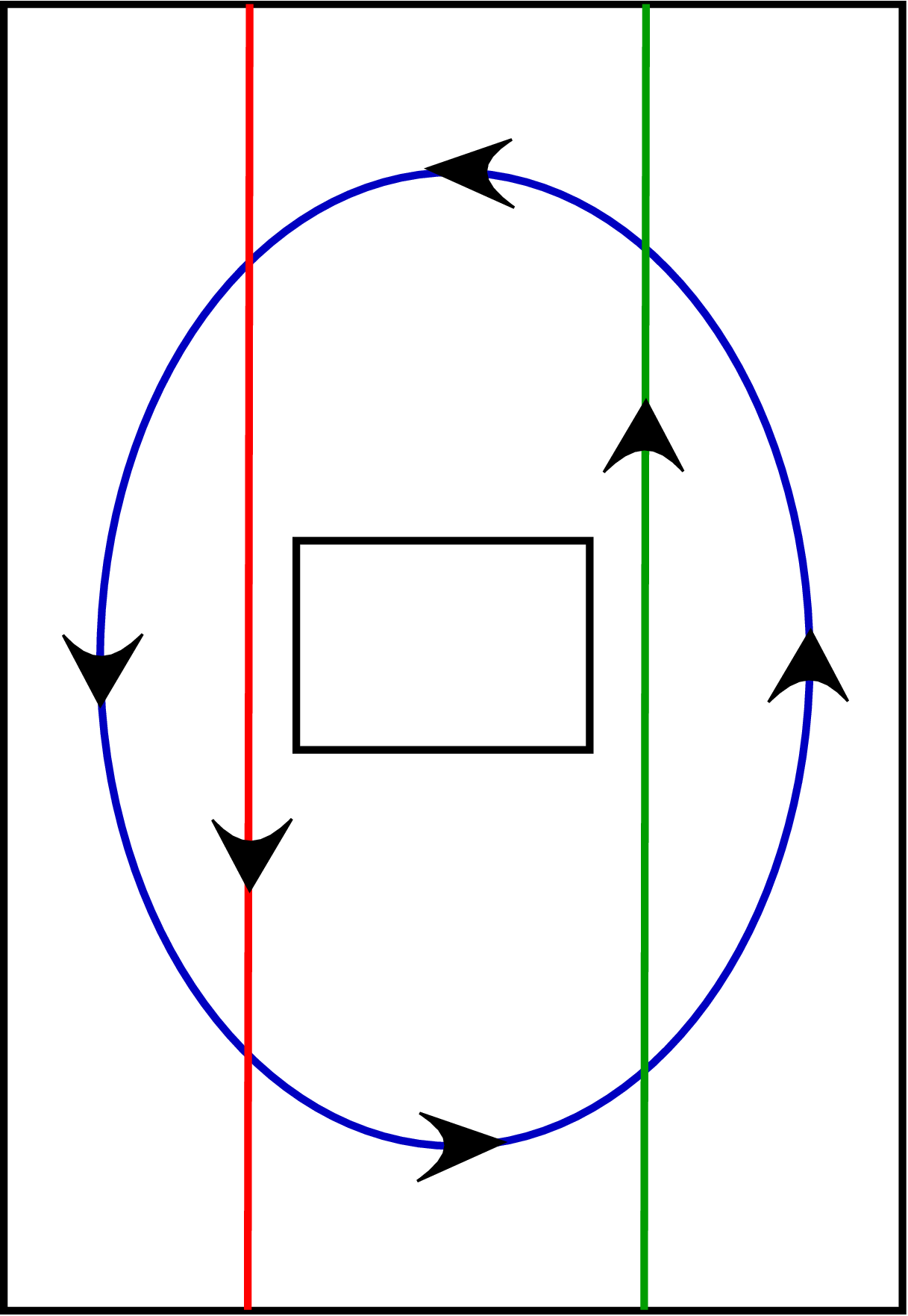} \endxy}=\vcenter{\xy (0,0)*{\def\svgscale{0.15}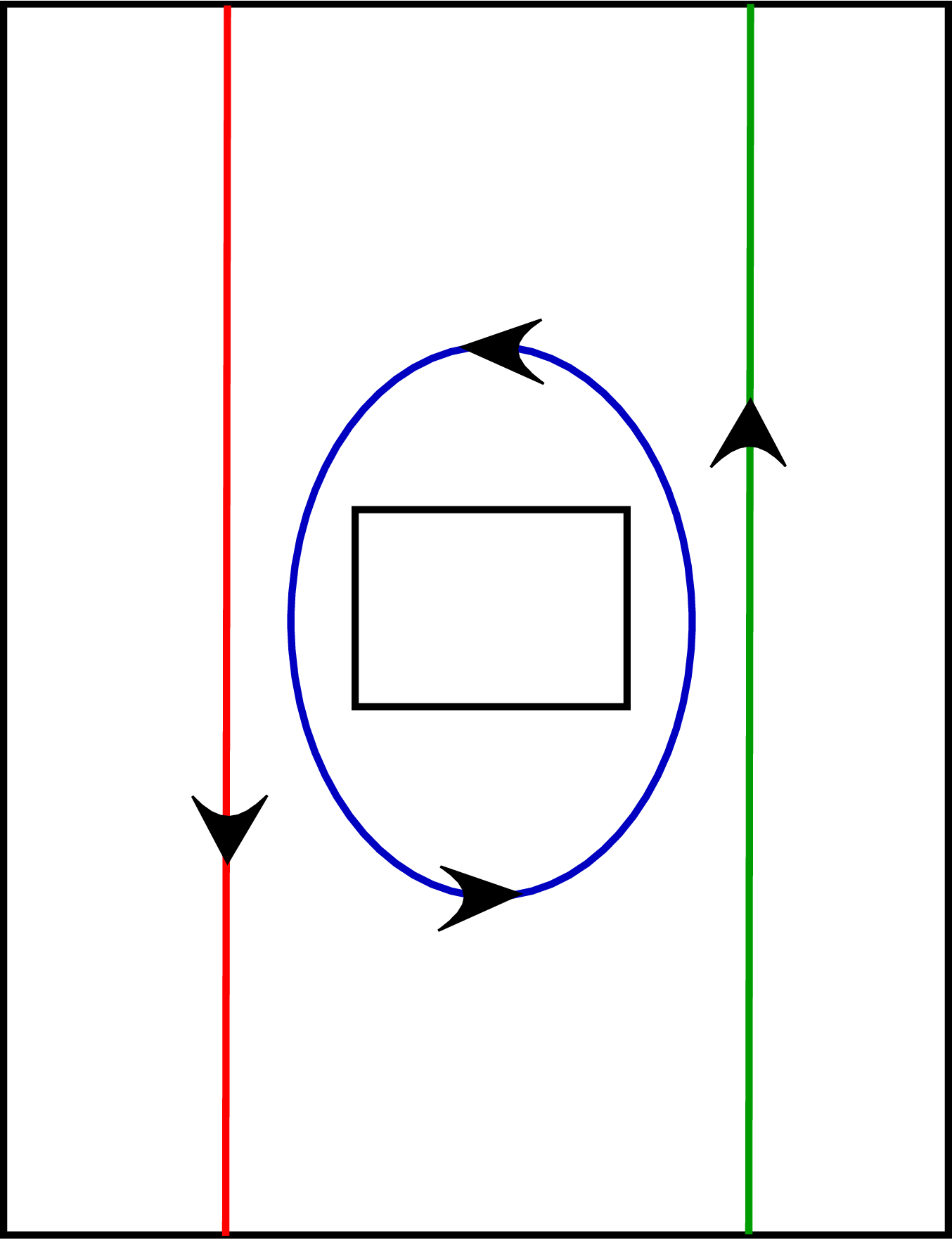} \endxy} = \vcenter{\xy (0,0)*{\def\svgscale{0.15}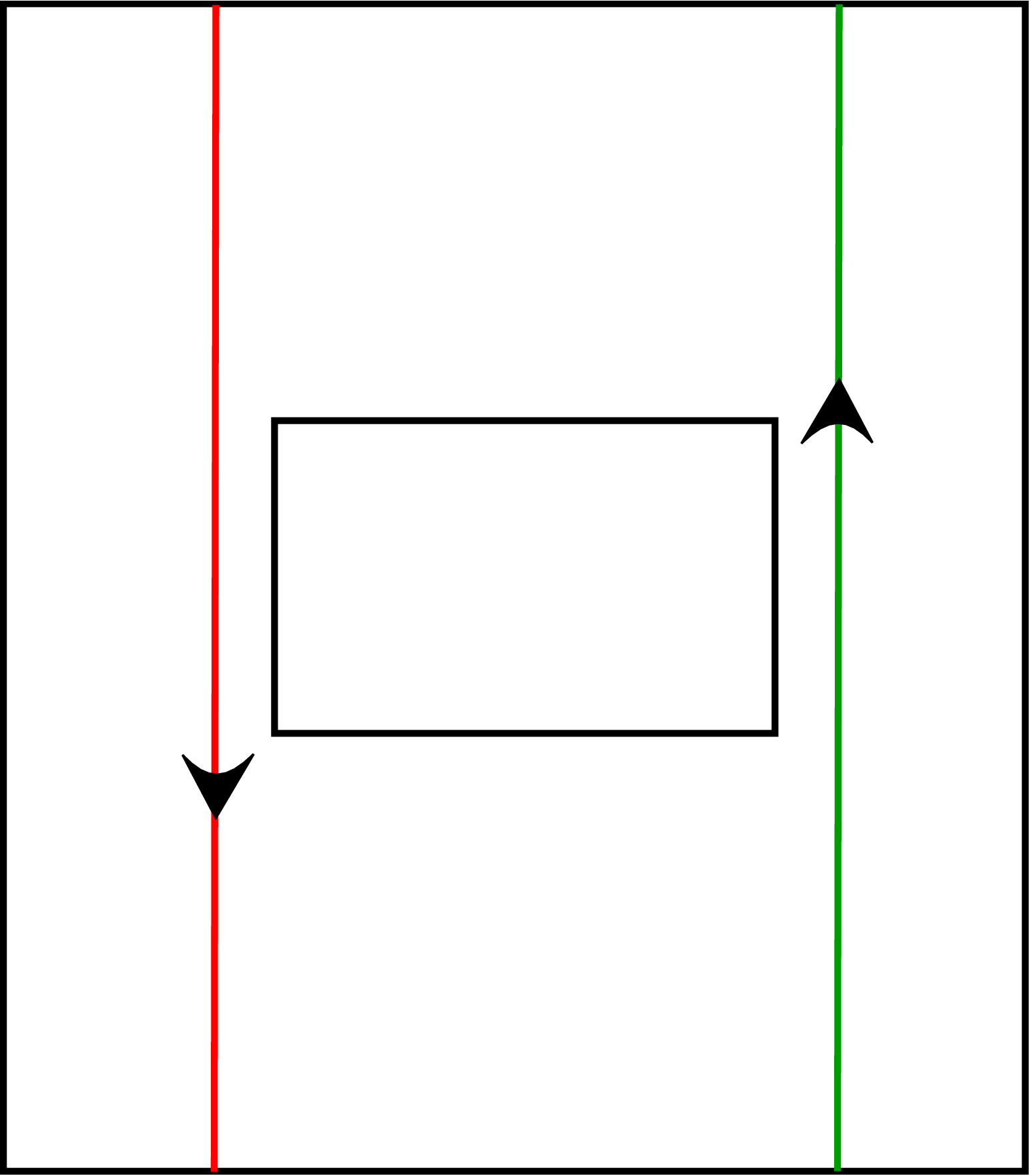} \endxy}.
\end{align}
By \cref{corollary- bigon bursting}, we know that $\dd_{\hh{0,k}}^{\hh{0,1,k}}\left(\mu_{\hh{1}}^{\hh{0,1},\hh{1,k}}\right)=(-1)^{k-1}\zeta^{-k(k-1)/2}q^{-(k-1)}[k]_q$.
The above discussion then shows that $\GS_{\zeta}$ sends the LHS of \cref{diag-digon bursting} to 
\[(-1)^{k-1}\zeta^{-k(k-1)/2}q^{-(k-1)}[k]_q\; \lambda_{k-1,1} \; \nu_{k-1,1} \vcenter{\xy (0,0)*{\def\svgscale{0.15}\input{arxiv-figures/image_of_tag_c_svg-tex.eps_tex}} \endxy}.\]
Comparing this with the image of the RHS of \cref{diag-digon bursting}, we have the conditions $$ \lambda_{k-1,1} \; \nu_{k-1,1}=(-1)^{k-1}\zeta^{k(k-1)/2}q^{k-1},$$ for $k\neq 1,n$ which is precisely \cref{cond-bigon killing}.

Now we consider the edge cases. 
For $k=1$, there is nothing to prove; for $k=n$, we use tag cancellation \cref{diag-tag cancellation}, and see that the resulting equality is true in $\dmSBSBim$ as follows.\\
Observe that $\GS_{\zeta}$ sends the the LHS to
\begin{align}\label{diag- bigon killing k=n}
  (-1)^{n-1}  \vcenter{\xy (0,0)*{\def\svgscale{0.1}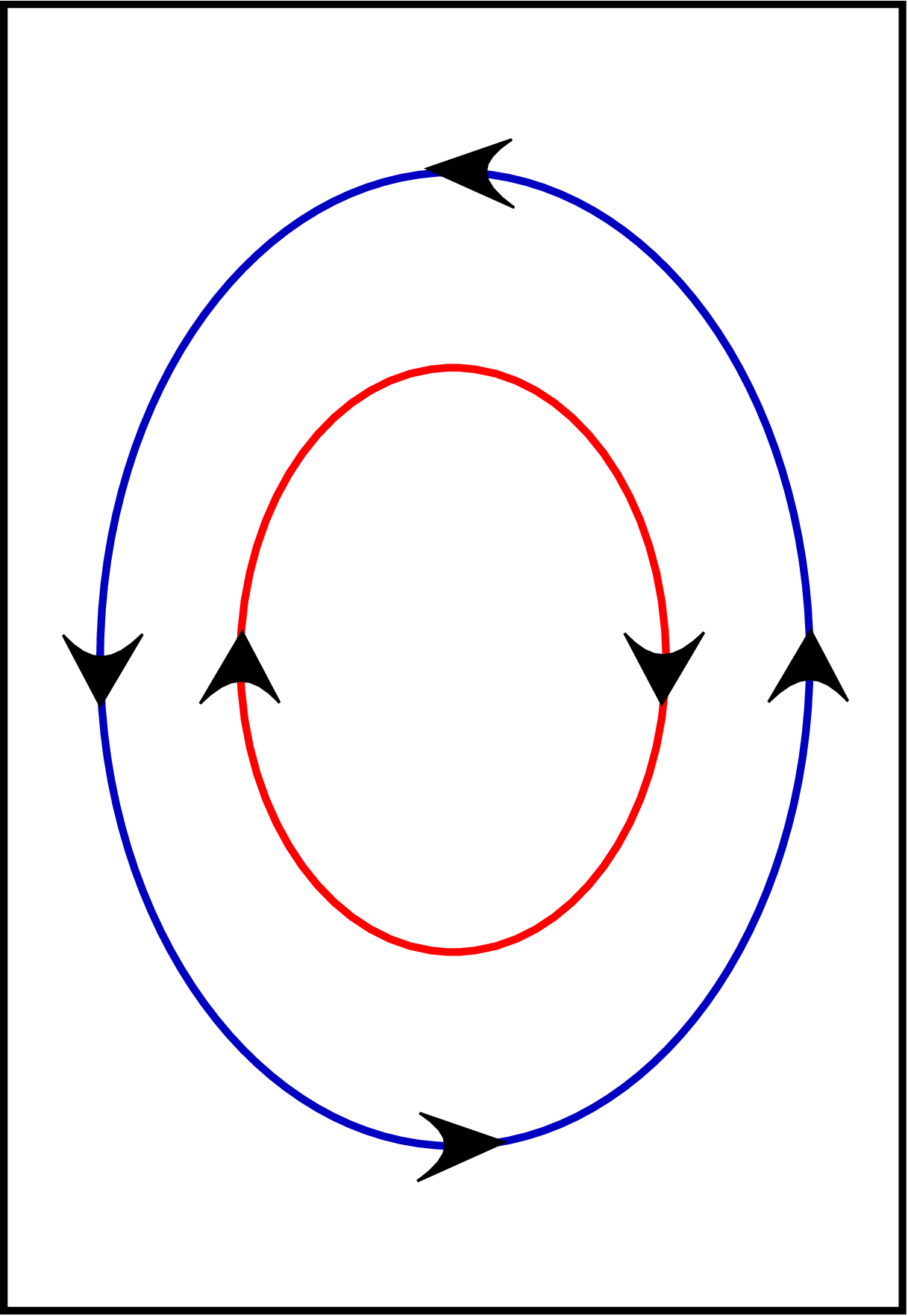} \endxy}.
\end{align}
By color symmetry, we may assume $a=0$. 
Using relations \cref{cccirc}, \cref{ccwcirc} we can simplify \cref{diag- bigon killing k=n} to $(-1)^{n-1} \dd_{\hh{0}}^{\hh{0,1}}\left(\mu^{\hh{0,1}}_{\hh{1}} \right)$.
From \cref{lemma demazure adjacent}, we have that $\dd_{\hh{0}}^{\hh{0,1}}\left(\mu^{\hh{0,1}}_{\hh{1}} \right)=(-1)^{n-1}[n]_q$, so that \cref{diag- bigon killing k=n} further reduces to $[n]_q$ as desired.

\subsection{Bursting a Square}\label{subsec-square flop}
For valid region labels, we had the following``bursting a square" relation in $\cwebs$:
\begin{align}\label{diag-square flop}
    {
      \phantom{ \tikz[x=1mm, y=1mm, baseline=-0.5ex] {\draw (0,-10) -- (0,10)}} 
      \tikz[x=1mm, y=1mm, baseline=-0.5ex] {
        \draw[->-] (45:4) arc (45:135:4) node [pos=0.5, above]{$1$};
        \draw[-<-] (135:4) arc (135:225:4);
        \draw[->-] (225:4) arc (225:315:4) node [pos=0.5, below]{$1$};
        \draw[->-] (315:4) arc (315:405:4);
        \draw[->-] (7,-7) node [right] {$1$} -- (315:4);
        \draw[->] (45:4) -- (7,7) node [right] {$1$};
        \draw[->] (135:4) -- (-7,7) node [left] {$k$};
        \draw[->-] (-7,-7) node [left] {$k$} -- (225:4);
      }
      =
      \tikz[x=1mm, y=1mm, baseline=-0.5ex] {
        \draw[->-] (-7,-7) node [left] {$k$}
          .. controls (-6,-6) and (-3,-5) .. (0,-4);
        \draw[->-] (7,-7) node [right] {$1$}
          .. controls (6,-6) and (3,-5) .. (0,-4);
        \draw[->-] (0,-4) -- (0,4);
        \draw[->] (0,4) .. controls (-3,5) and (-6,6)
          .. (-7,7) node [left] {$k$};
        \draw[->] (0,4) .. controls (3,5) and (6,6)
          .. (7,7) node [right] {$1$};
      }
      + [k - 1] \tikz[x=1mm, y=1mm, baseline=-0.5ex] {
        \draw[->] (-7,-7) .. controls (-4,-4) and (-4,4)
          .. (-7,7) node [left] {$k$};
        \draw[->] (7,-7) .. controls (4,-4) and (4,4)
          .. (7,7) node [right] {$1$};
      }}.
\end{align}

The edge cases $k=0,1$ are trivial, and we treat the edge cases $k \in \{n-1,n\}$ afterwards, so assume $1 < k < n-1$. Applying $\GS_{\zeta}$, we then desire the following equality in $\dmSBSBim$: 
\begin{align}\label{diag-ssbim square flop}
   \lambda_{k-1,1} \nu_{k-1,1} \lambda_{1,1} \nu_{1,1} \vcenter{\xy (0,0)*{\def\svgscale{0.15}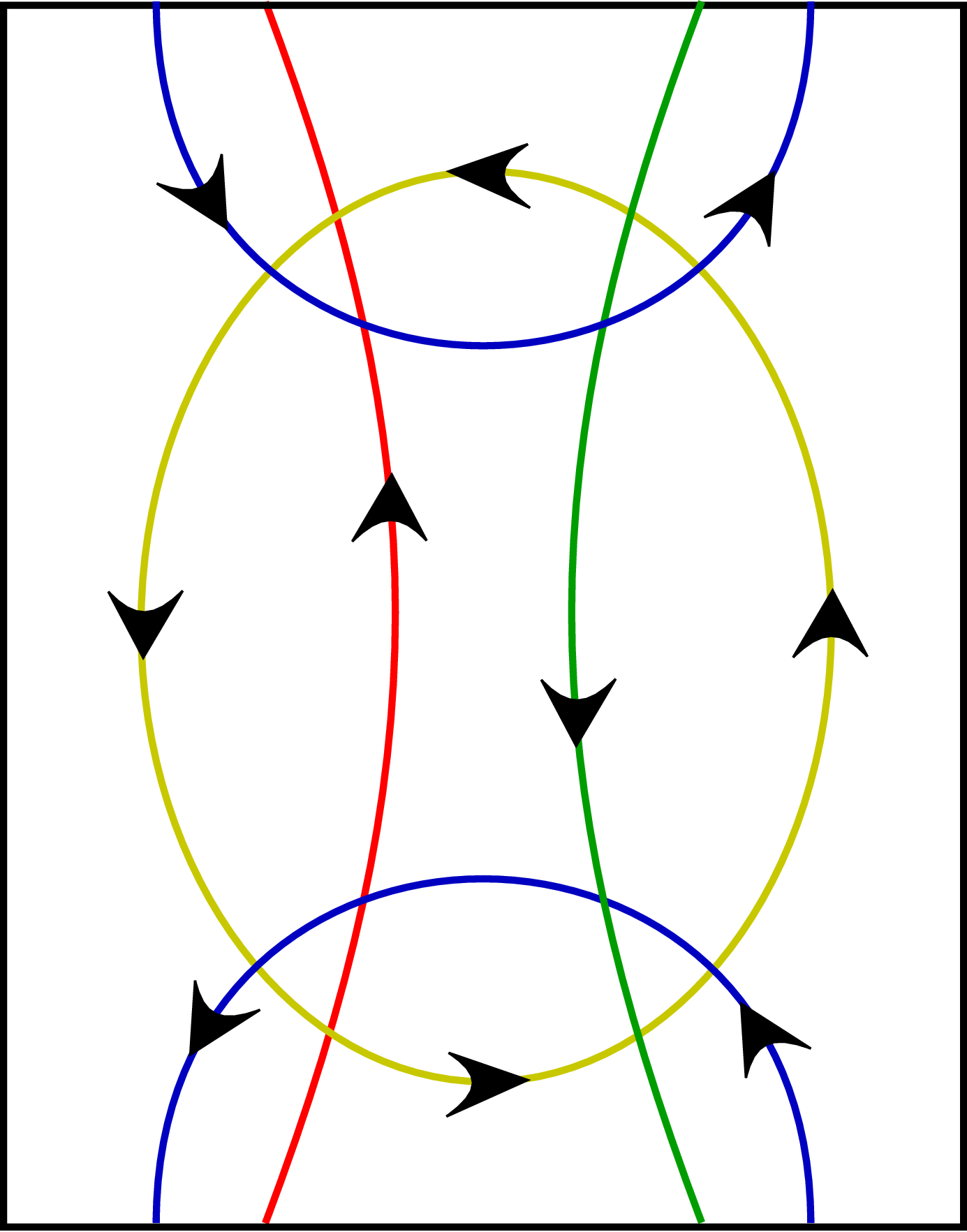} \endxy} = \lambda_{k,1}\nu_{k,1}\vcenter{\xy (0,0)*{\def\svgscale{0.15}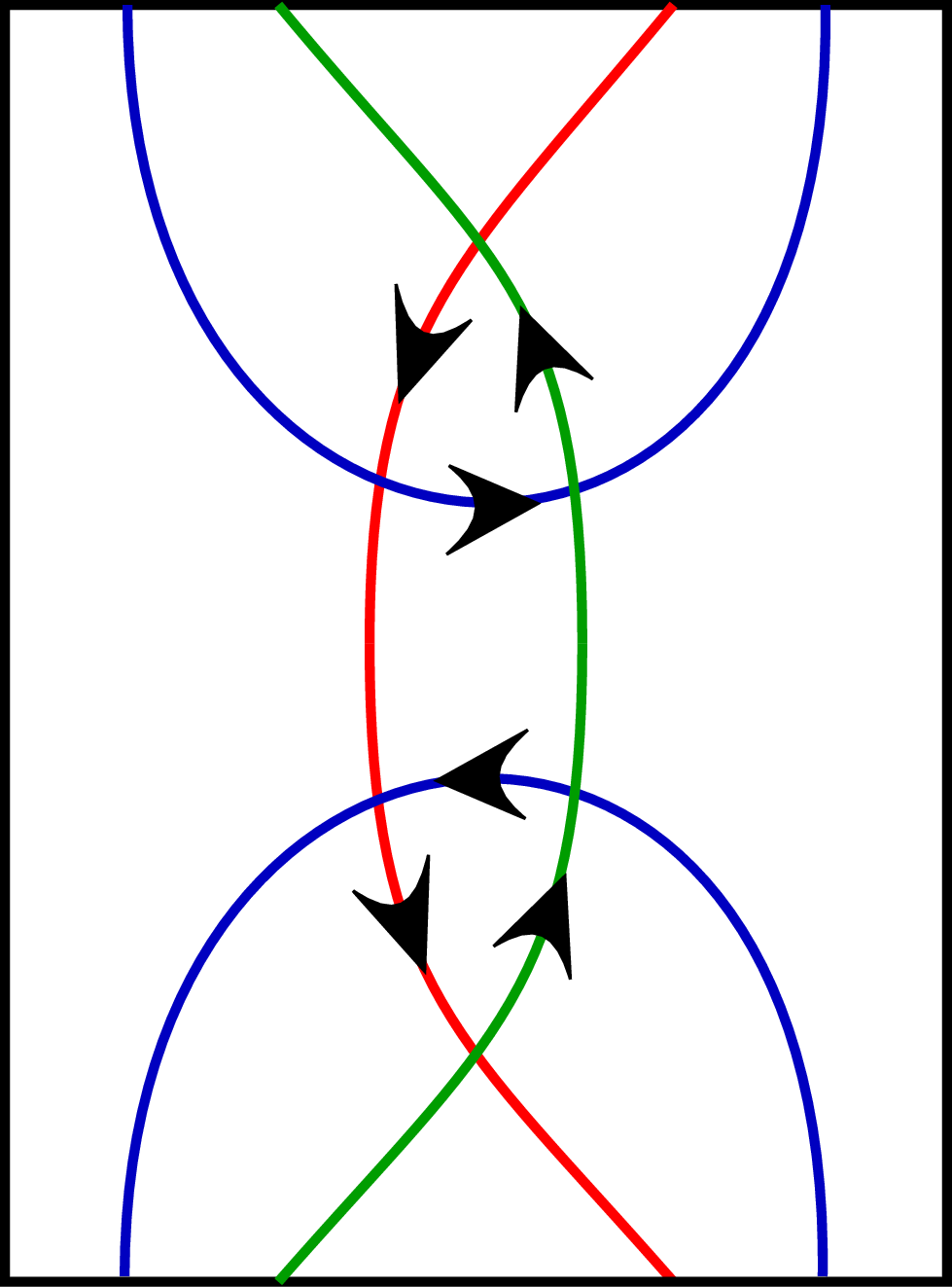} \endxy} + [k-1]_q\vcenter{\xy (0,0)*{\def\svgscale{0.15}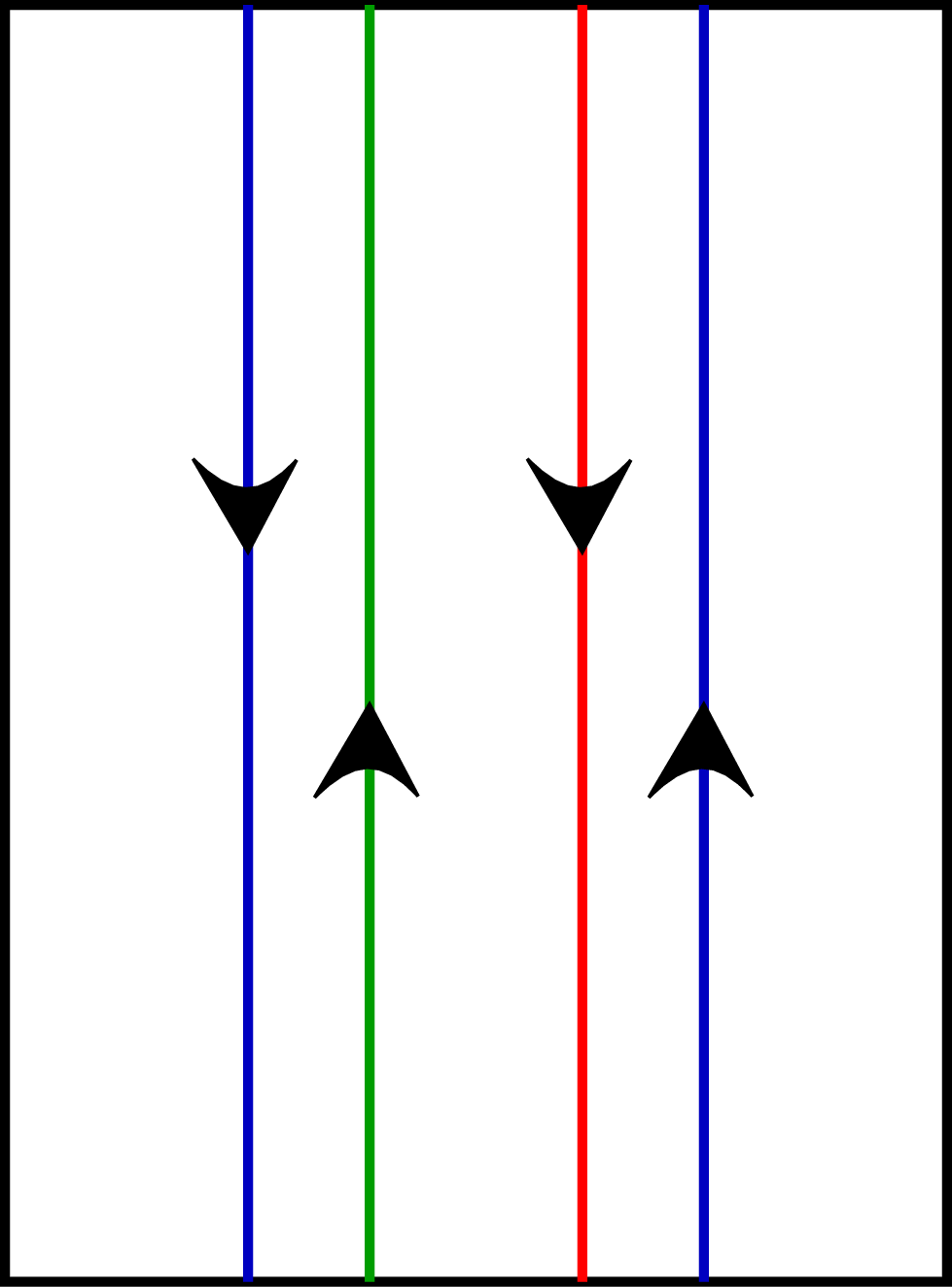} \endxy}.
\end{align}
The source and target of these morphisms is the object associated to the singlestep expression
\begin{equation} \label{flopsource} I_{\bullet} := [\hh{k+1} \supset \hh{1,k+1} \subset \hh{1} \supset \hh{0,1} \subset \hh{0}]. \end{equation}

From \cref{cond-bigon killing}, this simplifies to:
\begin{align}\label{diag-ssbim square flop simplified}
   (-1)^k\zeta^{\frac{k(k-1)}{2}+1}q^{k} \vcenter{\xy (0,0)*{\def\svgscale{0.14}\input{arxiv-figures/image_of_square_flop_LHS_svg-tex.eps_tex}} \endxy} = (-1)^k\zeta^{\frac{k(k+1)}{2}}q^{k} \vcenter{\xy (0,0)*{\def\svgscale{0.14}\input{arxiv-figures/image_of_square_flop_RHS_1_svg-tex.eps_tex}} \endxy} + [k-1]_q\vcenter{\xy (0,0)*{\def\svgscale{0.14}\input{arxiv-figures/image_of_square_flop_RHS_2_svg-tex.eps_tex}} \endxy}.
\end{align}
We may check this in the algebraic category $\mSBSBim$, as an equality in the (degree zero) endomorphism ring of 
\[ \BS(I_{\bullet}) \cong \Bimod{\hh{k+1}}{}{R^{\hh{1,k+1}}}\otimes_{\hh{1}}\Bimod{}{\hh{0}}{R^{\hh{1,0}}} (\ell). \]
The grading shift $\ell$ is irrelevant when considering endomorphisms, and we omit it below.

Explicitly, the LHS of \cref{diag-ssbim square flop} corresponds to the $(R^{\hh{k+1}}, R^{\hh{0}})$-bimodule morphism given by
\begin{align}\label{eq 1 sec square flop}
    f\otimes g \mapsto (-1)^k\zeta^{\frac{k(k-1)}{2}+1}q^{k}\;\ddm{1,2, k+1}{1,k+1}\left(\ddm{1,2,k+1}{2,k+1}(\Del{1,2}{2}{1}f)\Del{1,2}{1}{1}\right) \otimes \ddm{0,1,2}{0,1}\left(\ddm{0,1,2}{0,2}(\Del{1,2}{2}{2}g)\Del{1,2}{1}{2} \right).
\end{align}
The first term on the RHS corresponds to
\begin{align}\label{eq 1.1 sec square flop}
    f\otimes g \mapsto (-1)^k\zeta^{\frac{k(k+1)}{2}}q^{k}\;\ddm{0,1,k+1}{1,k+1}\left( \ddm{0,1,k+1}{0, k+1}(fg)\Del{0,1}{1}{1}\right) \otimes \Del{0,1}{1}{2},
\end{align}
and the second term to
\begin{equation}\label{eq 1.2 sec square flop}
f \ot g \mapsto [k-1]_q(f\otimes g). \end{equation}

Checking that the map in  \cref{eq 1 sec square flop} agrees with the sum of the maps in \cref{eq 1.1 sec square flop} and \cref{eq 1.2 sec square flop} may be done after forgetting the bimodule structures and just working within the category of $A$-modules. Since $\BS(I_\bullet)$ is a free module over $A$ (since it free as a right module over $R^{\hh{0}}$ which is free over $A$), it suffices to check the desired equality after base extension from $A$ to its fraction field. 

We now proceed with the check after base extension to the fraction field of $A$. Evaluating the first term on the RHS at $1 \ot 1$ yields zero for degree reasons, and the equality of the LHS with the second term on the RHS is \cref{lemma square flop 1 tensor 1}. Evaluating all three terms at $1 \ot x_1^k$ is done in \cref{lemma square flop 1 tensor x1k}, confirming \cref{diag-ssbim square flop simplified} as applied to this element. Note that all three terms are nonzero when applied to $1 \ot x_1^k$. These two calculations also suffice to show that the two terms on the RHS are linearly independent as bimodule endomorphisms.

In \cref{lem:dim2} we prove that $\dim \End^0\left(\BS(I_{\bullet})\right)=2$. (Our $k$ here agrees with $k-1$ in \cref{lem:dim2}.)
Thus the two terms on the RHS of \cref{diag-ssbim square flop} are a basis for $\End^0$. The LHS of \cref{diag-ssbim square flop} must be a unique linear combination of the terms in the RHS. Checking the equality on $1 \ot 1$ and $1 \ot x_1^k$ gave two independent constraints for two unknowns, verifying \cref{diag-ssbim square flop}.

Now we analyze the edge cases.

\subsubsection{When $k=n-1$}
In this case we desire the equality
\begin{align}\label{diag-square flop for k=n-1}
    \blkcoeff{\lambda_{n-2,1} \nu_{n-2,1} \lambda_{1,1} \nu_{1,1}} \vcenter{\xy (0,0)*{\def\svgscale{0.15}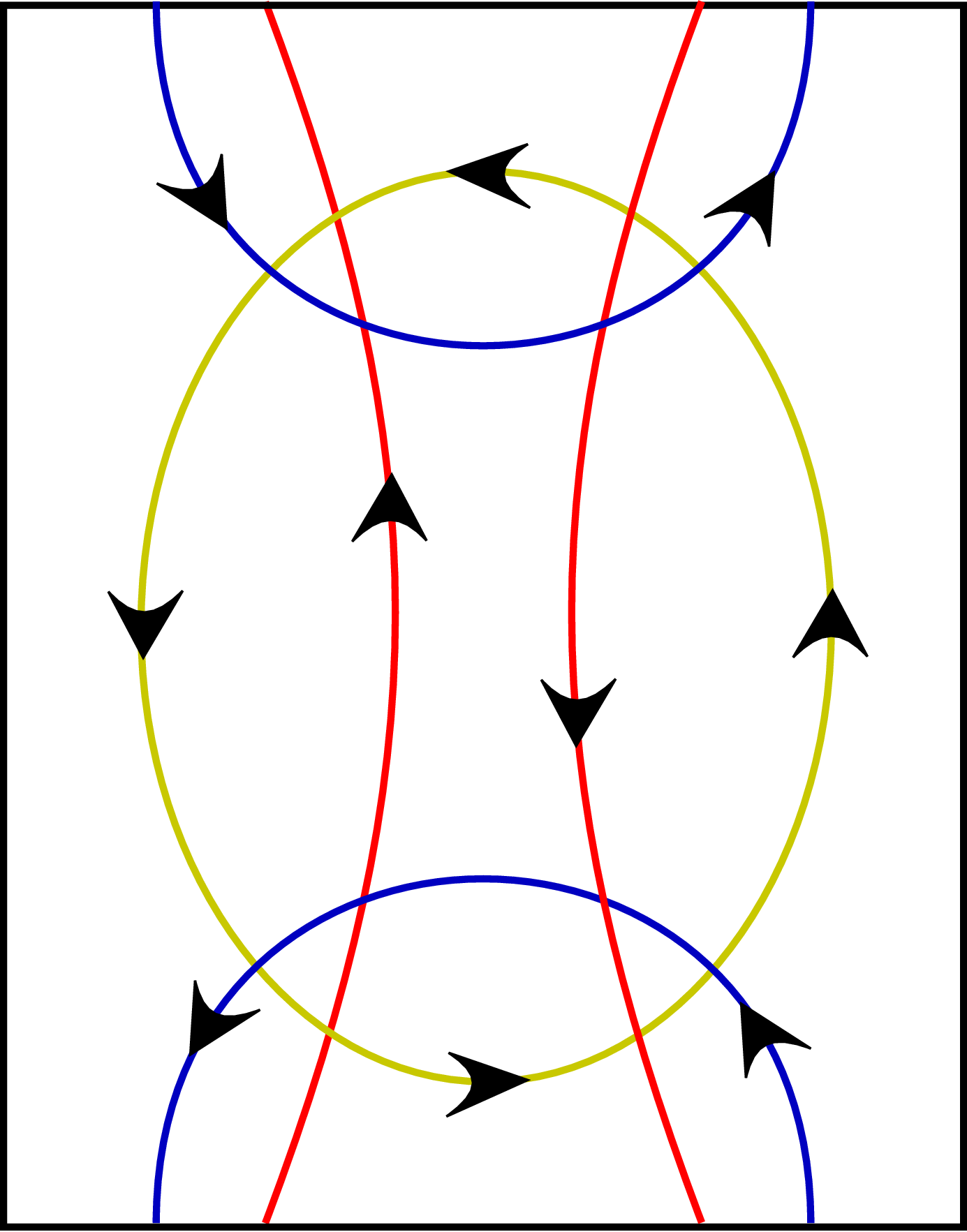} \endxy}=\blkcoeff{(-1)^{n-1}\lambda_{1}\nu_{n-1}}\vcenter{\xy (0,0)*{\def\svgscale{0.15}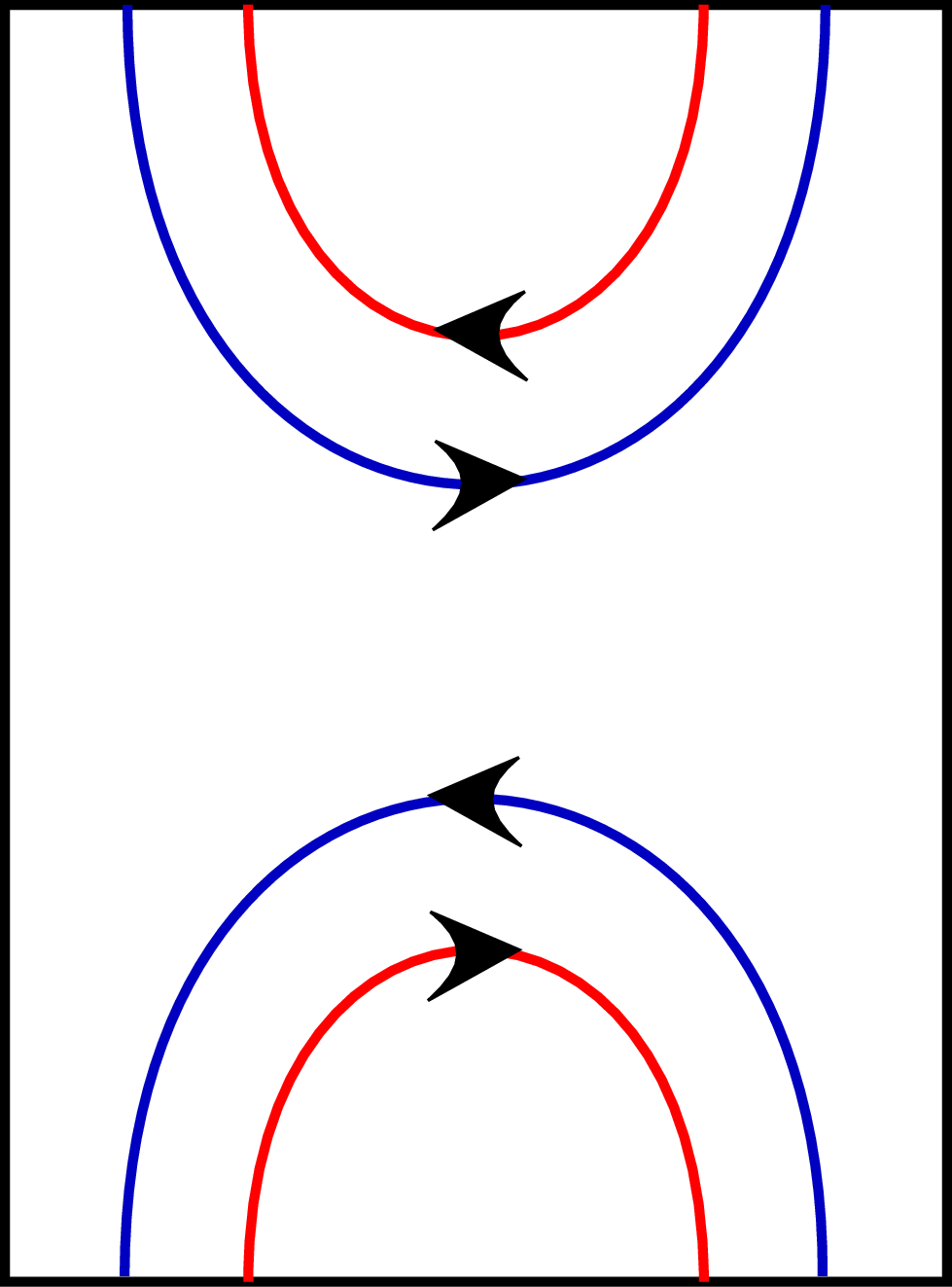} \endxy}+ \blkcoeff{[n-2]_q}\vcenter{\xy (0,0)*{\def\svgscale{0.15}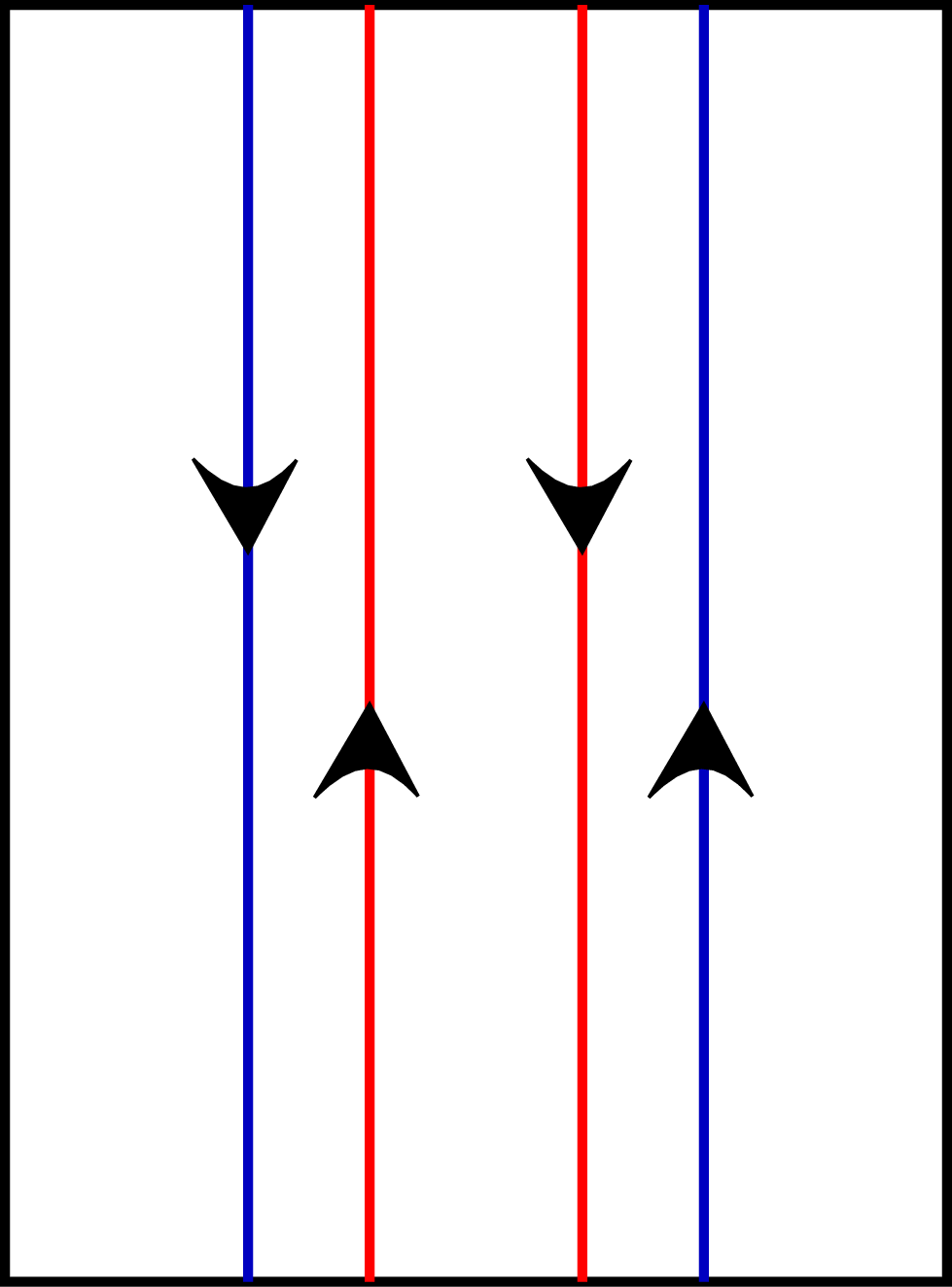} \endxy}.
\end{align}
From \cref{cond-bigon killing}, \cref{cond-tag switching} and \cref{cond-tag cancellation}, we have that $\lambda_{n-2,1} \nu_{n-2,1} \lambda_{1,1} \nu_{1,1} =(-1)^k\zeta^{\frac{k(k+1)}{2}}q^{k}$ and $\lambda_1 \nu_{n-1}=1$.
As before, it suffices to check this equality in the algebraic category $\mSBSBim$, in the ring $\End\left(\Bimod{\hh{0}}{}{R^{\hh{0,1}}}\otimes_{\hh{1}}\Bimod{}{\hh{0}}{R^{\hh{1,0}}}\right).$
The LHS corresponds to the bimodule homomorphism in \cref{eq 1 sec square flop}, given by
\begin{align}
    f\otimes g \mapsto (-1)^{n-1}\zeta^{\frac{(n-1)(n-2)}{2}+1}q^{n-1}\;\ddm{0,1,2}{0,1}\left(\ddm{0,1,2}{0,2}(\Del{1,2}{2}{1}f)\Del{1,2}{1}{1}\right) \otimes \ddm{0,1,2}{0,1}\left(\ddm{0,1,2}{0,2}(\Del{1,2}{2}{2}g)\Del{1,2}{1}{2} \right)
\end{align}
while the RHS corresponds to
\begin{align}
    f\otimes g\mapsto (-1)^{n-1}\ddm{0,1}{0}(fg)\Del{0,1}{1}{1} \otimes \Del{0,1}{1}{2} \;+\: [n-2]_q(f\otimes g).
\end{align}

The proof is almost the same as the previous case, and we need only check that the desired equality holds after base extension to the fraction field of $A$. Evaluation of the LHS on $1 \ot 1$ is still done using \cref{lemma square flop 1 tensor 1}. Evaluation of all terms at $1 \ot x_1^{n-1}$ is done in \cref{lemma square flop 1 tensor x1k for k=n-1}. This checks \cref{diag-square flop for k=n-1} in two linearly-independent cases.
As before, we showed in \cref{lem:dim2} that $\dim \End\left(\Bimod{\hh{0}}{}{R^{\hh{0,1}}}\otimes_{\hh{1}}\Bimod{}{\hh{0}}{R^{\hh{1,0}}}\right)=2$, and that the terms on the right are linearly independent, so they should form a basis. Thus equality holds in \cref{diag-square flop for k=n-1}. 

\subsubsection{When $k=n$:}
In this case, the first term on the RHS of $\cref{diag-square flop}$ is $0$ and the term on the LHS may be simplified using tag cancellation \cref{diag-tag cancellation}, so that the relation becomes (with valid region labels) 
\begin{align}\label{diag-square flop k=n}
   {
        \tikz[x=1mm, y=1mm, baseline=-0.5ex] {
          \draw[->-] (0,-8) -- (0,-4);
          \draw[-<-] (0,-4) arc (270:90:4) node[pos=0.5, left]{$1$};
          \draw[->-] (0,-4) arc (-90:90:4) node[pos=0.5,right]{$2$};
          \draw[->] (0,4) -- (0,8) node [right] {$1$};
        }
        = [n-1] \: \tikz[x=1mm, y=1mm, baseline=-0.5ex] {
          \draw[->] (0,-8) -- (0,8) node [right] {$1$};
        }
      }.
\end{align}
The 2-functor $\GS_{\zeta}$ sends the LHS to 
\begin{align}\label{diag-image of square flop k=n}
    (-1)^{n-1}\lambda_{1,1}\nu_{1,1} \vcenter{\xy (0,0)*{\def\svgscale{0.1}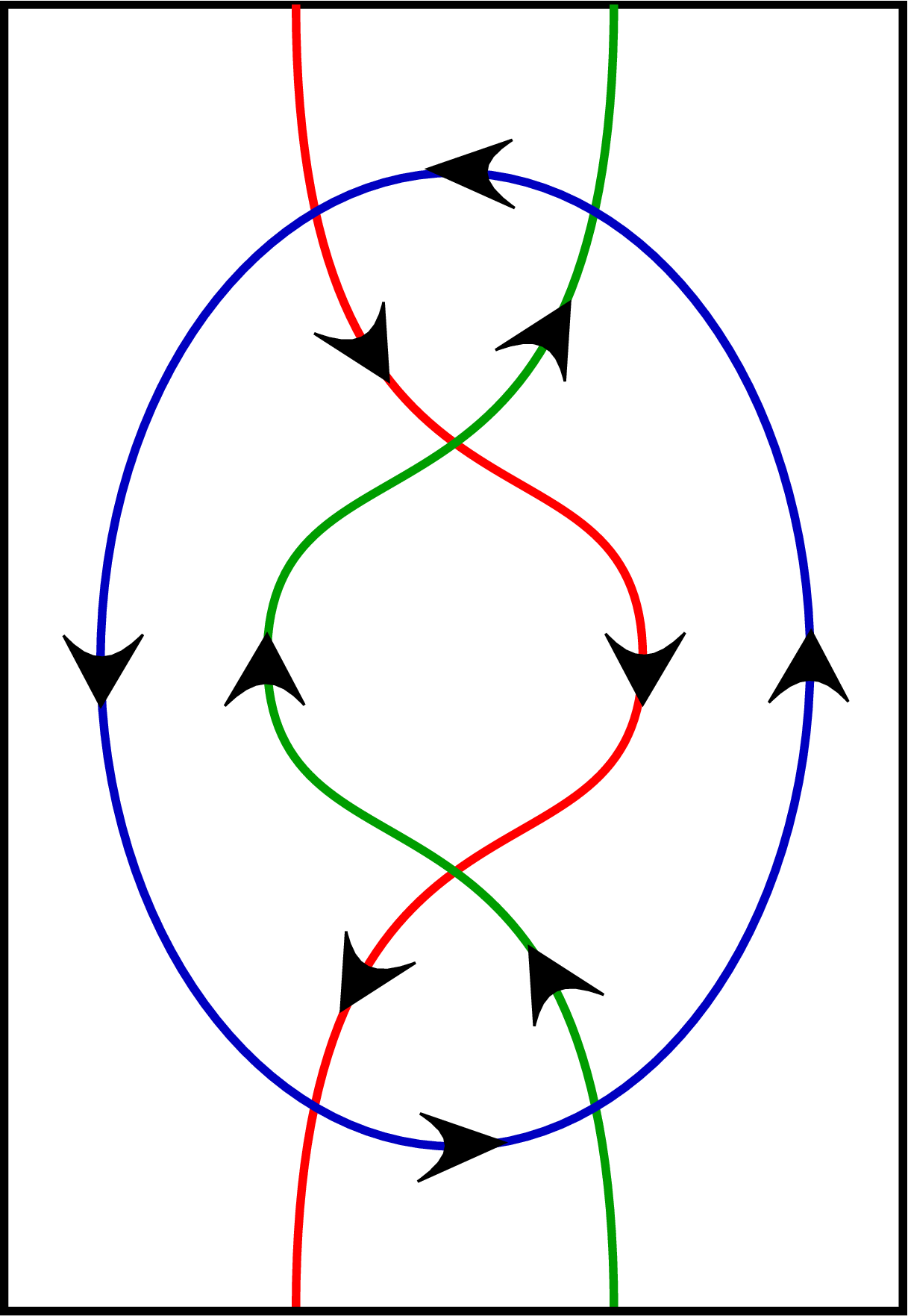} \endxy} \;\;.
\end{align}
By color symmetry, we may assume that $a=0$. 
Using \cref{R2unoriented2}, \cref{R2oriented} and \cref{ccwcirc} as in \cref{eq-bigon burst simplification},
the diagram above (without the coefficients) simplifies to
\begin{align}\label{diag-square flop k=n simplify}
    \vcenter{\xy (0,0)*{\def\svgscale{0.15}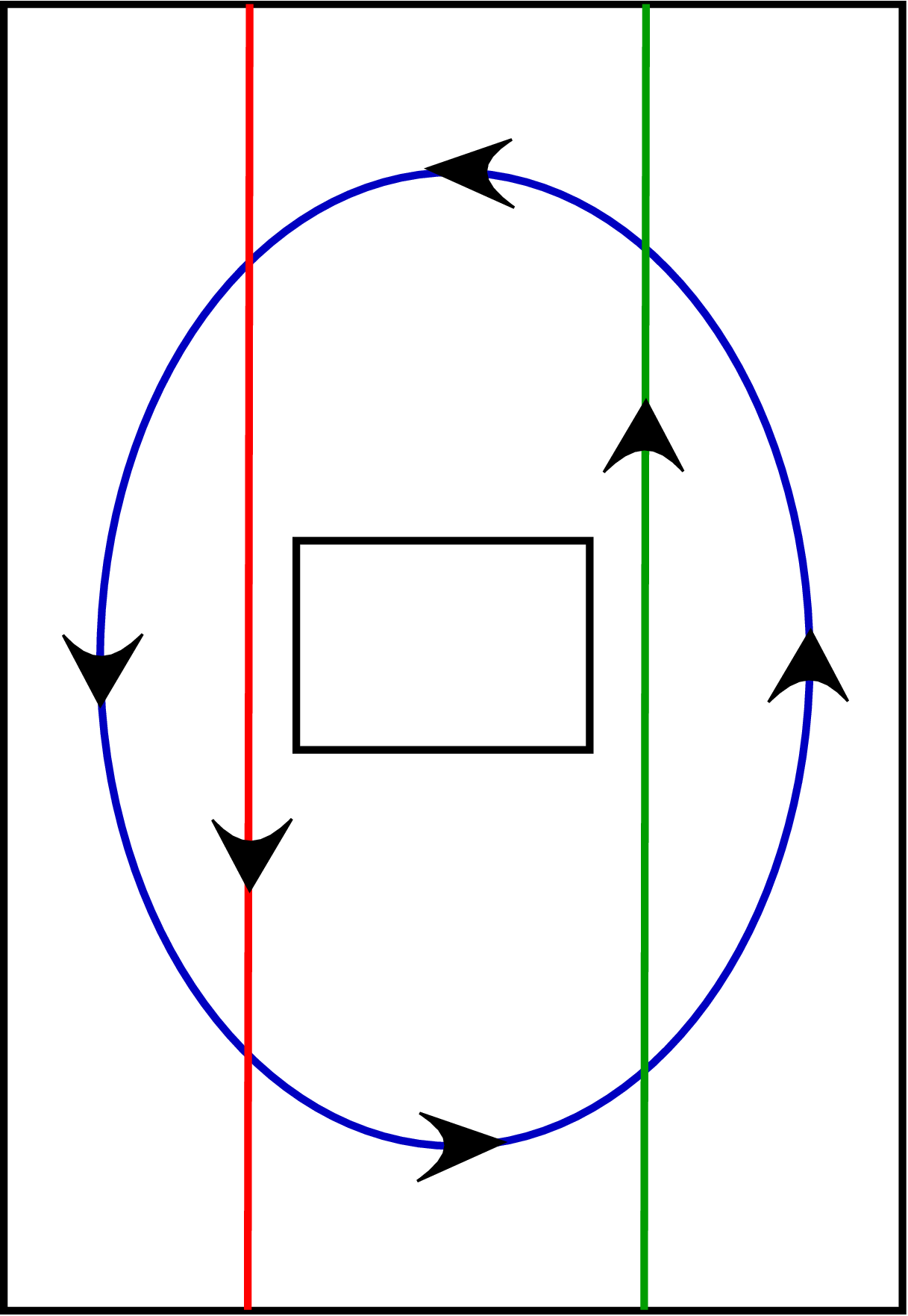} \endxy}=\vcenter{\xy (0,0)*{\def\svgscale{0.15}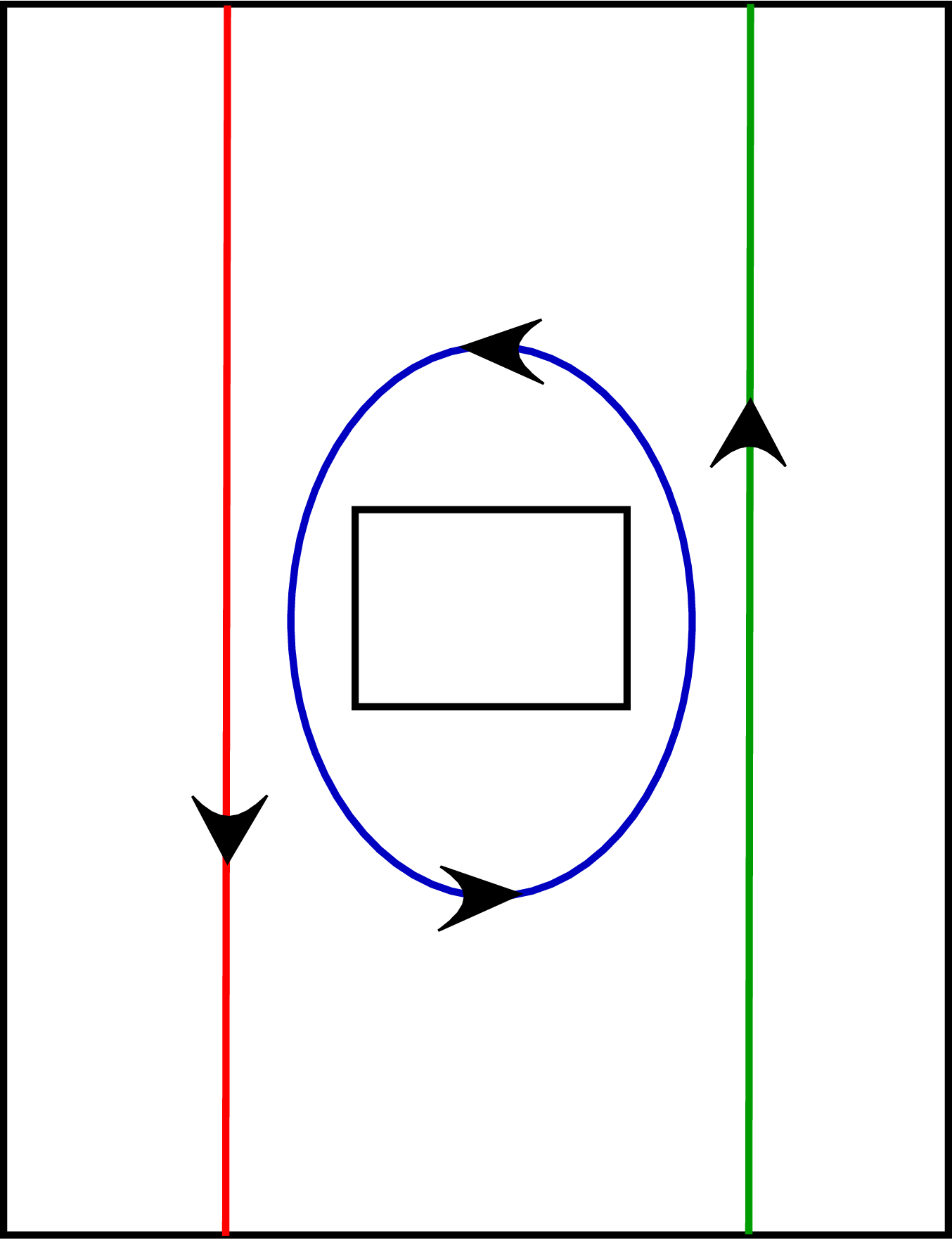} \endxy}= \vcenter{\xy (0,0)*{\def\svgscale{0.15}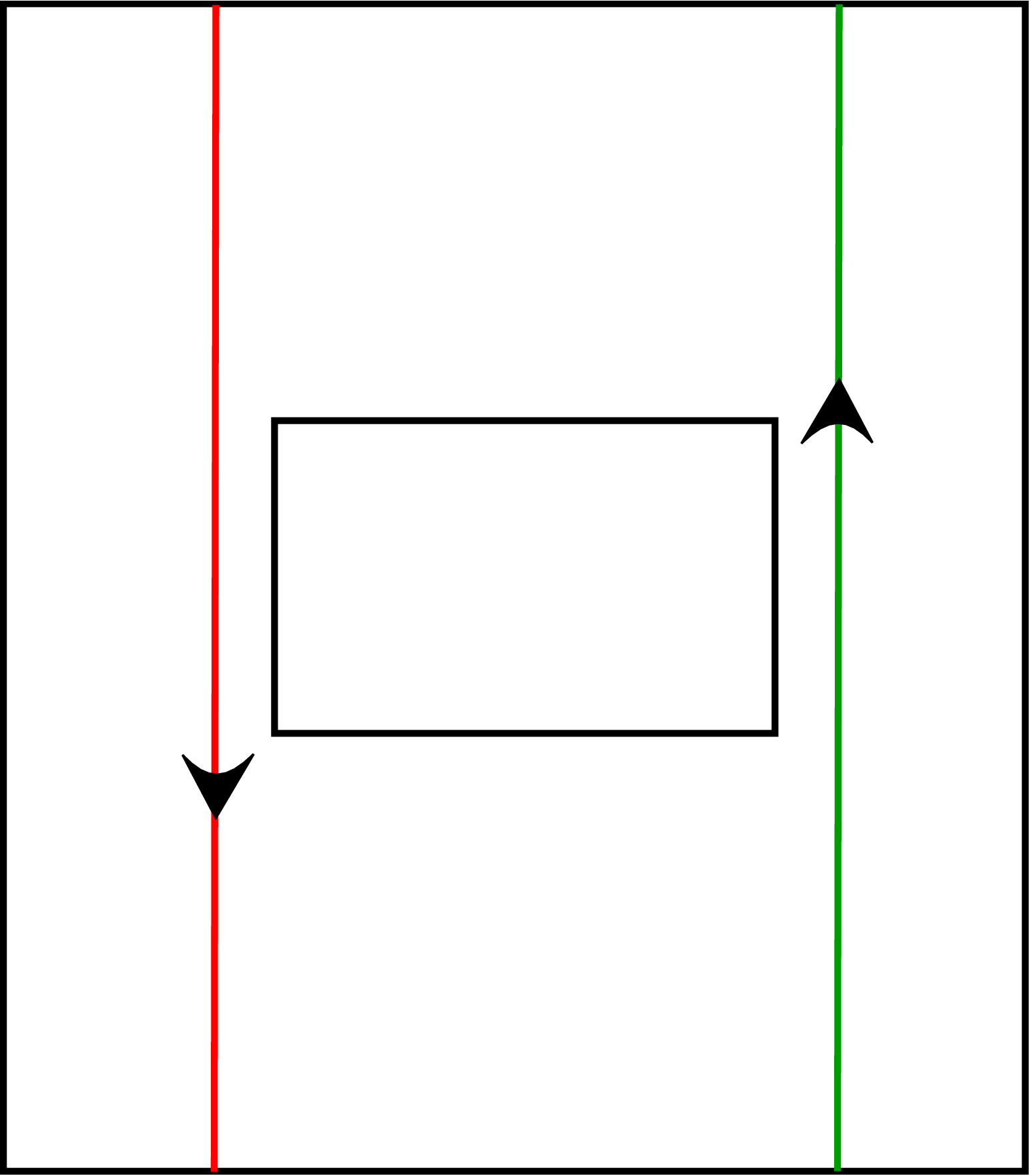} \endxy}.
\end{align}

We may simplify the polynomial in the box as follows.
\begin{align}\label{eq square flop k=n polynomial}
\notag \ddm{0,1,2}{0,1}\left(\mu^{\hh{0,2},\hh{1,2}}_{\hh{2}} \right) &= \zeta^{-\frac{(n-2)(n-3)}{2}}\dd_{n-1}\ldots \dd_3 \dd_2\left((x_3-\zeta^{n-1}x_2)(x_4- \zeta^{n-2}x_2)\ldots (x_n-\zeta^2x_2)\right)\\
\notag
    &\begin{array}{l}=\zeta^{-\frac{(n-2)(n-3)}{2}}(-1)^{n-2} \zeta^{(n-1)+      (n-2)+\ldots +2}\vspace{5 pt}\\
    \quad \quad\dd_{n-1}\ldots \dd_3\dd_2 \left((x_2-\zeta^{-2}x_n)(x_2-\zeta^{-3}x_{n-1})\ldots(x_2-\zeta^{-(n-1)}x_3)  \right)\end{array}\\
    &\begin{array}{l}=\zeta^{-\frac{(n-2)(n-3)}{2}}\zeta^{\frac{n(n-1)}{2}-1}(-1)^{n-2}\vspace{5 pt}\\
    \quad \quad\dd_{n-1}\ldots \dd_3\dd_2 \left((x_2-\zeta^{-2}x_n)(x_2-\zeta^{-3}x_{n-1})\ldots(x_2-\zeta^{-(n-1)}x_3)  \right).\end{array}
\end{align}
\cref{lemma demazure for bigon} and color symmetry tells us that 
$$\dd_{n-1}\ldots \dd_3\dd_2 \left((x_2-\zeta^{-2}x_n)(x_2-\zeta^{-3}x_{n-1})\ldots(x_2-\zeta^{-(n-1)}x_3)  \right)=q^{n-2}\zeta^{\frac{(n-2)(n-3)}{2}}[n-1]_q$$
so that \cref{eq square flop k=n polynomial} further simplifies to
\begin{align*}
   \ddm{0,1,2}{0,1}\left(\mu^{\hh{0,2},\hh{1,2}}_{\hh{2}} \right)  
   &= \zeta^{-\frac{(n-2)(n-3)}{2}}\zeta^{\frac{n(n-1)}{2}-1}(-1)^{n-2}q^{n-2}\zeta^{\frac{(n-2)(n-3)}{2}}[n-1]_q\\
   &=(-1)^{n-2}q^{-(n-1)}\zeta^{-1}q^{n-2}[n-1]_q =(-1)^{n-2}\zeta^{-1}q^{-1}[n-1]_q.
\end{align*}
Hence \cref{diag-image of square flop k=n} simplifies to
\[(-1)^{n-1} \lambda_{1,1}\nu_{1,1}\;(-1)^{n-2}\zeta^{-1}q^{-1}[n-1]_q=-\lambda_{1,1}\nu_{1,1}\zeta^{-1}q^{-1}[n-1]_q.\]
Comparing with the image of the RHS of \cref{diag-square flop k=n}, we have the conditions $\lambda_{1,1}\nu_{1,1}=-\zeta q.$
This is already captured in \cref{cond-bigon killing} for $k=2$.

%% file: arxiv-figures/tag_left_d_rotated_svg-tex.eps_tex
\begingroup%
  \makeatletter%
  \providecommand\color[2][]{%
    \errmessage{(Inkscape) Color is used for the text in Inkscape, but the package 'color.sty' is not loaded}%
    \renewcommand\color[2][]{}%
  }%
  \providecommand\transparent[1]{%
    \errmessage{(Inkscape) Transparency is used (non-zero) for the text in Inkscape, but the package 'transparent.sty' is not loaded}%
    \renewcommand\transparent[1]{}%
  }%
  \providecommand\rotatebox[2]{#2}%
  \newcommand*\fsize{\dimexpr\f@size pt\relax}%
  \newcommand*\lineheight[1]{\fontsize{\fsize}{#1\fsize}\selectfont}%
  \ifx\svgwidth\undefined%
    \setlength{\unitlength}{730.76518299bp}%
    \ifx\svgscale\undefined%
      \relax%
    \else%
      \setlength{\unitlength}{\unitlength * \real{\svgscale}}%
    \fi%
  \else%
    \setlength{\unitlength}{\svgwidth}%
  \fi%
  \global\let\svgwidth\undefined%
  \global\let\svgscale\undefined%
  \makeatother%
  \begin{picture}(1,0.86542458)%
    \lineheight{1}%
    \setlength\tabcolsep{0pt}%
    \put(0,0){\includegraphics[width=\unitlength]{tag_left_d_rotated_svg-tex.eps}}%
    \put(2.12603626,-0.9015041){\color[rgb]{0,0,0}\makebox(0,0)[lt]{\begin{minipage}{0.30656706\unitlength}\end{minipage}}}%
    \put(0.6012437,0.48508123){\color[rgb]{0,0,0}\makebox(0,0)[lt]{\begin{minipage}{0.10925356\unitlength}\end{minipage}}}%
    \put(2.05721508,-0.68252765){\color[rgb]{0,0,0}\makebox(0,0)[lt]{\begin{minipage}{0.4504657\unitlength}\end{minipage}}}%
    \put(0.54338646,0.14703318){\color[rgb]{0,0,0}\makebox(0,0)[lt]{\begin{minipage}{0.08276784\unitlength}\raggedright \gr{a-k}\end{minipage}}}%
    \put(0.31648533,0.72187127){\color[rgb]{0,0,0}\makebox(0,0)[lt]{\lineheight{1.25}\smash{\begin{tabular}[t]{l}\gr{a}\end{tabular}}}}%
    \put(0.535392,0.34426734){\color[rgb]{0,0,0}\makebox(0,0)[lt]{\begin{minipage}{0.05792111\unitlength}\raggedright \blk{n-k}\end{minipage}}}%
    \put(0.52948508,0.56574299){\color[rgb]{0,0,0}\makebox(0,0)[lt]{\begin{minipage}{0.26683146\unitlength}\raggedright \blk{k}                           \end{minipage}}}%
    \put(0.50508546,0.3417535){\color[rgb]{0,0,0}\makebox(0,0)[lt]{\begin{minipage}{0.05792111\unitlength}\end{minipage}}}%
  \end{picture}%
\endgroup%

%% file: arxiv-figures/image_of_tag_rotated_svg-tex.eps_tex
\begingroup%
  \makeatletter%
  \providecommand\color[2][]{%
    \errmessage{(Inkscape) Color is used for the text in Inkscape, but the package 'color.sty' is not loaded}%
    \renewcommand\color[2][]{}%
  }%
  \providecommand\transparent[1]{%
    \errmessage{(Inkscape) Transparency is used (non-zero) for the text in Inkscape, but the package 'transparent.sty' is not loaded}%
    \renewcommand\transparent[1]{}%
  }%
  \providecommand\rotatebox[2]{#2}%
  \newcommand*\fsize{\dimexpr\f@size pt\relax}%
  \newcommand*\lineheight[1]{\fontsize{\fsize}{#1\fsize}\selectfont}%
  \ifx\svgwidth\undefined%
    \setlength{\unitlength}{589.51936895bp}%
    \ifx\svgscale\undefined%
      \relax%
    \else%
      \setlength{\unitlength}{\unitlength * \real{\svgscale}}%
    \fi%
  \else%
    \setlength{\unitlength}{\svgwidth}%
  \fi%
  \global\let\svgwidth\undefined%
  \global\let\svgscale\undefined%
  \makeatother%
  \begin{picture}(1,0.95830354)%
    \lineheight{1}%
    \setlength\tabcolsep{0pt}%
    \put(0,0){\includegraphics[width=\unitlength]{image_of_tag_rotated_svg-tex.eps}}%
    \put(0.6046879,0.14019604){\color[rgb]{0,0,0}\makebox(0,0)[lt]{\lineheight{1.25}\smash{\begin{tabular}[t]{l}\dgr{a}\end{tabular}}}}%
    \put(0.26406039,0.83376259){\color[rgb]{0,0,0}\makebox(0,0)[lt]{\lineheight{1.25}\smash{\begin{tabular}[t]{l}\dgr{a-k}\end{tabular}}}}%
  \end{picture}%
\endgroup%

%% file: arxiv-figures/I=H_LHS_svg-tex.eps_tex
\begingroup%
  \makeatletter%
  \providecommand\color[2][]{%
    \errmessage{(Inkscape) Color is used for the text in Inkscape, but the package 'color.sty' is not loaded}%
    \renewcommand\color[2][]{}%
  }%
  \providecommand\transparent[1]{%
    \errmessage{(Inkscape) Transparency is used (non-zero) for the text in Inkscape, but the package 'transparent.sty' is not loaded}%
    \renewcommand\transparent[1]{}%
  }%
  \providecommand\rotatebox[2]{#2}%
  \newcommand*\fsize{\dimexpr\f@size pt\relax}%
  \newcommand*\lineheight[1]{\fontsize{\fsize}{#1\fsize}\selectfont}%
  \ifx\svgwidth\undefined%
    \setlength{\unitlength}{769.87426195bp}%
    \ifx\svgscale\undefined%
      \relax%
    \else%
      \setlength{\unitlength}{\unitlength * \real{\svgscale}}%
    \fi%
  \else%
    \setlength{\unitlength}{\svgwidth}%
  \fi%
  \global\let\svgwidth\undefined%
  \global\let\svgscale\undefined%
  \makeatother%
  \begin{picture}(1,0.73761988)%
    \lineheight{1}%
    \setlength\tabcolsep{0pt}%
    \put(0,0){\includegraphics[width=\unitlength]{I=H_LHS_svg-tex.eps}}%
    \put(0.7273098,0.53772226){\color[rgb]{0,0,0}\makebox(0,0)[lt]{\lineheight{1.25}\smash{\begin{tabular}[t]{l}\gr{a}\end{tabular}}}}%
    \put(0.63431946,0.08310254){\color[rgb]{0,0,0}\makebox(0,0)[lt]{\lineheight{1.25}\smash{\begin{tabular}[t]{l}\gr{a+m}\end{tabular}}}}%
    \put(0.36863262,0.27055942){\color[rgb]{0,0,0}\makebox(0,0)[lt]{\lineheight{1.25}\smash{\begin{tabular}[t]{l}\gr{a+m+l}\end{tabular}}}}%
    \put(0.51863607,0.14633947){\color[rgb]{0,0,0}\makebox(0,0)[lt]{\lineheight{1.25}\smash{\begin{tabular}[t]{l}\blk{l}\end{tabular}}}}%
    \put(0.20479239,0.28236774){\color[rgb]{0,0,0}\makebox(0,0)[lt]{\lineheight{1.25}\smash{\begin{tabular}[t]{l}\blk{k}\end{tabular}}}}%
    \put(0.80758844,0.15133576){\color[rgb]{0,0,0}\makebox(0,0)[lt]{\lineheight{1.25}\smash{\begin{tabular}[t]{l}\blk{m}\end{tabular}}}}%
    \put(0.05087915,0.53825433){\color[rgb]{0,0,0}\makebox(0,0)[lt]{\lineheight{1.25}\smash{\begin{tabular}[t]{l}\gr{a+k+l+m}\end{tabular}}}}%
  \end{picture}%
\endgroup%

%% file: arxiv-figures/I=H_RHS_svg-tex.eps_tex
\begingroup%
  \makeatletter%
  \providecommand\color[2][]{%
    \errmessage{(Inkscape) Color is used for the text in Inkscape, but the package 'color.sty' is not loaded}%
    \renewcommand\color[2][]{}%
  }%
  \providecommand\transparent[1]{%
    \errmessage{(Inkscape) Transparency is used (non-zero) for the text in Inkscape, but the package 'transparent.sty' is not loaded}%
    \renewcommand\transparent[1]{}%
  }%
  \providecommand\rotatebox[2]{#2}%
  \newcommand*\fsize{\dimexpr\f@size pt\relax}%
  \newcommand*\lineheight[1]{\fontsize{\fsize}{#1\fsize}\selectfont}%
  \ifx\svgwidth\undefined%
    \setlength{\unitlength}{769.87426195bp}%
    \ifx\svgscale\undefined%
      \relax%
    \else%
      \setlength{\unitlength}{\unitlength * \real{\svgscale}}%
    \fi%
  \else%
    \setlength{\unitlength}{\svgwidth}%
  \fi%
  \global\let\svgwidth\undefined%
  \global\let\svgscale\undefined%
  \makeatother%
  \begin{picture}(1,0.73761988)%
    \lineheight{1}%
    \setlength\tabcolsep{0pt}%
    \put(0,0){\includegraphics[width=\unitlength]{I=H_RHS_svg-tex.eps}}%
    \put(0.74242945,0.52591393){\color[rgb]{0,0,0}\makebox(0,0)[lt]{\lineheight{1.25}\smash{\begin{tabular}[t]{l}\gr{a}\end{tabular}}}}%
    \put(0.44574586,0.28384371){\color[rgb]{0,0,0}\makebox(0,0)[lt]{\lineheight{1.25}\smash{\begin{tabular}[t]{l}\gr{a+m}\end{tabular}}}}%
    \put(0.17415485,0.0432495){\color[rgb]{0,0,0}\makebox(0,0)[lt]{\lineheight{1.25}\smash{\begin{tabular}[t]{l}\gr{a+l+m}\end{tabular}}}}%
    \put(0.71826047,0.28656304){\color[rgb]{0,0,0}\makebox(0,0)[lt]{\lineheight{1.25}\smash{\begin{tabular}[t]{l}\blk{m}\end{tabular}}}}%
    \put(0.40884491,0.15395245){\color[rgb]{0,0,0}\makebox(0,0)[lt]{\lineheight{1.25}\smash{\begin{tabular}[t]{l}\blk{l}\end{tabular}}}}%
    \put(0.1466827,0.1572399){\color[rgb]{0,0,0}\makebox(0,0)[lt]{\lineheight{1.25}\smash{\begin{tabular}[t]{l}\blk{k}\end{tabular}}}}%
    \put(0.06009461,0.52349397){\color[rgb]{0,0,0}\makebox(0,0)[lt]{\lineheight{1.25}\smash{\begin{tabular}[t]{l}\gr{a+k+l+m}\end{tabular}}}}%
  \end{picture}%
\endgroup%

%% file: arxiv-figures/image_of_I=H_LHS_svg-tex.eps_tex
\begingroup%
  \makeatletter%
  \providecommand\color[2][]{%
    \errmessage{(Inkscape) Color is used for the text in Inkscape, but the package 'color.sty' is not loaded}%
    \renewcommand\color[2][]{}%
  }%
  \providecommand\transparent[1]{%
    \errmessage{(Inkscape) Transparency is used (non-zero) for the text in Inkscape, but the package 'transparent.sty' is not loaded}%
    \renewcommand\transparent[1]{}%
  }%
  \providecommand\rotatebox[2]{#2}%
  \newcommand*\fsize{\dimexpr\f@size pt\relax}%
  \newcommand*\lineheight[1]{\fontsize{\fsize}{#1\fsize}\selectfont}%
  \ifx\svgwidth\undefined%
    \setlength{\unitlength}{677.65956526bp}%
    \ifx\svgscale\undefined%
      \relax%
    \else%
      \setlength{\unitlength}{\unitlength * \real{\svgscale}}%
    \fi%
  \else%
    \setlength{\unitlength}{\svgwidth}%
  \fi%
  \global\let\svgwidth\undefined%
  \global\let\svgscale\undefined%
  \makeatother%
  \begin{picture}(1,1.05060582)%
    \lineheight{1}%
    \setlength\tabcolsep{0pt}%
    \put(0,0){\includegraphics[width=\unitlength]{image_of_I=H_LHS_svg-tex.eps}}%
    \put(0.33988245,0.44531598){\color[rgb]{0,0,0}\makebox(0,0)[lt]{\lineheight{1.25}\smash{\begin{tabular}[t]{l}\dgr{a+l+m}\end{tabular}}}}%
    \put(0.61284979,0.03989832){\color[rgb]{0,0,0}\makebox(0,0)[lt]{\lineheight{1.25}\smash{\begin{tabular}[t]{l}\dgr{a+m}\end{tabular}}}}%
    \put(0.78330071,0.81070158){\color[rgb]{0,0,0}\makebox(0,0)[lt]{\lineheight{1.25}\smash{\begin{tabular}[t]{l}\dgr{a}\end{tabular}}}}%
    \put(0.02638172,0.80904963){\color[rgb]{0,0,0}\makebox(0,0)[lt]{\lineheight{1.25}\smash{\begin{tabular}[t]{l}\dgr{a+k+l+m}\end{tabular}}}}%
  \end{picture}%
\endgroup%

%% file: arxiv-figures/image_of_I=H_RHS_svg-tex.eps_tex
\begingroup%
  \makeatletter%
  \providecommand\color[2][]{%
    \errmessage{(Inkscape) Color is used for the text in Inkscape, but the package 'color.sty' is not loaded}%
    \renewcommand\color[2][]{}%
  }%
  \providecommand\transparent[1]{%
    \errmessage{(Inkscape) Transparency is used (non-zero) for the text in Inkscape, but the package 'transparent.sty' is not loaded}%
    \renewcommand\transparent[1]{}%
  }%
  \providecommand\rotatebox[2]{#2}%
  \newcommand*\fsize{\dimexpr\f@size pt\relax}%
  \newcommand*\lineheight[1]{\fontsize{\fsize}{#1\fsize}\selectfont}%
  \ifx\svgwidth\undefined%
    \setlength{\unitlength}{677.65956526bp}%
    \ifx\svgscale\undefined%
      \relax%
    \else%
      \setlength{\unitlength}{\unitlength * \real{\svgscale}}%
    \fi%
  \else%
    \setlength{\unitlength}{\svgwidth}%
  \fi%
  \global\let\svgwidth\undefined%
  \global\let\svgscale\undefined%
  \makeatother%
  \begin{picture}(1,1.05060582)%
    \lineheight{1}%
    \setlength\tabcolsep{0pt}%
    \put(0,0){\includegraphics[width=\unitlength]{image_of_I=H_RHS_svg-tex.eps}}%
    \put(0.46565907,0.43375294){\color[rgb]{0,0,0}\makebox(0,0)[lt]{\lineheight{1.25}\smash{\begin{tabular}[t]{l}\dgr{a+m}\end{tabular}}}}%
    \put(0.21994037,0.02172775){\color[rgb]{0,0,0}\makebox(0,0)[lt]{\lineheight{1.25}\smash{\begin{tabular}[t]{l}\dgr{a+l+m}\end{tabular}}}}%
    \put(0.79675042,0.78261991){\color[rgb]{0,0,0}\makebox(0,0)[lt]{\lineheight{1.25}\smash{\begin{tabular}[t]{l}\dgr{a}\end{tabular}}}}%
    \put(0.02000906,0.78261978){\color[rgb]{0,0,0}\makebox(0,0)[lt]{\lineheight{1.25}\smash{\begin{tabular}[t]{l}\dgr{a+k+l+m}\end{tabular}}}}%
  \end{picture}%
\endgroup%

%% file: arxiv-figures/image_of_I=H_LHS_k+l+m=n_svg-tex.eps_tex
\begingroup%
  \makeatletter%
  \providecommand\color[2][]{%
    \errmessage{(Inkscape) Color is used for the text in Inkscape, but the package 'color.sty' is not loaded}%
    \renewcommand\color[2][]{}%
  }%
  \providecommand\transparent[1]{%
    \errmessage{(Inkscape) Transparency is used (non-zero) for the text in Inkscape, but the package 'transparent.sty' is not loaded}%
    \renewcommand\transparent[1]{}%
  }%
  \providecommand\rotatebox[2]{#2}%
  \newcommand*\fsize{\dimexpr\f@size pt\relax}%
  \newcommand*\lineheight[1]{\fontsize{\fsize}{#1\fsize}\selectfont}%
  \ifx\svgwidth\undefined%
    \setlength{\unitlength}{677.65956526bp}%
    \ifx\svgscale\undefined%
      \relax%
    \else%
      \setlength{\unitlength}{\unitlength * \real{\svgscale}}%
    \fi%
  \else%
    \setlength{\unitlength}{\svgwidth}%
  \fi%
  \global\let\svgwidth\undefined%
  \global\let\svgscale\undefined%
  \makeatother%
  \begin{picture}(1,1.05060582)%
    \lineheight{1}%
    \setlength\tabcolsep{0pt}%
    \put(0,0){\includegraphics[width=\unitlength]{image_of_I=H_LHS_k+l+m=n_svg-tex.eps}}%
    \put(0.33988245,0.44531598){\color[rgb]{0,0,0}\makebox(0,0)[lt]{\lineheight{1.25}\smash{\begin{tabular}[t]{l}\dgr{a+l+m}\end{tabular}}}}%
    \put(0.61284979,0.03989832){\color[rgb]{0,0,0}\makebox(0,0)[lt]{\lineheight{1.25}\smash{\begin{tabular}[t]{l}\dgr{a+m}\end{tabular}}}}%
    \put(0.41989013,0.8883393){\color[rgb]{0,0,0}\makebox(0,0)[lt]{\lineheight{1.25}\smash{\begin{tabular}[t]{l}\dgr{a}\end{tabular}}}}%
  \end{picture}%
\endgroup%

%% file: arxiv-figures/image_of_I=H_RHS_k+l+m=n_svg-tex.eps_tex
\begingroup%
  \makeatletter%
  \providecommand\color[2][]{%
    \errmessage{(Inkscape) Color is used for the text in Inkscape, but the package 'color.sty' is not loaded}%
    \renewcommand\color[2][]{}%
  }%
  \providecommand\transparent[1]{%
    \errmessage{(Inkscape) Transparency is used (non-zero) for the text in Inkscape, but the package 'transparent.sty' is not loaded}%
    \renewcommand\transparent[1]{}%
  }%
  \providecommand\rotatebox[2]{#2}%
  \newcommand*\fsize{\dimexpr\f@size pt\relax}%
  \newcommand*\lineheight[1]{\fontsize{\fsize}{#1\fsize}\selectfont}%
  \ifx\svgwidth\undefined%
    \setlength{\unitlength}{677.65956526bp}%
    \ifx\svgscale\undefined%
      \relax%
    \else%
      \setlength{\unitlength}{\unitlength * \real{\svgscale}}%
    \fi%
  \else%
    \setlength{\unitlength}{\svgwidth}%
  \fi%
  \global\let\svgwidth\undefined%
  \global\let\svgscale\undefined%
  \makeatother%
  \begin{picture}(1,1.05060582)%
    \lineheight{1}%
    \setlength\tabcolsep{0pt}%
    \put(0,0){\includegraphics[width=\unitlength]{image_of_I=H_RHS_k+l+m=n_svg-tex.eps}}%
    \put(0.46565907,0.43375294){\color[rgb]{0,0,0}\makebox(0,0)[lt]{\lineheight{1.25}\smash{\begin{tabular}[t]{l}\dgr{a+m}\end{tabular}}}}%
    \put(0.21994037,0.02172775){\color[rgb]{0,0,0}\makebox(0,0)[lt]{\lineheight{1.25}\smash{\begin{tabular}[t]{l}\dgr{a+l+m}\end{tabular}}}}%
    \put(0.50767381,0.88999122){\color[rgb]{0,0,0}\makebox(0,0)[lt]{\lineheight{1.25}\smash{\begin{tabular}[t]{l}\dgr{a}\end{tabular}}}}%
  \end{picture}%
\endgroup%

%% file: arxiv-figures/associativity_ssbim_LHS_svg-tex.eps_tex
\begingroup%
  \makeatletter%
  \providecommand\color[2][]{%
    \errmessage{(Inkscape) Color is used for the text in Inkscape, but the package 'color.sty' is not loaded}%
    \renewcommand\color[2][]{}%
  }%
  \providecommand\transparent[1]{%
    \errmessage{(Inkscape) Transparency is used (non-zero) for the text in Inkscape, but the package 'transparent.sty' is not loaded}%
    \renewcommand\transparent[1]{}%
  }%
  \providecommand\rotatebox[2]{#2}%
  \newcommand*\fsize{\dimexpr\f@size pt\relax}%
  \newcommand*\lineheight[1]{\fontsize{\fsize}{#1\fsize}\selectfont}%
  \ifx\svgwidth\undefined%
    \setlength{\unitlength}{677.65956526bp}%
    \ifx\svgscale\undefined%
      \relax%
    \else%
      \setlength{\unitlength}{\unitlength * \real{\svgscale}}%
    \fi%
  \else%
    \setlength{\unitlength}{\svgwidth}%
  \fi%
  \global\let\svgwidth\undefined%
  \global\let\svgscale\undefined%
  \makeatother%
  \begin{picture}(1,1.05060582)%
    \lineheight{1}%
    \setlength\tabcolsep{0pt}%
    \put(0,0){\includegraphics[width=\unitlength]{associativity_ssbim_LHS_svg-tex.eps}}%
    \put(0.38943844,0.44366414){\color[rgb]{0,0,0}\makebox(0,0)[lt]{\lineheight{1.25}\smash{\begin{tabular}[t]{l}\dgr{l+m}\end{tabular}}}}%
    \put(0.65745017,0.03989832){\color[rgb]{0,0,0}\makebox(0,0)[lt]{\lineheight{1.25}\smash{\begin{tabular}[t]{l}\dgr{m}\end{tabular}}}}%
    \put(0.78330071,0.81070158){\color[rgb]{0,0,0}\makebox(0,0)[lt]{\lineheight{1.25}\smash{\begin{tabular}[t]{l}\dgr{0}\end{tabular}}}}%
    \put(0.07593771,0.80904963){\color[rgb]{0,0,0}\makebox(0,0)[lt]{\lineheight{1.25}\smash{\begin{tabular}[t]{l}\dgr{k+l+m}\end{tabular}}}}%
  \end{picture}%
\endgroup%

%% file: arxiv-figures/associativity_ssbim_1_svg-tex.eps_tex
\begingroup%
  \makeatletter%
  \providecommand\color[2][]{%
    \errmessage{(Inkscape) Color is used for the text in Inkscape, but the package 'color.sty' is not loaded}%
    \renewcommand\color[2][]{}%
  }%
  \providecommand\transparent[1]{%
    \errmessage{(Inkscape) Transparency is used (non-zero) for the text in Inkscape, but the package 'transparent.sty' is not loaded}%
    \renewcommand\transparent[1]{}%
  }%
  \providecommand\rotatebox[2]{#2}%
  \newcommand*\fsize{\dimexpr\f@size pt\relax}%
  \newcommand*\lineheight[1]{\fontsize{\fsize}{#1\fsize}\selectfont}%
  \ifx\svgwidth\undefined%
    \setlength{\unitlength}{677.65956526bp}%
    \ifx\svgscale\undefined%
      \relax%
    \else%
      \setlength{\unitlength}{\unitlength * \real{\svgscale}}%
    \fi%
  \else%
    \setlength{\unitlength}{\svgwidth}%
  \fi%
  \global\let\svgwidth\undefined%
  \global\let\svgscale\undefined%
  \makeatother%
  \begin{picture}(1,1.05060582)%
    \lineheight{1}%
    \setlength\tabcolsep{0pt}%
    \put(0,0){\includegraphics[width=\unitlength]{associativity_ssbim_1_svg-tex.eps}}%
    \put(0.33823058,0.21901026){\color[rgb]{0,0,0}\makebox(0,0)[lt]{\lineheight{1.25}\smash{\begin{tabular}[t]{l}\dgr{l+m}\end{tabular}}}}%
    \put(0.64423524,0.02998717){\color[rgb]{0,0,0}\makebox(0,0)[lt]{\lineheight{1.25}\smash{\begin{tabular}[t]{l}\dgr{m}\end{tabular}}}}%
    \put(0.78330071,0.81070158){\color[rgb]{0,0,0}\makebox(0,0)[lt]{\lineheight{1.25}\smash{\begin{tabular}[t]{l}\dgr{0}\end{tabular}}}}%
    \put(0.08089331,0.80904963){\color[rgb]{0,0,0}\makebox(0,0)[lt]{\lineheight{1.25}\smash{\begin{tabular}[t]{l}\dgr{k+l+m}\end{tabular}}}}%
  \end{picture}%
\endgroup%

%% file: arxiv-figures/associativity_ssbim_2_svg-tex.eps_tex
\begingroup%
  \makeatletter%
  \providecommand\color[2][]{%
    \errmessage{(Inkscape) Color is used for the text in Inkscape, but the package 'color.sty' is not loaded}%
    \renewcommand\color[2][]{}%
  }%
  \providecommand\transparent[1]{%
    \errmessage{(Inkscape) Transparency is used (non-zero) for the text in Inkscape, but the package 'transparent.sty' is not loaded}%
    \renewcommand\transparent[1]{}%
  }%
  \providecommand\rotatebox[2]{#2}%
  \newcommand*\fsize{\dimexpr\f@size pt\relax}%
  \newcommand*\lineheight[1]{\fontsize{\fsize}{#1\fsize}\selectfont}%
  \ifx\svgwidth\undefined%
    \setlength{\unitlength}{677.65956526bp}%
    \ifx\svgscale\undefined%
      \relax%
    \else%
      \setlength{\unitlength}{\unitlength * \real{\svgscale}}%
    \fi%
  \else%
    \setlength{\unitlength}{\svgwidth}%
  \fi%
  \global\let\svgwidth\undefined%
  \global\let\svgscale\undefined%
  \makeatother%
  \begin{picture}(1,1.05060582)%
    \lineheight{1}%
    \setlength\tabcolsep{0pt}%
    \put(0,0){\includegraphics[width=\unitlength]{associativity_ssbim_2_svg-tex.eps}}%
    \put(0.33823058,0.21901026){\color[rgb]{0,0,0}\makebox(0,0)[lt]{\lineheight{1.25}\smash{\begin{tabular}[t]{l}\dgr{l+m}\end{tabular}}}}%
    \put(0.64423524,0.02998717){\color[rgb]{0,0,0}\makebox(0,0)[lt]{\lineheight{1.25}\smash{\begin{tabular}[t]{l}\dgr{m}\end{tabular}}}}%
    \put(0.78330071,0.81070158){\color[rgb]{0,0,0}\makebox(0,0)[lt]{\lineheight{1.25}\smash{\begin{tabular}[t]{l}\dgr{0}\end{tabular}}}}%
    \put(0.08089331,0.80904963){\color[rgb]{0,0,0}\makebox(0,0)[lt]{\lineheight{1.25}\smash{\begin{tabular}[t]{l}\dgr{k+l+m}\end{tabular}}}}%
  \end{picture}%
\endgroup%

%% file: arxiv-figures/associativity_ssbim_3_svg-tex.eps_tex
\begingroup%
  \makeatletter%
  \providecommand\color[2][]{%
    \errmessage{(Inkscape) Color is used for the text in Inkscape, but the package 'color.sty' is not loaded}%
    \renewcommand\color[2][]{}%
  }%
  \providecommand\transparent[1]{%
    \errmessage{(Inkscape) Transparency is used (non-zero) for the text in Inkscape, but the package 'transparent.sty' is not loaded}%
    \renewcommand\transparent[1]{}%
  }%
  \providecommand\rotatebox[2]{#2}%
  \newcommand*\fsize{\dimexpr\f@size pt\relax}%
  \newcommand*\lineheight[1]{\fontsize{\fsize}{#1\fsize}\selectfont}%
  \ifx\svgwidth\undefined%
    \setlength{\unitlength}{677.65956526bp}%
    \ifx\svgscale\undefined%
      \relax%
    \else%
      \setlength{\unitlength}{\unitlength * \real{\svgscale}}%
    \fi%
  \else%
    \setlength{\unitlength}{\svgwidth}%
  \fi%
  \global\let\svgwidth\undefined%
  \global\let\svgscale\undefined%
  \makeatother%
  \begin{picture}(1,1.05060582)%
    \lineheight{1}%
    \setlength\tabcolsep{0pt}%
    \put(0,0){\includegraphics[width=\unitlength]{associativity_ssbim_3_svg-tex.eps}}%
    \put(0.33823058,0.21901026){\color[rgb]{0,0,0}\makebox(0,0)[lt]{\lineheight{1.25}\smash{\begin{tabular}[t]{l}\dgr{l+m}\end{tabular}}}}%
    \put(0.64423524,0.02998717){\color[rgb]{0,0,0}\makebox(0,0)[lt]{\lineheight{1.25}\smash{\begin{tabular}[t]{l}\dgr{m}\end{tabular}}}}%
    \put(0.78330071,0.81070158){\color[rgb]{0,0,0}\makebox(0,0)[lt]{\lineheight{1.25}\smash{\begin{tabular}[t]{l}\dgr{0}\end{tabular}}}}%
    \put(0.08089331,0.80904963){\color[rgb]{0,0,0}\makebox(0,0)[lt]{\lineheight{1.25}\smash{\begin{tabular}[t]{l}\dgr{k+l+m}\end{tabular}}}}%
  \end{picture}%
\endgroup%

%% file: arxiv-figures/associativity_ssbim_4_svg-tex.eps_tex
\begingroup%
  \makeatletter%
  \providecommand\color[2][]{%
    \errmessage{(Inkscape) Color is used for the text in Inkscape, but the package 'color.sty' is not loaded}%
    \renewcommand\color[2][]{}%
  }%
  \providecommand\transparent[1]{%
    \errmessage{(Inkscape) Transparency is used (non-zero) for the text in Inkscape, but the package 'transparent.sty' is not loaded}%
    \renewcommand\transparent[1]{}%
  }%
  \providecommand\rotatebox[2]{#2}%
  \newcommand*\fsize{\dimexpr\f@size pt\relax}%
  \newcommand*\lineheight[1]{\fontsize{\fsize}{#1\fsize}\selectfont}%
  \ifx\svgwidth\undefined%
    \setlength{\unitlength}{677.65956526bp}%
    \ifx\svgscale\undefined%
      \relax%
    \else%
      \setlength{\unitlength}{\unitlength * \real{\svgscale}}%
    \fi%
  \else%
    \setlength{\unitlength}{\svgwidth}%
  \fi%
  \global\let\svgwidth\undefined%
  \global\let\svgscale\undefined%
  \makeatother%
  \begin{picture}(1,1.05060582)%
    \lineheight{1}%
    \setlength\tabcolsep{0pt}%
    \put(0,0){\includegraphics[width=\unitlength]{associativity_ssbim_4_svg-tex.eps}}%
    \put(0.33823058,0.21901026){\color[rgb]{0,0,0}\makebox(0,0)[lt]{\lineheight{1.25}\smash{\begin{tabular}[t]{l}\dgr{l+m}\end{tabular}}}}%
    \put(0.64423524,0.02998717){\color[rgb]{0,0,0}\makebox(0,0)[lt]{\lineheight{1.25}\smash{\begin{tabular}[t]{l}\dgr{m}\end{tabular}}}}%
    \put(0.78330071,0.81070158){\color[rgb]{0,0,0}\makebox(0,0)[lt]{\lineheight{1.25}\smash{\begin{tabular}[t]{l}\dgr{0}\end{tabular}}}}%
    \put(0.08089331,0.80904963){\color[rgb]{0,0,0}\makebox(0,0)[lt]{\lineheight{1.25}\smash{\begin{tabular}[t]{l}\dgr{k+l+m}\end{tabular}}}}%
  \end{picture}%
\endgroup%

%% file: arxiv-figures/associativity_ssbim_5_svg-tex.eps_tex
\begingroup%
  \makeatletter%
  \providecommand\color[2][]{%
    \errmessage{(Inkscape) Color is used for the text in Inkscape, but the package 'color.sty' is not loaded}%
    \renewcommand\color[2][]{}%
  }%
  \providecommand\transparent[1]{%
    \errmessage{(Inkscape) Transparency is used (non-zero) for the text in Inkscape, but the package 'transparent.sty' is not loaded}%
    \renewcommand\transparent[1]{}%
  }%
  \providecommand\rotatebox[2]{#2}%
  \newcommand*\fsize{\dimexpr\f@size pt\relax}%
  \newcommand*\lineheight[1]{\fontsize{\fsize}{#1\fsize}\selectfont}%
  \ifx\svgwidth\undefined%
    \setlength{\unitlength}{677.65956526bp}%
    \ifx\svgscale\undefined%
      \relax%
    \else%
      \setlength{\unitlength}{\unitlength * \real{\svgscale}}%
    \fi%
  \else%
    \setlength{\unitlength}{\svgwidth}%
  \fi%
  \global\let\svgwidth\undefined%
  \global\let\svgscale\undefined%
  \makeatother%
  \begin{picture}(1,1.05060582)%
    \lineheight{1}%
    \setlength\tabcolsep{0pt}%
    \put(0,0){\includegraphics[width=\unitlength]{associativity_ssbim_5_svg-tex.eps}}%
    \put(0.33823058,0.21901026){\color[rgb]{0,0,0}\makebox(0,0)[lt]{\lineheight{1.25}\smash{\begin{tabular}[t]{l}\dgr{l+m}\end{tabular}}}}%
    \put(0.64423524,0.02998717){\color[rgb]{0,0,0}\makebox(0,0)[lt]{\lineheight{1.25}\smash{\begin{tabular}[t]{l}\dgr{m}\end{tabular}}}}%
    \put(0.78330071,0.81070158){\color[rgb]{0,0,0}\makebox(0,0)[lt]{\lineheight{1.25}\smash{\begin{tabular}[t]{l}\dgr{0}\end{tabular}}}}%
    \put(0.08089331,0.80904963){\color[rgb]{0,0,0}\makebox(0,0)[lt]{\lineheight{1.25}\smash{\begin{tabular}[t]{l}\dgr{k+l+m}\end{tabular}}}}%
  \end{picture}%
\endgroup%

%% file: arxiv-figures/associativity_ssbim_6_svg-tex.eps_tex
\begingroup%
  \makeatletter%
  \providecommand\color[2][]{%
    \errmessage{(Inkscape) Color is used for the text in Inkscape, but the package 'color.sty' is not loaded}%
    \renewcommand\color[2][]{}%
  }%
  \providecommand\transparent[1]{%
    \errmessage{(Inkscape) Transparency is used (non-zero) for the text in Inkscape, but the package 'transparent.sty' is not loaded}%
    \renewcommand\transparent[1]{}%
  }%
  \providecommand\rotatebox[2]{#2}%
  \newcommand*\fsize{\dimexpr\f@size pt\relax}%
  \newcommand*\lineheight[1]{\fontsize{\fsize}{#1\fsize}\selectfont}%
  \ifx\svgwidth\undefined%
    \setlength{\unitlength}{677.65956526bp}%
    \ifx\svgscale\undefined%
      \relax%
    \else%
      \setlength{\unitlength}{\unitlength * \real{\svgscale}}%
    \fi%
  \else%
    \setlength{\unitlength}{\svgwidth}%
  \fi%
  \global\let\svgwidth\undefined%
  \global\let\svgscale\undefined%
  \makeatother%
  \begin{picture}(1,1.05060582)%
    \lineheight{1}%
    \setlength\tabcolsep{0pt}%
    \put(0,0){\includegraphics[width=\unitlength]{associativity_ssbim_6_svg-tex.eps}}%
    \put(0.33823058,0.21901026){\color[rgb]{0,0,0}\makebox(0,0)[lt]{\lineheight{1.25}\smash{\begin{tabular}[t]{l}\dgr{l+m}\end{tabular}}}}%
    \put(0.64423524,0.02998717){\color[rgb]{0,0,0}\makebox(0,0)[lt]{\lineheight{1.25}\smash{\begin{tabular}[t]{l}\dgr{m}\end{tabular}}}}%
    \put(0.78330071,0.81070158){\color[rgb]{0,0,0}\makebox(0,0)[lt]{\lineheight{1.25}\smash{\begin{tabular}[t]{l}\dgr{0}\end{tabular}}}}%
    \put(0.08089331,0.80904963){\color[rgb]{0,0,0}\makebox(0,0)[lt]{\lineheight{1.25}\smash{\begin{tabular}[t]{l}\dgr{k+l+m}\end{tabular}}}}%
  \end{picture}%
\endgroup%

%% file: arxiv-figures/image_of_bigon_svg-tex.eps_tex
\begingroup%
  \makeatletter%
  \providecommand\color[2][]{%
    \errmessage{(Inkscape) Color is used for the text in Inkscape, but the package 'color.sty' is not loaded}%
    \renewcommand\color[2][]{}%
  }%
  \providecommand\transparent[1]{%
    \errmessage{(Inkscape) Transparency is used (non-zero) for the text in Inkscape, but the package 'transparent.sty' is not loaded}%
    \renewcommand\transparent[1]{}%
  }%
  \providecommand\rotatebox[2]{#2}%
  \newcommand*\fsize{\dimexpr\f@size pt\relax}%
  \newcommand*\lineheight[1]{\fontsize{\fsize}{#1\fsize}\selectfont}%
  \ifx\svgwidth\undefined%
    \setlength{\unitlength}{583.56432731bp}%
    \ifx\svgscale\undefined%
      \relax%
    \else%
      \setlength{\unitlength}{\unitlength * \real{\svgscale}}%
    \fi%
  \else%
    \setlength{\unitlength}{\svgwidth}%
  \fi%
  \global\let\svgwidth\undefined%
  \global\let\svgscale\undefined%
  \makeatother%
  \begin{picture}(1,1.4489507)%
    \lineheight{1}%
    \setlength\tabcolsep{0pt}%
    \put(0,0){\includegraphics[width=\unitlength]{image_of_bigon_svg-tex.eps}}%
    \put(0.76145356,0.20273508){\color[rgb]{0,0,0}\makebox(0,0)[lt]{\lineheight{1.25}\smash{\begin{tabular}[t]{l}\dgr{a}\end{tabular}}}}%
    \put(0.01717356,0.19651644){\color[rgb]{0,0,0}\makebox(0,0)[lt]{\lineheight{1.25}\smash{\begin{tabular}[t]{l}\dgr{a+k}\end{tabular}}}}%
    \put(0.37370513,0.71474429){\color[rgb]{0,0,0}\makebox(0,0)[lt]{\lineheight{1.25}\smash{\begin{tabular}[t]{l}\dgr{a+1}\end{tabular}}}}%
    \put(0.75081674,1.23073522){\color[rgb]{0,0,0}\makebox(0,0)[lt]{\lineheight{1.25}\smash{\begin{tabular}[t]{l}\dgr{a}\end{tabular}}}}%
    \put(0.02709222,1.22451649){\color[rgb]{0,0,0}\makebox(0,0)[lt]{\lineheight{1.25}\smash{\begin{tabular}[t]{l}\dgr{a+k}\end{tabular}}}}%
  \end{picture}%
\endgroup%

%% file: arxiv-figures/ssbim_bigon_simplify_1_svg-tex.eps_tex
\begingroup%
  \makeatletter%
  \providecommand\color[2][]{%
    \errmessage{(Inkscape) Color is used for the text in Inkscape, but the package 'color.sty' is not loaded}%
    \renewcommand\color[2][]{}%
  }%
  \providecommand\transparent[1]{%
    \errmessage{(Inkscape) Transparency is used (non-zero) for the text in Inkscape, but the package 'transparent.sty' is not loaded}%
    \renewcommand\transparent[1]{}%
  }%
  \providecommand\rotatebox[2]{#2}%
  \newcommand*\fsize{\dimexpr\f@size pt\relax}%
  \newcommand*\lineheight[1]{\fontsize{\fsize}{#1\fsize}\selectfont}%
  \ifx\svgwidth\undefined%
    \setlength{\unitlength}{583.0181444bp}%
    \ifx\svgscale\undefined%
      \relax%
    \else%
      \setlength{\unitlength}{\unitlength * \real{\svgscale}}%
    \fi%
  \else%
    \setlength{\unitlength}{\svgwidth}%
  \fi%
  \global\let\svgwidth\undefined%
  \global\let\svgscale\undefined%
  \makeatother%
  \begin{picture}(1,1.44937129)%
    \lineheight{1}%
    \setlength\tabcolsep{0pt}%
    \put(0,0){\includegraphics[width=\unitlength]{ssbim_bigon_simplify_1_svg-tex.eps}}%
    \put(0.76169849,0.20245659){\color[rgb]{0,0,0}\makebox(0,0)[lt]{\lineheight{1.25}\smash{\begin{tabular}[t]{l}\dgr{0}\end{tabular}}}}%
    \put(0.1080148,0.20660639){\color[rgb]{0,0,0}\makebox(0,0)[lt]{\lineheight{1.25}\smash{\begin{tabular}[t]{l}\dgr{k}\end{tabular}}}}%
    \put(0.75105171,1.23141979){\color[rgb]{0,0,0}\makebox(0,0)[lt]{\lineheight{1.25}\smash{\begin{tabular}[t]{l}\dgr{0}\end{tabular}}}}%
    \put(0.10756849,1.23141979){\color[rgb]{0,0,0}\makebox(0,0)[lt]{\lineheight{1.25}\smash{\begin{tabular}[t]{l}\dgr{k}\end{tabular}}}}%
    \put(0.34782578,0.7252201){\color[rgb]{0,0,0}\makebox(0,0)[lt]{\lineheight{1.25}\smash{\begin{tabular}[t]{l}\blk{\mu^{\hh{0,1},\hh{1,k}}_{\hh{1}}}\end{tabular}}}}%
  \end{picture}%
\endgroup%

%% file: arxiv-figures/ssbim_bigon_simplify_2_svg-tex.eps_tex
\begingroup%
  \makeatletter%
  \providecommand\color[2][]{%
    \errmessage{(Inkscape) Color is used for the text in Inkscape, but the package 'color.sty' is not loaded}%
    \renewcommand\color[2][]{}%
  }%
  \providecommand\transparent[1]{%
    \errmessage{(Inkscape) Transparency is used (non-zero) for the text in Inkscape, but the package 'transparent.sty' is not loaded}%
    \renewcommand\transparent[1]{}%
  }%
  \providecommand\rotatebox[2]{#2}%
  \newcommand*\fsize{\dimexpr\f@size pt\relax}%
  \newcommand*\lineheight[1]{\fontsize{\fsize}{#1\fsize}\selectfont}%
  \ifx\svgwidth\undefined%
    \setlength{\unitlength}{649.81886894bp}%
    \ifx\svgscale\undefined%
      \relax%
    \else%
      \setlength{\unitlength}{\unitlength * \real{\svgscale}}%
    \fi%
  \else%
    \setlength{\unitlength}{\svgwidth}%
  \fi%
  \global\let\svgwidth\undefined%
  \global\let\svgscale\undefined%
  \makeatother%
  \begin{picture}(1,1.30037738)%
    \lineheight{1}%
    \setlength\tabcolsep{0pt}%
    \put(0,0){\includegraphics[width=\unitlength]{ssbim_bigon_simplify_2_svg-tex.eps}}%
    \put(0.82566073,0.17762733){\color[rgb]{0,0,0}\makebox(0,0)[lt]{\lineheight{1.25}\smash{\begin{tabular}[t]{l}\dgr{0}\end{tabular}}}}%
    \put(0.09397334,0.17576583){\color[rgb]{0,0,0}\makebox(0,0)[lt]{\lineheight{1.25}\smash{\begin{tabular}[t]{l}\dgr{k}\end{tabular}}}}%
    \put(0.82541625,1.1026757){\color[rgb]{0,0,0}\makebox(0,0)[lt]{\lineheight{1.25}\smash{\begin{tabular}[t]{l}\dgr{0}\end{tabular}}}}%
    \put(0.09357291,1.10081416){\color[rgb]{0,0,0}\makebox(0,0)[lt]{\lineheight{1.25}\smash{\begin{tabular}[t]{l}\dgr{k}\end{tabular}}}}%
    \put(0.38942771,0.65037439){\color[rgb]{0,0,0}\makebox(0,0)[lt]{\lineheight{1.25}\smash{\begin{tabular}[t]{l}\blk{\mu^{\hh{0,1},\hh{1,k}}_{\hh{1}}}\end{tabular}}}}%
  \end{picture}%
\endgroup%

%% file: arxiv-figures/ssbim_bigon_simplify_3_svg-tex.eps_tex
\begingroup%
  \makeatletter%
  \providecommand\color[2][]{%
    \errmessage{(Inkscape) Color is used for the text in Inkscape, but the package 'color.sty' is not loaded}%
    \renewcommand\color[2][]{}%
  }%
  \providecommand\transparent[1]{%
    \errmessage{(Inkscape) Transparency is used (non-zero) for the text in Inkscape, but the package 'transparent.sty' is not loaded}%
    \renewcommand\transparent[1]{}%
  }%
  \providecommand\rotatebox[2]{#2}%
  \newcommand*\fsize{\dimexpr\f@size pt\relax}%
  \newcommand*\lineheight[1]{\fontsize{\fsize}{#1\fsize}\selectfont}%
  \ifx\svgwidth\undefined%
    \setlength{\unitlength}{739.81886894bp}%
    \ifx\svgscale\undefined%
      \relax%
    \else%
      \setlength{\unitlength}{\unitlength * \real{\svgscale}}%
    \fi%
  \else%
    \setlength{\unitlength}{\svgwidth}%
  \fi%
  \global\let\svgwidth\undefined%
  \global\let\svgscale\undefined%
  \makeatother%
  \begin{picture}(1,1.14218466)%
    \lineheight{1}%
    \setlength\tabcolsep{0pt}%
    \put(0,0){\includegraphics[width=\unitlength]{ssbim_bigon_simplify_3_svg-tex.eps}}%
    \put(0.84850446,0.15562278){\color[rgb]{0,0,0}\makebox(0,0)[lt]{\lineheight{1.25}\smash{\begin{tabular}[t]{l}\dgr{0}\end{tabular}}}}%
    \put(0.08090624,0.16281724){\color[rgb]{0,0,0}\makebox(0,0)[lt]{\lineheight{1.25}\smash{\begin{tabular}[t]{l}\dgr{k}\end{tabular}}}}%
    \put(0.84828972,0.9681377){\color[rgb]{0,0,0}\makebox(0,0)[lt]{\lineheight{1.25}\smash{\begin{tabular}[t]{l}\dgr{0}\end{tabular}}}}%
    \put(0.08382472,0.96879173){\color[rgb]{0,0,0}\makebox(0,0)[lt]{\lineheight{1.25}\smash{\begin{tabular}[t]{l}\dgr{k}\end{tabular}}}}%
    \put(0.27697632,0.57478369){\color[rgb]{0,0,0}\makebox(0,0)[lt]{\lineheight{1.25}\smash{\begin{tabular}[t]{l}\blk{\dd^{\hh{0,1,k}}_{\hh{0,k}}\left(\mu^{\hh{0,1},\hh{1,k}}_{\hh{1}}\right)}\end{tabular}}}}%
  \end{picture}%
\endgroup%

%% file: arxiv-figures/image_of_bigon_k=n_svg-tex.eps_tex
\begingroup%
  \makeatletter%
  \providecommand\color[2][]{%
    \errmessage{(Inkscape) Color is used for the text in Inkscape, but the package 'color.sty' is not loaded}%
    \renewcommand\color[2][]{}%
  }%
  \providecommand\transparent[1]{%
    \errmessage{(Inkscape) Transparency is used (non-zero) for the text in Inkscape, but the package 'transparent.sty' is not loaded}%
    \renewcommand\transparent[1]{}%
  }%
  \providecommand\rotatebox[2]{#2}%
  \newcommand*\fsize{\dimexpr\f@size pt\relax}%
  \newcommand*\lineheight[1]{\fontsize{\fsize}{#1\fsize}\selectfont}%
  \ifx\svgwidth\undefined%
    \setlength{\unitlength}{583.0181444bp}%
    \ifx\svgscale\undefined%
      \relax%
    \else%
      \setlength{\unitlength}{\unitlength * \real{\svgscale}}%
    \fi%
  \else%
    \setlength{\unitlength}{\svgwidth}%
  \fi%
  \global\let\svgwidth\undefined%
  \global\let\svgscale\undefined%
  \makeatother%
  \begin{picture}(1,1.44937129)%
    \lineheight{1}%
    \setlength\tabcolsep{0pt}%
    \put(0,0){\includegraphics[width=\unitlength]{image_of_bigon_k=n_svg-tex.eps}}%
    \put(0.44632072,0.05929174){\color[rgb]{0,0,0}\makebox(0,0)[lt]{\lineheight{1.25}\smash{\begin{tabular}[t]{l}\dgr{a}\end{tabular}}}}%
    \put(0.37358681,0.71494547){\color[rgb]{0,0,0}\makebox(0,0)[lt]{\lineheight{1.25}\smash{\begin{tabular}[t]{l}\dgr{a+1}\end{tabular}}}}%
    \put(0.45019792,1.32271335){\color[rgb]{0,0,0}\makebox(0,0)[lt]{\lineheight{1.25}\smash{\begin{tabular}[t]{l}\dgr{a}\end{tabular}}}}%
  \end{picture}%
\endgroup%

%% file: arxiv-figures/image_of_square_flop_LHS_svg-tex.eps_tex
\begingroup%
  \makeatletter%
  \providecommand\color[2][]{%
    \errmessage{(Inkscape) Color is used for the text in Inkscape, but the package 'color.sty' is not loaded}%
    \renewcommand\color[2][]{}%
  }%
  \providecommand\transparent[1]{%
    \errmessage{(Inkscape) Transparency is used (non-zero) for the text in Inkscape, but the package 'transparent.sty' is not loaded}%
    \renewcommand\transparent[1]{}%
  }%
  \providecommand\rotatebox[2]{#2}%
  \newcommand*\fsize{\dimexpr\f@size pt\relax}%
  \newcommand*\lineheight[1]{\fontsize{\fsize}{#1\fsize}\selectfont}%
  \ifx\svgwidth\undefined%
    \setlength{\unitlength}{664.81888336bp}%
    \ifx\svgscale\undefined%
      \relax%
    \else%
      \setlength{\unitlength}{\unitlength * \real{\svgscale}}%
    \fi%
  \else%
    \setlength{\unitlength}{\svgwidth}%
  \fi%
  \global\let\svgwidth\undefined%
  \global\let\svgscale\undefined%
  \makeatother%
  \begin{picture}(1,1.27255604)%
    \lineheight{1}%
    \setlength\tabcolsep{0pt}%
    \put(0,0){\includegraphics[width=\unitlength]{image_of_square_flop_LHS_svg-tex.eps}}%
    \put(0.83867211,0.26488465){\color[rgb]{0,0,0}\makebox(0,0)[lt]{\lineheight{1.25}\smash{\begin{tabular}[t]{l}\dgr{0}\end{tabular}}}}%
    \put(0.04707195,0.26852387){\color[rgb]{0,0,0}\makebox(0,0)[lt]{\lineheight{1.25}\smash{\begin{tabular}[t]{l}\dgr{k+1}\end{tabular}}}}%
    \put(0.46010377,0.62333791){\color[rgb]{0,0,0}\makebox(0,0)[lt]{\lineheight{1.25}\smash{\begin{tabular}[t]{l}\dgr{2}\end{tabular}}}}%
    \put(0.83661353,0.99074486){\color[rgb]{0,0,0}\makebox(0,0)[lt]{\lineheight{1.25}\smash{\begin{tabular}[t]{l}\dgr{0}\end{tabular}}}}%
    \put(0.04850009,0.98892528){\color[rgb]{0,0,0}\makebox(0,0)[lt]{\lineheight{1.25}\smash{\begin{tabular}[t]{l}\dgr{k+1}\end{tabular}}}}%
    \put(0.45246163,1.19104044){\color[rgb]{0,0,0}\makebox(0,0)[lt]{\lineheight{1.25}\smash{\begin{tabular}[t]{l}\dgr{1}\end{tabular}}}}%
    \put(0.46738203,0.05272408){\color[rgb]{0,0,0}\makebox(0,0)[lt]{\lineheight{1.25}\smash{\begin{tabular}[t]{l}\dgr{1}\end{tabular}}}}%
  \end{picture}%
\endgroup%

%% file: arxiv-figures/image_of_square_flop_RHS_1_svg-tex.eps_tex
\begingroup%
  \makeatletter%
  \providecommand\color[2][]{%
    \errmessage{(Inkscape) Color is used for the text in Inkscape, but the package 'color.sty' is not loaded}%
    \renewcommand\color[2][]{}%
  }%
  \providecommand\transparent[1]{%
    \errmessage{(Inkscape) Transparency is used (non-zero) for the text in Inkscape, but the package 'transparent.sty' is not loaded}%
    \renewcommand\transparent[1]{}%
  }%
  \providecommand\rotatebox[2]{#2}%
  \newcommand*\fsize{\dimexpr\f@size pt\relax}%
  \newcommand*\lineheight[1]{\fontsize{\fsize}{#1\fsize}\selectfont}%
  \ifx\svgwidth\undefined%
    \setlength{\unitlength}{469.83444011bp}%
    \ifx\svgscale\undefined%
      \relax%
    \else%
      \setlength{\unitlength}{\unitlength * \real{\svgscale}}%
    \fi%
  \else%
    \setlength{\unitlength}{\svgwidth}%
  \fi%
  \global\let\svgwidth\undefined%
  \global\let\svgscale\undefined%
  \makeatother%
  \begin{picture}(1,1.35115455)%
    \lineheight{1}%
    \setlength\tabcolsep{0pt}%
    \put(0,0){\includegraphics[width=\unitlength]{image_of_square_flop_RHS_1_svg-tex.eps}}%
    \put(0.73212244,0.66740694){\color[rgb]{0,0,0}\makebox(0,0)[lt]{\lineheight{1.25}\smash{\begin{tabular}[t]{l}\dgr{0}\end{tabular}}}}%
    \put(0.07643495,0.66483219){\color[rgb]{0,0,0}\makebox(0,0)[lt]{\lineheight{1.25}\smash{\begin{tabular}[t]{l}\dgr{k+1}\end{tabular}}}}%
    \put(0.47205467,0.09325123){\color[rgb]{0,0,0}\makebox(0,0)[lt]{\lineheight{1.25}\smash{\begin{tabular}[t]{l}\dgr{1}\end{tabular}}}}%
    \put(0.4669053,1.20345723){\color[rgb]{0,0,0}\makebox(0,0)[lt]{\lineheight{1.25}\smash{\begin{tabular}[t]{l}\dgr{1}\end{tabular}}}}%
  \end{picture}%
\endgroup%

%% file: arxiv-figures/image_of_square_flop_RHS_2_svg-tex.eps_tex
\begingroup%
  \makeatletter%
  \providecommand\color[2][]{%
    \errmessage{(Inkscape) Color is used for the text in Inkscape, but the package 'color.sty' is not loaded}%
    \renewcommand\color[2][]{}%
  }%
  \providecommand\transparent[1]{%
    \errmessage{(Inkscape) Transparency is used (non-zero) for the text in Inkscape, but the package 'transparent.sty' is not loaded}%
    \renewcommand\transparent[1]{}%
  }%
  \providecommand\rotatebox[2]{#2}%
  \newcommand*\fsize{\dimexpr\f@size pt\relax}%
  \newcommand*\lineheight[1]{\fontsize{\fsize}{#1\fsize}\selectfont}%
  \ifx\svgwidth\undefined%
    \setlength{\unitlength}{469.83444011bp}%
    \ifx\svgscale\undefined%
      \relax%
    \else%
      \setlength{\unitlength}{\unitlength * \real{\svgscale}}%
    \fi%
  \else%
    \setlength{\unitlength}{\svgwidth}%
  \fi%
  \global\let\svgwidth\undefined%
  \global\let\svgscale\undefined%
  \makeatother%
  \begin{picture}(1,1.35115455)%
    \lineheight{1}%
    \setlength\tabcolsep{0pt}%
    \put(0,0){\includegraphics[width=\unitlength]{image_of_square_flop_RHS_2_svg-tex.eps}}%
    \put(0.81966189,0.6699816){\color[rgb]{0,0,0}\makebox(0,0)[lt]{\lineheight{1.25}\smash{\begin{tabular}[t]{l}\dgr{0}\end{tabular}}}}%
    \put(0.02236649,0.67255626){\color[rgb]{0,0,0}\makebox(0,0)[lt]{\lineheight{1.25}\smash{\begin{tabular}[t]{l}\dgr{k+1}\end{tabular}}}}%
    \put(0.45145716,0.67307134){\color[rgb]{0,0,0}\makebox(0,0)[lt]{\lineheight{1.25}\smash{\begin{tabular}[t]{l}\dgr{1}\end{tabular}}}}%
  \end{picture}%
\endgroup%

%% file: arxiv-figures/image_of_square_flop_LHS_for_n-1_svg-tex.eps_tex
\begingroup%
  \makeatletter%
  \providecommand\color[2][]{%
    \errmessage{(Inkscape) Color is used for the text in Inkscape, but the package 'color.sty' is not loaded}%
    \renewcommand\color[2][]{}%
  }%
  \providecommand\transparent[1]{%
    \errmessage{(Inkscape) Transparency is used (non-zero) for the text in Inkscape, but the package 'transparent.sty' is not loaded}%
    \renewcommand\transparent[1]{}%
  }%
  \providecommand\rotatebox[2]{#2}%
  \newcommand*\fsize{\dimexpr\f@size pt\relax}%
  \newcommand*\lineheight[1]{\fontsize{\fsize}{#1\fsize}\selectfont}%
  \ifx\svgwidth\undefined%
    \setlength{\unitlength}{664.81888336bp}%
    \ifx\svgscale\undefined%
      \relax%
    \else%
      \setlength{\unitlength}{\unitlength * \real{\svgscale}}%
    \fi%
  \else%
    \setlength{\unitlength}{\svgwidth}%
  \fi%
  \global\let\svgwidth\undefined%
  \global\let\svgscale\undefined%
  \makeatother%
  \begin{picture}(1,1.27255604)%
    \lineheight{1}%
    \setlength\tabcolsep{0pt}%
    \put(0,0){\includegraphics[width=\unitlength]{image_of_square_flop_LHS_for_n-1_svg-tex.eps}}%
    \put(0.83867211,0.26488465){\color[rgb]{0,0,0}\makebox(0,0)[lt]{\lineheight{1.25}\smash{\begin{tabular}[t]{l}\dgr{0}\end{tabular}}}}%
    \put(0.07982402,0.26488477){\color[rgb]{0,0,0}\makebox(0,0)[lt]{\lineheight{1.25}\smash{\begin{tabular}[t]{l}\dgr{0}\end{tabular}}}}%
    \put(0.46010377,0.62333791){\color[rgb]{0,0,0}\makebox(0,0)[lt]{\lineheight{1.25}\smash{\begin{tabular}[t]{l}\dgr{2}\end{tabular}}}}%
    \put(0.83661353,0.99074486){\color[rgb]{0,0,0}\makebox(0,0)[lt]{\lineheight{1.25}\smash{\begin{tabular}[t]{l}\dgr{0}\end{tabular}}}}%
    \put(0.09216952,0.98892528){\color[rgb]{0,0,0}\makebox(0,0)[lt]{\lineheight{1.25}\smash{\begin{tabular}[t]{l}\dgr{0}\end{tabular}}}}%
    \put(0.45246163,1.19104044){\color[rgb]{0,0,0}\makebox(0,0)[lt]{\lineheight{1.25}\smash{\begin{tabular}[t]{l}\dgr{1}\end{tabular}}}}%
    \put(0.46738203,0.05272408){\color[rgb]{0,0,0}\makebox(0,0)[lt]{\lineheight{1.25}\smash{\begin{tabular}[t]{l}\dgr{1}\end{tabular}}}}%
  \end{picture}%
\endgroup%

%% file: arxiv-figures/image_of_square_flop_RHS_1_for_n-1_svg-tex.eps_tex
\begingroup%
  \makeatletter%
  \providecommand\color[2][]{%
    \errmessage{(Inkscape) Color is used for the text in Inkscape, but the package 'color.sty' is not loaded}%
    \renewcommand\color[2][]{}%
  }%
  \providecommand\transparent[1]{%
    \errmessage{(Inkscape) Transparency is used (non-zero) for the text in Inkscape, but the package 'transparent.sty' is not loaded}%
    \renewcommand\transparent[1]{}%
  }%
  \providecommand\rotatebox[2]{#2}%
  \newcommand*\fsize{\dimexpr\f@size pt\relax}%
  \newcommand*\lineheight[1]{\fontsize{\fsize}{#1\fsize}\selectfont}%
  \ifx\svgwidth\undefined%
    \setlength{\unitlength}{469.83444011bp}%
    \ifx\svgscale\undefined%
      \relax%
    \else%
      \setlength{\unitlength}{\unitlength * \real{\svgscale}}%
    \fi%
  \else%
    \setlength{\unitlength}{\svgwidth}%
  \fi%
  \global\let\svgwidth\undefined%
  \global\let\svgscale\undefined%
  \makeatother%
  \begin{picture}(1,1.35115455)%
    \lineheight{1}%
    \setlength\tabcolsep{0pt}%
    \put(0,0){\includegraphics[width=\unitlength]{image_of_square_flop_RHS_1_for_n-1_svg-tex.eps}}%
    \put(0.44890667,0.62363717){\color[rgb]{0,0,0}\makebox(0,0)[lt]{\lineheight{1.25}\smash{\begin{tabular}[t]{l}\dgr{0}\end{tabular}}}}%
    \put(0.44888245,0.09582598){\color[rgb]{0,0,0}\makebox(0,0)[lt]{\lineheight{1.25}\smash{\begin{tabular}[t]{l}\dgr{1}\end{tabular}}}}%
    \put(0.45145716,1.17771036){\color[rgb]{0,0,0}\makebox(0,0)[lt]{\lineheight{1.25}\smash{\begin{tabular}[t]{l}\dgr{1}\end{tabular}}}}%
  \end{picture}%
\endgroup%

%% file: arxiv-figures/image_of_square_flop_RHS_2_for_n-1_svg-tex.eps_tex
\begingroup%
  \makeatletter%
  \providecommand\color[2][]{%
    \errmessage{(Inkscape) Color is used for the text in Inkscape, but the package 'color.sty' is not loaded}%
    \renewcommand\color[2][]{}%
  }%
  \providecommand\transparent[1]{%
    \errmessage{(Inkscape) Transparency is used (non-zero) for the text in Inkscape, but the package 'transparent.sty' is not loaded}%
    \renewcommand\transparent[1]{}%
  }%
  \providecommand\rotatebox[2]{#2}%
  \newcommand*\fsize{\dimexpr\f@size pt\relax}%
  \newcommand*\lineheight[1]{\fontsize{\fsize}{#1\fsize}\selectfont}%
  \ifx\svgwidth\undefined%
    \setlength{\unitlength}{469.83444011bp}%
    \ifx\svgscale\undefined%
      \relax%
    \else%
      \setlength{\unitlength}{\unitlength * \real{\svgscale}}%
    \fi%
  \else%
    \setlength{\unitlength}{\svgwidth}%
  \fi%
  \global\let\svgwidth\undefined%
  \global\let\svgscale\undefined%
  \makeatother%
  \begin{picture}(1,1.35115455)%
    \lineheight{1}%
    \setlength\tabcolsep{0pt}%
    \put(0,0){\includegraphics[width=\unitlength]{image_of_square_flop_RHS_2_for_n-1_svg-tex.eps}}%
    \put(0.81966189,0.6699816){\color[rgb]{0,0,0}\makebox(0,0)[lt]{\lineheight{1.25}\smash{\begin{tabular}[t]{l}\dgr{0}\end{tabular}}}}%
    \put(0.09445777,0.67255626){\color[rgb]{0,0,0}\makebox(0,0)[lt]{\lineheight{1.25}\smash{\begin{tabular}[t]{l}\dgr{0}\end{tabular}}}}%
    \put(0.45145716,0.67307134){\color[rgb]{0,0,0}\makebox(0,0)[lt]{\lineheight{1.25}\smash{\begin{tabular}[t]{l}\dgr{1}\end{tabular}}}}%
  \end{picture}%
\endgroup%

%% file: arxiv-figures/image_of_square_flop_k=n_svg-tex.eps_tex
\begingroup%
  \makeatletter%
  \providecommand\color[2][]{%
    \errmessage{(Inkscape) Color is used for the text in Inkscape, but the package 'color.sty' is not loaded}%
    \renewcommand\color[2][]{}%
  }%
  \providecommand\transparent[1]{%
    \errmessage{(Inkscape) Transparency is used (non-zero) for the text in Inkscape, but the package 'transparent.sty' is not loaded}%
    \renewcommand\transparent[1]{}%
  }%
  \providecommand\rotatebox[2]{#2}%
  \newcommand*\fsize{\dimexpr\f@size pt\relax}%
  \newcommand*\lineheight[1]{\fontsize{\fsize}{#1\fsize}\selectfont}%
  \ifx\svgwidth\undefined%
    \setlength{\unitlength}{583.56432731bp}%
    \ifx\svgscale\undefined%
      \relax%
    \else%
      \setlength{\unitlength}{\unitlength * \real{\svgscale}}%
    \fi%
  \else%
    \setlength{\unitlength}{\svgwidth}%
  \fi%
  \global\let\svgwidth\undefined%
  \global\let\svgscale\undefined%
  \makeatother%
  \begin{picture}(1,1.4489507)%
    \lineheight{1}%
    \setlength\tabcolsep{0pt}%
    \put(0,0){\includegraphics[width=\unitlength]{image_of_square_flop_k=n_svg-tex.eps}}%
    \put(0.76145356,0.20273508){\color[rgb]{0,0,0}\makebox(0,0)[lt]{\lineheight{1.25}\smash{\begin{tabular}[t]{l}\dgr{a}\end{tabular}}}}%
    \put(0.01717356,0.19651644){\color[rgb]{0,0,0}\makebox(0,0)[lt]{\lineheight{1.25}\smash{\begin{tabular}[t]{l}\dgr{a+1}\end{tabular}}}}%
    \put(0.37370513,0.71474429){\color[rgb]{0,0,0}\makebox(0,0)[lt]{\lineheight{1.25}\smash{\begin{tabular}[t]{l}\dgr{a+2}\end{tabular}}}}%
    \put(0.75081674,1.23073522){\color[rgb]{0,0,0}\makebox(0,0)[lt]{\lineheight{1.25}\smash{\begin{tabular}[t]{l}\dgr{a}\end{tabular}}}}%
    \put(0.02709222,1.22451649){\color[rgb]{0,0,0}\makebox(0,0)[lt]{\lineheight{1.25}\smash{\begin{tabular}[t]{l}\dgr{a+1}\end{tabular}}}}%
  \end{picture}%
\endgroup%

%% file: arxiv-figures/ssbim_square_flop_n=k_simplify_1_svg-tex.eps_tex
\begingroup%
  \makeatletter%
  \providecommand\color[2][]{%
    \errmessage{(Inkscape) Color is used for the text in Inkscape, but the package 'color.sty' is not loaded}%
    \renewcommand\color[2][]{}%
  }%
  \providecommand\transparent[1]{%
    \errmessage{(Inkscape) Transparency is used (non-zero) for the text in Inkscape, but the package 'transparent.sty' is not loaded}%
    \renewcommand\transparent[1]{}%
  }%
  \providecommand\rotatebox[2]{#2}%
  \newcommand*\fsize{\dimexpr\f@size pt\relax}%
  \newcommand*\lineheight[1]{\fontsize{\fsize}{#1\fsize}\selectfont}%
  \ifx\svgwidth\undefined%
    \setlength{\unitlength}{583.0181444bp}%
    \ifx\svgscale\undefined%
      \relax%
    \else%
      \setlength{\unitlength}{\unitlength * \real{\svgscale}}%
    \fi%
  \else%
    \setlength{\unitlength}{\svgwidth}%
  \fi%
  \global\let\svgwidth\undefined%
  \global\let\svgscale\undefined%
  \makeatother%
  \begin{picture}(1,1.44937129)%
    \lineheight{1}%
    \setlength\tabcolsep{0pt}%
    \put(0,0){\includegraphics[width=\unitlength]{ssbim_square_flop_n=k_simplify_1_svg-tex.eps}}%
    \put(0.76169849,0.20245659){\color[rgb]{0,0,0}\makebox(0,0)[lt]{\lineheight{1.25}\smash{\begin{tabular}[t]{l}\dgr{0}\end{tabular}}}}%
    \put(0.1080148,0.20660639){\color[rgb]{0,0,0}\makebox(0,0)[lt]{\lineheight{1.25}\smash{\begin{tabular}[t]{l}\dgr{1}\end{tabular}}}}%
    \put(0.75105171,1.23141979){\color[rgb]{0,0,0}\makebox(0,0)[lt]{\lineheight{1.25}\smash{\begin{tabular}[t]{l}\dgr{0}\end{tabular}}}}%
    \put(0.10756849,1.23141979){\color[rgb]{0,0,0}\makebox(0,0)[lt]{\lineheight{1.25}\smash{\begin{tabular}[t]{l}\dgr{1}\end{tabular}}}}%
    \put(0.34782578,0.7252201){\color[rgb]{0,0,0}\makebox(0,0)[lt]{\lineheight{1.25}\smash{\begin{tabular}[t]{l}\blk{\mu^{\hh{0,2},\hh{1,2}}_{\hh{2}}}\end{tabular}}}}%
  \end{picture}%
\endgroup%

%% file: arxiv-figures/ssbim_square_flop_n=k_simplify_2_svg-tex.eps_tex
\begingroup%
  \makeatletter%
  \providecommand\color[2][]{%
    \errmessage{(Inkscape) Color is used for the text in Inkscape, but the package 'color.sty' is not loaded}%
    \renewcommand\color[2][]{}%
  }%
  \providecommand\transparent[1]{%
    \errmessage{(Inkscape) Transparency is used (non-zero) for the text in Inkscape, but the package 'transparent.sty' is not loaded}%
    \renewcommand\transparent[1]{}%
  }%
  \providecommand\rotatebox[2]{#2}%
  \newcommand*\fsize{\dimexpr\f@size pt\relax}%
  \newcommand*\lineheight[1]{\fontsize{\fsize}{#1\fsize}\selectfont}%
  \ifx\svgwidth\undefined%
    \setlength{\unitlength}{649.81886894bp}%
    \ifx\svgscale\undefined%
      \relax%
    \else%
      \setlength{\unitlength}{\unitlength * \real{\svgscale}}%
    \fi%
  \else%
    \setlength{\unitlength}{\svgwidth}%
  \fi%
  \global\let\svgwidth\undefined%
  \global\let\svgscale\undefined%
  \makeatother%
  \begin{picture}(1,1.30037738)%
    \lineheight{1}%
    \setlength\tabcolsep{0pt}%
    \put(0,0){\includegraphics[width=\unitlength]{ssbim_square_flop_n=k_simplify_2_svg-tex.eps}}%
    \put(0.82566073,0.17762733){\color[rgb]{0,0,0}\makebox(0,0)[lt]{\lineheight{1.25}\smash{\begin{tabular}[t]{l}\dgr{0}\end{tabular}}}}%
    \put(0.09397334,0.17576583){\color[rgb]{0,0,0}\makebox(0,0)[lt]{\lineheight{1.25}\smash{\begin{tabular}[t]{l}\dgr{1}\end{tabular}}}}%
    \put(0.82541625,1.1026757){\color[rgb]{0,0,0}\makebox(0,0)[lt]{\lineheight{1.25}\smash{\begin{tabular}[t]{l}\dgr{0}\end{tabular}}}}%
    \put(0.09357291,1.10081416){\color[rgb]{0,0,0}\makebox(0,0)[lt]{\lineheight{1.25}\smash{\begin{tabular}[t]{l}\dgr{1}\end{tabular}}}}%
    \put(0.38942771,0.65037439){\color[rgb]{0,0,0}\makebox(0,0)[lt]{\lineheight{1.25}\smash{\begin{tabular}[t]{l}\blk{\mu^{\hh{0,2},\hh{1,2}}_{\hh{2}}}\end{tabular}}}}%
  \end{picture}%
\endgroup%

%% file: arxiv-figures/ssbim_square_flop_n=k_simplify_3_svg-tex.eps_tex
\begingroup%
  \makeatletter%
  \providecommand\color[2][]{%
    \errmessage{(Inkscape) Color is used for the text in Inkscape, but the package 'color.sty' is not loaded}%
    \renewcommand\color[2][]{}%
  }%
  \providecommand\transparent[1]{%
    \errmessage{(Inkscape) Transparency is used (non-zero) for the text in Inkscape, but the package 'transparent.sty' is not loaded}%
    \renewcommand\transparent[1]{}%
  }%
  \providecommand\rotatebox[2]{#2}%
  \newcommand*\fsize{\dimexpr\f@size pt\relax}%
  \newcommand*\lineheight[1]{\fontsize{\fsize}{#1\fsize}\selectfont}%
  \ifx\svgwidth\undefined%
    \setlength{\unitlength}{739.81886894bp}%
    \ifx\svgscale\undefined%
      \relax%
    \else%
      \setlength{\unitlength}{\unitlength * \real{\svgscale}}%
    \fi%
  \else%
    \setlength{\unitlength}{\svgwidth}%
  \fi%
  \global\let\svgwidth\undefined%
  \global\let\svgscale\undefined%
  \makeatother%
  \begin{picture}(1,1.14218466)%
    \lineheight{1}%
    \setlength\tabcolsep{0pt}%
    \put(0,0){\includegraphics[width=\unitlength]{ssbim_square_flop_n=k_simplify_3_svg-tex.eps}}%
    \put(0.84850446,0.15562278){\color[rgb]{0,0,0}\makebox(0,0)[lt]{\lineheight{1.25}\smash{\begin{tabular}[t]{l}\dgr{0}\end{tabular}}}}%
    \put(0.08090624,0.16281724){\color[rgb]{0,0,0}\makebox(0,0)[lt]{\lineheight{1.25}\smash{\begin{tabular}[t]{l}\dgr{1}\end{tabular}}}}%
    \put(0.84828972,0.9681377){\color[rgb]{0,0,0}\makebox(0,0)[lt]{\lineheight{1.25}\smash{\begin{tabular}[t]{l}\dgr{0}\end{tabular}}}}%
    \put(0.08382472,0.96879173){\color[rgb]{0,0,0}\makebox(0,0)[lt]{\lineheight{1.25}\smash{\begin{tabular}[t]{l}\dgr{1}\end{tabular}}}}%
    \put(0.27697632,0.57478369){\color[rgb]{0,0,0}\makebox(0,0)[lt]{\lineheight{1.25}\smash{\begin{tabular}[t]{l}\blk{\dd^{\hh{0,1,2}}_{\hh{0,1}}\left(\mu^{\hh{0,2},\hh{1,2}}_{\hh{2}}\right)}\end{tabular}}}}%
  \end{picture}%
\endgroup%

%% file: Supplementary_Material.tex
\subsection{From extended affine to affine.}\label{subsec extended affine to affine}
Most, if not all, of the content in this subsection is well-known to experts, but this subsection is warranted by the lack of a consolidated reference in the literature (as far as the authors are aware).
We retain the notation as in \cref{subsec extended affine Weyl groups}:
$W$ denotes the affine Weyl group, $W_{\ext}=\Omega \ltimes W $ denotes the extended affine Weyl group of $\PGL_n$, $\omega$ is a generator of $\Om \subset W_{\ext}$, with $\omega s_i \omega^{-1}=s_{i+1}$, and $\omega^n=1$, $\varpi_k$ denotes the fundamental weight for $1\leq k \leq n-1$, $\rot: \Lambda_{wt} \rightarrow \Om$ denotes the quotient map with respect to the root lattice.
While $W_{\ext}$ is not a Coxeter group, the length function on $W$ extends to a length function on $W_{\ext}$ such that $\ell(\omega^k)=0$ for any $k$.


Since $W_{\ext}= W_{\hh{0}}\ltimes\Lambda_{\wt}$, double cosets in
$W_{\hh{0}}\backslash W_{\ext}/W_{\hh{0}}$ have a unique representative in $\Lambda_{\wt}/W_{\hh{0}}$ which is naturally in bijection with $\Lambda_{\wt}^+$ . 
We refer to elements of $W_{\hh{0}}\backslash W_{\ext}/W_{\hh{0}}$  as \emph{spherical double cosets}, and the double coset $W_{\hh{0}}t_{\lambda}W_{\hh{0}}$ as the \emph{spherical double coset of $\lambda$}.
Note that when $\lambda \in \Lambda_{\rt}^+$, the spherical double  coset of $\lambda$ is an element of $W_{\hh{0}}\backslash W/W_{\hh{0}} \subset W_{\hh{0}}\backslash W_{\ext}/W_{\hh{0}}$.

\begin{defn}\label{defn the map Xi}
    We define the map 
    \begin{equation}\Xi: \Lambda_{\wt}^+ \rightarrow \coprod_{0 \leq k\leq n-1} W_{\hh{0}}\backslash W/W_{\hh{k}}\end{equation}
    by setting 
    \begin{equation}\Xi(\lambda):=W_{\hh{0}}t_{\lambda}\omega^{-\rot(\lambda)}W_{\hh{\rot(\lambda)}}
    \end{equation}   
    so that 
    \begin{equation}
        W_{\hh{0}}t_{\lambda}W_{\hh{0}}= \Xi(\lambda)\omega^{\rot(\lambda)}.
    \end{equation}
\end{defn}

\begin{lemma}\label{lemma Xi vs psi}
The map $\Xi$ is a bijective correspondence
which relates to the correspondence $\psi_0$ in \cref{def:psia} by the equation
$$\inv\circ\psi_0(\lambda)=\Xi(\lambda),$$
where $\inv: \displaystyle \coprod_{k, l \in \Om} W_{\hh{l}}\backslash W/W_{\hh{k}} \xrightarrow{\sim} \coprod_{l,m \in \Om}W_{\hh{k}}\backslash W/W_{\hh{l}}$ is the isomorphism induced by the map $w \mapsto w^{-1}$.
\end{lemma}
\begin{proof}
    This is an immediate consequence of \cref{lemma coset correspondence for extended affine spherical cosets}.
\end{proof}

\begin{remark}\label{rmk coset mismatch}
    For $\lambda \in \Lambda_{\wt}^+$ with $\rot(\lambda)=k$, the spherical double cosets $\omega^{-k}\psi_0(\lambda)$ and $\Xi(\lambda)\omega^k$ coincide with $W_{\hh{0}}t_{-\lambda}W_{\hh{0}}$ and $W_{\hh{0}}t_{\lambda}W_{\hh{0}}$ respectively. 
    Since the dominant weight representative in the $W$-orbit of $-\lambda$ is $-w_{\hh{0}}(\lambda)$,  the dominant representative of the spherical double coset $\omega^{-k}\psi_0(\lambda)$ is $-w_0(\lambda)$.
\end{remark}

Let $\HH_{\ext}$ denote the extended affine Hecke algebra associated to $W_{\ext}$, which has the affine Hecke algebra $\HH_{\aff}$ as a subalgebra. 
To make the formulation of the Satake isomorphism in terms of the Hecke algebroid precise, we compare the presentations of the Hecke algebra in \cite{WillSingular}  and \cite{Ach} (which offers a consolidated alternative to \cite{LusztigGS}).

The presentation we follow in \cref{subsec-Hecke} is from \cite{WillSingular}, so we start by recalling the explicit conventions from the same. 
The Hecke algebra in \cite{WillSingular} is defined using the quadratic relation
\begin{equation}H_s^2=(v^{-1}-v)H_s + 1.
\end{equation}
For a double coset $p \in W_I\backslash W/W_J$ with maximal element $\bar{p}$, the corresponding explicit standard basis element and Kazhdan-Lusztig basis element are given by
\begin{equation}\label{eq Williamson std basis and KL basis Hecke algebroid}
 \HAB{I}{J}{p}:=  \sum_{x\in p}v^{\ell(\bar{p})-\ell(x)}H_x, \qquad
 \KLB{I}{J}{p}:=\KLB{}{}{\bar{p}}.
\end{equation}
Note that the standard basis here is normalized precisely so that the transformation matrix from the standard basis to the Kazhdan-Lusztig basis (of $\Hecke{I}{J}$) is unipotent.

There is an action of $\Omega$ on $\HH_{\aff}$, corresponding to the action of $\Omega$ on the affine Dykin diagram, and $\HH_{\ext}$ may be thought of as the twisted group ring
\begin{equation}\label{eq H-ext as a semidirect product of H-aff and omega}
\HH_{\ext}=\HH_{\aff} \rtimes \Omega.
\end{equation}
Indeed, since $\ell(\omega^k)=0$, the elements $H_{\omega^k}\in \HH_{\ext}$ satisfies the equations 
\begin{equation}\label{eq Hecke multiplication length zero elements}
  H_x H_{\omega^k}=H_{x\omega^k}, \quad H_{\omega^k} H_x=H_{\omega^k x}  , \quad H_{\omega^k}H_x=H_{\omega^k(x)}H_{\omega^k}
\end{equation} for any $x\in W$ where $\omega^k(x)=\omega^k x \omega^{-k}$.

Now we recall the presentation for the extended affine Hecke algebra and explicit formulas for the standard basis for the corresponding spherical Hecke algebra from \cite{Ach}.
This presentation for the extended affine Hecke algebra, as in \cite{Ach}, is given over $\Z[q^{\pm1/2}]$, with basis elements $T_w$ for $w \in W_{\ext}$, and the quadratic relation 
\begin{equation}\label{eq Hecke algebra presentation Williamson vs Lusztig}
    T_s^2=(q-1)T_s+q T_1.
\end{equation}
Comparing with the presentation  from \cite{WillSingular}, we have
\begin{equation}v = q^{-1/2}, \qquad H_w=v^{\ell(w)}T_w.
\end{equation}

The spherical Hecke algebra $\HH_{\sph}$ is then defined as a subalgebra of the extended affine Hecke algebra localized at the Coxeter polynomial 

\begin{equation}\label{eq coxeter polynomial achar to williamson}
W_q:=\sum_{w \in W_{\hh{0}}}q^{\ell(w)}=\sum_{w \in W_{\hh{0}}}v^{-\ell(w)}=\tilde{\pi}(\hh{0})=v^{-\ell(w_{\hh{0}})}\pi(\hh{0})
\end{equation}
(for $\tilde{\pi}, \pi$ is as in \cref{subsec-Hecke}).
Explicitly, this subalgebra is the $\Z[v^{\pm 1}]$-span of $\set{K_\lambda}_{\lambda \in \lambda_{\wt}^+}$, where
\begin{equation}\label{eq standard basis spherical Hecke}
    K_\lambda:=W_q^{-1} \sum_{w\in W_{\hh{0}}t_\lambda W_{\hh{0}}}T_w.
\end{equation}
This subalgebra also coincides with $K_0\HH_{\ext}\cap \HH_{\ext}K_0=K_0\HH_{\ext}K_0$ (\cite[Lemma 9.3.6]{Ach}).
Note here that $K_0$ is an idempotent, and is explicitly related to the quasi-idempotent $\un{H}_{w_{\hh{0}}}$ as
\begin{equation}
    K_0=q^{\ell(w_{\hh{0}})/2}W_{q}^{-1}\un{H}_{w_{\hh{0}}}=\frac{1}{\pi(\hh{0})}\un{H}_{w_{\hh{0}}}.
\end{equation}
The spherical Hecke algebra $\HH_{\sph}$ inherits the bar involution from $\HH_{\ext}$, and there is an invariant basis $\set{\un{S}_{\lambda}}_{\lambda \in \Lambda^+_{\wt}}$ unique up to certain properties, as in \cite[Theorem 9.3.7]{Ach}.

For $\lambda \in \Lambda_{\wt}^+$, let $L_{\lambda}$ denote the irreducible representation of highest weight $\lambda$ in $\Rep_{\mathbb{C}}(\mathfrak{sl}_n)$.
Lusztig, in \cite{LusztigGS}, showed that the Satake isomorphism precisely identifies $\ch(L_\lambda)$ with $\un{S}_\lambda$. We formally state these results below (see \cite[Theorem 9.4.1, 9.4.2]{Ach}).

\begin{theorem}[Satake]\label{thm Satake original}
    There is an isomorphism of $\Z[v^{\pm 1}]$-algebras 
    $$\gamma: \HH_{\sph} \xrightarrow{\sim} \Z[v^{\pm 1}][\Lambda_{\wt}]^W.$$
\end{theorem}

\begin{theorem}[Lusztig]\label{thm Lusztig Satake}
    The Satake isomorphism $\gamma$ identifies $\un{S}_\lambda$ with $\ch(L_\lambda)$.
\end{theorem}

Inside $\HH_{\sph}$, we have the subalgebra $\HH_{\sph, \aff}$ spanned by the $K_\lambda$ for $\lambda \in \Lambda_{\rt}^+$.
The decomposition in \cref{eq H-ext as a semidirect product of H-aff and omega} then induces a decomposition
\begin{equation}
    \HH_{\sph}=\HH_{\sph, \aff} \rtimes \Om.
\end{equation}

The following result is well known (cf. \cite[Remark 2.2.2 (3)]{WillThesis}), but we provide a proof for the sake of clarity.

\begin{lemma}\label{lemma affine spherical Hecke in Hecke algebroid}
 There is an isomorphism of $\Z[v^{\pm 1}]$-algebras   
 $$\HH_{\sph, \aff} \cong \Hecke{\hh{0}}{\hh{0}}$$
 which sends the basis elements $K_\lambda$ to $v^{-\langle2\rho^\vee, \lambda \rangle}\HAB{\hh{0}}{\hh{0}}{W_{\hh{0}}t_{\lambda}W_{\hh{0}}}$ and $\un{S}_\lambda$ to $\KLB{\hh{0}}{\hh{0}}{W_{\hh{0}}t_{\lambda}W_{\hh{0}}}$.
\end{lemma}
\begin{proof}
For any $\lambda \in \Lambda_{\wt}$, the maximal coset representative of $W_{\hh{0}}t_\lambda W_{\hh{0}}$ has length $\ell(w_{\hh{0}}) + \langle 2\rho^\vee, \lambda \rangle$ (see \cite[Chapter 2]{Macdonaldhecke}\footnote{This fact may also be seen as a consequence of \cref{lemma dist coset reps}, \cref{thm:weightsadd} and \cref{lemma distinguished reps for fundamental weights}.}).
From \cref{eq standard basis spherical Hecke}, \cref{eq coxeter polynomial achar to williamson} and \cref{eq Hecke algebra presentation Williamson vs Lusztig},
we have that for $\lambda \in \Lambda_{\rt}^+$,
$$K_\lambda= \frac{v^{\ell(w_{\hh{0}})}}{\pi(\hh{0})} \sum_{w\in W_{\hh{0}}t_\lambda W_{\hh{0}}}v^{-\ell(w)}H_w= \frac{v^{-\langle2 \rho^\vee, \lambda \rangle}}{\pi(\hh{0})}\; \HAB{\hh{0}}{\hh{0}}{W_{\hh{0}}t_\lambda W_{\hh{0}}}$$
in $\HH_{\aff}[W_q^{-1}]=\HH_{\aff}[\pi(\hh{0})^{-1}].$
Additionally, given $\mu \in \Lambda_{\rt}$, we have
$$K_\lambda K_\mu
=\frac{v^{-\langle 2 \rho^\vee, \lambda+\mu \rangle}}{\pi(\hh{0})^2} \HAB{\hh{0}}{\hh{0}}{W_{\hh{0}}t_\lambda W_{\hh{0}}} \HAB{\hh{0}}{\hh{0}}{W_{\hh{0}}t_\mu W_{\hh{0}}}
=\frac{1}{\pi(\hh{0})}(v^{-\langle 2\rho^\vee, \lambda \rangle} \HAB{\hh{0}}{\hh{0}}{\lambda})
\ast_{\hh{0}}(v^{-\langle 2\rho^\vee, \mu \rangle} \HAB{\hh{0}}{\hh{0}}{\mu}).$$
It follows that the $\Z[v^{\pm 1}]$-linear automorphism of $\HH_{\aff}[\pi(\hh{0})^{-1}]$ defined by multiplication by $\pi(\hh{0})$ identifies  $K_\lambda$ with $v^{-\langle 2\rho^\vee, \lambda \rangle} \HAB{\hh{0}}{\hh{0}}{\lambda}$, and induces an isomorphism of algebras $\HH_{\sph, \aff} \cong \Hecke{\hh{0}}{\hh{0}}$.
This isomorphism preserves the bar involution, and the statement for the invariant basis may be deduced using the defining properties, or by using \cite[Theorem 6.12]{LusztigGS}, which states that
$\un{S}_\lambda$ (denoted by $C'_\lambda$ in \cite{LusztigGS}) equals 
$$\frac{1}{\pi(\hh{0})}\KLB{\hh{0}}{\hh{0}}{W_{\hh{0}}t_\lambda W_{\hh{0}}}$$
in $\HH_{\ext}[\pi(\hh{0})^{-1}]$.
\end{proof}

Now consider an algebra $\HHHH$ made using certain $\Hom$-spaces in the Hecke algebroid $\HH$ as defined below. We follow the notation from \cref{subsec-Hecke}.
\begin{defn}
  Let  $\HHHH:= \bigoplus_{0 \leq k \leq n-1}\Hecke{\hh{0}}{\hh{k}}$. 
  Define an $\Om-$graded $\Z[v,v^{-1}]$-algebra structure 
  $\bullet:\HHHH \times \HHHH \rightarrow \HHHH$ such that
  for $f \in \Hecke{\hh{0}}{\hh{k}}$, $g \in \Hecke{\hh{0}}{\hh{l}}$,
\begin{equation} \label{eq multiplication in HHHH}
    f\bullet g:= f\ast_{\hh{k}}\omega^k(g) \in \Hecke{\hh{0}}{\hh{k+l}},
\end{equation}
extended linearly to $\HHHH$. 
\end{defn}

\begin{lemma}\label{lemma HHHH isomorphic to spherical Hecke algebra}
  There is an isomorphism of $\Om$-graded $\Z[v^{\pm 1}]$-algebras
  $\HH_{\sph} \cong \HHHH$
  which sends the basis elements $K_\lambda$ to $v^{-\langle2 \rho^\vee, \lambda \rangle} \HAB{\hh{0}}{\hh{k}}{\Xi(\lambda)}$ and $\un{S}_\lambda$ to $\KLB{\hh{0}}{\hh{k}}{\Xi(\lambda)},$ for $\lambda \in \Lambda_{\wt}^+$ with $\rot(\lambda)=k$. In particular, $\HHHH$ is a commutative algebra.
\end{lemma}
\begin{proof}
    We proceed similarly as in the proof of \cref{lemma affine spherical Hecke in Hecke algebroid}, and use the equations in \cref{eq Hecke multiplication length zero elements} to see that  
    given $\lambda \in \Lambda_{\wt}^+$ with $\rot(\lambda)=k$, we have
    $$K_\lambda
    = \frac{v^{\ell(w_{\hh{0}})}}{\pi(\hh{0})} \sum_{w\in W_{\hh{0}}t_\lambda W_{\hh{0}}}v^{-\ell(w)}H_w 
    = \frac{v^{\ell(w_{\hh{0}})}}{\pi(\hh{0})} \sum_{w\in \Xi(\lambda)}v^{-\ell(w)}H_w H_{\omega^k}
    = \frac{v^{-\langle2 \rho^\vee, \lambda \rangle}}{\pi(\hh{0})} \HAB{\hh{0}}{\hh{k}}{\Xi(\lambda)} H_{\omega^k}$$
    in $\HH_{\ext}[\pi(\hh{0})^{-1}]$.
    Additionally, given $\mu \in \Lambda_{\wt}^+$ with $\rot(\mu)=l$, we have
    \begin{align*}
    K_\lambda K_\mu
    =\frac{v^{-\langle 2 \rho^\vee, \lambda+\mu \rangle}}{\pi(\hh{0})^2} \HAB{\hh{0}}{\hh{k}}{\Xi(\lambda)}H_{\omega^k} &\HAB{\hh{0}}{\hh{l}}{\Xi(\mu)}H_{\omega^l}
    =\frac{v^{-\langle 2 \rho^\vee, \lambda+\mu \rangle}}{\pi(\hh{0})\pi(\hh{k})} \HAB{\hh{0}}{\hh{k}}{\Xi(\lambda)}\HAB{\hh{k}}{\hh{k+l}}{\omega^k({\Xi(\mu)})}H_{\omega^{k+l}}\\
    &=\frac{1}{\pi(\hh{0})}(v^{-\langle 2\rho^\vee, \lambda \rangle} \HAB{\hh{0}}{\hh{k}}{\Xi(\lambda)})
    \ast_{\hh{k}}\omega^k \left(v^{-\langle 2\rho^\vee, \mu \rangle} \HAB{\hh{0}}{\hh{l}}{\Xi(\mu)}\right)H_{\omega^{k+l}}\\
    &=\frac{1}{\pi(\hh{0})}\left((v^{-\langle 2\rho^\vee, \lambda \rangle} \HAB{\hh{0}}{\hh{k}}{\Xi(\lambda)})
    \bullet (v^{-\langle 2\rho^\vee, \mu \rangle} \HAB{\hh{0}}{\hh{l}}{\Xi(\mu)})\right)H_{\omega^{k+l}}.
\end{align*}
It follows that the $\Z[v^{\pm 1}]$-linear isomorphism sending $K_\lambda$ to $v^{-\langle2 \rho^\vee, \lambda \rangle} \HAB{\hh{0}}{\hh{k}}{\Xi(\lambda)}$ is an isomorphism of $\Om$-graded algebras.
The statement for $\un{S}_{\lambda}$ then follows from \cite[Theorem 6.12]{LusztigGS} which implies that for $\lambda \in \Lambda_{\wt}^+$,
$$\un{S}_{\lambda}= \frac{1}{\pi(\hh{0})}\KLB{\hh{0}}{\hh{k}}{\Xi(\lambda)}H_{\omega^k}$$
in $\HH_{\ext}[\pi(\hh{0})^{-1}]$.
The statement regarding commutativity follows from the commutativity of $\HH_{\sph}.$
\end{proof}

To match with the conventions in this paper (in particular, the conventions in the coset correspondence in \cref{def:psia} and reading diagrams from right to left), we use the anti-involution of the Hecke algebroid, induced by that of the Hecke algebra which sends $H_w$ to $H_{w^{-1}}$.
Note that this lifts the involution of $W$ given by $w \mapsto w^{-1}$,  and identifies $\HAB{\hh{k}}{\hh{l}}{p}$ with $\HAB{\hh{l}}{\hh{k}}{\iota(p)}$ (where iota is as in \cref{lemma Xi vs psi}).
This anti-involution further identifies the commutative algebra $\HHHH^{\op}=\HHHH$ with $\HHH$ (see \cref{defn the algebra HHH}), identifying $\HAB{\hh{0}}{\hh{k}}{\Xi(\lambda)}$ with $\HAB{\hh{k}}{\hh{0}}{\Psi_0(\lambda)}$.
Combining \cref{thm Satake original}, \cref{thm Lusztig Satake} and \cref{lemma HHHH isomorphic to spherical Hecke algebra} with this observation, we then have the following theorem.
\begin{theorem}\label{thm Satake isomorphism to HHH}
 There is an isomorphism of $\Om$-graded $\Z[v^{\pm 1}]$-algebras
 \begin{gather}
     \notag \Z[v^{\pm 1}][\Lambda_{\wt}]^W \xrightarrow{\sim} \HH_{\sph} \xrightarrow{\sim} \HHHH \xrightarrow{\sim}\HHH \\
     \notag  \ch(L_\lambda) \mapsto \un{S}_\lambda \mapsto \KLB{\hh{0}}{\hh{\rot(\lambda)}}{\Xi(\lambda)} \mapsto \KLB{\hh{\rot(\lambda)}}{\hh{0}}{\Psi_0(\lambda)}. 
 \end{gather}   
\end{theorem}


\PTv2{Leaving the following comment for v2 modifications...}
\PTv2{I think I now understand the reason for the mismatch. The Soergel version of geometric Satake we construct is actually a composition of usual geometric Satake followed by the inversion anti-automorphism on the loop group/affine Kac Moody group. Normal geometric Satake at the affine Kac moody group level would take $$L_\lambda \mapsto IC_\lambda \mapsto IC_{W_{\hat{0}}1W_{\hat{k}}} \in \Perv_{P_{\hat{0}}\times P_{\hat{\hat{k}}}}(G_{KM}).$$
Composing with the isomorphism $$\Perv_{P_{\hat{0}}\times P_{\hat{\hat{k}}}}(G_{KM}) \xrightarrow{\sim} \Perv_{P_{\hat{k}}\times P_{\hat{0}}}(G_{KM})$$
induced by the inverse map of $G_{KM}$, we get our version. 
At the bimodule level, this is likely swapping the left and right actions of the rings, which I think decategorifies to the anti-involution obtained by composing the bar involution with the KL anti-involution $\omega$.
Also, $\ch(IC_{\omega_1)}$ canonically lives (via pulling back to $G_{KM}$, taking hypercohomology, and then decategoirfying) in $\Hecke{\hh{-1}}{\hh{0}}$ or $\Hecke{\hh{0}}{\hh{1}}$ which via the above anti-involution gets identified with $\Hecke{\hh{1}}{\hh{0}}$}

\subsection{Comparing the \texorpdfstring{$\zeta$}{zeta}-deformed and \texorpdfstring{$q$}{q}-deformed realizations} \label{appendix:comparisonqz}
In this section, we explicitly compare the deformed affine (Frobenius) realization \cref{defn deformed affine frobenius realization} with the following $q-$deformed version.
We assume $n \geq 3$, though there are similar statements in the $n=2$ case (see \cite{EJY1}).
Throughout this section, $\mc{A}_q=\kk[q^{\pm1}]$, $\mc{A}_\zeta=\kk[\zeta^{\pm1}]$ and $\mc{A}_{\zeta,q}=A_{\Z}\otimes_{\Z}\kk$, where $\kk$ is a commutative domain and $A$ is as in \cref{eq coefficient ring with q and zeta}.
Define $V_q$ as follows, over $\mc{A}_q=\kk[q^{\pm1}]$.
\begin{defn}\cite[Definition 2.4]{EJY1}\label{defn q-deformed affine realization}
    Let $V_q$ be the free module over $\mc{A}_q$ with basis $\set{v_i}_{1\leq i \leq n}$ and $\set{v_i^*}$ the corresponding dual basis.
    Let $$\alpha_i=v_i-v_{i+1},\;\; \alpha_{i}^\vee=v_i^*-v_{i+1}^*$$ for $i \neq 0$, and $$\al_0=qv_n-q^{-1}v_1, \;\;\al_0^\vee=q^{-1}v_n^*-qv_1^*.$$
    Then $W$ acts on $V_q$ as follows:
    for the subgroup $W_{\hh{0}}=S_n$ generated by $\set{s_i: i\neq 0}$, we get the usual permutation action, permuting the basis elements $\set{v_k}$.
    The action of $s_0$ is explicitly given by:
    $$s_0(v_n)=q^{-2}v_1,\;\; s_0(v_1)=q^{2} v_n.$$
\end{defn}

Let $V_\zeta$ denote the deformed affine realization over $\mc{A}_\zeta$, as in \cref{defn-deformed affine realization}. 
We then have:
\begin{lemma} \cite[Lemma 2.12]{EJY1}\label{lemma q-deformed to zeta-deformed affine realization}
    There is an isomorphism of $\mc{A}_{\zeta,q}$-modules
    $$V_q\otimes_{\mc{A}_q}\mc{A}_{\zeta,q} \xrightarrow{\sim} V_{\zeta}\otimes_{\mc{A}_\zeta}\mc{A}_{\zeta,q}: v_k \mapsto \zeta^k x_k, \;\text{$1\leq k \leq n$}.$$
    This isomorphism sends the simple roots $\alpha_k, \; 1\leq k \leq n-1$ in $V_q$ to $\zeta^k\alpha_k$ in the target, and sends $\alpha_0$ in $V_q$ to $q^{-1}\alpha_0$.
    The coroots are scaled by the inverse scalars.
\end{lemma}

Let us explicitly describe the Frobenius structure on $V_q$, which generalizes the Frobenius structure used for the $q$-deformed realization for quantum Satake functor in type $A_2$ in \cite{EQuantumI}.
For $s \in S$, we denote by $\ddq_s$ the corresponding Demazure operator in the $q$-deformed affine Frobenius realization, and for an expression $\un{w}$, we denote by $\ddq_{\un{w}}$ the corresponding Demazure operator.

\begin{defn}\label{defn q-deformed affine Frobenius realization}
The \emph{$q$-deformed affine Frobenius realization} is defined as follows.
\begin{itemize}
\item As in \cref{defn deformed affine frobenius realization}, we fix a system of positive roots $(\Phi_q)_I^+$ (or simply $\Phi_I^+$ when there's no ambiguity in the realization being referred to) for each finitary $I \subset S$.
We only define $\Phi_I^+$ for $I$ connected, since we may take a disjoint over the connected components for a more general $I$.

\item For $I \subset S \setminus \set{0}$, we have the standard set of positive roots from the permutation realization of $S_n=W_{\hh{0}}$.
In particular, for connected $I \subset S \setminus\set{0}$, say for $I=\set{i,i+1, \ldots, k}$, we have the positive roots to be $v_j-v_l$ for $i\leq j < l\leq k+1$.

\item Now consider a connected subset $I$ containing $n$. Then $I=\set{i, i+1, \ldots, n, 1,\ldots,k}.$ 
We then define the positive roots to be 
$$\Phi_I^+=\set{v_j-v_l}_{i\leq j<l\leq n}\sqcup \set{v_j-v_i}_{0<j<l\leq k+1}\sqcup \set{qv_j-q^{-1}v_l}_{i \leq j\leq n,\; 0<l\leq k+1}.$$
\item The $\procopdtq_I^J$, $\ddq_I^J$ for each $J \subset I$ are then defined exactly as in \cref{defn deformed affine frobenius realization}.
\end{itemize}
\end{defn}

Note that via the identification in \cref{lemma q-deformed to zeta-deformed affine realization}, the $\procopdtq_I^J$, $\ddq_I^J$ in $V_q$ are scalar multiples of their respective counterparts in $V_\zeta$.
One may show that this gives a Frobenius hypercube as desired, either using the fact above (using the results for $V_\zeta$ and comparing the Frobenius trace maps), or directly using the techniques that we used to show the statement for $V_\zeta$.
We formalize this in the following lemma.
\begin{lemma}
    Let $R=\Sym_{\mc{A}_q}V_q$ For finitary $I\subset J \subset S$, the extension $R^I \hookrightarrow R^J$ in \cref{defn q-deformed affine Frobenius realization} is a Frobenius extension, and $\ddq_I^J$ is a Frobenius trace, with product-coproduct element $\procopdtq_I^J$.
\end{lemma}

Here is the analogue of \cref{lemma demazure operators computation deformed affine}.

\begin{lemma}\label{lemma q-demazure operators computation}
    For $i \in \Om$, let $\ddq_i$ denote the corresponding Demazure operator on $R=\Sym_{\mc{A}_q}V_q$.     
We then have:
\begin{enumerate}
    \item For $i \in \Om$, $\ddq_i(v_j)=\begin{cases}
        1 &\text{ if $j=i<n$}\\
        -1 &\text{ if $j-1=i<n$}\\
        q^{-1} &\text{ if $j=i=n$}\\
        -q &\text{ if $j-1=i=0$}\\
        0  &\text{ otherwise}
    \end{cases}$
    \item For $i,j \in \Om$ such that $s_is_j=s_js_i$, $\ddq_i \ddq_j= \ddq_j \ddq_i$.
    \item For $i \in \Om$, 
    $$\ddq_{i}\ddq_{i+1}\ddq_{i}=
    \begin{cases}
        \ddq_{i+1}\ddq_i\ddq_{i+1} &\text{ if $1\leq i \leq n-2$}\\
        q\ddq_{i+1}\ddq_i\ddq_{i+1} &\text{ if $i=n-1, n$}
    \end{cases}$$
    \end{enumerate}
\end{lemma}



We now compare the Frobenius trace maps for the $q$-deformed realization with those of the $\zeta$-deformed realization.
\begin{lemma}\label{lemma comparing Frobenius structures q vs zeta}
As before, let $\dd_s$ denote the Demazure operators  for the $V_\zeta$, and $\ddq_s$ the operators for $V_q$. After base change to $\mc{A}_{\zeta,q}$, we have:
    \begin{enumerate}
        \item $\ddq_i=\zeta^{-k}\dd_i$ for $1 \leq i \leq n-1$.
        \item $\ddq_n=q\dd_n$.
        \item For $I=\set{i,i+1, \ldots,k} \subset S \setminus \set{0},$ we 
        have $\procopdtq_I= \kappa_I \mu_I$ and $\ddq_I=\kappa_I^{-1}\dd_I$, where
        $$\kappa_I= \displaystyle \prod_{i\leq j \leq k}\zeta^{j(k+1-j)}$$
        \item For $I = \set{i, i+1, \ldots, n, 1, \ldots, k}$, we have
        $\procopdtq_I= \kappa_I \mu_I$ and $\ddq_I=\kappa_I^{-1}\dd_I$, where
        \begin{align}\label{eqn 1 lemma comparing Frobenius structures q vs zeta}
        \notag \kappa_I &= q^{(n+1-i)(k+1)}\left(\prod_{i\leq j \leq n}\zeta^{j(n+k+1-j)}\right)\left( \prod_{1\leq j\leq k}\zeta^{j(k+1-j)}\right)\\
        &=q^{(n+1+k-i)(k+1)}\left(\prod_{i \leq j \leq n+k} \zeta^{j(n+k+1-j)}\right).
        \end{align}
        \item For an arbitrary finitary subset $I \subset S$, let $I_1, \ldots , I_k$ be the connected components of $I$.
        Then $\procopdtq_I=\kappa_I\mu_I$ and $\ddq_I=\kappa_I^{-1}\mu_I$, where $$\kappa_I=\kappa_{I_1}\kappa_{I_2}\ldots \kappa_{I_k}.$$
        \item For $J\subset I \subset S$ finitary subsets, we have $\procopdtq_I^J= \kappa_I^J \mu_I^J$ and $\ddq_I^J=(\kappa_I^J)^{-1}\dd_I^J$ where $\kappa_I^J=\kappa_J^{-1}\kappa_I.$
    \end{enumerate}
\end{lemma}
\begin{proof}
These are all straightforward computations.
\end{proof}

Let $\qdiag$ denote the diagrammatic category associated to the Frobenius hypercube from \cref{defn q-deformed affine Frobenius realization}, and $\zdiag$ for the one associated to \cref{defn deformed affine frobenius realization} (see \cref{subsec-diagrammatic singular bott-samelson bimodules}).
Let $\Phi_q: \qdiag \rightarrow \SBSBim_q$, $\Phi_\zeta: \zdiag \rightarrow \SBSBim_\zeta$ be the corresponding functors to the algebraic categories associated to $V_q$ and $V_\zeta$ respectively. Over $\mc{A}_{q,\zeta},$ we have an equivalence of 2-categories $\Theta^\zeta_q:\SSBim_q\xrightarrow{\sim}\SSBim_\zeta$ using \cref{lemma q-deformed to zeta-deformed affine realization}.
Working over $\mc{A}_{q,\zeta}$, we then have

\begin{equation}
\begin{tikzcd}
    \qdiag \arrow[r, hook, "\Phi_q"] \arrow[d, "\sim" sloped, "\vartheta_q^\zeta"'] & \SSBim_q \arrow[d, "\sim" sloped , "\Theta^\zeta_q"'] \\
    \zdiag \arrow[r, hook, "\Phi_\zeta"] & \SSBim_\zeta
\end{tikzcd}
\end{equation}
where $\vartheta$ is explicitly described below on the generators of $\qdiag$. This will be essential to produce a $q$-form $\GS_q$ for the $\GS_\zeta$ that we have constructed.

Recall from \cref{subsec-diagrammatics for frob hypercube} that the clockwise cap and counter-clockwise cup are mapped to the multiplication maps and inclusion maps for the respective rings, and is independent of the particular choice of the Frobenius structure we use.
Hence $\vartheta$ sends the clockwise caps (resp. counter-clockwise cups) in $\qdiag$ to those in $\zdiag$.

The counter-clockwise caps and clockwise cups are sent to the corresponding Frobenius trace maps and the co-product maps respectively in the algebraic category. 
The images of these diagrams in $\SSBim_q \simeq \SSBim_\zeta$ then differ by scalars specified by the $\kappa_I$ as in \cref{lemma comparing Frobenius structures q vs zeta}.
More precisely, we have:

\begin{equation}\label{diag- comparing q cap cc with zeta cap cc}
\vartheta^\zeta_q 
\left( \;
{\begin{tikzpicture}[baseline=(current bounding box.center),
    scale=1,
    line width=0.8pt,
    arrow_style/.style={
        decoration={
            markings,
            mark=at position 0.3 with {\arrow{stealth[black]}}, 
            mark=at position 0.8 with {\arrow{stealth[black]}}  
        },
        postaction={decorate}
    }
]
\draw[color={rgb,255:red,0; green,0; blue,192}, arrow_style] (0.8,0) arc (0:180:0.8);
\node at (0,1.5) {$\gr{I}$};
\node at (0,0.3) {$\gr{I\setminus \set{s}}$};
\node at (1,1) {}; 
\node at (-1,1) {}; 
\node at (0,1.9) {}; 

\node[draw, inner sep=0pt, fit=(current bounding box)] {};
\end{tikzpicture}_{\qdiag}}
\right)
= (\kappa_I^{I \setminus \set{s}})^{-1}
\;
{\begin{tikzpicture}[baseline=(current bounding box.center),
    scale=1,
    line width=0.8pt,
    arrow_style/.style={
        decoration={
            markings,
            mark=at position 0.3 with {\arrow{stealth[black]}}, 
            mark=at position 0.8 with {\arrow{stealth[black]}}  
        },
        postaction={decorate}
    }
]
\draw[color={rgb,255:red,0; green,0; blue,192}, arrow_style] (0.8,0) arc (0:180:0.8);
\node at (0,1.5) {$\gr{I}$};
\node at (0,0.3) {$\gr{I\setminus \set{s}}$};
\node at (1,1) {}; 
\node at (-1,1) {}; 
\node at (0,1.9) {}; 

\node[draw, inner sep=0pt, fit=(current bounding box)] {};
\end{tikzpicture}_{\zdiag}}
\end{equation}

\begin{equation}\label{diag- comparing q cup c with zeta cup c}
\vartheta^\zeta_q 
\left( \;
{\begin{tikzpicture}[baseline=(current bounding box.center),
    scale=1,
    line width=0.8pt,
    arrow_style/.style={
        decoration={
            markings,
            mark=at position 0.3 with {\arrow{stealth[black]}}, 
            mark=at position 0.8 with {\arrow{stealth[black]}}  
        },
        postaction={decorate}
    }
]
\draw[color={rgb,255:red,0; green,0; blue,192}, arrow_style] (0.8,1.9) arc (360:180:0.8);
\node at (0,1.5) {$\gr{I}$};
\node at (0,0.4) {$\gr{I\setminus \set{s}}$};
\node at (1,0) {}; 
\node at (-1,0) {}; 
\node[draw, inner sep=0pt, fit=(current bounding box)] {};
\end{tikzpicture}_{\qdiag}}
\right)
= \kappa_I^{I \setminus \set{s}}
\;
{\begin{tikzpicture}[baseline=(current bounding box.center),
    scale=1,
    line width=0.8pt,
    arrow_style/.style={
        decoration={
            markings,
            mark=at position 0.3 with {\arrow{stealth[black]}}, 
            mark=at position 0.8 with {\arrow{stealth[black]}}  
        },
        postaction={decorate}
    }
]
\draw[color={rgb,255:red,0; green,0; blue,192}, arrow_style] (0.8,1.9) arc (360:180:0.8);
\node at (0,1.5) {$\gr{I}$};
\node at (0,0.4) {$\gr{I\setminus \set{s}}$};
\node at (1,0) {}; 
\node at (-1,0) {}; 
\node[draw, inner sep=0pt, fit=(current bounding box)] {};
\end{tikzpicture}_{\zdiag}} \quad .
\end{equation}
\vspace{5 pt}

The generating upward and downward crossings (see \cref{diag generating crossings}) in both $\qdiag$ and $\zdiag$ go to the exact same 2-morphism in $\SBSBim_q \simeq \SBSBim_\zeta$ so that $\vartheta^\zeta_q$ sends these crossings in $\qdiag$ to those in $\zdiag$ (without any coefficients).

\begin{remark}
From here one could explicitly describe a functor $\GS_q': \cwebs \rightarrow \qdiag$ over $\mc{A}_{q, \zeta}$ when $\kk=\Q$, by composing $\GS_{\zeta}$ with $(\vartheta_q^\zeta)^{-1}$. However, this functor does not deserve to be called $\GS_q$, as it is not defined over $\mc{A}_q$ but only over an extension. 
One needs to make a careful, non-symmetric choice of the scalars $\lambda$ and $\nu$ from \cref{diag-dQGS on 2-morphisms} in order for the result to land in $\qdiag$ (over $\mc{A}_q$). We hope to provide this choice in an updated version of the paper. \PTv2{Add it.}
\end{remark}

\subsection{Symmetries of Diagrammatic and Algebraic Bott-Samelson bimodules.}\label{subsec reversal duality and phi}

This subsection is dedicated to explaining precisely how the diagrammatic symmetries of horizontal and vertical flipping (with suitable arrow reversals) connect to certain standard symmetries of $\SBSBim$ (and hence $\SSBim$). Nothing in this section should surprise the experts.
While we only discuss this in the context of $\SBSBim$ from the deformed affine Frealization, the contents of this section generalizes to any (partial, graded or ungraded) Frobenius hypercube $\Gamma$, with $\SBSBim$ replaced by the algebraic category $C(\Gamma)$ from \cite{EWSFrob}.

We retain notation from \cref{subsec-diagrammatic singular bott-samelson bimodules}.

\begin{defn}\label{defn reverse of singular multistep expression}
  Given a singular multistep expression  
  $$IK_{\bullet}=[[I_0 \subset K_1 \supset I_1 \subset \ldots \subset K_{d} \supset I_d]],$$
  we define the \emph{reverse of the multistep expression} $IK_\bullet$ to be
  $$\Rev(IK_\bullet):=[[I_d \subset K_{d} \supset I_{d-1} \subset \ldots \supset I_1 \subset K_1 \supset I_0 ]].$$
\end{defn}

\begin{defn}
    Let $\Rev$ denote the 2-functor on $\mc{D}^{\pre}$ which fixes regions, flips diagrams horizontally (i.e., reflects diagrams across a vertical axis) and  fixes floating polynomials. This 2-functor also necessarily changes some strand orientations to make the diagram consistent with the region labels. 
    Note that for a 1-morphism in $\mc{D}^{\pre}$ given by a singlestep expression $I_\bullet$, $\Rev(I_\bullet)$ is just the reverse of the singlestep expression $I_\bullet$.
\end{defn}

\begin{defn}
    Let $\Dual$ denote the 2-functor on $\mc{D}^{\pre}$ which fixes regions and flips diagrams vertically (i.e., reflects diagrams across a horizontal axis) and fixes floating polynomials. This 2-functor also necessarily changes some strand orientations to make the diagram consistent with the region labels. Note that for a 1-morphism in $\mc{D}^{\pre}$ given by a singlestep expression $I_\bullet$, $\Dual(I_\bullet)=I_\bullet$.
\end{defn}

There are a few conventions in the literature for the corresponding functors on $\SBSBim$, we follow those in \cite[(3.26), (3.48)]{BBE}.

\begin{defn}\label{def reversal symmetry R algebraic}
    Let $\Rev$ be the 2-functor on $\SBSBim$ defined as follows.
    $\Rev$ fixes objects, and acts on the 1-morphism categories by
    $$\Hom_{\SBSBim}(I,J) \xrightarrow{\sim} \Hom_{\SBSBim}(J,I): \Bimod{I}{J}{M} \mapsto \Bimod{J}{I}{M}(\ell(J)-\ell(I)).$$
    That is, $\Rev$ swaps the module actions on the left and right (using commutativity of the rings $R^I, R^J$) and adds a grading shift; $\Rev$ does not change the underlying functions of $2$-morphisms.
    This functor is a weak 2-functor which is contravariant on 1-morphisms and covariant on 2-morphisms, and as such includes the coherence data given by the 2-morphisms
    \begin{gather}
       \Rev(\Bimod{R^J}{R^K}{N})\otimes_{R^J}\Rev(\Bimod{R^I}{R^J}{M}) 
       \xrightarrow{\sim}\Rev(\Bimod{R^I}{R^J}{M} \otimes_{R^J} \Bimod{R^J}{R^K}{N}): n\otimes m \mapsto m\otimes n\\
       \Rev(\Bimod{R^I}{R^I}{R^I})=\Bimod{R^I}{R^I}{R^I} \xrightarrow[\id]{\sim} \Bimod{R^I}{R^I}{R^I}
    \end{gather}
    for $\Bimod{R^I}{R^J}{M} \in \Hom_{\SBSBim}(J,I)$, $\Bimod{R^J}{R^K}{N} \in \Hom_{\SBSBim}(K,J)$, $I,J\subset S$ finitary.
\end{defn}

\begin{defn}\label{def duality symmetry D algebraic}
    Let $\Dual$ denote the 2-functor on $\SBSBim$ defined as follows. 
    $\Dual$ fixes objects, and acts on the 1-morphism categories by
    $$\Hom(I,J)^{\op}\to \Hom(I,J): M \mapsto \Hom^\bullet_{(-,R^J)}(M, R^J)(2\ell(I)-2\ell(J)),$$ 
     where $\Hom^\bullet_{(-,R^J)}(B, R^J)$ denotes the space of right graded $R^J$-module homomorphisms into $R^J$, with $((r \otimes s)\cdot f)(v)=s(f(rv))$ for $r \in R^I$, $s\in R^J$, $f \in B^\vee$.
    The 2-functor $\Dual$ is a weak 2-functor which is covariant on 1-morphisms and contravariant on 2-morphisms, and as such includes the coherence data given by the 2-morphisms 
    \begin{gather}
    \Dual(M)\otimes_{R^J} \Dual(N) \xrightarrow[\epsilon_{M,N}]{\sim} \Dual(M\otimes_{R^J} N): f\otimes g \mapsto \bigl( m \otimes n \mapsto g(f(m)\cdot n) \bigr),\label{eq coherence data for duality}\\
    \Dual(\Bimod{J}{J}{R^J}) \xrightarrow[\varepsilon_J]{\sim} \Bimod{J}{J}{R^J}: f \mapsto f(1), \label{eq coherence data duality identity 1-morphism}
    \end{gather}
    for $M \in \Hom_{\SBSBim}(J,I), N \in \Hom_{\SBSBim}(K,J)$ for finitary $I,J \subset S$. 
    Note that the map in \cref{eq coherence data for duality} is an isomorphism because $M$, being a 1-morphism in $\SBSBim$, is a finitely generated projective (in fact, free) graded right $R^J$-module.
    We refer to $\Dual$ as the \emph{(right) duality functor} on $\SBSBim$.
\end{defn}

\begin{remark}
    By definition, the 2-functors $\Rev$ and $\Dual$ defined on $\mc{D}^{\pre}$ and $\SBSBim$ preserve objects. As such, natural transformations between 2-functors appearing in this section will only involve the data of 2-morphisms satisfying certain compatibility conditions (we default to the identity 1-morphism at the object level).
\end{remark}

\begin{remark}\label{rmk sbsbim larger for weel definedness of dual and rev}
    Strictly speaking, for $\Rev$ and $\Dual$ to be well-defined, we should allow $\SBSBim$ to be slightly bigger then (but equivalent to) the one in \cref{def sbsbim}: namely, we should allow 1-morphism categories to be closed under isomorphisms in the corresponding category of graded bimodules, so that 1-morphisms are those isomorphic (not just equal) to singular Bott-Samelson bimodules. That this category is closed under $\Rev$ and $\Dual$ then follows from \cref{lemma self duality of Bott-Samelsons} and \cref{lemma reversal on bott samelsons and multistep expressions} below, or by direct checks on the monoidal generators.
\end{remark}

The following lemma is easy and left as an exercise to the reader.
\begin{lemma}\label{lemma dual and rev are involutions}
    The 2-functors $\Dual$ and $\Rev$ on $\mc{D}^{\pre}$ and $\SBSBim$ are involutions.
    More precisely, there are equalities
    $$\Dual \circ \Dual = \Id_{\mc{D}^{\pre}}= \Rev \circ \Rev$$
    of 2-functors on $\mc{D}$, and an equality and an isomorphism
    $$\Id_{\SBSBim}=\Rev\circ \Rev, \qquad \Id_{\SBSBim} \Rightarrow \Dual \circ \Dual$$ of 2-functors on $\SBSBim$, where the latter is given by the evaluation maps $$B \xrightarrow{\sim} \Dual\circ \Dual(B): m \mapsto(f\mapsto f(m))$$ for 1-morphisms $B$ in $\SBSBim$.
\end{lemma}

Recall that 1-morphisms in $\SBSBim$ are\footnote{up to an isomorphism, see \cref{rmk sbsbim larger for weel definedness of dual and rev}.} given by singular Bott-Samelson bimodules
\begin{equation*}
    \BS(IK_{\bullet}):= R^{I_0}\otimes_{R^{K_1}} R^{I_1}\otimes_{R^{K_2}} \cdots \otimes_{R^{K_d}}R^{I_d}(\sum_{i=1}^d \ell(K_i)-\ell(I_i)) \in \Hom_{\SBSBim}(I_d,I_0),
\end{equation*}
for a singular multistep expression
\begin{equation}\label{eq singular multistep expression supplementary}
    IK_{\bullet}=[[I_0 \subset K_1 \supset I_1 \subset \ldots \subset K_{d} \supset I_d]].
\end{equation}

The following lemma is immediate from the definitions.
\begin{lemma}\label{lemma reversal on bott samelsons and multistep expressions}
    For any  multistep expression $IK_\bullet$ as above, there is an isomorphism 
    $$\etadot_{IK_\bullet}: \BS(\Rev(IK_\bullet)) \xrightarrow{\sim} \Rev(\BS(IK_\bullet)):r_d\otimes \cdots \otimes r_0 \mapsto r_0 \otimes \cdots \otimes r_d$$ in $\Hom_{\SBSBim}(I_0, I_d)$, for $r_j \in R^{I_j}$, $0 \leq j \leq d$.
\end{lemma}

The following lemma uses Frobenius structures to relate Bott-Samelson bimodules with their duals.

\begin{lemma}\label{lemma self duality of Bott-Samelsons}
    Fix a singular multistep expression $IK_\bullet$ as in \cref{eq singular multistep expression supplementary}. 
    The map 
    $$\eta_{IK_\bullet}: \BS(IK_{\bullet}) \to \Hom^\bullet_{R^{I_d}}(\BS(IK_{\bullet}), R^{I_d})(2(\ell(I_0)-\ell(I_d)))=D(\BS(IK_\bullet))$$ 
    given by
    \begin{equation}
        r_0 \otimes \dots \otimes r_d \mapsto \left(s_0\otimes \cdots \otimes s_d \mapsto \dd^{I_{d-1}}_{K_d}\left(\cdots \dd^{I_1}_{K_2}(\dd^{I_0}_{K_1}(r_0s_0)r_1s_1) \cdots r_{d-1}s_{d-1} \right)r_ds_d \right) \label{eq self duality of bott samelsons}
    \end{equation}
   for $r_j, s_j \in R^{I_j}$, $0 \leq j \leq d$, is an isomorphism in $\Hom_{\SBSBim}(I_d,I_0)$.
   Hence singular Bott-Samelson bimodules are self-dual (canonically, after choosing Frobenius structures). 
\end{lemma}
\begin{proof}
    The given map is a homomorphism of graded $(R^{I_d}, R^{I_0})$-bimodules, so we just have to show that the given map is an isomorphism of sets. We suppress the grading shifts and grading for notational simplicity.
    
    We proceed by induction on $d$. For $d=0$, the result is trivial.
    Now suppose the map 
    \begin{gather}
      \notag \eta': R^{I_0}\otimes_{R^{K_1}} \cdots \otimes_{R^{K_{d-1}}}R^{I_{d-1}} \rightarrow  \Hom^\bullet_{R^{I_{d-1}}}(R^{I_0}\otimes_{R^{K_1}} \cdots \otimes_{R^{K_{d-1}}}R^{I_{d-1}}, R^{I_{d-1}})\\
      \notag  r_0 \otimes \cdots \otimes r_d 
      \mapsto (s_0\otimes \cdots \otimes s_d \mapsto \dd^{I_{d-2}}_{K_{d-1}}\left(\cdots \dd^{I_1}_{K_2}(\dd^{I_0}_{K_1}(r_0s_0)r_1s_1) \cdots r_{d-2}s_{d-2} \right)r_{d-1}s_{d-1})
    \end{gather}
    is an isomorphism.
    Let us denote $ R^{I_0}\otimes_{R^{K_1}} R^{I_1}\otimes_{R^{K_2}} \cdots \otimes_{R^{K_{d-1}}}R^{I_{d-1}}$ by $M$.
     The map that we'd like to show is an isomorphism is then
    \begin{gather*}
       \eta_{IK_\bullet}:M\otimes_{R^{K_d}}R^{I_d} \rightarrow \Hom^\bullet_{R^{I_d}}(M \otimes_R^{K_d}R^{I_d}, R^{I_d}) \\
       m \otimes r \mapsto \left(m' \otimes s \mapsto \dd^{I_{d-1}}_{K_d}\left((\eta'(m))(m')\right)r s\right).
    \end{gather*}
The map $\eta_{IK_\bullet}$ may be factored as

\[
\begin{tikzcd}
   M \otimes_{R^{K_d}} R^{I_d} \arrow[r, "\eta_{IK_\bullet}"] \arrow[d, "\eta'\otimes \id"] &\Hom^\bullet_{R^{I_d}}(M \otimes_{R^{K_d}} R^{I_d}, R^{I_d}) \\
   \Hom^\bullet_{R^{I_{d-1}}}(M, R^{I_{d-1}})\otimes_{R^{K_d}} R^{I_d} \arrow[r, "\upsilon \otimes \id"] & \Hom^\bullet_{R^{K_d}}(M ,R^{K_d})\otimes_{R^{K_d}} R^{I_d} \arrow[u, "\beta"],
\end{tikzcd}\]
    where
    \begin{align*}
        \upsilon: f\mapsto  (m \mapsto \dd^{I_{d-1}}_{K_d}(f(m))) \qquad &f \in \Hom^\bullet_{R^{I_{d-1}}}(M,R^{I_{d-1}}),\; m \in M\\
        \beta: g \otimes r \mapsto (m'\otimes s\mapsto g(m')rs) \qquad&g \in \Hom^\bullet_{R^{K_d}}(M,R^{K_d}), \; m' \in M,\;\; r,s \in R^{I_d}.
    \end{align*}
    Both $\upsilon$ and $\beta$ are isomorphisms since $M$ is a finitely generated projective (in fact, free) graded module over the rings $R^{I_d}$ and $R^{K_d}$; $\eta'$ is an isomorphism by induction hypothesis.
    It follows that $\eta_{IK_\bullet}$ is an isomorphism as desired.
\end{proof}

Recall that we had the 2-functor $\Phi: \mc{D}^{\pre} \rightarrow \SBSBim$ which sends diagrams to corresponding morphisms between bimodules.
\begin{theorem}\label{thm phi commutes with dual and rev}
    The 2-functor $\Phi$ intertwines the 2-functors $\Rev$ and $\Dual$ on $\mc{D}^{\pre}$ with those on $\SBSBim$. 
    More precisely, there are isomorphisms of 2-functors 
    \begin{equation}
        \Phi \circ \Rev \xrightarrow{\sim} \Rev \circ \Phi, 
        \qquad
        \Phi \circ \Dual \xrightarrow{\sim} \Dual \circ \Phi \label{eq phi commutes with rev and dual}
    \end{equation} 
    given by the data
    \begin{gather}
    \etadot_{I_\bullet}:\Phi \circ \Rev (I_\bullet)=\BS(\Rev(I_\bullet)) \xrightarrow{\sim} \Rev(\BS(I_\bullet))= \Rev\circ \Phi(I_\bullet)
    \\
    \eta_{I_\bullet}:\Phi \circ \Dual (I_\bullet)=\BS(I_\bullet) \xrightarrow{\sim} \Dual(\BS(I_\bullet))=\Dual \circ \Phi (I_\bullet)
    \end{gather}
    for singlestep expressions $I_\bullet$, where $\etadot$ and $\eta$  are as in  \cref{lemma reversal on bott samelsons and multistep expressions} and \cref{lemma self duality of Bott-Samelsons}.
\end{theorem}
\begin{proof}
    Checking \cref{eq phi commutes with rev and dual} involves checking the compatibility of $\eta$ and $\etadot$ with the coherence 2-morphisms in \cref{def duality symmetry D algebraic} and \cref{def reversal symmetry R algebraic} respectively at the 1-morphism level, as well as checking certain equalities on (generating) 2-morphisms.
    Both are routine checks, and follows from the explicit formulas for the coherence data in \cref{def reversal symmetry R algebraic}, \cref{def duality symmetry D algebraic}, and the formulas for the isomorphisms in \cref{lemma self duality of Bott-Samelsons} and \cref{lemma reversal on bott samelsons and multistep expressions}.
    
    Here's an example to illustrate the check on the generating 2-morphisms. To check \cref{eq phi commutes with rev and dual} for $\Dual$ on a clockwise cap with region labels $I$ and $Ii$, we need to show that the diagram
    \begin{equation*}
    \begin{tikzcd}
        R^I \arrow[r, "\eta_{[I]}"] \arrow[d, "\Delta"] & \Dual(R^I) \arrow[d, "\Dual(m)"] \\
        R^I \otimes_{R^{Ii}} R^I(\ell(Ii)-\ell(I)) \arrow[r, "\eta_{[I\subset Ii \supset I]}"] &\Dual\left(R^I\otimes_{R^{Ii}}R^I(\ell(Ii)-\ell(I))\right)
    \end{tikzcd}            
    \end{equation*}
    commutes. 
    This is straightforward to check; we have
    \begin{equation*}
        \Dual(m)\circ \eta_{[I]}: r \mapsto (s\mapsto sr) \mapsto \bigl( s_1 \otimes s_2 \mapsto s_1 s_2 r \bigr)
    \end{equation*}
    which agrees with
    \begin{equation}
      \eta_{[I \subset Ii \supset I]}\circ \Delta: r\mapsto \sum_i e_i \otimes f_i r\mapsto \left(s_1 \otimes s_2 \mapsto \sum_i\dd^I_{Ii}(e_is_1)f_i rs_2 = s_1 s_2 r\right)  
    \end{equation}
    where $\{e_i\}, \{f_i\}$ are dual bases for $R^{Ii}\hookrightarrow R^I$ (so that $\sum_i\dd^{I}_{Ii}(se_i)f_i=s$ for any $s \in R^I$).
\end{proof}

\begin{defn}\label{def rotation by 180 degrees}
    Let $\rotation$ be the 2-functor on $\mc{D}^{\pre}$, given by rotating a diagram by 180 degrees. This is a strict 2-functor which coincides with the composition $\Rev \circ \Dual=\Dual \circ \Rev$.
    On $\SBSBim$, we define $\rotation$ to be the (weak) 2-functor given by the composition $\Dual \circ \Rev$.
\end{defn}

The following then easily follows from \cref{thm phi commutes with dual and rev}.
\begin{corollary}\label{corollary phi intertwines rotation}
    The 2-functor $\Phi: \mc{D}^{\pre} \rightarrow \SBSBim$ intertwines $\rotation$ with the 2-functor $\rotation= \Dual \circ \Rev$ (as well as with the 2-functor $\Rev \circ \Dual$) on $\SBSBim$. 
\end{corollary}


\begin{remark}\label{rmk rotation and flips and adjunctions}
In $\mc{D}^{\pre}$, the rotation $2$-functor $\rotation$ is more than just the composition of two flip symmetries; it is also realized by composition with caps and cups. That is, $\rotation$ agrees with the $2$-functor arising from the biadjunction of $[I \subset Ii]$ and $[Ii \supset I]$. After applying $\Phi$ one obtains the $2$-functor arising from the bijadunction of induction and restriction, with units and counits of adjunction determined by the Frobenius extension structures. In particular, this corollary implies that biadjunction on $\SBSBim$ agrees\footnote{Technically, this argument only implies agreement on the image of the functor $\Phi$. Recall that $\Phi$ was proven to be $2$-full for balanced realizations in \cite{EKLP}, and is only conjectured to be so for unbalanced realizations. The result does hold more generally though. Using the identifications $\vartheta_{IK_\bullet}:=\eta_{\Rev(IK_\bullet)} \circ \etadot_{IK_\bullet}:\BS(\Rev(IK_\bullet))\xrightarrow{\sim} \Dual \circ \Rev(\BS(IK_\bullet))$ 
for multistep expressions $IK_\bullet$, one can directly check that the resulting isomorphisms 
 \begin{gather*}
 \Hom_{\SBSBim}(\BS(IK_\bullet), \BS(IK'_\bullet)) \xrightarrow{\sim} \Hom_{\SBSBim}(\Rev (\BS(IK'_\bullet)), \Rev(\BS(IK_\bullet)))\\
 f \mapsto \vartheta_{IK_\bullet}^{-1}\circ \rotation(f) \circ \vartheta_{IK'_\bullet}
 \end{gather*}
  coincide precisely with the isomorphisms obtained by iteratively using the induction-restriction bi-adjunctions.} coherently with the composition of two other symmetries $\Dual$ and $\Rev$, which is not a priori obvious. 
\end{remark}


\begin{remark}\label{rmk left duals and right duals are naturally isomorphic}
    One may also define a \emph{left duality} 2-functor $\widetilde{\Dual}$ (see \cite[Remark 3.10]{BBE}) on $\SBSBim$ which fixes objects and act on 1-morphism categories via 
    \PTv2{Update from v1: removed incorrect grading shift.}
    \begin{equation}
    \Hom(I,J)^{\op}\to \Hom(I,J): M \mapsto \Hom^\bullet_{R^I}(M, R^I)
    \end{equation}
    with suitable changes to the coherence 2-morphisms in \cref{eq coherence data for duality} and \cref{eq coherence data for duality}.
    This 2-functor is related to $\Dual$ by the equality (not just an isomorphism) of 2-functors
    \begin{equation}\label{eq left duality is conjugate to right duality by reverse}
        \tilde{D}=\Rev \circ \Dual \circ \Rev.
    \end{equation}
    Analogous to \cref{lemma self duality of Bott-Samelsons}, we have an isomorphism (canonical, after choosing Frobenius structures)
    \begin{gather}
        \notag \tilde{\eta}_{IK_\bullet}: \BS(IK_\bullet) \xrightarrow{\sim} \widetilde{\Dual}(\BS(IK_\bullet)) \\
        r_0\otimes \cdots r_d \mapsto \Bigl(s_0 \otimes \cdots \otimes s_d \mapsto s_0r_0 \;\dd^{I_1}_{K_1}\bigl(r_1s_1\cdots \dd^{I_{d-1}}_{K_{d-1}}(r_{d-1}s_{d-1}\dd^{I_d}_{K_d}(r_d s_d)) \cdots \bigr)\Bigr) \label{eq self duality of Bott Samelsons for left dual}
    \end{gather}
    for any singular multistep expression $IK_\bullet$.
    One can then check that the isomorphisms 
    \begin{equation}
    \eta_{IK_\bullet} \circ (\eta_{IK_\bullet})^{-1}: \Dual(\BS(IK_\bullet)) \rightarrow \widetilde{\Dual}(\BS(IK_\bullet))\end{equation}
    provide the data of an isomorphism of (weak) 2-functors $\Dual \Rightarrow \widetilde{\Dual}$.
   Together with the equality in \cref{eq left duality is conjugate to right duality by reverse} and \cref{lemma dual and rev are involutions}, this also gives us an isomorphism of (weak) 2-functors $\Rev \circ \Dual \Rightarrow \Dual \circ \Rev$ (cf. \cref{corollary phi intertwines rotation}). 
\end{remark}

We now introduce certain non-degenerate bilinear forms on singular Bott-Samelson bimodules to compare with other approaches to duality in the literature.
Let $IK_\bullet$ be a singular multistep expression as in \cref{eq singular multistep expression supplementary}.
The self duality isomorphisms of singular Bott-Samelson bimodules (\cref{lemma self duality of Bott-Samelsons}) can be used to give a (right) $R^{I_d}$-bilinear form on $\BS(IK_\bullet)$ valued in $R^{I_d}$, using the canonical pairing
\begin{equation}
    \BS(IK_\bullet) \times \Dual(\BS(IK_\bullet)) \rightarrow R^{I_d}: (b, \beta)\mapsto \beta(b).
\end{equation}

\begin{defn}\label{defn invariant pairing on singular Bott Samelsons}
    Define the \emph{right (global) intersection form} on $\BS(IK_\bullet)$ to be the the (right) $R^{I_d}$-bilinear form given by 
    \begin{equation}
       \langle -,-\rangle_{IK_\bullet}: \BS(IK_\bullet)\times \BS(IK_\bullet) \rightarrow R^{I_d}: \langle b,b'\rangle_{IK_\bullet} := \bigl(\eta_{IK_\bullet}(b')\bigr)(b).
    \end{equation}
    From the explicit formula for $\eta_{IK\bullet}$ in \cref{eq self duality of bott samelsons}, 
    we have that when 
    $$b= r_0 \otimes \cdots \otimes r_d,\; b'=s_0\otimes \cdots \otimes s_d$$
    for $r_j, s_j \in R^{I_j}$, $0\leq j \leq d$,
    \begin{equation}
        \bigl(\eta_{IK_\bullet}(b')\bigr)(b) =\dd^{I_{d-1}}_{K_d} \left(\cdots \dd^{I_1}_{K_2}(\dd^{I_0}_{K_1}(r_0s_0)r_1s_1) \cdots r_{d-1}s_{d-1} \right)r_ds_d= \bigl(\eta_{IK_\bullet}(b)\bigr)(b')
    \end{equation}
    so that $\langle-,-\rangle_{IK_\bullet}$ is a symmetric bilinear form. 
    Moreover, it is $(R^{I_0}, R^{I_d})$-invariant, i.e., for any $f\in R^{I_0}$ and $g \in R^{I_d}$, we have
    \begin{equation}
        \langle fmg,m' \rangle_{IK_\bullet}=\langle m, fm'g\rangle_{IK_\bullet}.
    \end{equation}
    This form is also non-degenerate, since $\eta_{IK_\bullet}$ is an isomorphism and $\BS(IK_\bullet)$ is a free right $R^{I_d}$-module.
\end{defn}

\begin{remark}\label{rmk left intersection form}
    One can similarly also define a \emph{left (global) intersection form} $\widetilde{\langle -, -\rangle}_{IK_\bullet}$ on $\BS(IK_\bullet)$, which will be a symmetric $(R^{I_0},R^{I_d})$-invariant non-degenerate (left) $R^{I_0}$-bilinear form valued in $R^{I_0}$, using the self duality isomorphisms \cref{eq self duality of Bott Samelsons for left dual} for the left dual $\widetilde{D}$ in \cref{rmk left duals and right duals are naturally isomorphic}. 
    This form is related to the right intersection form by the equation
    \begin{equation}
        \widetilde{\langle f, g \rangle}_{IK_\bullet}= \langle \etadot_{IK_\bullet}^{-1}(f), \etadot_{IK_\bullet}^{-1}(g)\rangle_{\Rev (IK_\bullet)}
    \end{equation}
    where $\Rev(IK_\bullet)$ is the reverse of $IK_\bullet$ as in \cref{defn reverse of singular multistep expression} and $\etadot$ is as in \cref{lemma reversal on bott samelsons and multistep expressions}.
\end{remark}

\begin{remark}\label{rmk intersection forms in literature}
In Williamson's thesis \cite{WillThesis} where singular Soergel bimodules are originally defined, it is shown that the functor of tensoring with an induction or restriction bimodule commutes with duality, and the natural isomorphism involved is defined using the Frobenius structure, precisely as  $\eta_{I_\bullet}$ for the cases when $I_\bullet$ equals $[I \subset Ii]$ or $[Ii \supset I]$.  This happens in \cite[\S 3.2]{WillThesis} and \cite[Proposition 4.3.4]{WillThesis}, though the details are omitted from the published article \cite{WillSingular}. Iterating this natural isomorphism over iterated tensor products, one obtains $\eta_{I_\bullet}$ for any singlestep expression $I_\bullet$. This approach is also taken in \cite{BBE}.

Meanwhile, the subsequent literature tended to define the isomorphism between a Bott-Samelson bimodule and its dual using the global intersection form. References include \cite{EWHodge} and \cite{EMTW} for ordinary Soergel bimodules, and \cite[Appendix A]{patimo2021basesintersectioncohomologygrassmannian} for one-sided-singular Soergel bimodules, though we are unaware of a reference which treats global intersection forms for singular Soergel bimodules in general. These papers define the intersection form inductively, having chosen a particular top-degree element of each induction bimodule. In \cite[Appendix A]{patimo2021basesintersectioncohomologygrassmannian}, the Frobenius structure is explicitly used to choose a top-degree element.



\end{remark}

\begin{remark}
    While we explicitly defined the 2-functors $\Dual$ and $\Rev$ only on $\SBSBim$, it is straightforward to extend these 2-functors, as well as the results involving them, to the Karoubi envelope $\SSBim$.
\end{remark}

\subsection{Light ladders and light leaves: an example} \label{subsec:LLappendix}

As we illustrate a particular example, we briefly recall the purpose of elementary light leaves and ladders.

We consider representations of $\slf_n$, and let $V_k := \Lambda^k(\C^n)$. Let $\mu$ be a weight in $V_k$ for some $k$. There is a minimal dominant weight $\lambda_{\start}$ such that $\lambda_{\finish} := \lambda_{\start} + \mu$ is also dominant. If $L_{\lambda}$ denotes the simple representation of highest weight $\lambda$, then $L_{\lambda_{\finish}}$ is a direct summand of $L_{\lambda_{\start}} \ot V_k$, and we are interested in constructing the projection map.

Now $L_{\lambda_{\start}}$ is a summand of a tensor product $V_{i_1} \ot \cdots \ot V_{i_d}$ of fundamental representations, when $\sum_{j=1}^d \varpi_{i_j} = \lambda_{\start}$. Similarly, $L_{\lambda_{\finish}}$ is a summand of $V_{j_1} \ot \cdots \ot V_{j_e}$. The \emph{elementary light ladder} $\ELL(\mu)$ is a particular morphism
\[ \ELL(\mu) \colon (V_{i_1} \ot \cdots \ot V_{i_d}) \ot V_k \to (V_{j_1} \ot \cdots \ot V_{j_e}) \]
which remains nonzero after pre- and post-composition with the appropriate idempotents, so that it descends to a nonzero (projection) map
\[ L_{\lambda_{\start}} \ot V_k \to L_{\lambda_{\finish}}.\]

\begin{remark} \label{rmk:whenlambdanotminimal} If $\lambda$ is any dominant weight such that $\lambda + \mu$ is dominant, then $\lambda - \lambda_{\start}$ is also dominant. The projection map $L_{\lambda} \ot V_k \to L_{\lambda + \mu}$ can be obtained from $\id \ot \ELL(\mu)$ by pre- and post-composition with the appropriate idempotents. Here $\id$ represents the identity map of a tensor product of fundamental representations associated with $\lambda - \lambda_{\start}$. \end{remark}

Here is a particular example when $n \ge 5$ and $k=2$. Let $\mu = \varpi_4 - \varpi_3 + \varpi_2 - \varpi_1$. Then $\lambda_{\start} = \varpi_1 + \varpi_3$, and $\lambda_{\finish} = \varpi_2 + \varpi_4$. Note that $L_{\lambda_{\start}}$ is a summand of $V_1 \ot V_3$ and $L_{\lambda_{\finish}}$ is a summand of $V_2 \ot V_4$. Then 
\begin{equation} \ELL(\mu) \colon (V_1 \ot V_3) \ot V_2 \to (V_2 \ot V_4), \qquad \ELL(\mu) = \vcenter{\xy (0,0)*{\def\svgscale{0.15}\input{arxiv-figures/ELLforBen2_svg-tex.eps_tex}} \endxy}. \end{equation}
We have chosen, arbitrarily, to label the rightmost region with the label $0 \in \Omega$. After applying $\GS$ or $\GS_{\zeta}$ we obtain
\begin{equation} \label{eq:GSELLdiagram} \GS(\ELL(\mu)) = \textrm{scalar} \cdot {
\labellist
\small\hair 2pt
 \pinlabel {$\blk{5}$} [ ] at 8 -5
 \pinlabel {$\blk{6}$} [ ] at 24 -5
 \pinlabel {$\blk{2}$} [ ] at 40 -5
 \pinlabel {$\blk{5}$} [ ] at 56 -5
 \pinlabel {$\blk{0}$} [ ] at 72 -5
 \pinlabel {$\blk{2}$} [ ] at 88 -5
 \pinlabel {$\blk{4}$} [ ] at 8 80
 \pinlabel {$\blk{6}$} [ ] at 24 80
 \pinlabel {$\blk{0}$} [ ] at 72 80
 \pinlabel {$\blk{4}$} [ ] at 88 80
 \pinlabel {$\dgr{0}$} [ ] at 90 45
\endlabellist
\centering
\ig{1}{NastyELLimage}
}. \end{equation}

\vspace{.1cm}

We note that the colors $1, 3 \in \Omega$ remain present in every region of this diagram.

Now we consider elementary light leaves for singular Soergel bimodules. Whenever $I$ and $J \cup s$ are finitary, for $s \notin J$, and whenever $p \in W_I \backslash W / W_J$ and $q \in W_I \backslash W / W_{Js}$ are such that $p \subset q$, we call $[p \subset q]$ or $[q \supset p]$ a \emph{coset pair}. Observe that $B_q$ is a summand of $B_p \ot \BS([J \subset Js])$, and $B_p$ is a summand of $B_q \ot \BS([Js \supset J])$. We are interested in constructing the projection maps.

Choosing a reduced expressions $IK_{\bullet}$ for $p$, we have that $B_p$ is a summand of $\BS(IK_{\bullet})$. Similarly for a reduced expression $IK'_{\bullet}$ of $q$. The elementary light ladder $\ELL([p \subset q])$ is a particular morphism
\[ \ELL([p \subset q]) \colon \BS(IK_{\bullet}) \ot \BS(J \subset Js) \to \BS(IK'_{\bullet})\]
which descends to a nonzero projection map $B_p \ot \BS(J \subset Js) \to B_q$. Similarly, 
\[ \ELL([q \supset p]) \colon \BS(IK'_{\bullet}) \ot \BS(Js \supset J) \to \BS(IK_{\bullet}) \] descends to a projection map $B_q \ot \BS(Js \supset J) \to B_p$.

Soon enough we will wish to illustrate double cosets, for which it helps to fix a particular value of $n$, and we choose $n=8$. Below we will be using the symbol $n$ to indicate a particular double coset rather than the size of $\Omega$, following the notation in \cite{EKLP}. 

Just as in Remark \ref{rmk:whenlambdanotminimal}, $[p \subset q]$ might not be ``minimal.'' In this case there is a coset $z$ (think of $z$ as the analogue of $\lambda - \lambda_{\start}$) and a \emph{Grassmannian coset pair} $[m \subset n]$ such that $z.m = p$ and $z.n = q$ (i.e. a reduced expression for $z$ concatenated with a reduced expression for $m$ is a reduced expression for $p$). One has
\begin{equation} \ELL([p \subset q]) := \id_z \ot \ELL([m \subset n]), \end{equation}
where by $\id_z$ we mean the identity map of a reduced expression for $z$.

The analogue of $(-)\ot V_k$ after applying $\GS$ will be $(-) \ot [\hh{k} \supset \hh{0,k} \subset \hh{0}] = (-) \ot [-0+k]$. This involves two steps, not one, so we are interested in a composition of two elementary light leaves. Starting with a coset $q \in W_I \backslash W / W_{\hh{k}}$ we would find $[q \supset p \subset q']$ where
\[ p \in W_I \backslash W / W_{\hh{0,k}}, \qquad q' \in W_I \backslash W / W_{\hh{0}}\]
and take the composition $\ELL([p \subset q']) \circ (\ELL([q \supset p]) \ot \id_{[+k]})$.

Let 
\[ q_{\finish} = \psi_0(\lambda_{\finish}) \in W_{\hh{6}} \backslash W / W_{\hh{0}}, \qquad \text{resp. } \; q_{\start} = \psi_2(\lambda_{\start}) \in W_{\hh{6}} \backslash W / W_{\hh{2}}.\]
Following the bijection of \cref{def:psia} we have
\begin{equation} \overline{q}_{\finish} w_{\hh{0}} = {
\labellist
\tiny\hair 2pt
 \pinlabel {$1$} [ ] at 120 0
 \pinlabel {$6$} [ ] at 160 45
\endlabellist
\centering
\ig{1}{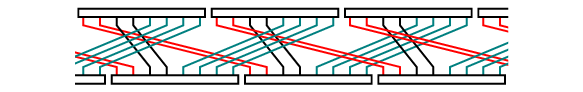}
}, \end{equation}
\begin{equation} \overline{q}_{\start} w_{\hh{2}} = {
\labellist
\tiny\hair 2pt
 \pinlabel {$1$} [ ] at 120 0
 \pinlabel {$6$} [ ] at 160 45
\endlabellist
\centering
\ig{1}{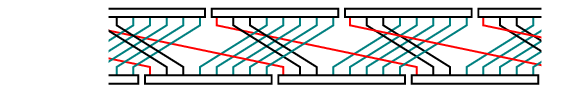}
}. \end{equation}
Let us also provide the minimal coset representatives:
\begin{equation} \label{eq:qendmin} \underline{q}_{\finish} = {
\labellist
\tiny\hair 2pt
 \pinlabel {$1$} [ ] at 120 0
 \pinlabel {$6$} [ ] at 160 45
\endlabellist
\centering
\ig{1}{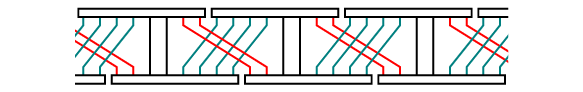}
}, \end{equation}
\begin{equation} \underline{q}_{\start} = {
\labellist
\tiny\hair 2pt
 \pinlabel {$1$} [ ] at 120 0
 \pinlabel {$6$} [ ] at 160 45
\endlabellist
\centering
\ig{1}{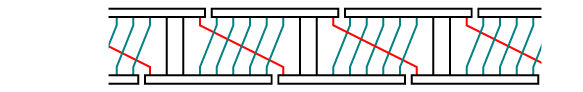}
}. \end{equation}
(Recall that $\overline{q}$ is the maximal coset representative and $\underline{q}$ is the minimal coset representative, in the Bruhat order.)

There is a unique double coset $p \in W_{\hh{6}} \backslash W / W_{\hh{0,2}}$ for which $p \subset q_{\start} \cap q_{\finish}$, and we have
\begin{equation} \underline{p} = {
\labellist
\tiny\hair 2pt
 \pinlabel {$1$} [ ] at 120 0
 \pinlabel {$6$} [ ] at 160 45
\endlabellist
\centering
\ig{1}{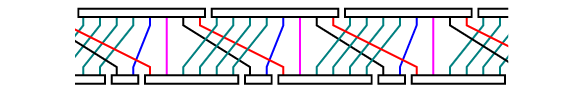}
}.\end{equation}

We wish to construct $\ELL([p \subset q_{\finish}])$, following the algorithm in \cite[\S 7]{EKLP}.

The right redundancy of $q_{\finish}$ is $Q = Q_{\finish} = \hh{024}$. Then $z = z_{\finish} \in W_{\hh{6}} \backslash W / W_Q$ satisfies $\underline{z} = \underline{q}_{\finish}$. A reduced expression for $z$ (in two different notations) is
\begin{equation} z \expr [[\hh{6} \supset \hh{246} \subset \hh{24} \supset \hh{024}]] = [\hh{6} - 2 - 4 + 6 - 0]. \end{equation}
One should think that $[\hh{6} - 2 - 4]$ breaks the parabolic subgroup into blocks associated to the different colors in \eqref{eq:qendmin}, and $[+6 - 0]$ is the cabled crossing of the red and blue strands.

Continuing, we have $m = m_{\finish} \in W_Q \backslash W / W_{\hh{02}}$ with $\underline{q}_{\finish} \cdot \underline{m} = \underline{p}$. Here is the picture:
\begin{equation} \underline{m} = {
\labellist
\tiny\hair 2pt
 \pinlabel {$1$} [ ] at 120 0
 \pinlabel {$6$} [ ] at 160 45
\endlabellist
\centering
\ig{1}{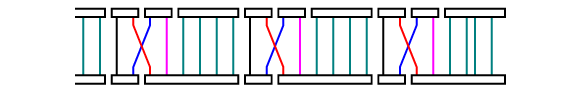}
}. \end{equation}
The reader should be able to see that stacking $\underline{q}_{\finish}$ above $\underline{m}_{\finish}$ yields $\underline{p}$. A reduced expression for $m$ is
\begin{equation} m \expr [Q - 1 - 3 + 2 - 2 + 1 + 3 + 4]. \end{equation}
The portion $[+2 - 2]$ of this expression is responsible for the crossing.
Now $n = n_{\finish} \in W_Q \backslash W / W_{\hh{0}}$ has reduced expression
\begin{equation} n \expr [Q + 2 + 4]. \end{equation}
The pair $[m \subset n]$ is a Grassmannian coset pair for $I = \hh{6}$ and $J = \hh{02}$ and $s = s_2$.

Tensoring $\ELL([m \subset n])$ (a ``sinister'' diagram with lots of left crossings and leftwards caps) with $\id_z$ yields
\begin{equation} \ELL([p \subset q_{\finish}]) = {
\labellist
\small\hair 2pt
 \pinlabel {$\blk{4}$} [ ] at 8 -4
 \pinlabel {$\blk{2}$} [ ] at 24 -4
 \pinlabel {$\blk{6}$} [ ] at 40 -4
 \pinlabel {$\blk{0}$} [ ] at 56 -4
 \pinlabel {$\blk{1}$} [ ] at 72 -4
 \pinlabel {$\blk{3}$} [ ] at 88 -4
 \pinlabel {$\blk{2}$} [ ] at 104 -4
 \pinlabel {$\blk{2}$} [ ] at 120 -4
 \pinlabel {$\blk{3}$} [ ] at 136 -4
 \pinlabel {$\blk{1}$} [ ] at 152 -4
 \pinlabel {$\blk{4}$} [ ] at 168 -4
 \pinlabel {$\blk{2}$} [ ] at 184 -4
 \pinlabel {$\blk{4}$} [ ] at 8 71
 \pinlabel {$\blk{2}$} [ ] at 24 71
 \pinlabel {$\blk{6}$} [ ] at 40 71
 \pinlabel {$\blk{0}$} [ ] at 56 71
 \pinlabel {$\blk{2}$} [ ] at 72 71
 \pinlabel {$\blk{4}$} [ ] at 120 71
 \pinlabel {$\dgr{0}$} [ ] at 184 45
\endlabellist
\centering
\ig{1}{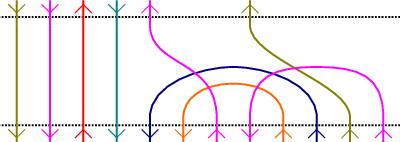}
}. \end{equation}

\vspace{.1cm}

Forgive us as we reuse the notation $Q, z, m, n$ for a new case. The right redundancy of $q_{\start}$ is $Q = Q_{\start} = \hh{235}$. Then $z = z_{\start} \in W_{\hh{6}} \backslash W / W_Q$ satisfies $\underline{z} = \underline{q}_{\start}$. A reduced expression for $z$ is
\begin{equation} z \expr [[\hh{6} \supset \hh{3,5,6} \subset \hh{35} \supset \hh{235}]] = [\hh{6} - 5 - 3 + 6 - 2]. \end{equation}
Again, $[\hh{6} - 5 - 3]$ breaks the parabolic subgroup into blocks associated to the different colors in the picture for $\underline{q}_{\start}$, and $[+6 - 2]$ is the cabled crossing of the red and blue strands.

Continuing, we have $m = m_{\start} \in W_Q \backslash W / W_{\hh{02}}$ with $\underline{q}_{\start} \cdot \underline{m} = \underline{p}$. Here is the picture:
\begin{equation} \underline{m} = {
\labellist
\tiny\hair 2pt
 \pinlabel {$1$} [ ] at 120 0
 \pinlabel {$6$} [ ] at 160 45
\endlabellist
\centering
\ig{1}{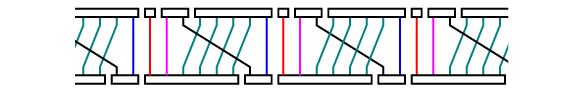}
}. \end{equation}
The reader should be able to see that stacking $\underline{q}_{\start}$ above $\underline{m}_{\start}$ yields $\underline{p}$. A reduced expression for $m$ is
\begin{equation} m \expr [[Q - 1 - 4 + 5 - 0 + 1 + 3 + 4]]. \end{equation}
The portion $[+5 - 0]$ of this expression is responsible for the cabled crossing.
Now $n = n_{\start} \in W_Q \backslash W / W_{\hh{2}}$ has reduced expression
\begin{equation} n \expr [Q + 3 + 5]. \end{equation}
The pair $[m \subset n]$ is a Grassmannian coset pair for $I = \hh{6}$ and $J = \hh{02}$ and $s = s_0$.

One obtains $\ELL([n \supset m])$ by flipping $\ELL([m \subset n])$ upside-down and using adjunction. Tensoring $\ELL([n \supset m])$ with $\id_z$ yields
\begin{equation} \ELL([q_{\start} \supset p]) = {
\labellist
\small\hair 2pt
 \pinlabel {$\blk{5}$} [ ] at 8 -2
 \pinlabel {$\blk{3}$} [ ] at 24 -2
 \pinlabel {$\blk{6}$} [ ] at 40 -2
 \pinlabel {$\blk{2}$} [ ] at 56 -2
 \pinlabel {$\blk{5}$} [ ] at 72 -2
 \pinlabel {$\blk{3}$} [ ] at 88 -2
 \pinlabel {$\blk{0}$} [ ] at 168 -2
 \pinlabel {$\blk{5}$} [ ] at 8 71
 \pinlabel {$\blk{3}$} [ ] at 24 71
 \pinlabel {$\blk{6}$} [ ] at 40 71
 \pinlabel {$\blk{2}$} [ ] at 56 71
 \pinlabel {$\blk{1}$} [ ] at 72 71
 \pinlabel {$\blk{4}$} [ ] at 88 71
 \pinlabel {$\blk{5}$} [ ] at 104 71
 \pinlabel {$\blk{0}$} [ ] at 120 71
 \pinlabel {$\blk{3}$} [ ] at 136 71
 \pinlabel {$\blk{4}$} [ ] at 152 71
 \pinlabel {$\blk{1}$} [ ] at 168 71
\pinlabel {$\dgr{02}$} [ ] at 168 35
\endlabellist
\centering
\ig{1}{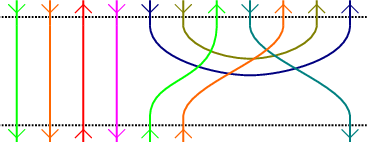}
}. \end{equation}

Now we wish to compose $\ELL([q_{\finish} \supset p])$ with $\ELL([p \subset q_{\start}]) \otimes \id_{[+2]}$ to obtain a morphism comparable with $\GS(\ELL(\mu))$ from \eqref{eq:GSELLdiagram}. However, the sources and targets of these morphisms do not match, because they use different reduced expressions for $q_{\start}$ and $q_{\finish}$ and $p$. These different reduced expressions were essential aspects of the algorithm, since at various stages we were required to use reduced expressions starting with $z_{\start}$ or $z_{\finish}$. Different reduced expressions for the same double coset are related by the (singular) braid relations, introduced in \cite{EKo}, and to each braid relation we have a corresponding diagrammatic morphism called a rex move, introduced in \cite[\S 6]{EKLP}.
To compose these elementary light leaves we should insert rex moves in between.

Examples of rex moves are: crossings of upward-oriented strands (the up-up move), crossings of downward oriented strands (down-down), sideways crossings where the colors commute in the ambient parabolic subgroup (commuting switchback), and the type $A$ switchback maps which look like a cup or cap straddling two other strands, see \cite[Example 6.1]{EKLP}. See also \cref{rmk:switchback}.

A rex move from the target of $\ELL([p \subset q_{\finish}])$ to the target of $\GS(\ELL(\mu))$ is a single type $A$ switchback
\begin{equation} \rex_3 = \ig{1}{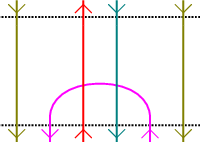}. \end{equation}
As a reminder:
\[ \ELL([p \subset q_{\finish}]) = \ig{1}{ELLend}. \]
A rex move from the target of $\ELL([q_{\start} \supset p]) \ot \id_{[+2]}$ to the source of $\ELL([p \subset q_{\finish}])$ is 
\begin{equation} \rex_2 = \ig{1}{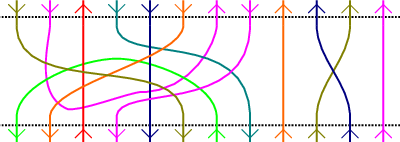}. \end{equation}
This is a composition of all the rex moves mentioned elsewhere.
As a reminder:
\[ \ELL([q_{\start} \supset p]) \ot \id_{[+2]} = 
\ig{1}{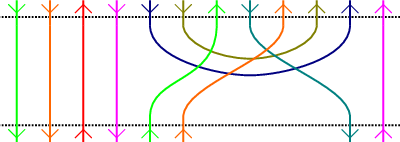}
. \]
A rex move from the source of $\GS(\ELL(\mu))$ to the source of $\ELL([q_{\start} \supset p]) \ot \id_{[+2]}$ is
\begin{equation} \rex_1 = \ig{1}{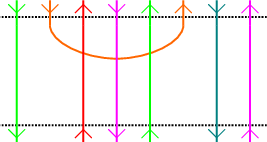}. \end{equation}
The orange cup straddling the red and fuschia strands is a type $A$ switchback, and the orange-green crossing is up-up.

\begin{remark} \label{rmk:switchback} Let $\Dual$ flip diagrams upside-down and reverse orientations without changing region labels. The composition $\rex_3 \circ \Dual(\rex_3)$ is equal to the identity map. This is called the \emph{switchback relation} and is proven in \cite[Theorem 6.7]{EKLP} as a consequence of the relations from \cite{EWSFrob}. A prototypical depiction of this relation is \cite[Lemma 3.4]{BBE}. \end{remark}

Composing all these morphisms together, we obtain
\begin{equation} \label{eq:finalnastydiagram} \rex_3 \circ \ELL([p \subset q_{\finish}]) \circ \rex_2 \circ (\ELL([q_{\start} \supset p]) \ot \id_{[+2]}) \circ \rex_1 = \ig{1.2}{NastyLL}. \end{equation}
From here, one can use a series of the relations: oriented R3 \eqref{R3oriented}, distant R3 \eqref{distantR3}, distant R2 \eqref{R2 distant}, and switchback (see \cref{rmk:switchback}), eventually deducing that the diagrams in \eqref{eq:finalnastydiagram} and \eqref{eq:GSELLdiagram} are equal. We leave this exercise to the reader.

We note that the colors $1, 3 \in \Omega$ (navy blue and orange) do not label strands in \eqref{eq:GSELLdiagram}. They arise in \eqref{eq:finalnastydiagram} because the light leaf must factor through a reduced expression for $p$.